\documentclass[a4paper,11pt]{report}

\usepackage[dvipdfx,cmyk]{xcolor}
\definecolor{bookColor}{cmyk}{0 ,0 ,0 ,1} 
\color{bookColor} 
\usepackage{amsmath,amsthm,amssymb} 
\usepackage[hidelinks]{hyperref}
\usepackage{titlesec}
\titleformat{\section}
{\normalfont\Large\scshape\centering}{\thesection}{0.5em}{}

\titleformat{\subsection}
{\normalfont\large \scshape\centering}{\thesubsection}{0.3em}{}

\usepackage{xcolor}
\hypersetup{
	colorlinks,
	linkcolor={red},
	citecolor={blue}
}

\usepackage{graphicx}
\usepackage{caption}
\usepackage{float}
\usepackage{mathtools}
\usepackage{braket}
\usepackage{dsfont}

\usepackage{changepage}

\usepackage{enumitem}
\usepackage{geometry}

\usepackage[affil-it]{authblk}

\usepackage{tikz}
\usetikzlibrary{graphs,trees,snakes,automata,topaths,shapes,arrows}
\usetikzlibrary{positioning}

\usepackage{graphicx}
\newcommand{\downspoon}{\mathrel{\reflectbox{\rotatebox[origin=c]{270}{$\multimap$}}}}

\makeatletter
\DeclareRobustCommand\widecheck[1]{{\mathpalette\@widecheck{#1}}}
\def\@widecheck#1#2{%
	\setbox\z@\hbox{\m@th$#1#2$}%
	\setbox\tw@\hbox{\m@th$#1%
		\widehat{%
			\vrule\@width\z@\@height\ht\z@
			\vrule\@height\z@\@width\wd\z@}$}%
	\dp\tw@-\ht\z@
	\@tempdima\ht\z@ \advance\@tempdima2\ht\tw@ \divide\@tempdima\thr@@
	\setbox\tw@\hbox{%
		\raise\@tempdima\hbox{\scalebox{1}[-1]{\lower\@tempdima\box
				\tw@}}}%
	{\ooalign{\box\tw@ \cr \box\z@}}}
\makeatother

\newcommand*\circled[1]{\tikz[baseline=(char.base)]{
		\node[shape=circle,draw,inner sep=2pt] (char) {#1};}}

\newcommand{\maxx}[1]{\underset{#1}{\mbox{max}}}

\newcommand{\supp}[1]{\underset{#1}{\mbox{sup}}}

\DeclareMathOperator{\Aut}{Aut}
\DeclareMathOperator{\Hom}{Hom}

\DeclareMathOperator{\Sgn}{Sgn}
\DeclareMathOperator{\sgn}{sgn}
\DeclareMathOperator{\Ind}{Ind}
\DeclareMathOperator{\Res}{Res}

\DeclareMathOperator{\Sym}{Sym}
\DeclareMathOperator{\Alt}{Alt}
\DeclareMathOperator{\id}{id}
\DeclareMathOperator{\proj}{proj}
\DeclareMathOperator{\Sl}{SL}
\DeclareMathOperator{\PSl}{PSL}
\DeclareMathOperator{\Stab}{Stab}
\DeclareMathOperator{\sphe}{Sphe}
\DeclareMathOperator{\spe}{Spe}
\DeclareMathOperator{\cusp}{Cusp}
\DeclareMathOperator{\Rep}{Rep}
\DeclareMathOperator{\Irr}{Irr}
\DeclareMathOperator{\Ker}{Ker}
\DeclareMathOperator{\Prim}{Prim}
\DeclareMathOperator{\Sp}{Sp}
\DeclareMathOperator{\cl}{cl}
\DeclareMathOperator{\Gl}{GL}

\newcommand{\Hr}[1]{\mathcal{H}_{#1}}

\newcommand{\N}{\mathbb N}
\newcommand{\Z}{\mathbb Z}
\newcommand{\Q}{\mathbb Q}
\newcommand{\R}{\mathbb R}
\newcommand{\C}{\mathbb C}

\newcommand{\tg}[1]{\textbf{#1}}
\newcommand{\ub}[1]{\overline{#1}}

\newcommand{\es}{\varnothing}

\newcommand{\norm}[2]{\lVert #1 \lVert_{#2}}

\newcommand{\modu}[1]{\lvert#1\lvert}
\newcommand{\lb}{\lbrack}
\newcommand{\rb}{\rbrack}

\newcommand{\s}[2]{\sum\limits_{#1}^{#2}}
\newcommand{\li}[2]{\xrightarrow[#1\rightarrow#2]{}}

\newcommand{\prods}[2]{\langle #1,#2\rangle}

\newcommand{\restr}[2]{{
		\left.\kern-\nulldelimiterspace 
		#1 
		\vphantom{\big|} 
		\right|_{#2} 
}}

\newcommand{\fct}[4]{\qq:\qq #1\qq\longrightarrow\qq #2\qq:\qq #3\qq \mapsto\qq #4}

\newcommand{\q}{\quad}
\newcommand{\qq}{\mbox{ }}

\newcommand{\diff}{\text{\rm d}}

\newcommand{\hr}[1]{\mathcal{H}_{#1}}

\newcommand{\Ch}{\text{\rm Ch}}

\newcommand{\Fix}{\text{\rm Fix}}

\newcommand{\chiup}{\raisebox{\depth}{\(\chi\)}}

\theoremstyle{plain}
\newtheorem{theorem}{Theorem}
\newtheorem{proposition}[theorem]{Proposition}
\newtheorem{lemma}[theorem]{Lemma}

\newtheorem{corollary}[theorem]{Corollary}

\newtheorem{theoremletter}{Theorem}

\theoremstyle{definition}
\newtheorem{definition}[theorem]{Definition}
\newtheorem{remark}[theorem]{Remark}

\newtheorem{example}[theorem]{Example}
\newtheorem{nonexample}[theorem]{Non-example}

\theoremstyle{plain}
\newtheorem*{theorem*}{Theorem}
\newtheorem*{conjecture*}{Conjecture}
\newtheorem*{proposition*}{Proposition}
\newtheorem*{lemma*}{Lemma}
\newtheorem*{corollary*}{Corollary}
\newtheorem*{problem*}{Problem}

\theoremstyle{definition}
\newtheorem*{definition*}{Definition}
\newtheorem*{remark*}{Remark}
\newtheorem*{example*}{Example}

\numberwithin{theorem}{section}
\numberwithin{equation}{section}

\author{Lancelot Semal}
\date{April 2023}

\begin{document}
	\pagenumbering{roman}
	\newgeometry{left=0cm,bottom=0cm, top=0cm}
	\begin{figure}
			\includegraphics[scale=1]{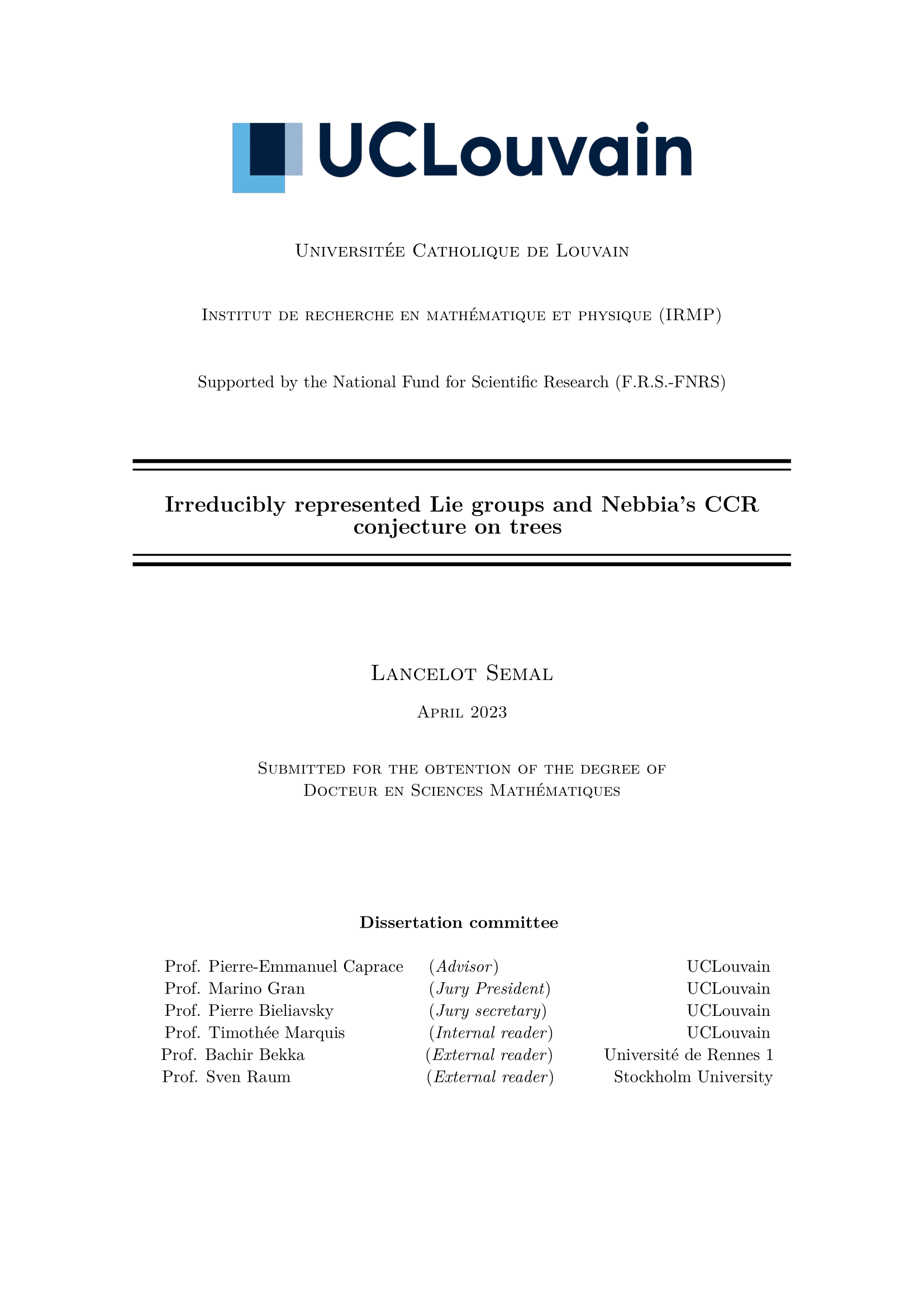}
	\end{figure}
	\restoregeometry
	\newpage
	\thispagestyle{empty}
	\mbox{}
	\newpage
	\section*{\textbf{\Large Remerciements}}
	\mbox{ }
	\newline
	Il est temps pour moi de conclure un des chapitres de ma vie. Mais avant de poser ma plume, de tourner la page et de voguer vers de nouveaux horizons, je me dois de remercier certaines personnes. 
	
	À bien des égards, le doctorat s'apparente à un voyage en mer. Vous voilà explorateur. Voguez maintenant. Visez l'horizon. Partez étendre les frontières de la connaissance. Je n'aurais pu espérer un meilleur guide que Pierre-Emmanuel. Il me semble clair que cette thèse n'aurait pu exister sans lui. Au-delà de ses indéniables qualités de mathématicien, j'aimerais le remercier pour le cadre de travail empreint de bienveillance dans lequel il m'a permis de m'épanouir pendant ces quelques années. Je ne peux que le remercier pour la façon avec laquelle il m'a partagé ses idées, écouté chacune des miennes, pour son humilité permanente et pour m'avoir toujours encouragé à faire les choses à ma manière. Tu as été et tu seras encore longtemps pour moi un modèle.
	
	J'aimerais également remercier les membres de mon jury pour le temps qu'ils m'ont consacré. Je remercie également Sven pour tout l'interêt qu'il a porté à mes travaux, Bachir pour ses ouvrages inspirants et Timothée pour son implication minutieuse dans la relecture de ce manuscrit. Je remercie bien entendu mes frères de galère. Le mousse Justin pour sa joie de vivre communicative, le flibustier Philippe pour nos exquises pauses cafés et le timonier Alex pour la passion qu'il partage au quotidien. Je remercie Cathy pour ses petites attentions, Alexandre, Abel et Benoît pour nos moments de franche camaraderie et Sandrine pour avoir été un peu comme une seconde maman. 
 
 	Ensuite et parce que la vie ne se réduit pas au monde du travail, j'aimerais remercier tous les gens qui ont empli mon quotidien d'amour et de légèreté. Je remercie mes parents pour leur soutien sans faille, pour tout l'amour qu'il m'ont donné et pour m'avoir encouragé à faire de ma vie une plaine de jeux. Je remercie Alix pour toute la douceur qu'elle a mise dans ma vie. Et je remercie Marie pour plus de raisons que je ne peux l'exprimer. 
	 
	Enfin, j'aimerais remercier une dernière personne à qui je dédie entièrement cet ouvrage. Pour m'avoir inspiré depuis l'enfance jusqu'à la dernière ligne de ce manuscrit et pour tous les moments que nous avons partagés. Simplement, merci ...
	
	\vfill
	{\begin{flushright}\textit{Pour mon grand Dada}	\end{flushright}}
	\newpage
	\tableofcontents

	\chapter*{Introduction}
	\addcontentsline{toc}{chapter}{Introduction}
	\pagenumbering{arabic}
	\setcounter{page}{1}

	The notion of symmetry is fundamental in mathematics and plays a major role in many physical models of nature. In modern mathematics, this concept is best apprehended by the abstract algebraic notion of a \tg{group} whose purpose is to encode only the rules of composition of a set of symmetries. As a result, a group is not attached with a specific object anymore. It thus becomes natural to study the actions of groups on various structures. The theory of representations aims at studying their linear actions on vector spaces.  
	The theory of linear representations of finite groups was founded by Frobenius in the late 1890's and further developed by Burnside and Schur among others. In this framework, a \tg{representation} of a finite group $G$ is a group homomorphism $\pi : G\rightarrow \Gl(V)$ to the group $\Gl(V)$ of invertible linear transformations of a finite dimensional complex vector space $V$ and a representation is said to be \tg{irreducible} if it does not have any proper invariant subspace.
	Furthermore, as $G$ is finite, the vector space $V$ can be equipped with a structure of Hilbert space by an averaging argument so that each of the $\pi(g)$ is a unitary operator (a linear map preserving the inner product). In particular, every linear representation of a finite group $G$ is equivalent to a \tg{unitary representation}. In addition, it can be shown that the set $\widehat{G}$ of equivalence classes of irreducible representations of $G$ is finite and that each representation $\pi$ of $G$ admits an essentially unique decomposition as a direct sum:
	$$\pi \simeq \bigoplus_{\sigma \in \widehat{G}} n_{\sigma} \sigma$$
	where $n_\sigma\sigma$ is a direct sum of $n_\sigma$ copies of $\sigma$. Hence, the theory of representations of finite groups reduces to the study of their irreducible representations. As a result, the theory of linear representations of finite groups became somewhat dominated for the next 30 years by the attempt to solve the following two problems:
	\begin{enumerate}[leftmargin=*, label=(\Roman*)]
		\item Given a finite group $G$ and a representation $\pi$ of $G$ determine the multiplicities $n_\sigma$ in the decomposition $\pi \simeq \bigoplus_{\sigma \in \widehat{G}} n_{\sigma}\sigma.$
		\item Given a finite group $G$, determine the dual space $\widehat{G}$, that is the set of irreducible representations of $G$ up to equivalence. 
	\end{enumerate}
	The work achieved in this direction led to the beautiful and profound theory of characters which to the knowledge of the author stays to the present days one of the most powerful tools in the study of representations of finite groups. According to A.~Kirillov \cite{Kirillov1994}, a fundamental series of paper by A.~Young devoted to the representation theory of the symmetric groups brought somewhat an end to this era of representation theory in the middle of the 1920's and, even if the subject was (and still is) far from being exhausted, the main interest moved into the representation theory of compact topological groups. In this framework, a representation of a compact topological group $G$ is a continuous group homomorphisms $\pi: G\rightarrow \mathcal{L}(\mathcal{H})$ to the group of continuous linear invertible operators of a (possibly infinite dimensional) complex vector space $\mathcal{H}$ and a representation is said to be \tg{irreducible} if it does not have any closed proper invariant subspace. One of the most fundamental achievements of this second period of representation theory is the proof by A.~Haar and J.~von Neumann of the existence of a finite invariant Radon measure on any compact topological group. In fact, locally compact groups are characterised among topological groups by the property of admitting an invariant Radon measure (which is finite precisely when the group is compact). As a result, averaging arguments become available on any compact group and most of Frobenius' theory extends to compact topological groups. To be more precise, every representation of a compact group is equivalent to a unitary representation in the sense that $\mathcal{H}$ can be equipped with a structure of complex Hilbert space so that $\pi$ ranges inside the unitary group $\mathcal{U}(\mathcal{H})$. Furthermore, it can be shown that:
	\begin{enumerate}[leftmargin=*, label=(\roman*)]
		\item Every representation $\pi$ of a compact group $G$ admits an essentially unique decomposition as a direct sum $$\pi \simeq \bigoplus_{\sigma \in \widehat{G}} n_{\sigma} \sigma$$ with $n_\sigma\in \N\cup \{\infty\};$
		\item Every irreducible representation of a compact group $G$ is finite dimensional and appears as a subrepresentation of the regular representation of $G$ with multiplicity equal to its dimension.
	\end{enumerate}
	In particular, just as for finite groups the set of equivalence classes of irreducible representations of a first-countable compact group is at most countable and the theory of representations of compact groups reduces nicely to the study of their irreducible representations. 
	Contemporarily to these discoveries, the necessity to consider non-compact groups and their infinite dimensional representations became clear in quantum physics at the beginning of the 1940's. As a striking example, the irreducible unitary representations of the Lorentz group are known to play a substantial role in the classification of elementary particles. On the other hand, the transition to the framework of locally compact groups needed a substantial overhaul of the theory and of its basic tools. For instance, the role played by purely algebraic constructions needed to be replaced by arguments relying on topology, measure theory and functional analysis. In this general framework where $G$ is a locally compact group without additional restrictions, the word ``representation'' usually refers to a continuous unitary representation and a \tg{unitary representation} is defined as a continuous group homomorphism $\pi: G\rightarrow \mathcal{U}(\mathcal{H})$ ranging in the unitary group $\mathcal{U}(\mathcal{H})$ of a (possibly infinite dimensional) complex Hilbert space $\mathcal{H}$. According to G.W.~Mackey \cite{Mackey1963}, who has been one of the most influential contributors of representation theory of locally compact groups, this domain of mathematics stayed somewhat dominated by the attempts to answer the following two questions:
	\begin{enumerate}[leftmargin=*, label=(\Roman*)]
		\item To what extent can an arbitrary unitary representation of a locally compact group $G$ be decomposed into irreducible representations ? 
		\item To what extent is it possible to determine all the irreducible unitary representations of a locally compact group $G$ ? 
	\end{enumerate}
	Simultaneously to the efforts made to solve these problems, Murray and von Neumann were laying the foundations of the tightly related theory of operator algebras which to the present days stands on its own as a fundamental discipline of  mathematics. The tools developed in this domain and the deep investigations of Segal, Mautner, Mackey, Fell, Dixmier and Glimm (to cite just a few) led to the conclusion that every unitary representation of a locally compact group $G$ can be decomposed into irreducible representations. However, by contrast with compact groups, certain representations do not even contain an irreducible representation and hence can not decompose as a direct sum of irreducible representations. In fact, every unitary representation can be decomposed as a direct integral of irreducible representations. On the other hand, this decomposition need not be (and is in general far from being) unique. In fact, the uniqueness feature is a characteristic property of the so-called \textbf{Type} \textbf{\rm I} locally compact groups (Bernstein and Kirillov have used the term ``tame'' to qualify Type {\rm I} groups, and ``wild'' for those which are not of Type {\rm I}). According to Auslander and Moore \cite{AuslanderMoore1966}, one of the fundamental questions one can raise about any locally compact group is that of determining when it is of Type {\rm I}. Loosely speaking, Type {\rm I} groups are characterised by the property that each of their factor representations (the general conceptual analogue of representations consisting of a unique isotypical component) contains (and hence is a multiple of) an irreducible representation. On the other hand, every unitary representation of a locally compact group admits a canonical decomposition as a direct integral of factor representations so that each unitary representation of a Type {\rm I} group admits an essentially unique decomposition as a direct integral of irreducible representations. In a sense, this shows that factor representations are somewhat the appropriate building blocks of unitary representations in the general case and that the theory reduces nicely to irreducible representations when the group is of Type {\rm I}. 
	
	Concerning the second question, the determination of the unitary dual of a locally compact group (its set of equivalence classes of irreducible unitary representations) is known to be an intractable problem unless the group is of Type {\rm I}. This puts a second emphasis on the importance of Type {\rm I} groups. This follows from a result of Glimm ensuring that the unitary dual of a non-Type {\rm I} group does not carry the structure of a standard Borel space. The first classification theorem for infinite-dimensional representations of non-compact groups were published in 1947 for groups such as the Lorentz group, its restricted version $\Sl_2(\R)$ and the group of affine transformations of the real line.
	In parallel, it became an important task to determine which locally compact groups are of Type {\rm I} without necessarily being able to identify their unitary duals. For instance, compact groups and abelian groups were known to be of Type {\rm I}. In 1953, Harish-Chandra proved that all semisimple Lie groups are of Type {\rm I}. A few years later, an alternative proof was found by Godement and then generalised by Bernstein to the setup of $p$-adic groups. In addition, Dixmier proved that real algebraic groups are of Type {\rm I} and Takenouchi proved that all exponential Lie groups are of Type {\rm I}. However, the property of being of Type {\rm I} should be seen more as an exception than the rule. For instance a fundamental result of Thoma shows that a countable discrete group is of Type {\rm I} if and only if it is virtually abelian, so that the algebraic structure of Type {\rm I} discrete groups is drastically restricted. 
	
	\subsection*{Key contributions}
	To the present day, despite the fact that many of the basic tools of representation theory apply to all locally compact groups without any restriction whatsoever, most of the advanced results concern discrete groups, Lie groups and algebraic groups over local fields. This thesis falls within Mackey's vision of representation  theory with a particular interest in Type {\rm I} groups. It aims at studying the relations between the algebraic structure of locally compact groups and the properties of their unitary representations with an emphasis on non-linear groups for which the current state of knowledge remains largely elusive. The work presented in these notes focuses on two independent problems and is separated into two parts accordingly. We provide a description of these problems together with a global overview of our personal contributions just below. However, as each of the chapters exposing our results already has an introduction, we do not provide all the details here and refer the reader to these introductions for a deeper overview.
	
	\subsection*{Problem I : Irreducibly represented Lie groups} 
	The first problem treated in these notes focuses on the concept of irreducibly represented groups with an emphasis on certain families of Type {\rm I} Lie groups. We recall that a locally compact group is \tg{irreducibly represented} if it admits a faithful irreducible unitary representation. A classical problem initially formulated by W.~Burnside \cite[Appendix F]{Burnside1955} in the context of finite groups and extending naturally to the context of locally compact groups is given by the following:
	\begin{problem*}
		Find a characterisation of the locally compact groups that are irreducibly represented purely in terms of their algebraic structure.
	\end{problem*}
	For instance, it can be deduced from Schur's Lemma that a finite abelian group is irreducibly represented if and only if it is cyclic. To the present day, algebraic characterisations were obtained for finite groups \cite{Gaschutz1954} and more generally for countable discrete groups \cite{BekkadeleHarpe2008} but not beyond this framework (at least to the knowledge of the author). Quite remarkably and despite the fact that the arguments used in the finite and in the countable case are very different from one another, the algebraic characterisation is the same in both setups. To be more precise, a countable discrete group $G$ is irreducibly represented if and only if every every finite normal subgroup of $G$ contained in the subgroup generated by all the minimal non-trivial finite normal subgroups of $G$ is generated by a single $G$-conjugacy class. One of our main achievements is an algebraic characterisation of the connected nilpotent Lie groups and of the Hausdorff connected component of the identity of real linear algebraic groups that are irreducibly represented. Similarly to the case of countable groups our reasoning relies on direct integral decompositions of representations and on measure theory. However, by contrast with the case of discrete groups, our reasoning relies heavily on the fact that our groups are of Type {\rm I}. Loosely speaking our characterisation reduces to some algebraic conditions on the center of the group and on a closed topologically characteristic subgroup isomorphic as a real vector space to a finite dimensional real vector space called the \tg{linear socle}, see Definition \ref{definition linear socle}. Among other things, we prove the following:
	\begin{theorem*}[Theorem \ref{Corollary letter amenable lie groups}]
		Let $G$ be the Hausdorff connected component of the identity of the $\R$-points of a linear algebraic group defined over $\R$. Then, the following assertions are equivalent: 
		\begin{enumerate}
			\item $G$ is irreducibly represented.
			\item The following assertions hold:
			\begin{enumerate}[leftmargin=*]
				\item The center $\mathcal{Z}(G)$ of $G$ is isomorphic to the circle group or is discrete countable and does not contain a subgroup isomorphic to $C_p\times C_p$ for any prime $p$.
				\item The linear dual $\mathcal{V}(G)^*$ of the linear socle $\mathcal{V}(G)$ of $G$ is spanned as a real vector space by a single coadjoint $G$-orbit.
			\end{enumerate}
			\item There exists a faithful irreducible representation of $G$ that is weakly contained in the regular representation $\lambda_G$.
		\end{enumerate}
	\end{theorem*}

	\subsection*{Problem II : Nebbia's conjecture}
	By contrast with the first part of these notes, the second part of this thesis is entirely centred around Problem II. This second problem is motivated by a conjecture of Nebbia formulated in \cite{Nebbia1999} and focuses on a surprising parallel existing between the regularity of the unitary duals of groups of tree automorphisms and the regularity of their actions on the tree and its boundary. We recall that a locally compact group $G$ is \textbf{CCR} if the operator $\pi(f)$ is compact for each irreducible representation $\pi$ of $G$ and for all $f\in L^{1}(G)$. An important feature of CCR locally compact groups is that they are of Type {\rm I}. In fact, the unitary dual $\widehat{G}$ of any locally compact group $G$ can be equipped with a topology called the Fell topology and both the property of being Type {\rm I} and the property of being CCR can be expressed in terms of the regularity of this topology. To be more precise, a locally compact group $G$ is of Type {\rm I} (resp.~CCR) if and only if the Fell topology of $\widehat{G}$ is $T_0$-separated (resp.~$T_1$-separated). 
	In addition, we recall from the beginning of the introduction that determining which locally compact group is of Type {\rm I} is a problem of fundamental importance in the representation theory of locally compact groups. At the end of the $80$'s, the classification of the irreducible representations of the full group of automorphisms $\Aut(T)$ of any thick regular tree $T$ led to the conclusion that those groups are all CCR and hence of Type {\rm I} \cite{Olshanskii1980}. A few years later, C.~Nebbia's work highlighted a surprising parallel between the regularity of the unitary dual of groups of tree automorphisms and the regularity of their action on the boundary of the tree. More precisely, he showed that any closed unimodular CCR vertex-transitive subgroup $G\leq \Aut(T)$ must act transitively on the boundary $\partial T$. Further progress going in this direction were recently achieved by Houdayer and Raum \cite{HoudayerRaum2019} and with even higher level of generality by Caprace, Kalantar and Monod \cite{CapraceKalantarMonod2022}. Among other things, they showed that any closed non-amenable Type {\rm I} subgroup acting minimally on a locally finite tree $T$ must act $2$-transitively on the boundary $\partial T$ \cite[Corollary D]{CapraceKalantarMonod2022}. On the other hand, Nebbia conjectured at the end of the $90$'s that conversely, the regularity of the action of a group of automorphisms of a tree on the boundary should imply the regularity of its unitary dual \cite{Nebbia1999}. Without loss of generality, his conjecture can be stated as follows.
	\begin{conjecture*}[Nebbia's conjecture]
		Let $T$ be a locally finite tree all of whose vertices have degree~$\geq 3$. Then, every closed subgroup $G\leq \Aut(T)$ acting $2$-transitively on the boundary $\partial T$ is CCR. 
	\end{conjecture*}
	Since compact groups are CCR and since every closed non-compact subgroup $G\leq \Aut(T)$ acting transitively the boundary $\partial T$ automatically acts $2$-transitively on the boundary (\cite[Lemma 3.1.1]{BurgerMozes2000}), notice that the claim of the conjecture is that every closed subgroup $G\leq \Aut(T)$ acting transitively on the boundary $\partial T$ is CCR. Prior to our work, this claim was already supported by the fact that rank one semi-simple algebraic groups over local-fields \cite{Bernshtein1974} and groups satisfying the Tits independence property \cite{Amann2003} are CCR. However, tackling Nebbia's conjecture in its full generality is an extremely ambitious task. Our personal contribution to its resolution focuses on the following restricted setup: 
	\begin{problem*}
		Let $T$ be a locally finite tree all of whose vertices have degree~$\geq 6$ and let $G \leq \Aut(T)$ be a closed subgroup acting $2$-transitively on the boundary $\partial T$ and whose local action at every vertex contains the alternating group. Determine whether $G$ is CCR.
	\end{problem*}
	The groups appearing in this problem were classified and explicitly described by N.~Radu in his thesis \cite{Radu2017} so that we refer to them as Radu groups. In addition, N.~Radu has highlighted that for some semi-regular trees $T$, all the closed groups of tree automorphisms acting $2$-transitively on the boundary $\partial T$ are Radu groups. In fact, this statement holds for all $(d_0,d_1)$-semi-regular tree $T$ with $d_0,d_1\in \Theta$ where 
	\begin{equation*}
		\Theta= \{d \geq 6 : \mbox{ each }\mbox{2-transitive subgroup of }\Sym(d)\mbox{ contains Alt}(d)\}.
	\end{equation*}
	In addition, $\Theta$ is known to be an asymptotically dense subset of $\N$ in the sense that 
	$$\lim_{n\rightarrow \infty} \frac{\modu{\Theta\cap \{1,2,...,n\}}}{n}=1.$$ Our main contribution to Nebbia's conjecture is given by the following result.
	\begin{theorem*}[Theorems \ref{Theorem A} and \ref{Thoerem letter ABCDE}]
		Radu groups are CCR so that Nebbia's conjecture is confirmed for any $(d_0,d_1)$-semi-regular tree with $d_0,d_1\in \Theta$.
	\end{theorem*} 
	The proof of that theorem has two major steps. The first one has a general scope on the representation theory of totally disconnected locally compact groups and is inspired by Ol'shanskii's work, while the second one restricts its attention to Radu groups and relies heavily on Radu's classification. To be more precise, we recall the irreducible representations of any group of tree automorphisms split in three disjoint subsets. Each of them is either spherical, special or cuspidal depending on the ``size'' of its ``biggest'' isotropy group. The classification of spherical and special representations of any closed subgroup $G\leq \Aut(T)$ acting $2$-transitively on the boundary $\partial T$ is a classical result that was achieved in the $70$'s due to Godement, Cartier and Matsumoto to cite just a few. This classification led to the conclusion that both spherical and special representations range into the set of compact operators and hence are CCR. A few years later, Ol'shanskii provided a classification of the cuspidal representations of the full group of automorphisms of any regular tree by exploiting the fact that these groups satisfy the Tits independence property. It appears that this independence property can be translated as a factorisation property on a basis of neighbourhoods of the identity consisting of compact open subgroups. In Chapter~\ref{Chapter Olshanskii's factor}, we develop an abstraction of Ol'shanskii's framework allowing one to describe the irreducible representations of totally disconnected locally compact groups all of whose isotropy groups are ``small'' provided that the group admits a basis of neighbourhoods of the identity consisting of compact open subgroups satisfying the same kind of factorisation property. For instance, Theorem \ref{la version paki du theorem de classification} provides a description of these irreducible representations and ensures that they are all square-integrable. We apply this machinery on certain families of groups of automorphisms of trees and right-angled buildings in Chapter~\ref{Chapter application olsh facto}. To be more precise, Theorems \ref{thm B1} and \ref{THM B} lead to a description of some of the cuspidal representations of groups of tree automorphisms satisfying the Property \ref{IPk} as defined in \cite{BanksElderWillis2015} and Theorem \ref{Theorem E} leads to a description of certain irreducible representations of universal groups of particular right-angled buildings as defined in \cite{Universal2018}. We apply this machinery to Radu groups in Chapter \ref{Chapter Radu groups} and obtain a description of the cuspidal representations of any simple Radu groups via Theorem \ref{theorem C}. Since these are all square-integrable, we deduce from Radu's classification that each Radu group is CCR, providing our contribution to Nebbia's conjecture. Finally, we describe the Fell topology of the known part of the unitary dual of closed groups of automorphism of trees acting $2$-transitively on the boundary of the tree in Chapter \ref{Chapter Fell topology of Aut(T)}. For instance Theorems \ref{theorem la topo de fell de Sphe spe G cas transitif}, \ref{theorem la topo de fell de Sphe spe G cas 2orbits} and \ref{thm C} ensure that the unitary dual of any Radu group is homeomorphic to the following topological space (where the left hand part of the picture represents the interval with double endpoints; the dotted lines link the non-Hausdorff points):
	\begin{figure}[H]
		\includegraphics[scale=0.15]{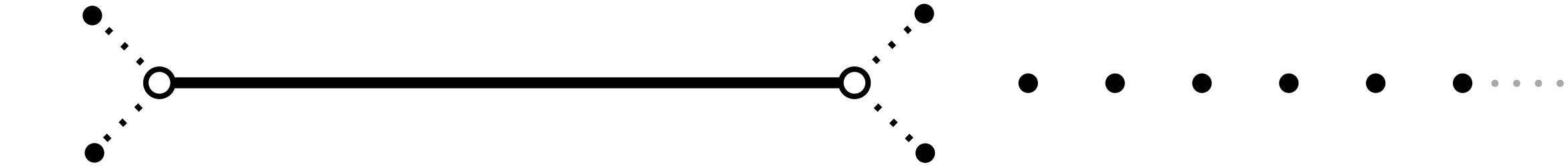}
	\end{figure}
	\subsection*{Structure of the thesis}
	This thesis is divided into two parts. The first part consists of Chapters \ref{Chapter I preliminaries} and \ref{Chapter Lie groups are irred faithfull} and contains the overall preliminaries required for the global understanding of the other chapters as well as our personal contributions to the representation theory of Lie groups. The aim of Chapter \ref{Chapter I preliminaries} is to provide a toolbox of notions that will be used throughout the text as well as a proper introduction to the theory of unitary representations of locally compact groups to anyone with a basic background in topology, measure theory and functional analysis. We expose our personal contribution to the representation theory of Lie groups and treat Problem I focusing on the concept of irreducibly represented groups in Chapter \ref{Chapter Lie groups are irred faithfull}. By contrast, the second part consisting of Chapters \ref{Preliminaries on totally disconnected locally compact groups and representations}, \ref{Chapter Olshanskii's factor}, \ref{Chapter application olsh facto}, \ref{Chapter Radu groups} and \ref{Chapter Fell topology of Aut(T)} focuses on the representation theory of totally disconnected locally compact groups and is entirely centred around Problem II. This second part should be seen as the core of the thesis and contains the vast majority of our personal contributions. The aim of Chapter \ref{Preliminaries on totally disconnected locally compact groups and representations} is to provide an introduction to the representation theory of totally disconnected locally compact groups and to the classification of the irreducible representations of the full groups of automorphisms of any semi-regular tree. In addition, this chapter serves as a toolbox for the advanced concepts that will be used throughout the rest of the notes. We expose our personal contributions on the representation theory of totally disconnected locally compact groups in Chapters \ref{Chapter Olshanskii's factor}, \ref{Chapter application olsh facto}, \ref{Chapter Radu groups} and \ref{Chapter Fell topology of Aut(T)}. The content of Chapter \ref{Chapter Olshanskii's factor} has a general scope and provides an axiomatic framework on totally disconnected locally compact groups under which we are able to describe some of their irreducible representations. We provide concrete applications of this machinery to certain families of groups of automorphisms of trees and of right-angled buildings in Chapter \ref{Chapter application olsh facto}. These two chapters are part of a preprint that was accepted for publication \cite{Semal2021O}. In Chapter \ref{Chapter Radu groups}, we treat the Problem II and provide a classification of the irreducible representations of Radu groups based on the machinery developed in Chapter \ref{Chapter Olshanskii's factor}. This part of our work is also the object of a preprint available on ArXiv \cite{SemalR2022}. Finally, we describe the Fell topology of the known part of the unitary dual of closed groups of automorphism of trees acting $2$-transitively on the boundary of the tree in Chapter \ref{Chapter Fell topology of Aut(T)}. The chapter dependencies are described by the following diagram:
	\begin{figure}[H]\label{drawingchapterrepart}
		\begin{center}
			\includegraphics[scale=0.15]{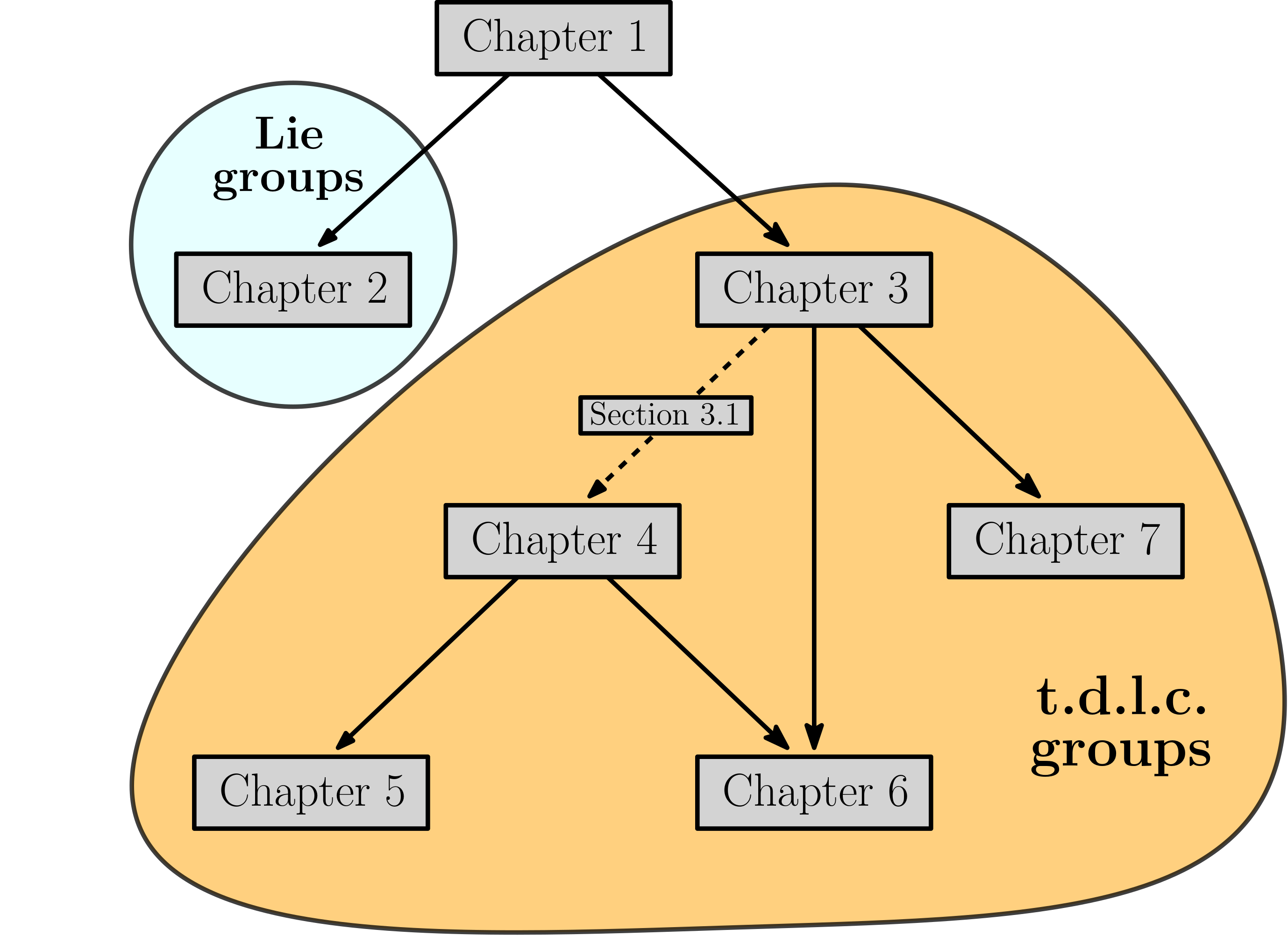}
		\end{center}
	\end{figure}
	
	\newpage
	\thispagestyle{empty}
	\mbox{}
	\newpage
	
	\part{General theory and irreducibly represented Lie groups}
	\newpage
	\thispagestyle{empty}
	\mbox{}
	\newpage
	\chapter{Representations of locally compact groups}\label{Chapter I preliminaries}
	The purpose of this chapter is to provide a general overview of the theory of unitary representations of locally compact groups and its basic tools. Most of the content presented here is based on \cite{Tao2014}, \cite{Folland2016} and \cite{BekkadelaHarpe2020} to which we refer any reader interested in a deeper understanding of these mathematical concepts. 
	\section{Locally compact groups}\label{Section locally comapct groups}
	
	The purpose of this section is to provide an introduction to locally compact groups with a special emphasis on Lie groups and totally disconnected locally compact groups. 

	\begin{definition}
		A \tg{topological group} $G$ is a group equipped with a topology such that both the multiplication map $$G\times G \rightarrow G:(g,h)\mapsto gh$$ and the inverse map $$G\rightarrow G:g\mapsto g^{-1}$$ are continuous. We denote the neutral element of $G$ by $1_G$.
	\end{definition}  
	A topological space $X$ is \tg{locally compact} if every element admits a compact neighbourhood, that is, is contained in an open set with compact closure.
	\begin{definition}
		A \tg{locally compact group} is a topological group whose topology is $T_2$ (Hausdorff) and locally compact.
	\end{definition}
	\begin{remark}
		In these notes, we will only consider groups whose topologies admit countable bases of open sets. A topological group with this property is said to be \tg{second-countable}. Loosely speaking, this hypothesis ensures that the group is not too large and that its topology can be apprehended by the convergence and divergence of sequences. 
	\end{remark} 
	\begin{remark}
		A second-countable locally compact group $G$ is $\sigma$-compact, meaning that it can be written as a union of countably many compact sets. Indeed, by second-countability there exists a countable basis of open subsets $\mathcal{B}$. On the other hand, since $G$ is locally compact, each element $g\in G$ is contained in an open neighbourhood $V_g$ with compact closure $\overline{V_g}$.	By definition of basis of a topology, there exists $B_g\in \mathcal{B}$ such that $g\in B_g\subseteq V_g$. and hence $\overline{B_g}\subseteq \overline{V_g}$. It follows that $\overline{B_g}$ is compact. Applying this for each point $g\in G$ and since $\mathcal{B}$ is countable we obtain a countable cover of $G$ by compact subsets. 
	\end{remark}
	\begin{remark}
		Some authors do not require that the topology of a locally compact group is Hausdorff. However, if $G$ is a topological group whose topology is locally compact, the closure of the identity $\overline{\{1_G\}}=H$ defines a closed normal subgroup of $G$ and the quotient group $G/H$ (equipped with the quotient topology) is a  topological group whose topology is both $T_2$ and locally compact, see \cite[Exercise 4.1.2]{Tao2014}. 
	\end{remark}
	We now provide examples of locally compact groups.
	\begin{example}
		Every group with the discrete topology is locally compact. Such a group is second-countable if and only if it is countable. 
	\end{example}
	\begin{example}
		Let $\mathcal{A}$ be a Banach space. The additive group $(\mathcal{A},+)$ equipped with the norm topology is a topological group. Notice that this group is  locally compact if and only if every closed ball of $\mathcal{A}$ is compact. It follows that $(\mathcal{A},+)$ is a locally compact group if and only if $\dim(\mathcal{A})<\infty$. 
	\end{example}
 	In this thesis, we are interested into two somewhat opposite families of locally compact groups, namely real Lie groups and totally disconnected locally compact groups. In the first part of these notes most of our original results concern mostly real Lie groups and in the second part of these notes our results concern solely totally disconnected locally compact groups. We now provide a definition of both families and explore the differences between them. The first family consists of those locally compact groups that admit a smooth structure. We recall that a $d$-dimensional \tg{smooth manifold} is a topological space $M$ admitting a smooth atlas; that is a family  $\{\phi_\alpha :\alpha\in \Omega\}$ of homeomorphisms $\phi_\alpha : U_\alpha \rightarrow V_\alpha$ from open subsets $U_\alpha$ of $M$ to open subsets $V_\alpha$ of $\R^d$, such that $\bigcup_{\alpha\in\Omega} U_\alpha$ is an open cover of $M$ and such that the map $\phi_\beta\circ \phi_\alpha^{-1}$ is smooth on its domain $\phi_\alpha(U_\alpha\cap U_\beta)$ for all $\alpha,\beta\in \Omega$. 
 	\begin{definition}
 		A \tg{real Lie group} is a group equipped with a structure of smooth manifold for which both group operations are smooth. 
 	\end{definition}
	
 	One of the most common interpretation of  Hilbert's fifth problem asks whether any locally Euclidean group admits a compatible smooth structure for which it is a real Lie group and more generally what sort of obstructions prevents a locally compact group from being a Lie group. This problem was initially solved by D.~Montgomery, L.~Zippin \cite{MontgomeryZippin1952} and A.~Gleason \cite{Gleason1952} for locally Euclidean groups and lead to the so called Gleason-Yamabe's theorem which greatly clarified the structure of general locally compact groups. We now state the theorem and refer the reader to \cite[Section 6]{Tao2014} for a detailed exposition.
 	\begin{theorem*}[Gleason-Yamabe's theorem {\cite[Theorem 6.0.11]{Tao2014}}]\label{Gleason Yamabe theorem }
 		Let $G$ be a locally compact group. For any open neighbourhood $U$ of the identity in the connected component of the identity $G^\circ$ of $G$, there is a compact normal subgroup $K$ of $G^\circ$ contained in $U$ for which $G^\circ/K$ is isomorphic to a Lie group.  
 	\end{theorem*}
 	In particular, one recovers a Lie group from any locally compact group after ignoring the large scale (by restricting to an open subgroup) and the small scale (by allowing a quotient by a compact normal subgroup). However, the Lie groups are far from being the only family of locally compact groups. In fact, the Lie groups are characterised among the locally compact groups by the property of not having small subgroups. We now formalise this concept and clarify our statement.
 	\begin{definition}
 		A topological group $G$ has \tg{no small subgroups} if there exists a neighbourhood of the identity which does not contain any subgroup other than the trivial group $\{1_G\}$. 
 	\end{definition}
 	 \begin{theorem}[{\cite[Corollary 5.3.3]{Tao2014}}]\label{Theorem Lie groups are the locally comapct groups without small subgroups}
 		A locally compact group $G$ has no small subgroups if and only if it is isomorphic to a Lie group. 
 	\end{theorem}
  	We provide a proof of one of the implications for instructive purposes.
 	\begin{lemma}
 		Let $G$ be a Lie group, then $G$ has no small subgroups.
 	\end{lemma}
 	\begin{proof}
 		If $G$ is discrete the statement is trivial so we might as well suppose that $\dim(G)\geq1$. Replacing $G$ by its connected component of the identity (which is an open subgroup of $G$) we might as well suppose that $G$ is connected. Now, let $\mathfrak{g}$ be the Lie algebra of $G$ and consider a neighbourhood $V$ of $0$ in $\mathfrak{g}$ and a neighbourhood $U$ of the identity $1_G$ in $G$ such that the exponential map 
 		$$\exp_G: V\rightarrow U$$ 
 		is a diffeomorphism. Shrinking $U$ and $V$ further if needed we suppose that $V$ is an open ball $B(0,\varepsilon)$ of radius $\varepsilon$ and that the Campbell-Baker-Hausdorff formula holds on the product of any two element of $U$. Since the exponential map is a local diffeomorphism the set $U'=\exp_G(B(0,\frac{\epsilon}{2}))$ is an open neighbourhood of the identity in $G$. We claim that $U'$ does not contain any subgroup of $G$ other than the trivial group. To see this, suppose that $H\subseteq U'$ is a subgroup and let $h\in H-\{1_G\}$. Since the exponential map is injective on $V$ there exists a non zero vector $v\in B(0,\frac{\epsilon}{2})-\{0\}$ such that $h=\exp_G(v)$. On the other hand, for all $n\in \N$ such that $nv\in V$, the Campbell-Baker-Hausdorff formula ensures that 
 		$$\exp_G(nv)=\exp_G(v)^n=h^n\in H.$$
 		As $v\in B(0,\frac{\epsilon}{2})-\{0\}$, there exists $n\in \N$ such that $nv\in B(0,\epsilon)-B(0,\frac{\epsilon}{2})$. Since $\exp_G:V\rightarrow U$ is bijective, one obtains a contradiction with the inclusion $H\subseteq  \exp_G(B(0,\frac{\epsilon}{2}))$.
 	\end{proof} 
	In the second part of this thesis we will be interested in a family of locally compact groups with somewhat opposite behaviour : the family of totally disconnected locally compact groups. As we will see in Chapter \ref{Preliminaries on totally disconnected locally compact groups and representations}, these groups are characterised among the locally compact groups by the property of admitting a basis of neighbourhoods of the identity consisting of compact open subgroups. We now provide a precise definition. We recall that a (non-empty) topological space is \tg{connected} if the only clopen subsets are the whole space and the empty set. At the other extreme, a topological space is \tg{totally disconnected} if the only connected subsets are the singleton sets. 
	\begin{definition}
		A \tg{totally disconnected locally compact group} is a topological group whose topology is both locally compact and totally disconnected.
	\end{definition}
	\begin{remark}
		Discrete groups are both Lie groups and totally disconnected locally compact groups. However, these are the only locally compact groups that belong to both families.
	\end{remark}
	We now provide an important family of examples by looking at the completions of the field $\Q$. We recall that for any prime $p$ and every $q\in \Q$, there is a unique decomposition $q=p^nr/d$ with $r$, $d$ relatively prime and relatively prime with respect to $p$.
	\begin{definition}
		For each prime $p$, the $p$-\tg{adic norm} $\norm{\cdot}{p}$ of $\Q$ is defined by $\norm{q}{p}=p^{-n}$ where $n$ is determined by the above decomposition.
	\end{definition}
	By contrast with the Euclidean norm, each of the $p$-adic norm satisfies a stronger version of the triangular identity called the \tg{ultrametric} inequality that is  $$\norm{q_1+q_2}{p}\leq \max\{\norm{q_1}{p},\norm{q_2}{p}\}\q \forall q_1,q_2\in \Q.$$ 
	A famous result from Ostrowski ensures that an absolute value on $\Q$ is either equivalent to the Euclidean norm or to a $p$-adic norm for some prime $p$. 
	\begin{example}
		It is well known that the completion of $\Q$ with respect to the Euclidean norm is the field of real numbers $\R$ and hence defines a Lie group. By contrast, for all prime $p$, the completion of $\Q$ with respect to the $p$-adic norm is a totally disconnected locally compact group (since the open balls for the $p$-adic norm are clopen sets). This completion is called the $p$\tg{-adic numbers} and denoted by $\Q_p$. 
	\end{example}
 	
 	We now make precise the statement that locally compact groups are characterised among topological groups by the property of admitting a regular invariant measure. 
	\begin{definition}
		Let $X$ be a Hausdorff topological space. The Borel $\sigma$-algebra $\mathcal{B}(X)$ of $X$ is the $\sigma$-algebra generated by the open subsets of $X$ and a \tg{Borel measure }on $X$ is a countably additive non-negative measure $\mu: \mathcal{B}(X)\rightarrow \lb 0,+\infty\rb$. A \tg{Radon measure} on $X$ is a Borel measure satisfying the three following properties:
		\begin{enumerate}[label=(\roman*)]
			\item (local finiteness) Every point has a neighbourhood with finite $\mu$-measure.
			\item (inner regularity) One has $\mu(B)=\sup_{K\subseteq B,\qq K\mbox{\tiny compact}}\mu(K)$ for every Borel measurable set $B\in \mathcal{B}(X)$.
			\item (outer regularity) One has $\mu(B)=\inf_{B\subseteq U,\qq U\mbox{\tiny open}}\mu(U)$ for every Borel measurable set $B\in \mathcal{B}(X)$.
		\end{enumerate}
	\end{definition}
	\begin{definition}
		Let $G$ be a $\sigma$-compact locally compact group and let $\mu$ be a Radon measure on $G$. The measure $\mu$ is \tg{left-invariant} (resp. right-invariant) if one has that $\mu(gB)=\mu(B)$ (resp. $\mu(Bg)=\mu(B)$) for all $g\in G$ and every $B\in \mathcal{B}(X)$. A \tg{left-invariant Haar measure }(resp. \tg{bi-invariant Haar measure}) on $G$ is a non-zero Radon measure that is left-invariant (resp. both left-invariant and right-invariant). 
	\end{definition}
	\begin{example}
		Let $G$ be a discrete group. Then, the counting point measure is a  bi-invariant Haar measure on $G$.
	\end{example}
	\begin{example}
		The Lesbegue measure on $\R^d$ is a bi-invariant Haar measure of the additive group $(\R^d,+)$ for the topology of the Euclidean norm.
	\end{example}
	The following classical result ensures the existence of a Haar measure for every second-countable locally compact group. 
	\begin{theorem}[{\cite[Theorem 4.1.6]{Tao2014}}]\label{Theorem exisetnce of Haar measure}
		Let $G$ be a $\sigma$-compact locally compact group. There exists  a left-invariant Haar measure $\mu$ on $G$. Furthermore, this measure is unique up to multiplication by scalars: if $\mu$, $\nu$ are two left-invariant Haar measures on $G$, $\nu=\lambda \mu$ for some scalar $\lambda >0$.
	\end{theorem} 
	The existence of a Haar measure characterises locally compact groups among topological groups. 
	\begin{proposition}[{\cite[443 E]{Fremlin2006}}]
		Let $G$ be a Hausdorff topological group and suppose that $G$ admits a left-invariant Haar measure $\mu$. Then $G$ is locally compact. 
	\end{proposition}
	\begin{proof}
		Since $\mu$ is non-zero, the inner regularity ensures the existence of a compact subset $K$ of $G$ with positive measure. Replacing $K$ by $k^{-1}K\cup K^{-1}k$ for any element $k\in K$, we can suppose without loss of generality that $K$ contains $1_G$ and is symmetric. We claim that $KK=\{kk': k,k'\in K\}$ has non-empty interior. It follows that $G$ contains an open subset with compact closure and hence is locally compact. The proof of the claim relies on an ingenious trick using the regularity of the convolution product. Since $K$ is compact, the local finiteness of $\mu$ ensures that $\mu(K)$ is finite and the convolution product $\mathds{1}_K*\mathds{1}_K:G\rightarrow \C$ defined by 
		$$\mathds{1}_K*\mathds{1}_K(g)=\int_K \mathds{1}_K(k^{-1}g)\diff \mu(k)$$
		is a continuous function \cite[Proposition 444R]{Fremlin2006} (the symmetry of $K$ is used here). On the other hand, the support of $\mathds{1}_K*\mathds{1}_K$ is contained inside $KK$ because $\mathds{1}_K(k^{-1}g)=\mathds{1}_{kK}(g)$ for all $k\in K$. As $\mathds{1}_K*\mathds{1}_K(1_G)=\mu(K)\not=0$, the continuity of $\mathds{1}_K*\mathds{1}_K$ ensures the existence of an open neighbourhood $U$ of $1_G$ such that $$U\subseteq \mbox{\rm supp}(\mathds{1}_K*\mathds{1}_K)\subseteq KK$$
		and the result follows.
	\end{proof}

	The last two statements stay true when replacing the left-invariance by right-invariance. However, most locally compact groups do not admit a bi-invariant Haar measure. In fact, the function quantifying the non-right invariance of a left-invariant Haar-measure defines a continuous group homomorphism called the modular function.  
	\begin{definition}\label{Definition modualr function}
		Let $G$ be a $\sigma$-compact locally compact group and, let $\mu$ be a left-invariant Haar measure on $G$. For all $t\in G$, one might consider the conjugate measure $\mu_t$ defined by $$\mu_t(B)=\mu(t^{-1}Bt)=\mu(Bt)\q \forall B\in \mathcal{B}(G).$$ It is not hard to see that $\mu_t$ is also a left-invariant Haar measure on $G$.  In particular, Theorem \ref{Theorem exisetnce of Haar measure} ensures the existence of a positive scalar $\Delta_G(t)$ such that $\mu_t=\Delta_G(t)\mu$. The \tg{modular function} is the corresponding map $$\Delta_G:G \rightarrow \R^*:t\mapsto\Delta_G(t)$$ 
		where $\R^*$ is the multiplicative group of $\R-\{0\}$.
		This map is a continuous group homomorphism that does not depend on the choice of $\mu$, see \cite[Chapter 12]{Robert1983}.
		The group $G$ is \tg{unimodular} if it admits a bi-invariant Haar measure or, equivalently, if the modular function of $G$ is trivial.
	\end{definition}
	\begin{remark}
		Since the measure of an open subset of $G$ with compact closure is both finite and non-zero with respect to any Haar measure on $G$, the modular function of $G$ can explicitly be computed as follows. Considering an open subset $U$ of $G$ with compact closure, one has that $\Delta_G(t)=\frac{\mu(Ut)}{\mu(U)}$. 
	\end{remark}
	We now provide a list of examples.
	\begin{example}
		Every discrete group is unimodular.
	\end{example}
	\begin{example}
		Every locally compact abelian group is unimodular.
	\end{example}
	\begin{example}
		Every locally compact group $G$ with $G=\overline{\lb G,G\rb}$ is unimodular. Indeed, since the modular function is a continuous group homomorphism, $\Delta_G(G)$ is isomorphic to an abelian quotient of $G$ and is therefore trivial. In particular, every simple locally compact group is unimodular.
	\end{example}
	\begin{example}\label{example compaclty generated is unimodular}
		Every locally compact group generated by compact subgroups is unimodular. This follows from the fact that the modular function is a continuous group homomorphism, as the image of a compact subgroup defines a compact subgroup of the multiplicative group of strictly positive real numbers $\R^*_{+}$ and is therefore trivial. 
	\end{example}
	\begin{nonexample}
		The group $G$ of invertible affine transformations of $\R$ is not unimodular. Indeed, this group can be realised as group of matrices
		$$G=\bigg\{\begin{pmatrix}
		a & b\\
		0 & 1
		\end{pmatrix}:\qq a\in \R^*_{+},\qq b\in \R\bigg\}$$ where $\R^*_{+}$ denote the multiplicative group of strictly positive real numbers. From this presentation, one can show that $\frac{1}{a^2}\diff a\diff b$ is a left-invariant Haar measure on $G$. This measure is clearly not right-invariant due to the definition of the multiplication. 
	\end{nonexample}

	\section{Unitary representations}\label{Section unitary representaions}
	The purpose of this thesis is to study the unitary representation theory of locally compact groups. That is, the continuous actions by linear isometries of those groups on Hilbert spaces. In Section \ref{Section locally comapct groups}, we have introduced the notion of locally compact groups and shown that those groups are naturally equipped with a rich analytic structure due to the existence of Haar measures. The purpose of the present section is to introduce unitary representations and related notions. Among other things, we show that the rich analytic structure of locally compact groups ensures the existence of a faithful unitary representation and of sufficiently many irreducible unitary representations to separate the points.
	
	In this thesis, we will consider only complex and separable Hilbert spaces unless specifically mentioned otherwise. We recall that a Hilbert space is separable if it admits a countable dense subset. Let $\mathcal{H}$ be a Hilbert space and denote by $\mathcal{L}(\mathcal{H})$ the involutive algebra of bounded linear operators on $\mathcal{H}$ where the involution $*$ is given by adjoint operators. We recall that the adjoint $V^*$ of an element $V\in \mathcal{L}(\mathcal{H})$ is the unique element of $\mathcal{L}(\mathcal{H})$ satisfying 
	$$\prods{V\xi}{\eta}_\mathcal{H}=\prods{\xi}{V^*\eta}_\mathcal{H} \q \forall \xi,\eta \in \mathcal{H}.$$
	We denote by  
	$$\mathcal{U}(\mathcal{H})=\{V\in \mathcal{L}(\mathcal{H}) : V^*V=VV^*=\id_{\mathcal{H}}\}$$
	the group of unitary operators on $\mathcal{H}$ and equip it with the \tg{weak topology} that is the smallest topology for which the maps $$\phi_{\xi,\eta}:\mathcal{U}(\mathcal{H})\rightarrow \C:V\mapsto \prods{V \xi}{\eta}$$ are continuous for all $\xi,\eta\in \mathcal{H}$. Now, let $G$ be a locally compact group.
	\begin{definition}
		A \tg{continuous unitary representation} of $G$ is a pair $(\pi,\Hr{\pi})$,  where $\Hr{\pi}$ is a Hilbert space called  \tg{representation space} of $\pi$ and a continuous group homomorphism $\pi: G \rightarrow \mathcal{U}(\Hr{\pi}).$ In other words we require that the functions 
		$$\varphi_{\xi,\eta}:G\rightarrow \C : g\mapsto \prods{\pi(g)\xi}{\eta}_{\Hr{\pi}}$$
		are continuous for all $\xi,\eta\in \hr{\pi}$. We also define the \tg{dimension} of $\pi$ to be the dimension of $\Hr{\pi}$. For brevity and when it leads to no confusion we will say ``\tg{representation}'' instead of ``continuous unitary representation'' and denote the pair $(\pi,\Hr{\pi})$ by $\pi$. 
	\end{definition}
	The following example shows that locally compact groups always admits a faithful unitary representation. 
	\begin{example}\label{Example permutation representation}
		Let $G$ be a locally compact group acting from the left on a $\sigma$-compact locally compact Hausdorff topological space $X$, let $\mathcal{B}(X)$ be the Borel $\sigma$-algebra generated by the open subsets of $X$, let $\nu$ be a Radon measure $\nu$ on $X$ and suppose that $G$ acts in a measure preserving way, that is $\nu(g B)=\nu(B)$ for all $B\in \mathcal{B} (X)$. As possible choices for $X$, the reader might consider a countable set with the counting point measure or $G$ itself with a left-invariant Haar measure. Now, consider the set $$L^2(X,\nu)=\bigg\{\varphi : X\rightarrow \C : \int_X\modu{\varphi(x)}^2\diff \mu(x)<\infty\bigg\}$$  of complex valued $\nu$-square integrable functions. This space is a complex Hilbert space when equipped with the inner product
		$$\prods{\varphi}{\psi}_{L^2(X,\nu)}=\int_X \prods{\varphi(x)}{\psi(x)}_\C\qq \diff \nu(x)\q \forall \varphi,\psi\in L^2(X,\nu).$$  
		The \tg{Koopman representation} $\lambda_{X,\nu}$ is defined on $L^2(X,\nu)$ as 
		$$\lambda_{X,\nu}(t)\varphi(x)=\varphi(t^{-1} x)\q \forall \varphi\in L^2(X,\nu), \qq\forall  t\in G,\qq \forall x\in X.$$
		Notice that $\lambda_{X,\nu}$ indeed ranges inside $\mathcal{U}(L^2(X,\nu))$, since 
		$$ \prods{\lambda_{X,\nu}(t) \varphi}{\lambda_{X,\nu}(t)\psi}_{L^2(X,\nu)}=\prods{\varphi}{\psi}_{L^2(X,\nu)}\q  \forall \varphi,\psi \in L^2(X,\nu), \qq  \forall t\in G.$$
		The Koopman representation $\lambda_G,\mu$ corresponding to the choice $X=G$ equipped with a left-invariant Haar measure $\mu$ is called the \tg{left-regular representation of $G$}. We usually denote this representation by $(\lambda_G,L^2(G))$ for brevity. Notice that the left-regular representation of $G$ is faithful. Indeed, let $g\in G-\{1_G\}$ and consider an open neighbourhood $U$ of $1_G$ with compact closure in $G$ such that $g^{-1}\not \in U$. Since $G$ is locally compact, the existence of such a neighbourhood is guaranteed and $\mathds{1}_U\in L^2(G,\mu)$. On the other hand, $\lambda_G(g)\mathds{1}_U(1_G)=\mathds{1}_U(g^{-1})=0$ so that $\lambda_G(g)\mathds{1}_U\not= \mathds{1}_U$. Similarly, it is possible to define the \tg{right-regular representation} $\rho_G$ of $G$ on the Hilbert space $L^2(G,\mu)$ by 
		$$\rho_G(t)\varphi(g)=\sqrt{\Delta_G(t)}\varphi(gt)\q \forall g\in G,\qq \forall t\in G,\qq \forall \varphi\in L^2(G,\mu)$$
		where $\Delta_G$ is the modular function of $G$ defined on page \pageref{Definition modualr function}. The modular function ensures that $\rho_G$ ranges inside $\mathcal{U}(L^2(G,\mu))$ since a straightforward computation shows for all $ \varphi,\psi \in L^2(G,\mu)$ and for every  $t\in G$ that
		\begin{equation*}
		\begin{split}
		\prods{\rho_G(t) \varphi}{\rho_G(t)\psi}_{L^2(G)}&=\int_G \Delta_G(t) \varphi(gt)\psi(gt) \diff \mu(g)\\
		&=\int_G \varphi(g')\psi(g') \diff \mu(g') =\prods{\varphi}{\psi}_{L^2(G)}.
		\end{split}
		\end{equation*}
	\end{example}
	\begin{example}
		Let $G$ be a locally compact group and let $\pi$ be a unitary representation of $G$. A closed subspace $\mathcal{K}$ of $\Hr{\pi}$ is called $\pi(G)$\tg{-invariant} (or simply $G$-invariant) if $\pi(g)\mathcal{K}\subseteq \mathcal{K}$ for all $g\in G$. Every such space gives rise to a \tg{subrepresentation} $G\rightarrow \mathcal{U}(\mathcal{K})$ of $\pi$ assigning to $g\in G$ the restriction of $\pi(g)$ to $\mathcal{K}$. We write $\sigma\leq \pi$ when $\sigma$ is a subrepresentation of $\pi$. 
	\end{example}
	Our next task will be to introduce the concept of sum of representations. We motive this concept with a situation arising naturally from the study of unitary representations. Let $G$ be a locally compact group, let $\pi$ be a representation of $G$ and suppose the existence of a proper closed $G$-invariant subspace $\mathcal{K}$ of $\Hr{\pi}$. In particular, this subspace gives rise to a subrepresentation $\pi_{\mathcal{K}}$ of $G$. Since $\pi$ is unitary, notice that the orthogonal complement
	$$\mathcal{K}^\perp=\{\xi\in \Hr{\pi}: \prods{\xi}{\eta}_{\Hr{\pi}}=0\qq \forall\eta \in \mathcal{K}\}$$
	is also a closed $G$-invariant subspace and hence gives rise to another subrepresentation $\pi_{\mathcal{K}^\perp}$ of $G$. Now, notice that $\Hr{\pi}$ decomposes as $\mathcal{K}\oplus \mathcal{K}^\perp$ and that $\pi_{\mathcal{K}}$ and $\pi_{\mathcal{K}^\perp}$ act as $\pi$ on their respective space. In a way, we would like to say that $\pi$ is obtained as a sum of those two subrepresentations. The following example formalises this concept.
	\begin{example}
		Let $G$ be a locally compact group, let $\mathcal{I}$ be a countable set and let $\pi_i$ be a representation of $G$ for all $i\in \mathcal{I}$. The direct sum $$\bigoplus_{i\in \mathcal{I}}\Hr{\pi_i}=\{(\xi_i)_{i\in\mathcal{I}}: \sum_{i\in \mathcal{I}} \norm{\xi_i}{\mathcal{H}_i}^2<\infty\}$$ is a separable Hilbert space with inner product defined as 
		$$\prods{(\xi_i)_{i\in\mathcal{I}}}{(\eta_i)_{i\in\mathcal{I}}}= \sum_{i\in \mathcal{I}} \prods{\xi_i}{\eta_i}_{\mathcal{H}_i}.$$
		The \tg{direct sum} $\bigoplus_{i\in \mathcal{I}}\pi_i$ is the representation of $G$ acting on $\bigoplus_{i\in \mathcal{I}}\Hr{\pi_i}$ defined by
		$$\bigg\lb\bigoplus_{i\in \mathcal{I}}\pi_i(g)\bigg\rb(\xi_i)_{i\in \mathcal{I}}=(\pi_i(g)\xi_i)_{i\in \mathcal{I}}\q \forall g\in G.$$
		For any given non-zero cardinal $\aleph\in \N\cup \aleph_0$ and any representation $\sigma$ of $G$ we denote by $\aleph\sigma$ the direct sum $\bigoplus_{n \leq  \aleph }\sigma$. 
	\end{example} 
	Another important operation on representations is provided by the tensor product. As a reminder, for every two Hilbert spaces $\mathcal{H}_1$, $\mathcal{H}_2$ the complex tensor product of the vector spaces underlying $\mathcal{H}_1$ and $\mathcal{H}_2$ can be equipped with an inner product given by the linear extension of the expressions
	$$\prods{\xi_1\otimes \xi_2}{\eta_1\otimes \eta_2}_{\mathcal{H}_1\otimes \mathcal{H}_2}=\prods{\xi_1}{\eta_1}_{\mathcal{H}_1}\prods{\xi_2}{\eta_2}_{\mathcal{H}_2}\q \forall \xi_1,\eta_1\in \mathcal{H}_1,\qq \forall \xi_2,\eta_2\in \mathcal{H}_2.$$
	The completion of this tensor product with respect to this inner-product is a Hilbert space called the \tg{tensor product} of $\mathcal{H}_1$ and $\mathcal{H}_2$. This one is denoted by $\mathcal{H}_1\otimes \mathcal{H}_2$. Another definition of the tensor product relying on Hilbert-Schmidt operators and traces can be found in {\cite[Section 3.1]{Mackey1976theory}. This construction is particularly useful in order to build representations of direct products of locally compact groups.
		\begin{example}
			Let $G$ and $H$ be two locally compact group. Every pair of representations  $\pi_1$ of $G$ and $\pi_2$ of $H$ gives rise to a representation $\pi_1\otimes \pi_2$ of the direct product $G\times H$ called the \tg{outer Kronecker product} of $\pi_1$ and $\pi_2$. Its representation space is $\Hr{\pi_1}\otimes \Hr{\pi_2}$ and the representation is defined as the continuous linear extension of the expressions
			$$\pi_1\otimes \pi_2(g,h) \xi_1\otimes \xi_2=\pi_1(g)\xi_1\otimes \pi_2(h)\xi_2\q \forall (g,h)\in G\times H, \qq \xi_1\in \Hr{\pi_1},\qq \xi_2\in \Hr{\pi_2}.$$ 
			Since $G$ embeds diagonally in $G\times G$, notice that any pair of representations $\pi_1,\pi_2$ of $G$ gives rise to a representation of $G$ by composition of the embedding and with the outer Kronecker product. This representation of $G$ is called the \tg{inner-Kronecker product} of $\pi_1$ and $\pi_2$ and is denoted in the same way by $\pi_1\otimes \pi_2$ when it leads to no confusion. 
		\end{example}
		There are several equivalence relations on the set of representations that we will consider in this thesis : the unitary equivalence, the quasi-equivalence and the weak equivalence. When it leads to no confusion we will denote an equivalence class of representations by any of its representative and when the distinction is not clear we will denote by $\lb \pi \rb$ the equivalence class of $\pi$. We start by introducing the notion of unitary equivalence and quasi-equivalence and refer to Section \ref{section Fell topology} for the notion of weak equivalence. Given two representations $\pi_1$ and $\pi_2$ of $G$, a bounded linear operator $V: \mathcal{H}_1\rightarrow \mathcal{H}_2$ is an \tg{intertwining operator} between $\pi_1$ and $\pi_2$ if $\pi_2(g)V = V\pi_1(g)$ for all $g\in G$. We denote by $\Hom_G(\pi_1,\pi_2)$ the set of intertwining operators of $G$. 
		\begin{definition}
			Two representations $\pi_1$ and $\pi_2$ of $G$ are:
			\begin{enumerate}[label=(\roman*)]
				\item \tg{unitary equivalent} (or \tg{equivalent} for short) which we denote by $\pi_1\simeq\pi_2$ if there exists an invertible unitary operator $V\in \Hom_G(\pi_1,\pi_2)$, that is, $V^*\pi_1(g)V = \pi_2(g)\qq \forall g\in G$.
				\item \tg{disjoint} which we denote by $\pi_1\downspoon \pi_2$ if there does not exist a subrepresentation of $\pi$ equivalent to a subrepresentation of $\pi_2$.
				\item \tg{quasi-equivalent} which we denote by $\pi_1 \approx\pi_2$ if there does not exist a subrepresentation of $\pi_1$ that is disjoint from $\pi_2$ and vice-versa.
			\end{enumerate} 
		\end{definition}
		\begin{remark}
			Both the relation of unitary equivalence and of quasi-equivalence are equivalence relations on the representations of $G$. Furthermore, unitary equivalent representations are automatically quasi-equivalent.
		\end{remark}
		Equipped with the notion of unitary equivalence and direct sums, it becomes natural to study the ``indecomposable" representations of $G$; that is those representations of $G$ that are not unitarily equivalent to a direct sum of subrepresentations. Those representations are the irreducible representations of $G$. When the group $G$ is compact the theory of unitary representations reduces to the study of these simple pieces. Indeed, every representation $\pi$ of a compact group $G$ admits a unique decomposition $\pi \simeq \bigoplus_{\sigma\in \widehat{G}} n_\sigma \sigma$ as a direct sum of irreducible representations. However, such a decomposition does not extend to the general case where $G$ is a locally compact group.  In fact, certain locally compact groups have representations that do not contain any irreducible subrepresentation. For instance, the left-regular representation of the additive group $(\R,+)$ does not contain any irreducible representation. For that reason, the general theory of unitary representations requires a finer notion of ``indecomposable'' pieces. This role is played by factor representations which happen to be the simple pieces provided by the relation of quasi-equivalence. 
	\begin{definition}
		Let $G$ be a locally compact group. A representation $\pi$ of $G$ is said to be:
		\begin{enumerate}[label=(\roman*)]
				\item \tg{cyclic} if if there exists a vector $\xi\in \Hr{\pi}$ such that the linear span of $\{\pi(g)\xi: g\in G\}$ is dense in $\Hr{\pi}$ ($\xi$ is said to be a \tg{cyclic} vector of $\pi$). 
				\item  \tg{irreducible} if it does not have a non-zero proper closed $G$-invariant subspace. 
				\item \tg{multiplicity free} if, for every decomposition $\pi\simeq \pi_1\oplus \pi_2$, the representations $\pi_1$ and $\pi_2$ are disjoint. 
				\item  \tg{factorial} (or a \tg{factor representation}) if $\pi$ cannot be written as a sum of two disjoint representations.
		\end{enumerate}
	\end{definition}
	In particular, notice that a representation $\pi$ is irreducible if and only if every of its non-zero vector $\xi \in \Hr{\pi}$ is cyclic and that the irreducible representations are exactly the representations which are both factorial and multiplicity free. In addition, notice that these four notions are preserved by unitary equivalence and that the notion of factor representation is even preserved by quasi-equivalence. 
	\begin{definition}
		Let $G$ be a second-countable locally compact group. We define $\widetilde{G}$ to be the set of representations of $G$ up to unitary equivalence and $\widehat{G}$ to be the set of irreducible representations of $G$ up to unitary equivalence. The set $\widehat{G}$ is called the \tg{unitary dual} of $G$.
	\end{definition}
	\begin{remark}\label{remark bound on the cardinal}
		The classes $\widetilde{G}$ and $\widehat{G}$ are well defined sets. Indeed, as the Hilbert spaces considered in this thesis are all separable, every representation of $G$ is unitary equivalent to a subrepresentation of a representation of $G$ on the Hilbert space $L^2(\N,\C)$. In particular, both $\widetilde{G}$ and $\widehat{G}$ can be be realised as a subset of the set of functions $G\rightarrow \mathcal{U}(L^2(\N,\C))$. Some authors do not require any separability condition in the definition of a representation. We highlight that the unitary dual $\widehat{G}$ of a second-countable locally compact group $G$ is the same set in their set up as \cite[Proposition 1.A.13]{BekkadelaHarpe2020} ensures that any irreducible representation of a second-countable locally compact group on a possibly non-separable Hilbert space is unitary equivalent to a representation on a separable Hilbert space.
	\end{remark}
	We now provide some examples. 
	\begin{example}\label{Example unitary character are irreducible}
			Let $G$ be a locally compact group. A \tg{unitary character} of $G$ is a unitary representation $\chiup :G\rightarrow \mathcal{U}(\C)$ on the one-dimensional Hilbert space $\C$. Notice that these representations are irreducible since the representation space does not contain any proper subspace. The \tg{trivial representation} of $G$ 
			$$1_{\widehat{G}}: G\rightarrow \mathcal{U}(\C): g\mapsto \id_\C$$
			always provides such a representation. However, in general $\widehat{G}$ does not need to contain any other finite dimensional representation (we refer to the second part of the thesis for a concrete example).
	\end{example}
	\begin{example}\label{Example of permutation represenatiton capturing the two transitivity}
			Let $G$ be a locally compact group, let $X$ be a finite set with $\modu{X}\geq 2$ and suppose that $G$ acts continuously and transitively on $X$. Consider the counting point measure $\nu$ on $X$ and the corresponding Koopman representation $\lambda_{X,\nu}$ provided by Example \ref{Example permutation representation}. For brevity we denote this representation by $\lambda_X$ and its representation space by $\ell^2(X)$. Consider the one-dimensional subspace $\C\mathds{1}_X$ of $\ell^2(X)$ spanned by the constant function $\mathds{1}_X$. This subspace is $G$-invariant under $\lambda_X$ and the corresponding subrepresentation of $\lambda_X$ is equivalent to the trivial representation of $G$. Denote by $\ell^2(X)_0$ the orthogonal complement of $\C\mathds{1}_X$ inside $\ell^2(X)$, that is 
			$$\ell^2(X)_0=\{f\in \ell^2(X): \sum_{x\in X}f(x)=\prods{f}{\mathds{1}_X}_{\ell^2(X)}=0\}.$$
			Since $\modu{X}\geq 2$, $\ell^2(X)_0$ is a non-trivial invariant subspace of dimension $\modu{X}-1$. We let $\rho$ the subrepresentation of $\lambda_X$ corresponding to $\ell^2(X)_0$ and notice that 
			$$\lambda_X \simeq 1_{\widehat{G}}\oplus \rho.$$
			This decomposition has a particular interest since the irreducibility or non-irreducibility of the representation $\rho$ captures the dynamics of the action of $G$ on $X$. To be more precise, $\rho$ \textit{is irreducible if and only if $G$ acts $2$-transitively on} $X$, that is $G$ is transitive on $X$ and $\Fix_G(x)=\{g\in G: gx=x\}$ is transitive on $X-\{x\}$ $\forall x\in X$. We provide a proof in light of Schur's Lemma in Example \ref{Example of permutation represenatiton capturing the two transitivity revisiter} below.
	\end{example}
		\noindent The following provides examples of multiplicity free and factor representations that are not necessarily irreducible representations.
	\begin{example}\label{Example irreductble factorial and mult free}
		Let $G$ be a second-countable locally compact group. Let $\mathcal{S}\subseteq \widehat{G}$ be a countable subset of equivalence classes of irreducible representation, for each $\sigma \in \mathcal{S}$ let $n_\sigma\in \N\cup \aleph_0$ be a non-zero cardinal and consider the representation $\pi=\bigoplus_{\sigma\in \mathcal{S}}n_\sigma \sigma$. Such a representation is called \tg{discretely decomposable} and it is:
		\begin{enumerate}[label=(\roman*)]
				\item multiplicity free if and only if $\pi=\bigoplus_{\sigma\in \mathcal{S}}\sigma.$
				\item factorial if and only if $\pi=n_\sigma \sigma$ for some $\sigma\in \widehat{G}$. 
		\end{enumerate}
	\end{example}
	Since every representation of a compact group is discretely decomposable, the above example is universal on these groups. However, certain locally compact groups admit unitary representations that are not discretely decomposable. To develop further the intuition of the reader, we would like to provide examples of multiplicity free and of factor representations that are not covered by Example \ref{Example irreductble factorial and mult free}. This is better done by considering the characterisation of factoriality, multiplicity freeness and irreducibility in terms of intertwining operators. We define the \tg{commutant} $\mathcal{S}'$ of a subset $\mathcal{S}\subseteq \mathcal{L}(\mathcal{H})$ of bounded linear operators on a Hilbert space $\mathcal{H}$ as the subalgebra of $\mathcal{L}(\mathcal{H})$ consisting of the bounded operators commuting with the elements of $\mathcal{S}$ that is 
	$$\mathcal{S}'=\{T\in \mathcal{L}(\mathcal{H}): TS=ST\q \forall S\in \mathcal{S}\}.$$  
	Note that for the choice $\mathcal{S}=\pi(G)=\{\pi(g): g\in G\}\subseteq \mathcal{L}(\Hr{\pi})$ where $\pi$ is a representation of $G$, the commutant $\pi(G)'$ coincides with the algebra of intertwining operators $\Hom_G(\pi,\pi)$. The characterisation of irreducible representations of $G$ by intertwining operators is known as Schur's Lemma and is part of the following more general picture:
	\begin{proposition}\label{Proposition characterisation de mult free facto and irred by commutatnt}
		A representation $\pi$ of a locally compact group $G$ is:
		\begin{enumerate}[label=(\roman*)]
				\item factorial if and only if $\pi(G)'\cap \pi(G)''=\C \id_{\Hr{\pi}}.$
				\item multiplicity free if and only if $\pi(G)'$ is commutative $(\pi(G)'\subseteq\pi(G)'')$.
				\item irreducible if and only if $\pi(G)'=\C\id_{\Hr{\pi}}$ (\tg{Schur's Lemma}). 
		\end{enumerate}
	\end{proposition}
	We illustrate the power of this characterisation with a few examples. 
	\begin{example}
			Let $G$ be a locally compact abelian group. Then, $\pi(G)$ is commutative so that $\pi(G)\subseteq \pi(G)'$. It follows from Schur's Lemma that every irreducible representation of $G$ is a unitary character. 
	\end{example}
	\begin{example}\label{example representations of the quotient provides representations of the group}
			Let $G$ be a locally compact group, let $H$ be a closed normal subgroup of $G$ and let $\rho$ be a representation of $G/H$. Then, $\rho$ gives rise to a representation of $G$ by composition with the quotient projection. The resulting representation is the \tg{inflation} of $\rho$ to $G$. It follows from Proposition \ref{Proposition characterisation de mult free facto and irred by commutatnt} that the inflation of $\rho$ is factorial, irreducible or multiplicity free if and only if $\rho$ is respectively factorial, irreducible or multiplicity free.  
	\end{example}
	\begin{example}\label{Example of permutation represenatiton capturing the two transitivity revisiter}
			Let $G$ be a locally compact group, let $X$ be a finite set with $\modu{X}\geq 2$, suppose that $G$ acts continuously and transitively on $X$ and consider the decomposition $\lambda_X\simeq 1_{\widehat{G}}\oplus \rho$ given in Example \ref{Example of permutation represenatiton capturing the two transitivity}. We prove that $\rho$ is irreducible if and only if $G$ acts $2$-transitively on $X$. In light of Schur's Lemma, $\rho$ is irreducible if and only if $\dim(\Hom_G(\rho,\rho))=1.$ Consider the basis of $\ell^2(X)$ provided by the Dirac functions $\{\delta_x:X\rightarrow \C:x\in X\}$ and denote by $V_{x,y}$ the matrix components of an operator $V\in \mathcal{L}(\ell^2(X))$ in this basis, that is, $V\delta_x=\sum_{y\in X}V_{x,y}\delta_y$. Notice from a  straightforward computation that $V\in \Hom_G(\lambda_X,\lambda_X)$ if and only if 
			\begin{equation*}\label{equation preuve Example of permutation represenatiton capturing the two transitivity revisiter 1}
			V_{gx,gy}=V_{x,y}\q \forall g\in G, \qq\forall x,y\in X. 
			\end{equation*}
			In particular, $\dim(\Hom_G(\lambda_X,\lambda_X))\geq 2$ and $\dim(\Hom_G(\lambda_X,\lambda_X))= 2$ if and only if $G$ acts $2$-transitively on $X$. We are going to show that $\Hom_G(\rho,\rho)$ is a subspace of codimension $1$ of $\Hom_G(\lambda_X,\lambda_X)$. This will prove as desired that $\dim(\Hom_G(\rho,\rho))= 1$ if and only if $G$ acts $2$-transitively on $X$. Now, notice that $\Hom_G(1_{\widehat{G}},1_{\widehat{G}})$ identifies as the subspace of elements of $\Hom_G(\lambda_X,\lambda_X)$ stabilising $\C \mathds{1}_X$ and mapping its orthogonal complement $\ell^2(X)_0$ to $0$. Similarly,
			$\Hom_G(\rho,\rho)$ identifies as the subspace of elements of $\Hom_G(\lambda_X,\lambda_X)$ stabilising $\ell^2(X)_0$ and mapping its orthogonal complement $\C \mathds{1}_X$ to $0$. On the other hand, the zero function is the only element of $\ell^2(X)_0$ that is fixed by $\lambda_X$. In particular, no subrepresentation of $\rho$ is isomorphic to the trivial representation and vice-versa. It follows as desired that 
			$$\Hom_G(\lambda_X,\lambda_X)\simeq \Hom_G(1_{\widehat{G}},1_{\widehat{G}}) \oplus \Hom_G(\rho,\rho)\simeq \C \oplus \Hom_G(\rho,\rho).$$
	\end{example}
	The following provides our first example of a multiplicity free representation that is not discretely decomposable.
	\begin{example}
			Let $G$ be a locally compact abelian group. Then, the regular representation of $G$ is multiplicity free. This follows from Proposition \ref{Proposition characterisation de mult free facto and irred by commutatnt} and from the fact that the commutant of the left-regular representation $\lambda_G(G)'$ is given by the \tg{bi-commutant} $\rho_G(G)''$ of the right-regular representation \cite[Section 13.10.2]{Dixmier1977} (and is therefore commutative by Von Neumann bi-commutant theorem). Notice from Theorem \ref{Theorem les differentes implication topo fell}, Example \ref{Example decomposition en integrale direct d'irred de la reguliere} and Theorem \ref{Theorem de decomposition des representation de type I en integral direct} below that this representation is not discretely decomposable when $G$ is not compact.
	\end{example}
	The following fundamental result due to Murray and von Neumann provides an example of a factor representation that does not contain any irreducible representation and is therefore not discretely decomposable. 
	\begin{example}[{\cite[pg. 60]{Mackey1976theory}}]\label{factor representation of ICC}
			A group is ICC when the conjugacy class of any non-trivial element is infinite. For instance every non-abelian free group is ICC. The left-regular representation of an infinite discrete ICC group is factorial and does not contain any irreducible representation. Indeed, a discrete group whose regular representation contains an irreducible representation must be finite \cite[Proposition 7.A.4]{BekkadelaHarpe2020}. On the other hand, we now prove that a discrete group $\Gamma$ is ICC if and only if its regular representation $\lambda_\Gamma$ is a factor representation. Since $\Gamma$ is discrete, it is unimodular and a basis of $\ell^2(\Gamma)$ is provided by the Dirac functions $\{\delta_\gamma:\gamma\in \Gamma\}$. Since the right and left-regular representations of $\Gamma$ commute with each other, we have that $\rho_\Gamma(\Gamma)\subseteq \lambda_\Gamma(\Gamma)'$ which implies that $\lambda_\Gamma(\Gamma)''\subseteq \rho_\Gamma (\Gamma)'$. It follows that every $T\in \lambda_\Gamma(\Gamma)'\cap \lambda_\Gamma(\Gamma)''$ commutes with both $\lambda_\Gamma$ and $\rho_\Gamma$. Now, let $T\in \lambda_\Gamma(\Gamma)'\cap \lambda_\Gamma(\Gamma)''$ and consider the function $\varphi=T\delta_{1_\Gamma}\in \ell^2(\Gamma)$. Since $\lambda_\Gamma(\gamma)\rho_\Gamma(\gamma)\delta_{1_\Gamma}=\delta_{1_\Gamma}$ for every $\gamma \in \Gamma$, we have that $$\lambda_\Gamma (\gamma)\rho_\Gamma (\gamma)\varphi= \varphi. $$ In particular, $\varphi$ is constant on the conjugacy class of $\Gamma$. Since $\varphi$ belongs to $\ell^2(\Gamma)$, this implies that $\varphi$ is supported on elements with finite conjugacy class. If $\Gamma$ is ICC, this ensures that $\varphi =T\delta_{1_\Gamma}= c\delta_{1_\Gamma}$ for some $c\in  \C$. Therefore, for every $\gamma \in \Gamma$, one has 
			$$T\delta_\gamma=T\lambda_\Gamma(\gamma )\delta_{1_\Gamma}= \lambda_\Gamma(\gamma )T\delta_{1_\Gamma}= c\lambda_\Gamma(\gamma )\delta_{1_\Gamma}=c\delta_{\gamma}.$$
			Since the span of $\{\delta_\gamma: \gamma \in \Gamma\}$ is dense in $\ell^2(\Gamma)$, it follows that $\lambda_\Gamma(\gamma )'\cap \lambda_\Gamma(\gamma )''=\C \id_{\ell^2(\Gamma)}$ and  $\lambda_\Gamma$ is a factor representation. 
			
			Now assume that $\Gamma$ is not ICC, that is, there exists a finite conjugacy class $C\subseteq \Gamma$, $C\not=\{1_\Gamma\}$. Then, the operator $T= \sum_{\gamma \in C}\lambda_\Gamma(\gamma)$ belongs to $\lambda_\Gamma(\Gamma)''$. On the other hand, this operator also belongs to $\lambda_\Gamma(\Gamma)'$ since for every $\omega \in \Gamma$ we have
			\begin{equation*}
			\begin{split}
			T\lambda_\Gamma(\omega)=\sum_{\gamma \in C}\lambda_\Gamma(\gamma)\lambda_\Gamma (\omega)= \lambda_\Gamma(\omega)\sum_{\gamma \in C}\lambda_\Gamma(\omega^{-1}\gamma \omega)=\lambda_\Gamma(\omega)T.
			\end{split}
			\end{equation*}
			On the other hand, $T$ is not a multiple of the identity since $T\delta_{1_\Gamma}=\mathds{1}_C.$ It follows that $\lambda_\Gamma$ is not factorial.
	\end{example}
	We end this section with an emphasis on the concept of factor representation and of Type {\rm I} groups. We recall that the representations of compact groups all admit an essentially unique decomposition as a direct sum of irreducible representations. However, the above example shows that this behaviour does not extend to general locally compact since certain representations do not even contain an irreducible representation. On the other hand, we will see in Section \ref{section direct integral decompostion} that the notion of direct integral of representations provides a suitable generalisation of the concept of direct sum of representations and that any representation admits decompositions as a direct integral of factor representations or even irreducible representations. This direct integral decomposition into irreducible representations is in general far from being unique. On the other hand, by contrast with the decomposition as a direct integral of irreducible representations, there is a canonical choice of decomposition as a direct integral of factors. This motivates the study of factor representations which tends to be more appropriate building blocks of representations in the general case. Now, notice that the factors provided by Example \ref{Example irreductble factorial and mult free} and Example \ref{factor representation of ICC} are quite different in nature. In fact, there are three types of factor representations: Types {\rm I}, {\rm II} and {\rm III}. These types can be defined independently from one another by relying on von Neumann algebras and traces see \cite[Definition $7.B.7$]{BekkadelaHarpe2020}. However, we will only be interested into the dichotomy between Type {\rm I} and non-Type {\rm I} representations in this thesis and hence adopt a more practical definition provided by the relation of quasi-equivalence as in \cite{Mackey1976theory}.
	\begin{definition}
		Let $G$ be a locally compact group. A factor representation $\pi$ of $G$ is:
		\begin{enumerate}[label=(\roman*)]
			\item Type \tg{\rm I} if it is quasi-equivalent to an irreducible representation. 
			\item Type \tg{\rm II} if it is neither type \tg{\rm I} nor \tg{\rm III}.
			\item Type \tg{\rm III} if $\pi'\approx \pi$ implies $\pi'\simeq \pi$ for every representation $\pi'$ of $G$.
		\end{enumerate}
	\end{definition}
	In particular, notice that a factor representation is of type \tg{\rm I} if and only if it is equivalent to a multiple $\aleph\sigma$ of an irreducible representation $\sigma$ of $G$ and hence is discretely decomposable. Example \ref{Example irreductble factorial and mult free} therefore provides an example of a type \tg{\rm I} factor while Example \ref{factor representation of ICC} provides an example of a factor representation of type \tg{\rm II} since the left-regular representation of a countable discrete group is never of type {\rm III}, see \cite[pg. 61]{Mackey1976theory}. 
	\begin{definition}
		A locally compact group $G$ is of \tg{Type} \tg{\rm I} if every factor representation of $G$ is of type \tg{\rm I}.
	\end{definition}
	Loosely speaking Type {\rm I} groups are the locally compact groups whose representation theory reduces suitably to the study of irreducible representations. To be more precise, Type {\rm I} groups are exactly the locally compact groups whose dual space $\widehat{G}$ is a $T_0$-separated topological space and a standard Borel space for the natural topology and Borel structure, see Theorem \ref{Glimm's theorem} below. Furthermore those are also the locally compact groups for which each representation admits an essentially unique decomposition as a direct integral of irreducible representations, see Theorem \ref{Theorem de decomposition des representation de type I en integral direct}. By contrast, when the group is not of Type {\rm I}, it has a representation with non-unique decomposition, see Theorem \ref{theorem non-uniqueness of direct integral decompostiion}. It is therefore important to identify which locally compact groups are of Type {\rm I}. A fundamental result from Thoma shows that countable discrete group are almost never of Type {\rm I}.
	\begin{theorem}[\cite{Thoma1968}]
		A countable discrete group is Type {\rm I} if and only if it is virtually abelian, that is, it admits a finite index abelian subgroup.
	\end{theorem}
	For non-discrete groups, the question has been solved for various family of locally compact groups. We gather some important contributions in the following statement.
	\begin{theorem}\label{Theorem example of type I groups}
		The following locally compact groups are Type {\rm I}:
		\begin{enumerate}[label=(\roman*)]
			\item Locally compact abelian groups.
			\item Compact groups.
			\item Connected nilpotent locally compact groups, and more generally connected-
			by-compact locally compact groups whose solvable radical is actually nilpotent \cite{Dixmier1959},\cite{Lipsman1972}.
			\item Reductive algebraic groups over a local field \cite{Harish1970},\cite{Bernshtein1974}.
			\item Linear algebraic groups over a local field of characteristic $0$ \cite{Dixmier1957},\cite{BekkaEchterhoff2021}.
			\item Some reductive adelic groups \cite{Clozel2002}.
			\item Full groups of automorphisms of regular trees \cite{Olshanskii1977},\cite{Olshanskii1980}.
		\end{enumerate}
	\end{theorem}
	\section{Induced representations}\label{Section induced representations}
	\subsection{Definition and properties}
	Let $G$ be a locally compact group and let $H$ be a closed subgroup. Given a representation $\pi$ of $G$, we can build a representation of $H$ with the same representation space $\Hr{\pi}$ as $\pi$ by considering the restriction $\Res_H^G(\pi)$ of $\pi$ to $H$. Similarly, given a representation $\sigma$ of $H$, it is possible to build a representation $\Ind_H^G(\sigma)$ of $G$ called the \tg{induced representation} of $\sigma$. The purpose of this section is to provide a proper definition of this concept and to recall some of its fundamental properties. The notion of induced representation can always be defined in the context where $G$ is a locally compact group and $H\leq G$ is a closed subgroup. However, most of the complexity vanishes when $H$ is open in $G$ (because the quotient space $G/H$ is discrete). Although the notion of induced representations from a closed subgroup of $H$ is used in this thesis, we provide an explicit description of the induced representation only when $H$ is open in this section. On the other hand, the results stated without the assumption that $H$ is open stay valid in full generality. We refer the reader to \cite[Chapter $2$]{KaniuthTaylor2013} for details in both setups. 
	
	Let $G$ be a locally compact group, $H\leq G$ be an open subgroup and $\sigma$ be a representation of $H$. The induced representation  $\Ind_{H}^G(\sigma)$ is a representation of $G$ with representation space given by
	$$\Ind_H^G(\Hr{\sigma})=\bigg\{\phi:G\rightarrow\Hr{\sigma}\Big\lvert  \phi(gh)=\sigma(h^{-1})\phi(g), \s{gH\in G/H}{}\prods{\phi(g)}{\phi(g)}_{\Hr{\sigma}}<+\infty\qq \bigg\}.$$
	Notice that the condition on the sum is well defined since $ \phi(gh)=\sigma(h^{-1})\phi(g)$ and since $\sigma$ is unitary. We define the inner product of any two element $\psi,\phi\in \Ind_H^G(\Hr{\sigma})$ as
	\begin{equation*}
	\prods{\psi}{\phi}_{\Ind_H^G(\Hr{\sigma})}= \s{gH\in G/H}{} \prods{\psi(g)}{\phi(g)}_{\Hr{\sigma}}.
	\end{equation*}
	Equipped with this inner product, $\Ind_H^G(\Hr{\sigma})$ is a separable complex Hilbert space. The \tg{induced representation} $\Ind_H^G(\sigma)$ is the representation of $G$ with representation space $\Ind_H^G(\Hr{\sigma})$ defined by
	\begin{equation*}
	\big\lb\Ind_H^G(\sigma)(h)\big\rb\phi(g)=\phi(h^{-1}g)\q \forall \phi \in \Ind_H^G(\Hr{\sigma}) \mbox{ and }\forall g,h\in G.
	\end{equation*}
	In the more general context where $H$ is closed rather than open, this sort of description fails since $G/H$ does not admit in general a $G$-invariant measure. However, as $H$ is closed, the quotient space $G/H$ admits a measure $\mu$ that is absolutely continuous with respect to each of its $G$-left-translates. The measure class of $\mu$ is essentially unique and is called the \tg{invariant measure} class of $G/H$. The definition of the induction in this setup is similar to the above but requires using the Radon-Nikodym derivatives associated with the invariant measure class of $G/H$ (depending on the modular functions of $G$ and $H$).
	
	It is natural to ask how the processes of induction and restriction of representations behave with respect to one another. When $G$ and $H$ are finite, strong relations between these operations are provided by the Frobenius reciprocity. This reciprocity ensures for every representation $\sigma$ of $H$ and $\pi$ of $G$ we have that
	\begin{equation*}
	{\rm I}(\Ind_H^G(\sigma),\pi)= {\rm I} (\sigma, \Res_H^G(\pi)) 
	\end{equation*}
	where ${\rm I}(\pi_1,\pi_2)$ is the dimension of the space of intertwining operators between the two representations $\pi_1$ and $\pi_2$. However, in the general context where $G$ is a locally compact group and $H\leq G$ is a closed subgroup, this reciprocity is generally nothing more than an inequality. 
	\begin{theorem}[{\cite[Corollary $1$ of Theorem $3.8$]{Mackey1976theory}} ]\label{criterion of mackey weak frobenius reciprocity}
		Let $G$ be a locally compact group, $H\leq G$ be a closed subgroup, $\pi$ be a representation of $G$ and $\sigma$ be a representation of $H$. Then, we have that
		\begin{equation*}
		{\rm I}(\Ind_H^G(\sigma),\pi)\leq {\rm I} (\sigma, \Res_H^G(\pi)).
		\end{equation*}
		Furthermore, if $\lb G: H\rb $ is finite, this relation becomes an equality.
	\end{theorem}
	In the rest of this section we provide criteria to check the irreducibility and equivalence of induced representations. Those results are classical and will be useful later in these notes. However, since we did not find any convenient references we present a proof relying on the above description of induced representations.
	\begin{lemma}\label{la correspondance Hsigma D}
		Let $G$ be a locally compact group, $H\leq G$ be an open subgroup and let $\sigma$ be a representation of $H$. There exists an isomorphism between the Hilbert spaces $\mathfrak{D}^\sigma_H=\{\phi\in \
		\Ind_H^G(\Hr{\sigma}): \phi(g)=0 \qq \forall g\in G-H\}$ and $\Hr{\sigma}$ that intertwines the representations $\restr{\Ind_H^G(\sigma)}{H}$ and $\sigma$. Furthermore, we have that
		\begin{equation*}\label{la densite des fonction a support dans H dans l induite}
		\Ind_H^G(\Hr{\sigma})=\ub{ \underset{gH\in G/H}{\bigoplus}\qq\big\lb\Ind_H^G(\sigma)(g)\big\rb\mathfrak{D}^\sigma_H}.
		\end{equation*}
	\end{lemma}
	\begin{proof} 
		First notice that $\mathfrak{D}^\sigma_H$ is a closed subspace of $\Ind_H^G(\Hr{\sigma})$ and therefore defines a Hilbert space. Now, for every $\xi \in \Hr{\sigma}$, let \begin{equation*}
		\phi_\xi(g)\qq=\qq \begin{cases}
		\qq\sigma(g^{-1})\xi&\mbox{ if } g\in H,\\
		\q \q 0 &\mbox{ if } g\not\in H.
		\end{cases}
		\end{equation*}
		Notice that $\phi_\xi \in \Ind_H^G(\Hr{\sigma})$ and that $\phi_\xi=0$ if and only if $\xi=0$. On the other hand, every function $\phi \in \mathfrak{D}^\sigma_H$ is uniquely determined by the value it takes on $1_G$ since $\phi(g)=0 \qq \forall g\in G -H$ and $\phi(g)=\sigma(g^{-1})\phi(1_G)$  $\forall g\in H$. Hence, we have that $\mathfrak{D}^\sigma_H=\{\phi_\xi: \xi \in \Hr{\sigma}\}$ and the map
		\begin{equation*}\label{equation the isomorphism for the induced representation between the subspace of the small rep and the set of coeff function supported inside the subgroup}
		\Phi\fct{\Hr{\sigma}}{\mathfrak{D}^\sigma_H}{\xi}{\phi_\xi}\qq
		\end{equation*}
		is a linear isomorphism between $\Hr{\sigma}$ and $\mathfrak{D}^\sigma_H$. Moreover, $\Phi$ is unitary since $\forall \xi,\eta \in \Hr{\sigma}$ we have
		\begin{equation*}
		\begin{split}
		\prods{\Phi(\xi)}{\Phi(\eta)}_{\Ind_H^G(\Hr{\sigma})}&=\s{gH\in G/H}{} \prods{\phi_\xi(g)}{\phi_\eta(g)}=\prods{\phi_\xi(1_G)}{\phi_\eta(1_G)}=\prods{\xi}{\eta}.
		\end{split}
		\end{equation*}
		Finally, $\Phi$ intertwines $\sigma$ and $\restr{\Ind_H^G(\sigma)}{H}$ since for all $ h\in H$ and every $g\in G$ we have that
		\begin{equation*}
		\begin{split}
		\Phi(\sigma(h)\xi)(g)&=\phi_{\sigma(h)\xi}(g)=\begin{cases}
		\sigma(g^{-1}h)\xi &\mbox{ if }g\in H\\
		\q \q 0&\mbox{ if } g\not\in H
		\end{cases}\\
		&=\phi_{\xi}(h^{-1}g)=\big\lb\Ind_H^G(\sigma)(h)\big\rb\phi_{\xi}(g).
		\end{split}
		\end{equation*}
		This proves the first part of the claim. To prove the second part of the claim, let $\phi\in \Ind_H^G(\Hr{\sigma})$ and let us show that 
		$$\phi=\s{gH\in G/H}{}\lb \Ind_H^G(\sigma)(g)\rb\Phi(\phi(g)).$$ 
		First, notice that $\Ind_H^G(\sigma)(g)\mathfrak{D}^\sigma_H$ is the set of functions of $\Ind_{H}^G(\Hr{\sigma})$ that are supported on $gH$. On the other hand, $\forall g,t\in G$, $\forall h\in H$ we have that 
		\begin{equation*}
		\begin{split}
		\big\lb\big\lb\Ind_{H}^G(\sigma)(gh)\big\rb\Phi(\phi(gh))\big\rb(t)&=\Phi(\phi(gh))(h^{-1}g^{-1}t)\\
		&=\begin{cases}
		\sigma(t^{-1}gh)\phi(gh) &\mbox{ if }h^{-1}g^{-1}t\in H\\
		\q \q 0&\mbox{ if }h^{-1}g^{-1}t\not\in H
		\end{cases}\\
		&=\begin{cases}
		\sigma(t^{-1}g)\phi(g) &\mbox{ if }g^{-1}t\in H\\
		\q \q 0&\mbox{ if }g^{-1}t\not \in H
		\end{cases}\\
		&=\Phi(\phi(g))(g^{-1}t)=\big \lb\lb\Ind_{H}^G(\sigma)(g)\rb\Phi(\phi(g))\big\rb(t).
		\end{split}
		\end{equation*}
		This proves that the map $G/H\rightarrow \Ind_H^G(\Hr{\sigma}): gH\mapsto \big\lb\Ind_{H}^G(\sigma)(g)\big\rb\Phi(\phi(g))$ is well defined.
		On the other hand, $\forall g\in G$, $\forall h\in H$, we have
		\begin{equation*}
		\big\lVert\big\lb\big\lb \Ind_{H}^G(\sigma)(g)\big\rb\Phi(\phi(g))\big\rb(gh)\big\rVert=\norm{\Phi(\phi(g))(h)}{}= \norm{\sigma(h^{-1})\phi(g)}{}= \norm{\phi(g)}{}.
		\end{equation*} 
		Hence, we obtain that
		\begin{equation*}
		\begin{split}
		&\Bigg\lVert\s{gH\in G/H}{} \big\lb\Ind_{H}^G(\sigma)(g)\big\rb\Phi(\phi(g)) \Bigg\rVert_{\Ind^G_H(\Hr{\sigma})}^2 = \s{tH\in G/H}{} \Bigg \lVert \s{gH\in G/H}{} \big\lb \big\lb\Ind_{H}^G(\sigma)(g)\big\rb\Phi(\phi(g))\big \rb(t)\Bigg \rVert^2\\
		&\q\q= \s{tH\in G/H}{} \big \lVert \big\lb \big\lb\Ind_{H}^G(\sigma)(t)\big\rb\Phi(\phi(t))\big \rb(t)\big \rVert^2= \s{tH\in G/H}{}\norm{\phi(t)}{}^2=\norm{\phi}{\Ind_H^G(\Hr{\sigma})}^2.
		\end{split}
		\end{equation*}
		This proves that $\s{gH\in G/H}{} \big\lb\Ind_{H}^G(\sigma)(g)\big \rb\Phi(\phi(g))$ belongs to $\ub{ \underset{gH\in G/H}{\bigoplus}\qq\Ind_H^G(\sigma)(g)\mathfrak{D}^\sigma_H}$. Finally, a direct computation shows that
		\begin{equation*}
		\begin{split}
		&\Bigg\lVert\phi-\s{gH\in G/H}{} \Ind_{H}^G(\sigma)(g)\Phi(\phi(g))\Bigg\lVert^2_{\Ind_{H}^G(\sigma)}=\s{tH\in G/H}{}\Bigg\lVert\phi(t)-\s{gH\in G/H}{}  \big\lb\Ind_{H}^G(\sigma)(g)\big\rb\Phi(\phi(g))(t)\Bigg\lVert^2\\
		&\q\q =\s{tH\in G/H}{} \big\lVert\phi(t)-\Ind_{H}^G(\sigma)(t)\Phi(\phi(t))(t)\big\lVert^2= \s{tH\in G/H}{} \norm{\phi(t)- \phi (t)}{}^2=0.
		\end{split}
		\end{equation*}
	\end{proof}
	The following provides a useful criterion to check the irreducibility of an induced representation.  
	\begin{proposition}\label{a criterion so that the induced of an irreducible is irreducible}
		Let $G$ be a locally compact group, $H\leq G$ be an open subgroup and let $\sigma$ be an irreducible representation of $H$. Then, $\Ind_H^G(\sigma)$ is irreducible if and only if every non-zero closed invariant subspace of $\Ind_H^G(\Hr{\sigma})$ contains a non-zero function supported in $H$. 
	\end{proposition}
	\begin{proof}
		The forward implication is trivial. For the other implication, let $M$ be a non-zero closed invariant subspace of $\Ind_H^G(\Hr{\sigma})$. By the hypothesis, $M\cap \mathfrak{D}_{H}^{\sigma}\not=\es$. In particular, since $\sigma$ is irreducible, the correspondence of Lemma \ref{la correspondance Hsigma D} implies that $\mathfrak{D}_H^{\sigma}\subseteq M$. We conclude from that same lemma that $M=\Ind_{H}^G(\Hr{\sigma})$ since $M$ is a closed invariant subspace of $\Ind_H^G(\Hr{\sigma})$ and since \begin{equation*}
		\Ind_H^G(\Hr{\sigma})=\ub{ \underset{gH\in G/H}{\bigoplus}\qq\big\lb\Ind_H^G(\sigma)(g)\big\rb\mathfrak{D}^\sigma_H}.
		\end{equation*}
	\end{proof}
	On the other hand, it is possible to determine whether an irreducible representation of $G$ is induced from a representation of one of its open subgroups by looking at the matrix coefficient functions. To be more precise, we have the following result.
	\begin{lemma}\label{lemma determine wether a representation is induced from an open subgroup of G}
		Let $G$ be a locally compact group, $H\leq G$ be an open subgroup and $\pi$ be an irreducible representation of $G$ with a non-zero matrix coefficient supported on $H$. Then, there exists an irreducible representation $\sigma$ of $H$ such that $\pi\simeq \Ind_H^G(\sigma)$.
	\end{lemma} 
	\begin{proof}
		The proof of this result relies on the notions of function of positive type (Definition \ref{definition functions of positive type}) and the GNS construction (Theorem \ref{Theorem GNS construction}) that will be introduced in Section \ref{section Fell topology} below.
		Let $\xi\in \Hr{\pi}$ be such that the matrix coefficient function $$\varphi_{\xi,\xi}:G\rightarrow \C : g\mapsto \prods{\pi(g)\xi}{\xi}_{\Hr{\pi}}$$ is not identically zero and is supported on $H$. Since $\varphi_{\xi,\xi}$ is a function of positive type on $G$, notice that its restriction $\restr{\varphi_{\xi,\xi}}{H}$ defines a function of positive type on $H$. In particular, the GNS construction ensures the existence of a cyclic representation $\sigma$ of $H$ and of a cyclic vector $\eta\in \Hr{\sigma}$ such that $$\varphi_{\xi,\xi}(h)=\prods{\sigma(h)\eta}{\eta}_{\Hr{\sigma}}.$$
		We claim that $\pi\simeq \Ind_H^G(\sigma)$. Indeed, as shown in the proof of Lemma \ref{la correspondance Hsigma D} the function
		\begin{equation*}
		\phi_\eta(g)\qq=\qq \begin{cases}
		\qq\sigma(g^{-1})\eta&\mbox{ if } g\in H\\
		\q \q 0 &\mbox{ if } g\not\in H
		\end{cases}
		\end{equation*}
		defines an element of $\Ind_H^G(\Hr{\sigma})$. On the other hand, notice that the function of positive type 
		$$\varphi_{\phi_{\eta},\phi_{\eta}}:G\rightarrow \C : g\mapsto\prods{\Ind_H^G(\sigma)(g)\phi_{\eta}}{\phi_\eta}_{\Ind_H^G(\Hr{\sigma})}$$
		is supported on $H$ and that for every $h\in H$ we have that
		\begin{equation*}
		\begin{split}
		\varphi_{\phi_{\eta},\phi_{\eta}}(h)= \prods{\Ind_H^G(\sigma)(h)\phi_{\eta}}{\phi_\eta}_{\Ind_H^G(\Hr{\sigma})}=\prods{\sigma(h)\eta}{\eta}_{\Hr{\sigma}}=\varphi_{\xi,\xi}(h).
		\end{split}
		\end{equation*}
		This proves that $\varphi_{\xi,\xi}=\varphi_{\phi_{\eta},\phi_{\eta}}$. On the other hand,  $\xi$ is cyclic since $\pi$ is irreducible. It follows from the uniqueness of the GNS construction that $\pi\simeq\Ind_H^G(\Hr{\sigma})$. The irreducibility of $\sigma$ follows from this equivalence as $\pi$ is irreducible.
	\end{proof}
	Finally, the following result provides a sufficient condition for induced representations to be inequivalent. 
	\begin{lemma}\label{les induites dinequivalente sont inequivalente}
		Let $G$ be a locally compact group, $H\leq G$ be an open subgroup and let $\sigma_1,\sigma_2$ be two inequivalent irreducible representations of $H$. Suppose that there exists a subgroup $K\leq H$ such that $\mathfrak{D}^{\sigma_i}_H=(\Ind_H^G(\Hr{\sigma_i})\big)^{{K}}$. Then, $\Ind_H^G(\sigma_1)$ and  $\Ind_H^G(\sigma_2)$ are inequivalent. 
	\end{lemma}
	\begin{proof}
		The proof is by contradiction. Suppose that there exists a unitary operator ${\rm U}$ intertwining $\Ind_H^G(\sigma_1)$ and $\Ind_H^G(\sigma_2)$. In particular, we have ${\rm U}(\Ind_H^G(\Hr{\sigma_1})\big)^{{K}}=(\Ind_H^G(\Hr{\sigma_2})\big)^{{K}}$. Let $\Phi_i: \Hr{\sigma_i}\rightarrow \mathfrak{D}_H^{\sigma_i}$ be the correspondences given by Lemma \ref{la correspondance Hsigma D} and notice that $\Phi_i(\Hr{\sigma_i})=\mathfrak{D}_H^{\sigma_i}=(\Ind_H^G(\Hr{\sigma_i})\big)^{{K}}.$ In particular, since $\Phi_i$ is a unitary operator intertwining $\sigma_i$ and $\Res_H^G\big(\Ind_H^G(\sigma_i)\big)$, the unitary operator $\Phi_2^{-1}{\rm U}\Phi_1$ intertwines $\sigma_1$ and $\sigma_2$ which leads to a contradiction.
	\end{proof}
	
	\subsection{Dynamics for index two subgroups}\label{section induction depuis sous group indice 2}

	In the case of finite groups, Clifford theory describes the restriction induction dynamics between the representation of a group $G$ and any of its normal subgroup $H$ \cite[Chapter 11]{CurtisReiner1981}. In the more general context where $G$ is locally compact group and $H\leq G$ is a closed normal subgroup, the dynamics becomes far more subtle to analyse and is provided by Mackey's Theory, see \cite[Chapter 3]{Mackey1976theory} and \cite[Chapter 4]{KaniuthTaylor2013}. The purpose of this section is to explain the relations between the irreducible representations of a locally compact group $G$ and the irreducible representations of its closed subgroups $H\leq G$ of index $2$. Theorem \ref{les rep dun group loc compact par rapport a celle d'un de ses sous groupes} describes these dynamics. 
	The results presented in this section are part of the Folklore and follow directly from  Mackey's theory. However, the author was unable to find appropriate references to the present setup and decided to provide elementary proofs for convenience of the reader. We refer to \cite[Chapter 8, Section 5.1]{Procesi2007} and \cite[Section 2]{LabesseSchwermer2019} for similar discussions and partial overlaps.
	
	We start with some preliminaries. Let $G$ be a locally compact group and let $H\leq G$ be a closed normal subgroup of finite index (in particular $H$ is open in $G$). We recall two procedures that can be applied to the representations of these groups. The first procedure is the conjugation. To be more precise, for every representation $\sigma$ of $H$ and for all $g\in G$ we have a \tg{conjugate representation} $\sigma^g$ with representation space $\Hr{\sigma}$ defined as $$\sigma^g(h)=\sigma(ghg^{-1}) \q \forall h\in H.$$ Notice that this action preserves irreducibility as a vector $\xi\in \Hr{\sigma}$ is cyclic for $\sigma$ if and only if it is cyclic for $\sigma^g$. Since the representation $\sigma^g$ does not depend up to equivalence on the representative of the coset $gH$, we obtain an action of $G/H$ on $\widehat{H}$. The second procedure consists of twisting the representations of $G$ with characters. To be more precise, consider a representation $\pi$ of $G$ and a unitary character $\chiup:G\rightarrow \C$, we define the \tg{twisted representation} of $\pi$ by $\chiup$ by $$\pi\otimes \chiup: G\rightarrow\mathcal{U}(\Hr{\pi}):g \mapsto \chiup (g)\pi(g).$$
	Notice that $\pi\otimes\chiup$ is irreducible if and only if $\pi$ is irreducible, as a vector is cyclic for one representation if and only if it is cyclic for the other. 
	
	We are finally able to state properly the correspondence. Consider a locally compact group $G$ and a closed subgroup $H\leq G$ of index $2$. In particular, $H$ is normal, $G/H$ is isomorphic to the cyclic group of order two and there is a unitary character $$\tau: G\rightarrow \C:g\mapsto \begin{cases}
	+1\qq &\mbox{if }g\in H\\
	-1\qq &\mbox{if }g\in G-H
	\end{cases}$$ associated with the quotient $G/H$. Let $t\in G-H$ and for every representation $\pi$ of $G$, let $\pi^\tau$ denote twisted representation $\pi\otimes \tau$ of $G$. The purpose of this section is to prove the following result.
	\begin{theorem}\label{les rep dun group loc compact par rapport a celle d'un de ses sous groupes}
		Let $G$ be a locally compact group, $H\leq G$ be a closed subgroup of index $2$ and $t\in G-H$. For every irreducible representation $\pi$ of $G$ the following assertions hold: 
		\begin{itemize}
			\item $\pi\not\simeq \pi^\tau$ if and only if $\Res_H^G(\pi)$ is an irreducible representation of $H$ and in that case $\Res_H^G(\pi)\simeq\Res_H^G(\pi)^t$.
			\item $\pi\simeq \pi^\tau$ if and only if $\Res_H^G(\pi)\simeq \sigma \oplus \sigma^t$ for some irreducible representation  $\sigma$ of $H$ and in that case $\sigma\not \simeq \sigma^t$.
		\end{itemize}
		For every irreducible representation $\sigma$ of $H$, the following assertions hold: 
		\begin{itemize}
			\item $\sigma\not\simeq \sigma^t$ if and only if $\Ind_H^G(\sigma)$ is an irreducible representation of $G$ and in that case $\Ind_H^G(\sigma)\simeq\Ind_H^G(\sigma)^\tau$.
			\item $\sigma \simeq \sigma^t$ if and only if $\Ind_H^G(\sigma)\simeq \pi \oplus \pi^\tau$ for some irreducible representation $\pi$ of $G$ and in that case $\pi \not \simeq \pi^\tau$.
		\end{itemize}
		Furthermore every irreducible representation of $G$ is a subrepresentation of $\Ind_H^G(\sigma)$ for some irreducible representation $\sigma$ of $H$ and every irreducible representation of $H$ is a subrepresentation of $\Res_H^G(\pi)$ for some irreducible representation $\pi$ of $G$.
	\end{theorem}
	The proof of this theorem is split in the following few results. 
	\begin{lemma}\label{Lemma les rep de G se decompose ou pas en irred si H est d'indice 2}
		Let $\pi$ be an irreducible representation of $G$. Then, exactly one of the following happens:
		\begin{itemize}
			\item $\Res_H^G(\pi)$ is an irreducible representation of $H$ and $\Res_H^G(\pi)\simeq \Res_H^G(\pi)^t$.
			\item $\Res_H^G(\pi)\simeq\sigma\oplus \sigma^t$ for some irreducible representation $\sigma$ of $H$.
		\end{itemize}
	\end{lemma} 
	\begin{proof}
		If $\Res_H^G(\pi)$ is an irreducible representation of $H$, notice for every $h\in H$ and every $\xi \in \Hr{\pi}$ that
		\begin{equation*}
		\begin{split}
		\pi(t)\lb \Res_H^G(\pi)(h)\rb \xi&=\pi(t)\pi(h)\xi=\pi(tht^{-1})\pi(t)\xi\\
		&=\lb \Res_H^G(\pi)(tht^{-1})\rb\pi(t)\xi=\lb\Res_H^G(\pi)^t(h)\rb\pi(t)\xi.
		\end{split}			
		\end{equation*}
		In particular, $\pi(t):\Hr{\pi}\rightarrow\Hr{\pi}$ is an intertwining operator between $\Res_H^G(\pi)$ and $\Res_H^G(\pi)^t$ which settles the first case.
		
		Now, suppose that $\Res_H^G(\pi)$ is not an irreducible representation of $G$. Since $\pi$ is irreducible, any non-zero $\xi\in \Hr{\pi}$ is a cyclic vector. Hence, the subspace spanned by $\{\pi(g)\xi : g\in G\}$ is dense in $\Hr{\pi}$. On the other hand, since $\Res_H^G(\pi)$ is not irreducible, there exists a non-zero vector $\xi \in \Hr{\pi}$ such that the subspace spanned by $\{\pi(h)\xi : h\in H\}$ is not dense in $\Hr{\pi}$. Let $M$ denote the closure of this space. First, let us show that $(\Res_H^G(\pi),M)$ is irreducible. Let $N$ be a $\pi(H)$-invariant subspace of $M$ and let us prove that $N$ cannot be a proper subspace of $M$. Notice that for every  closed $\pi(H)$-invariant subspace $L$ of $\Hr{\pi}$, the subspace $\pi(t)L$ is also $\pi(H)$-invariant since $\pi(H)\pi(t)L=\pi(Ht)L=\pi(tH)L=\pi(t)L$. Now, since $\xi$ is a cyclic vector for $\pi$, notice that $\Hr{\pi}=M+ \pi(t)M$ (where the sum is a priori not a direct sum). On the other hand, since $\Hr{\pi }\not=M$ and since $t^2\in H$, notice that $M\not \subseteq \pi(t)M$. In particular, replacing $N$ by $N^\perp\cap M$ if needed, we can suppose that $\Hr{\pi}\not=N +\pi(t) M$ and therefore that $\Hr{\pi}\not= N+\pi(t)N $. On the other hand, $N+ \pi(t)N$ is a closed $\pi(G)$-invariant subspace of $\pi$. Since $\pi$ is irreducible, this implies that $N+\pi(t)N=\{0\}$ and therefore that $N=\{0\}$ which proves that $(\Res_H^G(\pi),M)$ is irreducible. Since $\pi(t)\pi(t)M=M$, the same reasoning shows that $(\Res_H^G(\pi),\pi(t)M)$ is an irreducible representation of $H$. In particular, since $M^\perp$ is a closed $\pi(H)$-invariant subspace of $\pi(t)M$ and since $(\pi(t)M)^{\perp}$ is a closed $\pi(H)$-invariant subspace of $M$ this proves that $\Hr{\pi}=M\oplus \pi(t)M$. Now let $(\sigma, \Hr{\sigma})=(\Res_H^G(\pi), M)$ and notice that the unitary operator $\pi(t):\pi(t)M\rightarrow M$ intertwines $(\Res_H^G(\pi),\pi(t)M)$ and $(\sigma^t,\Hr{\sigma})$ since for every $\xi \in \pi(t)M$ and for every $h\in H$, we have
		\begin{equation*}
		\begin{split}
		\pi(t)\lb \Res_H^G(\pi)(h)\rb\xi&=\pi(t)\pi(h)\pi(t^{-1})\pi(t)\xi=\pi(tht^{-1})\pi(t)\xi \\
		&=\lb \Res_H^G(\pi)(tht^{-1})\rb \pi(t)\xi =\sigma^t(h)\pi(t)\xi.
		\end{split}
		\end{equation*} 
		This proves as desired that $\Res_H^G(\pi)\simeq \sigma \oplus \sigma^t$ for some irreducible representations $\sigma$ of $H$. 
	\end{proof}
	\begin{lemma}\label{Lemma equivalentce criterion between sigma et sigma g et pi et pitau}
		Let $\pi$ be an irreducible representation of $G$ such that $\pi \simeq \pi^\tau$. Then, $\Res_H^G(\pi)\simeq\sigma \oplus \sigma^t$ for some irreducible representation $\sigma$ of $H$.
	\end{lemma}
	\begin{proof}
		Lemma \ref{Lemma les rep de G se decompose ou pas en irred si H est d'indice 2} ensures that $\Res_H^G(\pi)$ is either irreducible or split as desired. Suppose for a contradiction that $\sigma= \Res_H^G(\pi)$ is irreducible and let $\mathcal{U}: \Hr{\pi}\rightarrow \Hr{\pi}$ be the unitary operator intertwining $\pi$ and $\pi^\tau$. Notice that $\pi$ and $\sigma$ have the same representation space $\Hr{\pi}$. Furthermore, for every $h\in H$, we have that $\pi(h)=\pi^\tau(h)=\sigma(h)$. In particular, $\mathcal{U}$ is a unitary operator that intertwines $\sigma$ with itself. Since $\sigma$ is irreducible this implies that $\mathcal{U}$ is a scalar multiple of the identity. However, this is impossible since for every $h\in H$ and every $\xi \in \Hr{\pi}$ we have 
		\begin{equation*}
		\mathcal{U}\pi(th)\xi=\pi^{\tau}(th)\mathcal{U}\xi=-\pi(th)\mathcal{U}\xi.
		\end{equation*} 
		We obtain as desired that $\Res_H^G(\pi)\simeq\sigma \oplus\sigma^t$ for some irreducible representation $\sigma$ of $H$ when $\pi\simeq \pi^\tau$.
	\end{proof}
	The rest of the proof relies on the fact that the Frobenius reciprocity holds for finite index subgroup by Theorem \ref{criterion of mackey weak frobenius reciprocity}.
	\begin{proposition}\label{prop des induite qui se split ou pas}
		Let $\sigma$ be an irreducible representation of $H$. Then, we have that $\Res_H^G(\Ind_H^G(\sigma))\simeq\sigma \oplus \sigma^t$ and the following hold:
		\begin{itemize}
			\item $\sigma\not\simeq \sigma^t$ if and only if $\Ind_H^G(\sigma)$ is an irreducible representation of $G$ and in that case $\Ind_H^G(\sigma)\simeq\Ind_H^G(\sigma)^\tau$.
			\item $\sigma \simeq \sigma^t$ if and only if $\Ind_H^G(\sigma)\simeq \pi \oplus \pi^\tau$ for some irreducible representation $\pi$ of $G$.
		\end{itemize}
	\end{proposition}
	\begin{proof}
		We start by showing that $\Res_H^G(\Ind_H^G(\sigma))\simeq\sigma \oplus \sigma^t$. We set $$\mathcal{L}=\{\varphi\in \Ind_H^G(\Hr{\sigma}) : \mbox{supp}(\varphi)\subseteq H\}\mbox{ and }\mathcal{L}^t=\{\varphi\in \Ind_H^G(\Hr{\sigma}) : \mbox{supp}(\varphi)\subseteq Ht\}.$$ By the definition of $\Ind_H^G(\Hr{\sigma})$ and since $G= H\sqcup Ht$, $\Ind_H^G(\Hr{\sigma})=\mathcal{L}\oplus \mathcal{L}^t.$
		Now, notice that $\mathcal{U}: \mathcal{L}\rightarrow \Hr{\sigma}: \varphi \mapsto \varphi(1_G)$ is a unitary operator and that 
		\begin{equation*}
		\begin{split}
		\sigma(h)\mathcal{U}\varphi=\sigma(h)\varphi(1_G)=\lb\Ind_H^G(\sigma)(h)\rb\varphi(1_G)=\mathcal{U}\lb \Ind_H^G(\sigma)(h)\rb\varphi\q \forall h\in H,\forall \varphi\in \mathcal{L}. 
		\end{split}
		\end{equation*}
		In particular, one has $(\Res_H^G(\Ind_H^G(\sigma)), \mathcal{L})\simeq (\sigma,\Hr{\sigma})$.	Similarly, notice that $\mathcal{U}^t: \mathcal{L}^t\rightarrow \Hr{\sigma}: \varphi \mapsto \varphi(t^{-1})$ is a unitary operator and that 
		\begin{equation*}
		\begin{split}
		\sigma^t(h)\mathcal{U}^t\varphi&=\sigma(tht^{-1})\varphi(t^{-1})=\varphi(t^{-1}th^{-1}t^{-1})=\varphi(h^{-1}t^{-1})\\
		&=\lb \Ind_H^G(\sigma)(h)\rb\varphi(t^{-1})=\mathcal{U}^t\lb \Ind_H^G(\sigma)(h)\rb\varphi\q \forall h\in H,\forall \varphi\in \mathcal{L}^t.
		\end{split}
		\end{equation*}
		This proves that $(\Res_H^G(\Ind_H^G(\sigma)), \mathcal{L}^t)\simeq (\sigma^t,\Hr{\sigma^t})$ and we obtain as desired that $\Res_H^G(\Ind_H^G(\sigma))\simeq\sigma \oplus \sigma^t$. By Frobenius reciprocity
		$${\rm I}(\Ind_H^G(\sigma),\Ind_H^G(\sigma))={\rm I}(\sigma,\Res_H^G\big(\Ind_H^G(\sigma)\big)$$
		this implies that $\Ind_H^G(\sigma)$ is irreducible (that is ${\rm I}(\Ind_H^G(\sigma),\Ind_H^G(\sigma))=1$) if and only if $\sigma\not \simeq \sigma^t$. Furthermore, in that case, we obtain that
		\begin{equation*}
		\begin{split}
		{\rm I}\big(\Ind_H^G(\sigma),\Ind_H^G(\sigma)^\tau\big)={\rm I}\big(\sigma,\Res_H^G\big(\Ind_H^G(\sigma)^\tau\big)\big)&={\rm I}\big(\sigma,\Res_H^G\big(\Ind_H^G(\sigma)\big)\big)\\
		&={\rm I}\big(\Ind_H^G(\sigma),\Ind_H^G(\sigma)\big)=1
		\end{split}
		\end{equation*}
		which proves that $\Ind_H^G(\sigma)\simeq\Ind_H^G(\sigma)^\tau$ and settles the first case.
		
		For the second case, Frobenius reciprocity ensures that $\sigma\simeq \sigma^t$ if and only if $\Ind_H^G(\sigma)$ is not irreducible. In that case, $\Ind_H^G(\Hr{\sigma})$ must split as a sum of two non-zero closed $G$-invariant subspaces $M$ and $M'$. On the other hand, since $\Res_H^G\big(\Ind_H^G(\sigma)\big)$ splits as a sum of two irreducible representations of $H$, and since every $G$-invariant subspace is $H$-invariant, $M$ and $M'$ do not admit any proper invariant subspaces. This proves that $\Ind_H^G(\sigma)\simeq \pi\oplus\pi'$ for some irreducible representations $\pi$ and $\pi'$ of $G$. On the other hand, since $\Res_H^G(\pi)=\Res_H^G(\pi^\tau)$, Frobenius reciprocity ensures that
		$${\rm I}(\Ind_H^G(\sigma),\pi)={\rm I}(\sigma,\Res_H^G(\pi))={\rm I}(\sigma,\Res_H^G(\pi^\tau))={\rm I}(\Ind_H^G(\sigma),\pi^\tau)$$
		for every irreducible representation $\pi$ of $G$. In particular, if $\pi\not \simeq \pi^\tau$, we obtain that $\Ind_H^G(\sigma)\simeq \pi \oplus \pi^\tau$. On the other hand, if $\pi\simeq \pi^\tau$, notice from Lemma \ref{Lemma equivalentce criterion between sigma et sigma g et pi et pitau} that $\Res_H^G(\pi)\simeq \sigma\oplus \sigma^t$. Hence, since $\sigma\simeq \sigma^t$, Frobenius reciprocity implies that
		\begin{equation*}
		\begin{split}
		{\rm I}(\Ind_H^G(\sigma),\pi)={\rm I}(\sigma,\Res_H^G(\pi))={\rm I}(\sigma,\sigma \oplus \sigma^t)>1
		\end{split}
		\end{equation*}
		which proves that $\Ind_H^G(\sigma)\simeq \pi \oplus \pi \simeq \pi \oplus \pi^\tau$.
	\end{proof}
	The first part of Theorem \ref{les rep dun group loc compact par rapport a celle d'un de ses sous groupes} follows from the above and the impossibility to have simultaneously $\pi\simeq \pi^\tau$ and $\sigma\simeq \sigma^t$. Indeed, if $\pi\simeq \pi^\tau$, Lemma \ref{Lemma equivalentce criterion between sigma et sigma g et pi et pitau} ensures that $\Res_H^G(\pi)\simeq \sigma\oplus \sigma^t$. However, if $\sigma \simeq \sigma^t$, Proposition \ref{prop des induite qui se split ou pas} ensures that $\Ind_H^G(\sigma)\simeq \pi \oplus \pi^\tau$. In particular, if those conditions were satisfied simultaneously one would obtain that $$\Res_H^G\big(\Ind_H^G(\sigma)\big)\simeq \sigma \oplus \sigma^t \oplus \sigma \oplus \sigma^t\simeq 4 \sigma.$$
	This is a contradiction since $\sigma$ is irreducible and Proposition \ref{prop des induite qui se split ou pas} ensures that $$\Res_H^G(\Ind_H^G(\sigma))\simeq\sigma \oplus \sigma^t\simeq 2\sigma.$$
	
	\noindent The following result completes the proof of Theorem \ref{les rep dun group loc compact par rapport a celle d'un de ses sous groupes}.
	\begin{lemma}\label{Uniformly admissible if and only if subgroup of index 2 is}
		Every irreducible representation $\pi$ of $G$ is a subrepresentation of $\Ind_H^G(\sigma)$ for some irreducible representation $\sigma$ of $H$ and every irreducible representation $\sigma$ of $H$ is a subrepresentation of $\Res_H^G(\pi)$ for some irreducible representation $\pi$ of $G$. 
	\end{lemma}
	\begin{proof}
		Let $\pi$ be an irreducible representation of $G$ and notice that
		$${\rm I}(\Ind_H^G\big(\Res_H^G(\pi)\big), \pi)={\rm I}(\Res_H^G(\pi),\Res_H^G(\pi))\geq1$$
		which proves that $\pi\leq \Ind_H^G\big(\Res_H^G(\pi)\big)$. If $\Res_H^G(\pi)$ is not irreducible, Lemma \ref{Lemma les rep de G se decompose ou pas en irred si H est d'indice 2} ensures that $\Res_H^G(\pi)\simeq \sigma\oplus \sigma^t$ for some irreducible representation $\sigma$ of $H$ and we have either that $\pi \leq \Ind_H^G(\sigma)$ or that $\pi \leq \Ind_H^G(\sigma^t)$. 
		
		Now, let $\sigma$ be an irreducible representation of $H$ and notice that
		$${\rm I}(\sigma,\Res_H^G\big(\Ind_H^G(\sigma)\big))={\rm I}(\Ind_H^G(\sigma), \Ind_H^G(\sigma))\geq1$$
		which proves that $\sigma\leq \Res_H^G\big(\Ind_H^G(\sigma)\big)$. If $\Ind_H^G(\pi)$ is not irreducible, Proposition \ref{prop des induite qui se split ou pas} ensures that $\Ind_H^G(\sigma)\simeq \pi\oplus \pi^\tau$ for some irreducible representation $\pi$ of $G$ and we have either that $\pi \leq \Ind_H^G(\sigma)$ or that $\pi \leq \Ind_H^G(\sigma^t)$.
	\end{proof}
	
	\section{Fell topology}\label{section Fell topology}
	Let $G$ be a second-countable locally compact group. The purpose of this section is to associate a topology to $\widetilde{G}$ and $\widehat{G}$. The topology that we consider here is the \tg{Fell topology}. This one comes directly from the relation of weak containment and can also be defined algebraically at the level of $C^*$-algebras. As we will see, the Fell topology carries a lot of information on the group and is well behaved with respect to natural operations such as direct sum, restriction, induction,... However, if $\pi$ is a representation of $G$, the Fell topology is not able to separate $\pi$ from $n\pi$ and does not separates the points of $\widehat{G}$ in full generality. In particular, the natural Borel structure associated to the Fell topology will usually not be suitable to support direct integrals of representations. We refer to Section \ref{section direct integral decompostion} below for more details. 
	
	We start with the notion of functions of positive type  and the relation of weak containment. 
	\begin{definition}\label{definition functions of positive type}
		Let $G$ be a locally compact group. A continuous function $\varphi : G \rightarrow \C$ is a \tg{function of positive type on $G$} if for all $n\in \N$, for every $c_1,...,c_n\in \C$ and for all $g_1,...,g_n\in G$ one has that
		$$\sum_{i=1}^n\sum_{j=1}^n c_i\overline{c_j}\varphi(g_{j}^{-1}g_i)\geq 0.$$ 
	\end{definition}
	This notion plays a central role in the theory of unitary representations due to the following example as well as Theorem \ref{Theorem GNS construction} below. 
	\begin{example}\label{Example function of positive type universal example}
		Let $G$ be a locally compact group and $\pi$ be a unitary representation of $G$. A function $\varphi:G\rightarrow \C$ is a \textbf{matrix coefficient} of $\pi$ if there exists a vector $\xi \in \Hr{\pi}$ such that 
		$$\varphi(g)=\prods{\pi(g)\xi}{\xi}_{\Hr{\pi}}\q \forall g\in G.$$
		Every such function is of positive type, see \cite[Proposition C.4.3]{BekkaHarpeValette2008}. 
	\end{example}
	A remarkable result due to  I. M.~Gelfand, M. A. ~Neimark and I. E.~Segal ensures that this example is universal and provides a correspondence between functions of positive type and unitary representations. 
	\begin{theorem}[\tg{GNS Construction} {\cite[Theorem C.4.10]{BekkaHarpeValette2008}}]\label{Theorem GNS construction}
		 Let $G$ be a locally compact group and $\varphi:G\rightarrow \C$ be a function of positive type. Then, there exists a triple $(\pi,\mathcal{H}, \xi)$ consisting of a unitary representation $\pi$ of $G$ with representation space $\mathcal{H}$ such that the space spanned by $\{\pi(g)\xi:g\in G\}$ is dense in $\mathcal{H}$ and $\varphi(g)=\prods{\pi(g)\xi}{\xi}_{\mathcal{H}}$. Furthermore, this triple is unique in the following sense: if $(\pi',\mathcal{H}',\xi')$ is another such triple, then there exists a Hilbert space isomorphism  $\mathcal{U}: \mathcal{H}\rightarrow \mathcal{H}'$
		intertwining $\pi$ and $\pi'$ and such that $\mathcal{U}\xi=\xi'$. 
	\end{theorem}
	 Let $\pi$ be a representation of $G$ and $\mathcal{S}\subseteq \widetilde{G}$ be a family of representations of $G$. A \tg{function of positive type} $\varphi:G\rightarrow \C$ is said to be \tg{associated} with $\pi$ when $\varphi(g)=\prods{\pi(g)\xi}{\xi}_{\Hr{\pi}}\qq \forall g\in G$ for some vector $\xi \in \hr{\pi}$. We say that $\varphi$ is \tg{associated} with $\mathcal{S}$ when it is associated with some $\rho\in \mathcal{S}$ (this definition makes sense since equivalent representations of $G$ are associated with the same functions of positive type).
	\begin{definition}
		Let $G$ be a locally compact group and $\mathcal{S}, \mathcal{P}$ be two subsets of $\widetilde{G}$. We say that $\mathcal{S}$ is \tg{weakly contained} in $\mathcal{P}$ which we denote by $\mathcal{S}\preceq \mathcal{P}$ if every function of positive type associated with $\mathcal{S}$ is the limit of a sum of functions of positive type associated with $\mathcal{P}$ for the topology of uniform convergence on compacta. We say that $\mathcal{S}$ and $\mathcal{P}$ are \tg{weakly equivalent}, which we denote by $\mathcal{S}\sim\mathcal{P}$, if $\mathcal{S}\preceq \mathcal{P}$ and $\mathcal{P}\preceq \mathcal{S}$. Similarly, given two representations $\sigma$ and $\pi$ of $G$, we say that $\sigma$ is weakly contained in $\pi$ if $\{\sigma\}\preceq \{\pi\}$ and we say that $\sigma$ is weakly equivalent to $\pi$ if $\{\sigma\}\preceq \{\pi\}$ and $\{\pi\}\preceq \{\sigma\}$. By abuse of notation, we denote these relations respectively by $\sigma\preceq \pi$ and $\sigma\sim \pi$.
	\end{definition}
	\begin{remark}
		If $\pi$ is irreducible, we recall that every vector $\xi$ of $\Hr{\pi}$ is cyclic, meaning that the linear span of $\{\pi(g)\xi:g\in G \}$ is dense inside $\Hr{\pi}$. Therefore, in order to prove that $\pi\preceq \mathcal{S}$ it is equivalent to prove that some function of positive type associated with $\pi$ is the  limit of a sum of functions of positive type associated with $\mathcal{S}$ for the topology of uniform convergence on compacta. 
	\end{remark}
	For every subset $\mathcal{S}\subseteq \widetilde{G}$, we define the closure of $\mathcal{S}$ in $\widetilde{G}$ as the set
	$$\cl(\mathcal{S})=\{\pi \in \widetilde{G}: \pi \preceq \mathcal{S}\}.$$
	The topology of $\widetilde{G}$ generated by these closed sets is called the \tg{Fell topology}. We now provide an explicit description of a basis of open sets for this topology. For an equivalence class of representation $\pi \in \widetilde{G}$, a finite sequence of functions of positive type $\varphi_1,...,\varphi_t$ associated with $\pi$, a compact subset $K$ of $G$ and a positive real $\varepsilon>0$ we let 
	\begin{equation*}\label{voisinage pour la topo de Fell}
		\mathcal{W}(\pi, \varphi_1,...,\varphi_t, K,\varepsilon)
	\end{equation*} be the set of all the equivalence classes of representations $\sigma\in \tilde{G}$ for which there exist sums of functions of positive type $\varphi_1',..,\varphi_t'$ associated with $\sigma$ such that  $$\modu{\varphi_i(k)-\varphi_i'(k)}<\varepsilon\q \forall k\in K,\qq \forall i=1,...,t.$$ 
	The sets $\mathcal{W}(\pi, \varphi_1,...,\varphi_t, K,\varepsilon)$ form a basis of the Fell topology on $\widetilde{G}$. 
	The \tg{Fell topology} on $\widehat{G}$ is defined as the restriction of the Fell topology on $\widetilde{G}$ provided by the inclusion $\widehat{G}\subseteq \widetilde{G}$. This topological space is in general far from being Hausdorff but the following result ensures nevertheless a certain amount of regularity. 	
	\begin{theorem}[{\cite[Proposition 1.D.1, 1.D.3 and 8.A.15]{BekkadelaHarpe2020}}]\label{theorem regularilte de la topologie de fell secodn countble bair}
		Let $G$ be a second-countable locally compact group and equip $\widehat{G}$ with the Fell topology. Then, $\widehat{G}$ is a second-countable locally quasi-compact Baire space.
	\end{theorem} 
	\begin{definition}\label{definition of support in the Fell topo}
		For every subset $\mathcal{S}\subseteq \widetilde{G}$, we define the \tg{support} $\Sp(\mathcal{S})$ of $\mathcal{S}$ as the subset of all those representations $\sigma \in\widehat{G}$ that are weakly contained in $\mathcal{S}$. The \textbf{support} of a representation $\pi$ of $G$ is the set $\Sp(\{\pi\})$ and is denoted by abuse of notation by $\Sp(\pi)$.
	\end{definition}
	\begin{remark}\label{remark the support is weakly equivalent to the set of rep in the fell topo}
		Let $\mathcal{S}\subseteq \widetilde{G}$. From the definition of the Fell topology, it is clear that $\Sp(\mathcal{S})$ is a closed subset of $\widehat{G}$. In fact, $\Sp(\mathcal{S})$ is the unique closed subset of $\widehat{G}$ that is weakly equivalent to $\mathcal{S}$, see \cite[Theorem 1.6]{Fell1960}.
	\end{remark}
	The following result shows that the Fell topology of the dual of a locally compact group contains a lot of information on the group itself. 
	\begin{theorem}\label{Theorem les differentes implication topo fell}
		Let $G$ be a second-countable locally compact group and consider the Fell topology on $\widehat{G}$. Then, the following hold:
		\begin{enumerate}[leftmargin=*,label=(\roman*)]
			\item\label{amenable} $G$ is \tg{amenable} if and only if $\Sp(\lambda_G)=\widehat{G}$ or equivalently $1_{\widehat{G}}\preceq\lambda_G$.
			\item\label{Compact} $G$ is \tg{compact} if and only if $\widehat{G}$ is discrete.
			\item\label{PropT} $G$ has \tg{Kazhdan's property (T)} if and only if $1_{\widehat{G}}$ is isolated in $\widehat{G}$.
			\item\label{T0} $G$ is \tg{Type {\rm I}} if and only if the Fell topology on $\widehat{G}$ is $T_0$.
			\item\label{T1} $G$ is \tg{CCR} if and only if the  Fell topology on $\widehat{G}$ is $T_1$.
		\end{enumerate} 
	\end{theorem}
	All those statement can be found in \cite{BekkadelaHarpe2020}, see Theorem 1.C.3 for \ref{amenable}, Proposition 1.D.7 and Remark 1.D.8 for \ref{Compact}, Remark 1.D.16 for \ref{PropT} and Theorem 1.D.13 for \ref{T0} and \ref{T1}. In addition, the Fell topology is well behaved with respect to natural operations. For instance, when $G$ is a locally compact group and $H$ is a closed subgroup, it is quite easy to see that the operation of restriction is a continuous map $\widetilde{G}\rightarrow \widetilde{H}: \pi \mapsto \Res_H^G(\pi)$ for the respective Fell topologies.  The same holds for the process of induction.
	\begin{theorem}[{\cite[Theorem 4.1]{Fell1962}}]\label{Theorem continity of induction}
		Let $G$ be a second-countable locally compact group and let $H$ be a closed subgroup of $G$, then the map $$\widetilde{H}\rightarrow \widetilde{G}:\pi\mapsto \Ind_{H}^G(\pi)$$ is continuous  for the respective Fell topologies. 
	\end{theorem}
	We now provide a description of the Fell topology of the unitary dual of some locally compact groups. The following provides an example where the Fell topology is super regular. 
	\begin{example}\label{Example duel of abelian group fell topo}
		Let $G$ be a second-countable abelian locally compact group. Schur's Lemma ensures that the irreducible representations of $G$ are all one-dimensional, so that the dual $\widehat{G}$ of $G$ is given by the set of unitary characters 
		$$\widehat{G}=\{\chiup: G\rightarrow \mathcal{U}(\C)\mid \chiup \mbox{ is a continuous group homomorphism}\}.$$
		Now, for every unitary character $\chiup$ of $G$, notice that every function of positive $\varphi: G\rightarrow \C$ associated to $\chiup$ is a scalar multiple of $\chiup$. The Fell topology is therefore given by uniform convergence on compacta of unitary characters. To be more precise, a sequence $(\chiup_n)_{n\in \N}$ of irreducible representations of $G$ converges to $\chiup\in \widehat{G}$ in the Fell topology if and only if the sequence of functions $(\chiup_n)_{n\in \N}$ converges uniformly on compacta to $\chiup$. In particular, $\widehat{G}$ is a Hausdorff second-countable locally compact topological space for the Fell topology. In fact, $\widehat{G}$ is itself a second-countable abelian locally compact group for the multiplication given by pointwise multiplication of unitary characters. Under the same process, the unitary dual of $\widehat{G}$ is isomorphic to $G$ as a locally compact group. This is the so called \tg{Pontryagin duality}. 
	\end{example}
	By contrast the following example shows that the dual frequently contains pathologies.
	\begin{example}\label{Example group ax+b Fell topology}
		Let $G$ be the affine group $\R\rtimes \R^*$ where $\R^*$ is the multiplicative group of real numbers and let $H=\{(r,1): r\in \R\}$. This subgroup $H$ is both normal and closed inside $G$. Furthermore, the dual of $G$ can easily be described by Mackey's theory of induction. To see this, notice that $\widehat{H}=\{\chiup_s: s\in \R\}$ where
		$$\chiup_s:H\rightarrow \C^* :(r,1) \mapsto e^{2\pi i sr}.$$
		Since $H$ is normal in $G$, notice that $G$ acts by conjugation on the unitary dual $\widehat{H}$. To be more precise, for every $g\in G$ and each irreducible representation $\chiup$ of $H$, we define the conjugate representation $\chiup^g$ by $$\chiup^g(h)=\chiup (gh g^{-1}) \qq \forall h\in H.$$ For each $t\in \R$, notice that $\chiup_s^{(0,t)}= \chiup_{st}$. In particular, there are exactly two $G$-orbits in $\widehat{H}$. The $G$-orbit $G\cdot 1_{\widehat{H}}=\{1_{\widehat{H}}\}$ of the trivial representation of $H$ and the $G$-orbit $G\cdot \chiup_1=\{\chiup_s : s\in \R-\{0\}\}$ of $\chiup_1$. For each $G$-orbit $\mathcal{O}$ in $\widehat{H}$, chose an element $\chiup \in \mathcal{O}$ and consider the subgroup $H_{\mathcal{O}}$ of all the elements $g$ of $G$ for which $\chiup^g\simeq \chiup$. The group $H_{\mathcal{O}}$ is closed, contains $H$ and does not depend on our choice of element $\chiup$ of $\mathcal{O}$. In our cases $H_{G\cdot  1_{\widehat{H}}}=G$ and $H_{G\cdot \chi_1}=H$. Now, \cite[Theorem 3.12]{Mackey1976theory} ensures that each irreducible representation $\pi$ of $G$ comes from an irreducible representation of $H_{\mathcal{O}}/H$ for some $G$-orbit $\mathcal{O}$ in the following way: if $\pi$ comes from $G\cdot 1_{\widehat{H}}$, it is the inflation of an irreducible representation of $G/H\simeq \R^*$ and if $\pi$ comes from $G\cdot \chiup_1$ it is isomorphic to $\Ind_H^G(\chiup_1)$. In particular, we have that
		$$\widehat{G}=\{\pi_t: t\in \R\}\sqcup \{\Ind_H^G(\chiup_1)\} $$
		where $\pi_t: \R\rtimes \R^*\rightarrow \C^*: (a,b)\mapsto b^{2\pi i t}$. We refer to \cite[pg. 194]{Mackey1976theory} for a similar discussion on the group $\R\rtimes \R^*_{+}$. Since each of the $\pi_t$ is a unitary characters of $G$, the map $\R\mapsto \widehat{G}-\{\Ind_H^G(\chiup_1)\}$ is an homeomorphism between $\R$ and the restriction of the Fell topology to $\widehat{G}-\{\Ind_H^G(\chiup_1)\}$. On the other hand, the closure of $G\cdot \chiup_1$ for the Fell topology of $\widehat{H}$ contains $1_{\widehat{H}}$. It follows from \cite[Theorem $4.3$]{Fell1962} that $\Ind_H^G(\chiup_1)$ weakly contains each of the $\pi_t$ (because $\Res^G_H(\pi_t)\sim1_{\widehat{H}}$). In particular, $\{\Ind_H^G(\chiup_1)\}$ is a dense point in $\widehat{G}$ for the Fell topology. We refer to \cite[pg. 263]{Fell1962} for a similar discussion on $\R\rtimes \R^*_{+}$.
	\end{example}
	Quite remarkably, the Fell topology can also be defined algebraically on the level of $C^*$-algebras. To be more precise, to every locally compact group is naturally associated an operator algebra encoding the topological dynamic of the group. The Fell topology of its dual can be derived from the inclusion of ideals of this algebraic structure. We give a precise version of this statement below. We start with some preliminaries concerning $C^*$-algebras. 
	\begin{definition}
		A \tg{Banach-algebra} $A$ is an associative complete normed algebra over the complex numbers whose norm $\norm{\cdot}{A}$ satisfies the followings: 
		\begin{enumerate}
			\item $\norm{\lambda a}{A}=\modu{\lambda}\norm{a}{A}$ for all $\lambda\in \C$ and $a\in A.$
			\item $\norm{ab}{A}\leq \norm{a}{A}\norm{b}{A}$ for all $a,b\in A$.
		\end{enumerate}
		A $*$-\tg{algebra} is an associative complex algebra $A$ equipped with a map $*:A\rightarrow A$ satisfying the following axioms:
		\begin{enumerate}
			\item $(ab)^*=b^*a^*$ for all $a,b\in A.$
			\item $(\lambda1_A)^*=\overline{\lambda}1_A$ for all $\lambda\in \C$.
			\item $a^{**}=a$ for all $a\in A.$
			\item $(a+b)^*=a^*+b^*$ for all $a,b\in A$. 
		\end{enumerate}
		A $C^*$-\tg{algebra} is a Banach $*$-algebra $A$ satisfying the $C^*$-\tg{identity} that is $$\norm{a^*a}{}=\norm{a}{}^2 \q \forall a\in A.$$
	\end{definition}
	We now provide two important examples.
	\begin{example}\label{example operators on Hilbert spaces are Cstar alegrbas}
		Let $\mathcal{H}$ be a Hilbert space and consider the associative algebra $\mathcal{L}(\mathcal{H})$ of bounded operators on $\mathcal{H}$. This is a $C^*$-algebra when equipped with the operator norm 
		$$\norm{V}{op(\mathcal{H})}\qq =\sup_{\xi \in \mathcal{H}: \norm{\xi}{\mathcal{H}}\leq 1} \norm{V\xi}{\mathcal{H}}$$ and with involution $*: \mathcal{L}(\mathcal{H})\rightarrow \mathcal{L}(\mathcal{H})$ given by adjunction.
	\end{example}
	\begin{example}\label{example compactoperators on Hilbert spaces are Cstar alegrbas}
		Let $\mathcal{H}$ be a Hilbert space. A bounded operator $V\in \mathcal{L}(\mathcal{H})$ is said to be \tg{compact} if the image of the unit ball of $\mathcal{H}$ under $V$ has compact closure in $\mathcal{H}$. It is not hard to prove that the set $\mathcal{K}(\mathcal{H})$ of compact operators on $\mathcal{H}$ is a self-adjoint two sided ideal of $\mathcal{L}(\mathcal{H})$ and that the operator norm limit of compact operators is itself a compact operator. In particular, $\mathcal{K}(\mathcal{H})$ is a $C^*$-algebra for the operator norm and the involution given by the adjunction. This algebra is even simple in the sense that the only closed two-sided ideals of $\mathcal{K}(\mathcal{H})$ are $\{0\}$ and $\mathcal{K}(\mathcal{H})$ itself.
	\end{example}
	Now, let $G$ be a locally compact group and $\mu$ be a left-invariant Haar measure of $G$. We consider the space $L^1(G)$ of complex valued measurable functions $\varphi:G\rightarrow \C$ that are integrable with respect to $\mu$ in the sense that $\int_G\modu{\varphi(g)}\diff \mu(g)<+\infty$. The convolution product of any two functions $\varphi, \psi\in L^1(G)$ is defined by 
	\begin{equation}\label{equation convolution product def}
	\varphi*\psi (g)=\int_{G}\varphi(h)\psi(h^{-1}g)\diff \mu(h)\q \forall g\in G.
	\end{equation}
	The space $L^1(G)$ is a Banach $*$-algebra for the product given by the convolution, the $L^1$-norm $$\norm{\varphi}{L^1(G)}=\int_G\modu{\varphi(g)}\diff \mu(g)$$ and the involution given by 
	\begin{equation}\label{equation the involution in C^*(G)}
	\varphi^*(g)=\Delta_G(g^{-1})\overline{\varphi(g^{-1})}\q \forall g\in G
	\end{equation}
	However, in general, this algebra does not satisfy the $C^*$-identity and is therefore not a $C^*$-algebra. For instance, we have the following result.
	\begin{lemma}
		Let $G$ be a non-trivial discrete group. Then, $\ell^1(G)$ is not a $C^*$-algebra for the convolution and the natural involution.
	\end{lemma}
	\begin{proof}
		For each $g\in G$, let $\delta_g$ be the Dirac function on $g$ and notice that for all $g_1,g_2\in G$ that
		$$\delta_{g_1}*\delta_{g_2}=\delta_{g_1g_2}.$$ Now, let $h\in G$ be a non-trivial element of $G$. A straightforward computation shows (even when $h$ has order $2$) that 
		$$\norm{\delta_{1_G}+i\delta_{h}+\delta_{h^2}}{\ell^1(G)}^2=\Big(\sum_{g\in G} \modu{\delta_{1_G}(g)+i\delta_{h}(g)+\delta_{h^2}(g)}\Big)^2=9.$$
		However, another computation shows that 
		\begin{equation*}
			\begin{split}
				(\delta_{1_G}+i\delta_{h}+\delta_{h^2})^**(\delta_{1_G}+i\delta_{h}+\delta_{h^2})&=(\delta_{1_G}-i\delta_{h^{-1}}+\delta_{h^{-2}})*(\delta_{1_G}+i\delta_{h}+\delta_{h^2})\\
				&=3 \delta_{1_G} + \delta_{h^{2}}+\delta_{h^{-2}}.
			\end{split}
		\end{equation*}
	It follows that 
	$$\norm{(\delta_{1_G}+i\delta_{h}+\delta_{h^2})^**(\delta_{1_G}+i\delta_{h}+\delta_{h^2})}{\ell^1(G)}=5.$$
	\end{proof}
	In order to obtain a $C^*$-algebra from $L^1(G)$ in full generality we therefore need to consider another norm. This one will take advantage from the fact that the $C^*$-identity is satisfied on $\mathcal{L}(\mathcal{H})$. We start with some preliminaries. 
	 \begin{definition}
	 	Let $A$ be a Banach $*$-algebra. A $*$-\tg{representation} of $A$ is a continuous homomorphism of algebras $\rho:A\rightarrow \mathcal{L}(\mathcal{\mathcal{H}})$ such that $\mathcal{H}$ is a separable Hilbert space and  $\rho(a^*)=\rho(a)^*$ for all $a\in A$. Such a representation is \tg{non-degenerate} if $\bigcap_{a\in A}\Ker(\rho(a))=\{0\}$ and is said to be \tg{irreducible} if the only invariant subspaces are $\{0\}$ and $\mathcal{H}$. Furthermore, two such representations $\rho_i: A\rightarrow\mathcal{L}(\mathcal{H}_i)$ ($i=1,2$) are \tg{equivalent} if there exists a unitary intertwining operator $V: \mathcal{H}_1 \rightarrow \mathcal{H}_2$ such that $\rho_1(a)V = V\rho_2(a)$ for all $a\in A$. We denote by $\widetilde{A}$ and $\widehat{A}$ the respective sets of equivalence classes of non-degenerate and of irreducible non-degenerate $*$-representations of $A$.
	 \end{definition} Now, let $G$ be a locally compact group. Each unitary representation $\pi:G\rightarrow \mathcal{U}(\Hr{\pi})$ gives rise to a non-degenerate $*$-representation  $$\tilde{\pi}:L^1(G)\rightarrow \mathcal{L}(\Hr{\pi})$$ defined by 
	\begin{equation*}\label{equation extention des represetations de G a L1G}
	\tilde{\pi}(\varphi)=\int_G\varphi(g)\pi(g)\diff \mu (g)\q \forall \varphi\in L^1(G).
	\end{equation*}
	In other words, for every $\xi \in \Hr{\pi}$, $\tilde{\pi}(\varphi)\xi$ is defined as the unique vector of $\Hr{\pi}$ satisfying that
	$$\prods{\tilde{\pi}(\varphi)\xi}{\eta}_{\Hr{\pi}}= \int_G\varphi(g)\prods{\pi(g)\xi}{\eta}_{\Hr{\pi}}\diff \mu (g)\q \forall \eta \in \Hr{\pi}.$$
	The correspondence $\pi\mapsto \tilde{\pi}$ sets up a bijection between the unitary representations of $G$ and the non-degenerate $*$-representation of $L^1(G)$ see \cite[Section 13.3]{Dixmier1977}. Furthermore, notice that this correspondence preserves both the notion of irreducibility and the respective relations of equivalence. On the other hand, each non-degenerate $*$-representation $\tilde{\pi}$ of $L^1(G)$ provides a semi-norm on $L^1(G)$ defined by $$\norm{\varphi}{\tilde{\pi}}=\norm{\tilde{\pi}(\varphi)}{op(\Hr{\tilde{\pi}})}.$$ The $C^*$-identity is satisfied for this semi-norm as for all $\varphi\in L^1(G)$ one has that
	$$\norm{\varphi^**\varphi}{\tilde{\pi}}= \norm{\tilde{\pi}(\varphi^**\varphi)}{op(\Hr{\tilde{\pi}})}=\norm{\big(\tilde{\pi}(\varphi)^*\big)\tilde{\pi}(\varphi)}{op(\Hr{\tilde{\pi}})}= \norm{\tilde{\pi}(\varphi)}{op(\Hr{\tilde{\pi}})}^2=\norm{\varphi}{\tilde{\pi}}^2.$$
	Now, consider the semi-norm on $L^1(G)$ defined by
	$$\norm{\varphi}{C^*}= \sup_{\tilde{\pi}\in \widetilde{L^1(G)}}\norm{\tilde{\pi}(\varphi)}{op(\Hr{\tilde{\pi}})} \q \forall \varphi \in L^1(G).$$
	Since the $*$-representation of $L^1(G)$ corresponding to the regular representation of $G$ is injective by \cite[Section 13.3.6]{Dixmier1977} $\norm{\cdot}{C^*}$ defines a norm satisfying the $C^*$-identity on $L^1(G)$. 
	\begin{definition}
		The completion of $L^1(G)$ with respect to $\norm{\cdot}{C^*}$ is a $C^*$-algebra called the \tg{maximal} $C^*$-algebra of $G$ that we denote by $C^*(G)$. 
	\end{definition}
	Notice by construction that every $*$-representation of $C^*(G)$ comes from a $*$-representation of $L^1(G)$ and is completely determined by its restriction to $L^1(G)$. In particular the correspondence 
	\begin{equation}\label{the one ot one correspondence between rep of G and of the C star completion}
	\widetilde{G}\rightarrow \widetilde{L^1(G)}: \pi\mapsto \tilde{\pi}
	\end{equation} 
	extends in a bijective correspondence between $\widetilde{G}$ and $\widetilde{C^*(G)}$ (the same holds for the bijective correspondence between $\widehat{G}$ and $\widehat{C^*(G)}$). In light of this remark, for every representation of $G$, we denote by $\tilde{\pi}$ the non-degenerate $*$-representation of $C^*(G)$ corresponding to $\pi$. 
	\begin{definition}
		The \tg{$C^*$-kernel} $C^*\Ker(\pi)$ of a representation $\pi$ of $G$ is the kernel of the $*$-representation $\tilde{\pi}$ inside $C^*(G)$. 
	\end{definition}
	\begin{remark}
		The $C^*$-kernel of a representation $\pi$ of $G$ is a closed two-sided ideal of $C^*(G)$. However, it does not coincide in general with the closure in $C^*(G)$ of the kernel of the restriction of $\tilde{\pi}$ to $L^1(G)$. For instance, $\widetilde{\lambda_G}$ is always faithful on $L^1(G)$ but $C^*\Ker(\lambda_G)=\{0\}$ only when $G$ is amenable. 
	\end{remark}
	A two-sided-ideal of a $C^*$-algebra is \tg{primitive} if it is the kernel of an irreducible $*$-representation. We denote by $\Prim(G)$ the set of primitive ideals of $C^*(G)$ and recall that set of primitive ideals of a $C^*$-algebra is naturally equipped with a topology called the \tg{Jacobson topology}. For this topology, the closure of a set $\mathcal{S}\subseteq \mbox{\rm Prim}(G)$ is the set of all ideals $I\in \mbox{\rm Prim}(G)$ such that $\bigcap_{J\in \mathcal{S}}J\subseteq I$. Equipped with this topology, $\Prim(G)$ is a $T_0$-separated topological space, see \cite[Chapter 3.1.]{Dixmier1977}. However, this space is not always $T_1$-separated and is therefore usually not Hausdorff. Now, notice that $C^*(G)$ and hence $\widehat{G}$ can be equipped with the pullback of the Jacobson topology of $\Prim(G)$ by the surjection $$\kappa:\widehat{G}\twoheadrightarrow \mbox{\rm Prim}(G): \pi \mapsto C^*\Ker(\tilde{\pi}).$$ 
	This topology on $\widehat{G}$ is called the \tg{hull-kernel topology}. Since the map $\kappa$ is not injective in full generality notice that the hull-kernel topology of $\widehat{G}$ is not guaranteed to be $T_0$ separated and hence does not separate the points of $\widehat{G}$. This provides an algebraic description of the Fell topology since, as shown in \cite{Fell1962}, the hull-kernel and the Fell topology coincide on $\widehat{G}$.
	
	A natural observation coming from this algebraic description of the Fell topology is that $\widehat{G}$ is a priory not $T_1$-separated. In particular, nothing ensures that the points of $\widehat{G}$ are closed. However, as we are going to see some representations of $G$ are closed in $\widehat{G}$. To be more precise, given a representation $\pi$ of $G$ we recall from the Example \ref{example compactoperators on Hilbert spaces are Cstar alegrbas} that the associative algebra $\mathcal{K}(\Hr{\pi})$ of compact operators is a simple $C^*$-algebra contained inside $\mathcal{L}(\Hr{\pi})$. Since the $*$-representation $\tilde{\pi}$ of $C^*(G)$ corresponding to $\pi$ ranges inside $\mathcal{L}(\Hr{\pi})$ it becomes natural to ask whether $\tilde{\pi}(C^*(G))$ intersects $\mathcal{K}(\Hr{\pi})$. The following result ensures that this intersection is either trivial or contains all the compact operators of $\Hr{\pi}$.
	\begin{proposition}[{\cite[Corollary $4.1.10$ and $4.1.11$]{Dixmier1977}}]\label{Proposotion Dixmier sur les rep CCR}
		Let $G$ be a locally compact group and let $\pi$ be an irreducible representation of $G$ such that $\tilde{\pi}(C^*(G))\cap \mathcal{K}(\Hr{\pi})\not=0$. Then, the following holds:
		\begin{enumerate}[label=(\roman*)]
			\item $\mathcal{K}(\Hr{\pi})\subseteq \tilde{\pi}(C^*(G))$ and, if $\mathcal{K}(\Hr{\pi})=\tilde{\pi}(C^*(G))$, one has that $C^*\Ker(\pi)$ is a maximal closed two-sided ideal of $C^*(G)$.
			\item Every irreducible representation $\pi'$ of $G$ with $C^*\Ker(\pi)= C^*\Ker(\pi')$ is unitary equivalent to $\pi$. 
		\end{enumerate}
	\end{proposition}
	This result motivates the following definition.
	\begin{definition}
		Let $G$ be a locally compact group. An irreducible representation $\pi$ of $G$ is a \tg{completely continuous representation} or \tg{CCR} for short if $\tilde{\pi}(C^*(G))= \mathcal{K}(\Hr{\pi}).$ A group $G$ is said to be \tg{CCR} if all its irreducible representations is CCR. 
	\end{definition}
	The following result follows directly from Proposition \ref{Proposotion Dixmier sur les rep CCR} and from the definition of the Fell topology in terms of inclusion of ideals in $C^*(G)$. 
	\begin{lemma}\label{lemma Dixmier sur les rep CCR sont fermee}
		Let $G$ be a locally compact group and let $\pi$ be CCR irreducible representation of $G$. Then, $\{\pi\}$ is closed in $\widehat{G}$.
	\end{lemma}
	In particular, if $G$ is a CCR locally compact group, every point of $\widehat{G}$ is closed. Theorem \ref{Theorem les differentes implication topo fell} ensures that $\widehat{G}$ is $T_1$ when $G$ is CCR so that CCR groups are of Type {\rm I}. In fact, Type {\rm I} locally compact groups are exactly the locally compact groups all of whose non-trivial irreducible representations $\pi$ satisfy that $\mathcal{K}(\Hr{\pi})\subseteq \tilde{\pi}(C^*(G))$, see \cite[Theorems $4.3.7$ and $9.1$]{Dixmier1977}.
	\begin{remark}
		Some locally compact groups are of Type {\rm I} but are not CCR. The affine group $\R\rtimes \R^*$ considered in the Example \ref{Example group ax+b Fell topology} provides a concrete example. To see this, observe the Fell topology of its unitary dual is $T_0$-separated but can be not $T_1$-separated as it  contains a dense point. It follows from Theorem \ref{Theorem les differentes implication topo fell} that this group is of Type {\rm I} but is not CCR. 
	\end{remark}
	
	We end this section by explaining how the Fell topology of $\widetilde{G}$ can be derived from the hull-kernel topology of $\widehat{G}$. Given a topological space $X$, the family $\mathcal{C}(X)$ of closed subsets of $X$, can be equipped with a topology called the \tg{inner topology} derived from $X$. This topology is defined as follows. For each finite set $\{U_1,U_2,...,U_n\}$ of non-empty open subsets of $X$, let 
	$$U(U_1,U_2,...U_n)=\{C\in \mathcal{C}(X): C\cap U_i\mbox{ is not empty }\forall i=1,...,n\}.$$ A subset of $\mathcal{C}(X)$ is open for the inner topology derived from $X$ if it is a union of sets of the form $U(U_1,...,U_n)$. Now, consider a locally compact group $G$. For every subset $\mathcal{S}\subseteq \widetilde{G}$ we define the closure of $\mathcal{S}$ as the set of all the representations $\pi$ of $G$ whose support $\Sp(\pi)$ belongs to the closure of $\{\Sp(\sigma): \sigma\in \mathcal{S}\}$ with respect to the inner topology  derived from $\mathcal{C}(\widehat{G})$.  This topology is called the \tg{inner-hull kernel topology} of $\widetilde{G}$ and coincides with the Fell topology, see \cite{Fell1962}.
	
	\section{Direct integral of representations}\label{section direct integral decompostion}
	
	For compact groups, each representation decomposes in an essentially unique way as a direct sum of irreducible representations. On the other hand, this result does not extend to general locally compact groups. However, a similar decomposition can be achieved relying on the more general notion of direct integral of representations. The purpose of this section is to introduce this concept. We begin the section with a disguised example to provide some intuition. The following example shows that the direct integral decomposition of the regular representation of a second-countable abelian locally compact group is provided by the Fourier transform.
	\begin{example}\label{Example decomposition en integrale direct d'irred de la reguliere}
		Let $G$ be a second-countable abelian locally compact group. For a concrete example with a non-compact group we advise the reader to think of $\Z$ or $\R$ and for a concrete example with a compact group we advise the reader to think of the circle group $\mathbb{T}$ (however, in that case, the direct integral will become a direct sum). We recall from Example \ref{Example duel of abelian group fell topo} that the unitary dual $\widehat{G}$ consists of set of unitary characters 
		$$\widehat{G}=\{\chiup: G\rightarrow \mathcal{U}(\C)\mid \chiup \mbox{ is a continuous group homomorphism}\}$$ and that it is a second-countable abelian locally compact group for the Fell topology and the multiplication given by pointwise multiplication of unitary characters. Let $\mu_G$ be a Haar measure on $G$ and let $C_c(G)$ be the set of complex valued continuous compactly supported function on $G$. By \cite[Theorem A.7.6]{BekkadelaHarpe2020}, there exists a renormalisation of the Haar measure $\mu_{\widehat{G}}$ on $\widehat{G}$ for which the Fourier transform defined by 
		$$\mathcal{F}(\varphi)(\chiup)=\int_G \chiup(g)\varphi(g)\diff\mu_G(g) \q \forall \varphi\in C_c(G),\qq \forall \chiup \in \widehat{G}$$
		extends to an isometry $\mathcal{F}:L^2(G,\mu_G)\rightarrow L^2(\widehat{G},\mu_{\widehat{G}})$. 
		As we will see below, the Hilbert space $L^2(\widehat{G},\mu_{\widehat{G}})$ can naturally be seen as a ``direct integral'' 
		$$L^2(\widehat{G},\mu_{\widehat{G}})=\int_{\widehat{G}}^\oplus \C\diff \mu_{\widehat{G}}(\chiup)$$ of a constant field of Hilbert spaces $\chiup\mapsto\C$ indexed by $\widehat{G}$ (because each $\chiup$ has dimension $1$). We can therefore define a ``direct integral'' of representations
		$$\pi=\int_{\widehat{G}}^\oplus \chiup \diff \mu_{\widehat{G}}(\chiup)$$
		on $L^2(\widehat{G},\mu_{\widehat{G}})$ by
		$$(\pi(g)\xi)(\chiup)=\chiup(g) \xi(\chiup)\q \forall \xi\in L^2(\widehat{G},\mu_{\widehat{G}}), \qq \chiup\in \widehat{G},\qq g\in G.$$ 
		Now, notice for all $\varphi\in L^2(G,\mu_G)$ and for every $g\in G$ that
		$$\mathcal{F}(\lambda_G(g)\varphi)=\pi(g)\mathcal{F}(\varphi).$$
		so that the Fourier transform is a unitary equivalence between the regular representation of $G$ and the direct integral $\pi$. In particular, this provides a direct integral decomposition of $\lambda_G$ into irreducible representations as 
		$$\lambda_G\simeq \int_{\widehat{G}}^\oplus \chiup \diff \mu_{\widehat{G}}(\chiup).$$
		This is the canonical decomposition of the regular representation of a second-countable abelian group into a direct integral of irreducible representations. Furthermore, notice that if $G$ is compact, $\widehat{G}$ is discrete and $\mu_{\widehat{G}}$ is the counting point measure, so that this direct integral decomposition coincides with the direct sum decomposition provided by Peter-Weyl's theorem \cite[Theorem 18.0.4]{Tao2014}.
	\end{example}
	We now formalise the concepts of direct integral of Hilbert spaces and representations. Our formalism is based on \cite{Dixmier1977} and \cite{BekkadelaHarpe2020} and we refer to these books for details. 
	
	Consider a measured space $(\Omega,\mathcal{B}, \nu)$ with a $\sigma$-finite measure $\nu$.
	\begin{definition}
		A \textbf{field of Hilbert spaces} over $\Omega$  is a family $(\mathcal{H}_\omega)_{\omega \in \Omega}$, where $\mathcal{H}_\omega$ is a Hilbert space for each $\omega \in \Omega$. The elements in $\prod_{\omega\in \Omega} \mathcal{H}_\omega$ are called the \textbf{vector fields} over $\Omega$. 
	\end{definition}
	A \tg{fundamental family of measurable vector fields} of $(\mathcal{H}_\omega)_{\omega \in \Omega}$ is a sequence $(e^n)_{n\in \N}$ of vector fields over $\Omega$ with the following properties:
	\begin{enumerate}[label=(\roman*)]
		\item the map $\Omega\rightarrow \C : \omega\mapsto \prods{e^n_\omega}{e^m_\omega}_{\mathcal{H}_\omega}$ is measurable $\forall n,m\in \N$; 
		\item the linear span of the $\{e^n_\omega:n\in \N\}$ is dense in $\mathcal{H}_\omega$ for all $\omega\in \Omega$.
	\end{enumerate}  
	Given a fundamental family of measurable vector fields $(e^n)_{n\in \N}$, a vector field $\xi$ over $\Omega$ is said to be \tg{measurable} if the function $$\Omega\rightarrow \C: \omega\mapsto \prods{\xi_\omega}{e^n_\omega}_{\mathcal{H}_\omega}$$ is measurable for all $n\in \N$. In addition, two measurable vector fields are said to be \textbf{equivalent} if they coincide $\nu$-almost everywhere. From now, we abusively write ``measurable vector fields'' instead of ``equivalence classes of measurable vector fields'' for this equivalence relation. The set $\mathcal{M}$ of measurable vector fields form a complex vector space. Moreover, for every $\xi,\eta\in \mathcal{M}$ the function $ \omega\mapsto \prods{\xi_\omega}{\eta_\omega}_{\mathcal{H}_\omega}$ is measurable. The pair $((\mathcal{H}_\omega)_{\omega \in \Omega},\mathcal{M})$ is a called a \tg{measurable field of Hilbert spaces} over $\Omega$. By abuse of notation, we usually say that $\omega \mapsto \mathcal{H}_\omega$ is a measurable field of Hilbert spaces over $\Omega$. 
	A measurable vector field $\xi$ over $\Omega$ is said to be \tg{square-integrable} if it satisfies that $$\int_\Omega \norm{\xi_\omega}{\mathcal{H}_\omega}^2\diff \nu(\omega)< \infty.$$ 
	We denote by $\mathcal{H}$ the vector space of square-integrable vector fields and equip it with an inner product defined by
	$$\prods{\xi}{\eta}_\mathcal{H}=\int_\Omega \prods{\xi_\omega}{\eta_\omega}_{\mathcal{H}_\omega}\diff \nu(\omega)\q \forall \xi,\eta\in \mathcal{H}.$$ 
	\begin{definition}
		The Hilbert space $(\mathcal{H}, \prods{\cdot}{\cdot}_{\mathcal{H}})$ is the \tg{direct integral} of $\omega\mapsto \mathcal{H}_\omega$ and is denoted by 
		$$\mathcal{H}=\int_\Omega^\oplus \mathcal{H}_\omega\diff \nu (\omega).$$
	\end{definition} 
	The following example shows that the notion of direct integral is nothing more than a generalisation of the notion of direct sum. 
	\begin{example}
		Let $\Omega$ be a countable set, let $\nu$ be the counting point measure and consider any field of Hilbert spaces $\omega\mapsto  \mathcal{H}_\omega$ on $\Omega$. Any sequence  $(e^n)_{n\in \N}$ of vector fields over $\Omega$ where $\{e^n_\omega: n \in \N\}$ is dense in $\mathcal{H}_\omega$ provides a fundamental family of measurable vector fields. Given such a family, every vector field is measurable and 
		$$\int^\oplus_\Omega \mathcal{H}_\omega \diff \nu( \omega )\simeq \bigoplus_{\omega \in \Omega}\mathcal{H}_\omega.$$
	\end{example}
	\begin{example}
		Let $(\Omega,\mathcal{B}, \nu)$ be a measured space with a $\sigma$-finite measure $\nu$, $\mathcal{H}_\aleph$ be the standard Hilbert space of dimension $\aleph\leq \aleph_0$ and consider the constant field of Hilbert spaces $\omega \mapsto \mathcal{H}_\aleph$. Choose an orthonormal basis $\{\xi_n:n\leq \aleph\}$ of  $\mathcal{H}_\aleph$ and define a fundamental family of measurable vector fields $(e^n)_{n\in \N}$ by
		$$e^n_\omega= \begin{cases}
		\qq\xi_n\qq&\forall n\leq \aleph,\\
		\qq 0 &\forall n> \aleph 
		\end{cases}\q \forall \omega \in \Omega.$$
		Then $\int_\Omega^\oplus\mathcal{H}_\omega\diff \nu (\omega)$ coincides with the Hilbert space 
		$$L^2(\Omega,\nu, \mathcal{H}_\aleph)=\Big\{\varphi: \Omega \rightarrow \mathcal{H}_\aleph\Big\lvert\int_\Omega \norm{\varphi(\omega)}{\mathcal{H}_\aleph}^2\diff \nu (\omega)<+\infty\Big\}$$
		with inner product defined by 
		$$\prods{\varphi}{\psi}_{L^2(\Omega,\nu, \mathcal{H}_\aleph)}=\int_\Omega \prods{\varphi(\omega)}{\psi(\omega)}_{\mathcal{H}_\aleph}\diff \nu(\omega).$$
	\end{example}
	\begin{definition}
		A \tg{measurable field of operators} over $\Omega$ is a map $$\Omega\rightarrow \mathcal{L}(\mathcal{H}_\omega): \omega \mapsto T_\omega$$ such that $$\Omega\rightarrow \C: \omega\mapsto\prods{T_\omega\xi_\omega}{\eta_\omega}_{\mathcal{H}_\omega}$$ is measurable for all $\xi,\eta \in \mathcal{H}$.
	\end{definition}
	Now, let $\mathcal{H}=\int^\oplus_\Omega \mathcal{H}_\omega \diff \nu( \omega )$ be a direct integral of Hilbert spaces over $\Omega$. Given a measurable field of operators $\omega\rightarrow T_\omega$ over $\Omega$ such that the map $$\Omega\rightarrow \C: \omega\mapsto\norm{T_\omega}{\mathcal{L}(\mathcal{H}_\omega)}$$ is $\nu$-essentially bounded, we define the \textbf{direct integral of operators} $\int_\Omega^\oplus T_\omega \diff \nu (\omega)$ as the operator over $\mathcal{H}$ given by
	$$\bigg(\int_\Omega^\oplus T_\alpha \diff \nu (\alpha)\xi\bigg)({\omega})=T_\omega \xi_\omega\q \forall \xi \in \mathcal{H}, \qq \forall \omega\in \Omega.$$
	An operator $T\in \mathcal{L}(\mathcal{H})$ of that form is said to be \tg{decomposable} over $\Omega$.
	\begin{remark}
		The operator norm $\norm{\int_\Omega^\oplus T_\omega \diff \nu (\omega)}{\mathcal{L}(\mathcal{H})}$ coincides with the $\nu$-essential supremum of the function $\omega\mapsto \norm{T_\omega}{\mathcal{L}(\mathcal{H}_\omega)}$. 
	\end{remark}
	\begin{example}
		Let $(\Omega,\mathcal{B}, \nu)$ be a measure space with a $\sigma$-finite measure $\nu$ and let $\omega \mapsto \mathcal{H}_\omega$ be a measurable field of Hilbert spaces over $\Omega$. To each function $\varphi\in L^\infty (\Omega, \nu)$ we can associate a decomposable operator $T_\varphi$ on $\mathcal{H}=\int^\oplus_\Omega \mathcal{H}_\omega \diff \nu (\omega)$ defined by 
		$$T_\varphi=\int_\Omega^\oplus \varphi(\omega )\id_{\mathcal{H}_\omega}\diff \nu (\omega).$$
		In other words, for all $\xi \in \mathcal{H}$ and for every $\omega \in \Omega$, one has 
		$$(T_\varphi\xi)(\omega)=\varphi(\omega)\xi(\omega).$$
		An operator $T\in \mathcal{L}(\mathcal{H})$ of the form is called \tg{diagonalisable}.
	\end{example}
	
	We come back to the main goal of this chapter that is the decomposition of representations as a direct integral of representations. Let $G$ be a locally compact group. A \tg{measurable field of representations} $((\pi_\omega,\mathcal{H}_\omega))_{\omega\in \Omega}$ of $G$ is a set of representations $(\pi_\omega, \mathcal{H}_\omega)$ such that the map $\omega \mapsto\pi_\omega(g)$ is a measurable field of operators over $\Omega$ for all $g\in G$. By abuse of notations, we denote such a measurable field by $\omega \mapsto \pi_\omega$. Given a measurable field of representations $\omega \mapsto \pi_\omega$ of $G$, we define a unitary representation $\pi$ on the Hilbert space $\mathcal{H}=\int_\Omega^\oplus \mathcal{H}_\omega \diff \nu (\omega)$ by 
	$$\pi(g)=\int_\Omega^\oplus\pi_\omega(g)\diff \nu(\omega)\q \forall g\in G.$$
	\begin{definition}
		The representation $(\pi,\mathcal{H})$ is the \tg{direct integral} of the measurable field of representations $\omega \mapsto \pi_\omega$ and is denoted by $\int_\Omega^\oplus\pi_\omega\diff \nu(\omega)$.
	\end{definition} 
	\begin{remark}\label{remark that enables one to consider the measurable fields of Hilbert spaces as measiravle on the Mackey-Borel}
		Direct integrals of representations can always be reduced to a direct sum of direct integrals of constant fields of Hilbert spaces. To be more precise, given a measurable field of representations $((\pi_\omega,\mathcal{H}_\omega))_{\omega\in \Omega}$ of $G$ and a cardinal $\aleph\leq \aleph_0$, the set $\Omega_\aleph=\{\omega\in \Omega: \dim(\pi_\omega)=\aleph\}$ is a measurable subset of $\Omega$ \cite[Proposition 1.G.2]{BekkadelaHarpe2020}. In particular, each direct integral $\pi=\int_\Omega^\oplus \pi_\omega\diff \nu(\omega)$ admits a canonical decomposition as a direct sum 
		$$\pi=\bigoplus_{\aleph\leq \aleph_0} \int_{\Omega_\aleph}^\oplus \pi_\omega\diff \nu(\omega)$$
		with representation space 
		$$\bigoplus_{\aleph\leq \aleph_0} \int_{\Omega_\aleph}^\oplus \mathcal{H}_\omega\diff \nu(\omega)\simeq \bigoplus_{\aleph\leq \aleph_0} L^2(\Omega_\aleph, \nu, \mathcal{H}_\aleph).$$
	\end{remark}
	Let $G$ be a locally compact group and $\pi$ be a representation of $G$. A direct integral $\int_\Omega^\oplus\pi_\omega\diff \nu(\omega)$ of representations of $G$ is a \tg{direct integral decomposition} of $\pi$ if 
	$$\pi \simeq \int_\Omega^\oplus\pi_\omega\diff \nu(\omega).$$
	This concept of direct integral decomposition can easily be seen as a generalisation of the concept of direct sum decomposition. Furthermore, such a decomposition is well behaved with respect to natural operations. For instance, if $H\leq G$ is a closed subgroup, it follows from the definition of direct integral of representations that 
	$$\Res_H^G(\pi)\simeq \int_{\Omega}^\oplus\Res_H^G(\pi_\omega) \diff \nu (\omega)$$
	when $\int_{\Omega}^\oplus\pi_\omega \diff \nu (\omega)$ is a direct integral decomposition of $\pi$. The same holds for the process of induction.
	\begin{theorem}[{\cite[Theorem 10.1]{Mackey1952induced}}]
		Let $G$ be a locally compact group $G$, $H$ be a closed subgroup and $\sigma\simeq \int_{\Omega}^\oplus\sigma_\omega \diff \nu (\omega)$ be a direct integral decomposition of a representation $\sigma$ of $H$. Then, we have that $$\Ind_H^G(\sigma)\simeq \int_{\Omega}^\oplus\Ind_H^G(\sigma_\omega) \diff \nu (\omega).$$
	\end{theorem}
	In general, the representations of a second-countable locally compact groups admit several decompositions as direct integrals of representations. The following theorem provides a way to obtain such decompositions but requires some preliminaries. We recall that the commutant of  $\mathcal{S}\subseteq \mathcal{L}(\mathcal{H})$ is the subalgebra of $\mathcal{L}(\mathcal{H})$ defined by
	$$\mathcal{S}'=\{T\in \mathcal{L}(\mathcal{H}): TS=ST\q \forall S\in \mathcal{S}\}$$  
	and that $\pi(G)'$ is the algebra $\Hom_G(\pi,\pi)$ of intertwining operators of $\pi$. 
	\begin{theorem}[{\cite[Theorem 1.G.7]{BekkadelaHarpe2020}}]\label{Theorem decomposition en integrale direct}
		Let $G$ be a second-countable locally compact group and let $\pi$ be a representation of $G$ on a separable Hilbert space $\mathcal{H}$. Let $\mathcal{Z}$ be an abelian subalgebra of $\pi(G)'$ stable by double commutant. Then, there exists a standard Borel space $\Omega$, a $\sigma$-finite measure $\nu$ on $\Omega$, a measurable field of Hilbert space $\omega\mapsto  \mathcal{H}_\omega$ a measurable field of representation $\omega \mapsto \pi_\omega$ of $G$ with representation spaces $\mathcal{H}_\omega$ and an isomorphism of Hilbert spaces 
		$$\mathcal{U}: \mathcal{H}\rightarrow \int_\Omega^\oplus \mathcal{H}_\omega\qq \diff \nu(\omega)$$
		with the following properties :
		\begin{enumerate}[label=(\roman*)]
			\item $\mathcal{U}\pi(g)\mathcal{U}^{-1}=\int_\Omega^\oplus \pi_\omega(g)\diff \nu ( \omega)\q \forall g\in G.$
			\item $\mathcal{U}\mathcal{Z}\mathcal{U}^{-1}$ is the algebra of diagonalisable operators on $ \int_\Omega^\oplus \mathcal{H}_\omega\diff \nu ( \omega).$
		\end{enumerate}
	\end{theorem}
	It is clear that such a direct integral decomposition is not unique as for instance if two measures $\nu$ and $\nu'$ on $\Omega$ are equivalent, the direct integrals $\int_\Omega^\oplus\pi_\omega \diff \nu(\omega)$ and $\int_\Omega^\oplus\pi_\omega \diff \nu'(\omega)$ are unitary equivalent.
	However, {\cite[Theorem 2.7]{Mackey1976theory}} ensures that on a measure-theoretic point of view, the triple $(\Omega,\nu,(\pi_\omega)_{\omega\in \Omega})$ is essentially unique given $\mathcal{Z}$. To be more precise for any other triple  $(\Omega',\nu',(\pi_{\omega'})_{\omega'\in \Omega'})$ with the above properties, there exist null Borel sets $\Omega_0$ and $\Omega_0'$ of $(\Omega,\nu)$ and $(\Omega',\nu')$ respectively and a Borel isomorphism $\Phi: \Omega\rightarrow \Omega'$ such that $\nu\circ \Phi^{-1}$ and $\nu'$ have the same null sets and $\pi_{\omega}\simeq \pi_{\Phi(\omega)}$ for all $\omega\in \Omega-\Omega_0$. On the other hand, for any representation $\pi$ of $G$ several choices of $\mathcal{Z}$ are usually possible. The following two important choices are omnipresent and show that every representation of a locally compact group admits a decomposition as a direct integral decomposition of factors of $G$ as well as decompositions into irreducible representations: 
	\begin{enumerate}[label=(\roman*)]
		\item Take $\mathcal{Z}$ to be the center of $\pi(G)'$ that is $\pi(G)'\cap \pi(G)''$.
		\item Take $\mathcal{Z}$ to be a maximal abelian subalgebra of $\pi(G)'$. 
	\end{enumerate}
	In the first case, the representations $\pi_\omega$ appearing in the direct decomposition given by Theorem \ref{Theorem decomposition en integrale direct} are $\nu$-almost all factors and are mutually disjoint. Furthermore, this decomposition is canonical in the sense that  $\pi(G)'$ has a unique centre. In the second case, the representations $\pi_\omega$ the representations are $\nu$-almost all irreducible \cite[Proposition 1.G.8]{BekkadelaHarpe2020}.  
	However, as there exists usually several maximal abelian subalgebras of $\pi(G)'$, there is no canonical choice of decomposition in a direct integral of irreducible representations. In fact, Theorem \ref{theorem non-uniqueness of direct integral decompostiion} below ensures that this direct integral decomposition is never unique for non-Type {\rm I} groups. 
	
	Notice, that the measure spaces appearing in each of the direct integral decomposition are a priori different. It would be desirable to base those direct integral decompositions on a common Borel space and it is natural to ask whether $\widetilde{G}$ and $\widehat{G}$ could be equipped with an appropriate Borel structure to provide such Borel spaces. However, the answer to this question involves fairly delicate measure-theoretic arguments and is only well behaved for Type {\rm I} groups, see Theorem \ref{Glimm's theorem}. As highlighted in Section \ref{section Fell topology}, the Borel structure coming from the Fell topology is not suitable to support direct integral decompositions since it does not separate a representation from any of its multiples and sometimes does not even separate the elements of $\widehat{G}$. We now provide a description of a finer Borel structure on $\widetilde{G}$ and $\widehat{G}$ called the Mackey-Borel structure. For every cardinal $\aleph\leq\aleph_0$, we fix a Hilbert space $\mathcal{H}_{\aleph}$ of dimension $\aleph$. Since the construction is similar for both $\widetilde{G}$ and $\widehat{G}$ we present a uniform description. For every $\aleph\leq \aleph_0$ we let $\mathcal{R}_\aleph$ be a set of representations of $G$ with representation space $\mathcal{H}_\aleph$. The two examples of interest in this thesis are:
	\begin{enumerate}
		\item the set $\Rep(G,\mathcal{H}_\aleph)$ of representations of $G$ on the space $\mathcal{H}_\aleph$. 
		\item the set $\Irr(G,\mathcal{H}_\aleph)$ of irreducible representations of $G$ on the space $\mathcal{H}_\aleph$.  
	\end{enumerate}
	We equip $\mathcal{R}_\aleph$ with the smallest $\sigma$-algebra for which the maps
	\begin{equation*}\label{equation les fonctions definissant la Mackey-Borel structure}
	\mathcal{R}_\aleph\rightarrow \C:\pi \mapsto\prods{\pi(g)\xi}{\eta}_{\mathcal{H}_\aleph}
	\end{equation*}
	are measurable for all $g\in G$, $\xi,\eta \in \mathcal{H}_\aleph$. We consider the direct union 
	$$\mathcal{R}=\bigsqcup_{\aleph\leq \aleph_0}\mathcal{R}_\aleph$$
	and equip this one with the sum $\sigma$-algebra $\mathcal{B}$. In other words, a subset $B\subseteq \mathcal{R}$ belongs to $\mathcal{B}$ if and only if $B\cap \mathcal{R}_\aleph$ is a Borel set of all $\aleph\leq \aleph_0$. Notice that for our respective choices of $\mathcal{R}_\aleph$, the quotient of $\mathcal{R}$ by the equivalence relation $\mathcal{R}/\simeq$ coincides as a set respectively with $\widetilde{G}$ and $\widehat{G}$. The respective quotient Borel structures are the \tg{Mackey-Borel structures} of $\widetilde{G}$ and $\widehat{G}$.
	\begin{remark}
		Since each measurable field of Hilbert spaces is essentially equivalent to a direct sum of constant measurable fields on standard Borel spaces by \cite[Proposition 1.G.2]{BekkadelaHarpe2020}, we notice for each measurable field of representations $\omega \mapsto \pi_\omega$ that the map 
		$$\Omega \rightarrow \widetilde{G}: \omega \mapsto \lb \pi_\omega \rb$$ is Borel measurable for the Mackey-Borel structure on $\widetilde{G}$.
	\end{remark}
	The Mackey-Borel structure of $\widehat{G}$ is always finer than the Borel structure coming from the Fell topology. Furthermore, by contrast with the Borel structure coming from the Fell topology, the Mackey-Borel structure always separate the points. A venerable result due to Glimm shows that both structures coincide exactly when the Fell topology is able to distinguish the elements of $\widehat{G}$. This provides the following deep equivalence.
	\begin{theorem}[Glimm's theorem {\cite[Theorem 8.F.3]{BekkadelaHarpe2020}}]\label{Glimm's theorem}
		Let $G$ be a $\sigma$-compact locally compact group. Then, the following are equivalent:
		\begin{enumerate}[leftmargin=*,label=(\roman*)]
			\item $G$ is type \tg{\rm I}.
			\item The map $\kappa : \widehat{G}\twoheadrightarrow \mbox{\rm Prim}(G)$ is a homeomorphism for the Fell topology. 
			\item The relation of equivalence and weak-equivalence coincide on $\widehat{G}$.
			\item The Mackey-Borel structure on $\widehat{G}$ is countably separated. 
			\item The Mackey-Borel structure on $\widehat{G}$ is standard. 
			\item The Mackey-Borel structure on $\widehat{G}$ coincides with the Borel structure coming from the Fell topology. 
		\end{enumerate}
	\end{theorem}
	The following result shows that the decomposition into a direct integral of irreducible representations is well behaved in that case. 
	\begin{theorem}[{\cite[Theorem 6.D.7]{BekkadelaHarpe2020}}]\label{Theorem de decomposition des representation de type I en integral direct}
		Let $G$ be a Type {\rm I} second-countable locally compact group and let $\pi$ be a representation of $G$. There exist a subset $\mathcal{N}\subseteq \N_0\cup\{\infty\}$ of extended positive integers and a family $(\nu_n)_{n\in \mathcal{N}}$ of positive measures with disjoint supports on $\widehat{G}$ such that 
		$$\pi\simeq \bigoplus_{n\in \mathcal{N}} n \int_{\widehat{G}}\sigma \diff \nu_n(\sigma).$$
		Furthermore, this decomposition is unique in the sense that $\mathcal{N}$ and the measure class of the $\nu_n$ is uniquely determined by $\pi$.
	\end{theorem}
	For non-Type {\rm I} groups on the other hand, this uniqueness does not hold.
	\begin{theorem}[{\cite[Theorem 1.G.11]{BekkadelaHarpe2020}}]\label{theorem non-uniqueness of direct integral decompostiion}
		Let $G$ be a non-Type {\rm I} locally compact group. Then, there exists a representation $\pi$ of $G$ admitting two direct integral decompositions $$\pi=\int^\oplus_\Omega \pi_\omega\diff \nu(\omega)= \int^\oplus_{\Omega'} \pi_{\omega'}\diff \nu'(\omega')$$ such that $\pi_\omega$ is disjoint from $\pi_{\omega'}$ for all $(\omega,\omega')\in \Omega\times \Omega'$.
	\end{theorem} 

	\chapter{Irreducibly represented Lie groups}\label{Chapter Lie groups are irred faithfull}
	\section{Introduction and main results}
	A locally compact group is \tg{irreducibly represented} if it admits a faithful irreducible unitary representation. A classical problem initially formulated by W.~Burnside \cite[Appendix F]{Burnside1955} in the context of finite groups and extending naturally to the context of locally compact groups is given by the following:
	\begin{problem*}
		Find a characterisation of the locally compact groups that are irreducibly represented purely in terms of their algebraic structure.
	\end{problem*}
	For instance, it can be deduced from Schur's Lemma that a finite abelian group is irreducibly represented if and only if it is cyclic. To the present day, algebraic characterisations were obtained for finite groups \cite{Gaschutz1954} and more generally for countable discrete groups \cite{BekkadeleHarpe2008} but not beyond this framework. Quite remarkably and despite the fact that the arguments used in the finite and in the countable case are very different from one another, the algebraic characterisation of the irreducibly represented groups is the same in both setups. To be more precise, a countable discrete group $G$ is irreducibly represented if and only if every every finite normal subgroup of $G$ contained in the subgroup generated by all the minimal non-trivial finite normal subgroups of $G$ is generated by a single $G$-conjugacy class. The purpose of this chapter is to provide an algebraic characterisation of the real Lie groups that are irreducibly represented among the connected real nilpotent Lie groups and the Hausdorff connected components of the identity of the $\R$-points of linear algebraic groups defined over $\R$. Similarly to the case of countable groups our reasoning relies on direct integral decompositions and on measure theory. However, by contrast with the case of discrete groups, our reasoning relies heavily on Type {\rm I} groups. 
	
	The statement of our results requires some preliminaries. Let $G$ be a locally compact group and $H\leq G$ be a closed normal subgroup. We recall that $G$ acts by conjugation on the unitary dual $\widehat{H}$ by defining the conjugate representation $\sigma^g$ of a representation $\sigma$ by $\sigma^g(h)=\sigma(ghg^{-1})$ $\forall h\in H$. If the subgroup $H$ is of Type {\rm I}, it is said to be \tg{regularly embedded} in $G$ if the quotient Borel structure of $\widehat{H}$ under this action is countably separated. For instance, closed central subgroups and compact normal subgroups are always regularly embedded. A representation $\sigma$ of $H$ is said to be $G$\tg{-faithful} if $\bigcap_{g \in G} \Ker(\sigma^g)=\{1_H\}$, or equivalently, if $\Ind_H^G(\sigma)$ is a faithful representation of $G$. The following theorems provide relations between the faithful irreducible representation of $G$ and the $G$-faithful irreducible representations of some of its closed Type {\rm I} regularly embedded normal subgroups.
	\begin{theoremletter}\label{theorem letter locally comapct groups equivalence of faitfull}
		Let $G$ be an irreducibly represented locally compact group and $N\leq G$ be a closed Type {\rm I} regularly embedded normal subgroup of $G$. Then, $N$ admits a $G$-faithful irreducible representation. 
	\end{theoremletter}
	\begin{theoremletter}\label{theorem letter locally comapct groups equivalence of faitfull part II}
		Let $G$ be a locally compact group, $K\times H\leq G$ be a closed subgroup such that $K$ is a compact normal subgroup, $H\leq G$ is a closed normal abelian regularly embedded subgroup 
		and both $K$ and $H$ admit a $G$-faithful irreducible representation. Then, there exists an irreducible representation $\pi$ of $G$ that is weakly contained in the regular representation $\lambda_G$ and whose kernel intersects $K\times H$ trivially. 
	\end{theoremletter}
	The strategy adopted in these notes to obtain an algebraic characterisation of certain families of Lie groups that are irreducibly represented is to rely on these theorems by identifying suitable closed normal subgroups of $G$. This approach is shown to be successful below for the Hausdorff connected component of the identity of the $\R$-points of a linear algebraic group defined over $\R$ and for connected nilpotent real Lie groups (which are ``almost algebraic" due to \cite[Theorem 3.6.3]{Varadarajan2013}).
	
	Let $G$ be a connected Lie group. Inspired by \cite{BekkadeleHarpe2008}, we make the following definitions. A \tg{semi-simple Lie foot} $K$ of $G$ is a minimal non-trivial compact connected semi-simple normal subgroup of $G$ with trivial center (minimal in the sense that the only subgroup of $K$ that is non-trivial connected and normal in $G$ is $K$ itself). We let $\mathcal{K}_G$ be the set of semi-simple Lie feet of $G$ and define the \tg{semi-simple Lie socle} $\mathcal{S}_s(G)$ as the closed subgroup of $G$ generated by the elements of $\mathcal{K}_G$. Similarly, an \tg{abelian Lie foot} $V$ of $G$ is a minimal non-trivial closed connected abelian normal subgroup of $G$ intersecting the center of $G$ trivially (minimal in the sense that the only subgroup of $V$ that is non-trivial closed connected and normal in $G$ is $V$ itself). We let $\mathcal{A}_G$ be the collection of abelian Lie feet of $G$ and define the \tg{abelian Lie socle} $\mathcal{S}_a(G)$ as the closed subgroup generated by the set of abelian Lie feet and the center $\mathcal{Z}(G)$ of $G$. In order to use Theorems \ref{theorem letter locally comapct groups equivalence of faitfull} and \ref{theorem letter locally comapct groups equivalence of faitfull part II}, we prove the following list of assertions:
	\begin{enumerate}[leftmargin=*, label=(\roman*)]\label{assertions intro Lie groups general} 
		\item\label{item the subgroup Sa times Ss} (Proposition \ref{Proposition la structure du group engendrer par Sa et Ss}) The closed subgroup of $G$ generated by $\mathcal{S}_s(G)$ and $\mathcal{S}_a(G)$ is a topologically characteristic subgroup of $G$ isomorphic to $\mathcal{S}_s(G)\times \mathcal{S}_a(G)$ (with the natural inclusions). 
		\item\label{item the subgroup Ss} (Proposition \ref{proposition the subgroup spanned by the semi simple or abelian feet is a closed connectedted normal ss subgroup} and Corollary \ref{corollary the semi-simple Lie group is irreducibly represented}) The semi-simple Lie socle $\mathcal{S}_s(G)$ is an irreducibly represented compact connected topologically characteristic semi-simple subgroup of $G$.
		\item\label{item the subgroup Sa} (Corollary \ref{Corollary the abelian Lie socle is a closed abelian topologically characteristic subgroup of G}) The abelian Lie socle $\mathcal{S}_a(G)$ is a closed connected topologically characteristic abelian subgroup of $G$.
		\item\label{item le noyau des rep faiblement reguliere intersectes les socles} (Proposition \ref{Proposition le noyau des rep faiblement reguliere intersectes les socles})
		Let $\pi$ be a non-faithful representation of $G$ that is weakly contained in the regular representation $\lambda_G$. Then, $\Ker(\pi)$ intersects either $\mathcal{S}_s(G)$ or $\mathcal{S}_a(G)$ non-trivially.
		\item\label{item alternative du centre} (Lemma \ref{lemma property Harpe}) If $G$ is irreducibly represented, its center $\mathcal{Z}(G)$ is either discrete and countable or connected and isomorphic to the circle group.
	\end{enumerate}
	In addition, if the center $\mathcal{Z}(G)$ is discrete and countable or connected and isomorphic to the circle group, the following assertions hold:
	\begin{enumerate}[leftmargin=*, label=(\roman*)]\label{assertions intro Lie groups under centre alternative} 
		\setcounter{enumi}{5}
		\item\label{item Sa G if G is nilpotent} (Lemma \ref{lemma in nilpotenet groups minimal connected subgroups are central}) If $G$ is nilpotent, we have that $\mathcal{S}_s(G)=\{1_G\}$, that $\mathcal{S}_a(G)=\mathcal{Z}(G)$ and that $\mathcal{Z}(G)$ is isomorphic to the circle group.
		\item\label{item if G is algebraix the abelian socle is reg embedded} (Proposition \ref{proposition the ab elian Lie socle is regularily embedded}) If $G$ is the Hausdorff connected component of the identity of the $\R$-points of a linear algebraic group defined over $\R$, then $\mathcal{S}_a(G)$ is regularly embedded.
		\item\label{item structure of the abelian lie socle}(Theorem \ref{Theorem LG est characteristic et isomoprhe a un espace vectoriel}) There exists a closed topologically characteristic subgroup $\mathcal{V}(G)$ of $G$ isomorphic as Lie groups to some finite dimensional real vector space such that the abelian Lie socle $\mathcal{S}_a(G)$ is isomorphic to $$\mathcal{V}(G)\times \mathcal{Z}(G)\mbox{ (with the natural inclusions)}.$$ The group $\mathcal{V}(G)$ is called the \tg{linear socle} of $G$.
		\item\label{item the moment Sa G admits a G faitfull representaiton }(Lemma \ref{lemma Sa has a G faith iff V (G) and Z are G irrep}) There exists a $G$-faithful irreducible representation of $\mathcal{S}_a(G)$ if and only if $\mathcal{V}(G)$ has a $G$-faithful irreducible representation and $\mathcal{Z}(G)$ is irreducibly represented.
		\item\label{item the moment LG admits a G faitfull representaiton }(Proposition \ref{Proposition, the lienar socle is G irred rep if if coadjoint orbit generates}) The group $\mathcal{V}(G)$ has a $G$-faithful irreducible representation if and only if the linear dual $\mathcal{V}(G)^*$ is spanned as a real vector space by a single coadjoint $G$-orbit. 
		\item\label{item the moment cenre faitfull representaiton }(\cite[Corollary $1.3$ and Corollary $1.9$]{CapracedelaHarpe2020}) The center $\mathcal{Z}(G)$ is not irreducibly represented if and only it is discrete and contains a subgroup isomorphic to $C_p\times C_p$ for some prime $p$.
	\end{enumerate}
	The following two theorems follow.
	\begin{theoremletter}\label{Corollary letter Nilpotentnt groups}
		Let $G$ be a connected real nilpotent Lie group. Then, $G$ is irreducibly represented if and only if its center $\mathcal{Z}(G)$ is isomorphic as a Lie group to the circle.  
	\end{theoremletter}
	\begin{proof}
		Suppose that $G$ is irreducibly represented. The assertions \ref{item alternative du centre} and \ref{item Sa G if G is nilpotent} ensure that $\mathcal{Z}(G)$ is isomorphic as a Lie group to the circle. Now, suppose that $\mathcal{Z}(G)$ is isomorphic as a Lie group to the circle. The assertion \ref{item Sa G if G is nilpotent} ensures that $S_s(G)=\{1_G\}$ and that $\mathcal{S}_a(G)$ is a central (hence regularly embedded) subgroup of $G$ isomorphic as a Lie group to the circle. Since the circle group is irreducibly represented, Theorem \ref{theorem letter locally comapct groups equivalence of faitfull part II} ensures the existence of an irreducible representation $\pi$ of $G$ that is weakly contained in the regular representation $\lambda_G$ and whose kernel intersects $S_s(G)\times S_a(G)$ trivially. The assertion \ref{item le noyau des rep faiblement reguliere intersectes les socles} ensures that such a representation of $G$ is faithful. 
	\end{proof}
	\begin{theoremletter}\label{Corollary letter amenable lie groups}
		Let $G$ be the Hausdorff connected component of the identity of the $\R$-points of a linear algebraic group defined over $\R$. Then, the following assertions are equivalent: 
		\begin{enumerate}
			\item\label{item irred rep} $G$ is irreducibly represented.
			\item\label{item algebraic charact of irred rep} The following assertions hold:
			\begin{enumerate}[leftmargin=*]
				\item The center $\mathcal{Z}(G)$ is isomorphic to the circle group or is discrete countable and does not contain a subgroup isomorphic to $C_p\times C_p$ for any prime $p$.
				\item The linear dual $\mathcal{V}(G)^*$ is spanned as a real vector space by a single coadjoint $G$-orbit.
			\end{enumerate}
			\item\label{item irred rep by weakly regular} There exists a faithful irreducible representation of $G$ that is weakly contained in the regular representation $\lambda_G$.
		\end{enumerate} 
	\end{theoremletter}
	\begin{proof}
		The implication \ref{item irred rep by weakly regular} $\Rightarrow$ \ref{item irred rep} follows from the definition of irreducibly represented groups. We now show that \ref{item irred rep} $\Rightarrow$ \ref{item algebraic charact of irred rep} $\Rightarrow$ \ref{item irred rep by weakly regular}. Suppose that $G$ is irreducibly represented. The assertion \ref{item alternative du centre} ensures that the center $\mathcal{Z}(G)$ is either is either discrete and countable or connected and isomorphic as a Lie group to the circle $\mathbb{T}$. Furthermore, the assertions \ref{item the subgroup Sa} and \ref{item if G is algebraix the abelian socle is reg embedded} ensures that $\mathcal{S}_a(G)$ is a closed regularly embedded normal subgroup of $G$.
		In particular, Theorem \ref{theorem letter locally comapct groups equivalence of faitfull} ensures that it admits a $G$-faithful irreducible representation. Furthermore, the assertion \ref{item the moment Sa G admits a G faitfull representaiton } ensures that this happens if and only if $\mathcal{V}(G)$ has a $G$-faithful irreducible representation and $\mathcal{Z}(G)$ is irreducibly represented. The statement \ref{item algebraic charact of irred rep} hence follows from the assertions \ref{item the moment LG admits a G faitfull representaiton } and \ref{item the moment cenre faitfull representaiton }. Now, the assertions \ref{item the subgroup Sa times Ss}, \ref{item the subgroup Ss}, \ref{item the subgroup Sa} and \ref{item the moment LG admits a G faitfull representaiton } ensure that one satisfies the hypotheses of Theorem \ref{theorem letter locally comapct groups equivalence of faitfull part II} for $K=\mathcal{S}_s(G)$ and $H=\mathcal{S}_a(G)$. This implies the existence of an irreducible representation $\pi$ of $G$ that is weakly contained in the regular representation of $G$ and whose kernel intersects $\mathcal{S}_s(G)\times\mathcal{S}_a(G)$ trivially. The assertion \ref{item le noyau des rep faiblement reguliere intersectes les socles} ensures that such a representation is faithful. 
	\end{proof}

	\subsection*{structure of the chapter}
	We start by giving a concrete application of Theorem \ref{Corollary letter amenable lie groups} in the Section \ref{section lie concrete application}. We provide a proof of  Theorems \ref{theorem letter locally comapct groups equivalence of faitfull} and \ref{theorem letter locally comapct groups equivalence of faitfull part II} in Section \ref{section proof of theorem I and II Lie} and we prove the assertions \ref{item the subgroup Sa times Ss} to \ref{item the moment cenre faitfull representaiton } in Section \ref{section Lie proof of all the assertiosn}.

	\section{A concrete application of Theorem \ref{Corollary letter amenable lie groups}}\label{section lie concrete application}
	
	The purpose of this section is to provide a concrete application of Theorem \ref{Corollary letter amenable lie groups} to determine if a certain Lie group is irreducibly represented. 
	
	Choose a field $\mathbb{K} \in \{\R, \C\}$ and two non-negative integers $m, n\in \N$. Let $V=\mathbb{K}^n$ and $\Gl(V)^{\circ}$ be the connected component of $\id_V$ in $\Gl(V)$. More precisely, $\Gl(V)^\circ$ is the set of invertible linear maps $A:V\rightarrow V$ with positive determinant if $\mathbb{K}=\R$ and it is $\Gl(V)$ if $\mathbb{K}=\C$. We consider the diagonal action of $\Gl(V)$ on $V^m$ and the associated semi-direct product $$G_{n,m}^+=V^m \rtimes \Gl(V)^\circ.$$
	
	\begin{lemma}
		$G_{n,m}^+$ is irreducibly represented if and only if $m\leq n$.
	\end{lemma}
	\begin{proof}
		By the definition of $G_{n,m}^+$, notice that every closed normal subgroup of $G_{n,m}^+$ is contained in $V^m$. From there, it is not hard to see that $G_{n,m}^+$ is center free. In particular, Theorem \ref{Corollary letter amenable lie groups} ensures that $G_{n,m}^+$ is irreducibly represented if and only if $\mathcal{V}(G_{n,m}^+)^*$ is spanned as a real vector space by a single coadjoint $G$-orbit. In addition, as $G_{n,m}^+$ is center free, one has that $\mathcal{V}(G_{n,m}^+)=\mathcal{S}_a(G_{n,m}^+)$. Now, let us show that $\mathcal{S}_a(G_{n,m}^+)=V^m$. Since $G_{n,m}^+$ is center free, the elements of $\mathcal{A}_{G_{n,m}^+}$ are just the non-trivial closed connected abelian normal subgroups of $G$ that are minimal for these properties. Since every closed normal subgroup of $G_{n,m}^+$ is contained in $V^m$ we obtain that $\mathcal{S}_a(G_{n,m}^+)\subseteq V^m$. To prove the other inclusion, we let $V_i$ be the $i^{th}$ copy of $V$ in the direct product $V^m=V\times \cdots \times V$. Each of these $V_i$ is a non-trivial closed connected abelian normal subgroup of $G$ and is minimal for these properties since the $G_{n,m}^+$-conjugacy class of an element $v\in V_i$ is the set 
		$$\Gl(V)^\circ v=\{Av:a\in \Gl(V)^\circ\}=V_i-\{0\}.$$ Since $V^m$ is generated by the $V_i$, we obtain as desired that $\mathcal{S}_a(G_{n,m}^+)= V^m.$ Now, let $f\in \mathcal{S}_a(G_{n,m}^+)^*$. Due to the above identification, $f$ can be realised as an $m$-tuple $(f_1,f_2,..., f_m)$ where $f_i\in V_i^*$ for each $i\in \{1,...,m\}$. Under this identification, the coadjoint $G_{n,m}^+$-orbit of $f$ is the set
		$$G_{n,m}^+\cdot f=\{(A^T f_1,A^T f_2,...,A^T f_m)  : A\in \Gl(V)^\circ\}.$$ 
		If $m\leq n$, we can take $f_1,...,f_m$ to be linearly independent so that $G_{n,m}^+\cdot f=(V^m)^*-\{0\}$ and $\mathcal{S}_a(G_{n,m}^+)^*$ is spanned by a single coadjoint $G_{n,m}^+$-orbit so that $G_{n,m}^+$ is irreducibly represented. On the other hand, if $m> n$, the vector space spanned by $G_{n,m}^+\cdot f$ has dimension at most $n^2$. Since $V^m$ has dimension $nm>n^2$, $\mathcal{S}_a(G_{n,m}^+)$ can not be spanned by a single coadjoint $G_{n,m}^+$-orbit. It follows that $G_{n,m}^+$ is not irreducibly represented. 
	\end{proof}

	\section{Regular embedding and irreducibly represented groups}\label{section proof of theorem I and II Lie}
	The purpose of this section is to provide a way to study the existence of a faithful irreducible representation of locally compact group $G$ from the existence of $G$-faithful irreducible representations of its closed regularly embedded normal subgroups. In particular, we prove Theorems \ref{theorem letter locally comapct groups equivalence of faitfull} and \ref{theorem letter locally comapct groups equivalence of faitfull part II}. 
	\subsection{Regular embedding}
	We start with an introduction on the concept of regular embedding playing a central role in both of these results. Let $G$ be a locally compact group and $H\leq G$ be a closed normal subgroup. As $H$ is normal in $G$, we have an action by conjugation of $G$ on $H$. This one provides an action of $G$ on $\widetilde{H}$ and on $\widehat{H}$. To be more precise, for every representation $\sigma$ of $H$ and for all $g\in G$, the \tg{conjugate representation} $\sigma^g$ of $\sigma$ by $g$ is the representation of $H$ with representation space $\Hr{\sigma}$ defined by $$\sigma^g(h)=\sigma(ghg^{-1}) \q \forall h\in H.$$ It is not hard to see that this action by conjugation preserves the unitary equivalences as well as the notion of irreducibility so that $G$ also acts by conjugation on both $\widetilde{H}$ and $\widehat{H}$. On the other hand, we recall from Section \ref{section Fell topology} that $\widehat{H}$ carries a natural Borel structure called the Mackey-Borel structure. Furthermore, we recall from Theorem \ref{Glimm's theorem} that this Borel structure is countably separated and coincides with the Borel structure generated by the Fell topology exactly when $H$ is of Type {\rm I}. 
	\begin{definition}\label{definition regular embedding}
		A closed Type {\rm I} normal subgroup $H\leq G$ is said to be \tg{regularly embedded} in $G$ if the quotient Borel structure of $\widehat{H}$ under the $G$-action by conjugation is countably separated in the sense that there exists a countable family $\mathcal{B}$ of $G$-invariant Borel sets such that each $G$-orbit in $\widehat{H}$ is the intersection of the elements of $\mathcal{B}$ in which it is contained.
	\end{definition}
	\begin{remark}\label{remark locally closed orbits implies regular embedding}
		Notice in light of Theorem \ref{theorem regularilte de la topologie de fell secodn countble bair} that a sufficient condition for $H$ to be regularly embedded in $G$ is that all the $G$-orbits of $\widehat{H}$ are locally closed (open in their closure) for the Fell topology. For certain families of groups such as abelian groups, these properties are even equivalent, see Proposition \ref{proposition for abelian groups regular embedding is equivalent to the fact taht every orbits is locally closed} below.
	\end{remark}
	The following results provide important examples of regularly embedded subgroups that will be used later in these notes.
	\begin{lemma}\label{lemma compact groups are regularly embedded}
		Let $G$ be a locally compact group and $K\leq G$ be a compact normal subgroup. Then, every $G$-orbit in $\widehat{K}$ is both open and closed for the Fell topology. In particular, $K$ is regularly embedded in $G$.
	\end{lemma}
	\begin{proof}
		As $G$ is compact, Theorem \ref{Theorem les differentes implication topo fell}\ref{Compact} ensures that $\widehat{K}$ is a countable discrete space for the Fell topology. In particular, every $G$-orbit in $\widehat{K}$ is both open and closed for the Fell topology. It follows from the Remark \ref{remark locally closed orbits implies regular embedding} that $K$ is regularly embedded in $G$.
	\end{proof}
	\begin{lemma}\label{lemma a central subgroup is regularily embedded}
		Let $G$ be a connected locally compact group and $\mathcal{Z}\leq G$ be a closed central subgroup. Then, $\mathcal{Z}$ is regularly embedded. 
	\end{lemma}
	\begin{proof}
		Since $\mathcal{Z}$ is central, it is an abelian normal subgroup of $G$ and the action of $G$ on $\widehat{\mathcal{Z}}$ is trivial. In particular, $\widehat{\mathcal{Z}}$ is a second-countable locally compact abelian group and Theorem \ref{Glimm's theorem} 
		ensures that the Mackey-Borel structure coincide with Borel structure generated by the Fell topology on $\widehat{\mathcal{Z}}$. The result therefore follows from the fact that $\widehat{\mathcal{Z}}$ has a countable basis for the topology which separated the points. 
	\end{proof}
	The following result due to Dixmier provides one of the key step of the proof of Theorem \ref{Corollary letter Nilpotentnt groups}.
	\begin{theorem}[{\cite[Theorem 1]{Dixmier1959}}]\label{Theorem dixmier nilp reg embedded}
		Let $G$ be a connected nilpotent Lie group and $H\leq G$ be a closed connected abelian normal subgroup. Then, $H$ is regularly embedded in $G$. 
	\end{theorem}
	In light of the Remark \ref{remark locally closed orbits implies regular embedding}, it is natural to ask whether the property of being regularly embedded is equivalent to the property of having locally closed orbits in the unitary dual. The next result shows that this is the case for abelian subgroups.
	\begin{proposition}\label{proposition for abelian groups regular embedding is equivalent to the fact taht every orbits is locally closed}
		Let $G$ be a locally compact group and $H\leq G$ be a closed normal abelian subgroup. Then, the following are equivalent:
		\begin{enumerate}[label=(\roman*)]
			\item $H$ is regularly embedded in $G$.
			\item Every $G$-orbit is locally closed in $\widehat{H}$.
		\end{enumerate}
	\end{proposition}
	This follows from the following classical result.
	\begin{theorem}[{\cite[Theorem $2.1.14$]{Zimmer1984}}]\label{Theorem de zimmer etre reg embedded c'est avoir des orbites loca closed}
		Let $G$ be a locally compact group acting continuously on a complete separable metrisable space $\Omega$. Then, the following are equivalent :
		\begin{enumerate}[label=(\roman*)]
			\item The quotient Borel structure of $\Omega/G$ is countably separated.
			\item All $G$-orbits are locally closed in $\Omega$.
			\item For every $\omega \in \Omega$, the natural map $G/\Fix_G(\omega)\rightarrow G\cdot \omega$ is an homeomorphism, where the $G$-orbit $G\cdot \omega$ has the relative topology as a subset of $\Omega$.
		\end{enumerate}
	\end{theorem} 
	\begin{proof}[Proof of Proposition \ref{proposition for abelian groups regular embedding is equivalent to the fact taht every orbits is locally closed}]
		We recall that all the locally compact groups considered in these notes are second-countable and hence are complete separable metrisable space. Furthermore, since $H$ is an abelian locally compact group, Pontryagin duality ensures that its unitary dual $\widehat{H}$ is also a locally compact group. In particular, in order to apply Theorem \ref{Theorem de zimmer etre reg embedded c'est avoir des orbites loca closed}, we are left to prove that the conjugation map $$\phi: G\times \widehat{H} \rightarrow \widehat{H}: (g,\chiup)\mapsto \chiup^g$$ is continuous. We recall that a basis of neighbourhoods of $\chiup\in \widehat{H}$ given by sets of the form 
		\begin{equation*}
			U_{\widehat{H}}(\chiup, K, \varepsilon)=\{\rho\in \widehat{H}\qq : \qq \modu{\chiup(t)-\rho(t)}< \varepsilon\qq \forall t\in K\}
		\end{equation*} 
		with $K\subseteq \widehat{H}$ compact and $\varepsilon>0$. In particular, it suffices to show that $\phi^{-1}(U_{\widehat{H}}(\chiup, K, \varepsilon))$ is open in $G\times \widehat{H}$ for all $(g,\chiup)\in G\times \widehat{H}$ and every $\varepsilon>0$. Let $\chiup\in \widehat{H}$, $\varepsilon>0$ and $(g_0,\chiup_0)\in \phi^{-1}(U_{\widehat{H}}(\chiup, K, \varepsilon))$. By continuity, the map
		$$K \rightarrow \R : k \mapsto \modu{\chiup(k)-\chiup_0^{g_0}(k)}$$ 
		attains its maximum $\lambda\in \R_{\geq 0}$. In particular, one has that 
		\begin{equation}\label{equation regular embedding in compactly gen 1}
			\modu{\chiup(k)-\chiup_0^{g_0}(k)}\leq \lambda\q \forall k\in K,
		\end{equation}
		and the definition of $U_{\widehat{H}}(\chiup, K, \varepsilon)$ ensures $\lambda<\varepsilon$. Now, choose a strictly positive real number $\delta<\frac{\varepsilon-\lambda}{2}$. By continuity of the map
		$$G\times K \rightarrow \R : (g,k)\rightarrow \modu{ \chiup_0^{g_0}(k)-\chiup_0^g(k)}$$
		and by compactness of $K$, there exists an open neighbourhood $V$ of $g_0\in G$ with compact closure $\overline{V}$ for which one has that 
		\begin{equation}\label{equation regular embedding in compactly gen 2}
			\modu{\chiup_0^{g_0}(k)-\chiup_0^g(k)}< \delta\q \forall k\in K, \qq \forall g\in\overline{V}.
		\end{equation}
		By continuity of the conjugation map $c_G:G\times H \rightarrow H : (g,h)\mapsto ghg^{-1}$, the set $c_G(\overline{V},K)$ is compact in $H$.
		In particular, $V\times U_{\widehat{H}}(\chiup_0,c_G(\overline{V},K),\delta)$ is an open neighbourhood of $(g_0,\chiup_0)$ in $G\times \widehat{H}$.
		Now, notice that $$V\times  U_{\widehat{H}}(\chiup_0,c_G(\overline{V},K),\delta) \subseteq \phi^{-1}(U_{\widehat{H}}(\chiup, K, \varepsilon)).$$ Indeed, for all $ \rho\in U_{\widehat{H}}(\chiup_0,c_G(\overline{V},K),\delta)$, $ g\in V$ and every $k\in K$ we have that
		\begin{equation}\label{equation regular embedding in compactly gen 3}
			\modu{\chiup_0^{g}(k)-\rho^g(k)}=\modu{\chiup_0(gkg^{-1})-\rho(gkg^{-1})} < \delta.
		\end{equation}
		In particular, in light of \eqref{equation regular embedding in compactly gen 1}, \eqref{equation regular embedding in compactly gen 2} and \eqref{equation regular embedding in compactly gen 3} we obtain for all $g\in V$, $k\in K$ and $\rho\in U_{\widehat{H}}(\chiup_0,c_G(\overline{V},K),\delta)$ that
		\begin{equation*}
			\begin{split}
				\modu{\chiup(t)-\rho^{g}(k)}&\leq \modu{\chiup(k)-\chiup_0^{g_0}(k)}+ \modu{\chiup_0^{g_0}(k)-\chiup_0^g(k)}+ \modu{\chiup_0^g(k)-\rho^g(k)}\\
				&<\lambda+ \delta + \delta \leq \lambda + (\varepsilon-\lambda)=\varepsilon.
			\end{split}
		\end{equation*}
		This proves that $\phi$ is continuous and the result follows.
	\end{proof}
	
	\subsection{Quasi-orbits and $G$-faithful representations}\label{subsection quasi orbit}
	The purpose of this section is to recall the notion of quasi-orbit and to provide a proof of Theorem \ref{theorem letter locally comapct groups equivalence of faitfull}. We start with an important observation.
	\begin{lemma}[{\cite[Section $3.8$, Corollary $2$]{Mackey1976theory}}]\label{lemma Mackey les restriction sont of uniform multiplicity}
		Let $G$ be a locally compact group, $H\leq G$ be a closed normal Type {\rm I} subgroup and $\pi$ be a factor representation of $G$. Then, $\Res_H^G(\pi)$ is a multiple of a multiplicity free representation.
	\end{lemma}
	Now, let $G$ be a locally compact group, $H\leq G$ be a closed normal Type {\rm I} subgroup and $\pi$ be a factor representation of $G$. We recall that the decomposition of a representation of a Type {\rm I} group as a direct integral decomposition of irreducible representations is essentially unique and that for a multiplicity free representation, this decomposition can be achieved with a measure on the unitary dual equipped with the Mackey-Borel structure, see Theorem \ref{Theorem de decomposition des representation de type I en integral direct} and \cite[Theorem $2.15$]{Mackey1976theory}. In particular, in light of Lemma \ref{lemma Mackey les restriction sont of uniform multiplicity}, there exists (up to equivalence) a unique Borel measure $\nu$ on $\widehat{H}$ and a unique $n\in \N\cup \{\infty\}$ such that
	$$\Res_H^G(\pi)\simeq n\int_{\widehat{H}}^\oplus \sigma \qq \diff \nu (\sigma)$$
	is a decomposition of $\Res_H^G(\pi)$ as a direct integral decomposition of irreducible representations of $H$. Furthermore, notice that the measure $\nu$ is $G$-invariant. To be more precise, for all $g\in G$ and every Borel set $B$ of $\widehat{H}$, we let $$g\cdot B=\{\sigma^g: \sigma \in B\}.$$ We define the conjugate measure $\nu^g$ on $\widehat{H}$ by $\nu^g(B)=(g^{-1}\cdot B)$ for each Borel subset $B$ of $\widehat{H}$. The measures $\nu$ and $\nu^g$ are equivalent. Indeed, for each $g\in G$, one has that $\pi\simeq \pi^g$ so that 
	$$\Res_H^G(\pi) \simeq \Res_H^G(\pi^g)\simeq \Res_H^G(\pi)^g\simeq n\int_{\widehat{H}}^\oplus \sigma^g \qq \diff \nu (\sigma)\simeq n \int_{\widehat{H}}^\oplus \sigma \qq \diff \nu^g (\sigma).$$ 
	The uniqueness of the decomposition of $\Res_H^G(\pi)$ as a direct integral of irreducible representations of $H$ ensures that $\nu$ and $\nu^g$ are equivalent.
	\begin{definition}
		The measure class of the measure $\nu$ appearing in the essentially unique decomposition of $\Res_H^G(\pi)$ as direct integral of irreducible representations of $H$ is called the $H$-\tg{quasi-orbit} of $\pi$. 
	\end{definition}
	Our next task is to show that  if $H$ is regularly embedded, the $H$-quasi-orbit of any irreducible representation $\pi$ of $G$ is essentially supported in a single $G$-orbit in $\widehat{H}$. This requires some preliminaries.
	\begin{definition}
		Let $\pi$ and $\rho$ be two representations of a locally compact group $G$. These representations are said to be \tg{disjoint} if no subrepresentation of $\pi$ is equivalent to a subrepresentation of $\rho$. In addition, $\pi$ and $\rho$ are said to be \tg{orbitally disjoint} if $\pi^g$ and $\sigma$ are disjoint for all $g\in G$. 
	\end{definition}
	\begin{lemma}[{\cite[Chapter 3.8, Corollary 1]{Mackey1976theory}}]\label{Lemma the rest of a primary cannot be written as orbitally disjoint rep}
		Let $G$ be a locally compact group, $H\leq G$ be a closed normal subgroup and $\pi$ be a factor representation of $G$. Then, $\Res_H^G(\pi)$ is not a sum of two orbitally disjoint representations.  
	\end{lemma}
	This leads as announced to the following result. 
	\begin{proposition}[{\cite[Chapter 3.8, pg. 186]{Mackey1976theory}}]\label{Proposition the decomp of the restriction is essentially in one orbit}
		Let $G$ be a locally compact group, $H\leq G$ be a closed Type {\rm I} regularly embedded normal subgroup and $\pi$ be a factor representation of $G$. Then, the quasi-orbit of $\pi$ on $H$ is essentially supported in a single $G$-orbit of $\widehat{H}$.
	\end{proposition}
	\begin{proof}
		Since $H$ is regularly embedded in $G$, there exists a countable family $\mathcal{O}$ of $G$-stable Borel sets separating the orbits in $\widehat{H}$. Appending their complement to $\mathcal{B}$ if needed, we can suppose that $\mathcal{B}$ is stable by complements. Now, let $B\subseteq \widehat{G}$  be a $G$-stable Borel set and notice 
		$$\Res_H^G(\pi)= n \int_B^\oplus \sigma\qq  \diff \nu (\sigma)\oplus n\int_{B^c}^\oplus \sigma \qq \diff \nu (\sigma).$$
		On the other hand, Lemma \ref{Lemma the rest of a primary cannot be written as orbitally disjoint rep} ensures that $\Res_H^G(\pi)$ is not the sum of two orbitally disjoint representations of $G$. However, as both $B$ and $B^c$ are $G$-stable measurable sets and since the decomposition of representations of $H$ as a direct integral of irreducible representations are essentially unique, no subrepresentation of $n\int_B^\oplus \sigma\qq  \diff \nu (\sigma)$ is equivalent to the conjugate of a subrepresentation of $n\int_{B^c}^\oplus \sigma \qq \diff \nu (\sigma)$. This implies that $$\Res_H^G(\pi)=n\int_B^\oplus \sigma \diff \nu (\sigma)$$ or that $\nu(B)=0$. Now, consider the set $$\mathcal{P}=\bigg\{B\in \mathcal{B} : \Res_H^G(\pi)\simeq n\int_{B}^\oplus \sigma \qq \diff \nu(\sigma)\bigg\}=\{B\in \mathcal{B} : \nu (B^c)=0\}$$ and notice that each $B\in \mathcal{B}$ either belongs to $\mathcal{P}$ or has a complement $B^c$ which belongs to $\mathcal{P}$. Since $\mathcal{B}$ is countable, we have that $\nu(\bigcup_{B\in \mathcal{P}} B^c)=0$ and hence that $$ \Res_H^G(\pi)\simeq n \int_{\bigcap_{B\in \mathcal{P}}B}^\oplus \sigma \qq \diff \nu (\sigma).$$ 
		In particular, $\nu$ is supported on $\bigcap_{B\in \mathcal{P}} B$. We are left to show that $\bigcap_{B\in \mathcal{P}} B$ consist of a single $G$-orbit in $\widehat{H}$. Suppose for a contradiction that this set contains two disjoint orbits $\mathcal{O}_1,\mathcal{O}_2$. Since $\mathcal{B}$ separates the orbits, there exists a set $B'\in \mathcal{B}$ such that $\mathcal{O}_1 \subseteq B'$ and $\mathcal{O}_2 \subseteq B'^c$. If $B'\in \mathcal{P}$, this provides a contradiction with the inclusion $\mathcal{O}_2\subseteq \bigcap_{B\in \mathcal{P}}B$ and if $B'^c\in \mathcal{P}$, this provides a contradiction with the inclusion $\mathcal{O}_1\subseteq \bigcap_{B\in \mathcal{P}}B.$ 
	\end{proof}
	Our purpose is to use this result in order to prove Theorem \ref{theorem letter locally comapct groups equivalence of faitfull}. We start by recalling the definition of a $G$-faithful representation. 
	\begin{definition}
		Let $G$ be a locally compact group and $H\leq G$ be a closed normal subgroup. A representation $\sigma$ of $H$ is said to be $G$-\tg{faithful} if $\bigcap_{g \in G} \Ker(\sigma^g)=\{1_H\}.$ 
	\end{definition}
	It is clear from this definition that every faithful representation of $H$ is $G$-faithful. The following example exposes a more interesting situation.
	\begin{example}\label{Example G faitfull rep}
		Let $G$ be the affine group $\R\rtimes \R^*_+$ where $\R^*_+$ is the multiplicative group of strictly positive real numbers and let $$H=\{(r,1): r\in \R\}.$$ We claim that not a single irreducible representation of $H$ is faithful but that each non-trivial irreducible representation of $H$ is $G$-faithful. Indeed, the dual of $H$ is easily identified to $$\widehat{H}=\{\chiup_s:H\rightarrow \C^* :(r,1) \mapsto e^{2\pi i sr}\mid s\in \R\}.$$
		Now, let $s\in \R-\{0\}$ and notice that $\Ker(\chiup_s)=(\frac{1}{s}\Z,1)$. It follows that not a single irreducible representation of $H$ is faithful. On the other hand, for each $t\in \R_+^*$ we have that $\chiup_s^{(0,t)}= \chiup_{st}$. It follows that $\bigcap_{g\in G}\Ker(\chiup_s^g)=\{(0,1)\}$ so that each non-trivial irreducible representation of $H$ is $G$-faithful.  
	\end{example}
	The following result provides a relation between the faithful irreducible representations of $G$ and the $G$-faithful irreducible representation of $H$ if this subgroup is regularly embedded in $G$.
	\begin{proposition}
		Let $G$ be a locally compact group, $H\leq G$ be a closed Type {\rm I} regularly embedded normal subgroup and $\pi$ be a faithful irreducible representation of $G$ with $H$-quasi-orbit $\nu$. Then, every irreducible representation of $H$ in the support of $\nu$ is $G$-faithful. 
	\end{proposition}
	\begin{proof}
		Proposition \ref{Proposition the decomp of the restriction is essentially in one orbit} ensures that $\nu$ is essentially supported on a single $G$-orbit $\mathcal{O}$ in $\widehat{H}$. Furthermore, as $H$ is regularly embedded in $G$, the $G$-orbit $\mathcal{O}$ is Mackey-Borel measurable. It follows for each $\sigma\in \mathcal{O}$, that one has that
		$$\bigcap_{g\in G}\Ker (\sigma^g)\subseteq\Ker\bigg(n\int_{G\cdot \sigma}^\oplus\sigma \qq\diff \nu (\sigma)\bigg)=\Ker(\Res_H^G(\pi))\subseteq \Ker(\pi)=\{1_G\}.$$
	\end{proof}
	In particular, if $G$ is irreducibly represented, each closed Type {\rm I} regularly embedded normal subgroup $H$ of $G$ admits a $G$-faithful irreducible representation. This proves Theorem \ref{theorem letter locally comapct groups equivalence of faitfull}. 
	
	\subsection{A new look at Gelfand-Raikov's theorem}
	The purpose of this section is to provide a proof of Theorem \ref{theorem letter locally comapct groups equivalence of faitfull part II}. This one relies on delicate measure theoretical arguments involving direct integral decompositions into irreducible representations just as in \cite{BekkadeleHarpe2008}. To provide an intuition of the type of information that can be deduced from this kind of reasoning we start with a proof of the famous Gelfand-Raikov theorem due to Mackey. 
	\begin{theorem}[Gelfand Raikov's theorem {\cite[pg.110]{Mackey1976theory}}]\label{Theorem Gelfand Raikov's Theorem}
		Let $G$ be a second-countable locally compact group and $g\in G-\{1_G\}$. Then, there exists an irreducible representation $\pi$ of $G$ that is weakly contained in the regular representation of $G$ for which $\pi(g)\not= \id_{\Hr{\pi}}$
	\end{theorem}
	\begin{proof}
		Suppose for a contradiction that $G$ does not admit any irreducible representation such that $\pi(g)\not=\pi(1_G)$ and consider the left-regular representation $\lambda_G$ of $G$. This is a faithful representation and Theorem \ref{Theorem decomposition en integrale direct} ensures the existence of a decomposition  $\lambda_G\simeq\int_\Omega \pi_\omega \diff \nu (\omega)$ as a direct integral of irreducible representations of $G$. However, the above hypothesis implies that $$\lambda(g)=\bigg(\int_\Omega \pi_\omega \qq \diff \nu(\omega)\bigg)(g)= \id_{L^2(G)}= \lambda_G(1_G).$$ This is a contradiction. 
	\end{proof} 
	The following corollary provides an important source of examples of irreducibly represented locally compact groups.
	\begin{corollary}\label{corollary example de groupes irred faithfull}
		Let $G$ be a locally compact group without any closed non-trivial proper amenable normal subgroup. Then, $G$ is irreducibly represented. 
	\end{corollary}
	This follows from the following observation. 
	\begin{lemma}\label{Lemma the kernel of of a representation that is weakly contained in the regular representation is amenable}
		Let $G$ be a locally compact group and $\pi$ be a representation of $G$ that is weakly contained in the regular representation $\lambda_G$. Then, one has that $$\Ker(\pi)=\{g\in G : \pi(g)=\id_{\Hr{\pi}}\}$$ is a closed normal amenable subgroup of $G$.  
	\end{lemma}
	\begin{proof}
		Let $H=\Ker(\pi)$ and recall from our hypothesis that $\pi\prec \lambda_G$. By continuity of the restriction map and since $\pi$ is trivial on $H$, this implies that 
		$$1_{\widehat{H}}\sim \Res_H^G(\pi)\preceq \Res_H^G(\lambda_G)\sim \lambda_H.$$
		It follows that $H$ is amenable.
	\end{proof}
	\begin{proof}[Proof of Corollary \ref{corollary example de groupes irred faithfull}]
		The trivial group is irreducibly represented. On the other hand, if there exists $g\in G-\{1_G\}$, Gelfand-Raikov theorem ensures the existence of an irreducible representation $\pi$ of $G$ that is weakly contained in the regular representation $\lambda_G$ for which $\pi(g)\not=\id_{\Hr{\pi}}$. Now, Lemma \ref{Lemma the kernel of of a representation that is weakly contained in the regular representation is amenable} ensures that $\Ker(\pi)$ is a closed normal amenable subgroup. It follows from our hypotheses that $\Ker(\pi)=\{1_G\}$ so that $G$ is irreducibly represented.
	\end{proof}
	The proof of Theorem \ref{theorem letter locally comapct groups equivalence of faitfull part II} is obtained from a refinement of the above reasoning based on Mackey's theory of induction and regular embeddings. Let $G$ be a locally compact group and $H$ be a closed Type {\rm I} regularly embedded subgroup. A natural question to ask in light of Section \ref{subsection quasi orbit} is if for each $G$-orbit $\mathcal{O}$ in $\widehat{H}$ there exists an irreducible representation $\pi$ of $G$ with quasi-orbit essentially supported on $\mathcal{O}$. We are going to show that this is the case if the $G$-orbit $\mathcal{O}$ is locally closed in $\widehat{H}$. We denote by $\widetilde{G}$ the equivalence classes of representations of $G$ and, for every subset $\mathcal{S}\subseteq \widetilde{G}$, we denote by $\cl(\mathcal{S})$ the closure of $\mathcal{S}$ in $\widetilde{G}$ for the Fell topology. Furthermore, we recall from the Remark \ref{remark the support is weakly equivalent to the set of rep in the fell topo}, the existence of a closed subset $\Sp(\mathcal{S})$ of $\widehat{G}$ called the \tg{support} of $\mathcal{S}$ such that $\mathcal{S}$ is weakly equivalent to $\Sp(\mathcal{S})$. In other words, the closure of $\mathcal{S}$ and $\Sp(\mathcal{S})$ coincide for the Fell topology on $\widetilde{G}$ so that $$\Sp(\mathcal{S})=\cl({\mathcal{S}})\cap \widehat{G}.$$
	\begin{remark}
		If $H\leq G$ is a closed Type {\rm I} regularly embedded normal subgroup, a $G$-orbit $\mathcal{O}$ in $\widehat{H}$ is locally closed if and only if $$\mathcal{O}\cap \Sp(\Sp(\mathcal{O})-\mathcal{O})=\es.$$
	\end{remark} 
	\begin{theorem}\label{le theorem technique de l'existence d'une irred qui se restreint bien dans les orbits}
		Let $G$ be a locally compact group, $K\times H\leq G$ be a closed subgroup such that $K$ is a compact normal subgroup of $G$ 
		and $H$ is a closed normal abelian regularly embedded subgroup of $G$. Then, for every $\sigma\in \widehat{K}$ and each $\chiup\in \widehat{H}$, there exists an irreducible representation $\pi$ of $G$ with the three following properties:
		\begin{enumerate}[label=(\roman*)]
			\item $\pi$ is weakly contained in the regular representation $\lambda_G$.
			\item The $K$-quasi-orbit of $\pi$ is supported in the $G$-orbit of $\sigma$ in $\widehat{K}$.
			\item The $H$-quasi-orbit of $\pi$ is supported in the $G$-orbit of $\chiup$ in $\widehat{H}$.
		\end{enumerate}
	\end{theorem}
	The proof of this theorem requires some preliminary results. We start with two classical results that will be of fundamental importance.
	\begin{theorem}[{\cite[Theorem 3.1]{Fell1962}}]\label{Theorem fell sur les integrale direct}
		Let $G$ be a locally compact group and $\int^\oplus_{\Omega}\pi_\omega \qq \diff \nu (\omega)$ be the direct integral of a measurable field $\omega \mapsto \pi_\omega$ of representations $G$. Then, $\int^\oplus_{\Omega}\pi_\omega \qq \diff \nu (\omega)$ is weakly equivalent to the set of all those representations $\pi_\omega$ with $\omega \in \Omega$ such that $\nu(\Phi^{-1}(U))>0$ for all open neighbourhoods $U$ of $\pi_\omega$ in the Fell topology of $\widetilde{G}$.
	\end{theorem}
	\begin{theorem}[{\cite[Theorem 4.5]{Fell1962}}]\label{Fell la restr de l'in est weak equi a l'orbite}
		Let $G$ be a locally compact group, $H\leq G$ be a closed normal subgroup and $\sigma$ be a representation of $H$, then $\Res_H^G(\Ind_H^G(\sigma))$ is weakly equivalent to the $G$-orbit of $\sigma$ in $\widetilde{H}$.
	\end{theorem}
	Finally, we have the following technical result.
	\begin{lemma}\label{lemma technique pour la mesurabilité de l'ensemble des rep qui se restreigne}
		Let $G$ be a locally compact group, $H\leq G$ be a closed subgroup, $\pi$ be a representation of $G$ and $\int_{\Omega}^\oplus\pi_\omega \diff \nu( \omega)$ be a direct integral decomposition of $\pi$ into irreducible representations. Then, for each subset $\mathcal{S}\subseteq\widetilde{H}$, the set
		$$\Omega_\mathcal{S}=\{\omega\in\Omega: \Res_H^G(\pi_\omega)\preceq \mathcal{S}\}$$
		is a measurable subset of $\Omega$. 
	\end{lemma}
	\begin{proof}
		We recall from \cite[Lemma $1.11$ and Theorem $2.3$]{Fell1962} that the map $\Res_H^G: \widetilde{G}\rightarrow \widetilde{H}$ is continuous for the respective Fell topologies. In particular, $$\big(\Res_H^G\big)^{-1}(\mbox{cl}(\mathcal{S}))\cap \widehat{G}$$ defines a closed subset of $\widehat{G}$. On the other hand, 
		the map $$\Phi:\Omega \rightarrow \widehat{G}: \omega \mapsto \pi_\omega $$ is Borel measurable for the Mackey-Borel structure on $\widehat{G}$. 
		It follows that $\Omega_\mathcal{S}$ is measurable since the Mackey-Borel structure is finer than the Borel structure coming from the Fell topology on $\widehat{G}$ and since 
		$$\Omega_\mathcal{S}= \Phi^{-1}\bigg(\big(\Res_H^G\big)^{-1}(\mbox{cl}(\mathcal{S}))\cap \widehat{G}\bigg).$$
	\end{proof}
	We are finally ready to prove Theorem \ref{le theorem technique de l'existence d'une irred qui se restreint bien dans les orbits}.
	\begin{proof}[Theorem \ref{le theorem technique de l'existence d'une irred qui se restreint bien dans les orbits}]
		Consider the representation $ \sigma\otimes \chiup $ of $K\times H$ and the decomposition of the induced representation
		$$\Ind_{K\times H}^G( \sigma\otimes \chiup)\simeq \int^\oplus_{\Omega}\pi_\omega\qq \diff \nu(\omega)$$
		as a direct integral of irreducible representations of $G$. We are going to prove the existence of a measurable subset $\Lambda\subseteq \Omega$ of positive measure $\nu(\Lambda)>0$ such that the three following assertions hold for all $\omega$ in $\Lambda$:
		\begin{enumerate}[label=(\roman*)]
			\item $\pi_\omega$ is weakly contained in the regular representation $\lambda_G$.
			\item The $K$-quasi-orbit of $\pi_\omega$ is supported on the $G$-orbit of $\sigma$ in $\widehat{K}$.
			\item The $H$-quasi-orbit of $\pi_\omega$ is supported on the $G$-orbit of $\chiup$ in $\widehat{H}$.
		\end{enumerate}
		We start by showing that $\pi_\omega$ is weakly contained in the regular representation $\lambda_G$ for $\nu$-almost all $\omega$ in $\Omega$. Indeed, since $K$ is compact and $H$ is abelian, the group $K\times H$ is amenable. It follows that $ \sigma\otimes \chiup \preceq \lambda_{K\times H}$ so that
		$$\int^\oplus_{\Omega}\pi_\omega\qq \diff \nu(\omega)\simeq \Ind_{K\times H}^G(\sigma \otimes \chiup)\preceq \Ind_{K\times H}^G(\lambda_{K\times H})\simeq\lambda_G.$$ 
		In light of Theorem \ref{Theorem fell sur les integrale direct}, this implies that $\pi_\omega$ is weakly contained in the regular representation $\lambda_G$ for $\nu$-almost all $\omega$ in $\Omega$. 
		
		Now, let us prove that the $K$-quasi-orbit of $\pi_\omega$ is supported on the $G$-orbit of $\sigma$ in $\widehat{K}$ for $\nu$-almost all $\omega$ in $\Omega$. Since $K$ is a compact group, Lemma \ref{lemma compact groups are regularly embedded} ensures that $K$ is regularly embedded in $G$ and that each $G$-orbit in $\widehat{K}$ is both open and closed for the Fell topology. We denote by $G\cdot\sigma$ the $G$-orbit of $\sigma$ in $\widehat{K}$. Lemma \ref{lemma technique pour la mesurabilité de l'ensemble des rep qui se restreigne} ensures that the set  
		$$\Omega_{G\cdot \sigma}=\{\omega\in\Omega: \Res_K^G(\pi_\omega)\preceq G\cdot\sigma\}$$
		is measurable. On the other hand, as $G\cdot\sigma$ is closed in the Fell topology of $\widehat{K}$, every $\pi_\omega$ with $\omega \in \Omega_{G\cdot \sigma}$ has its $K$-quasi-orbit supported on $G\cdot \sigma$. We are going to show that $\Omega- \Omega_{G\cdot \sigma}$ has measure $0$.
		Since $\Omega_{G\cdot \sigma}$ and $\Omega-\Omega_{G\cdot \sigma}$ are Mackey-Borel sets, both direct integrals $$\rho_0=\int^\oplus_{\Omega_{G\cdot \sigma}}\Res_{K\times H}^G(\pi_\omega) \qq \diff \nu(\omega)\mbox{ and }\rho_1= \int^\oplus_{\Omega-\Omega_{G\cdot \sigma}}\Res_{K\times H}^G(\pi_\omega) \qq \diff \nu(\omega)$$ are well defined on a possibly zero-dimensional Hilbert space. Furthermore, one has that 
		$$ \Res_{K}^G(\Ind_{K\times H}^G(\chiup\otimes\sigma))\simeq\Res_{K}^{K\times H}(\rho_0)\oplus \Res_{K}^{K\times H}(\rho_1).$$
		On the other hand, Theorem \ref{Fell la restr de l'in est weak equi a l'orbite} ensures that 
		\begin{equation}\label{equation weak equivalence of res of ind}
			\Res_{K\times H}^G(\Ind_{K\times H}^G(\sigma\otimes\chiup))\sim G\cdot (\sigma \otimes \chiup),
		\end{equation}so that
		\begin{equation*}
			\begin{split}
				\Res_{K}^G(\Ind_{K\times H}^G(\sigma\otimes\chiup ))\sim \Res_{K}^{K\times H}(G\cdot (\sigma \otimes \chiup))\sim G\cdot (\Res_K^{K\times H}(\sigma \otimes \chiup))\sim G\cdot\sigma.
			\end{split}
		\end{equation*}
		Since $G\cdot \sigma$ is closed, it follows that $\Res_{K}^G(\Ind_{K\times H}^G(\sigma\otimes\chiup))$ does not weakly contain any element of $\widehat{K}-G\cdot\sigma$. 
		However, by the definition of $\Omega_{G\cdot \sigma}$, Theorem \ref{Theorem fell sur les integrale direct} ensures that 
		$$\Res_{K}^{K\times H}(\rho_1)\simeq \int_{\Omega-\Omega_{G\cdot \sigma}} \Res_K^G(\pi_{\omega})\diff \nu (\omega)$$
		is weakly equivalent to a subset of $\widehat{K}-G\cdot\sigma$. Since $\rho_1$ is weakly contained in $\Res_{K}^G(\Ind_{K\times H}^G(\sigma\otimes\chiup))$, this implies that $\nu(\Omega-\Omega_{G\cdot \sigma})=0$.
		
		We are left to show the existence of a measurable subset $\Lambda\subseteq \Omega$ of positive measure $\nu(\Lambda)>0$ such that the $H$-quasi-orbit of $\pi_\omega$ is supported on the $G$-orbit of $\chiup$ in $\widehat{H}$ for all $\omega$ in $\Lambda$. We let $\mathcal{S}= \Sp(G\cdot \chiup)-G\cdot \chiup$. Lemma \ref{lemma technique pour la mesurabilité de l'ensemble des rep qui se restreigne} ensures that the set   
		$$\Omega_{\mathcal{S}}=\{\omega\in\Omega: \Res_H^G(\pi_\omega)\preceq \Sp(G\cdot \chiup)-G\cdot \chiup\}$$
		is measurable. We are going to show that $\Lambda=\Omega-\Omega_{\mathcal{S}}$ is as desired. Since $\Omega_{\mathcal{S}}$ and $\Omega- \Omega_{\mathcal{S}}$ are  Mackey-Borel sets, both direct integrals $$\tau_0=\int^\oplus_{\Omega_{\mathcal{S}}}\Res_{K\times H}^G(\pi_\omega) \qq \diff \nu(\omega)\mbox{ and }\tau_1= \int^\oplus_{\Omega-\Omega_{\mathcal{S}}}\Res_{K\times H}^G(\pi_\omega) \qq \diff \nu(\omega)$$ are well defined on a possibly zero-dimensional vector space. Furthermore, one has that
		\begin{equation*}
			\begin{split}
				\Res_{H}^G(\Ind_{K\times H}^G(\sigma\otimes\chiup))\simeq\Res_{H}^{K\times H}(\tau_0)\oplus \Res_{H}^{K\times H}(\tau_1).
			\end{split}
		\end{equation*}
		On the other hand, Theorem \ref{Theorem fell sur les integrale direct} ensures that $$\Res_{H}^{K\times H}(\tau_0)\simeq\int^\oplus_{\Omega_{\mathcal{S}}}\Res_H^G(\pi_{\omega})\diff \nu (\omega)$$ is weakly equivalent to a subset of $\Sp(\Sp(G\cdot \chiup)-G\cdot \chiup)$.
		Furthermore, $G\cdot \chiup$ is locally closed due to Proposition \ref{proposition for abelian groups regular embedding is equivalent to the fact taht every orbits is locally closed}. In particular, one has that $\Sp(\Sp(G\cdot \chiup)-G\cdot \chiup)$ is a closed subset of $\widehat{H}$ that does not weakly contain $G\cdot \chiup$. Nevertheless, the weak equivalence \eqref{equation weak equivalence of res of ind} implies that 
		\begin{equation}\label{equation weak equivalence}
			\Res_{H}^G(\Ind_{K\times H}^G(\sigma\otimes\chiup))\sim \Res_{H}^{K\times H}(G\cdot (\sigma \otimes \chiup))\sim G\cdot (\Res_H^{K\times H}( \sigma\otimes \chiup))\sim G\cdot\chiup.
		\end{equation}
		It follows that $\Res_{H}^{K\times H}(\tau_1)$ is not zero-dimensional, so that $\nu(\Omega- \Omega_{\mathcal{S}})>0$. On the other hand, the weak equivalence \eqref{equation weak equivalence} also ensures that the $H$-quasi-orbit of $\pi_\omega$ is essentially supported on a $G$-orbit inside $\Sp(G\cdot \chiup)$ for $\nu$-almost all $\omega$ in $\Omega$. Since $\{\pi_\omega:\omega \in \Omega_\mathcal{S}\}$ is exactly the set of representations $\pi_\omega$ ($\omega\in \Omega$) whose $H$-quasi-orbit is not supported on the $G$-orbit of $\chiup$ in $\widehat{H}$, we deduce that the $H$-quasi-orbit of $\pi_\omega$ is essentially supported on the $G$-orbit of $\chiup$ in $\widehat{H}$ for all $\omega \in \Omega-\Omega_{\mathcal{S}}$. The result follows.
	\end{proof}
	\noindent This theorem is particularly relevant in light of the following result.
	\begin{proposition}\label{proposition if the orbit is G faith and se rep has quasi orbit essentially supported there then the kernel intersects the subgroup H triv}
		Let $G$ be a locally compact group, $H\leq G$ be a closed Type {\rm I} regularly embedded normal subgroup, $\sigma$ be a $G$-faithful irreducible representation of $H$ and $\pi$ be an irreducible representation of $G$ with $H$-quasi-orbit essentially supported on the $G$-orbit of $\sigma$. Then, $\Ker(\pi)$ intersects $H$ trivially.
	\end{proposition}
	\begin{proof}
		Our hypotheses on $\pi$ ensures that 
		$$\Res_H^G(\pi)\simeq n \int_{\widehat{H}}^\oplus \rho \qq \diff \nu (\rho)$$
		for some measure $\nu$ on $\widehat{H}$ that is essentially supported on $G\cdot \sigma$. Since the measure $\nu$ is $G$-invariant, this implies that each representation $\rho \in G\cdot \sigma$ is in the support of $\nu$. Now, notice that $\Ker(\Res_H^G(\pi))=\Ker(\pi)\cap H$ and let us prove that $\Ker(\Res_H^G(\pi))\subseteq \bigcap_{g\in G}\Ker(\sigma^g)$. The result will follow as $\sigma$ is $G$-faithful. Taking the complement, let $h\in \bigcup_{g\in G}(H-\Ker(\sigma^g))$ and let us prove that $h\not\in \Ker(\Res_H^G(\pi))$. By choice of $h$, there exists some $g\in G$ such that $\sigma^g(h)\not=\id_{\Hr{\sigma}}$. In particular, there exists a function of positive type $\varphi$ associated with $\sigma^g$ such that $\varphi(h)\not=1$. By definition of the Fell topology, there exists an open neighbourhood $\mathcal{O}$ of $\sigma^g$ in $\widehat{H}$ such that every representation $\rho \in \mathcal{O}$ shares that property (for instance consider the set $\mathcal{O}=\mathcal{W}(\sigma^g, \varphi, h, \modu{1-\varphi(ghg^{-1})}/2)$ as defined on page \pageref{voisinage pour la topo de Fell}). In particular, for each $\rho\in \mathcal{O}$ one has that $\rho(h)\not=\id_{\Hr{\rho}}$. On the other hand, as $\sigma^g$ is in the support of $\nu$, the representation $\int_{\mathcal{O}}^\oplus \rho \qq \diff \nu (\rho)$ is well defined on a non-zero dimensional Hilbert space and the above discussion implies that  
		$$h\not \in \Ker\bigg(\int_{\mathcal{O}}^\oplus \rho \qq \diff \nu (\rho)\bigg).$$ 
		It follows that $h\not\in \Ker(\Res_H^G(\pi))$ because
		$$\Res_H^G(\pi)\simeq n \int_{\mathcal{O}}^\oplus \rho \qq \diff \nu (\rho)\oplus n \int_{\widehat{H}-\mathcal{O}}^\oplus \rho \qq \diff \nu (\rho).$$
	\end{proof}
	Theorem \ref{theorem letter locally comapct groups equivalence of faitfull part II} follows from Theorem \ref{le theorem technique de l'existence d'une irred qui se restreint bien dans les orbits} and Proposition \ref{proposition if the orbit is G faith and se rep has quasi orbit essentially supported there then the kernel intersects the subgroup H triv}.

	\section{The center alternative and the Lie socles}\label{section Lie proof of all the assertiosn}
	
	We recall that the purpose of this chapter is to provide an algebraic characterisation of the connected nilpotent Lie groups and of the Hausdorff connected components of the identity of the $\R$-points of a linear algebraic groups defined over $\R$ that are irreducibly represented. Our strategy is to rely on Theorems \ref{theorem letter locally comapct groups equivalence of faitfull} and \ref{theorem letter locally comapct groups equivalence of faitfull part II} for a suitable choice of subgroups of $G$. The purpose of this section is to define these subgroups and to prove the assertions \ref{item the subgroup Sa times Ss} to \ref{item the moment cenre faitfull representaiton }. 
	\subsection{The role played by the Lie socles}
	Let $G$ be a connected Lie group. The purpose of this section is to define a subgroup of $G$ intersecting non-trivially the kernels of non-faithful representations of $G$ that are weakly contained in the regular representation $\lambda_G$. Inspired by \cite{BekkadeleHarpe2008}, we make the following definitions.
	\begin{definition}\label{definition of the semi-simple Lie socle}
		Let $G$ be a connected Lie group. A \tg{semi-simple Lie foot} $K$ of $G$ is a minimal non-trivial compact connected semi-simple normal subgroup of $G$ with trivial center (minimal in the sense that the only subgroup of $K$ that is non-trivial connected and normal in $G$ is $K$ itself). We let $\mathcal{K}_G$ be the set of semi-simple Lie feet of $G$ and define the \tg{semi-simple Lie socle} $\mathcal{S}_s(G)$ as the closed subgroup of $G$ generated by the elements of $\mathcal{K}_G$. 
	\end{definition}
	\begin{definition}\label{definition of the abelian Lie socle}
		Let $G$ be a connected Lie group. An \tg{abelian Lie foot} $V$ of $G$ is a minimal non-trivial closed connected abelian normal subgroup of $G$ intersecting the center of $G$ trivially (minimal in the sense that the only subgroup of $V$ that is non-trivial closed connected and normal in $G$ is $V$ itself). We let $\mathcal{A}_G$ be the collection of abelian Lie feet of $G$ and define the \tg{abelian Lie socle} $\mathcal{S}_a(G)$ as the closed subgroup generated by the set of abelian Lie feet and the center of $G$. 
	\end{definition} 
	In order to obtain Theorems \ref{Corollary letter Nilpotentnt groups} and \ref{Corollary letter amenable lie groups} from Theorems \ref{theorem letter locally comapct groups equivalence of faitfull} and \ref{theorem letter locally comapct groups equivalence of faitfull part II}, we are left to prove that the groups $N=K\times H$, $K=\mathcal{S}_s(G)$ and $H=\mathcal{S}_a(G)$ satisfy the desired hypotheses. This is done by proving the assertions \ref{item the subgroup Sa times Ss} to \ref{item the moment cenre faitfull representaiton }. 
	
	If the group $G$ is nilpotent, both $\mathcal{S}_s(G)$ and $\mathcal{S}_a(G)$ are easily described. This is the purpose of the following result which implies the assertion \ref{item Sa G if G is nilpotent}.
	\begin{lemma}\label{lemma in nilpotenet groups minimal connected subgroups are central}
		Let $G$ be a connected nilpotent Lie group then $\mathcal{Z}(G)$ is positive dimensional, $\mathcal{S}_a(G)= \mathcal{Z}(G)$ and $\mathcal{S}_s(G)=\{1_G\}$.
	\end{lemma}
	\begin{proof}
		Let $H$ be a closed connected normal subgroup of $G$. We are going to show that $H$ contains a positive dimensional closed central subgroup of $G$. Our hypotheses ensures that $H$ corresponds to a non-trivial ideal $\mathfrak{h}\subseteq \mathfrak{g}$ of the Lie algebra $\mathfrak{g}$ of $G$ under the Lie correspondence. Now, let $\mathfrak{h}^{0}=\mathfrak{h}$ and for each positive integer $r\geq 1$ set $\mathfrak{h}^{r}=\lb\mathfrak{h}^{r-1},\mathfrak{g}\rb$. Since $\mathfrak{g}$ is nilpotent and since $\mathfrak{h}$ is non-trivial, there exists a positive integer $n\in \N$ such that $\mathfrak{h}^{n}\not=\{0\}$ and $\lb\mathfrak{h}^{n},\mathfrak{g}\rb= \mathfrak{h}^{n+1}=\{0\}$. Since the exponential map is surjective on connected nilpotent Lie groups, the Campbell-Baker-Hausdorff formula ensures that the elements of $\exp(\mathfrak{h}^n)$ commute with all the elements of $G$. In particular, the closed subgroup generated by $\exp(\mathfrak{h}^n)$ is a positive dimensional closed connected central subgroup of $G$ contained in $H$. It follows directly that $\mathcal{A}_G=\es$ and that $\mathcal{K}_G=\es$. In particular, one has that $\mathcal{S}_a(G)= \mathcal{Z}(G)$ and that $\mathcal{S}_s(G)=\{1_G\}$. Furthermore, taking $H=G$, we obtain that $\mathcal{Z}(G)$ is positive dimensional.
	\end{proof}
	Now, let us prove the assertion \ref{item le noyau des rep faiblement reguliere intersectes les socles}. We start with an important observation.
	\begin{lemma}\label{lemma discrete groups are central}
		Let $G$ be a connected topological group and let $\Gamma\leq G$ be a closed discrete normal subgroup. Then, $\Gamma$ is central. 
	\end{lemma}
	\begin{proof}
		Let $k\in \Gamma$. Since $\Gamma$ is a closed normal subgroup of $G$, notice that the map $$\phi_k:G \rightarrow \Gamma :g\mapsto gkg^{-1}k^{-1}$$ is well defined and continuous. On the other hand, one has that $\phi_k(1_G)=1_G$. As $G$ is connected and since $\Gamma$ is discrete we deduce that $\phi_k(G)=\{1_G\}$. It follows that $k\in \mathcal{Z}(G)$.
	\end{proof}
	\begin{proposition}\label{Proposition le noyau des rep faiblement reguliere intersectes les socles}
		Let $G$ be a connected Lie group and let $\pi$ be a non-faithful representation of $G$ that is weakly contained in the regular representation $\lambda_G$. Then, $\Ker(\pi)$ intersects either $\mathcal{S}_s(G)$ or $\mathcal{S}_a(G)$ non-trivially.
	\end{proposition}
	\begin{proof}
		Since $\pi$ is weakly contained in the regular representation $\lambda_G$, Lemma \ref{Lemma the kernel of of a representation that is weakly contained in the regular representation is amenable} ensures that $\Ker(\pi)$ is a closed normal amenable subgroup of $G$. If $\Ker(\pi)$ is discrete, Lemma \ref{lemma discrete groups are central} ensures it is central an hence intersects $\mathcal{S}_a(G)$ non-trivially. Now, suppose that $\Ker(\pi)$ is not discrete. The connected component of the identity $\Ker(\pi)^\circ$ is a closed topologically characteristic subgroup of $\Ker(\pi)$. It follows that $H=\Ker(\pi)^\circ$ is a non-trivial closed connected amenable Lie group that is normal in $G$. On the other hand, \cite[Corollary 4.1.9]{Zimmer1984} ensures that a connected Lie group is amenable if and only if it is the compact extension of its radical. We treat several cases separately. If the radical $\mathcal{R}$ of $H$ is trivial, $H$ is a compact connected semi-simple subgroup of $G$. In that case, $H$ either intersects the center of $G$ (and hence $\mathcal{S}_a(G)$) non-trivially or contains a minimal non-trivial compact connected semi-simple normal subgroup $K$ of $G$ whose center must be trivial because discrete normal subgroups of $G$ are central. In particular, in that last case, $H$ contains a semi-simple Lie foot of $G$ and hence intersects $\mathcal{S}_s(G)$ non-trivially. If on the other hand the radical $\mathcal{R}$ of $H$ is not trivial, it corresponds to a non-trivial solvable ideal $\mathfrak{r}$ of the Lie algebra $\mathfrak{g}$ of $G$ under the Lie correspondence (since $\mathcal{R}$ is a topological characteristic subgroup of $H$ and since $H$ is normal in $G$). As $\mathfrak{r}$ is solvable, there exists a positive integer $n$ such that $\lb\mathfrak{r}^{(n)},\mathfrak{r}^{(n)}\rb\not=0$ but $\lb\mathfrak{r}^{(n+1)},\mathfrak{r}^{(n+1)}\rb\not=\{0\}$ where $\mathfrak{r}^{(r)}=\lb\mathfrak{r}^{(r-1)},\mathfrak{r}^{(r-1)}\rb$ for each positive integer $r\geq 1$. In particular, $\mathfrak{r}^{(n)}$ is an abelian ideal of $\mathfrak{g}$ contained in $\mathfrak{r}$ and hence corresponds to a non-trivial connected abelian normal Lie group $N\subseteq H$ under the exponential map. As $H$ is closed, one has that $\overline{N}\subseteq H$. On the other hand, $\overline{N}$ either intersects the center of $G$ or contains an abelian Lie foot, so that $H$ intersects $\mathcal{S}_a(G)$ non-trivially. 
	\end{proof}
	In the rest of this section, we explore the properties of $\mathcal{S}_s(G)$, $\mathcal{S}_a(G)$ and of the closed subgroup of $G$ generated by these two. We start with a series of basic observations.
	\begin{lemma}\label{lemm athe semisimple Lie feet are center free}
		Let $G$ be a connected Lie group and $K\in \mathcal{K}_G$. Then, $K$ is center free. 
	\end{lemma}
	\begin{proof}
		First, notice that the center $\mathcal{Z}(K)$ of $K$ is discrete. Indeed, the connected component $\mathcal{Z}(K)^\circ$ of $\mathcal{Z}(K)$ is a closed connected abelian subgroup of $K$. In particular, it corresponds to an abelian Lie subalgebra of the Lie algebra of $K$ under Lie correspondence. However, $K$ is semi-simple so that its Lie algebra does not contain any non-trivial abelian Lie algebra. It follows that $K$ does not have any non-trivial connected abelian subgroup and $\mathcal{Z}(K)$ is discrete. On the other hand, $\mathcal{Z}(K)$ is a characteristic subgroup of $K$. Since $K$ is normal in $G$, it follows that $\mathcal{Z}(K)$ is normal in $G$. Since $G$ is connected, Lemma \ref{lemma discrete groups are central} ensures that $\mathcal{Z}(K)$ is central in $G$. However, by definition $K$ intersects the center of $G$ trivially so that $K$ is center free. 
	\end{proof}
	\begin{lemma}\label{Lemma group gen par minimal and abelian is dirrect product for the connected comp ULTIMATE}
		Let $G$ be a connected Lie group, $M\leq G$ be a closed normal subgroup and $H$ be a minimal non-trivial closed connected normal subgroup of $G$ intersecting $\mathcal{Z}(G)$ trivially. Then, $H\subseteq M$ or, $H\cap M=\{1_G\}$ and the elements of $H$ and $M$ commute with each other. 
	\end{lemma}
	\begin{proof}
		The connected component of the identity $(H\cap M)^\circ$ of $H\cap M$ is a closed connected normal subgroup of $G$ (by normality of both $H$ and $M$). Furthermore, by minimality of $H$, we have either that $(H\cap M)^\circ=H$ and $H\subseteq M$, or that $(H\cap M)^\circ=\{1_G\}$. In that last case, $H\cap M$ is discrete and Lemma \ref{lemma discrete groups are central} ensures that $H\cap M$ is central. Since $H\cap \mathcal{Z}(G)=\{1_G\}$ this implies that $H\cap M=\{1_G\}$. On the other hand, the commutator subgroup $\lb H,M\rb$ belongs to $H\cap M$ by normality of both subgroups. This implies that $\lb H,M\rb=\{1_G\}$ which proves that the elements of $H$ and $M$ commutes with each other.
	\end{proof} 
	\begin{lemma}\label{Lemma the closed subgroup generated by AG of KG is topologically characteristic}
		Let $G$ be a connected Lie group. The closed subgroup generated by the elements of $\mathcal{A}_G$ (resp. $\mathcal{K}_G$) is topologically characteristic.
	\end{lemma}
	\begin{proof}
		Denote by $\Aut(G)$ the set of topological group automorphisms of $G$. For all $\phi\in \Aut(G)$ and every closed non-trivial connected normal subgroup $H\leq G$, notice that $\phi(H)$ is still closed non-trivial connected and normal in $G$. In particular, if $H$ is minimal for these properties so is $\phi(H)$. Since the properties of being central or abelian is also preserved by any automorphism $\phi\in \Aut(G)$, it follows that 
		$$\phi(\mathcal{A}_G)=\{\phi(V):V\in \mathcal{A}_G\}=\mathcal{A}_G.$$
		In particular, the closed subgroup generated by the elements of $\mathcal{A}_G$ is topologically characteristic. Similarly, the properties of being semi-simple and center free are preserved by any automorphism $\phi\in \Aut(G)$ so that $\phi(K)$ is semi-simple and center free for all $K\in \mathcal{K}_G$ in light of Lemma \ref{lemm athe semisimple Lie feet are center free}. In particular, $\phi(K)$ intersects the center of $G$ trivially and hence belongs to $\mathcal{K}_G$. It follows that 
		$$\phi(\mathcal{K}_G)=\{\phi(K):K\in \mathcal{K}_G\}=\mathcal{K}_G.$$
		In particular, the closed subgroup generated by the elements of $\mathcal{K}_G$ is topologically characteristic.
	\end{proof}
	\begin{proposition}\label{proposition the subgroup spanned by the semi simple or abelian feet is a closed connectedted normal ss subgroup}
		Let $G$ be a connected Lie group. Then, the following assertions hold:
		\begin{enumerate}[label=(\roman*)]
			\item The closed subgroup generated by the elements of $\mathcal{A}_G$ is a connected abelian subgroup of $G$. 
			\item The semi-simple Lie socle $\mathcal{S}_s(G)$ is a compact connected semi-simple subgroup of $G$ with trivial center. 
		\end{enumerate}
		Furthermore, both of these groups are topologically characteristic in $G$. 
	\end{proposition}
	\begin{proof}
		We treat both cases in parallel since the proofs are very similar for both situations. Let $\mathcal{F}$ be a finite subset of $\mathcal{A}_G$ (resp. of $\mathcal{K}_G$). We start by showing that the closed subgroup generated by $\mathcal{F}$ is a closed connected abelian subgroup of $G$ (resp. a compact connected semi-simple subgroup of $G$ with trivial center). The proof is made by induction on $l=\modu{F}$. If $l=1$, the result is trivial in the abelian case and follows from Lemma \ref{lemm athe semisimple Lie feet are center free} in the semi-simple case. Now let $\mathcal{F}=\{H_1,...,H_l\}$ for some $l>1$ and let $H$ be the closed subgroup generated by the elements of $\mathcal{F}$. By our induction hypothesis, the closed subgroup $M$ of $G$ generated by $\{H_2,...,H_{l}\}$ is a closed connected abelian subgroup of $G$ (resp. a compact connected semi-simple subgroup of $G$ with trivial center). Since both $H_1$ and $M$ are connected Lie groups, they are path-wise connected. In particular, the subgroup of $G$ generated by $H_1$ and $M$ is path-wise connected. Similarly, the subgroup of $G$ generated by $H_1$ and $M$ is normal in $G$ by normality of both $H_1$ and $M$. These properties are inherited by the closure so that the closed subgroup generated by the elements of $\mathcal{A}_G$ (reps. $\mathcal{K}_G$) is connected and normal in $G$. On the other hand, by Lemma \ref{Lemma group gen par minimal and abelian is dirrect product for the connected comp ULTIMATE}, we either have that $H_1\subseteq M$ or that $H_1\cap M=\{1_G\}$ and the elements of $H_1$ and $M$ commute with each other. In the first case, we have that $H=M$ and the result follows. Now, suppose that $H_1\cap M=\{1_G\}$ and that the elements of $H_1$ and $M$ commute with each other. In the abelian case, this implies that the subgroup generated by $H_1$ and $M$ is abelian, so that the claim follows by taking the closure. On the other hand, in the compact semi-simple case, the subgroup generated by $H_1$ and $M$ is already compact and hence closed. It follows that the closed subgroup generated by the elements of $\mathcal{F}$ is isomorphic to the direct product $H_1\times M$ and is thus semi-simple and center free. Now, consider the set $\mathcal{P}$ of all closed subgroups of $G$ generated by finitely many elements of $\mathcal{A}_G$ (resp. $\mathcal{K}_G$) and equip $\mathcal{P}$ with the order given by inclusion. Since the length of a totally ordered chain $H_1\lneq H_2\lneq ...\lneq H_n$ of elements of $\mathcal{P}$ is bounded by $\dim(G)$, $\mathcal{P}$ admits a maximal element, say $N$. By the first part of the proof, $N$ is a closed connected abelian subgroup of $G$ (resp. a compact connected semi-simple subgroup of $G$ with trivial center). Furthermore, by maximality, the closed subgroup generated by $N$ and any element of $\mathcal{A}_G$ (resp. $\mathcal{K}_G$) coincides with $N$. This proves that $N$ contains every of the elements of $\mathcal{A}_G$ (resp. $\mathcal{K}_G$) and hence is the closed subgroup generated by the elements of $\mathcal{A}_G$ (resp. $\mathcal{K}_G$). Together with Lemma \ref{Lemma the closed subgroup generated by AG of KG is topologically characteristic}, the result follows. 
	\end{proof}
	The following result follows directly and proves the assertion \ref{item the subgroup Sa}.
	\begin{corollary}\label{Corollary the abelian Lie socle is a closed abelian topologically characteristic subgroup of G}
		Let $G$ be a connected Lie group. Then, the abelian Lie socle $\mathcal{S}_a(G)$ is a closed abelian topologically characteristic subgroup of $G$.
	\end{corollary}
	\begin{proof}
		Proposition \ref{proposition the subgroup spanned by the semi simple or abelian feet is a closed connectedted normal ss subgroup} ensures that the closed subgroup of $G$ generated by the elements of $\mathcal{A}_G$ is a closed abelian topologically characteristic subgroup of $G$. Since the center of $G$ shares all of these properties and since $\mathcal{S}_a(G)$ is the closed subgroup of $G$ generated by these two groups, the result follows. 
	\end{proof}
	The following result proves the assertion \ref{item the subgroup Sa times Ss}.
	\begin{proposition}\label{Proposition la structure du group engendrer par Sa et Ss}
		Let $G$ be a connected Lie group. Then, the closed subgroup of $G$ generated by $\mathcal{S}_a(G)$ and $\mathcal{S}_s(G)$ is a topologically characteristic subgroup of $G$ isomorphic to $\mathcal{S}_a(G)\times \mathcal{S}_s(G)$ (with the natural inclusions). 
	\end{proposition}
	\begin{proof}
		The abelian Lie socle $\mathcal{S}_a(G)$ is a closed normal subgroup of $G$. On the other hand, Proposition \ref{proposition the subgroup spanned by the semi simple or abelian feet is a closed connectedted normal ss subgroup} also ensures that semi-simple Lie socle $\mathcal{S}_s(G)$ is a topologically characteristic compact connected normal subgroup of $G$ intersecting the center trivially. In particular, the closed subgroup generated by $\mathcal{S}_a(G)$ and $\mathcal{S}_s(G)$ is topologically characteristic. Furthermore, applying Lemma \ref{Lemma group gen par minimal and abelian is dirrect product for the connected comp ULTIMATE} to the Lie feet of $G$, one obtains that $\mathcal{S}_a(G)\cap\mathcal{S}_s(G)=\{1_G\}$ and that the elements of $\mathcal{S}_a(G)$ and $\mathcal{S}_s(G)$ commute with each other. It follows that the subgroup generated by $\mathcal{S}_a(G)$ and $\mathcal{S}_s(G)$ is isomorphic to $\mathcal{S}_a(G)\times \mathcal{S}_s(G)$ (with the natural inclusions). On the other hand, this group is closed since $\mathcal{S}_a(G)$ is closed and $\mathcal{S}_s(G)$ is compact. The result follows.
	\end{proof}
	
	\subsection{The semi-simple Lie socle is irreducibly represented}\label{Section semisimple Lie socle} 
	Let $G$ be a connected Lie group. The purpose of this section is to prove that the semi-simple Lie socle $\mathcal{S}_s(G)$ of $G$ is irreducibly represented. We recall from Proposition \ref{proposition the subgroup spanned by the semi simple or abelian feet is a closed connectedted normal ss subgroup} that $\mathcal{S}_s(G)$ is a compact connected semi-simple subgroup of $G$ with trivial center. We show below that every such group is irreducibly represented. We start with an analysis of their structure based on the following classical result.
	\begin{theorem}[{\cite[Chapter $10$, Section $7.2$, Theorem $4$]{Procesi2007}}]\label{theorem claudio structure of semi simple connecte compact Lie groups}
		Let $K$ be a connected compact Lie group. Then, there exists connected simply connected compact simple Lie groups $K_1,...,K_n$, a non-negative integer $r\in \N$ and a finite normal subgroup $\Gamma\subseteq K_1\times K_2\times ... \times K_n\times \mathbb{T}^r$ such that 
		$$K\simeq (K_1\times \cdots \times K_n\times \mathbb{T}^r)/\Gamma.$$
	\end{theorem}
	The following result follows directly.
	\begin{proposition}\label{proposition structure of compact connected center free Lie groups}
		Let $K$ be a connected center free compact Lie group. Then, there exists connected simply connected compact simple Lie groups $K_1,...,K_n$ such that 
		$$K\simeq K_1/_{\mathcal{Z}(K_1)}\times \cdots \times K_n/_{\mathcal{Z}(K_n)}.$$
	\end{proposition}
	\begin{proof}
		Since $K$ is a connected compact Lie group, Theorem \ref{theorem claudio structure of semi simple connecte compact Lie groups} ensures the existence of connected simply connected compact simple Lie groups $K_1,...,K_n$, a non negative integer $r\in \N$ and a finite subgroup $\Gamma\subseteq K_1\times K_2\times ... \times K_n\times \mathbb{T}^r$ such that 
		$$K\simeq (K_1\times \cdots \times K_n\times \mathbb{T}^r)/\Gamma.$$
		Since $K_1\times K_2\times ...\times K_n\times \mathbb{T}^r$ is connected, Lemma \ref{lemma discrete groups are central} ensures that $\Gamma$ is central. In particular, 
		$\mathcal{Z}(K_1\times K_2\times ...\times K_n\times  \mathbb{T}^r)/\Gamma$ is a central subgroup of $K$. Since $K$ is center free, this implies that 
		$$\Gamma=\mathcal{Z}(K_1\times K_2\times ...\times K_n\times  \mathbb{T}^r) = \mathcal{Z}(K_1)\times ... \times \mathcal{Z}(K_n)\times  \mathbb{T}^r \q (r=0).$$
		It follows as desired that 
		$$K\simeq K_1/_{\mathcal{Z}(K_1)}\times \cdots \times K_n/_{\mathcal{Z}(K_n)}.$$
	\end{proof}
	The dual of such a direct product is easily understood via the Kronecker product. To be more precise, we have the following classical result.
	\begin{theorem}[{\cite[Theorem $3.1$ and Theorem $3.2$]{Mackey1976theory}}]\label{Theorem Mackey rep of the product}
		Let $H_1$, $H_2$ be two Type {\rm I} locally compact groups. Then, the following holds:
		$$\widehat{H_1\times H_2}=\widehat{H_1}\otimes \widehat{H_2}=\{\sigma_1\otimes \sigma_2:\sigma_1\in \widehat{H_1}, \sigma_2\in \widehat{H_2}\}.$$
	\end{theorem} 
	In addition, the kernel of a Kronecker product of representations can be apprehended via the following result.
	\begin{lemma}[{\cite[Lemma 12]{BekkadeleHarpe2008}}]\label{Lemma Bekka le produit tensoriel de 2 operateur de Hilbert est mult de l'identité ssi les 2 le sont}
		Let $\mathcal{H}_1$ and $\mathcal{H}_2$ be two Hilbert spaces. Let $S_1\in \mathcal{L}(\mathcal{H}_1)$ and $S_2\in \mathcal{L}(\mathcal{H}_2)$ be such that $S_1\otimes S_2\in \mathcal{L}(\mathcal{H}_1\otimes \mathcal{H}_2)$ is a non-zero multiple of the identity. Then $S_1$ and $S_2$ are multiples of the identity.
	\end{lemma}
	We are finally able to prove our claim.
	\begin{proposition}
		Let $K$ be a center free compact Lie group. Then, $K$ is irreducibly represented.
	\end{proposition}
	\begin{proof}
		Proposition \ref{proposition structure of compact connected center free Lie groups} ensures the existence of connected simply connected compact simple Lie groups $K_1,...,K_n$ such that 
		$$K\simeq K_1/_{\mathcal{Z}(K_1)}\times \cdots \times K_n/_{\mathcal{Z}(K_n)}.$$
		For each $i=1,...,n$, choose a non-trivial irreducible representation $\sigma_i$ of  $K_i/_{\mathcal{Z}(K_i)}$ and let us check that, for each $k_i\in K_i/\mathcal{Z}(K_i)$, $k_i\not=1_{K_i/\mathcal{Z}(K_i)}$, the operator $\sigma_i(k_i)$ is not a multiple of the identity operator. Consider the natural projection $p_i:K_i\rightarrow K_i/\mathcal{Z}(K_i)$, let $\rho_i=\sigma_i\circ p_i$ be the inflation of $\sigma_i$ to $K_i$ and observe that the 
		$$T_i=\{t\in K_i:\rho(t)\mbox{ is a multiple of the identity operator}\}$$
		is a closed normal subgroup of $K_i$ containing $\mathcal{Z}(K_i)$. In particular, as $K_i$ is a simple Lie group and since $\sigma_i$ is non-trivial, $T_i$ must be discrete. Since $K_i$ is connected, it follows that $T_i$ is central so that $\sigma_i(k_i)$ is not a multiple of the identity operator for each $k_i\in K_i/\mathcal{Z}(K_i)$, $k_i\not=1_{K_i/\mathcal{Z}(K_i)}$. Now, Lemma \ref{Lemma Bekka le produit tensoriel de 2 operateur de Hilbert est mult de l'identité ssi les 2 le sont} ensures that the irreducible representation $\sigma_1\otimes \cdots \otimes\sigma_n$ of $$ K_1/_{\mathcal{Z}(K_1)}\times \cdots \times K_n/_{\mathcal{Z}(K_n)}$$ is faithful.
	\end{proof}
	\begin{corollary}\label{corollary the semi-simple Lie group is irreducibly represented}
		Let $G$ be a connected Lie group. Then, $\mathcal{S}_s(G)$ has a faithful and irreducible representation. 
	\end{corollary}
	The assertion \ref{item the subgroup Ss} follows from Proposition \ref{proposition the subgroup spanned by the semi simple or abelian feet is a closed connectedted normal ss subgroup} and Corollary \ref{corollary the semi-simple Lie group is irreducibly represented}.
	
		\subsection{The center alternative and the abelian Lie socle}\label{Section abelian Lie socle}
	
	Let $G$ be a connected Lie group. We recall that an abelian Lie foot $V$ of $G$ is a minimal non-trivial closed connected abelian normal subgroup of $G$ intersecting the center of $G$ trivially (minimal in the sense that the only subgroups of $V$ that is non-trivial connected and normal in $G$ is $V$ itself). In addition, we recall that the abelian Lie socle $\mathcal{S}_a(G)$ of $G$ is defined as the closed subgroup of $G$ generated by the set of abelian Lie feet of $G$ and its center. The purpose of this section is to describe its structure for groups that have a chance to be irreducibly represented. 
	
	We start with an important observation concerning the center of an irreducibly represented Lie group. We recall, as explained in \cite{BekkadeleHarpe2008}, that as a consequence of Schur’s lemma the center of an irreducibly represented group embeds in the circle group. In the context of Lie groups, this provides a simple alternative. The following result proves the assertion \ref{item alternative du centre}. 
	\begin{lemma}[\tg{center alternative}]\label{lemma property Harpe}
		Let $G$ be an irreducibly represented Lie group. Then, the center $\mathcal{Z}(G)$ of $G$ is either discrete and countable or connected and isomorphic as a Lie group to the circle group $\mathbb{T}$.
	\end{lemma}
	\begin{proof}
		Since $\mathcal{Z}(G)=\bigcap_{g\in G} \Ker(c_g)$ where $c_g(h)=ghg^{-1}$, the center $\mathcal{Z}(G)$ is a closed subgroup of $G$ and therefore inherits a structure of a Lie group. If $\dim(\mathcal{Z}(G))=0$, then $\mathcal{Z}(G)$ is discrete and countable since $G$ is second-countable. On the other hand,  if $\dim(\mathcal{Z}(G))\geq 1$, let $\mathcal{Z}(G)^\circ$ be the connected component of the identity in $\mathcal{Z}(G)$ and let $\pi$ be an irreducible and faithful representation of $G$. By Schur's Lemma, $\pi(\mathcal{Z}(G))$ consists of scalar multiples of the identity $\id_{\Hr{\pi}}$ by a unit complex number. Since $\pi$ is a continuous group homomorphism, this implies that $\pi(\mathcal{Z}(G)^\circ)$ is a connected subgroup of the circle group $\mathbb{T}$. Furthermore, as $\pi$ is faithful, this subgroup is non-trivial. Since there is no non-trivial proper connected subgroup of $\mathbb{T}$, we obtain that $$\pi(\mathcal{Z}(G)^\circ)=\{s \id_{\Hr{\pi}}: s\in \C,\qq \modu{s}=1\}=\pi(\mathcal{Z}(G)).$$ It follows from the injectivity of $\pi$ that $\mathcal{Z}(G)=\mathcal{Z}(G)^\circ$. Finally, the classification of connected abelian Lie groups ensures that $\mathcal{Z}(G)$ is isomorphic as a Lie group to the circle group $\mathbb{T}$.
	\end{proof}
	\begin{definition}
		A Lie group is said to satisfy the \tg{center alternative} if its center is either discrete or isomorphic as a Lie group to the circle group.
	\end{definition} 
	The purpose of the rest of this section is to prove the following result implying the assertion \ref{item structure of the abelian lie socle}.
	\begin{theorem}\label{Theorem LG est characteristic et isomoprhe a un espace vectoriel}
		Let $G$ be a connected Lie group satisfying the center alternative. Then, there exists a closed connected topologically characteristic subgroup $\mathcal{V}(G)$ of $G$ isomorphic as a Lie group to a finite dimensional real vector space such that $$\mathcal{S}_a(G)\simeq \mathcal{V}(G)\times \mathcal{Z}(G) \q(\mbox{with the natural inclusions}).$$ 
	\end{theorem} 
	\begin{definition}\label{definition linear socle}
		The closed subgroup $\mathcal{V}(G)$ provided by the above theorem is called the \tg{linear socle} of $G$.
	\end{definition}
	The proof of Theorem \ref{Theorem LG est characteristic et isomoprhe a un espace vectoriel} requires some preliminaries. 
	\begin{lemma}\label{lemma les sous groups comp con norm sont centraux}
		Let $G$ be a connected Lie group and $K\leq G$ be a compact connected abelian normal subgroup. Then, $K$ is central. 
	\end{lemma}
	\begin{proof}
		Since $K$ is compact it is closed in $G$ and hence defines a Lie group. By the classification of connected abelian Lie groups, $K$ is isomorphic to the direct product $\mathbb{T}^n=\mathbb{T}\times \cdots \times \mathbb{T}$ of $n$ copies of the circle group for some positive integer $n\in \N$. In particular, its group of automorphisms is discrete as $\Aut(K)\simeq\Aut(\mathbb{T}^n)\simeq \Gl_n(\Z)$. On the other hand, since $K$ is normal, the map $\phi: G\rightarrow \Aut(K)$ sending $g$ to the restriction of the conjugation by $g$ to $K$ is a continuous group homomorphism. It follows from the connectedness of $G$ that $\phi$ has trivial image and hence that $K$ is central.
	\end{proof}
	By contrast, we have the following result.
	\begin{lemma}\label{Lemma les Rs intersecte le center triviallement}
		Let $G$ be a connected Lie group satisfying the center alternative and $H$ be a closed abelian subgroup isomorphic as a topological group to some finite dimensional real vector space. Then, $H\cap \mathcal{Z}(G)=\{1_G\}$.  
	\end{lemma}
	\begin{proof}
		The group of automorphism $\Aut(V)$ of a finite dimensional real vector space $V$ coincides with the linear group $\Gl(V)$. In particular, under the correspondence provided by our hypothesis, for all $g\in G$ the restriction of the conjugation by $g$ to $H$ is $\R$-linear. Suppose for a contradiction that there exists an element $h\in H\cap \mathcal{Z}(G)$. In particular, $ghg^{-1}=h$ for all $g\in G$ and every element of the real vector space $\R h$ spanned by $h$ in $H$ is fixed by the conjugation maps by linearity. The subgroup corresponding to $\R h$ is therefore contained in the center $\mathcal{Z}(G)$. On the other hand, $G$ satisfies the center alternative so that its center does not contain any non-compact non-discrete uncountable subgroup. This is a contradiction. 
	\end{proof} 
	The following two technical lemmas provide the last missing pieces of the proof of Theorem \ref{Theorem LG est characteristic et isomoprhe a un espace vectoriel}. 
	\begin{lemma}\label{techinal lemma ensureing that the group LG is characteristic}
		Let $G$ be a topological group, $H\leq G$ be a closed subgroup, $K$ be a compact group, $V$ be a topological group without non-trivial compact subgroups and 
		$$\Psi: V\times K\rightarrow H$$
		be an isomorphism of topological groups. Then, every topological group isomorphism $\phi\in \Aut(G)$ that stabilises $H$ (in the sense that $\phi(H)\subseteq H$) also stabilises $\Psi(V\times\{1_K\})$ and $\Psi(\{1_V\}\times K)$.
	\end{lemma}
	\begin{proof}
		Let $\phi$ be an isomorphism of topological group stabilising $H$ and denote by $p_V:V\times K\rightarrow V$ and $p_K:V\times K\rightarrow V$ the two natural projection maps. Since $K$ is compact and since $\phi$ is stabilises by $H$, $\phi\big( \Psi(\{1_V\}\times K)\big)$ defines a compact subgroup of $G$  contained in $H$. In particular, notice that $p_V\circ\Psi^{-1}\big( \phi\big( \Psi(\{1_V\}\times K)\big)\big)$ defines a compact subgroup of $V$. It follows from our hypothesis that this group is trivial and hence that $$\phi(\Psi(\{1_V\}\times K))\subseteq \Psi(\{1_V\}\times K).$$
		Similarly, $p_K\circ \Psi^{-1} \big(\phi\big( \Psi(V\times \{1_K\})\big)\big)$ defines a closed subgroup of $K$ that does not contain any non-trivial compact subgroup. It follows that this group is trivial and hence that $ \phi\big( \Psi(V\times \{1_K\})\big)\subseteq \Psi(V\times \{1_K\})$.
	\end{proof}
	
	\begin{lemma}\label{Lemma le sous group engendré par un ev et un center discret est un produit direct trop cool ferme}
		Let $G$ be a connected Lie group with discrete center $\mathcal{Z}(G)$ and $V\leq G$ be a closed normal subgroup that is isomorphic as a Lie group to a finite dimensional vector space. Then, the subgroup generated by $V$and $\mathcal{Z}(G)$ is closed in $G$ and isomorphic as a Lie group to $V\times \mathcal{Z}(G)$.
	\end{lemma}
	\begin{proof}
		Since $\mathcal{Z}(G)$ is discrete and since $V\leq G$ is a closed normal subgroup isomorphic as a Lie group to a finite dimensional vector space Lemma \ref{Lemma les Rs intersecte le center triviallement} ensures that $V\cap \mathcal{Z}(G)=\{1_G\}$. In particular, the smooth map
		$$V\times \mathcal{Z}(G)\rightarrow G : (v,z)\mapsto vz$$
		is injective. Since the elements of $V$ and $\mathcal{Z}(G)$ commute with each other, this implies that the subgroup of $G$ generated by $V$ and $\mathcal{Z}(G)$ is $$V\mathcal{Z}(G)=\{vz: v\in V, z\in \mathcal{Z}(G)\}.$$
		In particular, provided that this group is closed,  the subgroup generated by $V$ and $\mathcal{Z}(G)$ is closed in $G$ and isomorphic as a Lie group to $V\times \mathcal{Z}(G)$. We are left to show that $V\mathcal{Z}(G)$ is closed in $G$. Let $H$ be the closure of $V\cdot\mathcal{Z}(G)$ in $G$. The connected component of the identity $H^\circ$ of $H$ is a connected abelian Lie group. Hence, by the classification of connected abelian Lie groups, there exists a finite dimensional real vector space $W$ and a Lie group isomorphism 
		$$\psi: W\times \mathbb{T}^r \rightarrow H^\circ.$$
		On the other hand, since $H^\circ$ is normal in $G$, Lemma \ref{techinal lemma ensureing that the group LG is characteristic} ensures that $\psi(\{1_W\}\times\mathbb{T}^r )$ is a normal subgroup of $G$. It follows from Lemma \ref{Theorem LG est characteristic et isomoprhe a un espace vectoriel} that $\psi(\{1_W\}\times\mathbb{T}^r )$ is central and hence that $r=0$ because $\mathcal{Z}(G)$ is discrete and countable. It follows from Lemma \ref{Lemma les Rs intersecte le center triviallement} that $H^\circ\cap \mathcal{Z}(G)=\{1_G\}$. Since $V\subseteq H^\circ$ we therefore obtain that $H^\circ\cap Vz=\es$ for all $z\in\mathcal{Z}(G)-\{1_G\}$. On the other hand, the connected component of the identity $H^\circ$ is open in $H$ so that $H^\circ$ coincides with the closure of $H^\circ\cap (V \mathcal{Z}(G))$. It follows that $H^\circ=V$. In particular, $V\mathcal{Z}(G)$ is an open subgroup of $H$ and is therefore closed in $H$. The result follows.
	\end{proof}
	We are finally ready to prove the theorem. 
	\begin{proof}[Proof of Theorem \ref{Theorem LG est characteristic et isomoprhe a un espace vectoriel}]
		Corollary \ref{Corollary the abelian Lie socle is a closed abelian topologically characteristic subgroup of G} ensures that $\mathcal{S}_a(G)$ is a closed abelian topologically characteristic subgroup of $G$. Its connected component of the identity $\mathcal{S}_a(G)^\circ$ is therefore a closed connected topologically characteristic abelian subgroup of $G$. Now, the classification of connected abelian Lie groups, ensures the existence of a finite dimensional real vector space $V$, a non-negative integer $r\in \N$ and a Lie group isomorphism 
		$$\Psi:  V\times \mathbb{T}^r \rightarrow \mathcal{S}_a(G)^\circ.$$
		It follows from Lemma \ref{techinal lemma ensureing that the group LG is characteristic} that both $\Psi(V\times \{1_K\})$ and $\Psi(\{0\}\times \mathbb{T}^r)$ are topologically characteristic subgroups of $G$. In particular, $\Psi(\{0\} \times \mathbb{T}^r)$ is a compact normal subgroup of $G$ and Lemma \ref{lemma les sous groups comp con norm sont centraux} ensures that this group is central. Since $G$ satisfies the center alternative we have two cases. 
		
		If $\mathcal{Z}(G)\simeq \mathbb{T}$, we deduce that $r=1$, $\Psi(\{0\}\times \mathbb{T}^r)=\mathcal{Z}(G)$ and hence that
		$$\mathcal{S}_a(G)\simeq \Psi(V\times\{1_K\}) \times \mathcal{Z}(G).$$ In particular, the result follows by posing $\mathcal{V}(G)=\Psi(V\times \{1_K\})$. 
		
		If on the other hand, $\mathcal{Z}(G)$ is discrete, we deduce that $r=0$ so that $\mathcal{S}_a(G)^\circ$ is isomorphic as a Lie group to a real vector space. 
		It follows from \ref{Lemma le sous group engendré par un ev et un center discret est un produit direct trop cool ferme} that the subgroup generated by $\mathcal{S}_a(G)^\circ$ and $\mathcal{Z}(G)$ is closed in $G$ an isomorphic as a Lie group to 
		$$\mathcal{S}_a(G)^\circ \times \mathcal{Z}(G).$$
		On the other hand, Proposition \ref{proposition the subgroup spanned by the semi simple or abelian feet is a closed connectedted normal ss subgroup} ensures that the closed subgroup generated by the elements of $\mathcal{A}_G$ is connected and is therefore contained in $\mathcal{S}_a(G)^\circ$. It follows that the subgroup generated by $\mathcal{S}_a(G)^\circ$ and $\mathcal{Z}(G)$ is $\mathcal{S}_a(G)$. In particular, the result follows by posing $\mathcal{V}(G)=\mathcal{S}_a(G)^\circ$. 
	\end{proof}
	
	\subsection{A regular embedding of the abelian Lie socle}\label{subsection regular embedding of the abelian Lie socle}
	
	The purpose of this section is to prove the assertion \ref{item if G is algebraix the abelian socle is reg embedded}. More precisely, we prove the following result.
	\begin{proposition}\label{proposition the ab elian Lie socle is regularily embedded} 
		Let $G^\circ$ be the Hausdorff connected component of the identity of the $\R$-points $G$ of a linear algebraic group defined over $\R$ and suppose that $G^\circ$ satisfies the center alternative. Then, $\mathcal{S}_a(G^\circ)$ is regularly embedded in $G^\circ$.   
	\end{proposition}
	In light of Theorem \ref{Theorem LG est characteristic et isomoprhe a un espace vectoriel}, the abelian Lie socle $\mathcal{S}_a(G^\circ)$ of $G^\circ$ is a direct product of its linear socle $\mathcal{V}(G^\circ)$ and its center $\mathcal{Z}(G^\circ)$. The following result provides a sufficient criterion for $\mathcal{S}_a(G^\circ)$ to be regularly embedded.
	\begin{lemma}\label{Lemma the direct product of reg embedded is reg embedded if one is central}
		Let $G$ be a locally compact group, $H_1\times H_2\leq G$ be a closed normal subgroup such that $H_1\leq G$ is a closed abelian regularly embedded normal subgroup and $H_2\leq G$ is a closed central subgroup of $G$. Then, $H_1\times H_2$ is regularly embedded in $G$.
	\end{lemma}
	\begin{proof}
		Let $H=H_1\times H_2$ and denote by $$p_{H_i}:H\rightarrow H_i$$
		the natural projection map. Since abelian groups are of Type {\rm I}, Theorem \ref{Theorem Mackey rep of the product} ensures that the map
		$$\Phi : \widehat{H_1}\times \widehat{H_2}\rightarrow \widehat{H}: \chiup_1\times \chiup_2\mapsto \chiup_1\otimes \chiup_2.$$
		is bijective. We are going to show that this map is open (in fact it is a homeomorphism). Recall that for any locally compact abelian group $H$ a basis of the fell topology of $\widehat{H}$ is given by sets of the form 
		$$ U_{\widehat{H}}(\chiup, K,\varepsilon)=\big\{\sigma \in \widehat{H}: \sup_{k\in K}\modu{\chiup(k)-\sigma(k)}<\varepsilon\big\}$$
		with $\chiup \in \widehat{H}$, $K\subseteq H$ compact and $\varepsilon>0$. To show that $\Phi$ is open, we consider for all for $i\in \{1,2\}$ some $\chi_i\in \widehat{H}_i$, $K_i\subseteq H$ compact, $\varepsilon_i>0$ and we let $$U=\Phi( U_{\widehat{H}_1}(\chiup_1, K_1,\varepsilon_1)\times  U_{\widehat{H}_2}(\chiup_2, K_2,\varepsilon_2)).$$ Our purpose is to prove that $U$ is an open subset of $\widehat{H}$.
		Let $\chiup'_1\otimes \chiup'_2\in U$. There exists two scalars $\delta_1,\delta_2>0$ such that
		$$\max_{k_i\in K_i}\qq \modu{\chiup_i(k_i)-\chiup_i'(k_i)}+\delta_i<\varepsilon_i \q \forall  i\in \{1,2\}.$$
		Let $\delta=\min\{\delta_1,\delta_2\}$, let $K$ be a compact subset of $H$ containing $K_1\times \{1_{K_2}\} \cup\{1_{K_1}\}\times  K_2$ and let $\sigma_1\otimes \sigma_2\in U_{\widehat{H}}(\chiup'_1\otimes \chiup'_2, K,\delta).$ We notice from a  straight forward computation that for all $t\in K_1$ one has
		\begin{equation}
			\begin{split}
				\modu{\chiup_1(t)-\sigma_1(t)}&\leq \modu{\chiup_1(t)-\chiup_1'(t)}+\modu{\chiup_1'(t)-\sigma_1(t)}\\
				&\leq \max_{k\in K_1} \big(\modu{\chiup_1(k)-\chiup_1'(k)}\big) + \modu{\chiup_1'(k)\chiup_2'(1_G)-\sigma_1(k)\sigma_2(1_G)}\\
				&\leq \max_{k\in K_1} \big(\modu{\chiup_1(k)-\chiup_1'(k)}\big) + \delta \leq \varepsilon_1.
			\end{split}
		\end{equation}
		The same inequalities hold inverting the indices proving that $\sigma_1\otimes \sigma_2\in U$. It follows that $U_{\widehat{H}}(\chiup'_1\otimes \chiup'_2, K,\delta)\subseteq U$ and $\Phi$ is open. Now, notice that $H_2$ is regularly embedded in $G$ by Lemma \ref{lemma a central subgroup is regularily embedded}. In particular, for all $i\in \{1,2\}$ there is a countable family $\mathcal{B}_i$ of $G$-stable Borel subsets of $\widehat{H}_i$ such that each $G$-orbit of $\widehat{H}_i$ is the intersection of the elements of $\mathcal{B}_i$ in which it is contained. Since $\Phi$ is open, $\mathcal{B}=\{\Phi(B_1\times B_2):B_1\in \mathcal{B}_1, B_2\in \mathcal{B}_2\}$ is a countable family of $G$-stable Borel subsets of $\widehat{H}$. Finally, since $H_2$ is central the $G$-orbits in $\widehat{H}$ are of the form $\Phi(\mathcal{O}_1\times \{\sigma_2\})$ for some $G$-orbit $\mathcal{O}_1$ in $\widehat{H}_1$. It follows that every $G$-orbit in $\widehat{H}$ is the intersection of the element of $\mathcal{B}$ in which it is contained.
	\end{proof}
	In light of this result, the following theorem provides the last missing piece required to prove Proposition \ref{proposition the ab elian Lie socle is regularily embedded}.
	\begin{theorem}[{\cite[Corollary 9]{BekkaEchterhoff2021}}]\label{Thoerem bekka alg reg ebdedded}
		Let $k$ be a local field of characteristic $0$, $G$ be a linear algebraic group defined over $k$, $G\times V \mapsto V$ be a $k$-rational action of $G$ on a finite dimensional vector space $V$ over $k$ and $G(k)$ be the $k$-rational points of $G$. Then every $G(k)$-orbit in $\widehat{V}$ is locally closed for the Hausdorff topology. 
	\end{theorem}
	\begin{proof}[Proof of Proposition \ref{proposition the ab elian Lie socle is regularily embedded}]
		Theorem \ref{Theorem LG est characteristic et isomoprhe a un espace vectoriel} ensures that 
		$$\mathcal{S}_a(G^\circ)\simeq \mathcal{V}(G^\circ)\times \mathcal{Z}(G^\circ) \q(\mbox{with the natural inclusions}).$$
		In particular, in light of Lemma \ref{Lemma the direct product of reg embedded is reg embedded if one is central}, it is enough to show that $\mathcal{V}(G^\circ)$ is regularly embedded in $G^\circ$. Due to Theorem \ref{Theorem de zimmer etre reg embedded c'est avoir des orbites loca closed}, this is equivalent to prove that each $G^\circ$-orbit is locally closed in $\widehat{\mathcal{V}(G^\circ)}$. On the other hand, as $\mathcal{V}(G^\circ)$ is a topologically characteristic subgroup of $G^\circ$ and since $G^\circ$ is an open normal subgroup of $G$, $\mathcal{V}(G^\circ)$ is a closed normal subgroup of $G$. Now, Theorem \ref{Thoerem bekka alg reg ebdedded} ensures that every $G$-orbit in $\widehat{\mathcal{V}(G^\circ)}$ is locally closed. In particular, by Theorem \ref{Theorem de zimmer etre reg embedded c'est avoir des orbites loca closed}, for every $\chiup \in \widehat{\mathcal{V}(G^\circ)}$, the natural map $G/\Fix_G(\chiup)\mapsto G\cdot \chiup$ is a homeomorphism for the relative topology of $G\cdot \chiup$ as a subset of $\widehat{\mathcal{V}(G^\circ)}$. Since the quotient map $G \rightarrow G/\Fix_G(\chiup): g\mapsto g\Fix_G(\chiup)$ is open, this implies that $G^\circ \cdot \chiup$ is open in $G\cdot \chiup$. It follows that $G^\circ \cdot \chiup$ is locally closed in $\widehat{\mathcal{V}(G^\circ)}$ because $G\cdot \chiup$ is locally closed in $\widehat{\mathcal{V}(G^\circ)}$. The result follows.
	\end{proof}

		\subsection{Existence of a $G$-faithful representation of the abelian Lie socle}
	The purpose of this section is to provide necessary and sufficient criteria for the abelian Lie socle to admit a $G$-faithful representation. In particular, we prove the assertions \ref{item the moment Sa G admits a G faitfull representaiton }, \ref{item the moment LG admits a G faitfull representaiton } and \ref{item the moment cenre faitfull representaiton }. 
	
	We start with the following result implying the assertion \ref{item the moment Sa G admits a G faitfull representaiton }.
	\begin{lemma}\label{lemma Sa has a G faith iff V (G) and Z are G irrep}
		Let $G$ be a connected Lie group satisfying the center alternative. Then $\mathcal{S}_a(G)$ has a $G$-faithful irreducible representation if and only if $\mathcal{V}(G)$ has a $G$-faithful irreducible representation and $\mathcal{Z}(G)$ is irreducibly represented.
	\end{lemma}
	\begin{proof}
		Theorem \ref{Theorem LG est characteristic et isomoprhe a un espace vectoriel} ensures that the abelian Lie socle decomposes as a direct product $\mathcal{V}(G)\times\mathcal{Z}(G)$ of two abelian (hence Type {\rm I}) groups. In particular, Theorem \ref{Theorem Mackey rep of the product} ensures that
		$$\widehat{\mathcal{S}_a(G)}=\{\chiup_1\otimes \chiup_2: \chiup_1\in \widehat{\mathcal{V}(G)},\qq \chiup_2\in \widehat{\mathcal{Z}(G)}\}.$$
		On the other hand, since $\mathcal{V}(G)$ is normal and as $\mathcal{Z}(G)$ is central in $G$, notice that 
		$$(\chiup_1\otimes \chiup_2)^g=\chiup_1^g\otimes \chiup_2^g= \chiup_1^g\otimes \chiup_2\q \forall g\in G.$$
		In particular, realising $\mathcal{S}_a(G)$ as a direct product we have that 
		$$\bigg(\bigcap_{g\in G} \Ker(\chiup_1^g)\bigg)\times\Ker(\chiup_2)\subseteq \bigcap_{g\in G}\Ker((\chiup_1\otimes \chiup_2)^g).$$
		If $\chiup_1\otimes \chiup_2$ is $G$-faithful, it follows that $\chiup_1$ is a $G$-faithful irreducible representation of $\mathcal{V}(G)$ and that $\chiup_2$ is a faithful irreducible representation of $\mathcal{Z}(G)$. To prove the other implication, let $\chiup_1$ be a $G$-faithful irreducible representation of $\mathcal{V}(G)$ and $\chiup_2$ be a faithful irreducible representation of $\mathcal{Z}(G)$. Now, notice that for every $v\in \mathcal{V}(G)-\{0\}$ there exists an element $g\in G$ such that 
		$$\chiup_1^g(v)\not=\chiup_1(v).$$
		Suppose for a contradiction that this is not the case, that is, suppose that $gvg^{-1}-v\in \Ker(\chiup_1)$ for all $g\in G$. Since $\mathcal{V}(G)$ intersects the center trivially notice that $gvg^{-1}-v\not=0$ for at least one $g\in G$. In particular, the subgroup $H$ generated by $\{gvg^{-1}-v: g\in G\}$ is a non-trivial normal subgroup of $G$ contained in $\mathcal{V}(G)$. However, by construction we have that $H\subseteq \Ker(\chiup_1)$ and hence by normality that $H\subseteq \bigcap_{g\in G}\Ker (\chiup_1^g)$. This provides the desired contradiction since $\chiup_1$ is $G$-faithful. It follows that 
		$$\bigcap_{g\in G}\Ker((\chiup_1\otimes \chiup_2)^g)= \bigcap_{g\in G}\Ker(\chiup_1^g\otimes \chiup_2)\subseteq \bigg(\bigcap_{g\in G} \Ker(\chiup_1^g)\bigg)\times\Ker(\chiup_2)=\{1_G\}.$$
	\end{proof}
	
	Our next purpose is to determine when the linear socle $\mathcal{V}(G)$ admits a $G$-faithful irreducible representation. We recall that the linear socle $\mathcal{V}(G)$ of $G$ is isomorphic as a Lie group to a finite dimensional real vector space. In particular, the linear dual $$\mathcal{V}(G)^*=\{f: \mathcal{V}(G) \rightarrow \R: f \mbox{ is }\R\mbox{-linear}\}$$ of $\mathcal{V}(G)$ identifies with its unitary dual $\widehat{\mathcal{V}(G)}$ under the isomorphism of topological groups 
	$$\Phi: \mathcal{V}(G)^*\rightarrow \widehat{\mathcal{V}(G)}: f\mapsto e^{2\pi i f(\cdot)}.$$
	Under this correspondence, the conjugation map provides a continuous linear (non-unitary) representation $$\rho: G\rightarrow \Gl(\mathcal{V}(G)^*)$$
	where $$\rho(g)f(v)=f(g^{-1}vg)\q \forall g\in G,\qq  v\in \mathcal{V}(G)\mbox{ and }f\in \mathcal{V}(G)^*.$$
	This linear representation is called the \tg{coadjoint} representation of $G$ on $\mathcal{V}(G)^*$ and an orbit of this $G$-action on $\mathcal{V}(G)^*$ is called a \tg{coadjoint} $G$\tg{-orbit}. The following result proves the assertion \ref{item the moment LG admits a G faitfull representaiton }.
	\begin{proposition}\label{Proposition, the lienar socle is G irred rep if if coadjoint orbit generates}
		Let $G$ be a connected Lie group satisfying the center alternative. Then, $\mathcal{V}(G)$ admits a $G$-faithful irreducible representation if and only if $\mathcal{V}(G)^*$ is spanned as a real vector space by a single co-adjoint $G$-orbit. 
	\end{proposition}
	\begin{proof}
		We denote by $\Theta$ the inverse of the isomorphism of topological groups 
		$$\Phi: \mathcal{V}(G)^*\rightarrow \widehat{\mathcal{V}(G)}: f\mapsto e^{2\pi i f(\cdot)}.$$ Now, let $\chiup \in \widehat{\mathcal{V}(G)}$. Due to the form of $\chiup$ notice that $v\in \Ker(\chiup)$ if and only if $\Theta(\chiup)(v)\in \Z$ so that $\Ker(\chiup)^\circ=\Ker(\Theta(\chiup))$. On the other hand, we recall that the co-adjoint representation $\rho$ of $G$ is defined so that $$\rho(g)\Theta(\chiup)=\Theta(\chiup^g) \q \forall \chiup\in \widehat{\mathcal{V}(G)},\qq \forall g\in G.$$
		It follows that 
		\begin{equation}\label{egalitee entre le G noyou et le noyau de forme lineaire}
			\bigg(\bigcap_{g\in G}\Ker(\chiup^g)\bigg)^\circ=\bigcap_{g\in G} \Ker(\chiup^g)^\circ=\bigcap_{g\in G}\Ker(\rho(g) \Theta(\chiup)).
		\end{equation} 
		In particular, if $\chiup$ is $G$-faithful we obtain that $\bigcap_{g\in G}\Ker(\rho(g) \Theta(\chiup))=0$. In that case, $\mathcal{V}(G)^*$ is spanned as a real vector space by $\{\rho(g) \Theta(\chiup): g\in G\}$.
		
		On the other hand, if there exists a function $f\in\mathcal{V}(G)^*$ such that $\mathcal{V}(G)^*$ is spanned as a real vector space by $\{\rho(g) f: g\in G\}$, one has that $\bigcap_{g\in G}\Ker(\rho(g) f)=0$. Now, let $\chiup = \Theta^{-1}(f)$ and notice from \eqref{egalitee entre le G noyou et le noyau de forme lineaire} that 
		$$\Big(\bigcap_{g\in G}\Ker(\chiup^g)\Big)^\circ=0.$$
		It follows that $\bigcap_{g\in G}\Ker(\chiup^g)$ is a discrete subgroup of $G$ contained in $\mathcal{V}(G)$. On the other hand, Lemma \ref{lemma discrete groups are central} ensures that such a group is central and Lemma \ref{Lemma les Rs intersecte le center triviallement} ensures that $\mathcal{V}(G)\cap \mathcal{Z}(G)=\{1_G\}$. It follows that $\bigcap_{g\in G}\Ker(\chiup^g)$ is trivial so that $\chiup$ is $G$-faithful.
	\end{proof}
	We are left to determine when $\mathcal{Z}(G)$ is irreducibly represented. Since the circle group is irreducibly represented we are left to treat the case where $\mathcal{Z}(G)$ is discrete. In that case, we rely on the algebraic characterisation of countable discrete groups that are irreducibly represented:
	\begin{theorem}[{\cite[Corollary $1.3$ and Corollary $1.9$]{CapracedelaHarpe2020}}]
		A countable discrete abelian group $G$ is irreducibly represented if and only if it does not contain a subgroup isomorphic to $C_p\times C_p$ for any prime $p$.
	\end{theorem}
	The following result follows proving the assertion \ref{item the moment cenre faitfull representaiton }. 
	\begin{corollary}
		Let $G$ be a connected Lie group satisfying the center alternative. Then, $\mathcal{Z}(G)$ is not irreducibly represented if and only if it is discrete and contains a subgroup isomorphic to $C_p\times C_p$ for some prime $p$.
	\end{corollary}
	
	\part{Harmonic analysis of groups of tree automorphisms}
	\newpage
	\thispagestyle{empty}
	\mbox{}
	\newpage
	\chapter{Totally disconnected locally compact groups}\label{Preliminaries on totally disconnected locally compact groups and representations}
		
	\section{Representations of totally disconnected locally compact groups}\label{section tdlc groups aand rep}
	
	The purpose of this chapter is to highlight the structural particularities of totally disconnected locally compact groups and the direct consequences on their representation theory. We start with a proper introduction on the concept of totally disconnected locally compact group. We recall that a (non-empty) topological space is \tg{connected} if the only clopen subsets are the whole space and the empty set. A topological space is said to be \tg{totally disconnected} if its only connected components are the singletons. 
	\begin{definition}
		A \tg{totally disconnected locally compact group} is a topological group whose topology is Hausdorff, locally compact and totally disconnected.
	\end{definition} 
	We recall that the Lie groups are characterised among the locally compact groups by the properties of having no small subgroups, see Theorem \ref{Theorem Lie groups are the locally comapct groups without small subgroups}. By contrast, totally disconnected locally compact groups admit a large amount of ``small'' subgroups. To be more precise, they are characterised among the locally compact groups by the property of having a basis of neighbourhoods of the identity consisting of compact open subgroups. This theorem is best known as the fundamental theorem of Van Dantzig. 
	\begin{theorem}[Van Dantzig's Theorem]
		A locally compact group is totally disconnected if and only if it admits a basis of neighbourhoods of the identity consisting of compact open subgroups.
	\end{theorem}
	\begin{proof}
		Suppose that $G$ has a basis of neighbourhoods of the identity consisting of compact open subgroups. Since $G$ is Hausdorff, each compact subset of $G$ is also closed. In particular, this basis of neighbourhoods of the identity consists of clopen subsets of $G$. It follows that every non-singleton set contains a proper clopen subset for the relative topology given by the inclusion. Such a set is therefore not connected and the result follows.   
		
		For the other implication, we base our reasoning on \cite[Lemma 6.1.3]{Tao2014} which ensures that for every totally disconnected compact Hausdorff space $X$ and for every two distinct points $x,y\in X$ there exists a clopen subset $U\subseteq X$ containing $x$ but not $y$. Now, consider an open neighbourhood $K$ of the identity with compact closure and let us show the existence of a compact open subgroup $G'$ of $G$ contained in $K$. Since $K$ can be chosen arbitrarily small, the result will follow from this statement. From our choice of $K$, notice that $\partial K$ is a compact subset of $G$. In particular,  for every $k\in \partial K$, the Lemma ensures the existence of a clopen neighbourhood of the identity that avoids $k$. By compactness of $\partial K$, we obtain a clopen neighbourhood of $1_G$ that avoids $\partial K$. Intersecting this neighbourhood with $K$ we obtain a compact clopen neighbourhood $F$ of $1_G$ that avoids $\partial K$ and that is contained in $K$. On the one hand, since $F$ is open, the continuity of the multiplication map $m:G\times G \rightarrow G$ ensures the  existence of an open neighbourhood $U_f$ of the identity and an open neighbourhood $V_f$ of $f$ in $F$ such that $U_f\times V_f\subseteq m^{-1}(F)$ for all $f\in F$. On the other hand, as $F$ is compact, it can be covered by finitely many of those $V_f$ with $f\in F$. Now, notice that by construction, the intersection of the corresponding $U_f$'s provides an open neighbourhood $U$ of the identity such that $U F\subseteq F$. Furthermore, by continuity of the inverse map, $U$ can be chosen to be symmetric. The subgroup $G'$ generated by $U$ is therefore open in $G$ and stabilises $F$. In particular, $G'$ does not contain any element of $F^c$. This proves that $U\subseteq F$. The compactness of $G'$ follows as open subgroups are closed and since $F$ is compact.
	\end{proof}
	We now provide a list of examples. 
	\begin{example}
		Let $\Gamma$ be a discrete group. Such a group is clearly totally disconnected. Furthermore, notice that the set $\{\{1_\Gamma\}\}$ provides a basis of neighbourhoods of the identity consisting of compact open subgroups of $\Gamma$.
	\end{example}
	\begin{example}
		Let $G$ be a locally compact group.  The connected component $H$ of the identity $1_G$ is a closed normal subgroup of $G$. In particular, the quotient $G/H$ is a locally compact group for the quotient topology. On the other hand, the connected component of the identity is trivial in $G/H$ so that this group is totally disconnected. Notice further that $G/H$ is discrete if and only if $H$ is open.
	\end{example}
	\begin{example}\label{example profinite groups}
		Let $(\mathcal{I},\leq )$ be a countable directed set, consider a family of finite groups $\{G_{i}:i\in  \mathcal{I}\}$ equipped with the discrete topology, and a collection of homomorphisms $\{f_{i}^{j}:G_{j}\to G_{i}\mid i,j\in \mathcal{I},i\leq j\}$ satisfying the composition property $f_{i}^{j}\circ f_{j}^{k}=f_{i}^{k}$ for all $ i\leq j\leq k$ and such that $f_{i}^{i}$ is the identity on $G_i$. The projective limit of the $G_i$'s is the topological group
		$$\underset{\longleftarrow}{\mbox{\rm lim}} \qq G_i=\bigg\{(g_{i})_{i\in  \mathcal{I}}\in \prod _{i\in  \mathcal{I}}G_{i}:f_{i}^{j}(g_{j})=g_{i}{\text{ for all }}j\geq i\bigg\}$$
		equipped with the relative product topology. It is not hard to see that $\underset{\longleftarrow}{\mbox{\rm lim}} \qq G_i$ is closed inside $\prod _{i\in  \mathcal{I}}G_{i}$. It follows from Tychonoff's theorem that this group is compact and hence locally compact. A group $G$ of this form is called \tg{profinite}.  For each $j\in \mathcal{I}$, we let $p_j: G \rightarrow  \prod _{i\in  \mathcal{I}, i\leq j}G_{i}$ be the natural projection map. It is not hard to see that the sets $\{\Ker(p_i):i\in \mathcal{I} \}$ form a basis of neighbourhoods of the identity consisting of normal compact open subgroups of $G$. It follows that $G$ is totally disconnected. Concrete examples of profinite groups are given in Examples \ref{example p adic numbers are tdlc} and \ref{example automorphism group of graphs are tdlc} below.
	\end{example}
	\begin{example}\label{example p adic numbers are tdlc}
		Let $p$ be prime, consider the $p$-adic norm $\norm{\cdot}{p}$, the associated $p$-adic field $\Q_p$ and the corresponding set of $p$-adic integers $$\Z_p=\{q\in \Q_p:\norm{q}{p}\leq 1\}.$$ For every positive integer $n\in \N$, notice that $p^n\Z_p$ is both the closed ball of radius $p^{-n}$ and the open ball of radius $p^{-n+1}$ for the metric given by the $p$-adic norm.  It follows that $p^n\Z_p$ is a compact and open subset of $\Q_p$. On the other hand, due to the ultrametric equality each $p^n\Z_p$ is an additive subgroup of $\Q_p$. In particular, the $p$-adic groups $\Q_p$ and $\Z_p$ are both totally disconnected locally compact groups. Furthermore, notice that $\Z_p$ can be realised as the projective limit of the $\{\Z/p^i\Z: i\in \N\}$ for the natural projection maps and hence defines a profinite group.
	\end{example}
	\begin{example}\label{example automorphism group of graphs are tdlc}
		A \tg{graph} $\mathcal{G}$ is a couple of sets $(V,E)$ (called respectively the set of \tg{vertices} and the set of \tg{edges}) such that ${E \subseteq \binom{V}{2}}$ where $\binom{V}{2}$ denotes the set of unordered pairs of distinct vertices $\{v,w\}\subseteq V$. Two vertices $v,w\in V$ are \tg{adjacent} in $\mathcal{G}$ if $\{v,w\}\in E$. A \tg{connected path} in $\mathcal{G}$ is a finite sequence of successively adjacent vertices. The first and the last vertex of a connected path are respectively called the starting and the terminal vertex of the path. The \tg{length} of a connected path is equal the length of the sequence less one (so that a path consisting of a single vertex has length $0$). A \tg{geodesic} of $\mathcal{G}$ is a connected path of minimal length among all the connected paths with the same starting and terminal vertices. The set of vertices $V$ is naturally equipped with a metric $d_\mathcal{G}$ given by length of geodesics. An \tg{automorphism} of $\mathcal{G}$ is a bijection $\varphi: V \rightarrow V$
		such that $\{v,w\}\in E$ if and only if $\{\varphi(v),\varphi(w)\}\in E$. Equivalently, an automorphisms of $\mathcal{G}$ is an isometry of the metric space $(V, d_{\mathcal{G}})$. The set $\Aut(\mathcal{G})$ of all the automorphisms of $\mathcal{G}$ form a group for the composition. This group embeds as a subgroup of $\Sym(V)$ and is therefore naturally equipped with a topology called the \tg{permutation topology}. For this topology, a basis of neighbourhoods of the identity of $\Aut(\mathcal{\mathcal{G}})$ is given by the pointwise stabilisers
		$$\Fix_{\Aut(\mathcal{G})}(F)=\{g\in \Aut(\mathcal{G}): gv=v\qq \forall v\in F \}$$
		where $F$ runs through the entire collection of finite subsets of $V$. A graph is \tg{locally finite} when each of its vertex belongs to finitely many edges. When $\mathcal{G}$ is locally finite the $n$-sphere $$S_n(v)=\{w\in V :d_{\mathcal{G}}(v,w)=n\}$$ is finite for every $v\in V$. On the other hand, the stabiliser $\Fix_{\Aut(\mathcal{G})}(v)$ is a closed subgroup of $$\prod_{n\in \N}\Sym(S_n(v)).$$ It follows from Tychonoff's theorem that $\Fix_{\Aut(\mathcal{G})}(v)$ is compact and hence that $\Aut(\mathcal{G})$ is a locally compact group. On the other hand, since each of the $\Fix_{\Aut(\mathcal{G})}(v)$ is also open, each of the $\Fix_{\Aut(\mathcal{G})}(F)$
		with $F\subseteq V$ finite is a compact open subgroup of $\Aut(\mathcal{G})$. It follows form the easy implication of Van Dantzig's Theorem that $\Aut(\mathcal{G})$ and all of the $\Fix_{\Aut(\mathcal{G})}(F)$ are a totally disconnected locally compact groups. On the other hand, notice that $\Fix_{\Aut(\mathcal{G})}(v)$ is a profinite group for each $v\in V$ since it can be realised as the projective limit of the groups $\{\Aut_{\mathcal{G}}(B(v,i)):i\in \N\}$ where $B(v,i)$ denotes the ball of radius $i$ with center $v$ in $\mathcal{G}$ and $\Aut_{\mathcal{G}}(B(v,i))$ denotes the automorphism group of $B(v,i)$ induced by the action of $\Fix_{\Aut(\mathcal{G})}(v)$. 
	\end{example}
	It is natural to ask what sort of informations the existence of an open neighbourhood of the identity consisting of compact open subgroups provides on the representation theory of these groups. The following lemma gives a first answer to this question. 
	\begin{lemma}\label{existence of U invariant vector}
		Let $\pi$ be a representation of $G$ and let $\mathcal{S}$ be a basis of neighbourhoods of the identity consisting of compact open subgroups. Then,
		\begin{equation*}
		\bigcup_{K\in \mathcal{S}}\Hr{\pi}^{K}= \{\xi \in \Hr{\pi}: \exists K\in \mathcal{S}\mbox{ such that }\qq \pi(k)\xi=\xi\qq \forall k\in K\}
		\end{equation*} 
		is dense in $\Hr{\pi}$.
	\end{lemma}	
	\begin{proof} 
		Let $\xi \in \Hr{\pi}$ be a non-zero vector and $\varepsilon>0$ be a positive scalar. Our purpose is to exhibit an element of $\bigcup_{K\in \mathcal{S}}\Hr{\pi}^{K}$ with distance from $\xi$ bounded by $\varepsilon$. By continuity of the matrix coefficient $$\varphi_{\xi,\xi}:G\rightarrow\C :g\mapsto\prods{\pi(g)\xi}{\xi}_{\Hr{\pi}},$$ since $\varphi_{\xi,\xi}(1_G)=\norm{\xi}{}^2>0$ and as $\mathcal{S}$ is a basis of neighbourhoods of the identity, there exists $H\in \mathcal{S}$ such that $\norm{\xi}{}^2-\mbox{Re}(\varphi_{\xi,\xi}(h)) < \varepsilon/2\qq \forall h\in H$. Now, let $\mu$ be the left-invariant Haar measure on $G$ renormalised in such a way that $\mu(H)=1$ (such a renormalisation exists since $H$ is both open and compact). We claim that $\int_{H}\pi(h)\xi \qq \diff \mu (h)$ is as desired. Indeed, this vector belongs to $\Hr{\pi}^{H}$ and a  straightforward computation shows that 
		\begin{equation*}
		\begin{split}
		\Big\lVert\xi-\int_H \pi(h)\xi \qq \diff \mu(h) \Big\rVert^2&\leq  \int_H \norm{\xi -\pi(h)\xi}{}^2\qq\diff \mu (h)\\
		& =2\int_H\norm{\xi}{}^2- \mbox{Re}(\varphi_{\xi,\xi}(h))\qq \diff \mu (h)< \varepsilon.
		\end{split}
		\end{equation*}
	\end{proof}
	In particular, the above lemma shows that every representation $\pi$ of a totally disconnected locally compact group $G$ admits a non-zero invariant vector for some compact open subgroup. Notice, however, that the statement looses its relevance for discrete groups. It is natural to ask whether further information on the representation theory can be obtained if the basis of neighbourhoods of the identity has more properties. The answer to this question is positive. We end this section with an easy example and develop our personal contributions to this problematic in Chapters \ref{Chapter Olshanskii's factor}, \ref{Chapter application olsh facto} and \ref{Chapter Radu groups} below. We recall from Example \ref{example representations of the quotient provides representations of the group} that every irreducible representation of a quotient $G/H$ of a locally compact group by a closed subgroup gives rise to an irreducible representation of $G$. The following corollary shows that every irreducible representation of $G$ arises in this manner if $G$ admits a basis of neighbourhoods of the identity $\mathcal{S}$ consisting of \textbf{normal} compact open subgroups. 
	\begin{lemma}\label{lemma each irrep of G is a rep of the quotient}
		Let $G$ be a locally compact group, $\mathcal{S}$ be a basis of neighbourhoods of the identity consisting of normal compact open subgroups and $\pi$ be an irreducible representation of $G$. Then, $\pi$ is the inflation of an irreducible representation of $G/K$ for some $K$ in $\mathcal{S}$. 
	\end{lemma}
	\begin{proof}
		Since $\mathcal{S}$ consists of compact open subgroups of $G$, Lemma \ref{existence of U invariant vector} ensures the existence of a compact open subgroup $K$ such that the subspace $\Hr{\pi}^K$ of $K$-invariant vectors of $\Hr{\pi}$ is non-zero. On the other hand, since $K$ is normal in $G$, the space $\Hr{\pi}^K$ is $G$-invariant. Indeed, for every $k\in K$, $g\in G$ and $\xi\in \Hr{\pi}^K$ one has $$\pi(k)\pi(g)\xi=\pi(g)\pi(g^{-1}kg)\xi=\pi(g)\xi.$$
		By irreducibility and since $\Hr{\pi}^K$ is closed, this implies that $\Hr{\pi}^K=\Hr{\pi}$ and hence that $\pi$ is an irreducible representation of the discrete quotient $G/K.$
	\end{proof} 
	This statement has a particular flavour for profinite groups.
	\begin{example}\label{example representations of profinite tdlc groups}
		Let $G=\underset{\longleftarrow}{\mbox{\rm lim}} \qq G_i$ be a profinite group as defined in Example \ref{example profinite groups} and let $\pi$ be an irreducible representation of $G$. Consider the natural projection maps $p_j: G \rightarrow  \prod _{i\in  \mathcal{I}, i\leq j}G_{i}$ and the corresponding basis of neighbourhoods of the identity $\mathcal{S}=\{\Ker(p_i):i\in \mathcal{I} \}$. By Lemma \ref{lemma each irrep of G is a rep of the quotient}, there exists some $i\in \mathcal{I}$ such that $\pi$ is the inflation of an irreducible representation of $G/\Ker(p_i)\simeq G_i$. To provide concrete examples:
		\begin{enumerate}[label=(\roman*)]
			\item Let $p$ be a prime number. Every irreducible representation of $\Z_p$ is the inflation of a character $\chiup$ of some finite group $\Z_p/p^n\Z_p\simeq \Z/p^n\Z.$
			\item Let $\mathcal{G}$ be a locally finite graph, $G=\Aut(\mathcal{G})$ and $v\in V$ be a vertex. Every irreducible representation of $\Fix_{\Aut(\mathcal{G})}(v)$ is the inflation of some irreducible representation of $$\Aut_G(B(v,n))\simeq\Fix_{G}(v)/\Fix_{G}(B(v,n)) $$ where $\Aut_G(B(v,n))$ is the finite automorphism group of the ball of radius $n$ centred at $v$ induced by the action of $\Fix_G(v)$.   
		\end{enumerate}
	\end{example}
	\section{Groups of automorphisms of trees}\label{section groups of automorphisms of trees}
	 Groups of automorphisms of trees are one of the most important examples of totally disconnected locally compact groups. The work presented in this part of the thesis focuses mostly on the representation theory of these groups. The purpose of this section is to provide the basic related concepts and some of the details of the classical classification of the irreducible representations of the full group of automorphisms of a thick regular tree.
	
	We adopt the same formalism for graphs as in Example \ref{example automorphism group of graphs are tdlc}. However, if an emphasis needs to be made on the graph we denote by $V(\mathcal{G})$ and $E(\mathcal{G})$ the respective sets of vertices and edges of $\mathcal{G}$. We recall that a \tg{connected path} in $\mathcal{G}$ is a finite sequence of successively adjacent vertices. A \tg{simple cycle} is a connected path of $\mathcal{G}$ of length at least $3$ with the same starting and terminal vertex but without other repetitions. A graph $\mathcal{G}$ is said to be \tg{connected} if for every pair $x,y\in V(\mathcal{G})$ of vertices of $\mathcal{G}$, there exists a connected path in $\mathcal{G}$ starting at $x$ and ending at $y$.
	\begin{definition}
		 A \tg{tree} $T$ is connected graph without simple cycles. Such a graph is said to be:
		 \begin{enumerate}[label=(\roman*)]
		 	\item \textbf{thick} if the degree of each vertex is at least $3$.
		 	\item $d$-\tg{regular} if each vertex $v\in V$ has degree $d$.
		 	\item $(d_0,d_1)$-\tg{semi-regular} if there exists a bipartition $V=V_0\sqcup V_1$ such that each vertex of $V_t$ has degree $d_t$ and every edge of $T$ contains exactly one vertex in each $V_t$ ($t=0,1$).
		 \end{enumerate}
	\end{definition}
	In addition, a vertex $v\in V$ of a tree $T$ is called a \tg{leaf} of $T$ if it has degree $1$ in $T$. In these notes, we principally work with locally finite trees without leaves. These trees are naturally associated with a visual boundary. To formalise this concept we start another series of definitions. A \tg{simply infinite geodesic} (respectively \tg{doubly infinite geodesic}) on $T$ is a sequence $(v_n)_{n\in \N}$ (resp. a sequence $(v_n)_{n\in \Z}$) of distinct and successively adjacent vertices $v_n\in V$. We equip the set of simply infinite geodesics of $T$ with an equivalence relation where two sequences $(v_n)_{n\in \N}$ and $(w_n)_{n\in \N}$ are equivalent if and only if there exists some $t\in \Z$ such that $v_n=w_{n+t}$ for all big enough $n\in \N$.
 	\begin{definition}
 		Let $T$ be a tree without leaves. The \tg{visual boundary} $\partial T$ of $T$ is the set of equivalence classes of simply infinite geodesics of $T$.
 	\end{definition}
 	We recall from Example \ref{example automorphism group of graphs are tdlc} that the vertex set $V$ of a tree $T$ is naturally equipped with a metric $d_T$ given by length of geodesics and that an isometry $g : V\rightarrow V$ for this metric is called an \tg{automorphism} of $T$. Furthermore, we recall that the set $\Aut(T)$ of all the automorphisms of $T$ forms a group for the composition. This groups acts naturally on $V$ and hence on the set of simply infinite geodesics. On the other hand, this action preserves the equivalence relation so that $\Aut(T)$ acts on the visual boundary $\partial T$. We will see later in this thesis that there is a surprising parallel between the regularity of the representation theory of closed subgroups and the way these groups act on the visual boundary $\partial T$. The following classical result shows that there are only three kinds of automorphisms of trees. 
 	 \begin{theorem}[{\cite[Theorem 3.2]{FigaNebbia1991}}]\label{theorem definitions des trois types d'automorphisms}
 	 	Let $T$ be a tree and $g\in \Aut(T)$. Then, exactly one of the following occurs.
 	 		\begin{enumerate}[label=(\roman*)]
 	 		\item There exists a vertex $v\in V$ such that $gv=v$.
 	 		\item There exists an edge $e\in E$  such that $ge=e$ but $gv\not=v$ $\forall v\in V$. 
 	 		\item There exists a doubly infinite geodesic $(v_n)_{n\in \Z}$ and a strictly positive integer $t\in \N$ such that $gv_n=v_{n+t}$ $\forall n \in \Z$.
 	 	\end{enumerate}
 	 \end{theorem}
  	This motivates the following definition.
 	\begin{definition}
 		Let $T$ be a tree. An automorphism of $T$ is called a \tg{rotation} if it fixes a vertex, it is called an \tg{inversion} of an edge $e$ if it stabilises $e$ without fixing any vertex and it is called a \tg{translation of step $t$} if it satisfies the third condition of Theorem \ref{theorem definitions des trois types d'automorphisms}.
 	\end{definition}
	Let $T$ be a locally finite and $G$ be a closed subgroup of the full group of automorphisms $\Aut(T)$. We recall from Example \ref{example automorphism group of graphs are tdlc} that $G$ is a totally disconnected locally compact group with a basis of neighbourhoods of the identity given by the compact open subgroups
	$$\Fix_{G}(F)=\{g\in G:gv=v\qq \forall v\in F \}$$
	where $F\subseteq V$ runs over all the finite subsets of $V$. In these notes we will often be interested in the subgraphs of trees that have the property of being a tree rather than general finite subsets. 
	\begin{definition}
		Let $T$ be a tree. A \tg{subtree} $\mathcal{T}$ of $T$ is a connected graph that is a subgraph of $T$ in the sense that $V(\mathcal{T})\subseteq V(T)$ and $E(\mathcal{T})\subseteq E(T)$. 
	\end{definition}
	 Notice that such a subtree is entirely determined by its set of vertices. Therefore, when it leads to no confusion, we identify $\mathcal{T}$ with its set of vertices. A subtree $\mathcal{T}$ of $T$ is said to be \tg{complete} if the degree in $\mathcal{T}$ of each vertex $v\in V(\mathcal{T})$ that is not a leaf of $\mathcal{T}$ coincides with its degree in $T$. Since each finite subset of $T$ is contained in a complete finite subtree of $T$, a basis of neighbourhoods of the identity is equivalently provided by sets of the form
	$$\Fix_{G}(\mathcal{T})=\{g\in G:gv=v\qq \forall v\in V(\mathcal{T})\}$$
	where $\mathcal{T}$ runs over all the complete finite subtrees of $T$. In particular, in light of Lemma \ref{existence of U invariant vector} this implies that every representation of $G$ admits a non-zero $\Fix_G(\mathcal{T})$-invariant vector for some complete finite subtree $\mathcal{T}$. This motivates the following terminology. 
	\begin{definition}\label{definiton of spherical special cuspidal}
		An irreducible representation $\pi$ of $G$ is called:
		\begin{enumerate}[label=(\roman*)]\label{definition de spheric special cuspidal}
			\item \tg{spherical} if there exists a vertex $v\in V$ such that $\pi$ admits a non-zero $\Fix_G(v)$-invariant vector where $\Fix_G(v)=\{g\in G: gv=v\}$.
			\item \tg{special} if it is not spherical and there exists an edge $e\in E$ such that $\pi$ admits a non-zero $\Fix_G(e)$-invariant vector where $\Fix_G(e)=\{g\in G: gv=v\qq\forall v\in e\}$ is the pointwise stabiliser of the edge $e$.
			\item \tg{cuspidal} if it is not spherical nor special. In that case, there nevertheless exists a complete finite subtree $\mathcal{T}$ of $T$ such that $\pi$ admits a non-zero $\Fix_G(\mathcal{T})$-invariant vector where $\Fix_G(\mathcal{T})=\{g\in G: gv=v\qq\forall v\in V(\mathcal{T})\}$.
		\end{enumerate}
	\end{definition}
	Every irreducible representation of $G$ is of exactly one of the above types. The rest of this chapter is dedicated to recalling the classification of the irreducible representations of the full group of automorphisms of a semi-regular tree. 
	The classification of spherical and special representations of any closed subgroup $G\leq \Aut(T)$ acting $2$-transitively on the boundary $\partial T$ is a classical result that was achieved in the $70$'s due to numerous contributors such as Godement, Cartier, Matsumoto. This classification is presented in Sections \ref{Section spherical rep} and \ref{Section special rep} below. The classification of the cuspidal representations of the full group of automorphisms of a thick regular tree was later obtained by  Ol'shanskii. We present this classification  in Section \ref{cuspidal representations of the group of automorphism of a tree} below. 
		
	\subsection{Spherical representations}\label{Section spherical rep}
	Let $T$ be a thick semi-regular tree. The purpose of this section is to recall the classification of spherical representations of any closed subgroup  $G\leq \Aut(T)$ acting $2$-transitively on the boundary $\partial T$. We recall that a group $G$ acts $2$-transitively on a set $X$ if the action of $G$ on $X$ is transitive and if the restriction of the action of $\Fix_G(x)=\{g\in G : gx=x\}$ is transitive on $X-\{x\}$ for all $x\in X$. Various expositions of the classification of the spherical representations of such groups can be found in the literature, see \cite{Cartier1973}, \cite{Matsumoto1977}, \cite{Olshanskii1977}, \cite{FigaNebbia1991} and \cite{Choucroun1994} for instance. The exposition presented here is based on the approach adopted in \cite{FigaNebbia1991}. However, our presentation highlights slightly different intermediate results for the convenience of the writing of Chapter \ref{Chapter Fell topology of Aut(T)} where we describe the restriction of the Fell topology to spherical and special representations. The classification of spherical representations is gathered in Theorems \ref{thm la classification des spherique cas trans} and \ref{thm la classification des spherique cas 2 orbites} below. 
	
	We start with the following characterisation of closed subgroups $G\leq \Aut(T)$ acting $2$-transitively on the boundary. 
	\begin{theorem}[{\cite[Lemma 3.1.1]{BurgerMozes2000}}]\label{Theorem Burger mozes acting on the boundary trans}
		Let $T$ be a locally finite tree and let $G\leq \Aut(T)$ be a closed subgroup. Then, the following are equivalent:
		  \begin{enumerate}[label=(\roman*)]
		  	\item $G$ acts $2$-transitively on the boundary $\partial T$.
		  	\item $G$ is not compact and acts transitively on the boundary $\partial T$.
		  	\item $\Fix_G(v)$ is transitive on the boundary $\partial T$ for all $v\in V$.  
		  \end{enumerate}
	 If any of these properties is fulfilled, the tree $T$ is semi-regular, $\Fix_G(v)$ is $2$-transitive on the sphere of radius one around $v$: $$S_1(v)=\{w\in V: d_T(v,w)=1\}$$ and $G$ is either vertex-transitive or has two orbits of vertices (those of any pair of adjacent vertices). 
	\end{theorem}
	\noindent The classification of the spherical representations of these groups relies on the notion of Gelfand pair. We now provide reminders on this notion. 
	\begin{definition}
		Let $G$ be a locally compact group and $K\leq G$ be a compact subgroup $K\leq G$. We say that $(G, K)$ is a \textbf{Gelfand pair} if the convolution algebra $C_c(K \backslash  G \slash K)$ of continuous compactly supported $K$-bi-invariant functions on $G$ is commutative.
	\end{definition}
	\begin{example}
		Let $G$ be an abelian locally compact group. Then, $(G,\{1_G\})$ is a Gelfand pair as $C_c(K \backslash  G \slash K)$ is the convolution algebra $C_c(G)$.
	\end{example}
 	The following result provides a geometric criterion to prove that $(G,K)$ is a Gelfand pair even if the convolution might be complicated to compute. 
 	\begin{lemma}\label{lemma proof of symmetric gelfand pair}
 		Let $G$ be a locally compact group, $K\leq G$ be a compact subgroup and suppose that $g^{-1}\in KgK$ for all $g\in G$. Then, $(G,K)$ is a Gelfand pair. 
 	\end{lemma}
 	\begin{proof}
 		Let $\varphi,\psi\in C_c(K \backslash  G \slash K)$ and let $\mu$ be a Haar measure on $G$. We recall that the convolution of $\varphi$ and $\psi$ is defined by 
 		\begin{equation*}
 		\varphi * \psi (g)=\int_G \varphi (gh) \psi (h^{-1}) \qq \diff \mu(h)= \int_G \varphi(h) \psi(h^{-1}g) \diff \mu (h)\q \forall g\in G.
 		\end{equation*}
 		In particular, since both $\varphi$ and $\psi$ are $K$-bi-invariant, we deduce from our hypothesis that 
 		\begin{equation*}
 		\varphi * \psi (g) = \int_G \varphi (h) \psi (h^{-1}g) \qq \diff \mu(h) =  \int_G \varphi (h^{-1}) \psi (g^{-1}h) \qq \diff \mu(h)=\psi *\varphi(g^{-1}).
 		\end{equation*}
 		As $\psi*\varphi$ is $K$-bi-invariant, it follows that $\varphi*\psi=\psi*\varphi$. 
 	\end{proof}
 	\begin{remark}\label{Remark group of automorphism that are 2 trans are unimodular}
 		The hypothesis of unimodularity is a necessary condition since \cite[Proposition 6.1.2]{vanDijk2009} ensures that a locally compact group $G$ admitting a Gelfand pair $(G,K)$ is unimodular. 
 	\end{remark} 
 	\begin{lemma}\label{Lemma aut T fixator of a vertex is a Gelfand pair}
 		Let $T$ be a semi-regular tree and let $G\leq \Aut(T)$ be a closed subgroup of $G$ acting $2$-transitively on the boundary. Then, for every vertex $v\in V$, $(G,\Fix_G(v))$ is a Gelfand pair.
 	\end{lemma}
 	\begin{proof}
 		We first prove that $G$ is unimodular. Consider the subgroup $H$ of $G$ generated by any pair of stabilisers $\Fix_{G}(v)$, $\Fix_{G}(v')$ of adjacent vertices $v,v'\in V$. This subgroup $H$ is generated by compact groups. In particular, the restriction of the modular function $\Delta_G$ to $H$ is trivial just as in Example \ref{example compaclty generated is unimodular}. Note moreover that $H$ has index at most two in $G$. Since the modular function $\Delta_G$ is a continuous group homomorphism and as $\R^*_+$ is torsion free, it follows that $G$ is unimodular. 
 		
 		We now prove that for all $g\in G$, we have that $g^{-1}\in \Fix_G(v)g\Fix_G(v)$. The result will then follow from Lemma \ref{lemma proof of symmetric gelfand pair}. Since $G$ acts $2$-transitively on the boundary notice that $\Fix_G(v)$ is transitive on each of the sphere centred at $v$.  On the other hand, since $G$ acts by isometry on the tree, $g^{-1}v$ and $gv$ are equidistant from $v$ because $d_T(gv,v)=d_T(v,g^{-1}v)$. This proves the existence of an element $k\in \Fix_G(v)$ such that $kg^{-1}v=gv$.  Equivalently we have that $g^{-1}kg^{-1}\in \Fix_G(v)$ and it follows as desired that  $g^{-1}\in \Fix_G(v)g\Fix_G(v)$.
 	\end{proof}
 	Now, let $L^1(K \backslash  G \slash K)$ be the completion of $C_c(K \backslash  G \slash K)$ with respect to the $L^1$ norm. Just as $L^1(G)$, this convolution Banach algebra usually does not satisfy the $C^*$-identity but admits a $C^*$-completion with respect to a different norm. Let $L^1(K \backslash  G \slash K)^{\sim}$ be the set of equivalence classes of non-degenerate $*$-representations of $L^1(K \backslash  G \slash K)$. Each $\pi \in L^1(K \backslash  G \slash K)^{\sim}$ provides a semi-norm $\norm{\cdot}{\tilde{\pi}}$ satisfying the $C^*$-identity on $L^1(K \backslash  G \slash K)$ and defined by $$\norm{\psi}{\tilde{\pi}}=\norm{\tilde{\pi}(\psi)}{op(\Hr{\tilde{\pi}})}\q \forall \psi\in L^1(K \backslash  G \slash K).$$
 	Furthermore, the semi-norm $\norm{\cdot}{C^*_K}$ defined by
 	$$\norm{\varphi}{C^*_K}= \sup_{\tilde{\pi}\in L^1(K \backslash  G \slash K)^{\sim}}\norm{\tilde{\pi}(\varphi)}{op(\Hr{\tilde{\pi}})} \q \forall \varphi \in L^1(K \backslash  G \slash K)$$
 	is a norm satisfying the $C^*$-identity. We denote by $C^*(K \backslash  G \slash K)$ the completion of $L^1(K \backslash  G \slash K)$ with respect to $\norm{\cdot}{C^*_K}$ and notice that it defines a $C^*$-algebra containing $C_c(K \backslash  G \slash K)$. Furthermore, by construction each non-degenerate $*$-representation of $L^1(K \backslash  G \slash K)$ extends to a non-degenerate $*$-representation of $C^*(K \backslash  G \slash K)$ and is uniquely determined by the value it takes on $C_c(K \backslash  G \slash K)$. On the other hand, it is important to notice that $C^*(K \backslash  G \slash K)$ might be different from the completion of $C_c(K \backslash  G \slash K)$ inside $C^*(G)$. The interest of Gelfand pairs lays inside the fact that $C^*(K \backslash  G \slash K)$ is commutative. We recall from Section \ref{section Fell topology}  that every unitary representation $\pi$ of $G$ corresponds to a non-degenerate $*$-representation $\tilde{\pi}$ of $L^1(G)$ defined by 
 	$$\tilde{\pi}(\varphi)= \int_G \varphi(g)\pi(g)\diff \mu(g)\q \forall \varphi \in L^1(G).$$
 	Now, notice that the space $$\Hr{\pi}^K=\{\xi\in \Hr{\pi}: \pi(k)\xi= \xi\qq \forall k\in K\}$$ is $\tilde{\pi}(C_c(K \backslash  G \slash K))$-invariant and that $\mathds{1}_K$ acts trivially on $\Hr{\pi}^K$. In particular, when $\pi$ admits a non-zero $K$-invariant vector, $\tilde{\pi}$ gives rise to a non-degenerate $*$-representation of $L^1(K \backslash  G \slash K)$ with representation space $\Hr{\pi}^K$. Furthermore, this representation extends to a $*$-representation of $C^*(K \backslash  G \slash K)$ and it can be shown this representation is irreducible if and only if $\pi$ is irreducible {\cite[Chapter IV. \textsection{2}, Theorem $2$]{Serge1985}}. On the other hand, since $C^*(K\backslash G \slash K)$ is commutative when $(G,K)$ is a Gelfand pair, every of its irreducible representation is one-dimensional. The following result follows.
 	\begin{theorem}[{\cite[Chapter IV. \textsection{2}, Theorem $2$]{Serge1985}}]\label{Theorem dimension plus petite que 1}
 		Let $(G,K)$ be a Gelfand pair and let $\pi$ be an irreducible representation of $G$. Then, we have $\dim(\Hr{\pi}^K)\leq 1$.
 	\end{theorem}
 	In addition, the classification of the irreducible representations of $G$ admitting non-zero $K$-invariant vectors is equivalent to the classification of a particular family of functions on $G$ that are tightly related to the characters of $C^*(K \backslash  G \slash K)$. This is the purpose of the following definition.
 	\begin{definition}\label{defintion K spherical functions}
 		A function  $\varphi: G\rightarrow \C$ is $K$\textbf{-spherical} if it satisfies the following properties:
 		\begin{enumerate}[label=(\roman*)]
 			\item $\varphi$ is $K$-bi-invariant and continuous.
 			\item $\varphi(1_G)=1$.
 			\item $\varphi$ is an eigenfunction for the convolution on the right by elements of $C_c(K\backslash G \slash K)$. In other words, for every $\psi\in C_c(K\backslash G \slash K)$ we have that
 			$$\varphi*\psi= \lambda(\varphi,\psi) \varphi$$
 			for some complex number $\lambda(\varphi,\psi)\in \C$. 
 		\end{enumerate}
 	\end{definition}
 	The following result provides a lot of equivalent definitions and strong correspondence with the non-degenerate $*$-representations of $C^*(K \backslash  G \slash K)$.
 	\begin{theorem}[{\cite[Chapter IV. \textsection{3}, Theorem $5$, $6$ and $7$]{Serge1985}}]\label{Theorem equivalent defintoion of K spherical funcitons} Let $G$ be a unimodular locally compact group, $\varphi: G \rightarrow \C$ be a not identically $0$ continuous function and $\mu$ be a Haar measure on $G$ and $K\leq G$ be a compact subgroup such that $\mu(K)=1$. Then, the following are equivalent: 
 		\begin{enumerate}[label=(\roman*)]
 			\item $\varphi$ is a $K$-spherical function.
 			\item For all $g,h\in G$, we have that 
 			$$\int_K\varphi(gkh)\diff \mu(k)=\varphi(g)\varphi(h).$$
 			\item $\varphi\in C(K\backslash G \slash K)$ and the operator 
 			$$L_\varphi : C_c(K\backslash G \slash K)\rightarrow \C: \psi \mapsto \int_G \psi(g)\varphi(g)\diff \mu(g)$$
 			is an algebra homomorphism of $C_c(K\backslash G \slash K)$ onto $\C$. 
 		\end{enumerate}
 		Furthermore, if $\varphi$ is a bounded $K$-spherical function, $L_\varphi$ extends to a continuous algebra homomorphism of $L^1(K\backslash G \slash K)$ into $\C$ and every continuous algebra homomorphism of $L^1(K\backslash G \slash K)$ into $\C$ is of this form for some bounded $K$-spherical function.
 	\end{theorem} 
	Relying on this characterisation, we obtain a correspondence between the irreducible representations of $G$ with non-zero $K$-invariant vectors and the $K$-spherical functions of positive type. 
 	\begin{lemma}\label{Lemma la fonction de type positive est K spherique cas general}
 		Let $(G,K)$ be a Gelfand pair, $\pi$ be an irreducible representation of $G$ and $\xi \in \Hr{\pi}^K$ be a unit vector. Then, the function 
 		$$\varphi_{\xi,\xi}:G\rightarrow \C : g\mapsto \prods{\pi(g)\xi}{\xi}_{\Hr{\pi}}$$
 		is a $K$-spherical function of positive type.
 	\end{lemma} 
 	\begin{proof}
 		Since $\xi$ is a unit $K$-invariant vector, the function $\varphi_{\xi,\xi}$ clearly defines a continuous $K$-bi-invariant functions of positive type satisfying $\varphi_{\xi,\xi}(1_G)=1$. On the other hand, since $(G,K)$ is a Gelfand pair, Theorem \ref{Theorem dimension plus petite que 1} ensures that $\dim(\Hr{\pi}^K)=1$. It follows that $\Hr{\pi}^K$ is spanned by $\xi$. Now, denote by $\tilde{\pi}$ the $*$-representation of $L^1(G)$ corresponding to $\pi$. Notice for each $\psi\in C_c(K\backslash G \slash K)$ that the vector $\tilde{\pi}(\psi)\xi$ belongs to $\Hr{\pi}^K$. In particular, one has that  $\tilde{\pi}(\psi)\xi=\lambda(\psi)\xi$ for some complex number $\lambda(\psi)\in \C$  and the map $$\lambda:  C_c(K\backslash G \slash K)\rightarrow \C : \psi\mapsto \lambda(\psi)$$ 
 		is a $*$-homomorphism. A  straightforward computation shows that 
 		\begin{equation*}
 		\begin{split}
 		L_{\varphi_{\xi,\xi}}(\psi)&=\int_G \psi(g)\varphi_{\xi,\xi}(g)\diff \mu(g)\\
 		&= \int_G \prods{\psi(g)\pi(g)\xi}{\xi}_{\Hr{\pi}}\diff \mu(g)\\
 		&=\prods{\tilde{\pi}(\psi)\xi}{\xi}_{\Hr{\pi}}=\lambda(\psi).
 		\end{split}
 		\end{equation*}
 		It follows from Theorem \ref{Theorem equivalent defintoion of K spherical funcitons} that $\varphi_{\xi,\xi}$ is $K$-spherical.
 	\end{proof}
 	On the other hand, Theorem \ref{Theorem dimension plus petite que 1} ensures that the space of $K$-invariant vectors of an irreducible representation $\pi$ of $G$ admitting non-zero $K$-invariant vectors is one dimensional. In particular, in that case, the function $\varphi_{\xi,\xi}$ defined above does not depend on the choice of the unit vector $\xi \in \Hr{\pi}^K$ so that there is a unique $K$-spherical function of positive type on $G$ that is associated to $\pi$ by Lemma \ref{Lemma la fonction de type positive est K spherique cas general}. We denote by $\varphi_\pi$ this unique $K$-spherical function of positive type. The following result ensures that this correspondence is bijective. 
  	\begin{theorem}[{\cite[Chapter IV. \textsection{3}, Theorems $3$ and $9$]{Serge1985}}]\label{Thoerme la correspondence bijective entre spherical and spherical}
 		There is a bijective correspondence $\pi \rightarrow \varphi_\pi$ with inverse map given by the GNS construction between the equivalence classes of irreducible representations of $G$ with non-zero $K$-invariant vectors and the $K$-spherical functions of positive type on $G$.
 	\end{theorem}
 	In particular, the classification of irreducible representations of $G$ admitting non-zero $K$-invariant vectors reduces to the classification of $K$-spherical functions of positive type on $G$. Now, let $T$ be a thick semi-regular tree and let $G\leq \Aut(T)$ be a closed subgroup that is $2$-transitive on the boundary. Theorem \ref{Thoerme la correspondence bijective entre spherical and spherical} ensures the existence of a bijective correspondence between the equivalence classes of spherical representations of $G$ with non-zero $\Fix_G(v)$-invariant vectors and the $\Fix_G(v)$-spherical functions of positive type. We are left with the task of determining those. We start with a result focusing on the decomposition of $G$ into double $\Fix_G(v)$-cosets. 
 	\begin{lemma}\label{lemma decomposition of $G$ into double cosets}
 		Let $T$ be a thick semi-regular tree, $v$ be a vertex of $T$, $G\leq \Aut(T)$ be a closed subgroup acting $2$-transitively on the boundary $\partial T$ and let $\tau$ be a translation of minimal step in $G$ along a geodesic containing $v$. Then, the following decomposition holds
 		\begin{equation}\label{equation decomposition de G en K double coset avec translation n}
 		G = \bigsqcup_{n\in \N}\Fix_G(v)\tau^n\Fix_G(v).
 		\end{equation}
 	\end{lemma}
 	\begin{proof}
 		Denote by $t$ the step of $\tau$. Since $\tau$ is a translation of minimal step in $G$, we have that $t=1$ if $G$ is vertex-transitive and that $t=2$ if $G$ has two orbits of vertices. Now, notice that the elements of $\Fix_G(v)\tau^n\Fix_G(v)$ map $v$ to vertices at distance $nt$ from $v$ so that each of those double cosets are disjoint from one another. Consider any element $g\in G$. Since $\tau$ is a translation of minimal step $t$ in $G$, there exists a unique $n\in \N$ such that $nt=d(v,gv)$. In particular, $gv$ and $\tau^n v$ both belongs to the sphere $$S_{nt}(v)=\{w\in Gv: d_T(v,w)=nt\}.$$ On the other hand, since $G$ acts $2$-transitively on the boundary $\partial T$, Theorem \ref{Theorem Burger mozes acting on the boundary trans} ensures that $\Fix_G(v)$ is transitive on each such sphere. In particular, there exists an element $k\in \Fix_G(v)$ such that $gv= k\tau^nv$ which implies that $g^{-1}k\tau^n\in \Fix_G(v)$ and hence that $g\in \Fix_G(v)\tau^n\Fix_G(v)$. 
 	\end{proof}
 	As a direct consequence, every $K$-bi-invariant function $\varphi: G\rightarrow \C$ is entirely determined by the values $\varphi(\tau^n)$ it takes on $\tau^n$ with $n\in \N$. The following result shows further that every $K$-spherical function of $G$ is entirely determined by the value it takes on $\tau$. 
 	\begin{lemma}\label{lemma calcul des fonctions psheriques}
 		Let $T$, $G$, $v$ and $\tau$ be as in Lemma \ref{lemma decomposition of $G$ into double cosets} and let $\varphi$ be a $\Fix_G(v)$-spherical function of $G$. For each integer $n\geq2$ the following hold:
 		\begin{enumerate}[label=(\roman*)]
 			\item If $G$ is vertex-transitive: 
 			\begin{equation}\label{equation formule spherique cas transitive}
 			\varphi(\tau^n)= \frac{d}{d-1}\varphi(\tau)\varphi(\tau^{n-1}) - \frac{1}{d-1}\varphi(\tau^{n-2})
 			\end{equation}
 			where $d$ is the degree of $v$. 
 			\item If $G$ has two orbits of vertices: 
 			\begin{equation}\label{equation formule spherique cas 2 orbites}
 			\varphi(\tau^n)= \bigg(\frac{d}{d-1}\varphi(\tau)-\frac{d'-2}{(d'-1)(d-1)}\bigg)\varphi(\tau^{n-1}) - \frac{1}{(d'-1)(d-1)}\varphi(\tau^{n-2})
 			\end{equation}
 			where $d$ is the degree of $v$ and $d'$ is the degree of its adjacent vertex. 
 		\end{enumerate} 
 	\end{lemma}
 	\begin{proof}
 		Let $n\geq 2$ be a positive integer and renormalise the Haar measure $\mu$ of $G$ so that $\mu(\Fix_G(v))=1$. According to the second characterisation of $\Fix_G(v)$-spherical functions provided by Theorem \ref{Theorem equivalent defintoion of K spherical funcitons} we have
 		\begin{equation*}\label{equation pour determiner les psherique avec la charact de l'integrale}
 		\begin{split}
 		\varphi(\tau^{n-1})\varphi(\tau)&= \int_{\Fix_G(v)}\varphi(\tau^{n-1}k\tau)\diff \mu (k).
 		\end{split}
 		\end{equation*}
 		If $G$ is vertex-transitive, let $U=\{k\in \Fix_G(v): k\tau v=\tau^{-1} v\}$ and notice from the $K$-bi-invariance of $\varphi$ that 
 		\begin{equation*}
 		\begin{split}
 		\varphi(\tau^{n-1})\varphi(\tau)&= \int_{\Fix_G(v)-U}\varphi(\tau^{n-1}k\tau)\diff \mu (k)+ \int_{U}\varphi(\tau^{n-1}k\tau)\diff \mu (k)\\
 		&= \int_{\Fix_G(v)-U}\varphi(\tau^{n})\diff \mu (k)+ \int_{U}\varphi(\tau^{n-2}k)\diff \mu (k)\\
 		&=\frac{d-1}{d}\varphi(\tau^n)+ \frac{1}{d}\varphi(\tau^{n-2}).
 		\end{split}
 		\end{equation*}
 		This proves as desired that 
 		\begin{equation*}
 		\varphi(\tau^n)= \frac{d}{d-1}\varphi(\tau)\varphi(\tau^{n-1}) - \frac{1}{d-1}\varphi(\tau^{n-2}).
 		\end{equation*}
 		If on the other hand $G$ has two orbits of vertices, we set $$U=\{k\in \Fix_G(v):  k\tau v=\tau^{-1}v\}$$ and $$U'= \{k\in \Fix_G(v): d(k\tau v,v)\geq d(k\tau v,\tau^{-1} v)\},$$ and we obtain similarly that
 		\begin{equation*}
 		\begin{split}
 		\varphi(\tau^{n-1})\varphi(\tau)&=\int_{\Fix_G(v)} \varphi(\tau^{n-1}k\tau)\qq \diff\mu (k)\\
 		&=\int_{\Fix_G(v)-U'} \varphi(\tau^{n-1}k\tau)\qq \diff\mu (k)+\int_{U'-U} \varphi(\tau^{n-1}k\tau)\qq \diff\mu (k)\\
 		&\q \q \q \q \q \q \q \q \q \q \q \q \q \q \q \q +\int_{U} \varphi(\tau^{n-1}k\tau)\qq \diff\mu (k) \\
 		&=\frac{d-1}{d}\varphi(\tau^n)+\frac{d'-2}{d(d'-1)}\varphi(\tau^{n-1})+ \frac{1}{d(d'-1)}\varphi(\tau^{n-2})
 		\end{split}
 		\end{equation*}
 		This proves as desired that 
 		\begin{equation*}
 		\varphi(\tau^n)= \bigg(\frac{d}{d-1}\varphi(\tau)-\frac{d'-2}{(d'-1)(d-1)}\bigg)\varphi(\tau^{n-1}) - \frac{1}{(d'-1)(d-1)}\varphi(\tau^{n-2})
 		\end{equation*}
 	\end{proof}
 	On the one hand, it can be shown that the $\Fix_G(v)$-bi-invariant functions $\varphi: G\rightarrow \C$ satisfying \eqref{equation formule spherique cas transitive} if $G$ is vertex-transitive and \eqref{equation formule spherique cas 2 orbites} if $G$ has two orbits of vertices are all $\Fix_G(v)$-spherical functions. In particular, the above result provides a bijective correspondence between the $\Fix_G(v)$-spherical functions $\varphi$ of $G$ and their value $\varphi(\tau)$ at $\tau$. On the other hand, an argumentation at the level of the algebras shows $\varphi(\tau)$ and hence $\varphi$ is real for every $\Fix_G(v)$-spherical functions of positive type. To see this, notice in light of Lemma \ref{lemma decomposition of $G$ into double cosets} that the convolution algebra $C_c(\Fix_G(v)\backslash G \slash \Fix_G(v))$ is generated (as an algebra) by $\mathds{1}_{\Fix_G(v)\tau \Fix_G(v)}$. The $\Fix_G(v)$-spherical functions of $G$ are therefore the eigenfunctions for the convolution on the right by $\mathds{1}_{\Fix_G(v)\tau \Fix_G(v)}$. 
 	 \begin{remark}
 		If $G$ is vertex-transitive, this approach provides a deep connection between the representation theory of $G$ and the spectral theory of $T$. Indeed, when $G$ is vertex transitive there a bijective correspondence between the $\Fix_G(v)$-bi-invariant functions of $G$ and the set of functions on $V$ that are radial around $v$ (constant on the spheres centred at $v$). Furthermore, under this correspondence the operator of convolution on the right by $\mathds{1}_{\Fix_G(v)\tau \Fix_G(v)}$ can be interpreted as the adjacency operator of the tree. Indeed, if the Haar-measure $\mu$ of $G$ is renormalised so that $\mu(\Fix_G(v))=1$, we have for all $w\in V$ that 
 		$$\varphi*\mathds{1}_{\Fix_G(v)\tau \Fix_G(v)}(w)=\sum_{u\in V: d_T(w,u)=1}\varphi(u).$$
 		In particular, the $\Fix_G(v)$-spherical functions of $G$ are the radial eigenfunction of the adjacency operator of the tree.
 	\end{remark}

 	Since $G$ is unimodular, we recall from \eqref{equation the involution in C^*(G)} that the involution on $C^*(K\backslash G\slash K)$ is given by $\psi^*(g)=\overline{\psi(g^{-1})}$. 
 	In particular, $\mathds{1}_{\Fix_G(v)\tau \Fix_G(v)}$ is a self-adjoint element of $C^*(K\backslash G\slash K)$ because the proof of Lemma \ref{Lemma aut T fixator of a vertex is a Gelfand pair} ensures that $$\Fix_G(v)\tau \Fix_G(v)=\Fix_G(v)\tau^{-1} \Fix_G(v).$$ This implies that the spectrum of $\mathds{1}_{\Fix_G(v)\tau \Fix_G(v)}$ is real. On the other hand, since $C^*(\Fix_G(v)\backslash G\slash\Fix_G(v))$ is commutative, this spectrum coincides with the set of $L(\mathds{1}_{\Fix_G(v)\tau \Fix_G(v)})$ where $L$ runs over all characters of $C^*(\Fix_G(v)\backslash G\slash\Fix_G(v))$. Now, the third characterisation of spherical functions provided by Theorem \ref{Theorem equivalent defintoion of K spherical funcitons} ensures that each bounded $\Fix_G(v)$-spherical function $\varphi$ on $G$ gives rise to such character $L_\varphi$. Furthermore, a  straightforward computation shows that
 	\begin{equation*}
 	\begin{split}
 	L_\varphi\big(\mathds{1}_{\Fix_G(v)\tau \Fix_G(v)}\big)=\int_{\Fix_G(v)\tau \Fix_G(v)} \varphi (k)\diff \mu(k)  = \varphi(\tau).
 	\end{split}
 	\end{equation*}
 	Since every function of positive type on $G$ is bounded, this proves that $\varphi(\tau)$ is a real number for every $\Fix_G(v)$-spherical function $\varphi$ of positive type on $G$. To obtain a classification of the spherical representations of $G$ admitting non-zero $\Fix_G(v)$-invariant vectors we are left to find the real numbers $\alpha$ such that the $\Fix_G(v)$-spherical function $\varphi: G\rightarrow \C$ with $\varphi(\tau)=\alpha$ is of positive type. This part of the proof is quite technical and we do not present a full proof in these notes. This procedure leads to the following classifications. For clarity of the exposition, we present separately the classification in the vertex transitive case and in the case where $G$ has two orbits of vertices. 
	\begin{theorem}[{\cite[Chapter II]{FigaNebbia1991}}]\label{thm la classification des spherique cas trans}
		Let $T$ be a thick regular tree, $v\in V$ and let $G\leq \Aut(T)$ be a closed vertex-transitive subgroup acting $2$-transitively the boundary $\partial T$. Then, every spherical representation of $G$ admits a non-zero $\Fix_G(v)$-invariant vector. Furthermore, the equivalence classes of spherical representations of $G$ are in bijective correspondence with the interval $\mathcal{I}_v=\lb -1, 1\rb$ via the map $\phi_{v,G}:\pi\mapsto \varphi_\pi(\tau)$ where $\tau$ is any element of $G$ such that $d(\tau v,v)=1$ and $\varphi_\pi$ is the unique $\Fix_G(v)$-spherical function of positive type attached with $\pi$. Under this correspondence, $\phi_{v,G}^{-1}(1)$ is the trivial representation $1_{\widehat{G}}$.
	\end{theorem}
	The following theorem comes from \cite{Matsumoto1977} but is formulated differently for coherence of our exposition. 
	\begin{theorem}[{\cite[Chapter 5, Section 6]{Matsumoto1977}}]\label{thm la classification des spherique cas 2 orbites}
		Let $T$ be a thick semi-regular tree, $v,v'\in V$ be adjacent vertices and let $G\leq \Aut(T)^+$ be a closed subgroup of type-preserving automorphisms acting $2$-transitively on the boundary $\partial T$. Then, there is exactly one spherical representation $\pi_v$ of $G$ with a non-zero $\Fix_G(v)$-invariant vector but no non-zero $\Fix_G(v')$-invariant vector. Apart from the two exceptional representations $\pi_v$ and $\pi_{v'}$, every spherical representation of $G$ admits a non-zero $\Fix_G(w)$-invariant vector for all $w\in V$. In addition, the equivalence classes of spherical representations admitting a non-zero $\Fix_G(v)$-invariant vector are in bijective correspondence with the interval  $\mathcal{I}_v=\big\lb -\frac{1}{d'-1},1\big\rb$ via the map $\phi_{v,G}:\pi \mapsto \varphi_\pi(\tau)$ where $\tau$ is any element of $G$ such that $d(\tau v,v)=2$ and $d'$ is the degree of $v'$. Under this correspondence, $\phi_{v,G}^{-1}(-\frac{1}{d'-1})$ is the exceptional spherical representation $\pi_v$ and $\phi_{v,G}^{-1}(1)$ is the trivial representation $1_{\widehat{G}}$. Finally, if $\pi$ is a non-exceptional spherical representation of $G$ we have that 
		$$\phi_{v',G}(\pi)=\frac{d(d'-1)}{d'(d-1)}\phi_{v,G}(\pi) + \frac{d-d'}{d'(d-1)}.$$ 
	\end{theorem}
	\begin{remark}
		Certain extremal values of the real interval corresponding to spherical functions of positive types provided in the literature are not correct when $G$ has two orbits of vertices. The incorrect values provided by \cite{Bouaziz-Kellil1988} and \cite{Amann2003} for instance are based on earlier inaccuracies \cite{IozziPicardello1983} that were already pointed out in \cite[Remark $2$, pg. $243$]{CartwrightMlotowskiSteger1994}. In particular, while each of the bounded $\Fix_G(v)$-spherical functions are of positive type when $G$ is vertex transitive, there exists a bounded $\Fix_G(v)$-spherical functions that are not of positive type when $G$ has two orbits of vertices.
	\end{remark}
	We end this section by providing an explicit description of the spherical representations of $G$ that belong to the so called \tg{principal series} of $G$. In particular, this will provide a proof that the corresponding real numbers are indeed associated to spherical functions of positive type. The approach adopted here is to describe them as the induced representations form the stabiliser of an end of the tree. An equivalent approach relying on measures of the boundary of the tree is described in \cite[Chapter II]{Choucroun1994} to which we refer to for an alternative expository.
	
	Let $T$ be a thick semi-regular tree and $G$ be a closed subgroup $G\leq \Aut(T)$ that is $2$-transitive on the boundary. Consider an end of the tree $\omega \in \partial T$ and let $G_\omega=\{g\in G : g\omega = \omega\}$ be the stabiliser of $\omega$ in $G$. We recall for every vertex $w\in V$, the existence of a unique geodesic $\lb w,\omega\lb$ from $w$ to $\omega$ and that for all $w,w'\in V$, there is a unique $u\in V$ such that $\lb w,\omega\lb\cap \lb w',\omega\lb=\lb u,\omega\lb$. We define the $\omega$-\tg{horocycles distance} $\delta(w,w',\omega)$ between $w$ and $w'$ as 
	$$\delta(w, w', \omega) = d_T(w, u) - d_T(w', u).$$
	We start our construction of the representations of the principal series of $G$ with the following result.
	
	\begin{lemma}\label{Lemma les charracter de Ggamme}
		Let $T$ be a thick semi-regular tree, $\omega \in \partial T$, $v\in V$, $s\in \mathbb{C}$ be a complex number with $\modu{s}=1$ and $G\leq \Aut(T)$ be a closed subgroup acting $2$-transitively on the boundary. Then, the map $\chiup_{s}:G_\omega\rightarrow \mathbb{C}:  g\mapsto s^{\delta(v,gv,\omega)}$
		does not depend on $v\in V$ and defines a unitary character of $G_\omega$.
	\end{lemma} 
	The proof of this result relies on the following technical lemma.
	\begin{lemma}\label{Lemma la distance horocyclique va faire un morphisme}
		Let $T$ be a thick semi-regular tree, $\omega \in \partial T$, $v\in V$ and $G\leq \Aut(T)$ be a closed subgroup acting $2$-transitively on the boundary. Then, for all $g,h\in G_\omega$ one has that 
		\begin{equation*}\label{equation la formule qui montre qu'on va avoir un morphisme}
		\delta(v,hgh^{-1}v,\omega)=\delta (v,gv,\omega).
		\end{equation*}
	\end{lemma}
	\begin{proof}
		Let $\tau\in G$ be a translation of minimal step in $G$ along a geodesic containing $v$ and with $\omega$ as attracting point. The group $G_\omega$ is generated by $\Fix_{G_\omega}(v)$ and $\tau$. In particular, we just need to prove that the relation is satisfied for $g,h\in \Fix_{G_\omega}(v)\cup\{\tau,\tau^{-1}\}$. We treat these cases separately:
		\begin{itemize}[leftmargin=*]
			\item  If $g\in \Fix_{G_\omega}(v)$ and $h\in\Fix_{G_\omega}(v)\cup\{\tau^{-1}\}$, one has that $hgh^{-1}v=v$ and the equality follows trivially since $$\delta(v,hgh^{-1}v,\omega)=\delta (v,v,\omega)=\delta (v,gv,\omega).$$
			\item If $g\in  \Fix_{G_\omega}(v)$ and $h=\tau$, one has that $$\delta(v,hgh^{-1}v,\omega)=0=\delta(v,gv,\omega).$$ 
			\item If $g,h\in \{\tau, \tau^{-1}\}$ the result is trivial.
			\item If $g=\tau$ and $h\in  \Fix_{G_\omega}(v)$, one has that $hgh^{-1}v=h\tau v=\tau v=gv$ and the equality follows trivially.
			\item If $g=\tau^{-1}$ and $h\in  \Fix_{G_\omega}(v)$ one has that 
			\begin{equation*}
			\begin{split}
			\delta(v,hgh^{-1}v,\omega)&=\delta(v,h\tau^{-1}v,\omega)\\
			&=d_T(h\tau^{-1} v,v)-d_T(v,v)\\
			&=d_T(v,v)-d_T(\tau v,v)=-1=\delta(v,\tau^{-1} v,\omega).
			\end{split}
			\end{equation*}
		\end{itemize}
	\end{proof}
	\begin{proof}[Proof of Lemma \ref{Lemma les charracter de Ggamme}]
		To prove that $\chiup_s$ does not depend on $v\in V$, let $w\in V$, let $g\in G_\omega$ and let us show that $\delta(v,gv, \omega)=\delta(w, gw, \omega)$. This is easily verified as 
		\begin{equation*}
		\begin{split}
		\delta(v,gv,\omega)&=\delta(v,w,\omega)+ \delta(w,gw,\omega)+ \delta(gw,gv,\omega)\\
		&=\delta(v,w,\omega)+ \delta(w,gw,\omega)+ \delta(w,v,\omega)\\
		&=\delta(w,gw,\omega).
		\end{split}
		\end{equation*}
		On the other hand, $\chiup_s$ defines character, notice, in light of Lemma \ref{Lemma la distance horocyclique va faire un morphisme} that for all $g,h\in G_\omega$ one has that 
		\begin{equation*}
		\begin{split}
		\delta(v,ghv,\omega)&=\delta(v,hv,\omega)+\delta(hv,ghv,\omega)\\
		&=\delta(v,hv,\omega)+\delta(v,h^{-1}ghv,h^{-1}\omega)\\
		&=\delta(v,hv,\omega)+\delta(v,h^{-1}ghv,\omega)\\
		&=\delta(v,hv,\omega)+\delta(v,gv,\omega).
		\end{split}
		\end{equation*}
		It follows as desired that $\chiup_s(gh)=\chiup_s(g)\chiup_s(h)$.
	\end{proof}
	 We are going to show that the induced representation 
	$$\pi_s=  \Ind_{G_\omega}^G(\chiup_s)$$
	is a spherical representation of $G$ for every complex number $s\in \C$ with $\modu{s}=1$. 
	Since $G$ is transitive on the boundary for each $\zeta\in \partial T$, there exists an element $g_{\zeta}\in G$ such that $g_{\zeta}\omega=\zeta$. In particular, notice that the map $G/G_\omega$ is in bijective correspondence with $\partial T$ via the map 
	$$\partial T \rightarrow G/ G_\omega : \zeta\mapsto g_{\zeta}G_\omega$$ provides a bijective correspondence between $\partial T$ and $G/ G_\omega$.  By abuse of notation, if a function $\varphi: G\rightarrow \C$ is $G_\omega$-right invariant we will denote by $\varphi(\zeta)$ the value of $\varphi$ on any element of $g_{\zeta}G_\omega$. Since $G_\omega$ is a closed subgroup of $G$ notice that $\partial T\simeq G/G_\omega$ carries a unique quasi-invariant measure $\nu$ such that $\nu(\partial T)=1$. We consider the Hilbert space $\mathcal{H}_s$ of measurable functions $\psi: G\rightarrow \C$ satisfying:
	\begin{enumerate}
		\item $\psi(gh)=\chiup_s(h^{-1})\psi(g) \q \forall h\in G_\omega, \qq \forall g\in G,$ 
		\item ${\int_{\partial T} \modu{\psi(\zeta)}^2\qq \diff \nu(\zeta)}<+\infty,$ 
	\end{enumerate} 
	and where the inner product of two functions $\psi,\phi\in \mathcal{H}_s$ is given by 
	$$\prods{\psi}{\phi}_{\mathcal{H}_s}={\int_{\partial T} \prods{\psi(\zeta)}{\phi(\zeta)}_{\C}\qq \diff \nu(\zeta)}.$$
	Then, $\pi_s$ is defined for all $ \psi \in \mathcal{H}_s$ and for all $t,g\in G$ by 
	$$\pi_s(t)\psi(g)=\big(\sqrt{(d_0-1)(d_1-1)}\big)^{\delta(v,tv,g\omega)}\psi(t^{-1}g).$$
	\begin{theorem}\label{Theorem ce sont les spherique pour les valeurs}
		Let $T$, $\omega$, $G$ and $\chiup_s$ be as in Lemma \ref{Lemma les charracter de Ggamme} and let $v,v'\in V$ be adjacent vertices with respective degree $d$ and $d'$. Then, the representation $\pi_s=\Ind_{G_\omega}^G(\chiup_s)$ is a spherical representation of $G$ for every complex number $s\in \C$ with $\modu{s}=1$. Furthermore, we have that $$\phi_{v,G}(\pi_s)=2 \frac{ \sqrt{(d-1)} }{d} \mbox{\rm Re}(s)$$ when $G$ is vertex-transitive and we have that $$ \phi_{v,G}(\pi_s)=2\frac{\sqrt{d-1}}{d\sqrt{d'-1}} \mbox{\rm Re}(s^2)+   \frac{d'-2}{d (d'-1)}$$ when $G$ has two orbits of vertices.  
	\end{theorem}
	\begin{proof}
		First, notice that
		$$\psi_s: G\rightarrow \C: g\mapsto s^{\delta(v,g^{-1}v,\omega)}$$
		is a $\Fix_G(v)$-invariant unit vector of $\mathcal{H}_s$. On the other hand, since $G$ is $2$-transitive on the boundary $\partial T$, it can be shown that $\psi_s$ is a cyclic vector. In particular, the matrix coefficient $\varphi_{s}:t\mapsto \prods{\pi_s (t)\psi_s}{\psi_s}_{\mathcal{H}_s}$ corresponds to $\pi_s$ under the GNS construction (and not to one of its subrepresentation). To prove that $\pi_s$ is a spherical representation of $G$, we are left to show that $\varphi_{s}$ is a spherical function. Let $q=s\sqrt{(d_0-1)(d_1-1)}$ and notice from the definition of $\varphi_s$ that
		$$\varphi_{s}(t)= \int_{\partial T} q^{\delta(v,tv,\zeta)}\qq \diff \nu(\zeta)\q \forall t\in G.$$
		Now, let $\mu$ be the Haar measure on $G$ such that $\mu(\Fix_G(v))=1$. A  straightforward computation shows for all $g,h\in G$ that
		\begin{equation*}
		\begin{split}
		\int_{\Fix_G(v)}\varphi_s(gkh)\diff \mu(k)&=\int_{\Fix_G(v)}\bigg(\int_{\partial T} q^{\delta(v,gkhv,\zeta)}\qq \diff \nu(\zeta)\bigg)\diff \mu(k)\\
		&=\int_{\Fix_G(v)}\bigg(\int_{\partial T} q^{\delta(v,gv,\zeta)+\delta(gv,gkhv,\zeta)}\qq \diff \nu(\zeta)\bigg)\diff \mu(k)\\
		&=\int_{\partial T} q^{\delta(v,gv,\zeta)}\bigg(\int_{\Fix_G(v)}q^{\delta(v,hv,k^{-1}g^{-1}\zeta)}\qq \diff \mu(k)\bigg)\diff \nu(\zeta)\\
		&=\int_{\partial T} q^{\delta(v,gv,\zeta)}\bigg(\int_{\partial T}q^{\delta(v,hv,\zeta)}\qq \diff \nu(\zeta)\bigg)\diff \nu(\zeta)\\
		&=\int_{\partial T} q^{\delta(v,gv,\zeta)} \varphi_s(h)\diff \nu(\zeta)=\varphi_s(g)\varphi_s (h).
		\end{split}
		\end{equation*}
		This proves as desired that $\varphi_s$ is a $\Fix_G(v)$-spherical function of positive type on $G$. Now, let $\partial T(v,w)=\{\zeta\in \partial T : w\in \lb v,\zeta)\}$, let $\tau\in G$ be a translation of minimal step in $G$ along a geodesic containing $v$ and that 
		\begin{equation*}
		\begin{split}
		\varphi_{s}(\tau)&=\int_{\partial T(v,\tau v)} q^{\delta(v,\tau v , \zeta)} \diff \nu(\zeta) + \int_{\partial T(\tau v,v)} q^{\delta(v,\tau v , \zeta)} \diff \nu(\zeta)\\
		&\q\q \q\q \q\q  \q\q \q\q \q\q \q + \int_{\partial T - (\partial T(v,\tau v) \cup \partial T(\tau v,v))}q^{\delta(v,\tau v , \zeta)}\qq \diff \nu(\zeta).
		\end{split}
		\end{equation*}
		If $G$ is vertex-transitive this computation leads to
		\begin{equation*}
		\begin{split}
		\varphi_{s}(\tau)&=\int_{\partial T(v,\tau v)} q\qq \diff \nu(\zeta) + \int_{\partial T(\tau v,v)} q^{-1} \diff \nu(\zeta)+ 0\\
		&=\nu(\partial T(v,\tau v))q + \nu(\partial T(\tau v,v)) q^{-1}\\
		&= \frac{\sqrt{d-1}}{d}s+ \frac{(d-1)}{\sqrt{d-1}d}s^{-1}\\
		&= 2\frac{\sqrt{d-1}}{d}\mbox{\rm Re}(s).
		\end{split}
		\end{equation*}
		On the other hand, if $G$ has two orbits of vertices this leads to 
		\begin{equation*}
		\begin{split}
		\varphi_{s}(\tau)&=\int_{\partial T(v,\tau v)} q^{2}\diff \nu_v(\zeta) + \int_{\partial T(\tau v,v)} q^{-2} \diff \nu(\zeta)+ \int_{\partial T(v,\tau v)} q^0 \diff  \nu(\zeta)\\
		&=\frac{q}{d(d'-1)}s^2+ \frac{d-1}{d q}s^{-2}+ \frac{d'-2}{d (d'-1)}\\
		&= \frac{\sqrt{d-1}}{d\sqrt{d'-1}} s^2+  \frac{\sqrt{d-1}}{d\sqrt{d'-1}} s^{-2}+  \frac{d'-2}{d (d'-1)}\\
		&= 2\frac{\sqrt{d-1}}{d\sqrt{d'-1}} \mbox{\rm Re}(s^2)+  \frac{d'-2}{d (d'-1)}.
		\end{split}
		\end{equation*}
	\end{proof}
	In particular, if $G$ is vertex transitive the principal series of $G$ corresponds under $\phi_{v,G}$ to the real interval $$\bigg\lb -2\frac{\sqrt{d-1}}{d},  2\frac{\sqrt{d-1}}{d}\bigg\rb,$$ and if $G$ has two orbits of vertices, it corresponds to the real interval $$\frac{d'-2}{d(d'-1)}+\bigg\lb-2\frac{\sqrt{d-1}}{d\sqrt{d'-1}},2\frac{\sqrt{d-1}}{d\sqrt{d'-1}}\bigg\rb.$$ When $d=d'$ the last interval is equal to $$\bigg\lb -\frac{1}{d-1}, \frac{3 d-4}{d(d-1)}\bigg\rb.$$
	\subsection{Special representations}\label{Section special rep}

	 The purpose of this section is to recall the classification of special representations of any closed subgroup of automorphisms of a tree acting $2$-transitively on the boundary. The exposition presented here is based on the approach adopted in \cite{FigaNebbia1991}. However, our presentation highlights slightly different intermediate results for the convenience of the writing of Chapter \ref{Chapter Fell topology of Aut(T)} where we describe the restriction of the Fell topology to spherical and special representations. The classification of special representations is provided by Theorem \ref{thm la classification des speciale} below.

	Let $T$ be a thick semi-regular tree and $G\leq \Aut(T)$ be a closed subgroup acting $2$-transitively on the boundary $\partial T$. Notice from Theorem \ref{Theorem Burger mozes acting on the boundary trans} that $G$ is edge-transitive. Furthermore, we recall that $G$ is either vertex transitive or has two orbits of vertices (those of any pair of adjacent vertices). Let $\sigma$ be a special representation of $G$. The definition of special representations ensures the existence of an edge $e\in E$ and a unit $\Fix_G(e)$-invariant vector $\eta\in \Hr{\sigma}$. On the other hand, since $G$ acts transitively on $E$ and as $\sigma(g)\eta$ defines a unit $\Fix_G(ge)$-invariant vector, $\sigma$ admits $\Fix_G(f)$ invariant vectors for every edge $f\in E$. Now, fix an edge $e\in E$ and consider a unit ${\Fix_G(e)}$-invariant vector $\eta\in \Hr{\sigma}$. The matrix coefficient function 
	\begin{equation*}
	\varphi_{\eta,\eta}:G\rightarrow \C : g\mapsto \prods{\sigma(g)\eta}{\eta}_{\Hr{\sigma}}
	\end{equation*}
	obviously defines a continuous $\Fix_G(e)$-bi-invariant function of positive type. Furthermore, since $\sigma$ does not admit any non-zero invariant vector for the  stabiliser of a vertex, notice for all $g\in G$ and for every $ v\in V(e)$ that 
	\begin{equation*}\label{equation definissant les representations speciale}
	\int_{\Fix_G(v)}\varphi_{\eta,\eta}(gk)\diff \mu(k)=\prods{\sigma(g)\int_{\Fix_G(v)}\sigma(k)\eta\diff \mu(k)}{\eta}_{\Hr{\sigma}}=0.
	\end{equation*}
	The key of the classification of cuspidal representation is to determine the functions satisfying these properties. We start with a technical result concerning the coset dynamic.
	\begin{lemma}\label{Lemma the factorization of $G$ into FixGe double cosets}
		Let $T$ be a thick semi-regular tree, $e\in E$, $G\leq \Aut(T)$ be a closed subgroup acting $2$-transitively on the boundary and let $\tau\in G$ be a translation of minimal step in $G$ along a geodesic $\gamma$ containing $e$. Then, exactly one of the following holds:
		\begin{enumerate}[label=(\roman*)]
			\item $G$ is vertex-transitive, there exists an inversion $h\in G$ of $e$ and
			$$G= \bigsqcup_{n\in \Z} \Fix_G(e)\tau^n\Fix_G(e)\qq \sqcup\qq \bigsqcup_{n\in \Z}\Fix_G(e)\tau^nh\Fix_G(e).$$
			\item $G$ has two orbits of vertices and for every $v\in V(e)$ and each $k_v\in U_v$ where $U_v=\{k\in \Fix_G(v)-\Fix_G(e): \gamma \mbox{ contains }ke\}$ we have that
			$$G= \bigsqcup_{n\in \Z} \Fix_G(e)\tau^n\Fix_G(e)\qq \sqcup\qq \bigsqcup_{n\in \Z} \Fix_G(e)\tau^n k_v\Fix_G(e).$$
		\end{enumerate}
		In particular, in both cases, for every vertex $v\in V(e)$, we have that 
		$$G= \bigsqcup_{n\in \Z} \Fix_G(e)\tau^n\Fix_G(v)\mbox{ and that }G= \bigsqcup_{n\in \Z} \Fix_G(v)\tau^n\Fix_G(v).$$ 
	\end{lemma}
	\begin{proof}
		Since $G$ acts $2$-transitively on the boundary, Theorem \ref{Theorem Burger mozes acting on the boundary trans} ensures that $\Fix_G(w)$ acts transitively on the boundary $\partial T$ for every $w\in V$. In particular, for every element $g\in G$, there exists an element $k\in \Fix_G(e)$ such that $k^{-1}ge$ belongs to $\gamma$. We now treat two cases separately. 
		
		If $G$ is vertex-transitive, $\tau$ has step $1$ and since $k^{-1}ge$ belongs to $\gamma$, there exists an integer $n\in \Z$ such that $\tau^{-n}k^{-1}ge=e$ which proves that $g\in k\tau^n\Stab_G(e)\subseteq \Fix_G(e)\tau^n\Stab_G(e)$. Since $\Stab_G(e)=\Fix_G(e)\sqcup h\Fix_G(e)$ we obtain as desired that  
		$$G=\bigsqcup_{n\in \Z} \Fix_G(e)\tau^n\Fix_G(e)\qq\sqcup\qq\bigsqcup_{n\in \Z} \Fix_G(e)\tau^nh\Fix_G(e).$$
		Furthermore, for every vertex $v\in V(e)$, we have that $\Fix_G(e) \subseteq \Fix_G(v)$ and we either have that $\tau h \Fix_G(e)\subseteq \Fix_G(v)$ of that $\tau^{-1} h \Fix_G(e)\subseteq \Fix_G(v)$. In particular, we also obtain the following decompositions
		$$G= \bigsqcup_{n\in \Z} \Fix_G(e)\tau^n\Fix_G(v)\mbox{ and }G= \bigsqcup_{n\in \Z} \Fix_G(v)\tau^n\Fix_G(v).$$
		
		If on the other hand $G$ has two orbits of vertices, $\tau$ has step $2$. Since $k^{-1}ge$ belongs to $\gamma$ we have two cases depending on the distance between $k^{-1}ge$ and $e$. If there is an odd number of edges in the shortest geodesic containing both $e$ and $ge$, there exists an integer $n\in \Z$ such that $\tau^{-n}k^{-1}ge=e$ and we obtain that $g\in k\tau^{n}\Stab_G(e)\subseteq \Fix_G(e)\tau^n\Stab_G(e)$. On the other hand, if there is an even number of edges in the smallest geodesic containing both $e$ and $ge$, there exists an integer $n\in \Z$ such that $\tau^{-n}k^{-1}ge=k_ve$ and we obtain that $g\in k\tau^{n}k_v\Stab_G(e)\subseteq \Fix_G(e)\tau^nk_v\Stab_G(e)$. Since $\Stab_G(e)=\Fix_G(e)$, we obtain as desired that  
		$$G= \bigsqcup_{n\in \Z} \Fix_G(e)\tau^n\Fix_G(e)\qq \sqcup\qq \bigsqcup_{m\in \Z} \Fix_G(e)\tau^n k_v\Fix_G(e).$$ 
		In particular, since $\Fix_G(e)\subseteq \Fix_G(v)$ and $k_v\Fix_G(e)\subseteq \Fix_G(v)$, we also obtain that 
		$$G= \bigsqcup_{n\in \Z} \Fix_G(e)\tau^n\Fix_G(v)\mbox{ and that }G= \bigsqcup_{n\in \Z} \Fix_G(v)\tau^n\Fix_G(v).$$
	\end{proof}
	Let $\varphi: G\rightarrow \C$ be a $\Fix_G(e)$-bi-invariant function. The above decomposition of $G$ ensures that $\varphi$ is entirely determined by the values it takes on the set $$\{\tau^n:n\in \Z\}\cup \{\tau^nh\}$$ if $G$ is vertex-transitive and by the values it takes on the set $$\{\tau^n:n\in \Z\}\sqcup \{\tau^nk_v: n\in \Z\}$$ where $v\in V(e)$ and $k_v\in U_v$ if $G$ has two orbits of vertices. Now, let $\sigma$ be a special representation of $G$, $\eta\in \Hr{\sigma}$ be a ${\Fix_G(e)}$-invariant vector and $\varphi_{\eta,\eta}: G\rightarrow \C : g\mapsto \prods{ \pi(g)\eta}{\eta}_{\Hr{\sigma}}$. The following result shows that $\varphi_{\eta,\eta}$ is entirely determined by the values it takes on $1_G$ and on any inversion $h\in G$ of the edge $e$ if $G$ is vertex-transitive and is completely determined  by the value it takes on $1_G$ if $G$ has two orbits of vertices. 
	\begin{proposition}\label{proposition compute the special functions in both cases}
			Let $T$ be a thick semi-regular tree, $e\in E$, $v\in V(e)$, $G\leq \Aut(T)$ be a closed subgroup acting $2$-transitively on the boundary, $\tau\in G$ be a translation of minimal step in $G$ and such that $e$ belongs to $\lb v,\tau v\rb$ and $\varphi: G\rightarrow\C$ be a $\Fix_G(e)$-bi-invariant function satisfying 
		\begin{equation}\label{equation definitssant la fonction speciale avec les sommet}
		\int_{\Fix_G(w)}\varphi(gk)\diff \mu(k)=0\q\forall g\in G,\qq  \forall w\in V(e).
		\end{equation}
		If $G$ is vertex-transitive we have for all $n\in \Z$ that 
		\begin{equation*}
		 \varphi(\tau^n)=\bigg(\frac{-1}{d-1}\bigg)^{n}\varphi(h_{\lb n\rb_2})\q \mbox{and}\q\varphi(\tau^nh)=\bigg(\frac{-1}{d-1}\bigg)^{n}\varphi(h_{\lb n +1\rb_2 }).
		\end{equation*}
		where $d$ is the degree of any vertex of $T$, $h\in G$ is any inversion of $e$, $h_0=1_G$, $h_1=h$ and $\lb \cdot\rb_2$ is a notation for modulo $2$. On the other hand, if $G$ has two orbits of vertices, for any $k_v\in U_v=\{k\in \Fix_G(v)-\Fix_G(e): \gamma \mbox{ contains }ke\}$ we have that:
		 \begin{equation*}
			\varphi(\tau^n)=\frac{1}{(d-1)^{\modu{n}}(d'-1)^{\modu{n}}}\varphi(1_G) \q \forall n\in \Z,
		\end{equation*}
		\begin{equation*}
			\varphi(\tau^mk_v)=\frac{-1}{(d-1)^{m}(d'-1)^{m+1}}\varphi(1_G)\q \forall\qq m\geq 1 ,
		\end{equation*}
		\begin{equation*}
			\varphi(\tau^{-(m-1)}k_v)=\frac{-1}{(d-1)^{m}(d'-1)^{m-1}}\varphi(1_G) \q \forall\qq m\geq 1.
		\end{equation*}
	\end{proposition}
	\begin{proof}
		Suppose first that $G$ is vertex-transitive. Let $n\geq 1$ be a positive integer, $v'$ be the vertex of $e$ adjacent to $v$, $h\in G$ be an inversion of $e$ and renormalise the Haar measure $\mu$ of $G$ so that $\mu(\Fix_G(v))=1$. Computing the Equalities (\ref{equation definitssant la fonction speciale avec les sommet}) we obtain from the left $\Fix_G(e)$-invariance of $\varphi$ that  
				\begin{equation*}
				\begin{split}
				0= \int_{\Fix_G(v)}\varphi(\tau^nk)&=\int_{U_v} \varphi(\tau^{n}k)\diff \mu(k)+ \int_{\Fix_G(v)-U_v} \varphi(\tau^{n}k)\diff \mu(k)\\
				&=\int_{U_v}  \varphi(\tau^{n-1}h) \diff \mu(k)+ \int_{\Fix_G(v)-U_v}  \varphi(\tau^{n}) \diff \mu(k)\\
				&=\frac{1}{d}  \varphi(\tau^{n-1}h)  +  \frac{d-1}{d}\varphi(\tau^{n})
				\end{split}
				\end{equation*}
				and that
				\begin{equation*}
				\begin{split}
				0= \int_{\Fix_G(v)}\varphi(\tau^{-(n-1)}k)&=\int_{\Fix_G(e)} \varphi(\tau^{-(n-1)}k)\diff \mu(k)+ \int_{\Fix_G(v)-\Fix_G(e)} \varphi(\tau^{-(n-1)}k)\diff \mu(k)\\
				&=\int_{\Fix_G(e)}  \varphi(\tau^{-(n-1)}) \diff \mu(k)+ \int_{\Fix_G(v)-\Fix_G(e)}  \varphi(\tau^{-n}h) \diff \mu(k)\\
				&=\frac{1}{d} \varphi(\tau^{-(n-1)})  +  \frac{d-1}{d}\varphi(\tau^{-n}h).
				\end{split}
				\end{equation*}
				Applying the same reasoning on the stabiliser $\Fix_G(v')$ we end up with the four following equalities: 
				\begin{equation}\label{equation proposition compute the special functions in both cases 1}
				\varphi(\tau^{n})=-\frac{1}{d-1}  \varphi(\tau^{n-1}h),
				\end{equation}
				\begin{equation}\label{equation proposition compute the special functions in both cases 2}
				\varphi(\tau^{-n}h)=- \frac{1}{d-1}\varphi(\tau^{-(n-1)}),
				\end{equation}	 
				\begin{equation}\label{equation proposition compute the special functions in both cases 3}
				\varphi(\tau^{n}h)=-\frac{1}{d-1}  \varphi(\tau^{n-1}),
				\end{equation}
				\begin{equation}\label{equation proposition compute the special functions in both cases 4}
				\varphi(\tau^{-n})= -\frac{1}{d-1}\varphi(\tau^{-(n-1)}h).
				\end{equation}
			The result follows. Now, suppose that $G$ has two orbits of vertices. Let $n\geq 1$ be a positive integer, $v'$ be the vertex of $e$ adjacent to $v$, $d$ and $d'$ be the respective degree of $v$ and $v'$ and renormalise the Haar measure of $G$ in such a way that $\mu( \Fix_G(e))=1$. Computing the Equalities \eqref{equation definitssant la fonction speciale avec les sommet}, we obtain from the left $\Fix_G(e)$-invariance of $\varphi$ that 
			\begin{equation*}
			\begin{split}
			0=\int_{\Fix_G(v)}\varphi(\tau^n k)\diff \mu(k)&=\int_{\Fix_G(v)-U_v}\varphi(\tau^nk)\diff \mu(k)+ \int_{U_v}\varphi(\tau^nk)\diff \mu(k)\\
			&=\int_{\Fix_G(v)-U_v}\varphi(\tau^n)\diff \mu(k)+ \int_{U_v}\varphi(\tau^nk_v)\diff \mu(k)\\
			&=(d-1) \varphi(\tau^n)+ \varphi(\tau^n k_v),
			\end{split}
			\end{equation*}
			\begin{equation*}
			\begin{split}
			0=\int_{\Fix_G(v)}\varphi(\tau^{-(n-1)} k)\diff \mu(k)&=\int_{\Fix_G(v)-\Fix_G(e)}\varphi(\tau^{-(n-1)}k)\diff \mu(k)+ \int_{\Fix_G(e)}\varphi(\tau^{-(n-1)}k)\diff \mu(k)\\
			&=\int_{\Fix_G(v)-\Fix_G(e)}\varphi(\tau^{-(n-1)k_v})\diff \mu(k)+ \int_{\Fix_G(e)}\varphi(\tau^{-(n-1)}k_v)\diff \mu(k)\\
			&=(d-1) \varphi(\tau^{-(n-1)}k_v)+ \varphi(\tau^{-(n-1)}),
			\end{split}
			\end{equation*}
			\begin{equation*}
			\begin{split}
			0=\int_{\Fix_G(v')}\varphi(\tau^{n-1} k)\diff \mu(k)&=\int_{\Fix_G(v')-\Fix_G(e)}\varphi(\tau^{n-1}k)\diff \mu(k)+ \int_{\Fix_G(e)}\varphi(\tau^{n-1}k)\diff \mu(k)\\
			&=\int_{\Fix_G(v')-\Fix_G(e)}\varphi(\tau^{n}k_v)\diff \mu(k)+ \int_{\Fix_G(e)}\varphi(\tau^{n-1})\diff \mu(k)\\
			&=(d'-1) \varphi(\tau^nk_v)+ \varphi(\tau^{n-1}),
			\end{split}
			\end{equation*}
			and posing $U_{v'}=\{k\in \Fix_G(v'): kv=\tau v\}$ we also obtain that
			\begin{equation*}
			\begin{split}
			0=\int_{\Fix_G(v')}\varphi(\tau^{-n}k)\diff \mu(k)&=\int_{\Fix_G(v')-U_{v'}}\varphi(\tau^{-n}k)\diff \mu(k)+ \int_{U_{v'}}\varphi(\tau^{-(n-1)}k)\diff \mu(k)\\
			&=\int_{\Fix_G(v')-U_{v'}}\varphi(\tau^{-n})\diff \mu(k)+ \int_{U_{v'}}\varphi(\tau^{-(n-1)}k_v)\diff \mu(k)\\
			&=(d'-1) \varphi(\tau^{-n})+ \varphi(\tau^{-(n-1)} k_v).
			\end{split}
			\end{equation*}
			We end up for every integer $n\geq 1$ with the four following equalities: 
			\begin{equation}\label{equation proposition compute the special functions in both cases 1'}
			\varphi(\tau^n)= \frac{1}{(d-1)(d'-1)}\varphi(\tau^{n-1}),
			\end{equation}
			\begin{equation}\label{equation proposition compute the special functions in both cases 2'}
			\varphi(\tau^{-n})= \frac{1}{(d-1)(d'-1)}\varphi(\tau^{-(n-1)}),
			\end{equation}	 
			\begin{equation}\label{equation proposition compute the special functions in both cases 3'}
			\varphi(\tau^nk_v)= -\frac{1}{d'-1} \varphi(\tau^{n-1}),
			\end{equation}
			\begin{equation}\label{equation proposition compute the special functions in both cases 4'}
			 \varphi(\tau^{-(n-1)}k_v)= -(d'-1)\varphi(\tau^{-n}).
			\end{equation}
			The result follows. 
	\end{proof}
	It is now time to study the consequence of this result on the special representations of $G$. We start with a reminder on the notion of square-integrable representations. 
	\begin{definition}
		A representation $\pi$ of a unimodular locally compact group is \tg{square-integrable} if there exists a non-zero vector $\xi\in \Hr{\pi}$ such that the matrix coefficient 
		$$\varphi_{\xi,\xi}: G \rightarrow \C: g\mapsto \prods{ \pi(g)\xi}{\xi}_{\Hr{ \pi}}$$ 
		belongs to $L^2(G)$. 
	\end{definition}
	\noindent This concepts has a particular interest due to the following correspondence.
	\begin{lemma}[{\cite[Lemma 14.1.1 and 14.1.2]{Dixmier1977}}]\label{lemma dixmier sur les reps de carre integrable}
		Let $G$ be a unimodular locally compact group. Then, the following holds: 
		\begin{enumerate}[label=(\roman*)]
			\item Every subrepresentation of the left-regular representation of $G$ is square-integrable.
			\item Every irreducible square-integrable representation of $G$ is a subrepresentation of the left-regular representation of $G$. 
		\end{enumerate}  
	\end{lemma}
	Relying on the explicit description of $\varphi_{\eta,\eta}$ provided by Proposition \ref{proposition compute the special functions in both cases}, the following result ensures that every special representation of $G$ is square-integrable.
	\begin{lemma}\label{Lemma the fucntions of positive type are square integrable}
			Let $T$ be a thick semi-regular tree, $e\in E$, $G\leq \Aut(T)$ be a closed subgroup acting $2$-transitively on the boundary and let $\varphi: G\rightarrow\C$ be a $\Fix_G(e)$-bi-invariant function of $G$ satisfying 
			\begin{equation*}
			\int_{\Fix_G(v)}\varphi(gk)\diff \mu(k)=0\q\forall g\in G,\qq  \forall v\in V(e).
			\end{equation*}
			Then, $\varphi$ is square-integrable. 
	\end{lemma}
	\begin{proof}
		Suppose first that $G$ is vertex-transitive and renormalise the Haar measure $\mu$ of $G$ so that $\mu( \Fix_G(e))=1$. Proposition \ref{proposition compute the special functions in both cases} ensures for all $n\in \Z$ and for all $k,k'\in \Fix_G(e)$ that 
		\begin{equation*}
		\varphi(k\tau^nk')=\bigg(\frac{-1}{d-1}\bigg)^{\modu{n}}\varphi(h_{\lb n\rb_2})\q \mbox{and}\q\varphi(k\tau^nhk')=\bigg(\frac{-1}{d-1}\bigg)^{\modu{n}}\varphi(h_{\lb n +1\rb_2 }).
		\end{equation*}
		where $d$ is the degree of any vertex of $T$, $h\in G$ is any inversion of $e$, $h_0=1_G$, $h_1=h$.
		It follows from Lemma \ref{Lemma the factorization of $G$ into FixGe double cosets} that 
		\begin{equation*}
		\begin{split}
		\int_{G} \modu{\varphi (g)}^2\diff \mu(g)&=\sum_{n\in \Z}\int_{\Fix_G(e)\tau^n\Fix_G(e)}\modu{\varphi (g)}^2\diff \mu(g)\\
		&\q\q\q\q\q\q+ \sum_{n\in \Z}\int_{\Fix_G(e)\tau^nh\Fix_G(e)}\modu{\varphi (g)}^2\diff \mu(g)\\
		&= \sum_{n\in \Z}\bigg(\frac{-1}{d-1}\bigg)^{2\modu{n}} \mu(\Fix_G(e)\tau^n\Fix_G(e))\modu{\varphi (h_{\lb n\rb_2})}^2\\
		&\q\q\q\q\q+ \sum_{n\in \Z}\bigg(\frac{-1}{d-1}\bigg)^{2\modu{n}} \mu(\Fix_G(e)\tau^nh\Fix_G(e)) \modu{\varphi (h_{\lb n+1\rb_2})}^2\\
		&=(\modu{\varphi(1_G)}^2+\modu{\varphi (h)}^2)\sum_{n\in \Z}\bigg(\frac{1}{d-1}\bigg)^{\modu{n}} <\infty.
		\end{split}
		\end{equation*}
		The last equality follows from the fact that $$\mu(\Fix_G(e)\tau^n\Fix_G(e))= \mu(\Fix_G(e)\tau^nh\Fix_G(e))=(d-1)^n\mu(\Fix_G(e))$$
		since $G$ is $2$-transitive on the boundary and since there is exactly $(d-1)^n$ vertices at distance $n$ from $e$ on one or the other side of $e$.
		Suppose now that $G$ has two orbits of vertices. Proposition \ref{proposition compute the special functions in both cases} ensures for all $n\in \Z$, for all positive integer $m\geq 1$ and for all $k,k'\in \Fix_G(e)$ that 
		\begin{equation*}
		\varphi (k\tau^nk')=\frac{1}{(d-1)^{\modu{n}}(d'-1)^{\modu{n}}}\varphi(1_G),
		\end{equation*}
		\begin{equation*}
		\varphi (k\tau^mk_vk')=\frac{-1}{(d-1)^{m}(d'-1)^{m+1}}\varphi(1_G),
		\end{equation*}
		and
		\begin{equation*}
		\varphi (k\tau^{-(m-1)}k_vk')=\frac{-1}{(d-1)^{m}(d'-1)^{m-1}}\varphi(1_G).
		\end{equation*}
		It follows from Lemma \ref{Lemma the factorization of $G$ into FixGe double cosets} that 
		\begin{equation*}
		\begin{split}
		\int_{G} \modu{\varphi (g)}^2\diff \mu(g)&=\sum_{n\in \Z}\int_{\Fix_G(e)\tau^n\Fix_G(e)}\modu{\varphi (g)}^2\diff \mu(g)\\
		&\q\q\q\q\q\q+ \sum_{n\in \Z}\int_{\Fix_G(e)\tau^nk_v\Fix_G(e)}\modu{\varphi (g)}^2\diff \mu(g)\\
		&= \sum_{n\in \Z}\bigg(\frac{1}{(d-1)(d'-1)}\bigg)^{2{\modu{n}}} \modu{\varphi(1_G)}^2\mu(\Fix_G(e)\tau^n\Fix_G(e))\\
		&\q\q+ \frac{1}{(d-1)^2}\sum_{m\in \N}\bigg(\frac{1}{(d-1)(d'-1)}\bigg)^{2{m}}  \modu{\varphi(1_G)}^2\mu(\Fix_G(e)\tau^{-m}k_v\Fix_G(e)) \\	
		&\q\q + \frac{1}{(d'-1)^2}\sum_{m\in \N, m\geq 1}\bigg(\frac{1}{(d-1)(d'-1)}\bigg)^{2m}  \modu{\varphi(1_G)}^2 \mu(\Fix_G(e)\tau^{m}k_v\Fix_G(e)) \\		
		&=  \modu{\varphi(1_G)}^2\sum_{n\in \Z}\bigg(\frac{1}{(d-1)(d'-1)}\bigg)^{\modu{n}}+  \frac{\modu{\varphi(1_G)}^2}{d-1}\sum_{m\in\N}\bigg(\frac{1}{(d-1)(d'-1)}\bigg)^{m}\\
		&\q\q\q\q\q\q\q\q\q +   \frac{\modu{\varphi(1_G)}^2}{d'-1}\sum_{m\in\N,m\geq 1}\bigg(\frac{1}{(d-1)(d'-1)}\bigg)^{m+1}<\infty.
		\end{split}
		\end{equation*}
		The last equality follows from the fact that $$\mu(\Fix_G(e)\tau^n\Fix_G(e))= (d-1)^{\modu{n}}(d'-1)^{\modu{n}}\mu(\Fix_G(e))\q \forall n\in \Z,$$
		$$ \mu(\Fix_G(e)\tau^{-(m-1)}k_v\Fix_G(e))= (d-1)(d-1)^{m-1}(d'-1)^{m-1}\mu(\Fix_G(e))\q \forall m\geq 1,$$
		$$\mbox{ and } \mu(\Fix_G(e)\tau^{m}k_v\Fix_G(e))= (d'-1)(d-1)^{m-1}(d'-1)^{m-1}\mu(\Fix_G(e))\q \forall m\geq 1$$
		when $G$ is $2$-transitive on the boundary.
	\end{proof}
	Now, let $G\leq \Aut(T)$ be a closed subgroup acting $2$-transitively on the boundary, $\sigma$ be a special representation of $G$, $\eta$ be a unit $\Fix_G(e)$-invariant vector and consider the matrix coefficient
	$$\varphi_{\eta,\eta} : G\rightarrow \C: g\mapsto \prods{\sigma(g)\eta}{\eta}_{\Hr{\sigma}}.$$
	This function satisfies the hypothesis of the above Lemma and is therefore square-integrable. Since $G$ is unimodular, it follows from Lemma \ref{lemma dixmier sur les reps de carre integrable} that $\sigma$ is a subrepresentation of the left-regular representation of $G$. Our next task will be to provide an explicit description of this representation. We let $\mathcal{L}(e)$ be the closure in $L^2(G)$ of all functions $\varphi: G\rightarrow \C$ satisfying the following properties:
	\begin{enumerate}[label=(\roman*)]
		\item $\varphi$ is $\Fix_G(f)$-left-invariant for some $f\in E$.
		\item $\varphi$ is $\Fix_G(e)$-right-invariant.
		\item For all $v\in V(e)$ and every $g\in G$, we have that
		$$\int_{\Fix_G(v)}\varphi(gk)\diff \mu(k)=0.$$
	\end{enumerate}
	Since $G$ is edge-transitive, $\mathcal{L}(e)$ is exactly the closure of the subspace of $L^2(G)$ spanned by the left-translates of functions described by Proposition \ref{proposition compute the special functions in both cases}. In particular, this space is $\lambda_G$-invariant. We denote by $(\sigma,\mathcal{L}(e))$ the corresponding subrepresentation of the regular representation of $G$. The following result enures that this representation contains all the special representations of $G$ as discrete summand.
	\begin{proposition}\label{proposition, every special representation appears a s subrepresentation of sigma L e}
		Let $T$ be a thick semi-regular tree, $G\leq \Aut(T)$ be a closed subgroup acting $2$-transitively on the boundary $\partial T$ and let $e\in E$. Then, every special representation of $G$ is a subrepresentation of $(\sigma,\mathcal{L}(e))$.
	\end{proposition}
	\begin{proof}
		Let $\pi$ be a special representation of $G$ and let $\eta\in \Hr{\pi}$ be a unit $\Fix_G(e)$-invariant vector. For all $h\in G$ the function
		$$\overline{\varphi_{\eta,\pi(h) \eta}}: G\rightarrow \C: g\mapsto \prods{\pi(h)\eta}{\pi(g)\eta}_{\Hr{\pi}}$$
		coincides with $\sigma(h)\overline{\varphi_{\eta,\eta}}$ and therefore belongs to $\mathcal{L}(e)$. Let $\mathcal{H}$ be the subspace of $\Hr{\pi}$ spanned by $\{\pi(h)\eta: h\in G\}$ and let $\mathcal{L}$ be subspace of $\mathcal{L}(e)$ spanned by $\{\overline{\varphi_{\eta,\pi(h)\eta}}:h\in G\}$. The map $$\Phi:\mathcal{H}\rightarrow \mathcal{L}: \xi\mapsto \overline{\varphi_{\eta,\xi}}$$ defines a linear isomorphism between these spaces. On the other hand, $\Phi$ intertwines the representations $\pi$ and $\sigma$ since for all $g,t\in G$ and for every $\xi\in \mathcal{H}$ we have that 
		$$\sigma(t)\big\lb\Phi(\xi)\big\rb(g)=\overline{\varphi_{\eta,\xi}}(t^{-1}g)=\overline{\varphi_{\eta,\pi(t)\xi}}(g)=\big\lb\Phi(\pi(t)\xi)\big\rb(g).$$ Furthermore, \cite[Theorem 14.3.3.]{Dixmier1977} ensures the existence of a positive constant $d_\pi$ called the formal dimension of $\sigma$ such that 
		\begin{equation*}
		\begin{split}
		\prods{\Phi(\pi(h)\eta)}{\Phi(\pi(h')\eta)}_{L^2(G)}&= \int_{G}\overline{\varphi_{\eta,\pi(h)\eta}}(g)\varphi_{\eta,\pi(h')\eta}(g)\diff \mu(g)\\
		&= d_\pi^{-1}\prods{\eta}{\eta}_{\Hr{\pi}}\prods{\pi(h)\eta}{\pi(h')\eta}_{\Hr{\pi}}\\
		&= d_\pi^{-1}\prods{\pi(h)\eta}{\pi(h')\eta}_{\Hr{\pi}}.
		\end{split}
		\end{equation*}
		Thus, the map $\Psi= \sqrt{d_\pi}\Phi$ is a unitary operator on the pre-Hilbert space $\mathcal{H}$ intertwining $\pi$ and $\sigma$. Since $\pi$ is irreducible, the space $\mathcal{H}$ is dense in $\Hr{\pi}$ and $\Psi$ extends continuously in a unitary operator from $\Hr{\pi}$ to $\mathcal{L}(e)$ intertwining $\pi$ and a subrepresentation of $\sigma$.
		\end{proof}
	Our next task will be to show that each subrepresentation of $(\sigma, \mathcal{L}(e))$ is a special representation of $G$ and to identify explicitly these irreducible representations of $G$. We start by showing that each subrepresentation of $\sigma$ admits a non-zero $\Fix_G(e)$-invariant vector.
	\begin{lemma}\label{lemma every invariant subspace contains non trivial fix e invariant}
		Every non-trivial closed $\sigma$-invariant subspace of $\mathcal{L}(e)$ contains a non-trivial $\Fix_G(e)$-left-invariant function. 
	\end{lemma}
	\begin{proof}
		Consider a non-trivial closed $\sigma$-invariant subspace $M$ of $\mathcal{L}(e)$, $\psi\in M$ be a non-trivial function and $g\in G$ be such that $\psi(g)\not=0$. The function 
		$$\varphi=\int_{\Fix_G(e)} \sigma(kg^{-1})\psi\diff \mu(k)$$
		belongs to $M$, is $\Fix_G(e)$-left-invariant and satisfies
		$$\varphi(1_G)= \int_{\Fix_G(e)}\psi(gk^{-1})\diff \mu(k)=\mu(\Fix_G(e))\psi(g)\not=0.$$
	\end{proof}
	On the other hand, due to the properties of the $L^2$-convergence, each element $\varphi$ of $\mathcal{L}(e)$ is $\Fix_G(e)$-right invariant and satisfies 
	$$\int_{\Fix_G(v)}\varphi(gk)\diff \mu(k)=0\q \forall v\in V(e),\qq \forall g\in G.$$
	In particular, the non-trivial $\Fix_G(e)$-left-invariant functions of $\mathcal{L}(e)$ are completely determined by Proposition \ref{proposition compute the special functions in both cases}. In particular, this proposition ensures that the subspace of such functions has dimension two if $G$ is vertex transitive and has dimension $1$ is $G$ has two orbits of vertices. Thus, when $G$ has two orbits of vertices, it follows from Lemma \ref{lemma every invariant subspace contains non trivial fix e invariant} that $(\sigma,\mathcal{L}(e))$ is irreducible and is the unique special representation of $G$. We are left to treat the vertex transitive case. This is the purpose of the following result. For every $\epsilon\in \{-1,+1\}$ we let $\mathcal{L}(e)_{\epsilon}$ the subspace of functions of $\mathcal{L}(e)$ satisfying $\varphi(gh)=\epsilon\varphi(g)$ for all $g\in G$ and every inversion $h\in G$ of the edge $e$ and  let $\sigma^\epsilon$ denotes the restriction of the action of $\sigma$ to $\mathcal{L}(e)_{\epsilon}$. 
	\begin{lemma}\label{lemma if G trans L(e) decompose en somme direct de 2 rep}
		Let $G\leq \Aut(T)$ be a closed vertex-transitive subgroup acting $2$-transitively on the boundary and $e\in E$. Then,  $(\sigma^{-1},\mathcal{L}(e)_{-1})$ and $(\sigma^{+1},\mathcal{L}(e)_{+1})$ are inequivalent special representations and $\sigma \simeq \sigma^{-1}\oplus \sigma^{+1}$.
	\end{lemma}
	\begin{proof}
		 Choose an inversion $h\in G$ of the edge $e$ and consider the linear map $\nu:\mathcal{L}(e) \rightarrow \mathcal{L}(e)$ defined for all $\varphi\in \mathcal{L}(e)$ and every $g\in G$ by $\nu(\varphi)(g)=\varphi(gh)$. This map is well defined since for all $\varphi\in \mathcal{L}(e)$, $g\in G$ and $v\in V(e)$ we have that
		\begin{equation*}
		\begin{split}
		\int_{\Fix_G(v)} (\nu\varphi)(gk)\diff \mu(k)&=\int_{\Fix_G(v)} \varphi(gkh)\diff \mu(k)\\
		&=\int_{\Fix_G(v)} \varphi(ghh^{-1}kh)\diff \mu(k)\\
		&=\int_{\Fix_G(h^{-1}v)} \varphi(ghk)\diff \mu(k)=0.
		\end{split}
		\end{equation*} 
		Since every element of $\mathcal{L}(e)$ is $\Fix_G(e)$-right invariant, $\nu$ does not depend on our choice of inversion $e$ and defines an isometric involution on $\mathcal{L}(e)$. Notice for each $\epsilon\in \{-1,1\}$, that $\mathcal{L}(e)_{\epsilon}$ is the $\epsilon$-eigenspace of $\nu$. On the other hand, one has that
		$$\sigma(t)\nu(\varphi)(g)=\varphi(t^{-1}gh)=\nu(\sigma(t)\varphi)(g)\q \forall g,t\in G,\qq \varphi\in\mathcal{L}(e).$$
		This proves that each $\mathcal{L}(e)_{\epsilon}$ is a $\sigma$-left invariant subspace of $\mathcal{L}(e)$. 
		Our next task is to prove that
		$\mathcal{L}(e)=\mathcal{L}(e)_{-1}\oplus\mathcal{L}(e)_{+1}.$ Notice that $\mathcal{L}(e)_{-1}$ and $\mathcal{L}(e)_{+1}$ are orthogonal with respect to one another since for all $\varphi\in \mathcal{L}(e)_{-1},\qq \psi\in \mathcal{L}(e)_{+1}$ one has that
		$$\prods{\varphi}{\psi}_{L^2(G)}=\prods{\nu(\varphi)}{\nu(\psi)}_{L^2(G)}=\prods{-\varphi}{\psi}_{L^2(G)}=-\prods{\varphi}{\psi}_{L^2(G)}.$$
		Now, we recall from Proposition \ref{proposition compute the special functions in both cases} and  Lemma \ref{Lemma the fucntions of positive type are square integrable}, that for each $\epsilon\in \{-1,+1\}$, there is a unique element $\varphi^\epsilon$ of $\mathcal{\mathcal{L}}(e)^{\Fix_G(e)}$ satisfying by $\varphi^\epsilon(1_G)=1$, $\varphi^\epsilon(h)=\epsilon$ and notice that $\nu(\varphi^\epsilon)=\epsilon\varphi.$ We let $\mathcal{L}$ be the subspace of $\mathcal{L}(e)$ spanned by 
		$$\{\sigma(g)\varphi: g\in G,\qq \varphi\in \mathcal{L}(e)^{\Fix_G(e)}\}$$
		and recall that $\mathcal{L}$ is dense in $\mathcal{L}(e)$. Now, notice that $\mathcal{L}(e)^{\Fix_G(e)}$ is spanned by $\{\varphi^{-1}, \varphi^{+1}\}$, so that $\mathcal{L}$ is spanned by 
		$$\{\sigma(g)\varphi^-:g\in G\}\sqcup \{\sigma(g)\varphi^+: g\in G\}.$$ It follows that $\mathcal{L}(e)=\mathcal{L}(e_{-1})\oplus \mathcal{L}(e)_{+1}$ and hence that $\sigma \simeq \sigma^{-1}\oplus \sigma^{+1}$. On the other hand, the space of $\Fix_G(e)$-left-invariant functions of $\mathcal{L}(e)_\epsilon$ is one-dimensional for each $\epsilon\in \{-1,+1\}$. It follows from Lemma \ref{lemma every invariant subspace contains non trivial fix e invariant} that both $\sigma^{-1}$ and $\sigma^{+1}$ are irreducible and hence define special representations of $G$. Since $\varphi^\epsilon$ is the unique function of positive type associated with $\sigma^\epsilon$ corresponding to a unit $\Fix_G(e)$-invariant vector and since $\varphi^{-1}\not=\varphi^{+1}$ one sees that these irreducible representations are also inequivalent.
	\end{proof}
	In light of Proposition \ref{proposition, every special representation appears a s subrepresentation of sigma L e} this leads to the following classification of special representations of $G$.
	\begin{theorem}[{\cite[Chapter III, Section 2]{FigaNebbia1991}}, {\cite[Section 5.6]{Matsumoto1977}}]\label{thm la classification des speciale}
		Let $T$ be a thick semi-regular tree, $e\in E$ and $G\leq \Aut(T)$ be a closed subgroup acting $2$-transitively on the boundary $\partial T$. Every special representation of $G$ is square-integrable and admits a $\Fix_G(f)$-invariant vector for each $f\in E$. Furthermore:
		\begin{enumerate}[leftmargin=*, label=(\roman*)]
			\item If $G$ is vertex-transitive, there are two equivalence classes of special representations: these of $(\sigma^{-1},\mathcal{L}(e)_{-1})$ and $(\sigma^{+1},\mathcal{L}(e)_{+1})$.
			\item If $G$ has two orbits on $V$, there is a unique equivalence class of special representations of $G$: the equivalence class of $(\sigma,\mathcal{L}(e))$ 
		\end{enumerate}
	\end{theorem}
	\subsection{Cuspidal representations}\label{cuspidal representations of the group of automorphism of a tree}
	Let $T$ be a thick semi-regular tree. The purpose of this section is to recall the classification of the cuspidal representations of the full automorphism group of $T$. This classical result was achieved by Ol'shanskii \cite{Olshanskii1977}. By contrast with Sections \ref{Section spherical rep} and \ref{Section special rep} the classification of the cuspidal representations presented here does not apply to every closed subgroup $G\leq \Aut(T)$ acting $2$-transitively on the boundary. For instance, the key Lemma \ref{independence figa nebbia 3.1} does not apply to $\mbox{PSL}_2(\Q_p)$ (which can be realised as a closed subgroup acting $2$-transitively on the boundary of a $p+1$-regular tree \cite[Chapter II]{Serre1980}). The reasoning presented here applies to the family of closed subgroups $G\leq \Aut(T)$ satisfying the Tits independence property. We now recall the main ideas of this classification and refer to \cite{FigaNebbia1991} and \cite{Amann2003} for details. 
	
	\begin{definition}\label{equation prop ind de tits}
		A group $G\leq \Aut(T)$ has the \tg{Tits independence property} if for any two adjacent vertices $v,v'\in V$, the pointwise stabiliser of the edge $\{v,v'\}$ satisfies the following factorisation
		$$\Fix_G(\{v,v'\})= \Fix_G(T(v,v')) \Fix_G(T(v',v))$$ where $T(w,v)=\{u\in V: d_T(w,u)< d_T(v,u)\}$.
	\end{definition}
	In other words, a group  $G$ of automorphisms of tree satisfies the Tits independence property if every element of $G$ fixing an edge decomposes as a product of an automorphism of $G$ fixing all the vertices on one side of the edge and an automorphism fixing every vertices on the other side of the edge. On the other hand, the Tits independence property can be realised as a factorisation property of the pointwise stabilisers of complete finite subtrees of $T$. We start by formalising this idea. For each subtree $\mathcal{T}$ of $T$ we let
	\begin{equation*}\label{definition du set partial Eo(T)}
	E_o(\mathcal{T})=\{(v,w)\in  V(\mathcal{T})\times  V(\mathcal{T}): \qq v \mbox{ and } w \mbox{ are adjacent}\}
	\end{equation*}
	be the set of ordered pairs of distinct adjacent vertices of $\mathcal{T}$. The elements of $E_o(\mathcal{T})$ are called the \tg{oriented edges} of $\mathcal{T}$. For any oriented edge $f=(w,v)\in E_o(\mathcal{T})$ we say that $w$ is the \tg{origin} of $f$, that $v$ is the \tg{terminal} vertex of $f$,  that $f$ is \tg{terminal} in $\mathcal{T}$ if the terminal vertex of $f$ is a leaf of $\mathcal{T}$ and we denote by $\bar{f}$ the oriented edge $(v,w)$ with opposite orientation.  
	We denote by $\partial E_o(\mathcal{T})$ the set of all terminal edges of $\mathcal{T}$ and for every oriented edge $f=(w,v)\in E_o(T)$, we let $$Tf=T(w,v)=\{u\in V: d_T(w,u)< d_T(v,u)\}.$$ 
	As announced, the following result provides a characterisation of groups with the Tits independence property in terms of a factorisation property of the pointwise stabilisers of complete finite subtrees of $T$. Let $\mathfrak{T}$ be the set of complete finite subtrees for $T$.
	\begin{lemma}\label{lemma alternative definition of tites indep prop}
		Let $G\leq \Aut(T)$. Then $G$ satisfies the Tits independence property if and only if for every $\mathcal{T}\in \mathfrak{T}$ containing an edge we have that
		\begin{equation}\label{equation pour la prop d'indep de tits}
		\Fix_G(\mathcal{T})=\prod_{f\in \partial E_o(\mathcal{T})} \Fix_{G}(Tf)\cap \Fix_G(\mathcal{T}).
		\end{equation}
	\end{lemma}
	\begin{proof}
		The reverse implication is trivial as every edge of $T$ is a complete finite subtree of $T$ containing an edge. To prove the forward implication we apply an induction on the number of interior vertices of $\mathcal{T}$ and treat several cases. If $\mathcal{T}$ does not have an interior vertex, it is an edge and the group equality corresponds exactly to the Tits independence property. If $\mathcal{T}$ has exactly one interior vertex, there exists $v\in V$ such that $\mathcal{T}=B_T(v,1)$. Let $g\in \Fix_G(v)$, let $\partial E_o(\mathcal{T})=\{f_1,..., f_n\}$ and notice that  $\Fix_G(B_T(v,1))\subseteq \Fix_G(f)$ for all $f\in \partial E_o(\mathcal{T})$. Since $G$ satisfies the Tits independence property we have a decomposition $g=g_{f_1}g_{\bar{f_1}}$ where $g_{f_1}\in \Fix_{G}(Tf_1)$ and $g_{\bar{f_1}}\in  \Fix_{G}(T\bar{f_1})$. Furthermore, notice that $\Fix_G(Tf\cup f)= \Fix_G(Tf)$ for every $f\in E_o(\mathcal{T})$. On the other hand, since $f_1\in \partial E_o(\mathcal{T})$ we have that $B_T(v,1)\subseteq  T{f_1}\cup f_1$ which implies that $\Fix_{G}(T{f_1})\subseteq \Fix_G(B_T(v,1))$. This proves that $g_{f_1}\in \Fix_{G}(B_T(v,1))$. Since $g=g_{f_1}g_{\bar{f_1}}$ and as $g\in \Fix_G(B_T(v,1))$ we obtain that $g_{\bar{f_1}}\in  \Fix_{G}(T\bar{f_1})\cap \Fix_G(B_T(v,1))$. It follows that $g_{\bar{f_1}}\in \Fix_G(f_2)$ and since $G$ satisfies the Tits independence property we have a decomposition $g_{\bar{f_1}}=g_{f_2} g_{\bar{f_2}}$ where $g_{f_2}\in \Fix_{G}(T{f_2})$ and $g_{\bar{f_2}}\in \Fix_{G}(T\bar{f_2})$. Just as before, since $f_2\in \partial E_o(\mathcal{T})$, we have that $B_T(v,1)\subseteq Tf_2\cup f_2$. This implies that $\Fix_{G}(T\bar{f_2})\subseteq \Fix_G(B_T(v,1))$ and it follows that $g_{\bar{f_2}}\in\Fix_{G}(T{\bar{f_2}})$. On the other hand, $T{\bar{f_1}}\subseteq T{f_2}$ which implies that $g_{f_2}\in \Fix_G(T{\bar{f_1}})$. Since $g_{\bar{f_1}}\in \Fix_G(T{\bar{f_1}})$, the above decomposition of $g_{\bar{f_1}}$ implies that $g_{\bar{f_2}}\in \Fix_G(T{\bar{f_1}})$ and we obtain that $g_{\bar{f_2}}\in \Fix_{G}(T\bar{f_1}\cup T\bar{f_2})\cap \Fix_G(B_T(v,1))$. Proceeding iteratively, we obtain that 
		\begin{equation*}
			g= g_{f_1}g_{f_2}... g_{f_n}g_{\bar{f_n}}
		\end{equation*}
		for some $g_{f_i}\in \Fix_{}(T{f_i})\cap \Fix_G(B_T(v,1))$ and some $g_{\bar{f_n}}\in \Fix_{G}(\bigcup_{i=1}^nT{\bar{f_i}})\cap \Fix_G(B_T(v,1))$. In particular since $\{f_1,..., f_n\}=\partial E_o(\mathcal{T})$, notice that $g_{\bar{f_n}}\in \Fix_G(T)=\{1_G\}$. This proves as desired that
		\begin{equation*}
			g\in \prod_{f\in \partial E_o(\mathcal{T})}^{} \big\lb\Fix_{G}(T{f})\cap \Fix_G(B_T(v,1))\big\rb.
		\end{equation*}
		If $\mathcal{T}$ has at least two interior vertices, the reasoning is similar. We let $g\in \Fix_G(\mathcal{T})$, $v$ be an interior vertex of $\mathcal{T}$ at distance $1$ form the boundary $\partial\mathcal{T}$ and $\mathcal{R}$ be the unique maximal proper complete subtree of $\mathcal{T}$ for which $v$ is not an interior vertex. Notice that $\mathcal{R}$ is a complete subtree of $\mathcal{T}$ with one less interior vertex. In particular, our induction hypothesis, ensures that
		\begin{equation*}
			\Fix_G(\mathcal{R})= \prod_{f\in \partial E_o(\mathcal{R})}^{} \big\lb\Fix_{G}(Tf)\cap \Fix_G(\mathcal{R})\big\rb.
		\end{equation*}
		Now, notice that exactly one extremal edge $e\in \partial E_o(\mathcal{R})$ of $\mathcal{R}$ is not extremal in $\mathcal{T}$; the oriented edge $e$ of $\mathcal{R}$ with terminal vertex $v$. Furthermore, since $\mathcal{R}\subseteq \mathcal{T}$, we have that $g\in\Fix_G(\mathcal{R})$ which implies that
		\begin{equation*}
			g= g_e\prod_{f\in \partial E_o(\mathcal{T})\cap \partial E_o(\mathcal{R})}^{} g_f
		\end{equation*} 
		for some $g_e\in \Fix_{G}(Te)\cap \Fix_G(\mathcal{R})$ and some $g_f\in \Fix_{G}(Tf)\cap \Fix_G(\mathcal{R})$. Furthermore, for every $f\in \partial E_o(\mathcal{T})\cap \partial E_o(\mathcal{R})$ notice that $$\Fix_{G}(Tf)\cap \Fix_G(\mathcal{R})=\Fix_{G}(Tf)\cap\Fix_G(\mathcal{T})$$ which implies that $g_f \in \Fix_{G}(Tf)\cap \Fix_G(\mathcal{T})$. In particular, since $g\in \Fix_G(\mathcal{T})$, the above decomposition implies that $g_e\in\Fix_G(Te)\cap \Fix_G(\mathcal{T})$. On the other hand, $B_T(v,1)\subseteq T{e} \cup\mathcal{T}$ which implies that $g_e\in \Fix_G(B_T(v,1))$. Hence, by the first part of the proof, we have that
		\begin{equation*}
			g_e = g_{\bar{e}} \prod_{b\in \partial E_o(\mathcal{T}) \cap \partial E_o(B_T(v,1))} g_{b}
		\end{equation*}
		for some $g_{\bar{e}}\in \Fix_{G}(T{\bar {e}})\cap \Fix_G(B_T(v,1))$ and $g_{b}\in \Fix_{G}(Tb)\cap \Fix_G(B_T(v,1))$. Furthermore, for every $b\in  \partial E_o(\mathcal{T}) \cap \partial E_o(B_T(v,1))$, notice that $\mathcal{T}\subseteq Tb \cup B_T(v,1)$ and thus that $g_{b}\in \Fix_{G}(Tb)\cap \Fix_G(\mathcal{T})$. In particular, since $g_e\in \Fix_{G}(Te)\cap \Fix_G(\mathcal{T})$ the above decomposition of $g$ implies that $g_{\bar{e}}=1_G$. This proves as desired that $g$ belongs to 
		\begin{equation*}
			\prod_{f\in \partial E_o(\mathcal{T})\cap \lb \partial E_o(\mathcal{R})\cup\partial E_o(B_T(v,1))\rb }^{} \big\lb\Fix_{G}(Tf)\cap \Fix_G(\mathcal{T})\big\rb.
		\end{equation*}
	\end{proof}	 
	The following direct consequence is the key of the classification of the cuspidal representations (see \cite[Lemma $3.1$]{FigaNebbia1991} for the full group of automorphisms of a regular tree and \cite[Lemma $19$]{Amann2003} for the general case of groups satisfying the Tits independence property).
	\begin{lemma}\label{independence figa nebbia 3.1}
		Let $G\leq \Aut(T)$ be a closed subgroup satisfying the Tits independence property and let $\mathcal{T}$, $\mathcal{T}'\in \mathfrak{T}$ be such that $\mathcal{T}$ contains at least one interior vertex and $\mathcal{T}'$ does not contain $\mathcal{T}$. Then, there exists a complete proper subtree $\mathcal{R}$ of $\mathcal{T}$ such that $\Fix_{G}(\mathcal{R})\subseteq \Fix_{G}(\mathcal{T}') \Fix_{G}(\mathcal{T})$. 
	\end{lemma}
	\begin{proof}
		Since $\mathcal{T}$ and $\mathcal{T}'$ are complete and since $\mathcal{T}'$ does not contain $\mathcal{T}$ there exists a vertex $w_\mathcal{T}$ of $\mathcal{T}$ such that all the vertices adjacent to $w_\mathcal{T}$ but possibly one are leaves of $\mathcal{T}$ and none of these leaves is a vertex of $\mathcal{T}'$. Since $\mathcal{T}$ is complete it contains every of the oriented edges of $T$ with terminal vertex $w_\mathcal{T}$. Furthermore, for one of these oriented edges say $e$, the half-tree $Te\cup \{w_\mathcal{T}\}$ contains $\mathcal{T}'$. Let $\mathcal{R}=\big( Te\cup \{w_\mathcal{T}\}\big)\cap \mathcal{T}$. Since $\mathcal{R}$ is complete, contains an edge and since $G$ satisfies the Tits independence property, notice that $$\Fix_G(\mathcal{R})=\big\lb\Fix_{G}(Te)\cap  \Fix_G(\mathcal{R})\big\rb\prod_{f\in \partial E_o(\mathcal{R})-\{e\}} \big\lb\Fix_{G}(Tf)\cap \Fix_G(\mathcal{R})\big\rb.$$ On the other hand, $\forall f\in \partial E_o(\mathcal{R})-\{e\}$ we have that $\mathcal{T}\subseteq Tf$ and therefore that $\Fix_{G}(Tf)\cap \Fix_G(\mathcal{R})\subseteq \Fix_G(\mathcal{T})$. Furthermore, since $\mathcal{T}'\subseteq Te\cup \{w_\mathcal{T}\}$, we have that $\Fix_{G}(Te)\cap \Fix_G(\mathcal{R})\subseteq \Fix_G(\mathcal{T}')$. The desired inclusion follows. 
	\end{proof}
	It is now time to study the consequences of such a factorisation on the cuspidal representations of $G$. Let $T$ be a thick semi-regular tree, let $G\leq \Aut(T)$ be a closed subgroup satisfying the Tits independence property and let $\pi$ be a cuspidal representations of $G$.
	\begin{definition}
		 A complete finite subtree $\mathcal{T}$ of $T$ is called a \tg{minimal subtree} of $\pi$ if $\Hr{\pi}^{\Fix_G(\mathcal{T})}$ is not zero and no complete finite subtree of $T$ with this properties has strictly fewer interior vertices than $\mathcal{T}$.
	\end{definition} 
	 In light of the definition of cuspidal representations of $G$ we recall the existence of a complete finite subtree $\mathcal{T}$ such that $\pi$ admits a non-zero $\Fix_G(\mathcal{T})$-invariant vector. This ensures the existence of a minimal subtree for any cuspidal representation of $G$.
	\begin{remark}
		A cuspidal representation $\pi$ of $G$ is in general associated with many different minimal subtrees. For instance, if $\mathcal{T}$ is a minimal subtree of $\pi$, then the subtree $g\mathcal{T}$ is also a minimal subtree of $\pi$ for all $g\in G$. This follows from the observation that $$\Hr{\pi}^{\Fix_G(g\mathcal{T})}=\Hr{\pi}^{g\Fix_G(\mathcal{T})g^{-1}}=\pi(g)\Hr{\pi}^{\Fix_G(\mathcal{T})}\q \forall g\in G.$$
	\end{remark}
	\begin{lemma}
		Let $G\leq \Aut(T)$ be a closed subgroup satisfying the Tits independence property, $\pi$ be a cuspidal representation of $G$, $\mathcal{T}$ be a minimal subtree of $\pi$ and let $\xi\in \Hr{\pi}^{\Fix_G(\mathcal{T})}$ be a unit vector. Then, the matrix coefficient function 
		$$\varphi_{\xi,\xi}:G\rightarrow \C: g\mapsto \prods{\pi(g)\xi}{\xi}_{\Hr{\pi}}$$
		is compactly supported on 
		$$\Stab_G(\mathcal{T})=\{g\in G: g\mathcal{T}=\mathcal{T}\}.$$
	\end{lemma}
	\begin{proof}
		Since $\xi\in \Hr{\pi}^{\Fix_G(\mathcal{T})}$, notice that $\varphi_{\xi,\xi}$ is a $\Fix_G(\mathcal{T})$-bi-invariant functions. Now, let $g\in G-\Stab_G(\mathcal{T})$. In particular, $g^{-1}\mathcal{T}\not= \mathcal{T}$ and Lemma \ref{independence figa nebbia 3.1} ensures the existence of a complete proper subtree $\mathcal{R}$ of $\mathcal{T}$ such that 
		$$\Fix_{G}(\mathcal{R})\subseteq \Fix_{G}(g^{-1}\mathcal{T}) \Fix_{G}(\mathcal{T})=g^{-1}\Fix_{G}(\mathcal{T})g\Fix_{G}(\mathcal{T}).$$
		In particular, one has that 
		$$g\Fix_G(\mathcal{R})\subseteq \Fix_{G}(\mathcal{T})g\Fix_{G}(\mathcal{T}).$$
		It follows from the $\Fix_G(\mathcal{T})$-bi-invariance of $\varphi_{\xi,\xi}$ that 
		\begin{equation*}
		\begin{split}
		\mu(\Fix_G(\mathcal{R}))\varphi_{\xi,\xi}(g)=\int_{\Fix_G(\mathcal{R})}\varphi_{\xi,\xi}(gh)\diff \mu(h)=\prods{\pi(g)\int_{\Fix_G(\mathcal{R})}\pi(h)\xi\diff \mu(h)}{\xi}_{\Hr{\pi}}.
		\end{split}
		\end{equation*}
		However, the right-hand side of this equality is null by minimality of $\mathcal{T}$ since $\int_{\Fix_G(\mathcal{R})}\pi(h)\xi \diff \mu (h)$ defines a $\Fix_G(\mathcal{R})$-invariant vector of $\pi$. Furthermore, on has that $\mu(\Fix_G(\mathcal{R}))>0$ because $\Fix_G(\mathcal{R})$ is an open subset of $G$. It follows that $\varphi_{\xi,\xi}(g)=0$ for all $g\in G - \Stab_G(\mathcal{T})$.
	\end{proof}
	A great deal of information can be deduced from this result. Indeed, let $G\leq \Aut(T)$ be a closed subgroup satisfying the Tits independence property and let $\pi$ be a cuspidal representation of $G$. The above result ensures the existence of a non-zero function of positive type associated to $\pi$ and compactly supported on $\Stab_G(\mathcal{T})$. Since this function is square-integrable \cite[Theorem 2]{DufloMoore1976} ensures that $\pi$ is a subrepresentation of the left-regular representation of $G$ (this is a generalisation of Lemma \ref{lemma dixmier sur les reps de carre integrable} to non-unimodular groups). On the other hand, since this function is also integrable \cite[Corollary 1 pg.223]{DufloMoore1976} ensures that $\{\pi\}$ is open in the unitary dual $\widehat{G}$ for the Fell topology. In particular, the set of cuspidal representations of $G$ is a countable discrete open subset of $\widehat{G}$. Finally, since this function of positive type is supported on $\Stab_G(\mathcal{T})$, Lemma \ref{lemma determine wether a representation is induced from an open subgroup of G} ensures that $\pi$ is induced from an irreducible representation $\sigma$ of $\Stab_G(\mathcal{T})$. On the other hand, as $\Stab_G(\mathcal{T})$ is a profinite group, the Example \ref{example representations of profinite tdlc groups} shows that $\sigma$ is the inflation of an irreducible representation $\omega$ of a finite quotient  $\Stab_{G}(\mathcal{T})/\Ker(\sigma)$ of $\Stab_{G}(\mathcal{T})$. Furthermore, as we will show in Chapter \ref{Chapter Olshanskii's factor}, $\Ker{\sigma}$ can be identified with the pointwise stabiliser $$\Fix_G(\mathcal{T})=\{g\in G : g v= v\qq \forall v\in V(\mathcal{T})\}.$$  This procedure leads to a bijective correspondence between the equivalence classes of cuspidal representations $\pi$ of $G$ with minimal subtree $\mathcal{\mathcal{T}}$ and a family of equivalence classes of irreducible representations of the finite quotient $$\Aut_G(\mathcal{T})=\Stab_G(\mathcal{T})/\Fix_G(\mathcal{T}).$$ 
	To be more precise, we make the following definition.
	\begin{definition}
		A representation $\omega$ of $\Aut_G(\mathcal{T})$ is called \tg{standard} if it is irreducible and if, for any proper complete finite subtree $\mathcal{R}$ of $\mathcal{T}$, $\omega$ does not admit any non-zero $\Fix_{\Aut_G(\mathcal{T})}(\mathcal{R})$-invariant vector. 
	\end{definition}
	The bijective correspondence is then described by the following theorem.
	\begin{theorem}[\cite{FigaNebbia1991} and \cite{Amann2003}]
		Let $T$ be a thick semi-regular tree, $G\leq \Aut(T)$ be a closed unimodular subgroup satisfying the Tits-independence property and $\pi$ be a cuspidal representation of $G$ with minimal subtree $\mathcal{T}$. Then, $\pi$ admits a non-zero matrix coefficient compactly supported on $\Stab_{G}(\mathcal{T})$ and the set of minimal subtrees of $\pi$ is given by $\{g\mathcal{T}: g\in G\}$. On the other hand, if $\omega$ is a standard representation of $\Aut_{G}(\mathcal{T})$, the representation $\Ind_{\Stab_G(\mathcal{T})}^{G}(\omega \circ p_{\mathcal{T}})$ is a cuspidal representation of $G$ with minimal subtree $\mathcal{T}$ where $$p_\mathcal{T}: \Stab_{G}(\mathcal{T})\rightarrow \Aut_G(\mathcal{T})$$ denotes the natural projection map. Furthermore, every cuspidal representation of $G$ with minimal subtree $\mathcal{T}$ arises from this procedure and we have that $$\Ind_{\Stab_G(\mathcal{T})}^{G}(\omega_1 \circ p_{\mathcal{T}})\simeq \Ind_{\Stab_G(\mathcal{T})}^{G}(\omega_2 \circ p_{\mathcal{T}})$$ if and only if $\omega_1\simeq\omega_2$ when $\omega_1$ and $\omega_2$ are standard representations. 
	\end{theorem}

	\subsection{Induction and restriction dynamic for spherical and special representations}\label{Section induction dynamic}
	Let $T$ be a thick regular tree. Consider $T$ as a semi-regular tree with associated bipartition $V=V_0\sqcup V_1$. An automorphism $g\in \Aut(T)$ is said to be type preserving if $gv\in V_0$ and $gw\in V_1$ for all $v\in V_0, \qq w\in V_1$. We denote by $\Aut(T)^+$ be the subgroup of type-preserving automorphisms of $T$. Notice that $\Aut(T)^+$ is a closed index $2$ subgroup of $\Aut(T)$ and hence defines an open subgroup of $\Aut(T)$. The purpose of the present section is to describe the restriction-induction dynamic of spherical and special representations of any closed subgroup $G\leq \Aut(T)$ acting $2$-transitively on the boundary $\partial T$ and transitively on the vertices of $T$ and its index $2$ subgroup of type preserving automorphisms $G^+=G\cap \Aut(T)^+$. Notice that $G^+$ is still $2$-transitive on the boundary and hence that the classification of spherical and special representations of $G$ and $G^+$ are both handled by section \ref{Section special rep} and Section \ref{Section spherical rep}. We adopt here the same formalism as in those sections. This dynamic of restriction and induction is described by the following result. 
	\begin{theoremletter}\label{theoremletter induction restiction dynamic}
		Let $T$ be a $d$-regular tree, let $G\leq \Aut(T)$ be a closed subgroup acting $2$-transitively on the boundary $\partial T$ and transitively on the vertices and let $G^+=G\cap \Aut(T)^+$. Then, the following hold:
		\begin{enumerate}[leftmargin=*, label=(\roman*)]
			\item\label{item 1 theoremletter induction restiction dynamic} For all $\alpha\in \mathcal{I}_{v,G}-\{0\}$, we have $$\Res_{G^+}^G(\phi_{v,G}^{-1}(\alpha))\simeq \phi_{v,G^+}^{-1}\Big(\frac{d}{d-1}\alpha^2-\frac{1}{d-1}\Big).$$
			In particular, each non-exceptional spherical representation of $G^+$ is the restriction of exactly two non-equivalent spherical representations of $G$.
			\item\label{item 2 theoremletter induction restiction dynamic} $\Res_{G^+}^G(\phi_{v,G}^{-1}(0)) \simeq \pi_v \oplus \pi_{v'}.$
			\item\label{item 3 theoremletter induction restiction dynamic} For all $\gamma\in \mathcal{I}_{v,G^{+}}-\{ -\frac{1}{d-1}\}$ one has $$\Ind_{G^+}^G(\phi_{v,G^{+}}^{-1}(\gamma))\simeq \phi_{v,G^{+}}^{-1}\Big(-\sqrt{\frac{d-1}{d}\gamma+\frac{1}{d}}\Big)\oplus \phi_{v,G^{+}}^{-1}\Big(\sqrt{\frac{d-1}{d}\gamma+\frac{1}{d}}\Big).$$ 
			In particular, the representation induced from any non-exceptional spherical representation of $G$ splits as a direct sum of two non-equivalent spherical representations of $G$ corresponding to opposite elements in $\mathcal{I}_{v,G}$.
			\item\label{item 4 theoremletter induction restiction dynamic} $\Ind_{G^+}^G(\pi_v) \simeq \phi_{v,G}^{-1}(0)\simeq\Ind_{G^+}^G(\pi_{v'}).$
			\item\label{item 5 theoremletter induction restiction dynamic} $\Ind_{G^+}^G(\sigma)\simeq\sigma^{+1}\oplus \sigma^{-1}\mbox{ and}\qq \Res_{G^+}^G(\sigma^{\pm 1})\simeq\sigma.$
		\end{enumerate} 
	\end{theoremletter}

	The proof of this result relies on Theorem \ref{les rep dun group loc compact par rapport a celle d'un de ses sous groupes}. In particular, we adopt the same formalism as in Section \ref{section induction depuis sous group indice 2}. For clarity of the exposition, the proof of the theorem is spread into several results. For the rest of this section, we let $T$ be a thick $d$-regular tree,  $G\leq \Aut(T)$ be a closed subgroup acting $2$-transitively on the boundary $\partial T$ of $T$ and on the vertices if $T$. Let $e$ be an edge of $T$ with vertices $\{v,v'\}$, let $G^+=G\cap \Aut(T)^+$ be the index $2$ subgroup of type preserving automorphisms of $G$ and let $\chiup:G\rightarrow \C$ be the corresponding unitary character. We start with the following result. 
	\begin{lemma}\label{lemma equiavelence de pi et pi twister pour spherique}
		Let $\pi$ be a spherical representation of $G$. Then, $\pi^\chi$ is a spherical representation of $G$, and $\phi_{v,G}(\pi^\chi)=-\phi_{v,G}(\pi)$. In particular, $\pi\simeq \pi^\chi$ if and only if $\pi\simeq\phi_{v,G}^{-1}(0)$.
	\end{lemma}
	\begin{proof}
		Let $\xi$ be a unit $\Fix_G(v)$-invariant vector of $\pi$ and let $\tau\in G$ be a translation of step $1$ along a geodesic containing $v$. By the definition $\pi^\chi$ has representation space $\hr{\pi}$ and since $\Fix_G(v)\subseteq G^+$, $\xi$ is also a $\Fix_G(v)$-invariant vector of $\pi^\chi$. In particular, $\pi^\chi$ is a spherical representation of $G$ with corresponding $\Fix_G(v)$-spherical function of positive type $$\varphi_{\pi^\chi}:G\rightarrow \C:g\mapsto\prods{\pi^{\chi}(g)\xi}{\xi}.$$ 
		On the other hand, since  $\tau \in G-G^+$ we notice that $$\phi_{v,G}(\pi^\chi)=\varphi_{\pi^\chi}(\tau)=\prods{\pi^\chi(\tau)\xi}{\xi} =-\prods{\pi(\tau)\xi}{\xi}=-\varphi_\pi(\tau)=-\phi_{v,G}(\pi).$$
		This proves the first part of the lemma. Furthermore, since $\phi_{v,G}$ is a bijection between the set of unitary equivalence classes of spherical representations of $G$ and $\mathcal{I}_{v,G}$ notice that $\pi\simeq \pi^\chi$ if and only if $\phi_{v,G}(\pi^\chi)=\phi_{v,G}(\pi)$ which happens if and only if $\phi_{v,G}(\pi)=0$.
	\end{proof}
	
	\begin{lemma}\label{Lemma la restriction d'une spherique est une speherique}
		Let $\alpha\in \mathcal{I}_{v,G}$ and let $\varphi_\alpha$ be the spherical function of positive type on $G$ corresponding to $\phi_{v,G}^{-1}(\alpha)$. Then, $\restr{\varphi_{\alpha}}{G^+}$ is the $\Fix_{G^+}(v)$-spherical function of positive type of $G^+$ corresponding to $\phi_{v,G^+}^{-1}(\gamma_\alpha)$ where
		\begin{equation}\label{equation gamma en fonciton de alpha}
		\gamma_\alpha= \frac{d}{d-1}\alpha^2-\frac{1}{d-1}.
		\end{equation}
		In particular, the $\Fix_{G^+}(v)$-spherical function of positive type associated to any $\phi_{v,G^+}^{-1}(\gamma)$ with $\gamma\in \mathcal{I}_{v,G^+}$ is the restriction of a $\Fix_G(v)$-spherical function of positive type of $G$.
	\end{lemma}
	\begin{proof}
		Since $\Fix_G(v)\subseteq \Aut(T)^+$, notice that $\Fix_G(v)=\Fix_{G^+}(v)$. It follows from the definition of spherical function that $\restr{\varphi_{\alpha}}{G^+}$ is a spherical function of positive type. Now, let $\tau\in G$ be a translation of step $1$ along a geodesic containing $v$ and let $\gamma_\alpha\in \mathcal{I}_{v,G^+}$ correspond to the spherical representation of $G^+$ obtained by GNS construction to the function of positive type $\restr{\varphi_{\alpha}}{G^+}.$ Since $\tau^2\in G^+$ is a translation of step $2$ along a geodesic containing $v$, notice that $\gamma_\alpha=\varphi_{\alpha}(\tau^2)$. Let $U=\{k\in \Fix_G(v): k\tau(v)=\tau^{-1}(v)\}$ and notice since $\varphi_{\alpha}$ is a $\Fix_G(v)$-spherical function of $G$ that for all integers $n\geq 2$ we have
		\begin{equation*}
		\begin{split}
		\varphi_{\alpha}(\tau^{n-1})\varphi_{\alpha}(\tau)&=\int_{\Fix_G(v)}\varphi_{{\alpha}}(\tau^nk\tau)\diff \mu(k)\\
		&=\int_{\Fix_G(v)-U}\varphi_{{\alpha}}(\tau^nk\tau)\diff \mu(k)+ \int_{U}\varphi_{{\alpha}}(\tau^nk\tau)\diff \mu(k)\\
		&=\frac{d-1}{d}\varphi_{\alpha}(\tau^n)+ \frac{1}{d}\varphi_{\alpha}(\tau^{n-2}).
		\end{split}
		\end{equation*}
		In particular, for all integer $n\geq 2$ we have
		\begin{equation*}
		\varphi_{\alpha}(\tau^n)= \frac{d}{d-1}\varphi_{\alpha}(\tau)\varphi_\pi(\tau^{n-1})-\frac{1}{d-1}\varphi_{\alpha}(\tau^{n-2}).
		\end{equation*}
		Taking $n=2$ we obtain that
		$$\gamma_\alpha=\varphi_{\alpha}(\tau^2)=\frac{d}{d-1}\varphi_{\alpha}(\tau)^2-\frac{1}{d-1}=\frac{d}{d-1}\alpha^2-\frac{1}{d-1}.$$
		Now, notice by continuity that the $\gamma_\alpha$ exhaust all values in $\mathcal{I}_{v,G^{+}}=\big\lb -\frac{1}{d-1},1\big\rb$ as $\alpha$ varies in $\mathcal{I}_{v,G}=\lb -1; 1\rb$ since $\gamma_{0}=-\frac{1}{d-1}$ and $\gamma_{1}=1.$
	\end{proof}
	The following proposition proves Theorem \ref{theoremletter induction restiction dynamic}\ref{item 1 theoremletter induction restiction dynamic} and \ref{theoremletter induction restiction dynamic}\ref{item 3 theoremletter induction restiction dynamic}
	\begin{proposition}
		For all $\alpha\in \mathcal{I}_{v,G}-\{0\}$, $\Res_{G^+}^G(\phi_{v,G}^{-1}(\alpha))$ is a spherical representation of $G^+$ with non-zero $\Fix_G(v)$-invariant vectors and $$\Res_{G^+}^G(\phi_{v,G}^{-1}(\alpha))\simeq \phi_{v,G^+}^{-1}\Big(\frac{d}{d-1}\alpha^2-\frac{1}{d-1}\Big).$$
		Furthermore, for all $\gamma\in \mathcal{I}_{v,G^{+}}-\{ -\frac{1}{d-1}\}$ one has $$\Ind_{G^+}^G(\phi_{v,G^{+}}^{-1}(\gamma))\simeq \phi_{v,G^{+}}^{-1}\Big(-\sqrt{\frac{d-1}{d}\gamma+\frac{1}{d}}\Big)\oplus \phi_{v,G^{+}}^{-1}\Big(\sqrt{\frac{d-1}{d}\gamma+\frac{1}{d}}\Big).$$ 
	\end{proposition}
	\begin{proof}
		Let $\alpha\in \mathcal{I}_{v,G}-\{0\}$ and for brevity set $\pi_\alpha=\phi_{v,G}^{-1}(\alpha)$. As a consequence of Lemma \ref{lemma equiavelence de pi et pi twister pour spherique}, Theorem \ref{les rep dun group loc compact par rapport a celle d'un de ses sous groupes} ensures that $\Res_{G^+}^G(\pi_\alpha)$ is an irreducible representation of $G^+$ and that
		$$\Ind_{G^+}^G(\Res_{G^+}^G(\pi_\alpha)) \simeq \pi_\alpha \oplus \pi_\alpha^\chi\simeq \pi_\alpha \oplus \pi_{-\alpha}.$$
		On the other hand, Lemma \ref{Lemma la restriction d'une spherique est une speherique} ensures by irreducibility that $\Res_{G^+}^G(\pi_\alpha)$  coincide with $\phi_{v,G^+}^{-1}(\gamma_\alpha)$ where $\gamma_\alpha= \frac{d}{d-1}\alpha^2-\frac{1}{d-1}.$  The result follows from \eqref{equation gamma en fonciton de alpha} since the $\gamma_\alpha$ exhaust the interval $\mathcal{I}_{v,G^+}$ as $\alpha$ varies inside $\mathcal{I}_{v,G}$ and since $\pm\sqrt{\frac{d-1}{d}\gamma+\frac{1}{d}}$ are the two roots of the polynomial $P(x)=\frac{d}{d-1}x^2-\frac{1}{d-1}-\gamma$.
	\end{proof}
	The next proposition proves Theorem \ref{theoremletter induction restiction dynamic}\ref{item 2 theoremletter induction restiction dynamic} and \ref{theoremletter induction restiction dynamic}\ref{item 4 theoremletter induction restiction dynamic}.
	\begin{proposition}
		$\Res_{G^+}^G(\phi_{v,G}^{-1}(0)) \simeq \pi_v \oplus \pi_{v'}$ and $\Ind_{G^+}^G(\pi_v) \simeq \phi_{v,G}^{-1}(0)\simeq \Ind_{G^+}^G(\pi_{v'}).$
	\end{proposition}
	\begin{proof}
		For brevity set $\pi_0=\phi_{v,G}^{-1}(0)$ and let $h\in G$ be an inversion of the edge $e$. As a consequence of Lemma \ref{lemma equiavelence de pi et pi twister pour spherique}, Theorem \ref{les rep dun group loc compact par rapport a celle d'un de ses sous groupes} ensures that $\Res_{G^+}^G(\pi_\alpha)$ splits as a direct sum $\Res_{G^+}^G(\pi_\alpha)\simeq \rho \oplus \rho^h$ of non-equivalent irreducible representations of $G^+$ such that
		$\Ind_{G^+}^G(\rho) \simeq \pi_0 \simeq\Ind_{G^+}^G(\rho^t)$. 
		On the other hand, Lemma \ref{Lemma la restriction d'une spherique est une speherique} ensures that $\pi_{v}=\phi_{v,G^+}^{-1}(-\frac{1}{d-1})$ is a subrepresentation of $\Res_{G^+}^G(\pi_\alpha)$. We are left to show that $\pi_{v}^h\simeq\pi_{v'}$. Let $\varphi_{\pi_v}$ be the $\Fix_{G^+}(v)$-spherical function of positive type of $G$ associated with $\pi_v$. Then, notice that the function $\varphi_{\pi_v}^h:{G^+}\rightarrow \C : g\mapsto \varphi_{\pi_v}(hgh^{-1})$ is a $\Fix_{G^+}(v')$-spherical function of positive type on $G^+$ as for all $g,g'\in G$ we have
		\begin{equation*}
		\begin{split}
		\int_{\Fix_{G^+}(v')}\varphi_{\pi_v}^h(gk'g')\diff \mu(k')&= \int_{\Fix_{G^+}(v')}\varphi_{\pi_v}(hgh^{-1}hk'h^{-1}hg'h^{-1})\diff \mu(k')\\
		&= \int_{\Fix_{G^+}(v)}\varphi_{\pi_v}(hgh^{-1}khg'h^{-1})\diff \mu(k)\\
		&= \varphi_{\pi_v}(hgh^{-1}) \varphi_{\pi_v}(hg'h^{-1})= \varphi_{\pi_v}^h(g) \varphi_{\pi_v}^h(g')
		\end{split}
		\end{equation*}
		Now, let $\tau$ be a translation of step $2$ along a geodesic containing both $v$ and $v'$ and notice that $h\tau h^{-1}$ is also a translation of step $2$ along a geodesic containing both $v$ and $v'$ so that $\phi_{v',G^+}^{-1}(\pi_v^h)=\varphi_{\pi_v}^h(\tau)=\varphi_{\pi_v}(h\tau h^{-1})=\varphi_{\pi_v}(\tau)=\phi_{v,G^+}^{-1}(\pi_v)=-\frac{1}{d-1}.$ The result follows from the fact that $\phi_{v',G^+}^{-1}(\frac{-1}{d-1})=\pi_{v'}.$
	\end{proof}
	Finally, the following result treats the induction and restriction dynamic of special representations.
	\begin{proposition}
		$\Ind_{G^+}^G(\sigma)\simeq\sigma^{+1}\oplus \sigma^{-1}\mbox{ and}\qq \Res_{G^+}^G(\sigma^{\epsilon})\simeq\sigma$ for all $\epsilon\in \{-1,+1\}$.
	\end{proposition}
	\begin{proof}
		We apply directly Theorem \ref{les rep dun group loc compact par rapport a celle d'un de ses sous groupes}. Let $t\in G-G^+$ and notice that $\sigma^t$ is a special representation of $G^+$. By uniqueness of the special representation, this implies that $\sigma^t\simeq\sigma$. Thus, $\Ind_{G^+}^G(\sigma)$ splits as a sum of two inequivalent representations $\pi$, $\pi^\chi$ of $G$ such that $\sigma\simeq \Res_{G^+}^G(\pi)$ and $\sigma\simeq \Res_{G^+}^G(\pi^\chi)$. Since $\Fix_G(v)=\Fix_{G^+}(v)$ and $\Fix_G(e)=\Fix_{G^+}(e)$ for every vertex $v\in V$ and every edge $e\in E$, both $\pi$ and $\pi^\chi$ must be special representations of $G$ and the result follows from the fact that there exists only two such representations $\sigma^{-1}$ and $\sigma^{+1}$.
	\end{proof}
	\newpage
	\thispagestyle{empty}
	\mbox{}
	\newpage
	\chapter{Ol'shanskii's factorisation}\label{Chapter Olshanskii's factor}

	\section{Introduction and main results}\label{Section Olshanskii}
	We recall from Section \ref{section tdlc groups aand rep} that totally disconnected locally compact groups are characterised among locally compact groups by the property of admitting a basis of neighbourhoods of the identity consisting of compact open subgroups. The existence of such a basis ensures that each representation  admits a non-zero invariant vector for at least one of these compact open subgroup. Furthermore, we recall from  Lemma \ref{lemma each irrep of G is a rep of the quotient} and from Section \ref{cuspidal representations of the group of automorphism of a tree} that additional properties of such a basis of neighbourhoods can lead to a description of the irreducible representation having particular type of invariant vectors. The purpose of this chapter is to provide an abstraction of Ol'shanskii's framework presented in Section \ref{cuspidal representations of the group of automorphism of a tree}. To be more precise, in this chapter, we develop a machinery that aims at describing the irreducible representations admitting particular invariant vectors when the group has a basis of neighbourhoods of the identity satisfying a particular kind of factorisation property. We provide direct applications of this machinery to groups of automorphisms of trees and to universal groups of right-angled buildings in Chapters \ref{Chapter application olsh facto} and \ref{Chapter Radu groups}. 
	
	Now, let $G$ be a totally disconnected locally compact group, $U\leq G$ be a closed subgroup and $\pi$ be a representation of $G$ admitting a non-zero $U$-invariant vector $\xi\in \Hr{\pi}$. Notice that for each $g\in G$, that the vector $\pi(g)\xi$ defines a non-zero $gUg^{-1}$-invariant vector of $\Hr{\pi}$. In particular, the existence of a non-zero $U$-invariant vector is an invariant property of the conjugacy class of $U$. In light of this observation, we let $\mathcal{B}$ be the set of compact open subgroups of $G$, $P(\mathcal{B})$ be the power set of $\mathcal{B}$ and 
	$$\mathcal{C}: \mathcal{B}\rightarrow P(\mathcal{B})$$
	be the map sending a compact open subgroup to its conjugacy class in $G$. Let $\mathcal{S}$ be a basis of neighbourhoods of the identity consisting of compact open subgroups of $G$ and let $\mathcal{F}_{\mathcal{S}}=\{ \mathcal{C}(U): U\in \mathcal{S}\}$. In order to properly state our factorisation property, we require a notion of relative size for the elements of $\mathcal{S}$ that is well behaved with respect to these conjugacy classes. To this end, we equip $\mathcal{F}_\mathcal{S}$ with the partial order given by the reverse inclusion of representatives. In other words, we say that $\mathcal{C}(U)\leq \mathcal{C}(V)$ if there exist $\tilde{U}\in \mathcal{C}(U)$ and $\tilde{V}\in \mathcal{C}(V)$ such that $\tilde{V}\subseteq \tilde{U}$. For a poset $(P, \leq)$ and an element $x\in P$, we let $L_x$ is the maximal length of a strictly increasing chain in $P_{\leq x}=\{y\in P: y\leq x\}$. If $L_x$ is finite, we say that $x$ has \tg{height} $L_x-1$ in $(P,\leq)$. Otherwise, we say that $x$ has infinite height in $(P,\leq)$. 
	\begin{definition}
		A basis of neighbourhoods of the identity $\mathcal{S}$ consisting of compact open subgroups of $G$ is called a \tg{generic filtration} of $G$ if the height of every element in $\mathcal{F}_{\mathcal{S}}$ is finite.
	\end{definition} 
	The following lemma ensures that such a generic filtration exists for unimodular groups.
	\begin{lemma}\label{lemma bounded basis implies generic filtration}
		Let $G$ be a unimodular totally disconnected locally compact group, $\mu$ be a Haar measure on $G$ and $\mathcal{S}$ be a basis of neighbourhoods of the identity consisting of compact open subgroups such that $\mu(U)\leq 1$ for all $U\in \mathcal{S}$. Then, $\mathcal{S}$ is a generic filtration of $G$.
	\end{lemma}
	\begin{proof}
		Let $\mathcal{C}(U)\in \mathcal{F}_\mathcal{S}$. Since $G$ is unimodular, the measure $\mu(U)$ does not depend on the choice of representative $U\in \mathcal{C}(U)$. Now, let $\mathcal{C}(U_0)\leq\mathcal{C}(U_1)\leq ...\leq \mathcal{C}(U_{n-1})\leq \mathcal{C}(U)$ be a strictly increasing chain of elements of $\mathcal{F}_{\mathcal{S}}$. Changing representatives if needed, we can suppose that $U\subseteq U_{n-1}\subseteq ... \subseteq U_1\subseteq U_0$. In particular, notice that
		$$ \lbrack U_0:U\rbrack =\lbrack U_0 : U_1\rbrack \cdots \lbrack U_{n-1}: U\rbrack\geq 2^n.$$
		On the other hand, since $U$ and $U_0$ are both compact open subgroups of $G$ observe that $$\lbrack U_0:U\rbrack= \frac{\mu(U_0)}{\mu(U)}\leq \frac{1}{\mu(U)}.$$
		This proves that $n\leq -\log_2(\mu(U))$ and therefore that the height of $\mathcal{C}(U)$ in $\mathcal{F}_\mathcal{S}$ is finite. 
	\end{proof} 
	
	Every generic filtration $\mathcal{S}$ of $G$ splits as a disjoint union $\mathcal{S}=\bigsqcup_{l\in \N}\mathcal{S}\lb l \rb$ where $\mathcal{S}\lb l\rb$ denotes the set of elements $U\in \mathcal{S}$ such that $\mathcal{C}(U)$ has height $l$ in $\mathcal{F}_{\mathcal{S}}$. The element of $\mathcal{S}\lb l \rb$ are called the elements at \tg{depth} $l$. For every representation $\pi$ of $G$ there exists a smallest non-negative integer $l_\pi\in \N$ such that $\pi$ admits a non-zero $U$-invariant vector for some $U\in \mathcal{S}\lb l_\pi \rb$. This $l_\pi$ is called the \tg{depth} of $\pi$ with respect to $\mathcal{S}$ in analogy with the similar notion of depth for representations of reductive groups over non-Archimedean fields introduced by Moy-Prasad in \cite{Moy1996}.
	\begin{example}
		Let $T$ be a thick regular tree.  We recall that a subtree $\mathcal{T}\subseteq T$ is \tg{complete} if the degree in $\mathcal{T}$ of each vertex $v\in V(\mathcal{T})$ that is not a leaf of $\mathcal{T}$ coincides with its degree in $T$. Let $\mathfrak{T}$ be the set of all complete finite subtrees of $T$ and $\mathcal{S}=\{\Fix_{\Aut(T)}(\mathcal{T}): \mathcal{T}\in \mathfrak{T}\}$ be the basis of neighbourhoods of the identity consisting of the groups
		$$\Fix_{\Aut(T)}(\mathcal{T})=\{g\in \Aut(T) : \qq gv=v \qq \forall v\in V(\mathcal{T})\}.$$
		Then, $\mathcal{S}$ is a generic filtration (this is proved in Section \ref{Application to Aut T} below). Furthermore, Lemma \ref{Lemma la startification de S pour Aut(T)} ensures that the depth $l_\pi$ with respect to $\mathcal{S}$ of a representation $\pi$ of $\Aut(T)$ can be interpreted as follows:
		\begin{itemize}[leftmargin=*]
			\item $l_\pi=0$ if and only if $\pi$ is a spherical representation of $\Aut(T)$. 
			\item $l_\pi=1$ if and only if $\pi$ is a special representation of $\Aut(T)$.
			\item $l_\pi\geq 2$ if and only if $\pi$ is a cuspidal representation of $\Aut(T)$. Furthermore, $l_\pi$ is the smallest positive integer $l$ for which there exists a complete finite subtree $\mathcal{T}$ of $T$ with $l-1$ of interior vertices such that $\pi$ has a non-zero $\Fix_{\Aut(T)}(\mathcal{T})$-invariant vector.
		\end{itemize}
	\end{example}
	\noindent We now introduce the notion of \tg{Ol'shanskii's factorisation}.  
	\begin{definition}\label{definition olsh facto}
		Let $G$ be a non-discrete unimodular totally disconnected locally compact group, $\mathcal{S}$ be a generic filtration of $G$ and $l$ be a strictly positive integer. We say that $\mathcal{S}$ \tg{factorises at depth} $l$ if for all $U\in\mathcal{S}\lb l \rb$ the following conditions are satisfied:
		\begin{enumerate}[label=(\roman*)]
			\item\label{olshfacto definition condition 1} For every $V$ in the conjugacy class of an element of $\mathcal{S}$ with $V \not\subseteq U$, there exists a subgroup $W$ in the conjugacy class of an element of $\mathcal{S}\lb l -1\rb$ satisfying that $$U\subseteq W \subseteq V U=\{vu: v\in V, u\in U\}.$$
			\item\label{olshfacto definition condition 2} For every $V$ in the conjugacy class of an element of $\mathcal{S}$, the set
			\begin{equation*}
			N_G(U, V)= \{g\in G : g^{-1}Vg\subseteq U\}
			\end{equation*}
			is compact.
		\end{enumerate} 
		Furthermore, the generic filtration $\mathcal{S}$ of $G$ is said to \tg{factorise$^+$ at depth} $l$ if, in addition, for all $U\in\mathcal{S}\lb l\rb$ and every $W$ in the conjugacy class of an element of $\mathcal{S}\lb l-1\rb$ such that $U\subseteq W$ we have that
		\begin{equation*}
		W\subseteq N_G(U,U) =\{g\in G : g^{-1}Ug\subseteq U\}.
		\end{equation*}  
	\end{definition} 
	\begin{remark}
		The factorisation at depth $l$ defined here depends on the entire generic filtration $\mathcal{S}$ and not only on the elements of $\mathcal{S}\lb l\rb$ and $\mathcal{S}\lb l-1\rb$.
	\end{remark}
	\begin{remark}
		Since $G$ is unimodular, notice that the set $N_G(U,U)$ coincides with the normaliser $N_G(U)$ of $U$ in $G$. Furthermore, notice that the conditions \ref{olshfacto definition condition 1} and \ref{olshfacto definition condition 2} are satisfied for some $U\in \mathcal{S}\lb l\rb$ if and only if they are satisfied for each of its conjugates. In particular, a generic filtration $\mathcal{S}$ factorises at depth $l$ if and only if the conditions \ref{olshfacto definition condition 1} and  \ref{olshfacto definition condition 2} are satisfied for each subgroup $U$ that is conjugate to an element of $\mathcal{S}\lb l \rb$. The same remark holds for the notion of factorisation$^+$.
	\end{remark}
	The purpose of this chapter is to prove the following theorem.
	\begin{theoremletter}\label{la version paki du theorem de classification}
	Let $G$ be a non-discrete unimodular totally disconnected locally compact group and $\mathcal{S}$ be a generic filtration of $G$ factorising at depth $l$. Then, every irreducible representation $\pi$ of $G$ at depth $l$ satisfies the following:
	\begin{enumerate}[label=(\roman*)]
		\item\label{item 1 la version paki du theorem de classif avec seed} There exists a unique $C_\pi\in \mathcal{F}_\mathcal{S}=\{\mathcal{C}(U): U\in \mathcal{S}\}$ with height $l$ such that for all $U\in C_\pi$, $\pi$ admits a non-zero $U$-invariant vector.
		\item\label{item 2 la version paki du theorem de classif avec NGU} For each $U\in C_\pi$, there exists a function of positive type associated to $\pi$ that is supported in the compact open subgroup $N_G(U)$. In particular, $\pi$ is induced from an irreducible representation of $N_G(U)$, belongs to the discrete series of $G$ and its equivalence class is isolated in the unitary dual $\widehat{G}$ for the Fell topology.
	\end{enumerate}
	Furthermore, if $\mathcal{S}$ factorises$^+$ at depth $l$, there is a bijective correspondence between the equivalence classes of irreducible representations $\pi$ of $G$ at depth $l$ with $C_\pi=\mathcal{C}(U)$ and the equivalence classes of $\mathcal{S}$-standard representations of the finite group $N_G(U)/U$. This correspondence is explicitly given by Theorem \ref{the theorem of classification for cuspidal representations}.
	\end{theoremletter}

	\noindent We refer to Section \ref{Section the bijective correspondence of theorem A} for a proper definition of $\mathcal{S}$-standard representations and for the details of the correspondence given by Theorem \ref{the theorem of classification for cuspidal representations}. 
	\begin{remark}
		Theorem \ref{la version paki du theorem de classification} does not ensure the existence of irreducible representations at depth $l$. 
	\end{remark}
	\begin{remark}
		Different generic filtrations might factorise simultaneously and lead to a different description of the same representations. A concrete example of this phenomenon is given on the full group of automorphisms of a $(d_0,d_1)$-regular tree in Sections \ref{Application to Aut T}, \ref{application IPk} and \ref{application IPV1} if $d_0\not=d_1$. Furthermore, different generic filtrations might also describe different sets of irreducible representations. A concrete example of this phenomenon occurs for instance if we replace a generic filtration $\mathcal{S}$ of $G$ that factorises at all positive depths with the generic filtration $\mathcal{S}'=\mathcal{S}-\mathcal{S}\lb 0 \rb$. In that case, $\mathcal{S}'\lb l'\rb =\mathcal{S}\lb l+1\rb$ for all $ l\in \N$. In particular, when applied to $\mathcal{S}'$, Theorem \ref{la version paki du theorem de classification} does not describe the irreducible representations admitting a non-zero $U$-invariant vector if $U\in \mathcal{S}\lb 1\rb$ while it does when applied to $\mathcal{S}$.
	\end{remark} 
	\section{Proof of Theorem \ref{la version paki du theorem de classification}}\label{section generalisation of Ol'shanskii machinery}
	\subsection{Direct consequences of Ol'shanskii's factorisation}
	Let $G$ be a non-discrete unimodular totally disconnected locally compact group, $\mu$ be a Haar measure on $G$ and $\mathcal{S}$ be a generic filtration of $G$. This section explores the first consequences of a factorisation of $\mathcal{S}$ at depth $l$. We begin with the following key lemma.
	
	\begin{lemma}\label{les fonction continue de LSS sont a support compact}
		Suppose that $\mathcal{S}$ factorises at depth $l$. Let $U$ be conjugate to an element of $\mathcal{S}\lb l \rb$ and $V\leq G$ be conjugate to an element of $\mathcal{S}$. Suppose that $\varphi: G\rightarrow\C$ is a $U$-right-invariant, $V$-left-invariant function satisfying
		\begin{equation*}
		\int_{W} \varphi(gh)\qq \diff\mu (h)=0 \q \forall g\in G
		\end{equation*}
		for every $W$ that is conjugate to an element of $\mathcal{S}\lb l-1\rb$ and such that $U\subseteq W$. Then, $\varphi$ is compactly supported and we have that
		\begin{equation*}
		{\rm supp}{(\varphi)}\subseteq N_G(U, V )=\{g\in G: g^{-1}Vg \subseteq U\}.
		\end{equation*}
	\end{lemma}
	\begin{remark}
		Since $U$ is a compact open subgroup and since $\varphi$ is $U$-right-invariant, notice that $\varphi$ is continuous. In particular, the integrals $\int_{W} \varphi(gh)\qq \diff\mu (h)$ are all well-defined.
	\end{remark}
	\begin{proof} Since $\mathcal{S}$ factorises at depth $l$, $N_G(U, V)$ is a compact set. Let $g\not\in N_G(U, V )$ and notice that $g^{-1}Vg\not\subseteq U$. In particular, there exists $W$ in the conjugacy class of an element of $\mathcal{S}\lb l-1\rb $ such that $U\subseteq W\subseteq g^{-1}VgU$. Hence, $gW\subseteq VgU$ and we have by $U$-right-invariance and $V$-left-invariance that $\varphi(gh)= \varphi(g)$ for all $h\in W$. It follows that
		\begin{equation*}
		\varphi(g)=\frac{1}{\mu(W)}\int_{W}\varphi(gh)\qq \diff\mu(h)= 0,
		\end{equation*} 
		which proves as desired that	$${\rm supp}{(\varphi)}\subseteq N_G(U, V )=\{g\in G: g^{-1}Vg \subseteq U\}.$$
	\end{proof}
	\noindent The following lemma follows directly and proves Theorem \ref{la version paki du theorem de classification}\ref{item 1 la version paki du theorem de classif avec seed}.
	\begin{lemma}\label{Lemma existence of a seed}
		Let $\pi$ be an irreducible representation of $G$ at depth $l$ and suppose that $\mathcal{S}$ factorises at depth $l$. Then, there exists a unique $C_\pi\in \mathcal{F}_{\mathcal{S}}$ with height $l$ such that for all $U\in C_\pi$, $\pi$ has a non-zero $U$-invariant vector.   
	\end{lemma}
	\begin{proof}
		Since $\pi$ is at depth $l$, there exists a compact open subgroup $U\in \mathcal{S}\lb l \rb$ at depth $l$ and a non-zero vector $\xi \in \Hr{\pi}^U$. Now, let $V\leq G$ be conjugate to an element of $\mathcal{S}\lb l \rb$ and admitting a non-zero invariant vector $\xi'\in \Hr{\pi}$ and let us show that $\mathcal{C}(U)=\mathcal{C}(V)$. The function $$\varphi_{\xi, \xi'} :G \rightarrow \C : g\mapsto \prods{\pi(g)\xi}{\xi'}$$ is clearly $U$-right-invariant and $V$-left-invariant. On the other hand, since $\pi$ does not admit a non-zero $W$-invariant vector for any subgroup $W\leq G$ that is conjugate to an element of $\mathcal{S}\lb l-1\rb$ we have that
		\begin{equation*}
		\int_{W} \varphi_{\xi, \xi'}(gh)\qq \diff\mu (h)=\prods{\int_{W} \pi(g)\pi(h)\xi \diff \mu(h)}{\xi'} =0 \q \forall g \in G
		\end{equation*}
		for each such $W$.
		It follows from Lemma \ref{les fonction continue de LSS sont a support compact} that $\varphi_{\xi, \xi'}$ is supported on $N_G(U,V )=\{g\in G:g^{-1}Vg \subseteq U\}$. On the other hand, since $\pi$ is irreducible, $\xi$ is cyclic and the function $\varphi_{\xi, \xi'}$ is not identically zero. This implies the existence of an element $g\in G$ such that $g^{-1}Vg \subseteq U$. Considering now the function $\varphi_{\xi', \xi} :G \rightarrow \C : g\mapsto \prods{\pi(g)\xi'}{\xi}$ we obtain by symmetry the existence of an element $h\in G$ such that $h^{-1}Uh \subseteq V$. In particular, $h^{-1}Uh\subseteq V \subseteq g U g^{-1}$. Hence, $g^{-1}h^{-1}Uhg \subseteq U$. Since $U$ is a compact open subgroup of $G$ and since $G$ is unimodular, this implies that $g^{-1}h^{-1}Uhg=U$. In particular, we have that $h^{-1}Uh= V$ which implies that $\mathcal{C}(U)=\mathcal{C}(V)$. 
	\end{proof}
	\begin{definition}\label{definition seed}
		The unique element $C_\pi\in \mathcal{F}_\mathcal{S}$ with height $l$ such that for all $U\in C_\pi$, $\pi$ admits a non-zero $U$-invariant vector is called the \tg{seed} of $\pi$.  
	\end{definition}
	The following proposition proves Theorem \ref{la version paki du theorem de classification}\ref{item 2 la version paki du theorem de classif avec NGU}.
	\begin{proposition}
		Let $\pi$ be an irreducible representation of $G$ at depth $l$ and suppose that $\mathcal{S}$ factorises at depth $l$. Then, for every $U\in C_\pi$, there exists a function of positive type associated with $\pi$ that is supported on the compact open subgroup $N_G(U)$. In particular, $\pi$ is induced from an irreducible representation of $N_G(U)$, belongs to the discrete series of $G$ and its equivalence class is isolated in the unitary dual $\widehat{G}$ for the Fell topology. 
	\end{proposition}
	\begin{proof}
		Since $\pi$ is at depth $l$, there exists a compact open subgroup $U\in \mathcal{S}\lb l \rb$ and a non-zero vector $\xi \in \Hr{\pi}^U$. Now, Lemma \ref{Lemma existence of a seed} ensures that $$\varphi_{\xi,\xi}: G\rightarrow \C: g\mapsto \prods{\pi(g)\xi}{\xi}$$ is supported on $N_G(U,U)$. Furthermore, $N_G(U,U)=N_G(U)$ is a compact open subgroup of $G$. 
		Now, Lemma \ref{lemma determine wether a representation is induced from an open subgroup of G} ensures the existence of an irreducible representation $\sigma$ of $N_G(U)$ such that $\pi\simeq \Ind_{N_G(U)}^G(\sigma)$. Furthermore, since $\varphi_{\xi,\xi}$ is compactly supported, the representation $\pi$ is both square-integrable and integrable. Hence, $\pi$ belongs to the discrete series of $G$ and \cite[Corollary 1 pg.223 ]{DufloMoore1976} ensures that its equivalence class is open in the unitary dual $\widehat{G}$ for the Fell topology.
	\end{proof}
	\begin{remark}
		Notice that the irreducible representation $\sigma$ of $N_G(U)$ such that $\pi\simeq \Ind_{N_G(U)}^G(\sigma)$ is the inflation of an irreducible representation of a finite quotient of $N_G(U)$. Indeed, as $\sigma$ is an irreducible representation of a compact group it is finite dimensional. In particular, $\sigma(N_G(U))$ is a closed subgroup of the Lie group $\mathcal{U}(d)$ of unitary operators of the $d$ dimensional complex Hilbert space for some positive integer $d\in \N$. On the other hand, $\sigma(N_G(U))$ is a quotient of a totally disconnected compact group. In particular, $\sigma(N_G(U))$ is a totally disconnected compact Lie group and is therefore finite. This implies that ${\rm Ker}(\sigma)$ is an open subgroup of finite index of $N_G(U)$ and therefore that $\sigma$ is the lifted to $N_G(U)$ from an irreducible representation of the finite group $N_G(U)/{\rm Ker}(\sigma)$. The purpose of the rest of this chapter is identify the irreducible representations of $N_G(U)$ at depth $l$ that arises in this manner and to show that every irreducible representations of $G$ at depth $l$ can be obtained from this procedure if $\mathcal{S}$ factorises$^+$. 
	\end{remark}

	\subsection{A family of square-integrable representations}
	For the rest of this section, let $G$ be a non-discrete unimodular totally disconnected locally compact group, $\mathcal{S}$ be a generic filtration of $G$ factorising at depth $l$, $U\leq G$ be conjugate to an element of $\mathcal{S}\lb l\rb$, $\mu$ be a Haar measure of $G$ and let $\lambda_G$ and $\rho_G$ be respectively the left-regular and right-regular representations of $G$. The purpose of this section is to study the space of functions of positive type associated with irreducible representations of $G$ with seed $\mathcal{C}(U)$. Lemma \ref{les fonction continue de LSS sont a support compact} motivates the following definition.
	\begin{definition}\label{definition of LSS}
		We define $\mathcal{L}_{\mathcal{S}}(U)$ to be the closure in $L^2(G)$ of the set of all functions $\varphi:G\rightarrow \C$ satisfying the following properties:
		\begin{enumerate}\label{les 3 prop de LSHS}
			\item $\varphi$ is $U$-right-invariant.
			\item  $\varphi$ is $V$-left-invariant for some $V\leq G$ conjugate to an element of $\mathcal{S}$.
			\item For every $W\leq G$ containing $U$ and conjugate to an element of $\mathcal{S}\lb l-1\rb$ we have that
			\begin{equation*}
			\int_{W}\qq \varphi(gh)\qq \diff\mu(h)\qq=\qq 0\q \forall g\in G.
			\end{equation*}
		\end{enumerate}
		Equivalently these three properties can be formulated in terms of fixed point subspace and orthogonal complement as follows:
		\begin{enumerate}
			\item $\varphi\in L^2(G)^{\rho_G(U)}.$
			\item $\varphi\in \bigcup_{W\in\mathcal{S},V\in \mathcal{C}(W)}L^2(G)^{\lambda_G(V)}$.
			\item $\varphi\in \bigcap_{V\in \mathcal{S}\lb l-1\rb, U\subseteq W\in \mathcal{C}(V),g\in G} \big(\mathds{1}_{gW}\big)^\perp$.
		\end{enumerate}
	\end{definition}
	By definition $\mathcal{L}_{\mathcal{S}}(U)$ is a set of equivalence classes of functions up to negligible sets and not a set of functions. However, the following lemma ensures the existence of a canonical choice of representative. 
	\begin{lemma}\label{representatiove of elements in LSS}
		For every $\tilde{\varphi}\in \mathcal{L}_{\mathcal{S}}(U)$, there exists a unique $U$-right-invariant representative $\varphi$ of $\tilde{\varphi}$. Furthermore, for every $W\leq G$ containing $U$ and conjugate to an element of $\mathcal{S}\lb l\rb$  we have that
		\begin{equation*}
		\int_{W}\qq \varphi(gh)\qq \diff\mu(h)\qq=\qq 0\q \forall g\in G.
		\end{equation*}
	\end{lemma}
	\begin{proof}
		By the definition of $\mathcal{L}_{\mathcal{S}}(U)$, there exists a sequence $(\varphi_n)_{n\in \N}$ of complex valued functions such that: 
		\begin{enumerate}
			\item $\tilde{\varphi}_n\li{n}{\infty} \tilde{\varphi}$ in $L^2(G)$ where $\tilde{\varphi}_n$ denote the equivalence class of $\varphi_n$.
			\item $\varphi_n$ is $U$-right-invariant.
			\item For every $W\leq G$ containing $U$ and conjugate to an element of $\mathcal{S}\lb l-1\rb$ we have that
			\begin{equation*}
			\int_{W}\qq \varphi_n(gh)\qq \diff\mu(h)\qq=\qq 0 \q \forall g \in G.
			\end{equation*}
		\end{enumerate}
		Now, let $\varphi':G\rightarrow \C$ be any representative of $\tilde{\varphi}$. The above implies that
		\begin{equation*}
		\int_G \modu{\varphi_n (h)-\varphi'(h)}^2\qq\diff \mu (h)\li{i}{\infty}0.
		\end{equation*}
		On the other hand, since $G$ is a disjoint union of its $U$-left-cosets, we have 
		\begin{equation*}
		\int_G \modu{\varphi_n (h)-\varphi'(h)}^2\qq\diff \mu (h)= \s{gU\in G/U}{}\int_{gU}  \modu{\varphi_n (h)-\varphi'(h)}^2\qq\diff \mu (h).
		\end{equation*}
		Therefore, we obtain for all $g\in G$ that
		\begin{equation*}
		\int_{gU}  \modu{\varphi_n (h)-\varphi'(h)}^2\qq\diff \mu (h)=\int_{U}  \modu{\varphi_n (gh)-\varphi'(gh)}^2\qq\diff \mu (h)\li{n}{\infty}0.
		\end{equation*}
		In particular, for every $g\in G$ this implies that $\varphi_n (gh)$ converges to $\varphi'(gh)$ for almost all $h\in U$. Since the $\varphi_n$ are constant on $U$-left-cosets, there exists of a unique representative $\varphi:G\rightarrow \C$ of $\tilde{\varphi}$ such that $$\varphi(gh)=\varphi(g) \qq\forall g\in G, \forall h\in U.$$ On the other hand, since $U$ is a compact set, the convergence $\varphi_n\li{n}{\infty}\varphi$ is uniform on $U$-left-cosets. Now, let $g\in G$ and let $W$ be conjugate to an element of $\mathcal{S}\lb l-1\rb$ and such that $U\subseteq W$. Since $gW$ is compact, it can be covered by finitely many $U$-left-cosets and the convergence $\varphi_n\li{n}{\infty}\varphi$ is also uniform on $W$-left-cosets. This implies as desired that
		\begin{equation*}
		\begin{split}
		\bigg\lvert\int_{W}\varphi(gh)\diff \mu(h)\bigg\lvert &=\bigg\lvert \int_{W}\varphi(gh)-\varphi_n(gh)\diff \mu(h)\bigg\lvert\leq \int_{W}\modu{\varphi(gh)-\varphi_n(gh)}\diff \mu(h)\\
		&\leq \mu(W)\supp{k\in gW}\modu{\varphi(k)-\varphi_n(k)}\qq\li{n}{\infty}\qq 0\q \forall g \in G.
		\end{split}
		\end{equation*}	
	\end{proof}
	In light of this result, we identify each equivalence class $\tilde{\varphi}\in \mathcal{L}_{\mathcal{S}}(U)$ with its canonical continuous representative if it leads to no confusion. The following lemma shows that $\mathcal{L}_{\mathcal{S}}(U)$ is $G$-left-invariant and therefore defines a subrepresentation of the left-regular representation $(\lambda_G, L^2(G))$. 
	\begin{lemma}\label{LGSHS est G left invaraitn subspacede L2}
		$\mathcal{L}_{\mathcal{S}}(U)$ is a closed $G$-left-invariant subspace of $L^2(G)$.
	\end{lemma}
	\begin{proof}
		By the definition of $\mathcal{L}_{\mathcal{S}}(U)$, it is enough to prove that the subspace of functions satisfying the three properties of Definition \ref{definition of LSS} is $G$-left-invariant. Let $\varphi$ be such a function, $k\in G$ and notice that:
		\begin{enumerate}[leftmargin=*]
			\item $\lambda_G(k)\varphi\in L^2(G)^{\rho_G(U)}$.
			\item $\lambda_G(k)\varphi\in \bigcup_{W\in\mathcal{S},V\in \mathcal{C}(W)}L^2(G)^{\lambda_G(kVk^{-1})}=\bigcup_{W\in\mathcal{S},V\in \mathcal{C}(W)}L^2(G)^{\lambda_G(V)}$.
			\item $\lambda_G(k)\varphi\in \bigcap_{ V\in \mathcal{S}\lb l-1\rb, U\subseteq W\in \mathcal{C}(V),g\in G} \big(\mathds{1}_{kgW}\big)^\perp$\\
			$\mbox{ }$ $\mbox{ }$ $\mbox{ }$  $\mbox{ }$  $\mbox{ }$  $\mbox{ }$  $\mbox{ }$  $\mbox{ }$  $\mbox{ }$  $\mbox{ }$  $\mbox{ }$  $\mbox{ }$  $\mbox{ }$  $\mbox{ }$  $\mbox{ }$  $\mbox{ }$ $\mbox{ }$  $\mbox{ }$ $\mbox{ }$  $\mbox{ }$ $\mbox{ }$ $\mbox{ }$  $\mbox{ }$ $=\bigcap_{ V\in \mathcal{S}\lb l-1\rb, U\subseteq W\in \mathcal{C}(V),g\in G} \big(\mathds{1}_{gW}\big)^\perp.$
		\end{enumerate} 
	\end{proof} 
	We denote by $T_{\mathcal{S},U}$ the subrepresentation of $\lambda_G$ corresponding to $\mathcal{L}_\mathcal{S}(U)$. The following result shows that this representation depends, up to equivalence, only on the conjugacy class $\mathcal{C}(U)$ and not on his representative $U$. 
	\begin{lemma}
		Let $C\in \mathcal{F}_\mathcal{S}$ be a conjugacy class with height $l$ and let $U,U'\in C$. Then, the representations $(T_{\mathcal{S},U}, \mathcal{L}_{\mathcal{S}}(U))$ and $(T_{\mathcal{S},U'}, \mathcal{L}_{\mathcal{S}}(U'))$ are unitarily equivalent.  
	\end{lemma}
	\begin{proof}
		Since $U$ and $U'$ belongs to the same conjugacy class $C$, there exists an element $k\in G$ such that $U'= kUk^{-1}$. Now, notice that $\rho_G(k):L^2(G)\rightarrow L^2(G)$ is a unitary operator and that $\rho_G(k)\mathcal{L}_{S}(U)=\mathcal{L}_S(U')$. Indeed, for every function $\varphi: G\rightarrow \C$ such that $\varphi\in L^2(G)^{\rho_G(U)},$
		$\varphi\in \bigcup_{W\in\mathcal{S},V\in \mathcal{C}(W)}L^2(G)^{\lambda_G(V)}$
		and $\varphi\in \bigcap_{V\in \mathcal{S}\lb l-1\rb, U\subseteq W\in \mathcal{C}(V),g\in G} \big(\mathds{1}_{gW}\big)^\perp$ we have:
		\begin{enumerate}[leftmargin=*]
			\item $\rho_G(k)\varphi\in L^2(G)^{\rho_G(kUk^{-1})}=L^2(G)^{\rho_G(U')}.$
			\item $\rho_G(k)\varphi\in \bigcup_{W\in\mathcal{S},V\in \mathcal{C}(W)}L^2(G)^{\lambda_G(V)}$.
			\item $\rho_G(k)\varphi\in \bigcap_{V\in \mathcal{S}\lb l-1\rb, U\subseteq W\in \mathcal{C}(V),g\in G} \big(\mathds{1}_{gWk^{-1}}\big)^\perp$ \\
			$\mbox{ }$ $\mbox{ }$ $\mbox{ }$  $\mbox{ }$  $\mbox{ }$  $\mbox{ }$  $\mbox{ }$  $\mbox{ }$  $\mbox{ }$  $\mbox{ }$  $\mbox{ }$  $\mbox{ }$  $\mbox{ }$  $\mbox{ }$   $=\bigcap_{V\in \mathcal{S}\lb l-1\rb, U'\subseteq kWk^{-1}\in \mathcal{C}(V),gk^{-1}\in G} \big(\mathds{1}_{gk^{-1}(kWk^{-1})}\big)^\perp.$
		\end{enumerate}
		It follows by density and by continuity of $\rho_G(k)$ that $\rho_G(k)\mathcal{L}_{S}(U)\subseteq\mathcal{L}_S(U')$. Inverting the role of $U$ and $U'$, we obtain that $\rho_G(k^{-1})\mathcal{L}_{S}(U')\subseteq\mathcal{L}_S(U)$ and thus that $\rho_G(k)\mathcal{L}_S(U)=\mathcal{L}_S(U')$. Since $\lambda_G(g)\rho_G(k)= \rho_G(k)\lambda_G(g)$ for every $g\in G$, the restriction of $\rho_G(k)$ to $\mathcal{L}_{\mathcal{S}}(U)$ provides the desired intertwining operator between $T_{\mathcal{S},U}$ and $T_{\mathcal{S}, U'}$.
	\end{proof}
	Now, recall from Lemma \ref{Lemma existence of a seed} that for every irreducible representation $\pi$ of $G$ at depth $l$ with a non-zero $U$-invariant vector $\xi$, the matrix coefficient $\varphi_{\xi,\xi}$ is a non-zero $U$-bi-invariant element of $\mathcal{L}_\mathcal{S}(U)$. The following result shows that every subrepresentation of $\mathcal{L}_{\mathcal{S}}(U)$ contains such a function. 
	\begin{lemma}\label{Every invariant space in the regular contains a non trivial invariant vector}
		Every $G$-left-invariant subspace of $\mathcal{L}_{\mathcal{S}}(U)$ contains a non-zero $U$-bi-invariant function.
	\end{lemma}
	\begin{proof}
		Suppose that $M$ is a non-zero closed $G$-left-invariant subspace of $\mathcal{L}_{\mathcal{S}}(U)$, consider any non-zero class of functions $\varphi\in M$ and let $t\in G$ be such that  $\varphi(t)\not=0$. Notice that the function
		\begin{equation*}
			\psi= \int_{U} T_{\mathcal{S},U}(kt^{-1})\varphi\qq \diff  \mu(k)
		\end{equation*}
		is $U$-left-invariant, $U$-right-invariant and satisfies 
		\begin{equation*}
			\begin{split}
				\int_{W}\qq \psi(gh)\qq \diff\mu(h)\qq=\qq 0 \q \forall g\in G
			\end{split}
		\end{equation*} 
		for every $W$ that is conjugate to an element of $\mathcal{S}\lb l-1\rb$ and such that $U\subseteq W$.
		Since $M$ is a closed $G$-left-invariant subspace of $\mathcal{L}_{\mathcal{S}}(U)$, we have that $\psi\in M$. Furthermore, this function is not identically zero since the $U$-right invariance of $\varphi$ implies that
		\begin{equation*}
			\psi(1_G)= \int_{U} \varphi(tk^{-1})\diff \mu(k)= \int_{U}\varphi(t)\qq \diff \mu(k)= \mu(U) \varphi(t)\qq \not=0. 
		\end{equation*}
	\end{proof}
	The following proposition shows that the space of $U$-bi-invariant functions of $\mathcal{L}_{\mathcal{S}}(U)$ is finite-dimensional and therefore that $T_{\mathcal{S},U}$ decomposes as a finite sum of irreducible representations of $G$.
	\begin{proposition}\label{The space of invariant vectors is finite dimensional}
		The subspace of $U$-bi-invariant functions of $\mathcal{L}_{\mathcal{S}}(U)$ is finite dimensional.  
	\end{proposition}
	\begin{proof}
		Lemma \ref{les fonction continue de LSS sont a support compact} ensures that the $U$-bi-invariant functions of $\mathcal{L}_{\mathcal{S}}(U)$ are supported on the compact set $N_G(U)$. On the other hand every $U$-bi-invariant continuous function $\varphi:G\rightarrow \C$ is constant on the $U$-left-cosets and $N_G(U)$ can be covered by finitely many  such cosets. This proves that the subspace of $U$-bi-invariant functions of $\mathcal{L}_{\mathcal{S}}(U)$ is in the span of finitely many characteristic functions and is therefore finite dimensional. 
	\end{proof}
	\begin{corollary}\label{LSS is the finite sum of irreducible rep with seed S}
		$T_{\mathcal{S}, U}$ decomposes as a finite sum of irreducible square-integrable representations of $G$ with seed $\mathcal{C}(U)$. 
	\end{corollary}
	\begin{proof}
		Lemma \ref{Every invariant space in the regular contains a non trivial invariant vector} and Proposition \ref{The space of invariant vectors is finite dimensional} ensure that $T_{\mathcal{S},U}$ decomposes as a finite sum of irreducible square-integrable representations of $G$; each of these containing a non-zero $U$-invariant vector. Now, let $W\leq G$ be conjugate to an element of $\mathcal{S}\lb r\rb$ for some $r< l$. Lemma \ref{les fonction continue de LSS sont a support compact} ensures that the $W$-left-invariant functions of $\mathcal{L}_{\mathcal{S}}(U)$ are supported on $$N_G(U,W)=\{g\in G: gWg^{-1}\subseteq U\}.$$ On the other hand, since $\mathcal{C}(W)$ has height $r$, since $\mathcal{C}(U)$ has height $l$ and since $r< l$ the set $N_G(U,W)$ is empty. This proves that $\mathcal{L}_{\mathcal{S}}(U)$ does not contain any non-zero $W$-left invariant function and thus that every irreducible representation appearing in the decomposition of $T_{\mathcal{S},U} $ has seed $\mathcal{C}(U)$.
	\end{proof}
	It is a natural to ask whether the reversed implication holds. The following provides a positive answer to this question. 
	\begin{lemma}\label{irreducible representation with seed S are in the regular rep LSS}
		Every irreducible representation $\pi$ of $G$ with seed $\mathcal{C}(U)$ is equivalent to a subrepresentation of $T_{\mathcal{S},U}$.
	\end{lemma}
	\begin{proof}
		Let $\pi$ be an irreducible representation of $G$ with seed $\mathcal{C}(U)$ and let $\xi \in \Hr{\pi}^U$ be a non-zero vector. Since $\pi$ has seed $\mathcal{C}(U)$, notice with a similar reasoning as in the proof of Lemma \ref{Lemma existence of a seed} that $\widecheck{\varphi}_{\xi,\xi}$ is a non-zero $U$-bi-invariant function of $\mathcal{L}_\mathcal{S}(U)$ where $\widecheck{\varphi}_{\xi,\xi}: G\rightarrow \C : g\mapsto \prods{\pi(g^{-1})\xi}{\xi}$. Now consider the matrix coefficient $\psi(\cdot)=\prods{T_{\mathcal{S},U}(\cdot)\widecheck{\varphi}_{\xi,\xi}}{\widecheck{\varphi}_{\xi,\xi}}$ of $T_{\mathcal{S},U}$ and notice from \cite[Theorem 14.3.3]{Dixmier1977} that
		\begin{equation*}
		\begin{split}
		\psi(g)&=\int_G T_{\mathcal{S},U}(g)\widecheck{\varphi}_{\xi,\xi}(h)\overline{\widecheck{\varphi}_{\xi,\xi}(h)}\diff \mu (h)=\int_G {\varphi}_{\xi,\xi}(h^{-1}g)\overline{{\varphi}_{\xi,\xi}(h^{-1})}\diff \mu (h)\\
		&=\int_G {\varphi}_{\pi(g)\xi,\xi}(h^{-1})\overline{{\varphi}_{\xi,\xi}(h^{-1})}\diff \mu (h)=\int_G {\varphi}_{\pi(g)\xi,\xi}(h)\overline{{\varphi}_{\xi,\xi}(h)}\diff \mu (h)\\
		&=d_{\pi}^{-1} \norm{\xi}{}^2\prods{\pi(g)\xi}{\xi} = d_{\pi}^{-1} \norm{\xi}{}^2\varphi_{\xi,\xi}(g)
		\end{split}
		\end{equation*}
		where $d_\pi$ is the formal dimension of $\pi$. Renormalising $\xi$ if needed, the result follows from the uniqueness of the GNS construction see Theorem \ref{Theorem GNS construction}.
	\end{proof}
	\begin{corollary}
		Let $C\in \mathcal{F}_{\mathcal{S}}$ be a conjugacy class with height $l$. Then, there exists at most finitely many inequivalent classes of irreducible representations of $G$ with seed $C$. 
	\end{corollary}
	
	To improve the clarity of the exposition the author decided to gather the results of the previous two sections in a theorem.
	Up to this point, Theorem \ref{la version paki du theorem de classification}\ref{item 1 la version paki du theorem de classif avec seed} and \ref{item 2 la version paki du theorem de classif avec NGU} together with the following result.
	\begin{theorem}\label{theorem de correspondance entre les subre repre se LS U et les irrep avec seed CU}
		Let $G$ be a non-discrete unimodular totally disconnected locally compact group,  $\mathcal{S}$ be a generic filtration of $G$ factorising at depth $l$ and $U\leq G$ be conjugate to an element of $\mathcal{S}\lb l\rb$. Then, the following hold:
		\begin{enumerate}
			\item The square-integrable representation $(T_{\mathcal{S},U},\mathcal{L}_{\mathcal{S}}(U))$ depends, up to equivalence, only on the conjugacy class $\mathcal{C}(U)$. 
			\item Every $G$-left-invariant subspace of $\mathcal{L}_{\mathcal{S}}(U)$ contains a non-zero $U$-bi-invariant function and every such function is compactly supported on $N_G(U)$.
			\item $T_{\mathcal{S},U}$ decomposes as a finite sum of square-integrable irreducible representations of $G$ with seed $\mathcal{C}(U)$. 
			\item Every irreducible representation of $G$ with seed $\mathcal{C}(U)$ is a subrepresentation of $T_{\mathcal{S},U}$.
		\end{enumerate} 
	\end{theorem}
	In particular, if $\mathcal{S}$ factorises at depth $l$ we have a bijective correspondence between the equivalence classes of irreducible subrepresentations of $(T_{\mathcal{S},U}, \mathcal{L}_{\mathcal{S}}(U))$ and the equivalence classes of irreducible representations of $G$ with seed $\mathcal{C}(U)$. In the following sections, we introduce a family of irreducible representations of the finite group $N_G(U)/U$ which can be lifted to representations of $N_G(U)$ and show that the irreducible representations of $G$ with seed $\mathcal{C}(U)$ are induced from these when $\mathcal{S}$ factorises$^+$ at depth $l$. 
	
	\subsection{The bijective correspondence}\label{Section the bijective correspondence of theorem A}
	
	Let $G$ be a non-discrete unimodular totally disconnected locally compact group and let $\mathcal{S}$ be a generic filtration of $G$. The purposes of this section is to describe explicitly the bijective correspondence of Theorem \ref{la version paki du theorem de classification}. This requires some formalism that we now introduce.
	
	Let $U$ be conjugate to an element of $\mathcal{S}\lb l \rb$ and notice that if $\mathcal{S}$ factorises at depth $l$, $N_G(U)/U$ is a finite group. Furthermore, notice that $$N_G(g^{-1}Ug)=g^{-1}N_G(U)g \q \forall g\in G$$ which implies that $$N_G(U)/U\simeq N_G(V)/V\q \forall U,V\in C, \qq \forall C\in \mathcal{F}_\mathcal{S}.$$ This motivates the following definition.
	\begin{definition}
		For every conjugacy class $C\in \mathcal{F}_{\mathcal{S}}$, the \textbf{group of automorphisms} $\Aut_{G}(C)$ \textbf{of the seed} $C$ is the finite group $N_G(U)/U$ corresponding to any group $U\in C$. 
	\end{definition}
	\noindent For all $C\in \mathcal{F}_{\mathcal{S}}$ with height $l$ and every $U\in C$, we set
	\begin{equation*}
	\tilde{\mathfrak{H}}_{\mathcal{S}}(U)= \{W : \exists g\in G\mbox{ s.t. }gWg^{-1}\in \mathcal{S}\lb l-1 \rb \mbox{ and } U \subseteq  W \}.
	\end{equation*}
	and let $p_{U}: N_G(U)\mapsto N_G(U)/U$ denote the quotient map. If $\mathcal{S}$ factorises$^+$ at depth $l$, notice that $p_{U}(W)$ is a non-trivial (possibly not proper) subgroup of $\Aut_{G}(C)$ for every $W\in \tilde{\mathfrak{H}}_{\mathcal{S}}(U)$. Moreover, notice that $$\tilde{\mathfrak{H}}_{\mathcal{S}}(g^{-1}Ug)= g^{-1}\tilde{\mathfrak{H}}_{\mathcal{S}}(U)g \q \forall g\in G.$$ In particular, the subset of non-trivial subgroups of $\Aut_G(C)$ $$\mathfrak{H}_{\mathcal{S}}(C)=\{ p_U(W): W\in \tilde{\mathfrak{H}}_\mathcal{S}(U)\}$$ does not depend on our choice of representative $U\in C$. 
	\begin{definition}
		An irreducible representation $\omega$ of $\Aut_{G}(C)$ is a $\mathcal{S}$-\tg{standard} if it does not have any non-zero $H$-invariant vector for any $H\in \mathfrak{H}_{\mathcal{S}}(C)$.
	\end{definition}
	Our goal is to describe the irreducible representations of $G$ with seed $C$ from these $\mathcal{S}$-\tg{standard} representations. We recall that every representation $\omega$ of $N_G(U)/U$ can be lifted to a representation $\omega \circ p_U$ of $N_G(U)$ acting trivially on $U$ and with representation space $\Hr{\omega}$. Furthermore, notice that $\omega\circ p_U$ is irreducible if and only if $\omega$ is irreducible. We can now use the process of induction recalled in Section \ref{Section induced representations}. For brevity, we denote by $$T(U,\omega)=\Ind_{N_G(U)}^G(\omega\circ p_{U})$$ the resulting representation of $G$. Our purpose will be to show that $T(U,\omega)$ is an irreducible representation of $G$ with seed $\mathcal{C}(U)$ if $\omega$ is $\mathcal{S}$-standard. Conversely, if $\pi$ is an irreducible representation of $G$ with seed $C$, notice that $\Hr{\pi}^{U}$ is a non-zero $N_G(U)$-invariant subspace of $\Hr{\pi}$ for every $U\in C$. In particular, the restriction $(\restr{\pi}{N_G(U)}, \Hr{\pi}^U)$ defines a representation of $N_G(U)$ whose restriction to $U$ is trivial. This representation passes to the quotient group $N_G(U)/U$ and therefore provides a representation $\omega_\pi$ of $\Aut_{G}(C)$. The following theorem describes the bijective correspondence of Theorem \ref{la version paki du theorem de classification} using these constructions.
	\begin{theorem}\label{the theorem of classification for cuspidal representations}
		Let $G$ be a non-discrete unimodular totally disconnected locally compact group, $\mathcal{S}$ be a generic filtration of $G$ that factorises$^+$ at depth $l$ and $C\in \mathcal{F}_{\mathcal{S}}$ be a conjugacy class at height $l$. There exist a bijective correspondence between the equivalence classes of irreducible representations of $G$ with seed $C$ and the equivalence classes of $\mathcal{S}$-standard representations of $\Aut_G(C)$. More precisely, for every $U\in \mathcal{S}$ the following holds:
		\begin{enumerate}
			\item If $\pi$ is an irreducible representation of $G$ with seed $\mathcal{C}(U)$, the representation $(\omega_\pi, \Hr{\pi}^{U})$ of $\Aut_{G}(\mathcal{C}(U))$ is a $\mathcal{S}$-standard representation of $\Aut_{G}(\mathcal{C}(U))$ such that
			\begin{equation*}
			\pi \simeq T(U,\omega_\pi)=\Ind_{N_G(U)}^{G}(\omega_\pi \circ p_{U}).
			\end{equation*}
			\item If $\omega$ is a $\mathcal{S}$-standard representation of $\Aut_{G}(\mathcal{C}(U))$, the representation $T(U,\omega)$ is an irreducible representation of $G$ with seed $\mathcal{C}(U)$. 
		\end{enumerate}
		Furthermore, if $\omega_1$ and $\omega_2$ are $\mathcal{S}$-standard representations of $\Aut_{G}(\mathcal{C}(U))$, we have that $T(U,\omega_1)\simeq T(U,\omega_2)$ if and only if $\omega_1\simeq\omega_2$. In particular, the above two constructions are inverse of one another.
	\end{theorem}
	The structure of the proof is given as follows: 
	\noindent\begin{minipage}[t]{5\linewidth}%
		\centering
		\mbox{%
			\begin{minipage}[t]{\linewidth-2\fboxsep-2\fboxrule}
				\hspace*{-24em}\raisebox{-13em}{ 
					\tikzstyle{block} = [rectangle, text centered,
					minimum height=10mm,node distance=8em,  rounded corners=1mm]
					\begin{tikzpicture}
					\node[draw,text width=4.5cm] at (0,0) [block ] (n01) {Irreducible subrepresentations of $T_{\mathcal{S},U}$};
					\node[draw,text width=4.5cm] at (3.5,-2) [block ] (n02) {Representations $T(U,\omega)$};
					\node[draw,text width=4.5cm] at (-3.5,-2) [block ] (n03) {Irreducible representations with seed $\mathcal{C}(U)$};

					\path(n02.east)  edge[thick, -latex,out=30,in=0,looseness=1.5]node[anchor=center,fill=white] {\circled{2}}(n01.east) ;
					\path(n01.south) edge[ thick, -latex,out=330,in=0,looseness=1.5]node[anchor=center,fill=white] {\circled{1}}(n03.east);
					\path(n03.west)  edge[thick, -latex,out=150,in=180,looseness=1.5]node[anchor=center,fill=white] {\circled{1}}(n01.west);
					\path(n03) edge[thick, -latex,out=330,in=210,looseness=1]node[anchor=center,fill=white] {\circled{3}}(n02);
					\end{tikzpicture}  
				}
		\end{minipage}}
	\end{minipage}\hfill
	\begin{center}
		\circled{1} Theorem \ref{theorem de correspondance entre les subre repre se LS U et les irrep avec seed CU} $\mbox{ }$ \circled{2} Proposition \ref{les induite de S standard sont des sous rep de TS}  $\mbox{ }$ \circled{3} Lemma \ref{every irred rep with seed S is the induced rep}
	\end{center}
	The last part of the statement is then handled by Lemma \ref{unicitee des induite par res S standard} below.\\
	\noindent We start the proof with an intermediate result. From now on and for the rest of this section, let $U\leq G$ be conjugate to an element of $\mathcal{S}\lb l \rb$ and suppose that $\mathcal{S}$ factorises$^+$ at depth $l$.
	\begin{lemma}\label{les fonction U invarainte de LS U on inverse dans LsU}
		Let $\varphi$ be a $U$-bi-invariant function of $\mathcal{L}_{\mathcal{S}}(U)$ and let $\widecheck{\varphi}$ be the function defined by $\widecheck{\varphi}(g)=\varphi(g^{-1})$ $\forall g\in G$. Then, $\widecheck{\varphi}$ is a $U$-bi-invariant function of $\mathcal{L}_{\mathcal{S}}(U)$.
	\end{lemma}
	\begin{proof}
		The function $\widecheck{\varphi}$ is clearly $U$-bi-invariant. Furthermore, since $G$ is unimodular and since $\varphi\in L^2(G)$, notice that $\widecheck{\varphi}\in L^2(G)$. On the other hand, Lemma \ref{les fonction continue de LSS sont a support compact} ensures that $\varphi$ is supported on $N_G(U)$. The same holds for $\widecheck{\varphi}$ as $N_G(U)$ is symmetric.  Now let $W$ be conjugate to an element of $\mathcal{S}\lb l-1\rb$ and such that $U\subseteq W$. Since $\mathcal{S}$ factorises$^+$ at depth $l$, we have that $W\subseteq N_G(U)$. In particular, for every $g\not\in N_G(U)$, $gW\cap N_G(U)=\es$ and $\int_{W}\widecheck{\varphi}(gh)\diff \mu (h)=0.$ On the other hand, if $g\in N_G(U)$, notice that $U=gUg^{-1}\subseteq gWg^{-1}$ which implies that
		$$\int_W\widecheck{\varphi}(gh)\diff\mu(h)=\int_{W} \varphi(hg^{-1})\diff \mu(h)=\int_{gWg^{-1}} \varphi(g^{-1}h)\diff \mu(h)=0.$$
	\end{proof}
	
	\begin{proposition}\label{les induite de S standard sont des sous rep de TS}
		Let $\omega$ be a $\mathcal{S}$-standard representation of $\Aut_{G}(\mathcal{C}(U))$. Then, $\mathfrak{D}^{\omega \circ p_U}_{N_G(U)}=\Hr{T(U,\omega)}^U$ and $T(U,\omega)$ is equivalent to an irreducible subrepresentation of $T_{\mathcal{S},U}$.
	\end{proposition}
	\begin{proof}
		We start by showing that $T(U,\omega)$ is equivalent to a subrepresentation of $T_{\mathcal{S},U}$. Let $\xi\in \Hr{\omega}$ be non-zero and consider the function 
		$$\phi_\xi:G\rightarrow \C : g\mapsto \begin{cases}
		\omega\circ p_U(g^{-1})\xi &\mbox{ if }g\in N_G(U)\\
		\q \q \qq 0  &\mbox{ if }g\not\in N_G(U)
		\end{cases}.$$
		We have that $\phi_\xi \in \Hr{T(U,\omega)}$ and since $\omega$ is irreducible, Lemma \ref{la correspondance Hsigma D} ensures that $\phi_\xi$ is cyclic for $T(U,\omega)$. We consider the matrix coefficient $\varphi_{\phi_\xi,\phi_\xi}(\cdot)=\prods{T(U,\omega)(\cdot)\phi_\xi}{\phi_\xi}$ of $T(U,\omega)$ and let 
		\begin{equation*}
		\varphi_{\xi,\xi}:G\rightarrow \C : g\mapsto\begin{cases}
		\prods{\omega \circ p_U(g)\xi}{\xi}_{\Hr{\omega}} &\mbox{ if }g\in N_G(U)\\
		\q\q \q 0 &\mbox{ if }g\not\in N_G(U)
		\end{cases}.
		\end{equation*} A  straightforward computation shows for every $g\in G$ that
		\begin{equation}\label{equation piece 1 of GNS implication}
		\begin{split}
		\varphi_{\phi_\xi,\phi_\xi}(g)&=\prods{T(U,\omega)(g)\phi_\xi}{\phi_\xi}_{\Ind_{N_G(U)}^G(\Hr{\omega \circ p_U})}=\s{tN_G(U)\in G/N_G(U)}{}\prods{T(U,\omega)(g)\phi_\xi(t)}{\phi_\xi(t)}\\
		&=\s{tN_G(U)\in G/N_G(U)}{}\prods{\phi_\xi(g^{-1}t)}{\phi_\xi(t)} =\prods{\phi_\xi(g^{-1})}{\phi_\xi(1_G)}=\varphi_{\xi,\xi}(g).
		\end{split}
		\end{equation}
		We are going to show that $\varphi_{\xi,\xi}$ is a $U$-left-invariant function of $\mathcal{L}_{\mathcal{S}}(U)$. It is clear from the definition that $\varphi_{\xi,\xi}$ is $U$-bi-invariant, compactly supported and hence square integrable. Now, let $W$ be conjugate to an element of $\mathcal{S}\lb l-1\rb$ and such that $U\subseteq W$. Since $\mathcal{S}$ factorises$^+$ at depth $l$, notice that $W\subseteq N_G(U)$. In particular, for every $g\not\in N_G(U)$ we have that $gN_G(U)\cap N_G(U)=\es$ and thus $\int_W\varphi_{\xi,\xi}(gh)\diff\mu(h)=0$. On the other hand, for every $g\in N_G(U)$, since $\omega$ is a $\mathcal{S}$-standard representation of $\Aut_G(\mathcal{C}(U))$ we have 
		$$\int_W\varphi_{\xi,\xi}(gh)\diff\mu(h)=\prods{\int_W\omega\circ p_U(h)\xi\diff\mu(h)}{\omega\circ p_U(g^{-1})\xi}_{\Hr{\omega}}=0.$$
		This proves, as desired, that $\varphi_{\xi,\xi}$ is a $U$-left-invariant function of $\mathcal{L}_{\mathcal{S}}(U)$. In particular, Lemma \ref{les fonction U invarainte de LS U on inverse dans LsU} ensures that $\widecheck{\varphi}_{\xi,\xi}\in \mathcal{L}_\mathcal{S}(U)$. Now, consider the matrix coefficient $\psi(\cdot)=\prods{T_{\mathcal{S},U}(\cdot)\widecheck{\varphi}_{\xi,\xi}}{\widecheck{\varphi}_{\xi,\xi}}$ of $T_{\mathcal{S},U}$, let $\mu$ be the Haar measure of $G$ renormalised in such a way that $\mu(N_G(U))=1$ and notice from \cite[Theorem 14.3.3]{Dixmier1977} that
		\begin{equation}\label{equation piece 2 of GNS implication}
		\begin{split}
		\psi(g)&=\int_G T_{\mathcal{S},U}(g)\widecheck{\varphi}_{\xi,\xi}(h)\overline{\widecheck{\varphi}_{\xi,\xi}(h)}\diff \mu (h)=\int_{N_G(U)} {\varphi}_{\xi,\xi}(h^{-1}g)\overline{{\varphi}_{\xi,\xi}(h^{-1})}\diff \mu (h)\\&=\int_{N_G(U)} {\varphi}_{\xi,\xi}(hg)\overline{{\varphi}_{\xi,\xi}(h)}\diff \mu (h) =d_{\omega}^{-1} \norm{\xi}{}^2\varphi_{\xi,\xi}(g)
		\end{split}
		\end{equation}
		where $d_\omega$ is the dimension of $\omega$. In particular, renormalising $\xi$ if needed, \eqref{equation piece 1 of GNS implication}, \eqref{equation piece 2 of GNS implication} and the uniqueness of the GNS construction (Theorem \ref{Theorem GNS construction}) imply that $T(U,\omega)$ is a subrepresentation of $T_{\mathcal{S},U}$.
		
		Now, let us show that $\mathfrak{D}^{\omega \circ p_U}_{N_G(U)}=\big(\Ind_{N_G(U)}^G(\Hr{\omega\circ p_{U}})\big)^U=\Hr{T(U,\omega)}^U$. Let  $$\mathcal{U}: \Hr{T(U,\omega)}\rightarrow \mathcal{L}_\mathcal{S}(U)$$ be a unitary operator intertwining $T(U,\omega)$ and $T_{\mathcal{S},U}$ and notice that $$\mathfrak{D}^{\omega \circ p_U}_{N_G(U)}\subseteq\Hr{T(U,\omega)}^U.$$ Let $g\in G$ and $ \phi\in \lb T(U,\omega)\big\rb(g)\mathfrak{D}^{\omega\circ p_U}_{N_G(U)}$. Since $$\lb T(U,\omega)\big\rb(g)\mathfrak{D}^{\omega\circ p_U}_{N_G(U)}\subseteq\Hr{T(U,\omega)}^{gUg^{-1}},$$ the function $\mathcal{U}(\phi)$ is a $gUg^{-1}$-left invariant function of $\mathcal{L}_\mathcal{S}(U)$ and Lemma \ref{les fonction continue de LSS sont a support compact} ensures that it is supported inside $$N_G(gUg^{-1},U)=\{t\in G: t^{-1}gUg^{-1}t\subseteq U\}=g N_G(U).$$
		Now, let $\varphi\in \Big(\mathfrak{D}^{\omega \circ p_U}_{N_G(U)}\Big)^\perp\cap\Hr{T(U,\omega)}^U$. By Lemma \ref{la correspondance Hsigma D}, we have that
		$$  \Big(\mathfrak{D}^{\omega \circ p_U}_{N_G(U)}\Big)^\perp = \overline{\underset{gN_G(U)\in G/N_G(U)-\{N_G(U)\}}{\bigoplus}\qq\big\lb T(U,\omega)\big\rb(g)\mathfrak{D}^{\omega\circ p_U}_{N_G(U)}}.$$
		In particular, the above discussion implies that $\mbox{\rm supp}(\mathcal{U}(\varphi))\subseteq G-N_G(U).$ On the other hand, since $\varphi\in \Hr{T(U,\omega)}^U$, the function $\mathcal{U}(\varphi)$ is a $U$-left invariant function of $\mathcal{L}_\mathcal{S}(U)$ and is therefore supported inside $N_G(U)$. This implies that $\mathcal{U}(\varphi)=0$. Hence, $\varphi=0$ and $\mathfrak{D}^{\omega \circ p_U}_{N_G(U)}= \Hr{T(U,\omega)}^U$.
		
		Finally, we prove the irreducibility of $T(U,\omega)$ with Proposition \ref{a criterion so that the induced of an irreducible is irreducible}. Let $M$ be a non-zero closed invariant subspace of $\Hr{T(U,\omega)}$. Then $\mathcal{U}(M)$ is a non-zero closed invariant subspace of $\mathcal{L}_{\mathcal{S}}(U)$ and Lemma \ref{Every invariant space in the regular contains a non trivial invariant vector} ensures the existence of a non-zero $U$-bi-invariant function $\varphi\in \mathcal{U}(M)$. In particular, $\mathcal{U}^{-1}(\varphi)$ is a non-zero $U$-invariant function of $\Hr{T(U,\omega)}$ contained in $M$. The result follows from the fact that $\Hr{T(U,\omega)}^U=\mathfrak{D}^{\omega \circ p_U}_{N_G(U)}$.
	\end{proof}
	This proves \circled{2}. We now prove \circled{3}.
	\begin{lemma}\label{every irred rep with seed S is the induced rep}
		Let $\pi$ be an irreducible representation of $G$ with seed $\mathcal{C}(U)$. Then $\omega_\pi$ is a $\mathcal{S}$-standard representation of $\Aut_G(\mathcal{C}(U))$ (in particular it is irreducible) and
		\begin{equation*}
		\pi \qq \simeq \qq T(U,\omega_\pi)\qq =\qq \Ind_{N_G(U)}^G(\omega_\pi \circ p_{U}).
		\end{equation*}
	\end{lemma}
	\begin{proof}
		We start by showing that $(\omega_\pi, \Hr{\pi}^{U})$ is a $\mathcal{S}$-standard representation of $\Aut_G(\mathcal{C}(U))$. Let $M$ be a non-zero closed $\Aut_G(\mathcal{C}(U))$-invariant subspace of $\Hr{\pi}^{U}$ for $\omega_\pi$ and let $\xi\in \Hr{\pi}^U$ be a non-zero vector. Since $\xi$ is cyclic for $\pi$ the function $\varphi_{\xi, \eta }\fct{G}{\C}{g}{\prods{\pi(g)\xi}{\eta}}$ is not identically zero for every non-zero $\eta\in M$. On the other hand, Lemma \ref{les fonction continue de LSS sont a support compact} ensures that this function is supported inside $N_G(U)$. In particular, there exists an element $g\in N_G(U)$ such that $0\not= \varphi_{\xi,\eta}(g)=\prods{\xi}{\omega_{\pi}\circ p_U(g^{-1})\eta}$. Since $M$ is $\Aut_G(\mathcal{C}(U))$-invariant this proves the existence of a vector $\eta'\in M$ such that $\prods{\xi}{\eta'}\not=0$. It follows that the orthogonal complement of $M$ in $\Hr{\pi}^{U}$ is trivial. This proves that $M=\Hr{\pi}^{U}$ and $\omega_\pi$ is irreducible. Now, let $W\in \tilde{\mathfrak{H}}_{\mathcal{S}}(U)$ where we recall that 
		\begin{equation*}
		\tilde{\mathfrak{H}}_{\mathcal{S}}(U)= \{W : \exists g\in G \mbox{ s.t. } gWg^{-1}\in \mathcal{S}\lb l-1 \rb \mbox{ and } U \subseteq  W \}.
		\end{equation*}
		Since $\pi(g)=\omega_\pi\circ p_{U}(g)$ for every $g\in N_G(U)$ and since $\pi$ has seed $C$, there does not exist any non-zero $W$-invariant vector in $\Hr{\pi}^{U}$ for $\omega_\pi \circ p_{U}$. This proves that $\omega_\pi$ is $\mathcal{S}$-standard. 
		
		We now prove that $\pi\simeq T(U,\omega_\pi).$ Let $\xi\in \Hr{\pi}^U$ and  consider the function $\varphi_{\xi,\xi}:G\rightarrow\C: g\mapsto \prods{ \pi(g)\xi}{\xi}$. The proof of Lemma \ref{Lemma existence of a seed} ensures that $\varphi_{\xi,\xi}$ is a $U$-bi-invariant function of $\mathcal{L}_{\mathcal{S}}(U)$ and is therefore compactly supported inside $N_G(U)$. In particular, we have that
		\begin{equation*}
		\varphi_{\xi,\xi}(g)=\begin{cases}
		\prods{\omega_\pi \circ p_U(g)\xi}{\xi} &\mbox{ if }g\in N_G(U)\\
		\q\q \q 0 &\mbox{ if }g\not\in N_G(U)
		\end{cases}.
		\end{equation*}
		On the other hand, \eqref{equation piece 1 of GNS implication} ensures that $\varphi_{\xi,\xi}$ is a function of positive type associated with $\text{ \rm Ind}_{N_G(U)}^{G}(\omega_\pi\circ p_{U})$. Since $\text{ \rm Ind}_{N_G(U)}^{G}(\omega_\pi\circ p_{U})$ is irreducible by Proposition \ref{les induite de S standard sont des sous rep de TS}, the result follows from the
		uniqueness of the GNS construction (Theorem \ref{Theorem GNS construction}).
	\end{proof}
	\begin{lemma}\label{unicitee des induite par res S standard}
		Let $\omega_1,\qq \omega_2$ be $\mathcal{S}$-standard representations of $\Aut_G(\mathcal{C}(U))$. Then, $T(U,\omega_1)$ and $T(U,\omega_2)$ are equivalent if and only if $\omega_1$ and $\omega_2$ are equivalent.
	\end{lemma}
	\begin{proof}
		Proposition \ref{les induite de S standard sont des sous rep de TS} ensures that  $\mathfrak{D}^{\omega_i \circ p_U}_{N_G(U)}=\Hr{T(U,\omega_i)}^U$. The result thus follows from Lemma \ref{les induites dinequivalente sont inequivalente} applied with $H=N_G(U)$ and $K=U$.
	\end{proof}
	\subsection{Existence criteria}\label{section existence des supercuspidal}
	Let $G$ be a non-discrete unimodular totally disconnected locally compact group and $\mathcal{S}$ be a generic filtration of $G$. If $\mathcal{S}$ factorises$^+$ at depth $l$, Theorem \ref{la version paki du theorem de classification} provides a bijective correspondence between the equivalence classes of irreducible representations of $G$ at depth $l$ with seed $C\in \mathcal{F}_{\mathcal{S}}$ and the $\mathcal{S}$-standard representations of $\Aut_G(C)$. However, it does not guarantee the existence of an irreducible representations of $G$ at depth $l$. The purpose of this section is to provide some existence criteria that will be used in Chapter \ref{Chapter application olsh facto}. The following results were used in \cite{FigaNebbia1991} to prove the existence of cuspidal representations of the full group of automorphisms $\Aut(T)$ of a regular tree and their proofs are essentially covered by \cite[Lemma 3.6, Lemma 3.7, Lemma 3.8 and Theorem 3.9]{FigaNebbia1991} but we recall them for completeness of the argument. 
	\begin{lemma}\label{les rep de Qsur Qi}
		Let $Q$ be a finite group with $\modu{Q}\gneq 2$ acting $2$-transitively on a finite set $X=\{1,...,d\}$. There exists an irreducible representation of $Q$ without non-zero $\Fix_Q(i)$-invariant vector for all $i\in X$. 
	\end{lemma}
	\begin{proof}
		Since $Q$ acts transitively on $X$, notice that $\Fix_Q(i)$ and $\Fix_Q(j)$ are conjugate to one another for all $i,j\in X$. In particular, a representation $\pi$ of $Q$ admits a non-zero $\Fix_Q(i)$-invariant vector for all $i\in X$ if and only if it admits a non-zero $\Fix_Q(1)$-invariant vector. In light of these considerations we are going to prove the existence of an irreducible representation of $Q$ without non-zero $\Fix_Q(1)$-invariant vectors.
		We recall that the quasi-regular representation $\sigma$ of $Q/\Fix_{Q}(1)$ is the representation $\text{\rm Ind}_{\Fix_{Q}(1)}^{Q}(1_{\widehat{\Fix_{Q}(1)}})$ of $Q$ induced by the trivial representation of $\Fix_{Q}(1)$. On the other hand, for every representation $\pi$ of $Q$, the Frobenius reciprocity implies that 
		\begin{equation*}
		\prods{\text{\rm Res}_{\Fix_{Q}(1)}^Q(\pi)}{1_{\widehat{\Fix_{Q}(1)}}}_{\Fix_{Q}(1)}\qq=\qq \prods{\pi}{\text{\rm Ind}_{\Fix_{Q}(1)}^{Q}(1_{\widehat{\Fix_{Q}(1)}})}_{Q}\qq=\qq  \prods{\pi}{\sigma}_{Q}.
		\end{equation*}
		In particular, every irreducible representation $\pi$ of $Q$ with a non-zero $\Fix_{Q}(1)$-invariant vector is a subrepresentation of $\sigma$. Moreover, since $Q$ acts $2$-transitively on $X$, Example \ref{Example of permutation represenatiton capturing the two transitivity revisiter} ensures the existence of an irreducible representation $\psi$ of $Q$ such that $\sigma= 1_{\widehat{Q}} \oplus \psi$. Suppose for a contradiction that every irreducible representation of $Q$ has a non-zero $\Fix_{Q}(1)$-invariant vector and is therefore contained in $\sigma$. This implies that $Q$ has two conjugacy classes and is therefore isomorphic to the cyclic group of order two which contradicts our hypothesis that $\modu{Q}\gneq2$.
	\end{proof}
	The following result provides another useful criterion that we adapt below to the context of trees.
	\begin{lemma}[\cite{FigaNebbia1991}, Lemma 3.7]\label{critede de k non deg rep for finite group}
		Let $Q$ be a finite group, let $H\leq Q$ be a direct product $H_1\times H_2\times ...\times H_s$ of non-trivial subgroups $H_i$ of $Q$ and suppose that the group of inner automorphisms of $Q$ acts by permutation on the set $\{H_1,...,H_s\}$. Then, there exists an irreducible representation $\pi$ of $Q$ without non-zero $H_i$-invariant vectors for every $i=1,...,s$. 
	\end{lemma}
	\noindent  For every locally finite tree $T$ and each subtree $\mathcal{T}\subseteq T$ we set 
	$$\Stab_G(\mathcal{T})=\{g\in G: g\mathcal{T}\subseteq\mathcal{T}\}\mbox{ and }\Fix_G(\mathcal{T})=\{g\in G: gv=v\qq\forall v\in V(\mathcal{T})\}.$$
	We obtain the following proposition.
	\begin{proposition}\label{existence criterion}
		Let $T$ be a locally finite tree, $G\leq \Aut(T)$ be a closed subgroup, $\mathcal{T}$ be a finite subtree of $T$ and  $\{\mathcal{T}_1,\mathcal{T}_2,...,\mathcal{T}_s\}$ be a set of distinct finite subtrees of $T$ contained in $\mathcal{T}$ such that $\mathcal{T}_i\cup\mathcal{T}_j=\mathcal{T}$ for every $i\not=j$. Suppose that $\Stab_G(\mathcal{T})$ acts by permutation on the set $\{\mathcal{T}_1,\mathcal{T}_2,...,\mathcal{T}_s\}$ and that $\Fix_G(\mathcal{T})\subsetneq \Fix_G(\mathcal{T}_i)\subsetneq \Stab_G(\mathcal{T})$. Then, there exists an irreducible representation of $\Stab_G(\mathcal{T})/\Fix_G(\mathcal{T})$ without non-zero $\Fix_G(\mathcal{T}_i)/\Fix_G(\mathcal{T})$-invariant vector for every $i=1,...,s$.
	\end{proposition}
	\begin{proof}
		Since $\mathcal{T}$ is a finite subtree of $T$, $\Stab_G(\mathcal{T})$ and $\Fix_G(\mathcal{T})$ are compact open subgroups of $G$. Since $\Stab_G(\mathcal{T})$ is the normaliser of $\Fix_G(\mathcal{T})$ notice that $\Stab_G(\mathcal{T})/\Fix_G(\mathcal{T})$ is a finite group and our hypothesis ensures that every $H_i=\Fix_G(\mathcal{T}_i)/\Fix_G(\mathcal{T})$ is a non-trivial  subgroup of $\Stab_G(\mathcal{T})/\Fix_G(\mathcal{T})$. On the other hand, for every $i\not=j$, $\mathcal{T}_i\cup \mathcal{T}_j=\mathcal{T}$ which ensures that $H_i\cap H_j=\{1_{\Stab_G(\mathcal{T})/\Fix_G(\mathcal{T})}\}$ and that the supports of elements of $\Fix_G(\mathcal{T}_i)$ and $\Fix_G(\mathcal{T}_j)$ are disjoint from one another. This implies that the elements of $H_i$ and $H_j$ commute with one another and that the subgroup of $\Stab_G(\mathcal{T})/\Fix_G(\mathcal{T})$ generated by $\bigcup_{i=1}^s H_i$ is isomorphic to $H_1\times H_2\times ...\times H_s$. Since the elements of $\Stab_G(\mathcal{T})$ act by permutation on  $\{\mathcal{T}_1,...,\mathcal{T}_s\}$, notice that the group of inner automorphisms of $\Stab_G(\mathcal{T})$ acts by permutation on $\{\Fix_G(\mathcal{T}_1),...,\Fix_G(\mathcal{T}_s)\}$. The result follows from Lemma \ref{critede de k non deg rep for finite group}.
	\end{proof}
	\newpage
	\thispagestyle{empty}
	\mbox{}
	\newpage
	\chapter{Application of Ol'shanskii's factorisation}\label{Chapter application olsh facto}
	\section{Introduction and main results}
	The purpose of this chapter is to present concrete applications of the axiomatic framework developed in Chapter \ref{Chapter Olshanskii's factor} to groups of automorphisms of trees and right angled buildings. In the first section, we recover the classification of cuspidal representations of groups of automorphisms of trees satisfying the Tits-independence property described in Section \ref{cuspidal representations of the group of automorphism of a tree} but this time by relying on Ol'shanskii's factorisation.  In particular, this section is redundant from the point of view of new results. It serves instead as an application of Ol'shanskii's factorisation method on a familiar setup allowing the reader to develop an intuition of our new formalism.  
	
	The purpose of Section \ref{application IPk} is to apply our axiomatic framework to groups of automorphisms of thick semi-regular trees and satisfying the property \ref{IPk}(Definition \ref{definition IPk}) as defined in \cite{BanksElderWillis2015}. Loosely speaking, this property ensures that the pointwise stabiliser of any ball of radius $k-1$ around an edge decomposes as a direct product of the subgroup fixing all the vertices on one side of the edge and the subgroup fixing all the vertices on the other side. The main contributions of this section are Theorems \ref{thm B1} and \ref{THM B} below which provide two ways to build generic filtrations that factorise$^+$. To be more precise, let $T$ be a thick semi-regular tree with set of vertices $V$. For every finite subtree $\mathcal{T}$ of $T$ and each integer $r\geq 0$, we denote by $\mathcal{T}^{(r)}$ the ball of radius $r$ around $\mathcal{T}$ for the natural metric $d_T$ on $V$ that is
	$$\mathcal{T}^{(r)}=\{v\in V: \exists w\in V(\mathcal{T})\qq\mbox{s.t.}\qq d_T(v,w)\leq r\}.$$ 
	\begin{theoremletter}\label{thm B1}
		Let $T$ be a thick semi-regular tree, $G\leq \Aut(T)$ be a closed non-discrete unimodular subgroup satisfying the property \ref{IPk} for some integer $k\geq1$,  $ \mathcal{P}$ be a complete finite subtree of $T$ containing an interior vertex, $\Sigma_{\mathcal{P}}$ be the set of maximal complete proper subtrees of $ \mathcal{P}$ and 
		$$ \mathfrak{T}_{\mathcal{P}}=\{\mathcal{R}\in \Sigma_{\mathcal{P}}: \Fix_G((\mathcal{R}')^{(k-1)})\not\subseteq \Fix_G(\mathcal{R}^{(k-1)})\qq \forall \mathcal{R}'\in \Sigma_{\mathcal{P}}-\{\mathcal{R}\}\}.$$ 
		Suppose in addition that: \begin{enumerate}
			\item $\forall \mathcal{R},\mathcal{R}'\in \mathfrak{T}_{ \mathcal{P}}, \forall g\in G$, we do not have $\Fix_G(\mathcal{R}^{(k-1)})\subsetneq\Fix_G(g(\mathcal{R}')^{(k-1)})$.
			\item For all $\mathcal{R}\in \mathfrak{T}_{ \mathcal{P}}$, $\Fix_G( \mathcal{P}^{(k-1)})\not=\Fix_G(\mathcal{R}^{(k-1)})$. 
			Furthermore, if $\Fix_G( \mathcal{P}^{(k-1)})\subsetneq\Fix_G(g\mathcal{R}^{(k-1)})$ we have that $ \mathcal{P}\subseteq g\mathcal{R}^{(k-1)}$. 
			\item $\forall n\in \N,\forall v\in V$, $\Fix_G(v^{(n)})\subseteq \Fix_G( \mathcal{P}^{(k-1)})$ implies $ \mathcal{P}^{(k-1)}\subseteq v^{(n)}$.  
			\item For every $g\in G$ such that $g \mathcal{P}\not= \mathcal{P}$, $\Fix_G( \mathcal{P}^{(k-1)})\not= \Fix_G(g \mathcal{P}^{(k-1)})$.
		\end{enumerate}  Then, there exists a generic filtration $\mathcal{S}_{ \mathcal{P}}$ of $G$ factorising$^+$ at depth $1$ with 
		$$\mathcal{S}_{ \mathcal{P}}\lb 0\rb =\{\Fix_G(\mathcal{R}^{(k-1)}): \mathcal{R}\in \mathfrak{T}_{ \mathcal{P}}\}$$
		$$\mathcal{S}_{ \mathcal{P}}\lb 1\rb=\{\Fix_G( \mathcal{P}^{(k-1)})\}.$$
	\end{theoremletter}
	\noindent The author would like to underline how realistic these assumptions are. Indeed, any closed non-discrete unimodular subgroup $G\leq \Aut(T)$ satisfying the property \ref{IPk} and such that $\Fix_G(\mathcal{T})$ does not admit any fixed point other than the vertices of $\mathcal{T}$ for every complete finite subtree $\mathcal{T}$ of $T$ satisfies the hypothesis of Theorem \ref{thm B1}. In light of Theorem \ref{la version paki du theorem de classification} we obtain a description of the irreducible representations of $G$ admitting a non-zero $\Fix_G(\mathcal{P}^{(k-1)})$ invariant vector but which do not have any non-zero $\Fix_G(\mathcal{R}^{(k-1)})$-invariant vector for any $\mathcal{R}\in \mathfrak{T}_{\mathcal{P}}$. Notice furthermore that the theorem can be applied inductively with different $\mathcal{P}$. This is done for instance in Example \ref{example fixateur point au bords}. 
	
	Under a stronger hypothesis on $G$ (the hypothesis \ref{Hypothese Hq} (Definition \ref{definition de Hq})), we are even able to explicit a generic filtration that factorises$^+$ at all sufficiently large depths. To be more precise, we have the following result.
	\begin{theoremletter}\label{THM B}
		Let $T$ be a thick semi-regular tree, $G\leq \Aut(T)$ be a closed non-discrete unimodular subgroup satisfying the hypothesis \ref{Hypothese Hq} and the property \ref{IPk} for some integers $q\geq 0$ and $k\geq 1$. Then, there exists a generic filtration $\mathcal{S}_{q}$ of $G$ that factorises$^+$ at all depths $l\geq L_{q,k}$ where $$L_{q,k}=\begin{cases}
		\max\{1,2k-q-1\} \qq& \mbox{if } q\mbox{  is even.} \\
		\max\{1,2k-q\} & \mbox{if } q\mbox{  is odd.} 
		\end{cases}$$
	\end{theoremletter}
	\noindent The generic filtration $\mathcal{S}_q$ is defined on page \pageref{page Sq filtration} and we provide some existence criteria for the irreducible representations at depth $l$ with respect to $\mathcal{S}_{q}$ in Section \ref{existence of representaions ot depth l for IPk}. Furthermore, we provide concrete applications of these results  in Examples \ref{example autT pour changer} and \ref{example fixateur point au bords}. Furthermore, we 
	
	The purpose of Section \ref{application IPV1} is to apply our axiomatic framework to groups of type preserving automorphisms of locally finite trees satisfying the property \ref{IPV1}(Definition \ref{defintion IPV1}). This serves as a preamble for Section \ref{Application to universal groups of right-angled buildings} were we show that the universal groups of certain semi-regular right-angled buildings can be realised as such groups. The main result of this section is the following.
	\begin{theoremletter}\label{theorem Ipv1 letter}
		Let $T$ be a locally finite tree and $G$ be a closed non-discrete unimodular group of type preserving automorphisms of $T$ satisfying the property \ref{IPV1}(Definition \ref{defintion IPV1}) and the hypothesis \ref{Hypothese HV1}(Definition \ref{definition H v1}). Then, there exists a generic filtration $\mathcal{S}_{V_1}$ of $G$ that factorises$^+$ at all depths $l\geq 1$. 
	\end{theoremletter} \noindent The generic filtration $\mathcal{S}_{V_1}$ is defined on page \pageref{page de SV1 yesssouille} and an existence criterion for the irreducible representations at depth $l$ with respect to $\mathcal{S}_{V_1}$ is given in Section \ref{existence for IPV1}.\\
	
	Finally, the purpose of Section \ref{Application to universal groups of right-angled buildings} is to prove that the universal groups of certain semi-regular right-angled buildings as introduced in \cite{Universal2018} can be realised as groups of type-preserving automorphisms of a locally finite tree $T$ satisfying the hypothesis of Theorem \ref{theorem Ipv1 letter}. We refer to Section \ref{section application right angle preliminaire} for a proper reminder on the notion of universal groups of semi-regular right-angled buildings and state the main result of this section assuming that the reader is familiar with this notion. We let $(W,I)$ be a finitely generated right-angled Coxeter system and suppose that $I$ can be partitioned as $I=\bigsqcup_{k=1}^rI_k$ in such a way that $I_k=\{i\}\cup \{i\}^\perp$ for every $i\in I_k$ and for all $k=1,...,r$ where $\{i\}^\perp=\{j\in I: ij=ji\}$. In particular, $W$ is virtually free and isomorphic to a free product $W_1*W_2*....*W_r$ where each of the $W_k$ is a direct product of finitely many copies of the group of order $2$. Let $(q_i)_{i\in I}$ be a set of integers greater than $2$ and $\Delta$ be a semi-regular building of type $(W,I)$ and prescribed thickness $(q_i)_{i\in I}$. Let  $(h_i)_{i\in I}$ be a set of legal colorings of $\Delta$, $Y_i$ be a set of size $q_i$ and $G_i\leq \Sym(Y_i)$ for every $i\in I$. 
	\begin{theoremletter}\label{Theorem E}
		There exists a locally finite tree $T$ such that the universal group $\mathcal{U}((h_i,G_i)_{i\in I})$ embeds as a closed subgroup of the group $\Aut(T)^+$ of type-preserving automorphisms of $T$ satisfying the property \ref{IPV1}. Furthermore, if $G_i$ is $2$-transitive on $Y_i$ for each $i\in I$, this group of automorphisms of tree corresponding to $\mathcal{U}((h_i,G_i)_{i\in I})$ satisfies the hypothesis \ref{Hypothese HV1} and is unimodular.
	\end{theoremletter} 
	Together with Theorem \ref{theorem Ipv1 letter}, this provides a generic filtration factorising$^+$ at all depths $l\geq 1$ for every non-discrete universal group $\mathcal{U}((h_i,G_i)_{i\in I})$ with $2$-transitive $G_i$. \\ 
	
	\section{The full group of automorphisms of a semi-regular tree}\label{Application to Aut T}

	The purpose of this section is to present a classification of cuspidal representations of groups of automorphisms of trees satisfying the Tits independence property by relying on the formalism developed in Chapter \ref{Chapter Olshanskii's factor}. However, in order to apply our formalism uniformly in this section, we make the additional assumption that our group satisfy the hypothesis \ref{Hypothese H Tree} defined just below. Since the classification of the cuspidal representations of groups of automorphisms of trees satisfying the Tits independence property was already given in Section \ref{cuspidal representations of the group of automorphism of a tree}, the present section is redundant from the point of view of new results. It serves instead as an application of our axiomatic framework on a familiar setup allowing the reader to develop his intuition. We adopt the same notations and terminology as in Section \ref{cuspidal representations of the group of automorphism of a tree} and recall that a group of automorphisms of trees $G\leq \Aut(T)$ satisfies the Tits independence property if every element of $G$ fixing an edge decomposes as a product of an element of $G$ fixing all the vertices on one side of the edge and an element fixing all the vertices on the other side of the edge.

	Let $T$ be a thick semi-regular tree and let $\mathfrak{T}$ be the set of non-empty complete finite subtrees of $T$.
	\begin{definition}\label{definition Hypoth H}
		A closed subgroup $G\leq \Aut(T)$ is said to satisfy the hypothesis \ref{Hypothese H Tree} if for all $\mathcal{T},\mathcal{T}'\in \mathfrak{T}$ we have that
		\begin{equation}\tag{$H$}\label{Hypothese H Tree}
		\Fix_G(\mathcal{T}')\subseteq \Fix_G(\mathcal{T})\mbox{ if and only if }\mathcal{T}\subseteq \mathcal{T}'. 
		\end{equation}
	\end{definition}
	For every such group $G$, we consider the basis of neighbourhoods of the identity $\mathcal{S}$ provided by the pointwise stabilisers of these subtrees that is 
	$$\mathcal{S}=\{\Fix_G(\mathcal{T}): \mathcal{T}\in \mathfrak{T}\},$$
	where $\Fix_G(\mathcal{T})=\{g\in G: gv=v\qq \forall v\in V(\mathcal{T})\}$. We are going to prove that $\mathcal{S}$ is a generic filtration of $G$ and that:
	\begin{itemize}
		\item $\mathcal{S}\lb 0\rb=\{\Fix_{G}(v): v\in V\}$.
		\item $\mathcal{S}\lb 1\rb=\{\Fix_{G}(e): e\in E\}$.
		\item $\mathcal{S}\lb l\rb=\{\Fix_{G}(\mathcal{T}): \mathcal{T}\in \mathfrak{T} \mbox{ and } \mathcal{T}\mbox{ has }l-1 \mbox{ interior vertices}\}\q \forall l\geq 2.$
	\end{itemize}
	In particular, the spherical, special and cuspidal representations of $G$ correspond respectively to the irreducible representations at depth $0$, $1$ and at least $2$ with respect to $\mathcal{S}$. The purpose of this section is to prove the following theorem.
	\begin{theorem}\label{Theorem classif pour Aut(T)}
		Let $G\leq \Aut(T)$ be a closed non-discrete unimodular subgroup satisfying the hypothesis \ref{Hypothese H Tree} and the Tits independence property, then $\mathcal{S}$ is a generic filtration of $G$ factorising$^+$ at all depths $l\geq 2$. 
	\end{theorem}
	\noindent Together with Theorem \ref{la version paki du theorem de classification}, this provides a description of the equivalence classes of cuspidal representations in terms of the $\mathcal{S}$-standard representations of the group of automorphisms of the corresponding seed. In Section \ref{application IPk}, we will give a different generic filtration of $\Aut(T)$ that also factorises$^+$ at all depths. These will lead to different descriptions of the same cuspidal representations of $\Aut(T)$.
	\begin{remark}
		If $G$ is locally $2$-transitive, the generic filtration $\mathcal{S}$ does not factorises at depth $1$. Indeed, let $e$ and $f$ be two different edges containing a common vertex $v\in V$, let $U=\Fix_G(e)$, $V=\Fix_G(f)$ and let $W$ be in the conjugacy class of an element of $\mathcal{S}\lb 0 \rb$ such that $U\subseteq W$. We are going to show that $W \not\subseteq V U$. From the definition of $W$, Lemma \ref{Lemma la startification de S pour Aut(T)} below ensures the existence of a vertex $w\in V$ such that $W=\Fix_G(w)$. Since $U\subseteq W$ and since $G$ satisfies the hypothesis \ref{Hypothese H Tree}, we must have $w\in e$. However, if $w\not=v$, there exists an element $g\in W$ that does not fix $v$ (because $\Fix_G(e)\not=\Fix_G(v)$) but every element of $VU$ fixes $v$. Moreover, if $w=v$, there exists $g\in W$ such that $ge=f$ (because $G$ is locally $2$-transitive) but no element of $VU$ maps $e$ to $f$. In both cases, this proves that $W\not\subseteq VU$ and therefore that $\mathcal{S}$ does not factorise at depth $1$.
	\end{remark}
 	The following lemma ensures that $\mathcal{S}$ is a generic filtration of $G$ and identify the elements at all height if $G$ satisfies the hypothesis \ref{Hypothese H Tree}.
	\begin{lemma}\label{Lemma la startification de S pour Aut(T)}
		Let $G$ be a closed non-discrete unimodular subgroup of $\Aut(T)$ satisfying the hypothesis \ref{Hypothese H Tree}. Then, $\mathcal{S}$ is a generic filtration of $G$ and:
		\begin{itemize}
			\item $\mathcal{S}\lb 0\rb=\{\Fix_{G}(v): v\in V\}.$
			\item $\mathcal{S}\lb 1\rb=\{\Fix_{G}(e): e\in E\}.$
			\item $\mathcal{S}\lb l\rb=\{\Fix_{G}(\mathcal{T}): \mathcal{T}\in \mathfrak{T} \mbox{ has }l-1 \mbox{ interior vertices}\}\q \forall l\geq 2$.
		\end{itemize}
	\end{lemma}
	\begin{proof}
		For every $\mathcal{T}\in \mathfrak{T}$ and each $g\in \Aut(T)$, notice that $g\Fix_{G}(\mathcal{T})g^{-1}$ coincides with the pointwise stabiliser $\Fix_{G}(g\mathcal{T})$. In particular, the elements of $\mathcal{F}_\mathcal{S}=\{\mathcal{C}(U): U\in \mathcal{S}\}$ are of the form
		\begin{equation*}
		\mathcal{C}(\Fix_{G}(\mathcal{T}))=\{\Fix_{G}(g\mathcal{T}): g\in G\}
		\end{equation*}
		with $\mathcal{T}\in \mathfrak{T}$. Thus, for all $\mathcal{T},\mathcal{T}'\in \mathfrak{T}$ we have that $\mathcal{C}(\Fix_{G}(\mathcal{T}'))\leq \mathcal{C}(\Fix_{G}(\mathcal{T}))$ if and only if there exists $g\in G$ such that $\Fix_{G}(\mathcal{T})\subseteq \Fix_{G}(g\mathcal{T}')$. Since $G$ satisfies the hypothesis \ref{Hypothese H Tree}, this implies that $\mathcal{C}(\Fix_{G}(\mathcal{T}'))\leq \mathcal{C}(\Fix_{G}(\mathcal{T}))$ if and only if there exists $g\in G$ such that $g\mathcal{T}'\subseteq \mathcal{T}$. In particular, since $\mathfrak{T}$ is stable under the action of $\Aut(T)$, for every chain $C_0\lneq C_1\lneq...\lneq C_{n-1}\lneq C_{n}$ of elements of $\mathcal{F}_{\mathcal{S}}$, there exists a chain $\mathcal{T}_0\subsetneq\mathcal{T}_1\subsetneq...\subsetneq\mathcal{T}_{n-1}\subsetneq\mathcal{T}_n$ of elements of $\mathfrak{T}$ such that $C_r=\mathcal{C}(\Fix_G(\mathcal{T}_r))$. Reciprocally, for every chain $\mathcal{T}_0\subsetneq\mathcal{T}_1\subsetneq...\subsetneq\mathcal{T}_{n-1}\subsetneq\mathcal{T}$ of elements of $\mathfrak{T}$ contained in a subtree $\mathcal{T}\in \mathfrak{T}$ we obtain a chain $\mathcal{C}(\Fix_G(\mathcal{T}_0))\lneq\mathcal{C}(\Fix_G(\mathcal{T}_1))\lneq...\lneq\mathcal{C}(\Fix_G(\mathcal{T}_{n-1}))\lneq\mathcal{C}(\Fix_G(\mathcal{T}))$ of elements of $\mathcal{F}_{\mathcal{S}}$ with maximal element $\mathcal{C}(\Fix_G(\mathcal{T}))$. This proves that the height of $\mathcal{C}(\Fix_G(\mathcal{T}))$  is the maximal length of a strictly increasing chain of elements of $\mathfrak{T}$ contained in $\mathcal{T}$. The result then follows from the following observations. If $\mathcal{T}$ is a vertex, it does not contain any non-empty proper subtree. If $\mathcal{T}$ is an edge, every maximal strictly increasing chain of non-empty complete subtree of $\mathcal{T}$ is of the form $\mathcal{T}_0\subsetneq \mathcal{T}$ where $\mathcal{T}_0$ is a vertex. If $\mathcal{T}$ is a complete finite subtrees of $T$ with $n\geq 1$ interior vertices, every maximal strictly increasing chain of non-empty complete subtrees of $\mathcal{T}$ is of the form $\mathcal{T}_0\subsetneq\mathcal{T}_1\subsetneq...\subsetneq\mathcal{T}_{n}\subsetneq\mathcal{T}$ 
		where $\mathcal{T}_0$ is a vertex, $\mathcal{T}_1$ is an edge and $\mathcal{T}_i$ contains $i-1$ interior vertices for every $i\geq 2$.
	\end{proof}
	
	\begin{lemma}\label{la diffenrence entre les fixing group alors la difference entre les graph}
		$\Aut(T)$ satisfies the hypothesis \ref{Hypothese H Tree}. 
	\end{lemma}
	\begin{proof}
		It is clear that $\Fix_{\Aut(T)}(\mathcal{T}')\subseteq \Fix_{\Aut(T)}(\mathcal{T})$ if $\mathcal{T}\subseteq \mathcal{T}'$. Now, suppose that $\mathcal{T}\not\subseteq \mathcal{T}'$ and let us show that $\Fix_{\Aut(T)}(\mathcal{T}')\not\subseteq \Fix_{\Aut(T)}(\mathcal{T})$. Since $\mathcal{T}\not\subseteq\mathcal{T}'$ there exists a vertex $v_\mathcal{T}\in V(\mathcal{T})-V(\mathcal{T}')$. We claim the existence of an element of $\Fix_{\Aut(T)}(\mathcal{T}')$ that does not fix $v_{\mathcal{T}}$. If $\mathcal{T}'$ is a vertex, this is obvious since $\Fix_G(\mathcal{T}')$ acts transitively on the vertices at distance $n$ from $\mathcal{T}'$. On the other hand, if $\mathcal{T}'$ contains at least two vertices, there exists a unique vertex $v'_{\mathcal{T}}\in V(\mathcal{T'})$ that is closer to $v_{\mathcal{T}}$ than every other vertex of $\mathcal{T}'$ and there exists a unique vertex $w_\mathcal{T}'$ of $\mathcal{T}'$ that is at distance one from $v'_\mathcal{T}$. Furthermore, notice that $\mathcal{T}'$ is a subtree of the half-tree $$T(w_{\mathcal{T}}, v'_\mathcal{T})\cup \{v'_\mathcal{T}\}=\{v\in V: d_T(w_\mathcal{T},v)< d_T(v_\mathcal{T}',v)\}\cup \{v_{\mathcal{T}}'\}.$$ It follows that $\Fix_{\Aut(T)}(T(w_{\mathcal{T}}, v'_\mathcal{T}))\subseteq \Fix_{\Aut(T)}(\mathcal{T}')$. On the other hand, $\Fix_{\Aut(T)}(T(w_{\mathcal{T}}, v'_\mathcal{T}))$ acts transitively on the set  $$\{v\in V-T(w_{\mathcal{T}},v'_\mathcal{T}): d_T(v'_\mathcal{T},v)=d_T(v'_\mathcal{T},v_\mathcal{T})\}.$$ 
		Since $T$ is thick, this set of vertices contains $v_{\mathcal{T}}$ and at least one other vertex. This proves the existence of an element of $\Fix_{\Aut(T)}(T(w_{\mathcal{T}},v_\mathcal{T}'))$ which does not fix $v_\mathcal{T}$. It follows that $\Fix_{\Aut(T)}(\mathcal{T}')\not\subseteq \Fix_{\Aut(T)}(\mathcal{T})$. 
	\end{proof}
	Our next task is to prove that $\mathcal{S}$ factorises$^+$ at all depths $l\geq 2$ for groups satisfying the Tits independence property. The following lemma is the key of the proof.
	\begin{lemma}\label{independence figa nebbia 3.1 version Ols facto}
		Let $G\leq \Aut(T)$ be a subgroup satisfying the Tits independence property and the hypothesis \ref{Hypothese H Tree} and let $\mathcal{T}$, $\mathcal{T}'$ be complete finite subtrees of $T$ such that $\mathcal{T}$ contains at least one interior vertex and such that $\mathcal{T}'$ does not contain $\mathcal{T}$. Then, there exists a complete proper subtree $\mathcal{R}$ of $\mathcal{T}$ such that $\Fix_{G}(\mathcal{R})\subseteq \Fix_{G}(\mathcal{T}') \Fix_{G}(\mathcal{T})$ and such that $\mathcal{C}(\Fix_{G}(\mathcal{R}))$ has the height of $\mathcal{C}(\Fix_G(\mathcal{T}))$ less $1$ in $\mathcal{F}_\mathcal{S}$.
	\end{lemma}
	\begin{proof}
		Since $G$ satisfies the Tits independence property, Lemma \ref{independence figa nebbia 3.1} ensures the existence of a complete proper subtree $\mathcal{P}$ of $\mathcal{T}$ such that $$\Fix_{G}(\mathcal{P})\subseteq \Fix_{G}(\mathcal{T}') \Fix_{G}(\mathcal{T}).$$ If $\mathcal{T}$ has exactly one interior vertex, there exists an edge $\mathcal{R}$ in $\mathcal{T}$ such that $\mathcal{P}\subseteq \mathcal{R} \subsetneq \mathcal{T}$. On the other hand, if $\mathcal{T}$ has more than one interior vertex, there exists a complete subtree $\mathcal{R}$ of $\mathcal{T}$ with one less interior vertex than $\mathcal{T}$ such that $\mathcal{P}\subseteq \mathcal{R}\subsetneq \mathcal{T}$. In both cases, this implies that $\Fix_{G}(\mathcal{R})\subseteq \Fix_{G}(\mathcal{P})\subseteq \Fix_{G}(\mathcal{T}')\Fix_{G}(\mathcal{T})$ and $\mathcal{R}$ is as desired.
	\end{proof}
	We are finally able to prove our claim.
	\begin{proof}[Proof of Theorem \ref{Theorem classif pour Aut(T)}]
		To prove that $\mathcal{S}$ factorises$^+$ a depth $l\geq 2$, we shall successively verify the three conditions of the Definition \ref{definition olsh facto}. 
		
		First, we need to prove that for all $U$ in the conjugacy class of an element of $\mathcal{S}\lb l \rb$ and every $V$ in the conjugacy class of an element of $\mathcal{S}$ such that $V \not\subseteq U$, there exists a $W$ in the conjugacy class of an element of $\mathcal{S}\lb l -1\rb$ and $U\subseteq W \subseteq V U.$
		Let $U$ and $V$ be as above. Since $\mathfrak{T}$ is stable under the action of $G$, since $G$ satisfies the hypothesis \ref{Hypothese H Tree} and as a consequence of Lemma \ref{Lemma la startification de S pour Aut(T)} there exists subtrees $\mathcal{T},\mathcal{T}'\in \mathfrak{T}$ such that $\mathcal{T}$ has  $l-1$ interior vertices, $\mathcal{T}'$ does not contain $\mathcal{T}$, $U=\Fix_{G}(\mathcal{T})$ and $V=\Fix_{G}(\mathcal{T}')$. Now, Lemma \ref{independence figa nebbia 3.1 version Ols facto} ensures the existence of a proper subtree $\mathcal{R}$ of $\mathcal{T}$ such that $\Fix_{G}(\mathcal{R})$ is conjugate to an element of $\mathcal{S}\lb l-1 \rb$ and  $$\Fix_{G}(\mathcal{R})\subseteq \Fix_{G}(\mathcal{T}')\Fix_{G}(\mathcal{T}).$$ This proves the first condition. Next, we need to prove that 
		\begin{equation*}
		N_{G}(U, V)= \{g\in {G} : g^{-1}Vg\subseteq U\}
		\end{equation*}
		is compact for every $V$ in the conjugacy class of an element of $\mathcal{S}$. Just as before, notice that $V=\Fix_{G}(\mathcal{T}')$ for some $\mathcal{T}'\in \mathfrak{T}$. Furthermore, since $G$ satisfies the hypothesis \ref{Hypothese H Tree} we have 
		\begin{equation*}
		\begin{split}
		N_{G}(U, V)&= \{g\in {G} : g^{-1}Vg\subseteq U\}= \{g\in {G} : g^{-1}\Fix_{G}(\mathcal{T}')g\subseteq \Fix_{G}(\mathcal{T})\}\\
		&= \{g\in {G} : \Fix_G(g^{-1}\mathcal{T}')\subseteq \Fix_G(\mathcal{T})\}= \{g\in {G} : g\mathcal{T}\subseteq \mathcal{T}'\}
		\end{split}
		\end{equation*}
		and this set is compact since both $\mathcal{T}$ and $\mathcal{T}'$ are finite. This proves the second condition. Finally, we need to prove that for any subgroup $W$ in the conjugacy class of an element of $\mathcal{S}\lb l-1\rb$ such that $U\subseteq W$ we have
		\begin{equation*}
		W\subseteq N_{G}(U,U) =\{g\in {G} : g^{-1}Ug\subseteq U\}.
		\end{equation*} 
		For the same reasons as before, there exist a complete subtree $\mathcal{R}\in \mathfrak{T}$ such that $W = \Fix_{G}(\mathcal{R})$. Furthermore, since $U\subseteq W$ and since $G$ satisfies the hypothesis \ref{Hypothese H Tree}, we have that $\mathcal{R}\subseteq \mathcal{T}$. On the other hand, since $\mathcal{R}$ and $\mathcal{T}$ are both complete finite subtrees and since $\mathcal{R}$ has exactly one less interior vertex than $\mathcal{T}$, every interior vertex of $\mathcal{T}$ belongs to $\mathcal{R}$. Hence, $g\mathcal{T}\subseteq  \mathcal{T}$ for every $g\in \Fix_{G}(\mathcal{R})$. This implies that 
		\begin{equation*}
		\begin{split}
		W= \Fix_{G}(\mathcal{R})&\subseteq \{g\in G : g\mathcal{T}\subseteq \mathcal{T}\}= \{g\in {G} : \Fix_{G}(\mathcal{T})\subseteq \Fix_{G}(g\mathcal{T})\}\\
		&=\{g\in {G} : g^{-1}\Fix_{G}(\mathcal{T})g\subseteq \Fix_{G}(\mathcal{T})\}\\
		&=N_{G}(\Fix_{G}(\mathcal{T}),\Fix_{G}(\mathcal{T}))= N_{G}(U,U).
		\end{split}
		\end{equation*} 
	\end{proof}
	Now, Theorem \ref{la version paki du theorem de classification} provides a bijective correspondence between the equivalence classes of cuspidal representations of $G$ with seed $\mathcal{C}(\Fix_G(\mathcal{T}))$ and the equivalence classes of $\mathcal{S}$-standard representations of $\Aut_{G}(\mathcal{C}(\Fix_{G}(\mathcal{T})))$ for every complete finite subtree $\mathcal{T}$ containing an interior vertex. On the other hand, notice from the above computations that $\Aut_{G}(\mathcal{C}(\Fix_{G}(\mathcal{T})))$ can be be identified with the group of automorphisms of $\mathcal{T}$ coming from the action of $\Stab_{G}(\mathcal{T})=\{g\in G: g\mathcal{T}\subseteq \mathcal{T}\}$ on $\mathcal{T}$. Under this identification, the $\mathcal{S}$-standard representations of $\Aut_{G}(\mathcal{C}(\Fix_{G}(\mathcal{T})))$ are the irreducible representations which do not have any non-zero invariant vector for the pointwise stabiliser of any maximal proper complete subtree of $\mathcal{T}$. The existence of these $\mathcal{S}$-standard representations of $\Aut_{G}(C)$ is proved for any seed $C$ at height $l\geq 2$ in \cite[Theorem 3.9]{FigaNebbia1991}.

	\section{Groups of automorphisms of trees with the property ${\rm IP}_k$}\label{application IPk}
	In this section, we apply our machinery to groups of automorphisms of semi-regular trees satisfying the property \ref{IPk}. The purpose of the section is to prove Theorem \ref{thm B1} and Theorem \ref{THM B}. We use the same notations and terminology as in Section \ref{section groups of automorphisms of trees}. Let $T$ be a thick semi-regular tree. For every subtree $\mathcal{T}$ of $T$ we denote by $B_T(\mathcal{T},r)$ or simply by $\mathcal{T}^{(r)}$ the ball of radius $r\geq 0$ around $\mathcal{T}$ that is 
	$$\mathcal{T}^{(r)}=\{v\in V: \exists w\in V(\mathcal{T})\qq\mbox{s.t.}\qq d_T(v,w)\leq r\}.$$
	In addition, for each subtree $\mathcal{T}$ of $T$ we let
	\begin{equation*}
		E_o(\mathcal{T})=\{(v,w):v,w\in V(\mathcal{T}), \qq v\not=w\}
	\end{equation*}
	be the set of ordered pairs of distinct adjacent vertices of $\mathcal{T}$. We recall that the elements of $E_o(\mathcal{T})$ are the \tg{oriented edges} of $\mathcal{T}$. For every oriented edge $f=(w,v)\in E_o(\mathcal{T})$, we let $\bar{f}=(v,w)$ be the oriented edge with opposite orientation and we say that $w$ and $v$ are respectively the \tg{origin} and the \tg{terminal} vertex of $f$. Finally, for every oriented edge $f=(w,v)\in E_o(T)$, we let $$Tf=T(w,v)=\{u\in V: d_T(w,u)< d_T(v,u)\}.$$
	\begin{definition}\label{definition IPk} Let $k\geq 1$ be a positive integer. A group $G\leq \Aut(T)$ satisfies the \tg{property ${\rm IP}_k$} if for all $e\in E_o(T)$ we have 
		\begin{equation}\tag{${\rm IP}_k$}\label{IPk}
			\Fix_{G}(e^{(k-1)})=\big\lb\Fix_{G}(Te)\cap\Fix_{G}(e^{(k-1)}) \big\rb \big\lb \Fix_{G}(T\bar{e})\cap\Fix_{G}(e^{(k-1)})\big\rb.
		\end{equation}
	\end{definition} 
	In particular, under our convention, ${\rm IP}_1$ is the Tits independence property as defined in Section \ref{cuspidal representations of the group of automorphism of a tree}. Theorem \ref{thm B1} states the following.
	\begin{theorem*}\label{Theorem IPK avec un arbre}
		Let $G\leq \Aut(T)$ be a closed non-discrete unimodular subgroup satisfying the property \ref{IPk} for some integer $k\geq1$, $ \mathcal{P}$ be a complete finite subtree of $T$ containing an interior vertex, $\Sigma_{\mathcal{P}}$ be the set of maximal complete proper subtrees of $ \mathcal{P}$ and 
		$$ \mathfrak{T}_{\mathcal{P}}=\{\mathcal{R}\in \Sigma_{\mathcal{P}}: \Fix_G((\mathcal{R}')^{(k-1)})\not\subseteq \Fix_G(\mathcal{R}^{(k-1)})\qq \forall \mathcal{R}'\in \Sigma_{\mathcal{P}}-\{\mathcal{R}\}\}.$$ 
		Suppose in addition that: \begin{enumerate}
			\item $\forall \mathcal{R},\mathcal{R}'\in \mathfrak{T}_{ \mathcal{P}}, \forall g\in G$, we do not have $\Fix_G(\mathcal{R}^{(k-1)})\subsetneq\Fix_G(g(\mathcal{R}')^{(k-1)})$.
			\item For all $\mathcal{R}\in \mathfrak{T}_{ \mathcal{P}}$, $\Fix_G( \mathcal{P}^{(k-1)})\not=\Fix_G(\mathcal{R}^{(k-1)})$. 
			Furthermore, if $\Fix_G( \mathcal{P}^{(k-1)})\subsetneq\Fix_G(g\mathcal{R}^{(k-1)})$ we have $ \mathcal{P}\subseteq g\mathcal{R}^{(k-1)}$. 
			\item $\forall n\in \N,\forall v\in V$, $\Fix_G(v^{(n)})\subseteq \Fix_G( \mathcal{P}^{(k-1)})$ implies $ \mathcal{P}^{(k-1)}\subseteq v^{(n)}$.  
			\item For every $g\in G$ such that $g \mathcal{P}\not= \mathcal{P}$, $\Fix_G( \mathcal{P}^{(k-1)})\not= \Fix_G(g \mathcal{P}^{(k-1)})$.
		\end{enumerate}  Then, there exists a generic filtration $\mathcal{S}_{ \mathcal{P}}$ of $G$ factorising$^+$ at depth $1$ with
		$$\mathcal{S}_{ \mathcal{P}}\lb 0\rb =\{\Fix_G(\mathcal{R}^{(k-1)}): \mathcal{R}\in \mathfrak{T}_{ \mathcal{P}}\}$$
		$$\mathcal{S}_{ \mathcal{P}}\lb 1\rb=\{\Fix_G( \mathcal{P}^{(k-1)})\}.$$
	\end{theorem*}
	\noindent In particular, Theorem \ref{la version paki du theorem de classification} leads to a description of the irreducible representations of $G$ admitting a non-zero $\Fix_G( \mathcal{P}^{(k-1)})$ invariant vectors but which do not have any non-zero $\Fix_G(\mathcal{R}^{(k-1)})$-invariant vector for any $\mathcal{R}\in \Sigma_{\mathcal{P}}$.
	
	Under stronger hypothesis on $G$, Theorem \ref{THM B} explicit a generic filtration factorising$^+$ at all depths greater than a certain constant. To be more precise, let $q\in \N$ be a non-negative integer. If $q$ is even, let $$\mathfrak{T}_{q}=\bigg\{B_T(v,r)\Big\lvert v\in V, r\geq \frac{q}{2}+1\bigg\}\sqcup\bigg\{B_T(e,r)\Big\lvert e\in E, r\geq \frac{q}{2}\bigg\}.$$
	If $q$ is odd, let  $$\mathfrak{T}_{q}=\bigg\{B_T(v,r)\Big\lvert v\in V, r\geq \frac{q+1}{2}\bigg\}\sqcup\bigg\{B_T(e,r)\Big\lvert e\in E, r\geq \frac{q+1}{2}\bigg\}.$$ 
	For any closed subgroup $G\leq \Aut(T)$, we consider the set 
	\begin{equation*} \label{page Sq filtration}
	\mathcal{S}_{q}=\{\Fix_G(\mathcal{T}): \mathcal{T}\in \mathfrak{T}_{q}\}
	\end{equation*}
	of pointwise stabilisers of those subtrees. 
	\begin{definition}\label{definition de Hq}
		A group $G\leq\Aut(T)$ is said to satisfy the hypothesis \ref{Hypothese Hq} if for all $\mathcal{T},\mathcal{T'}\in \mathfrak{T}_{q}$ we have that
		\begin{equation}\tag{$H_q$}\label{Hypothese Hq}
		\Fix_G(\mathcal{T}')\subseteq \Fix_G(\mathcal{T})\mbox{ if and only if }\mathcal{T}\subseteq \mathcal{T}'. 
		\end{equation}
	\end{definition}
	\noindent If $G\leq \Aut(T)$ is a closed non-discrete unimodular subgroup of $\Aut(T)$ satisfying the hypothesis \ref{Hypothese Hq}, Lemma \ref{la forme des Sl pour IPk} below ensures that $\mathcal{S}_{q}$ is a generic filtration of $G$ and:
	\begin{itemize}
		\item  $\mathcal{S}_{q}\lb l\rb = \{\Fix_G(B_T(e,\frac{q+l}{2})): e\in E\}$ if $q+l$ is even.
		\item  $\mathcal{S}_{q}\lb l\rb = \{\Fix_G(B_T(v,\frac{q+l+1}{2})): v\in V\}$ if $q+l$ is odd.
	\end{itemize}
	Theorem \ref{THM B} states the following.
	\begin{theorem*}\label{le theorem pour IP_k}
		Let $T$ be a thick semi-regular tree and let $G\leq \Aut(T)$ be a closed non-discrete unimodular subgroup satisfying the hypothesis \ref{Hypothese Hq} and the property \ref{IPk} for some integers $q\geq 0$ and $k\geq 1$. Then $\mathcal{S}_{q}$ factorises$^+$ at all depths $l\geq L_{q,k}$ where $$L_{q,k}=\begin{cases}
		\max\{1,2k-q-1\} \qq& \mbox{if } q\mbox{  is even.} \\
		\max\{1,2k-q\} & \mbox{if } q\mbox{  is odd.} 
		\end{cases}$$
	\end{theorem*}
	Together with Theorem \ref{la version paki du theorem de classification} this provides a description of the irreducible representations of $G$ that do not admit any non-zero $U$-invariant vector for any $U\in \mathcal{S}_q$ with depth strictly less than $L_{q,k}$. For $G=\Aut(T)$, one can take $q=0$ and $k=1$ so that the generic filtration $\mathcal{S}_{0}$ factorises$^+$ at all positive depths. This provides a second description of the cuspidal representations of $\Aut(T)$. 
	
	\subsection{Preliminaries}\label{reminders prop IP k}
	The purpose of this section is to provide reminders and an equivalent characterisation of the property \ref{IPk} that will be useful to prove that the generic filtrations below factorise. Let $T$ be a thick semi-regular tree and $k\geq 1$ be an integer. We recall from the Definition \ref{definition IPk} that a subgroup $G\leq \Aut(T)$ has the property \ref{IPk} if for every oriented edge $e\in E_o(T)$ we have 
	\begin{equation*}
		\Fix_{G}(e^{(k-1)})=\big\lb\Fix_{G}(Te)\cap\Fix_{G}(e^{(k-1)}) \big\rb \big\lb \Fix_{G}(T\bar{e})\cap\Fix_{G}(e^{(k-1)})\big\rb.
	\end{equation*}
	
	This notion has been extensively studied in \cite{BanksElderWillis2015} and we now recall some of its properties. We start with the following result.
	\begin{lemma}[{\cite[Proposition 5.3.]{BanksElderWillis2015}}]\label{IPk then IPk' for k' geq k}
		A group $G\leq \Aut(T)$ satisfying \ref{IPk} satisfies ${\rm IP}_{k'}$ for every $ k'\geq k$.
	\end{lemma}
	In particular, groups satisfying the Tits independence property satisfy the property \ref{IPk} for all $k\geq 1$. Furthermore, natural examples of groups satisfying the property \ref{IPk} can be constructed by the operation of $k$-closure. 
	\begin{definition}\label{definition kclosure}
		Let $G\leq \Aut(T)$. The $k$-\tg{closure} $G^{(k)}$ of $G$ in $\Aut(T)$ is 
		\begin{equation*}
		G^{(k)}= \left\{h\in \Aut(T) : \forall v\in V, \qq \exists g\in G  \mbox{ with } \restr{h}{B_T(v,k)}=\restr{g}{B_T(v,k)}  \right\}.
		\end{equation*}
		A group is $k$-\tg{closed} group if it coincides with its $k$-closure.
	\end{definition}
	\begin{lemma}[{\cite[Lemma 3.2, Lemma 3.4 and Proposition 5.2]{BanksElderWillis2015}}]\label{lemme Burger mozes 2}
		Let $G\leq \Aut(T)$ and let $k\geq 1$ be an integer. The $k$-closure $G^{(k)}$ is a closed subgroup of $\Aut(T)$ containing $G$ and satisfying the property \ref{IPk}.
	\end{lemma}
	
	Our current purpose is to prove Proposition \ref{G(W) for Sk-1} which plays a similar role than Lemma \ref{independence figa nebbia 3.1 version Ols facto} in Section \ref{Application to Aut T}. To this end, we reformulate \ref{IPk}. For every subtree $\mathcal{T}$ of $T$ we denote by $\partial\mathcal{T}$ the boundary of $\mathcal{T}$ and by $\partial E_o(\mathcal{T})$ the set of terminal edges of $\mathcal{T}$ that is the set of oriented edges $e\in E_o(\mathcal{T})$ of $\mathcal{T}$ with terminal vertex in $\partial\mathcal{T}$. For every $H\leq G$, for every complete finite subtree $\mathcal{T}\subseteq  T$ and for all $ f,f'\in \partial E_o(\mathcal{T})$ notice that the elements of $\Fix_{H}(Tf)$ and $\Fix_{H}(Tf')$ commute with one another since their respective support are disjoint. This remark lifts the ambiguity in the group equality of the following proposition. 
	\begin{proposition}\label{alternative definition iof property IPk}
		A group $G\leq \Aut(T)$ satisfies the property \ref{IPk} if and only if for every complete finite subtree $\mathcal{T}$ of $T$ containing an edge we have 
		\begin{equation*}\label{formula for G(Bx,k-1) if IPk}
		\Fix_G(\mathcal{T}^{(k-1)})=\prod_{f\in \partial E_o(\mathcal{T})}^{}\big\lb \Fix_{G}(Tf)\cap \Fix_G(\mathcal{T}^{(k-1)})\big\rb.
		\end{equation*}
	\end{proposition}
	\begin{proof}
		
		The reverse implication is trivial as every edge of $T$ is a complete finite subtree of $T$ containing an edge. To prove the forward implication we apply an induction on the number of interior vertices of $\mathcal{T}$ and treat several cases. If $\mathcal{T}$ does not have an interior vertex, it is an edge and the group equality corresponds exactly to the property \ref{IPk}. If $\mathcal{T}$ has exactly one interior vertex, there exists $v\in V$ such that $\mathcal{T}=B_T(v,1)$. In particular, we have $\mathcal{T}^{(k-1)}=v^{(k)}$. Let $g\in \Fix_G(v^{(k)})$, let $\partial E_o(\mathcal{T})=\{f_1,..., f_n\}$ and notice that $f^{(k-1)} \subseteq v^{(k)}$ and thus that $\Fix_G(v^{(k)})\subseteq \Fix_G(f^{(k-1)})$ for all $f\in E_o(\mathcal{T})$. In particular, since $G$ satisfies the property \ref{IPk}, we have a decomposition $g=g_{f_1}g_{\bar{f_1}}$ where $g_{f_1}\in \Fix_{G}(Tf_1)\cap \Fix_G(f_1^{(k-1)})$ and $g_{\bar{f_1}}\in  \Fix_{G}(T\bar{f_1})\cap\Fix_G(f_1^{(k-1)})$. Now, notice that  $\Fix_{G}(Tf_1\cup f_1^{(k-1)})\subseteq \Fix_G(v^{{(k)}})$ so that $g_{f_1}\in  \Fix_G(v^{{(k)}})$. Since $g=g_{f_1}g_{\bar{f_1}}$ and as $g\in \Fix_G(v^{(k)})$,  we obtain that $g_{\bar{f_1}}\in  \Fix_{G}(T\bar{f_1})\cap \Fix_G(v^{(k)})$.  It follows that  $g_{\bar{f_1}}\in \Fix_G(f_2^{(k-1)})$ and since $G$ satisfies the property \ref{IPk} we have a decomposition $g_{\bar{f_1}}=g_{f_2} g_{\bar{f_2}}$ where $g_{f_2}\in \Fix_{G}(T{f_2})\cap \Fix_G(f_2^{(k-1)})$ and $g_{\bar{f_2}}\in \Fix_{G}(T\bar{f_2})\cap \Fix_G(f_2^{(k-1)})$. Just as before, since $f_2\in \partial E_o(\mathcal{T})$, we have that $v^{(k)}\subseteq Tf_2 \cup f_2^{(k-1)}$. This implies that $\Fix_{G}(T\bar{f_2})\cap \Fix_G(f_2^{(k-1)})\subseteq \Fix_G(v^{(k)})$ and it follows that $g_{\bar{f_2}}\in\Fix_{G}(T{\bar{f_2}})\cap \Fix_G(\mathcal{T}^{(k-1)})$. On the other hand, $T{\bar{f_1}}\subseteq T{f_2}$ which implies that $g_{f_2}\in \Fix_G(T{\bar{f_1}})$. Since $g_{\bar{f_1}}\in \Fix_G(T{\bar{f_1}})$, the above decomposition of $g_{\bar{f_1}}$ implies that $g_{\bar{f_2}}\in \Fix_G(T{\bar{f_1}})$  and we obtain that $g_{\bar{f_2}}\in \Fix_{G}(T\bar{f_1}\cup T\bar{f_2})\cap \Fix_G(v^{(k)})$. Proceeding iteratively, we obtain that 
		\begin{equation*}
			g= g_{f_1}g_{f_2}... g_{f_n}g_{\bar{f_n}}
		\end{equation*}
		for some $g_{f_i}\in \Fix_{}(T{f_i})\cap \Fix_G(v^{(k)})$ and some $g_{\bar{f_n}}\in \Fix_{G}(\bigcup_{i=1}^nT{\bar{f_i}})\cap \Fix_G(v^{(k)})$. In particular since $\{f_1,..., f_n\}=\partial E_o(\mathcal{T})$, notice that $g_{\bar{f_n}}\in \Fix_G(T)=\{1_G\}$. This proves as desired that
		\begin{equation*}
			g\in \prod_{f\in \partial E_o(\mathcal{T})}^{} \big\lb\Fix_{G}(T{f})\cap \Fix_G(v^{(k)})\big\rb.
		\end{equation*}
		If $\mathcal{T}$ has at least two interior vertices, the reasoning is similar. We let $g\in \Fix_G(\mathcal{T}^{(k-1)})$, $v$ be an interior vertex of $\mathcal{T}$ at distance $1$ form the boundary $\partial\mathcal{T}$ and $\mathcal{R}$ be the unique maximal proper complete subtree of $\mathcal{T}$ for which $v$ is not an interior vertex. Notice that $\mathcal{R}$ is a complete subtree of $\mathcal{T}$ with one less interior vertex. In particular, our induction hypothesis, ensures that
		\begin{equation*}
			\Fix_G(\mathcal{R}^{(k-1)})= \prod_{f\in \partial E_o(\mathcal{R})}^{} \big\lb\Fix_{G}(Tf)\cap \Fix_G(\mathcal{R}^{(k-1)})\big\rb.
		\end{equation*}
		Now, notice that exactly one extremal edge $e\in \partial E_o(\mathcal{R})$ of $\mathcal{R}$ is not extremal in $\mathcal{T}$; the oriented edge $e$ of $\mathcal{R}$ with terminal vertex $v$. Furthermore, since $\mathcal{R}^{(k-1)}\subseteq \mathcal{T}^{(k-1)}$, we have $g\in\Fix_G(\mathcal{R}^{(k-1)})$ which implies that
		\begin{equation*}
			g= g_e\prod_{f\in \partial E_o(\mathcal{T})\cap \partial E_o(\mathcal{R})}^{} g_f
		\end{equation*} 
		for some $g_e\in \Fix_{G}(Te)\cap \Fix_G(\mathcal{R}^{(k-1)})$ and some $g_f\in \Fix_{G}(Tf)\cap \Fix_G(\mathcal{R}^{(k-1)})$. Furthermore, for every $f\in \partial E_o(\mathcal{T})\cap \partial E_o(\mathcal{R})$ notice that $$\Fix_{G}(Tf)\cap \Fix_G(\mathcal{R}^{(k-1)})=\Fix_{G}(Tf)\cap\Fix_G(\mathcal{T}^{(k-1)})$$ which implies that $g_f \in \Fix_{G}(Tf)\cap \Fix_G(\mathcal{T}^{(k-1)})$. In particular, since $g\in \Fix_G(\mathcal{T}^{(k-1)})$, the above decomposition implies that $g_e\in\Fix_G(Te)\cap \Fix_G(\mathcal{T}^{(k-1)})$. On the other hand, $v^{(k)}\subseteq T{e} \cup\mathcal{T}^{(k-1)}$ which implies that $g_e\in \Fix_G(v^{(k)})$. Hence, by the first part of the proof, we have that
		\begin{equation*}
			g_e = g_{\bar{e}} \prod_{b\in \partial E_o(\mathcal{T}) \cap \partial E_o(v^{(1)})} g_{b}
		\end{equation*}
		for some $g_{\bar{e}}\in \Fix_{G}(T{\bar {e}})\cap \Fix_G(v^{(k-1)})$ and some $g_{b}\in \Fix_{G}(Tb)\cap \Fix_G(v^{(k-1)})$. Furthermore, for every $b\in  \partial E_o(\mathcal{T}) \cap \partial E_o(v^{(1)})$, notice that $\mathcal{T}^{(k-1)}\subseteq Tb \cup v^{(1)}$ and thus that $g_{b}\in \Fix_{G}(Tb)\cap \Fix_G(\mathcal{T}^{(k-1)})$. In particular, since $g_e\in \Fix_{G}(Te)\cap \Fix_G(\mathcal{T}^{(k-1)})$ the above decomposition of $g$ implies that $g_{\bar{e}}=1_G$. This proves as desired that $g$ belongs to 
		\begin{equation*}
			\prod_{f\in \partial E_o(\mathcal{T})\cap \lb \partial E_o(\mathcal{R})\cup\partial E_o(v^{(1)})\rb }^{} \big\lb\Fix_{G}(Tf)\cap \Fix_G(\mathcal{T}^{(k-1)})\big\rb.
		\end{equation*}
	\end{proof}	 
	\begin{proposition}\label{G(W) for Sk-1}
		Let $G\leq \Aut(T)$ be a subgroup satisfying the property \ref{IPk} and let $\mathcal{T}$, $\mathcal{T}'$ be complete finite subtrees of $T$ such that $\mathcal{T}$ contains at least one interior vertex and such that $\mathcal{T}'$ does not contain $\mathcal{T}$. Then, there exists a maximal complete proper subtree $\mathcal{R}$ of $\mathcal{T}$ such that $$\Fix_G(\mathcal{R}^{(k-1)})\subseteq \Fix_G((\mathcal{T}')^{(k-1)})\Fix_G(\mathcal{T}^{(k-1)}).$$
	\end{proposition}
	\begin{proof}
		Since $\mathcal{T}$ and $\mathcal{T}'$ are complete and since $\mathcal{T}'$ does not contain $\mathcal{T}$, there exists an extremal edge $b$ of $\mathcal{T}$ which does not belong to $\mathcal{T}'$ and such that $\mathcal{T}'\subseteq T{b}$. Let $\mathcal{R}$ be the maximal complete subtree of $\mathcal{T}$ such that $b\not \subseteq \mathcal{R}$. Notice that there is a unique $e\in \partial E_o(\mathcal{R})$ which is not extremal in $\mathcal{T}$. Furthermore, observe that $\mathcal{T}'\subseteq Te\cup e$, that $\mathcal{R}$ is a complete subtree of $T$ containing an edge and that $\mathcal{R}$ has one less interior vertex than $\mathcal{T}$. In particular, Proposition \ref{alternative definition iof property IPk} ensures that 
		\begin{equation*}
		\Fix_G(\mathcal{R}^{(k-1)})= \prod_{f\in \partial E_o(\mathcal{R})} \big \lb\Fix_{G}(Tf)\cap \Fix_G(\mathcal{R}^{(k-1)})\big\rb.
		\end{equation*} 
		However, by construction, we have that $(\mathcal{T}')^{(k-1)}\subseteq Te \cup \mathcal{R}^{(k-1)}$ which implies that $$\Fix_{G}(Te)\cap \Fix_G(\mathcal{R}^{(k-1)})\subseteq \Fix_G((\mathcal{T}')^{(k-1)}).$$ 
		On the other hand, notice that $\partial E_o(\mathcal{R})-\{e\}\subseteq \partial E_o(\mathcal{T})$. Therefore, for every $f\in \partial E_o(\mathcal{R})-\{e\}$ we have $\mathcal{T}^{(k-1)}\subseteq Tf \cup \mathcal{R}^{(k-1)}$ which implies that $$\Fix_{G}(Tf)\cap \Fix_G(\mathcal{R}^{(k-1)})\subseteq \Fix_G(\mathcal{T}^{(k-1)}).$$ 
		This proves as desired that
		\begin{equation*}
		\begin{split}
		\Fix_G(\mathcal{R}^{(k-1)})&= \big\lb \Fix_{G}(Te)\cap \Fix_G(\mathcal{R}^{(k-1)})\big\rb \prod_{f\in \partial E_o(\mathcal{R})-\{e\}} \big\lb \Fix_G(Tf)\cap \Fix_G(\mathcal{R}^{(k-1)})\big \rb\\
		& \subseteq  \Fix_G((\mathcal{T}')^{(k-1)})\Fix_G(\mathcal{T}^{(k-1)}).
		\end{split}
		\end{equation*}  
	\end{proof}
	\subsection{Factorisation of the generic filtrations}\label{generic filtration fertile IP k}
	Let $T$ be a thick semi-regular tree and $G\leq \Aut(T)$ be a closed non-discrete unimodular subgroup with bi-invariant Haar measure $\mu$. 
	\begin{lemma}\label{Lemme on peut etendre la la filtration generic de mathcalT}
		Let $ \mathcal{P}$ be a complete finite subtree of $T$ containing an interior vertex, $\Sigma_{\mathcal{P}}$ be the set of maximal complete proper subtrees of $ \mathcal{P}$ and 
		$$ \mathfrak{T}_{\mathcal{P}}=\{\mathcal{R}\in \Sigma_{\mathcal{P}}: \Fix_G((\mathcal{R}')^{(k-1)})\not\subseteq \Fix_G(\mathcal{R}^{(k-1)})\qq \forall \mathcal{R}'\in \Sigma_{\mathcal{P}}-\{\mathcal{R}\}\}.$$ Suppose that: 
		\begin{enumerate}
			\item $\forall \mathcal{R},\mathcal{R}'\in \mathfrak{T}_{ \mathcal{P}}, \forall g\in G$, we do not have $\Fix_G(\mathcal{R}^{(k-1)})\subsetneq\Fix_G(g(\mathcal{R}')^{(k-1)})$.
			\item For every $\mathcal{R}\in \mathfrak{T}_{ \mathcal{P}}$, $\Fix_G( \mathcal{P}^{(k-1)})\not=\Fix_G(\mathcal{R}^{(k-1)}).$ 
		\end{enumerate}
		Then, there exists a family $\mathfrak{R}\subseteq\mathfrak{T}_{ \mathcal{P}}\cup \{ \mathcal{P}\}\cup\{v^{(n)}: n\in \N, v\in V\}$  of complete finite subtrees of $T$ such that $\mathcal{S}_{ \mathcal{P}}=\{\Fix_G(\mathcal{T}^{(k-1)}): \mathcal{T}\in\mathfrak{R}\}$ is a generic filtration of $G$ such that:
		$$\mathcal{S}_{ \mathcal{P}}\lb 0\rb =\{\Fix_G(\mathcal{R}^{(k-1)}): \mathcal{R}\in \mathfrak{T}_{ \mathcal{P}}\},\qq \mathcal{S}_{ \mathcal{P}}\lb 1\rb=\{\Fix_G( \mathcal{P}^{(k-1)})\}$$
		and $\mu(\Fix_G( \mathcal{P}^{(k-1)})\not = \mu(\Fix_G(\mathcal{T}^{(k-1)}))$ for every $\mathcal{T}\in \mathfrak{R}-\{ \mathcal{P}\}$.
	\end{lemma}
	\begin{proof}
		Since $G$ is non-discrete, there exists a vertex $v\in V$ and an integer $N\geq k$ such that $\mathcal{P}\subsetneq v^{(N)}$ and $\Fix_G(v^{(N+k-1)})<\Fix_G( \mathcal{P}^{(k-1)}).$ 
		We set $\mathfrak{R}=\mathfrak{T}_{ \mathcal{P}}\sqcup\{ \mathcal{P}\}\sqcup\{ v^{(n)}: n\geq N\}$ and let $\mathcal{S}_{ \mathcal{P}}=\{\Fix_G(\mathcal{T}^{(k-1)}): \mathcal{T}\in\mathfrak{R}\}.$ Notice by construction, that for every $\mathcal{T}\in \mathfrak{R}$ there exists $\mathcal{R}\in \mathfrak{T}_{ \mathcal{P}}$ such that $\mathcal{R}\subseteq \mathcal{T}$. In particular, for every $U\in \mathcal{S}_{ \mathcal{P}}$ there exists $\mathcal{R}\in \mathfrak{T}_{ \mathcal{P}}$ such that $U\subseteq \Fix_G(\mathcal{R}^{(k-1)})$. On the other hand, since $ \mathcal{P}$ is a finite subtree of $T$ notice that $\mathfrak{T}_{ \mathcal{P}}$ is finite. This implies that 
		$$\mu(U)\leq \qq \maxx{\mathcal{R}\in \mathfrak{T}_{ \mathcal{P}}}\big \lb\mu( \Fix_G(\mathcal{R}^{(k-1)}))\big \rb<+\infty\q \q \forall U\in \mathcal{S}_{ \mathcal{P}}.$$ 
		Renormalising $\mu$ if needed, Lemma \ref{lemma bounded basis implies generic filtration} ensures that $\mathcal{S}_{ \mathcal{P}}$ is a generic filtration of $G$. Now, notice for every $\mathcal{T}\in \mathfrak{R}- (\mathfrak{T}_{ \mathcal{P}}\sqcup\{ \mathcal{P}\})$ and every $\mathcal{R}\in \mathfrak{T}_{ \mathcal{P}}$ we have that $\mathcal{R}\subseteq  \mathcal{P}\subseteq \mathcal{T}$. In particular, this implies that $ \Fix_G(\mathcal{T}^{(k-1)}) \subseteq \Fix_G( \mathcal{P}^{(k-1)})\subseteq \Fix_G(\mathcal{R}^{(k-1)})$. On the other hand, by the hypothesis none of those inclusion is an equality which implies that $$\mu(\Fix_G(\mathcal{T}^{(k-1)}))<\mu(\Fix_G( \mathcal{P}^{(k-1)}))< \mu(\Fix_G(\mathcal{R}^{(k-1)})).$$ 
		This proves, as desired, that $\mu(\Fix_G( \mathcal{P})^{(k-1)})\not = \mu(\Fix_G(\mathcal{T}^{(k-1)}))$ for every $\mathcal{T}\in \mathfrak{R}-\{ \mathcal{P}\}$. Furthermore, since $G$ is unimodular, the measure $\mu(U)$ is an invariant of the conjugacy class $\mathcal{C}(U)$ and one realises that $$\mathcal{C}(\Fix_G(\mathcal{R}^{(k-1)}))\lneq\mathcal{C}(\Fix_G( \mathcal{P}^{(k-1)}))\lneq \mathcal{C}(\Fix_G(\mathcal{T}^{(k-1)})).$$ In particular, the depth of every subgroup  $\Fix_G(\mathcal{T}^{(k-1)})\in \mathcal{S}_{ \mathcal{P}}$ for which $\mathcal{T}\in \mathfrak{R} -(\mathfrak{T}_{ \mathcal{P}}\sqcup \{  \mathcal{P}\})$ is at least $2$. Now, let us prove that $$\mathcal{S}_{ \mathcal{P}}\lb 0\rb =\{\Fix_G(\mathcal{R}^{(k-1)}): \mathcal{R}\in \mathfrak{T}_{ \mathcal{P}}\} \mbox{ and }\mathcal{S}_{ \mathcal{P}}\lb 1\rb=\{\Fix_G( \mathcal{P}^{(k-1)})\}.$$
		Let $\mathcal{T}\in \mathfrak{R}$ and $\mathcal{R}\in \mathfrak{T}_{ \mathcal{P}}$ be such that
		$\mathcal{C}(\Fix_G(\mathcal{T}^{(k-1)}))\leq \mathcal{C}(\Fix_G(\mathcal{R}^{(k-1)})).$ By the definition, this implies the existence of an element $g\in G$ such that $\Fix_G(\mathcal{R}^{(k-1)})\leq  \Fix_G(g\mathcal{T}^{(k-1)})$ and by the first part of the proof we have that $\mathcal{T}\in \mathfrak{T}_{ \mathcal{P}}$. Our hypothesis imply that $\Fix_G(g\mathcal{T}^{(k-1)})=\Fix_G(\mathcal{R}^{(k-1)})$ and therefore that $\mathcal{C}(\Fix_G(\mathcal{T}^{(k-1)}))= \mathcal{C}(\Fix_G(\mathcal{R}^{(k-1)}))$. In particular, this proves that $\{\Fix_G(\mathcal{R}^{(k-1)}): \mathcal{R}\in \mathfrak{T}_{ \mathcal{P}}\}\subseteq \mathcal{S}_{ \mathcal{P}}\lb 0\rb$. On the other hand, we have proved that the depth of every subgroup  $\Fix_G(\mathcal{T}^{(k-1)})$ with  $\mathcal{T}\in \mathfrak{R} -(\mathfrak{T}_{ \mathcal{P}}\sqcup \{  \mathcal{P}\})$ is at least $2$. Since there must exist an element at depth $1$, this implies that 
		$\mathcal{S}_{ \mathcal{P}}\lb 1\rb=\{\Fix_G( \mathcal{P}^{(k-1)})\}$
		and it follows that $\mathcal{S}_{ \mathcal{P}}\lb 0\rb =\{\Fix_G(\mathcal{R}^{(k-1)}): \mathcal{R}\in \mathfrak{T}_{ \mathcal{P}}\}.$
	\end{proof}
	\begin{proof}[Proof of Theorem \ref{thm B1}]
		Let $\mathfrak{R}$ be the family of subtrees of $T$ given by Lemma \ref{Lemme on peut etendre la la filtration generic de mathcalT} and consider the generic filtration $\mathcal{S}_{ \mathcal{P}}=\{\Fix_G(\mathcal{T}^{(k-1)}): \mathcal{T}\in \mathfrak{R}\}$ of $G$. In order to prove that $\mathcal{S}_{ \mathcal{P}}$ factorises$^+$ at depth $1$ we shall successively verify the three conditions of the Definition \ref{definition olsh facto}. 
		
		First, we need to prove that for all $U$ in the conjugacy class of an element of $\mathcal{S}_{\mathcal{P}}\lb l \rb$ and every $V$ in the conjugacy class of an element of $\mathcal{S}_\mathcal{P}$ such that $V \not\subseteq U$, there exists a $W$ in the conjugacy class of an element of $\mathcal{S}_\mathcal{P}\lb l -1\rb$ such that $U\subseteq W \subseteq V U.$
		Let $U$ and $V$ be as above. By the definition of $\mathcal{S}_{ \mathcal{P}}$ there exist $t,h\in G$ and some $\mathcal{T} \in \mathfrak{R}$ such that $U=\Fix_G(t\mathcal{P}^{(k-1)})$ and $V= \Fix_G(h\mathcal{T}^{(k-1)})$. Furthermore, since $V\not \subseteq U$, we have that $ \Fix_G(t^{-1}h\mathcal{T}^{(k-1)})\not \subseteq\Fix_G(\mathcal{P}^{(k-1)})$. In particular, we obtain that $\mathcal{P}\not \subseteq t^{-1}h\mathcal{T}$. This follows from hypothesis $3$ if $\mathcal{T}=v^{(n)}$ for some $n\geq N$, from hypothesis $4$ if $\mathcal{T}= \mathcal{P}$ and from the fact that $t^{-1}h\mathcal{T}^{(k-1)}$ can never contain $ \mathcal{P}^{(k-1)}$ if $\mathcal{T}\in \mathfrak{T}_{ \mathcal{P}}$. In particular, Proposition \ref{G(W) for Sk-1} ensures the existence of a complete finite subtree $\mathcal{Q}\in \Sigma_{ \mathcal{P}}$ such that
		$$\Fix_G(t\mathcal{Q}^{(k-1)})\subseteq \Fix_G(h\mathcal{T}^{(k-1)})\Fix_G( t\mathcal{P}^{(k-1)})= VU.$$
		On the other hand, by the definition of $\mathfrak{T}_\mathcal{P}$ there exists $\mathcal{R}\in \mathfrak{T}_\mathcal{P}$ such that $\Fix_G(t\mathcal{R}^{(k-1)})\subseteq \Fix_G(t\mathcal{Q}^{(k-1)}).$  
		Furthermore, by the definition of $\mathcal{S}_\mathcal{P}$, the group $\Fix_G(t\mathcal{R}^{(k-1)})$ is conjugate to an element of $\mathcal{S}_{ \mathcal{P}}\lb 0\rb$ which proves the first condition. 
		
		Next, we need to prove that $N_{G}(U, V)= \{g\in {G} : g^{-1}Vg\subseteq U\}$
		is compact for every $V$ in the conjugacy class of an element of $\mathcal{S}_\mathcal{P}$. Notice that $V= \Fix_G(h\mathcal{T}^{(k-1)})$ for some some $\mathcal{T} \in \mathfrak{R}$ and some $h\in G$ and that
		\begin{equation*}
		\begin{split}
		N_{G}(U, V)&= \{g\in G : g^{-1}Vg\subseteq U\}\\
		&= \{g\in G : \Fix_{G}(g^{-1}h\mathcal{T}^{(k-1)})\subseteq \Fix_{G}( t\mathcal{P}^{(k-1)})\}
		\end{split}
		\end{equation*}
		This leads to three cases. If $\mathcal{T}=v^{(n)}$ for some $n\in \N$, the hypothesis $3$ implies that $N_{G}(U, V)=\{g\in G: t\mathcal{P}\subseteq g^{-1}hv^{(n)}\}$ which is a compact set. If $\mathcal{T}=\mathcal{T}'$, the hypothesis $4$ ensures that $N_G(V,U)=\{g\in G: g t\mathcal{P}\subseteq  t\mathcal{P}\}$ which is a compact set. If $\mathcal{T}\in \mathfrak{T}_{ \mathcal{P}}$, notice from Lemma \ref{Lemme on peut etendre la la filtration generic de mathcalT} that $\mu(V)$ is strictly smaller than $\mu(U)$ which implies that $N_G(U,V)=\es$. In every cases, this proves the second condition.

		Finally, we need to prove that for any subgroup $W$ in the conjugacy class of an element of $\mathcal{S}_{ \mathcal{P}}\lb 0\rb$ such that $U\subseteq W$ we have
		\begin{equation*}
		W\subseteq N_G(U,U) =\{g\in G : g^{-1}Ug\subseteq U\}=\{g\in G : g^{-1}Ug=U\}.
		\end{equation*} 
		The hypothesis $4$ of the theorem implies that $N_G(U,U)=\{g\in G: g t\mathcal{P}= t\mathcal{P}\}$. On the other hand, by construction of $\mathcal{S}_{ \mathcal{P}}$, for every $W$ in the conjugacy class of an element of $\mathcal{S}_{ \mathcal{P}}\lb 0\rb$, there exists an element $h\in G$ and a subtree $\mathcal{R}\in \mathfrak{T}_{ \mathcal{P}}$ such that $W=\Fix_G(h\mathcal{R}^{(k-1)})$. Furthermore, if $U\subseteq W$ the hypothesis $2$ implies that $ \mathcal{P}\subseteq t^{-1}h\mathcal{R}^{(k-1)}$. In particular, we obtain that  $W=\Fix_G(h\mathcal{R}^{(k-1)})\subseteq \Fix_G(t\mathcal{P})\subseteq N_G(U,U)$ which proves the third condition. 
	\end{proof}
	
	Our next task is to prove Theorem \ref{THM B} concerning groups satisfying the hypothesis \ref{Hypothese Hq}(Definition \ref{definition de Hq}). Let $q\in \N$ be a non-negative integer. If $q$ is even, let $$\mathfrak{T}_{q}=\bigg\{B_T(v,r)\Big\lvert v\in V, r\geq \frac{q}{2}+1\bigg\}\sqcup\bigg\{B_T(e,r)\Big\lvert e\in E, r\geq \frac{q}{2}\bigg\}.$$
	If $q$ is odd, let  $$\mathfrak{T}_{q}=\bigg\{B_T(v,r)\Big\lvert v\in V, r\geq \frac{q+1}{2}\bigg\}\sqcup\bigg\{B_T(e,r)\Big\lvert e\in E, r\geq \frac{q+1}{2}\bigg\}.$$ 
	For any closed subgroup $G\leq \Aut(T)$, we consider the set $$\mathcal{S}_{q}=\{\Fix_G(\mathcal{T}): \mathcal{T}\in \mathfrak{T}_{q}\}$$
	of pointwise stabilisers of those subtrees.  
	\begin{lemma}\label{la forme des Sl pour IPk}
		Let $G\leq \Aut(T)$ be a closed non-discrete unimodular subgroup satisfying the hypothesis \ref{Hypothese Hq} for some integer $q\geq 0$. Then, $\mathcal{S}_{q}$ is a generic filtration of $G$ and:
		\begin{itemize}
			\item $\mathcal{S}_{q}\lb l\rb = \{\Fix_G(B_T(e,\frac{q+l}{2})): e\in E\}$ if $q+l$ is even.
			\item $\mathcal{S}_{q}\lb l \rb=\{ \Fix_G(B_T(v,(\frac{q+l+1}{2}))): v\in V\}$ if $q+l$ is odd.
		\end{itemize} 
	\end{lemma}
	\begin{proof}
		For brevity and readability, for all $v\in V$ and every $e\in E$ we denote by $Gv$ and $Ge$ their respective orbit under the action of $G$ on $V$ and $E$ and by $v^{(r)}$ and $e^{(r)}$ the balls of radius $r$ around $v$ and $e$ respectively. Since $g\Fix_G(v^{(r)})g^{-1}=\Fix_G(gv^{(r)})$ $\forall g\in G, \forall v\in V$ and $\forall r\in \N$, notice that $\mathcal{C}(\Fix_G(v^{(r)})=\{\Fix_G(w^{(r)}): w\in Gv\}$. Similarly, we have that  $\mathcal{C}(\Fix_G(e^{(r)}))=\{\Fix_G(f^{(r)}): f\in Ge\}$ $\forall e\in E$. Furthermore, for every $\mathcal{T}, \mathcal{T}'\in \mathfrak{T}_{q}$, we have that $\mathcal{C}(\Fix_G(\mathcal{T}'))\leq \mathcal{C}(\Fix_G(\mathcal{T}))$ if and only if there exist $g\in G$ such that $\Fix_G(\mathcal{T})\subseteq \Fix_G(g\mathcal{T}')$. Therefore, since $G$ satisfies the hypothesis \ref{Hypothese Hq},  we have $\mathcal{C}(\Fix_G(\mathcal{T}'))\leq \mathcal{C}(\Fix_G(\mathcal{T}))$ if and only if there exists $g\in G$ such that $g\mathcal{T}'\subseteq \mathcal{T}$. Furthermore, notice that  $\mathcal{T}_{q}$ is stable under the action of $G$. In particular, for every increasing chain $C_0\lneq C_1\lneq...\lneq C_{n-1}\lneq C_{n}$ of elements of $\mathcal{F}_{\mathcal{S}_{q}}$ there exists a strictly increasing chain $\mathcal{T}_0\subsetneq\mathcal{T}_1\subsetneq...\subsetneq\mathcal{T}_{n-1}\subsetneq\mathcal{T}_n$ of elements of $\mathfrak{T}_{q}$ such that $C_t=\mathcal{C}(\Fix_G(\mathcal{T}_t))$. It follows that the height of an element $\mathcal{C}(\Fix_G(\mathcal{T}))\in \mathcal{F}_{\mathcal{S}_{q}}$ is bounded above by the maximal length of a strictly increasing chain of elements of $\mathfrak{T}_{q}$ contained in $\mathcal{T}$. On the other hand, for every strictly increasing chain $\mathcal{T}_0\subsetneq\mathcal{T}_1\subsetneq...\subsetneq\mathcal{T}_{n}\subsetneq\mathcal{T}$ of elements of $\mathfrak{T}_{q}$ contained in $\mathcal{T}$, we can build a strictly increasing chain $\mathcal{C}(\Fix_G(\mathcal{T}_0))\lneq\mathcal{C}(\Fix_G(\mathcal{T}_1))\lneq...\lneq\mathcal{C}(\Fix_G(\mathcal{T}_{n}))\lneq\mathcal{C}(\Fix_G(\mathcal{T}))$ of elements of $\mathcal{F}_{\mathcal{S}_{q}}$. This proves that the height of $\mathcal{C}(\Fix_G(\mathcal{T}))$ is the maximal length of a strictly increasing chain of elements of $\mathfrak{T}_{q}$ contained in $\mathcal{T}$. The results therefore follows from the following observation.
		If $q$ is even, every such chain is of the form
		$$e_1^{(\frac{q}{2})}\subseteq v_{1}^{(\frac{q}{2}+1)}\subseteq e_2^{(\frac{q}{2}+1)}\subseteq v_{2}^{(\frac{q}{2}+2)}\subseteq...\subseteq\mathcal{T}$$
		where $v_t\in V$ and $e_t\in E$ for all $t$ and if $q$ is odd, every such chain is of the form
		$$v_1^{(\frac{q+1}{2})}\subseteq e_{1}^{(\frac{q+1}{2})}\subseteq v_2^{(\frac{q+1}{2}+1)}\subseteq e_{2}^{(\frac{q+1}{2}+1)}\subseteq...\subseteq\mathcal{T}$$
		where the $v_t\in V$ and $e_t\in E$ for every $t$. 
	\end{proof}
	\begin{lemma}\label{the lemma imply fertility for G S on IPk groups}
		Let $G\leq \Aut(T)$ be a closed unimodular subgroup satisfying the hypothesis \ref{Hypothese Hq} and the property \ref{IPk} for some integers $q\geq 0$, $k\geq 1$ and let $$L_{q,k}=\begin{cases}
		\max\{1,2k-q-1\} \qq& \mbox{if } q\mbox{  is even.} \\
		\max\{1,2k-q\} & \mbox{if } q\mbox{  is odd.} 
		\end{cases}$$ Suppose further that $l,l'\in \N$ are such that $l\geq L_{q,k}$ and $l \leq l'$. Then, for every $U$ in the conjugacy class of an element of $\mathcal{S}_q\lb l \rb$ and every $V$ in the conjugacy class of an element of $\mathcal{S}_q\lb l'\rb$ such that $V\not \subseteq U$, there exists $W\in\mathcal{S}_q\lb l-1\rb$ such that $U\subseteq W\subseteq VU$. 
	\end{lemma}
	\begin{proof} For every $t\in \N$, let $\mathfrak{T}_q\lb t\rb=\{\mathcal{T}\in \mathfrak{T}_q: \Fix_G(\mathcal{T})\in \mathcal{S}\lbrack t\rbrack\}$.
		Notice that $\mathcal{T}_q\lb t\rb$ is stable under the action of $G$. In particular, there exist $\mathcal{T}\in \mathfrak{T}_q\lbrack l\rbrack$ and $\mathcal{T}'\in \mathfrak{T}_{q}\lb l'\rb$ such that $U= \Fix_G(\mathcal{T})$ and $V= \Fix_G(\mathcal{T}')$. Since $V\not \subseteq U$, we have that $\mathcal{T}\not\subseteq \mathcal{T}'$.
		There is four cases to treat depending on the parity of $q$ and $l$. We suppose that $q$ is even (the reasoning with odd $q$ is similar). If $l$ is even, let $k'=\frac{l+q}{2}$. Lemma \ref{la forme des Sl pour IPk} implies the existence of an edge $e\in E$ such that $U=\Fix_G(B_T(e, k'))$. Furthermore, since $l\geq L_{q,k}$ and since $l$ is even, we have $k'\geq k$ and Lemma \ref{IPk then IPk' for k' geq k} implies that $G$ satisfies the property ${\rm IP}_{k'}$. Therefore, Proposition \ref{G(W) for Sk-1} ensures the existence of a vertex $v\in e$ such that 
		$$\Fix_G(B_T(v,k'))\subseteq VU$$
		and Lemma \ref{la forme des Sl pour IPk} ensures that $\Fix_G(B_T(v,k'))\in \mathcal{S}_q\lb l-1\rb$. Finally, notice that $U=\Fix_G(B_T(e,k'))\subseteq \Fix_G(B_T(v,k'))$ since $v\in e$. 
		
		If $l$ is odd, let $k'=\frac{l+q+1}{2}$. Lemma \ref{la forme des Sl pour IPk} implies the existence of some  $v\in V$  such that $U=\Fix_G(B_T(v, k'))$. Furthermore, since $l\geq L_{q,k}$, we have $k'\geq k$ and Lemma \ref{IPk then IPk' for k' geq k} implies that $G$ satisfies the property ${\rm IP}_{k'}$. Therefore, Proposition \ref{G(W) for Sk-1} ensures the existence of an edge $e\subsetneq B_T(v,1)$ such that
		$$\Fix_G(B_T(e,k'-1))\subseteq VU$$
		and Lemma \ref{la forme des Sl pour IPk} ensures that $\Fix_G(B_T(e,k'-1)\in \mathcal{S}_q\lb l-1\rb$. Finally, since $e\in E(B_T(v,1))$, notice that $\Fix_G(B_T(v,k'))\subseteq \Fix_G(B_T(e,k'-1))$. 
	\end{proof}
	\begin{proof}[Proof of Theorem \ref{THM B}]
		To prove that $\mathcal{S}_{q}$ factorises$^+$ at depth $l\geq L_{q,k}$ we shall successively verify the three conditions of Definition \ref{definition olsh facto}.
		
		First, we need to prove that for every $U$ in the conjugacy class of an element of $\mathcal{S}_q\lb l \rb$ and every $V$ in the conjugacy class of an element of $\mathcal{S}_q$ with $V \not\subseteq U$, there exists a $W$ in the conjugacy class of an element of $\mathcal{S}_q\lb l -1\rb$ such that $U\subseteq W \subseteq V U.$ Let $U$ and $V$ be as above. If $V$ is conjugate to an element of $\mathcal{S}_q\lb l'\rb$ for some $l'\geq l$ the result follows directly from Lemma \ref{the lemma imply fertility for G S on IPk groups}. Therefore, we suppose that $l'< l$. By the definition of $\mathcal{S}_q$ and since $\mathfrak{T}_q$ is stable under the action of $G$, there exist $\mathcal{T},\mathcal{T}'\in \mathfrak{T}_q$ such that $U=\Fix_G(\mathcal{T})$ and $V=\Fix_G(\mathcal{T}')$. We have two cases. Either, $\mathcal{T}'\subseteq \mathcal{T}$ and there exists a subtree $\mathcal{R}\in \mathfrak{T}_{q}$ such that $\mathcal{T}'\subseteq \mathcal{R}\subseteq \mathcal{T}$ and $\Fix_G(\mathcal{R})\in \mathcal{S}_q\lb l-1\rb$. In that case $$\Fix_G(\mathcal{T})\subseteq \Fix_G(\mathcal{R})\subseteq \Fix_G(\mathcal{T}')\subseteq \Fix_G(\mathcal{T}')\Fix_G(\mathcal{T}).$$ Or else, $\mathcal{T}'\not \subseteq \mathcal{T}$ and since $l'< l$, this implies the existence of a subtree $ \mathcal{P}\in \mathfrak{T}_q$ such that $\mathcal{T}'\subseteq  \mathcal{P}\not= \mathcal{T}$ and $\Fix_G( \mathcal{P})\in \mathcal{S}_q\lb l\rb$. In particular, Lemma \ref{the lemma imply fertility for G S on IPk groups} ensures the existence of a $W\in \mathcal{S}_q\lb l-1\rb$ such that $U\subseteq W \subseteq \Fix_G( \mathcal{P})U$. Since $\Fix_G( \mathcal{P}) \subseteq \Fix_G(\mathcal{T}')$, this proves the first condition. 
		
		Next, we need to prove that $N_{G}(U, V)= \{g\in {G} : g^{-1}Vg\subseteq U\}$ is compact for every $V$ in the conjugacy class of an element of $\mathcal{S}_q$. Just as before, notice that $V=\Fix_G(\mathcal{T}')$ for some $ \mathcal{T}'\in \mathfrak{T}_q$. Since $G$ satisfies the hypothesis \ref{Hypothese Hq} notice that
		\begin{equation*}
		\begin{split}
		N_{G}(U, V)&= \{g\in G : g^{-1}Vg\subseteq U\}=\{g\in G : g^{-1}\Fix_G(\mathcal{T}')g\subseteq \Fix_G(\mathcal{T})\}\\ 
		&=\{g\in G : \Fix_G(g^{-1}\mathcal{T}')\subseteq \Fix_G(\mathcal{T})\}= \{g\in G : g\mathcal{T}\subseteq \mathcal{T}'\}.
		\end{split}
		\end{equation*}
		In particular, since both $\mathcal{T}$ and $\mathcal{T}'$ are finite subtrees of $T$, $N_G(U,V)$ is a compact subset of $G$ which proves the second condition.
		
		Finally, we need to prove that for every $W$ in the conjugacy class of an element of $\mathcal{S}_q\lb l-1\rb$ with $U\subseteq W$ we have
		\begin{equation*}
		W\subseteq N_G(U,U) =\{g\in G : g^{-1}Ug\subseteq U\}.
		\end{equation*} 
		For the same reasons as before, there exists $\mathcal{R}\in \mathfrak{T}_q$ such that $W= \Fix_G(\mathcal{R})$. On the other hand, since $U\subseteq W$ and since $G$ satisfies the hypothesis \ref{Hypothese Hq}, notice that $\mathcal{R}\subseteq \mathcal{T}$. Furthermore, notice that $\Fix_G(\mathcal{R})$ has depth $l-1$ and therefore that $\mathcal{R}$ contains every interior vertex of $\mathcal{T}$. Since $G$ is unimodular and satisfies the hypothesis \ref{Hypothese Hq} this implies that
		\begin{equation*}
		\begin{split}
		\Fix_G(\mathcal{R})&\subseteq \{h\in G : h\mathcal{T}\subseteq \mathcal{T}\}= \{h\in G : \Fix_G(\mathcal{T})\subseteq\Fix_G(h\mathcal{T})\}\\
		&=\{h\in G : h^{-1}\Fix_G(\mathcal{T})h\subseteq\Fix_G(\mathcal{T})\}= N_G(U,U)
		\end{split}
		\end{equation*}
		which proves the third condition.
	\end{proof}
	\noindent We provide examples of groups satisfying the hypothesis of these Theorems. 
		\begin{example}\label{example autT pour changer}
		Let $T$ be a thick semiregular tree and consider the full group $\Aut(T)$. By Lemma \ref{la diffenrence entre les fixing group alors la difference entre les graph}, $\Aut(T)$ satisfies the hypothesis \ref{Hypothese H Tree} and hence the hypothesis $H_0$. Furthermore, $\Aut(T)$ coincides with its $1$-closure and therefore satisfies the property \ref{IPk} for every integer  $k\geq 1$ by Lemmas \ref{IPk then IPk' for k' geq k} and \ref{lemme Burger mozes 2}. Since $\Aut(T)$ is non-discrete and unimodular Theorem \ref{THM B} applies and the generic filtrations $\mathcal{S}_0$ factorises$^+$ at all depths $l\geq 1$. In particular, Theorem \ref{la version paki du theorem de classification} provides a description of all the irreducible representations at depth $l\geq 1$ for $\mathcal{S}_0$. Due to Lemma \ref{la forme des Sl pour IPk}, these are exactly the cuspidal representations of $\Aut(T)$. Notice however that $\mathcal{S}_{0}$ is quite different from the generic filtration $\mathcal{S}$ of $\Aut(T)$ considered in Chapter \ref{Application to Aut T} so that this procedure leads to different description of these representations.
	\end{example} 
	\begin{example}\label{example fixateur point au bords}
		Let $T$ be a thick semiregular tree and consider a $k$-closed group $G\leq \Aut(T)$ (Definition \ref{definition kclosure}). Let $\omega\in \partial T$ be an end $T$ and consider the stabiliser of the $\omega$-horocycles $$G_\omega^0 =\{g\in G: g\omega = \omega \mbox{ and }\exists v\in V\mbox{ s.t. }gv=v\}.$$ Notice that $G_\omega^0$ is still $k$-closed and hence satisfies the property \ref{IPk}. Now, consider an infinite geodesic $\gamma=(v_0,v_1,...)$ of $T$ with end $\omega$ and notice that 
		$$G_\omega^0=\bigcup_{n\in \N} \Fix_{G_\omega^0}(v_n).$$ In particular, $G_\omega^0$ is a union of compact groups and is therefore unimodular by Examples \ref{example compaclty generated is unimodular}. If $G_\omega^0$ is non-discrete and satisfies the hypothesis \ref{Hypothese Hq}, it satisfies the hypothesis of Theorem \ref{THM B}. In particular, in that case, Theorem \ref{la version paki du theorem de classification} provides a description of all the irreducible representations of $G_{\omega}^0$ at depth $l\geq L_{q,k}$ for $\mathcal{S}_q$. However, $G_{\omega}^0$ never satisfies the hypothesis $H_0$ since for any edges $e,f\in E$ along an infinite geodesic with end $\omega$ we have either that $\Fix_{G^0_\omega}(e)\subseteq \Fix_{G^0_\omega}(f)$ or that $\Fix_{G^0_\omega}(f)\subseteq \Fix_{G^0_\omega}(e)$. Nevertheless, in certain cases, a description of the remaining cuspidal representations of $G_{\omega}^0$ can be obtained using Theorem \ref{thm B1}. For instance let $G=\Aut(T)$. In that case, $G_\omega^0$ satisfies the hypothesis $H_1$ and the generic filtration $\mathcal{S}_1$ factorises$^+$ at all depths $l\geq 1$. In particular, by Theorem \ref{THM B} we obtain a description of the cuspidal representations admitting non-zero invariant vectors for the pointwise stabiliser of a ball of radius one around an edge or bigger but not for the pointwise stabiliser of a ball of radius one around a vertex. To obtain a description of the cuspidal representations admitting non-zero invariant vectors for the pointwise stabiliser of a ball $B_T(v,1)$ of radius $1$ around a vertex $v\in V$, we let $\mathcal{P}=B_T(v,1)$, notice that $\Sigma_{\mathcal{P}}=\{e\}$ where $e$ is the only edge of $B_T(v,1)$ contained in the geodesic $\lb v,\omega\rb$  and apply Theorem \ref{thm B1}. For $G=\Aut(T)$, the reaming irreducible representations of $G^0_\omega$ are all spherical and are classified in \cite{Nebbia1990}.
	\end{example}
	Other applications of Theorem \ref{thm B1} and Theorem \ref{THM B} could be made for instance on the $k$-closure of certain groups of automorphisms of trees and on the generalisation of Burger-Mozes groups described in \cite{Tornier2020}. 
	\subsection{Existence of $\mathcal{S}_q$-standard representations}\label{existence of representaions ot depth l for IPk}
	Let $T$ be a  thick semi-regular tree and let $q\in \N$ be a non-negative integer. If $q$ is even, let $$\mathfrak{T}_{q}=\bigg\{B_T(v,r)\Big\lvert v\in V, r\geq \frac{q}{2}+1\bigg\}\sqcup\bigg\{B_T(e,r)\Big\lvert e\in E, r\geq \frac{q}{2}\bigg\}.$$
	If $q$ is odd, let  $$\mathfrak{T}_{q}=\bigg\{B_T(v,r)\Big\lvert v\in V, r\geq \frac{q+1}{2}\bigg\}\sqcup\bigg\{B_T(e,r)\Big\lvert e\in E, r\geq \frac{q+1}{2}\bigg\}.$$ 
	For any closed non-discrete unimodular subgroup $G\leq\Aut(T)$ satisfying the hypothesis \ref{Hypothese Hq}, we have shown that  $$\mathcal{S}_q=\{\Fix_G(\mathcal{T}): \mathcal{T}\in \mathfrak{T}_q\}$$ is a generic filtration of $G$. Furthermore, if $G$ satisfies the property \ref{IPk} for some integer $k\geq 1$ we have shown that $\mathcal{S}_q$ factorises$^+$ at all depths $l\geq L_{q,k}$ where$$L_{q,k}=\begin{cases}
	\max\{1,2k-q-1\} \qq& \mbox{if } q\mbox{  is even.} \\
	\max\{1,2k-q\} & \mbox{if } q\mbox{  is odd.} 
	\end{cases}$$ 
	In particular, Theorem \ref{la version paki du theorem de classification} provides a bijective correspondence between the equivalence classes of irreducible representations of $G$ at depth $l\geq L_{q,k}$ with seed $C\in \mathcal{F}_{\mathcal{S}_q}$ and the $\mathcal{S}_q$-standard representations of $\Aut_G(C)$. However, no results so far ensures the existence of such representations of $G$. The purpose of the present section is to study the existence of those ${\mathcal{S}_q}$-standard representations. The following result ensures the existence of ${\mathcal{S}_q}$-standard representations of $\Aut_{G}(C)$ for all $C\in \mathcal{F}_{\mathcal{S}_q}$ with height $l\geq L_{q,k}$ if $q$ and $l$ have the same parity.
	\begin{proposition}\label{existence of standard for edges for IPK}
		Let $G\leq \Aut(T)$ be a closed non-discrete unimodular subgroup satisfying the hypothesis \ref{Hypothese Hq} and the property \ref{IPk} for some integers $q\geq 0$, $k\geq 1$ and let $l\geq L_{q,k}$. Suppose that one of the following happens:
		\begin{itemize}[leftmargin=*]
			\item $q$ and $l$ are even.
			\item $q$ and $l$ are odd but $l\not=1$.
			\item $q$ is odd, $l=1$ and $\Fix_G(v^{(\frac{q}{2}+1)})\not= \{g\in G: ge=e\}$  $\forall e\in E$, $\forall v\in e$.
		\end{itemize} Then, there exists a ${\mathcal{S}_q}$-standard representation of $\Aut_{G}(C)$ for every $C\in \mathcal{F}_{\mathcal{S}_q}$ at height $l$.
	\end{proposition} 
	\begin{proof}
		Let $C\in \mathcal{F}_{\mathcal{S}_q}$ be at height $l$. Since $q$ and $l$ have the same parity, Lemma \ref{la forme des Sl pour IPk} ensures the existence of an edge $e\in E$ and an integer $r\geq k$ and such that $B_T(e,r)\in \mathfrak{T}_q$ and $C=\mathcal{C}(\Fix_{G}(B_T(e,r))$. For brevity we let $\mathcal{T}$ denote the subtree $B_T(e,r)$. Since $G$ satisfies the hypothesis \ref{Hypothese Hq} and as a consequence of Lemma \ref{la forme des Sl pour IPk}, notice that $N_G(\Fix_G(\mathcal{T}))=\{g\in G : g\mathcal{T}\subseteq \mathcal{T}\}=\Stab_{G}(\mathcal{T})=\{g\in G : ge=e\}$, that $\Aut_G(C)\simeq \Stab_G(\mathcal{T})/\Fix_G(\mathcal{T})$ and that 
		\begin{equation*}
		\begin{split}
		\tilde{\mathfrak{H}}_{\mathcal{S}_q}(\Fix_G(\mathcal{T}))&= \{W : \exists g\in G \mbox{ s.t. } gWg^{-1}\in \mathcal{S}_q\lb l-1 \rb \mbox{ and } \Fix_G(\mathcal{T}) \subseteq  W \}\\ &=\{\Fix_G(B_T(v,r)): v\in e\}.
		\end{split}
		\end{equation*}
		Let $v_0,v_1$ denote the two vertices of $e$, let $\mathcal{T}_i=B_T(v_i,r)$ and notice that $\mathcal{T}_0\cup \mathcal{T}_1=\mathcal{T}$. In particular, the action of $N_G(\Fix_G(\mathcal{T}))$ on $\mathcal{T}$ permutes the subtrees $\{\mathcal{T}_0, \mathcal{T}_1\}$. On the other hand, since $G$ satisfies the hypothesis \ref{Hypothese Hq} our hypothesis imply that $$\Fix_G(\mathcal{T})\subsetneq \Fix_G(\mathcal{T}_i)\subsetneq \Stab_G(\mathcal{T}).$$  The result therefore follows from Proposition \ref{existence criterion}.
	\end{proof}
	The following two results ensures the existence of ${\mathcal{S}_q}$-standard representations of $\Aut_{G}(C)$ for all $C\in \mathcal{F}_{\mathcal{S}_q}$ with height $l\geq L_{q,k}$ if $q$ and $l$ have opposite parity. We start with the degenerate case $q=0,k=1$ and $l=1$ where Proposition \ref{existence criterion} does not apply.
	\begin{lemma}\label{existence of standard for vertices Tits}
		Let $G\leq \Aut(T)$ be a closed non-discrete unimodular subgroup satisfying the hypothesis $H_0$, the Tits independence property ${\rm IP}_1$ and such that $\Fix_G(v)$ is $2$-transitive on the set of edges of $B_T(v,1)$ for every $v\in V$. Then, there exists a $\mathcal{S}_{0}$-standard representation of $\Aut_{G}(C)$ for every $C\in \mathcal{F}_{\mathcal{S}_0}$ at height $1$. 
	\end{lemma}
	\begin{proof}
		Let $C\in \mathcal{F}_{\mathcal{S}_{0}}$ be at height $1$. Lemma \ref{la forme des Sl pour IPk} ensures the existence of a vertex $v\in V$ such that $C=\mathcal{C}(\Fix_{G}(B_T(v,1))$. Let $U=\Fix_{G}(B_T(v,1))$ and notice that
		$$N_G(U)= \{g\in G : gB_T(v,1)\subseteq B_T(v,1)\}= \{ g\in G : gv=v\}=\Fix_G(v).$$
		Furthermore, since $G$ satisfies the hypothesis $H_0$, Lemma \ref{la forme des Sl pour IPk} implies that 
		\begin{equation*}
		\begin{split}
		\tilde{\mathfrak{H}}_{\mathcal{S}_{0}}(U)&= \{W : \exists g\in G \mbox{ s.t. } gWg^{-1}\in \mathcal{S}_q\lb l-1 \rb \mbox{ and } \Fix_G(B_T(v,1)) \subseteq  W \}\\ &= \{\Fix_G(e): e\in E(B_T(v,1))\}
		\end{split}
		\end{equation*}
		where $E(B_T(v,1))$ denotes the set of edges of $B_T(v,1)$. Let $d$ be the degree of $v$ in $T$, let $X=E(B_T(v,1))$ and notice that our hypothesis imply that $\Aut_{G}(C)\simeq \Fix_G(v)/\Fix_G(B_T(v,1))$ is $2$-transitive $X$. In particular, Lemma \ref{les rep de Qsur Qi} implies the existence of an irreducible representation $\sigma$ of $\Aut_{G}(C)$ without non-zero $\Fix_{\Aut_{G}(C)}(e)$-invariant vectors for all $e\in X$. Since $\mathfrak{H}_{\mathcal{S}_{0}}(U)=\{p_U(\Fix_G(e)): e\in X\}=\{\Fix_{\Aut_{G}(C)}(e): e\in X\}$, this proves the existence of a $\mathcal{S}_{0}$-standard representation of $\Aut_{G}(C)$.  
	\end{proof}
	\noindent The following result treats the remaining cases.
	\begin{lemma}\label{existence of standard for odd depth}
		Let $G\leq \Aut(T)$ be a closed non-discrete unimodular subgroup satisfying the hypothesis \ref{Hypothese Hq} and the property \ref{IPk} for some integers $q\geq 0$, $k\geq 1$ and let $l\geq L_{q,k}$. Suppose further that  $$\Fix_G((B_T(v,1)-\{w\})^{(r)})\not= \Fix_G(B_T(v,r+1))\qq \forall v\in V,\forall w\in B_T(v,1)-\{v\}$$ for all $r\geq \frac{q}{2}$ if $q$ is even and for all $r\geq \frac{q-1}{2}$ if $q$ is odd and that one of the following happens:
		\begin{itemize}[leftmargin=*]
			\item $q$ is odd and $l$ is even.
			\item $q$ is even, $l$ is odd and $l\not=1$.
			\item $q$ is even, $q\not=0$, $l=1$  and $\forall v\in V$,  $\forall w\in  B_T(v,1)-\{v\}$ we have $$\Fix_G((B_T(v,1)-\{w\})^{(\frac{q}{2})})\not= \Fix_G(v).$$ 
		\end{itemize} 
		Then, there exists a $\mathcal{S}_q$-standard representation of $\Aut_{G}(C)$ for every $C\in \mathcal{F}_{\mathcal{S}_q}$ at height $l$. 
	\end{lemma}
	\begin{proof}
		Suppose that $C\in \mathcal{F}_{\mathcal{S}_q}$ is at height $l$. Since $q$ and $l$ have opposite parity, Lemma \ref{la forme des Sl pour IPk} ensures the existence of a vertex $v\in V$ and an integer $r\geq k-1$ such that $B_T(v,r+1)\in \mathfrak{T}_q$ and $C=\mathcal{C}(\Fix_{G}(B_T(v,r+1)))$. For brevity we let $\mathcal{T}$ denote the subtree $B_T(v,r+1)$. Since $G$ satisfies the hypothesis \ref{Hypothese Hq} and as a consequence of Lemma \ref{la forme des Sl pour IPk}, notice that $N_G(\Fix_G(\mathcal{T}))=\{g\in G : g\mathcal{T}\subseteq\mathcal{T}\}= \Stab_G(\mathcal{T})=\{g\in G : gv=v\}$, that $\Aut_G(C)\simeq \Stab_G(\mathcal{T})/\Fix_G(\mathcal{T})$ and that 
		\begin{equation*}
		\begin{split}
		\tilde{\mathfrak{H}}_{\mathcal{S}_q}(\Fix_G(\mathcal{T}))&= \{W: \exists g\in G \mbox{ s.t. }gWg^{-1}\in \mathcal{S}_q\lb l-1 \rb \mbox{ and }\Fix_G(\mathcal{T}) \subseteq W \}\\
		&= \{\Fix_G(B_T(e,r)): e\in E(B_T(v,1))\}.
		\end{split}
		\end{equation*}
		Now, let $\{w_1,...,w_d\}$ be the leaves of $B_T(v,1)$, let $\mathcal{T}_i= (B_T(v,1)-\{w_i\})^{(r-1)}$ $ i=1,...,d$ and notice that $\mathcal{T}_i\cup \mathcal{T}_j=\mathcal{T}$ $\forall i\not=j$. On the other hand, the action of $\Stab_G(\mathcal{T})$ on $T$ permutes the subtrees $\{\mathcal{T}_1,...,\mathcal{T}_d\}$ and since each $\mathcal{T}_i$ contains $v$ we have that $\Fix_G(\mathcal{T}_i)\subseteq \Stab_G(\mathcal{T})$  $\forall i=1,...,d$. Furthermore, the hypothesis on $G$ imply that $\Fix_G(\mathcal{T})\subsetneq \Fix_G(\mathcal{T}_i)\subsetneq \Stab_G(\mathcal{T})$. In particular, Proposition \ref{existence criterion} ensures the existence of an irreducible representation $\sigma$ of $\Aut_{G}(C)$ without non-zero $p_{\Fix_G(\mathcal{T})}(\Fix_G(\mathcal{T}_i))$-invariant vectors. Moreover, for every edge $e\in E(B_T(v,1))$ there exists some $i\in \{1,...,d\}$ such that $B_T(e,r)\subseteq \mathcal{T}_i$ which implies that $p_{\Fix_G(\mathcal{T})}(\Fix_G(\mathcal{T}_i))\subseteq p_{\Fix_G(\mathcal{T})}(\Fix_G(B_T(e,r))$. Hence, $\sigma$ is a $\mathcal{S}_q$-standard representation of $\Aut_{G}(C)$.
	\end{proof}

	\section{Groups of automorphisms of trees with the property ${\rm IP}_{V_1}$}\label{application IPV1}
	In this section, we apply our machinery to groups of automorphisms of locally finite trees satisfying the property \ref{IPV1}(Definition \ref{defintion IPV1}). We use the same notations and terminology as in Section \ref{section groups of automorphisms of trees}. Let $T$ be a locally finite tree and let $\Aut(T)^+$ be the group of type-preserving automorphisms of $T$. 
	\begin{definition}\label{defintion IPV1}
		A group $G\leq \Aut(T)^+$ is said to satisfy the \tg{property \ref{IPV1}}, if there exists a bipartition $V=V_0\sqcup V_1$ such that every edge of $T$ contains exactly one vertex in each $V_i$ and such that $\forall w\in V_1$ we have 
		\begin{equation}\tag{${\rm IP }_{V_1}$}\label{IPV1}
		\Fix_G(B_T(w,1))=\prod_{v\in B_T(w,1)-\{w\}}\Fix_G(T(w,v)) 
		\end{equation}
		where $T(w,v)=\{u\in V : d_T(w,u)< d_T(v,u)\}$.
	\end{definition}
	\begin{example}
		Let $T$ be a locally finite semi-regular tree and let $G\leq \Aut(T)^+$ satisfy the property ${\rm IP}_1$. Then, $G$ satisfies the property ${\rm IP}_{V_1}$. 
	\end{example}
	\noindent Other examples are given in Section \ref{Application to universal groups of right-angled buildings} where we show that the universal groups of certain semi-regular right-angled buildings can be realised as closed subgroups of $ \Aut(T)^+$ satisfying the property ${\rm IP}_{V_1}$ but where $T$ is in general not semi-regular (Theorem \ref{theorem pour les right angled buildings}).
	
	The purpose of the present section is to prove Theorem \ref{theorem Ipv1 letter} which provides an explicit generic filtration factorising$^+$ at all positive depths for subgroups $G\leq \Aut(T)^+$ satisfying the property \ref{IPV1} and the hypothesis \ref{Hypothese HV1}(Definition \ref{definition H v1}). This requires some formalism that we now introduce. Let $V=V_0\sqcup V_1$ be a bipartition of $T$ such that every edge of $T$ contains exactly one vertex in each $V_i$. For every subtree $\mathcal{T}\subseteq T$, we set
	$$Q_\mathcal{T}=\{v\in V_0: B_T(v,2) \subseteq \mathcal{T}\}$$
	and we define $\mathfrak{T}_{V_1}$ as follows:
	\begin{enumerate}
		\item $\mathfrak{T}_{V_1}\lb 0\rb=\{ B_T(v,1): v\in V_1\}.$
		\item For every $l\in \N$ such that $l\geq 0$, we define iteratively
		\begin{equation*}\label{page Tv1}
		\begin{split}
		\mathfrak{T}_{V_1}\lb l+1\rb=\{ \mathcal{T}\subseteq T: \exists \mathcal{R}\in \mathfrak{T}_{V_1}\lb l\rb, \exists w\in (V&(\mathcal{R})- Q_{\mathcal{R}})\cap V_0 \\
		&\mbox{ s.t. }  \mathcal{T}=\mathcal{R}\cup B_T(w,2)\}.
		\end{split}
		\end{equation*}
		\item We set $\mathfrak{T}_{V_1} =\bigsqcup_{l\in \N}\mathfrak{T}_{V_1}\lb l\rb$.
	\end{enumerate}
	For every closed subgroup $G\leq\Aut(T)^+$ we let
	\begin{equation*}\label{page de SV1 yesssouille}
	\mathcal{S}_{V_1}=\{\Fix_G(\mathcal{T}): \mathcal{T}\in \mathfrak{T}_{V_1}\}.
	\end{equation*}
	\begin{definition}\label{definition H v1}
		A group $G\leq\Aut(T)^+$ is said to satisfy the hypothesis \ref{Hypothese HV1} if for all $\mathcal{T},\mathcal{T'}\in \mathfrak{T}_{V_1}$ we have that
		\begin{equation}\tag{$H_{V_1}$}\label{Hypothese HV1}
		\Fix_G(\mathcal{T}')\subseteq \Fix_G(\mathcal{T})\mbox{ if and only if }\mathcal{T}\subseteq \mathcal{T}'. 
		\end{equation}
	\end{definition}
	\noindent If $G\leq \Aut(T)^+$ is a closed non-discrete unimodular subgroup of $\Aut(T)^+$ satisfying the hypothesis \ref{Hypothese HV1}, Lemma \ref{la forme des Sl pour IPV1} below ensures that $\mathcal{S}_{V_1}$ is a generic filtration of $G$ and that
	$$\mathcal{S}_{V_1}\lb l\rb =\{\Fix_G(\mathcal{T}): \mathcal{T}\in \mathfrak{T}_{V_1}\lb l \rb\}.$$
	Theorem \ref{theorem Ipv1 letter} states the following.
	\begin{theorem*}\label{Theorem pour IPV1}
		Let $T$ be a locally finite tree and $G\leq \Aut(T)^+$ be a closed non-discrete unimodular subgroup satisfying the hypothesis \ref{Hypothese HV1} and the property \ref{IPV1}. Then, the generic filtration $\mathcal{S}_{V_1}$ factorises$^+$ at all depths $l\geq 1$. 
	\end{theorem*}
	
	\subsection{Preliminaries}
	Let $T$ be a locally finite tree and let $V=V_0\sqcup V_1$ be a bipartition of $T$ such that every edge of $T$ contains exactly one vertex in each $V_i$. The purpose of the present section is to describe further the elements of $\mathfrak{T}_{V_1}$. For every subtree $\mathcal{T}\subseteq T$, we associated a set $Q_\mathcal{T}=\{v\in V_0: B_T(v,2) \subseteq \mathcal{T}\}$. The purpose of the following two lemmas is to give a characterisation of the elements $\mathcal{T}\in \mathfrak{T}_{V_1}$ in terms of their corresponding sets $Q_\mathcal{T}$. 
	\begin{lemma}\label{the form of trees in the family}
		The elements of $\mathfrak{T}_{V_1}$ satisfy the following:
		\begin{enumerate}[label=(\roman*)]
			\item\label{item 1 le lemme pour la famille d'indep}  Every $\mathcal{T}\in\mathfrak{T}_{V_1}$ is a complete finite subtree of $T$ with leaves in $V_0$.
			\item\label{item 2 le lemme pour la famille d'indep} For every $\mathcal{T}\in \mathfrak{T}_{V_1}- \mathfrak{T}_{V_1}\lb 0 \rb$ we have that $\mathcal{T}= \bigcup_{v\in Q_\mathcal{T}}B_T(v,2)$. 
			\item\label{item 3 le lemme pour la famille d'indep}  For every $\mathcal{T}\in \mathfrak{T}_{V_1}$, we have $\mathcal{T}\in \mathfrak{T}_{V_1}\lb l\rb$ if and only if $\modu{Q_\mathcal{T}}=l$.
		\end{enumerate}
	\end{lemma}
	\begin{proof}
		Since each element of $\mathfrak{T}_{V_1}$ belongs to some $\mathfrak{T}_{V_1}\lb l\rb$ for some $l\in \N$, in order to show \ref{item 3 le lemme pour la famille d'indep} it is enough to show that  $\forall \mathcal{T}\in \mathfrak{T}_{V_1}\lb l\rb$, $\modu{Q_\mathcal{T}}=l$. We prove \ref{item 1 le lemme pour la famille d'indep}, \ref{item 2 le lemme pour la famille d'indep} and that $\modu{Q_\mathcal{T}}=l$ for every $\mathcal{T}\in \mathfrak{T}_{V_1}\lb l\rb$ by induction on $l$. If $l=0$, $\mathcal{T}=B_T(v,1)$ for some $v\in V_1$. Hence, $\mathcal{T}$ is a complete finite subtree with leaves in $V_0$ and $\modu{Q_\mathcal{T}}=0$. Similarly, if $l=1$, $\mathcal{T}=B_T(v,2)$ for some $v\in V_0$. In particular, $\mathcal{T}$ is a complete finite subtree of $T$ with leaves in $V_0$ and since $Q_\mathcal{T}=\{v\}$ we have that $\mathcal{T}= \bigcup_{v\in Q_\mathcal{T}}B_T(v,2)$ and $\modu{Q_\mathcal{T}}=1$. If $l\geq 2$, by construction, there exist $\mathcal{R}\in \mathfrak{T}_{V_1}\lb l-1\rb$ and $w\in (V(\mathcal{R})-Q_\mathcal{R})\cap V_0$ such that $\mathcal{T}= \mathcal{R}\cup B_T(w,2)$. By the induction hypothesis we have:
		\begin{enumerate}
			\item[\rm (1)] $\mathcal{R}$ is a finite complete subtree of $T$ with leaves in $V_0$.
			\item[\rm(2)] $\mathcal{R}= \bigcup_{v\in Q_{\mathcal{R}}} B_T(v,2).$ 
			\item[\rm (3)] $\modu{Q_{\mathcal{R}}}=l-1.$
		\end{enumerate}
		Since $\mathcal{T}= \mathcal{R}\cup B_T(w,2)$, (1) implies that $\mathcal{T}$ is a complete finite subtree with leaves in $V_0$ which proves \ref{item 1 le lemme pour la famille d'indep}. On the other hand, (2) implies that $Q_{\mathcal{R}} \cup \{w\}\subseteq Q_\mathcal{T}$ and therefore that $\mathcal{T}\subseteq \bigcup_{v\in Q_\mathcal{T}}B_T(v,2)$. The reverse inclusion follows trivially from the definition of $Q_\mathcal{T}$ which proves \ref{item 2 le lemme pour la famille d'indep}. Now, let $w'\in Q_{\mathcal{T}} - Q_{\mathcal{R}}$. To prove that $\modu{Q_\mathcal{T}}=l$ we have to prove that $w'=w$. Since $w'\in Q_{\mathcal{T}} - Q_{\mathcal{R}}$, there exists $u\in B_T(w',2)\cap V_0$ such that $u\not \in V(\mathcal{R})$. Moreover, since the leaves of $\mathcal{T}$ belongs to $V_0$ and since the distance between two vertices of $V_0$ is even, notice that $d_T(u,\mathcal{R})=2$. On the other hand, there exists a unique vertex $x\in V(\mathcal{R})-Q_\mathcal{R}$ such that $u\in B_T(x,2)$. Since $\mathcal{T}= \mathcal{R}\cup B_T(w,2)$ this proves that $w=x=w'$, that $Q_\mathcal{T}=Q_{\mathcal{R}}\sqcup \{w\}$ and therefore that $\modu{Q_\mathcal{T}}=l$. 
	\end{proof}
	\begin{lemma}\label{alternative form of S}
		Let $\mathcal{T}= \bigcup_{v\in Q}B_T(v,2)$ for some finite set $Q\subseteq V_0$ of order $l\geq 1$. Then $\mathcal{T}\in \mathfrak{T}_{V_1}\lb l\rb$ if and only if $\mbox{Con}(Q)\cap V_0=Q$ where $\mbox{Con}(Q)$ denotes the convex hull of $Q$ in $T$.
	\end{lemma}
	\begin{proof}
		Suppose first that $\mathcal{T}\in \mathfrak{T}_{V_1}\lb l\rb$. The definition of $Q_{\mathcal{T}}$ implies that $Q\subseteq Q_\mathcal{T}$ and Lemma \ref{the form of trees in the family} ensures that $\modu{Q_\mathcal{T}}=l$ which proves that $Q=Q_\mathcal{T}$. We prove that $\mbox{Con}(Q_\mathcal{T})\cap V_0=Q_{\mathcal{T}}$ by induction on $l\geq 1$. If $l=1$ the result is trivial. Suppose that $l\geq 2$. By construction, there exists $\mathcal{R}\in \mathfrak{T}_{V_1}\lb l-1\rb$ and $w\in (V(\mathcal{R})\cap V_0)-Q_{\mathcal{R}}$ such that $\mathcal{T}= \mathcal{R}\cup B_T(w,2)$. Furthermore, since $l-1\geq 1$, Lemma \ref{the form of trees in the family} ensures that  $\mathcal{R}= \bigcup_{v\in Q_{\mathcal{R}}}B_T(v,2)$. Since $T$ is a tree and since $w\in \mathcal{R}$, there exists a unique $u\in Q_{\mathcal{R}}$ such that $d_T(u,w)=2$ and we have that $\mbox{Con}(Q_\mathcal{T})= \mbox{Con}(Q_{\mathcal{R}})\cup \lb u,w\rb$. Finally, notice that $\lb u,w\rb\cap V_0=\{u,w\}$ which proves that $\mbox{Con}(Q_\mathcal{T})\cap V_0= Q_{\mathcal{R}}\cup \{w\}= Q_\mathcal{T}$. 
		
		Now, we show by induction on $l$ that $\mathcal{T}= \bigcup_{v\in Q}B_T(v,2)\in \mathfrak{T}_{V_1}\lb l \rb$ if $Q\subseteq V_0$ is a set of order $l$ such that $\mbox{Con}(Q)\cap V_0=Q$. If $l=1$, the result is trivial. Suppose that $l\geq 2$, choose any vertex $v\in Q$ and let $Q^n=\{w\in Q\lvert d_T(w,v)=2n\}$. Since $Q$ is finite there exists $N\in \N$ such that $Q^N\not = \es$ but $Q^n=\es$ for every $n\gneq N$. In particular, notice that $Q=\bigsqcup_{n\leq N} Q^n$. For every $n\leq N$, we let $S_n=\bigcup_{w\in Q^n, s\leq n}B_T(w,2)$, $l_n=\big\lvert\bigsqcup_{s\leq n}Q^{s-1}\big\rvert$ and we notice by induction on $n$ that $S_n\in \mathfrak{T}_{V_1}\lb l_n\rb$. Notice that the result is trivial for $n=1$, so let $n\geq 2$ and let $Q^n=\{v_1,...,v_{r_n}\}$. Since $\mbox{Con}(Q)\cap V_0=Q$, for all $w\in Q^n$ there exists $v_w\in Q^{n-1}$ such that $d_T(w,v_w)=2$. In particular, starting from our induction hypothesis we obtain iteratively for every $0\leq t\leq r_n$ that 
		$S_n\cup \big( \bigcup_{i\leq t}B_T(v_i,2)\big)\in \mathfrak{T}_{V_1}\lb l_n+ t\rb$. The result follows since $l_N=l$ and $S_N=\mathcal{T}$.
	\end{proof}
	This description allows one to prove the following result.
	\begin{lemma}\label{invariance de TTT par rapport a l'action de AutT plus}
		Let $\mathcal{T}\in \mathfrak{T}_{V_1}\lb l \rb$ and $g\in \Aut(T)^+$, then we have $g\mathcal{T}\in \mathfrak{T}_{V_1}\lb l\rb$. 
	\end{lemma} 
	\begin{proof}
		If $l=0$, $\mathcal{T}=B_T(v,1)$ some $v\in V_1$. Furthermore, $gB_T(v,1)=B_T(gv,1)$ and since the element of $\Aut(T)^+$ are type-preserving, $gv\in V_1$ which proves that $g\mathcal{T}\in \mathfrak{T}_{V_1}\lb 0\rb$.  If $l\geq 1$, Lemma \ref{alternative form of S} ensures that $\mbox{Con}(Q_\mathcal{T})\cap V_0=Q_\mathcal{T}$. It is clear from the definition that $Q_{g\mathcal{T}}=gQ_{\mathcal{T}}$ and since $g$ is a type-preserving automorphism of a tree we have $\mbox{Con}(gQ_\mathcal{T})\cap V_0=gQ_\mathcal{T}$. In particular, Lemma \ref{alternative form of S} ensures that $g\mathcal{T}\in \mathfrak{T}_{V_1}\lb l \rb$.
	\end{proof}
	\subsection{Factorisation of the generic filtration $\mathcal{S}_{V_1}$}
	The purpose of this section is to prove Theorem \ref{theorem Ipv1 letter}. We adopt the same notations as in the above sections. 
	\begin{lemma}\label{la forme des Sl pour IPV1}
		Let $G\leq \Aut(T)^+$ be a closed non-discrete unimodular subgroup satisfying the hypothesis \ref{Hypothese HV1}. Then, $\mathcal{S}_{V_1}$ is a generic filtration of $G$ and
		$$\mathcal{S}_{V_1}\lb l\rb =\{\Fix_G(\mathcal{T}): \mathcal{T}\in \mathfrak{T}_{V_1}\lb l \rb\}\q \forall l\in \N.$$
	\end{lemma}
	\begin{proof}
		For every $\mathcal{T}\in \mathfrak{T}_{V_1}$ notice that $g\Fix_{G}(\mathcal{T})g^{-1}=\Fix_{G}(g\mathcal{T})$ $\forall g\in G$ and therefore that $\mathcal{C}(\Fix_{G}(\mathcal{T}))=\{\Fix_{G}(g\mathcal{T}): g\in G\}$. In particular, for every $\mathcal{T},\mathcal{T}'\in \mathfrak{T}_{V_1}$ we have that $\mathcal{C}(\Fix_{G}(\mathcal{T}'))\leq \mathcal{C}(\Fix_{G}(\mathcal{T}))$ if and only if there exists $g\in G$ such that $\Fix_{G}(\mathcal{T})\subseteq \Fix_{G}(g\mathcal{T}')$. Since $G$ satisfies the hypothesis \ref{Hypothese HV1}, this implies that $\mathcal{C}(\Fix_{G}(\mathcal{T}'))\leq \mathcal{C}(\Fix_{G}(\mathcal{T}))$ if and only if there exists some $g\in G$ such that $g\mathcal{T}'\subseteq \mathcal{T}$. On the other hand, Lemma \ref{invariance de TTT par rapport a l'action de AutT plus} ensures that $\mathfrak{T}_{V_1}$ is stable under the action of $G$. In particular, for every strictly increasing chain $C_0\lneq C_1\lneq ...\lneq C_{n-1}\lneq C_{n}$ of elements of $\mathcal{F}_{\mathcal{S}_{V_1}}$ there exists a strictly increasing chain $\mathcal{T}_0\subsetneq\mathcal{T}_1\subsetneq...\subsetneq\mathcal{T}_{n-1}\subsetneq\mathcal{T}_n$ of elements of $\mathfrak{T}_{V_1}$ such that $C_t=\mathcal{C}(\Fix_G(\mathcal{T}_t))$. On the other hand, for every strictly increasing chain $\mathcal{T}_0\subsetneq\mathcal{T}_1\subsetneq...\subsetneq\mathcal{T}_{n}\subseteq\mathcal{T}$ of elements of $\mathfrak{T}_{V_1}$ contained in $\mathcal{T}$ we can build a strictly increasing chain $\mathcal{C}(\Fix_G(\mathcal{T}_0))\lneq\mathcal{C}(\Fix_G(\mathcal{T}_1))\lneq...\lneq\mathcal{C}(\Fix_G(\mathcal{T}_{n}))\lneq\mathcal{C}(\Fix_G(\mathcal{T}))$ of elements of $\mathcal{F}_{\mathcal{S}_{V_1}}$. This proves that the height of $\mathcal{C}(\Fix_G(\mathcal{T}))$ is the maximal length of a strictly increasing chain of elements of $\mathfrak{T}_{V_1}$ contained in $\mathcal{T}$.
		The result therefore follows from the following observation: Lemma \ref{the form of trees in the family} ensures that every maximal strictly increasing chain of elements of $\mathfrak{T}_{V_1}$ contained in $\mathcal{T}$ is of the form $\mathcal{T}_0\subsetneq \mathcal{T}_1\subsetneq...\subsetneq\mathcal{T}_{l-1}\subsetneq\mathcal{T}$ where $\mathcal{T}_t\in \mathfrak{T}_{V_1}\lb t \rb$.
	\end{proof}
	The following lemma shows that independence properties such as ${\rm IP}_1$ or ${\rm IP}_{V_1}$ can be realised as factorisation properties for the pointwise stabilisers of particular families of complete finite subtrees. 
	\begin{lemma}\label{independence on trees}
		Let $T$ be a locally finite tree, $G\leq \Aut(T)$, $\mathcal{A}$ be a family of finite subtrees of $T$ with at least two vertices such that
		\begin{equation*}
		\Fix_G( \mathcal{P})=\prod_{v\in \partial  \mathcal{P}} \Fix_G(T( \mathcal{P},v))\mbox{ } \q\forall \mathcal{P}\in \mathcal{A}
		\end{equation*}  
		where $\partial  \mathcal{P}$ denotes the set of leaves of $ \mathcal{P}$ and $T(\mathcal{P},v)$ denotes the half-tree $T( \mathcal{P},v)=\{w\in V: d_T(w, \mathcal{P})< d_T(w,v)\}$ and let $\mathcal{T}$ be a non-empty complete finite subtree of $T$ such that for every $v\in \partial \mathcal{T}$, there exists a subtree $\mathcal{T}_v\in \mathcal{A}$ with $v\in \mathcal{T}_v\subseteq  \mathcal{T}$. Then, we have that
		\begin{equation*}
		\Fix_G(\mathcal{T})=\prod_{v\in \partial \mathcal{T}} \Fix_G(T(\mathcal{T},v)).
		\end{equation*}
	\end{lemma}
	\begin{proof}
		For every subset $F\subseteq V$, we denote by $F^c$ the complement of $F$ in $V$. First, notice that for every two distinct leaves $v,v'\in \partial \mathcal{T}$, the support of the elements of $\Fix_G(T(\mathcal{T},v))$ and $\Fix_G(T(\mathcal{T},v'))$ are disjoint. In particular, the elements of $\Fix_G(T(\mathcal{T},v))$ and of $\Fix_G(T(\mathcal{T},v'))$ commute with each other and $\prod_{v\in \partial \mathcal{T}} \Fix_G(T(\mathcal{T},v))$ is a well defined subgroup of $G$. On the other hand, $\forall v\in \partial \mathcal{T}$ we have $\Fix_G(\mathcal{T})\supseteq  \Fix_G(T(\mathcal{T},v))$ and therefore that $ \Fix_G(\mathcal{T})\supseteq  \prod_{v\in \partial \mathcal{T}}\Fix_G(T(\mathcal{T},v))$. 
		In order to prove the other inclusion let $g\in \Fix_G(\mathcal{T})$ and let $\partial\mathcal{T}=\{v_1,...,v_n\}$. 
		The hypothesis on $G$ imply the existence of a subtree $\mathcal{T}_{1}\in \mathcal{A}$ such that $v_1\in \mathcal{T}_{1}\subseteq\mathcal{T}$. Furthermore, since $v_1\in \partial \mathcal{T}$ and $\mathcal{T}_{1}\subseteq \mathcal{T}$, we observe that $v_1\in \partial \mathcal{T}_{1}$ and $T(\mathcal{T}_1,v_1)= T(\mathcal{T},v_1)$. Furthermore, since $\mathcal{T}_{1}\subseteq \mathcal{T}$ and $g\in \Fix_G(\mathcal{T})$, notice that $g\in \Fix_G(\mathcal{T}_1)$. Our hypothesis on $G$ ensures the existence of some $h_{1}\in \Fix_G(T(\mathcal{T}_1,v_1))$ and $$g_1\in \prod_{v\in \partial \mathcal{T}_1-\{v_1\}}\Fix_G(T(\mathcal{T}_1,v))\subseteq \Fix_G(T(\mathcal{T}_1,v_1)^c)$$
		such that $g=h_{1}g_1$. Since $g\in \Fix_G(\mathcal{T})$ and $\Fix_G(T(\mathcal{T}_1,v_1))\subseteq \Fix_G(\mathcal{T})$, this decomposition implies that $g_1\in \Fix_G(\mathcal{T})\cap \Fix_G(T(\mathcal{T},v_1)^c)$. Proceeding iteratively, we prove the existence of some $h_i\in \Fix_G(T(\mathcal{T},{v_i}))$ and some $g_i\in \Fix_G(\mathcal{T})\cap\big( \bigcap_{r\leq i}\Fix_G(T(\mathcal{T},v_r)^c)\big)$ such that $g_{i-1}=h_{i} g_i$. To see that $g_i \in \bigcap_{r\leq i}\Fix_G(T(\mathcal{T},v_r)^c)$, notice by induction that $$g_{i-1}\in \bigcap_{r\leq i-1}\Fix_G(T(\mathcal{T},v_r)^c),$$ that $h_{i} \in \Fix_G(T(\mathcal{T},{v_i}))$ and that $\Fix_G(T(\mathcal{T},v_j))\subseteq \Fix_G(T(\mathcal{T},v_i)^c)$  $\forall i\not=j$. This implies that $h_i\in \bigcap_{r\leq i-1}\Fix_G(T(\mathcal{T},v_r)^c)$ and therefore that $g_i \in \bigcap_{r\leq i}\Fix_G(T(\mathcal{T},v_r)^c)$. The result follows since we have by construction that $g= h_1h_2...h_ng_n$, that $h_i\in \Fix_G(T(\mathcal{T},v_i))$ and that $$g_n\in \Fix_G(\mathcal{T})\cap \bigg(\bigcap_{v\in \partial \mathcal{T}}\Fix_G(T(\mathcal{T},v)^c)\bigg)=\Fix_G(T)=\{1_{\Aut(T)}\}.$$
	\end{proof}
	This result provides an alternative proof of Lemma \ref{lemma alternative definition of tites indep prop} and allows one to prove the following result which is key to the proof of Theorem \ref{theorem Ipv1 letter}.
	\begin{proposition}
		Let $G\leq \Aut(T)^+$ be a subgroup satisfying the property \ref{IPV1}. Then, for every $\mathcal{T}\in \mathfrak{T}_{V_1}$, we have
		\begin{equation*}
		\Fix_G(\mathcal{T})=\prod_{v\in \partial \mathcal{T}} \Fix_G(T(\mathcal{T},v)).
		\end{equation*}  
	\end{proposition}
	\begin{proof}
		If $\mathcal{T}\in \mathfrak{T}_{V_1}\lb 0\rb$, there exists $w\in V_1$ such that $\mathcal{T}=B_T(w,1)$. Notice that $\partial \mathcal{T}=\{v\in V: d_T(v,w)=1\}$ and that $T(\mathcal{T},v)=T(w,v)$  $\forall v\in \partial \mathcal{T}$. Since $G$ satisfies the property \ref{IPV1}, we obtain, as desired, that
		\begin{equation*}
		\begin{split}
		\Fix_G(\mathcal{T})&= \prod_{v\in \partial \mathcal{T}}\Fix_G(T(\mathcal{T},v)).
		\end{split}
		\end{equation*}
		Now, let $\mathcal{A}=\mathfrak{T}_{V_1}\lb 0\rb$ and notice from Lemma \ref{the form of trees in the family}\ref{item 2 le lemme pour la famille d'indep} that the hypothesis of Lemma \ref{independence on trees} are satisfied for every $\mathcal{T}\in\mathfrak{T}_{V_1}$. The result follows.
	\end{proof}
	\begin{lemma}\label{the key lemma for fertility of groups with IP V1 (build)}
		Let $G\leq \Aut(T)^+$ and suppose that 
		\begin{equation*}
		\Fix_G(\mathcal{T})=\prod_{v\in \partial \mathcal{T}} \Fix_G(T(\mathcal{T},v))\q \forall \mathcal{T}\in \mathfrak{T}_{V_1}.
		\end{equation*}
		Then, for every integers $ l,l'\geq 1$ with $l'\geq l$, $\forall \mathcal{T}\in \mathfrak{T}_{V_1} \lb l \rb$ and $\forall \mathcal{T}'\in \mathfrak{T}_{V_1}\lb l'\rb$ such that $\mathcal{T}\not \subseteq \mathcal{T}'$, there exists a subtree $\mathcal{R}\subseteq\mathcal{T}$ such that $\mathcal{R}\in \mathfrak{T}_{V_1}\lb l-1\rb$ and
		\begin{equation*}
		\Fix_G(\mathcal{R})\subseteq \Fix_G(\mathcal{T}')\Fix_G(\mathcal{T}).
		\end{equation*}	 
	\end{lemma}
	\begin{proof}
		Let $l,l'\geq 1$ be such that $l'\geq l$, let $\mathcal{T}\in \mathfrak{T}_{V_1} \lb l \rb$ and let $\mathcal{T}'\in \mathfrak{T}_{V_1}\lb l'\rb$ be such that $\mathcal{T}\not \subseteq  \mathcal{T}'$. If $l=1$, $Q_\mathcal{T}=\{v_\mathcal{T}\}$ for some $v_\mathcal{T}\in V_0$. Since $\mathcal{T}\not \subseteq  \mathcal{T}'$, we have that $v_\mathcal{T}\not\in Q_{\mathcal{T}'}$. Hence $d_T(v_\mathcal{T},Q_{\mathcal{T}'})\geq 2$. In particular, there exists a unique vertex $w\in V_1\cap B_T(v_\mathcal{T},1)$ such that $\mathcal{T}'\subseteq T(w,v_\mathcal{T})\cup \{v_\mathcal{T}\}$. Let $\mathcal{R}= B_T(w,1)$ and notice that $\mathcal{R}\in \mathfrak{T}_{V_1}\lb 0\rb$. Our hypothesis on $G$ ensure that
		$$\Fix_G(\mathcal{R})= \prod_{v\in \partial \mathcal{R}} \Fix_G(T(\mathcal{R},v)).$$ 
		On the other hand, $\partial \mathcal{R}\subseteq \partial \mathcal{T} \cup \{v_\mathcal{T}\}$ and $T(\mathcal{R},v)=T(w,v)$ $\forall v\in \partial \mathcal{R}$. This proves that $\Fix_G(T(\mathcal{R},v))\subseteq \Fix_G(\mathcal{T})$ for every leaf $v\in \partial \mathcal{R} - \{v_\mathcal{T}\}$ and since $ \mathcal{T}'\subseteq T(w,v_\mathcal{T})\cup\{v_\mathcal{T}\} =T(\mathcal{R},v_\mathcal{T})\cup\{v_\mathcal{T}\}$ we also have that $\Fix_G(T(\mathcal{R},v_\mathcal{T}))\subseteq \Fix_G(\mathcal{T}')$. This proves, as desired, that $$\Fix_G(\mathcal{R})=\prod_{v\in \partial \mathcal{R}} \Fix_G(T(\mathcal{R},v))\subseteq \Fix_G(\mathcal{T}')\Fix_G(\mathcal{T}).$$ 
		
		If $l\geq 2$, we have that $\modu{Q_\mathcal{T}}\geq 2$ and since $\mathcal{T}\not \subseteq \mathcal{T}'$, $Q_\mathcal{T}\not \subseteq Q_{\mathcal{T}'}$. In particular, there exists $v_\mathcal{T}\in Q_\mathcal{T}$ such that $d_T(v_\mathcal{T},Q_{\mathcal{T}'})=\max\{d_T(v,Q_{\mathcal{T}'}): v\in Q_\mathcal{T}\}$. On the other hand, since $Q_\mathcal{T},Q_{\mathcal{T}'}\subseteq V_0$ the distance $d_T(v_\mathcal{T},Q_{\mathcal{T}'})$ must be even, hence $d_T(v_\mathcal{T},Q_{\mathcal{T}'})\geq 2$. Now, let $Q_{\mathcal{R}}=Q_\mathcal{T}- \{ v_\mathcal{T}\}$. Notice that $\mbox{Con}(Q_{\mathcal{R}})\cap V_0=Q_{\mathcal{R}}$. Indeed, suppose for a contradiction that there exists $w\in \mbox{Con}(Q_{\mathcal{R}})\cap V_0$ such that $w\not\in Q_{\mathcal{R}}$. Since $\mathcal{T}\in \mathfrak{T}_{V_1}\lb \modu{Q_\mathcal{T}}\rb$, Lemma \ref{alternative form of S} guarantees that $\mbox{Con}(Q_\mathcal{T})\cap V_0=Q_\mathcal{T}$. In particular, we observe that $w\in Q_\mathcal{T}-Q_{\mathcal{R}}$ and since $Q_{\mathcal{R}}= Q_\mathcal{T} -\{v_\mathcal{T}\}$ this implies that $w=v_\mathcal{T}$.  In particular, $v_\mathcal{T}\in\mbox{Con}(Q_{\mathcal{R}})\cap V_0$, there exists $w_1,w_2\in Q_\mathcal{T}- \{v_\mathcal{T}\}$ such that $(\lb w_1, w_2\rb\cap V_0)-\{w_1,w_2\} =\{v_\mathcal{T}\}$. If $l=2$, this leads to a contradiction since $Q_\mathcal{T}$ contains only two elements. On the other hand, if $l\geq 3$, we obtain a contradiction with our choice of $v_\mathcal{T}$. Indeed, for $i=1,2$ we have that $d_T(w_i, Q_{\mathcal{T}'})\leq d_T(v_\mathcal{T},Q_{\mathcal{T}'})$. Since $v_\mathcal{T}\not \in \mathcal{T}'$, this implies the existence of $\tilde{w}_i\in Q_{\mathcal{T}'}\cap T(w_i,v_\mathcal{T})$ and since there exists a unique simple path between $\tilde{w}_1$ and $\tilde{w}_2$, we obtain that $v_\mathcal{T}\in\lb \tilde{w}_1, \tilde{w}_2 \rb \cap V_0 \subseteq \mbox{Con}(Q_{\mathcal{T}'}) \cap V_0=Q_{\mathcal{T}'}\cap V_0$. This is a contradiction since $v_\mathcal{T}\not \in Q_{\mathcal{T}'}$ which proves proves that $\mbox{Con}(Q_{\mathcal{R}})\cap V_0=Q_{\mathcal{R}}$. In particular, Lemma \ref{alternative form of S} ensures that $\mathcal{R}=\bigcup_{w\in Q_{\mathcal{R}}}B_T(w,2)\in \mathfrak{T}_{V_1}\lb l-1\rb$. On the other hand, by choice of $v_\mathcal{T}$, we have $d_T(v_\mathcal{T}, w)\geq d_T(v,w)$ $\forall v\in Q_\mathcal{T}$, $\forall w\in Q_{\mathcal{T}'}$ which implies that  $\Fix_G(T(\mathcal{R},v_\mathcal{T}))\subseteq  \Fix_G(\mathcal{T}')$. On the other hand, $\partial \mathcal{R} \subseteq \partial \mathcal{T} \sqcup \{v_\mathcal{T}\}$ and $T(\mathcal{R},v)=T(\mathcal{T},v)$  $\forall v\in \partial \mathcal{R} \cap \partial \mathcal{T}$. In particular, since  $\Fix_G(\mathcal{T})=\prod_{v\in \partial \mathcal{T}} \Fix_G(T(\mathcal{T},v))$, we obtain that
		\begin{equation*}
		\begin{split}
		\Fix_G(\mathcal{R})&= \prod_{v\in \partial {\mathcal{R}}} \Fix_G(T(\mathcal{R},v))\subseteq \Fix_G(T(\mathcal{R},v_\mathcal{T}))\Fix_G(\mathcal{T})\subseteq \Fix_G(\mathcal{T}')\Fix_G(\mathcal{T}). 
		\end{split}
		\end{equation*}   
	\end{proof}
	\noindent The following result plays a similar role as Proposition \ref{G(W) for Sk-1} in Section \ref{application IPk}.
	\begin{proposition}\label{le lemme ultime de fertilité pour IPV1}
		Let $G\leq \Aut(T)^+$ be a closed subgroup satisfying the hypothesis \ref{Hypothese HV1} and the property \ref{IPV1}. Then, for ever integers $l,l'\geq 1$ such that $l'\geq l$, for every $U$ in the conjugacy class of an element of $\mathcal{S}_{V_1}\lb l \rb$ and every $V$ in the conjugacy class of an element of $\mathcal{S}_{V_1}\lb l'\rb$ such that $V\not \subseteq U$, there exists $W \in \mathcal{S}_{V_1}\lb l-1\rb$ such that $U\subseteq W\subseteq VU$.
	\end{proposition}
	\begin{proof}
		Lemma \ref{invariance de TTT par rapport a l'action de AutT plus} ensures that $\mathfrak{T}_{V_1}$ is stable under the action of $G$. In particular, by Lemma \ref{la forme des Sl pour IPV1} there exists $\mathcal{T}\in \mathfrak{T}_{V_1}\lb l\rb$ and $\mathcal{T}'\in \mathfrak{T}_{V_1}\lb l'\rb$ such that $U=\Fix_G(\mathcal{T})$ and $V=\Fix_G(\mathcal{T}')$. Furthermore, since $G$ satisfies the hypothesis \ref{Hypothese HV1} and since $V\not \subseteq U$ we have $\mathcal{T}\not \subseteq\mathcal{T}'$. In particular, since $G$ satisfies the property ${\rm IP}_{V_1}$, Lemma \ref{the key lemma for fertility of groups with IP V1 (build)} ensures the existence of $\mathcal{R}\in \mathcal{T}\lb l-1\rb$ such that $\mathcal{R}\subseteq \mathcal{T}$ and $\Fix_G(\mathcal{R})\subseteq \Fix_G(\mathcal{T}')\Fix_G(\mathcal{T})$. The result follows since, by  Lemma \ref{la forme des Sl pour IPV1},
		$W=\Fix_G(\mathcal{R})\in \mathcal{S}_{V_1}\lb l-1\rb$.
	\end{proof}
	We are finally ready to prove the main result of this section.
	\begin{proof}[Proof of Theorem \ref{theorem Ipv1 letter}]
		To prove that $\mathcal{S}_{V_1}$ factorises$^+$ at all depths $l\geq1$ we shall successively verify the three conditions of Definition \ref{definition olsh facto}. 
		
		First, we need to prove that for every $U$ in the conjugacy class of an element of $\mathcal{S}_{V_1}\lb l \rb$ and every $V$ in the conjugacy class of an element of $\mathcal{S}_{V_1}$ such that $V \not\subseteq U$, there exists a $W$ in the conjugacy class of an element of $\mathcal{S}_{V_1}\lb l -1\rb$ such that $U\subseteq W \subseteq V U$. Let $U$, $V$ be as above. If $V$ is conjugate to an element of $\mathcal{S}_{V_1}\lb l'\rb$ for some $l'\geq l$ the results follows directly from Proposition \ref{le lemme ultime de fertilité pour IPV1}. Therefore, we suppose that $l'< l$. Since $\mathfrak{T}_{V_1}$ is stable under the action of $G$ (Lemma \ref{invariance de TTT par rapport a l'action de AutT plus}), Lemma \ref{la forme des Sl pour IPV1} ensures the existence of $\mathcal{T}\in \mathfrak{T}_{V_1}\lb l \rb$ and $\mathcal{T}'\in \mathfrak{T}_{V_1}\lb l'\rb$ such that $U=\Fix_G(\mathcal{T})$ and $V=\Fix_G(\mathcal{T}')$. We have two cases.
		
		Either $\mathcal{T}'\subseteq \mathcal{T}$. In that case, we prove the existence of a finite subtree $\mathcal{R}\in \mathfrak{T}_{V_1}\lb l-1\rb$ such that $\mathcal{T}'\subseteq \mathcal{R}\subsetneq \mathcal{T}$ by induction on $l-l'$. If $l'-l=1$ we can take $\mathcal{R}=\mathcal{T}'$ and the result is trivial. On the other hand, if $l'-l\geq 2$, notice that $Q_\mathcal{T'}\subseteq Q_{\mathcal{T}'}$ and Lemma \ref{the form of trees in the family} ensures that $Q_{\mathcal{T}}- Q_\mathcal{T'}$ contains $l'-l$ vertices. If $Q_{\mathcal{T}'}$ is empty, there exists $v\in Q_\mathcal{T}$ such that $\mathcal{T}'\subseteq B_T(v,2)$. In particular, we let $\mathcal{P}=B_T(v,2)$ and notice that $\mathcal{P}\in\mathfrak{T}_{V_1}\lb l'+1\rb$ and that $\mathcal{T}'\subseteq \mathcal{P}\subsetneq \mathcal{T}$. If $Q_{\mathcal{T}'}$ is not empty, Lemma \ref{alternative form of S} ensures that $\mbox{Con}(Q_{\mathcal{T}})\cap V_0=Q_\mathcal{T}$ and at least one vertex $v\in Q_{\mathcal{T}}- Q_{\mathcal{T}'}$ satisfies that $d_T(v,Q_{\mathcal{T}'})=2$. We let $Q= Q_{\mathcal{T}'}\cup \{v\}$ and $\mathcal{P}=\bigcup_{w\in Q}B_T(w,2)$. Notice that $\mbox{Con}(Q)\cap V_0=Q$. In particular, Lemma \ref{alternative form of S} ensures that $\mathcal{P}\in \mathfrak{T}_{V_1}\lb l'+1\rb$ and we have, by construction, that $\mathcal{T}'\subseteq \mathcal{P}\subsetneq \mathcal{T}$. In both cases ($Q_\mathcal{T}$ is empty or not) our induction hypothesis ensures the existence of a finite subtree $\mathcal{R}\in \mathfrak{T}_{V_1}\lb l-1\rb$ such that $\mathcal{T}'\subseteq\mathcal{P}\subseteq \mathcal{R}\subsetneq \mathcal{T}$. In particular, we have $\Fix_G(\mathcal{R})\subseteq \Fix_G(\mathcal{T}')$ which implies, as desired, that
		$$\Fix_G(\mathcal{R})\subseteq \Fix_G(\mathcal{T}')\subseteq \Fix_G(\mathcal{T}')\Fix_G(\mathcal{T}).$$
		
		Or else, $\mathcal{T}'\not\subseteq \mathcal{T}$. If $Q_{\mathcal{T}'}=\es$, there exists a vertex $v\in V_0-Q_{\mathcal{T}}$ such that $\mathcal{T}'\subseteq B_T(v,2)$. In particular, we choose a set $Q\subseteq V_0$ of order $l$ containing $v$, such that $\mbox{Con}(Q)\cap V_0=Q$ and we set $\mathcal{P}=\bigcup_{w\in Q}B_T(w,2)$. Similarly, if $Q_{\mathcal{T}'}\not=\es$, since $\mbox{Con}(Q_{\mathcal{T}'})\cap V_0=Q_{\mathcal{T}'}$ by Lemma \ref{the form of trees in the family}, there exists a finite set $Q\subseteq V_0$ of order $l$ containing $Q_{\mathcal{T}'}$ and such that $\mbox{Con}(Q)\cap V_0=Q$. We set $\mathcal{P}=\bigcup_{w\in Q}B_T(w,2)$. In both cases ($Q_\mathcal{T}$ is empty or not), Lemma \ref{alternative form of S} ensures that $\mathcal{P}\in \mathfrak{T}_{V_1}\lb l \rb$ and we have by construction that  $\mathcal{T}'\subseteq \mathcal{P}$. In particular, Proposition \ref{le lemme ultime de fertilité pour IPV1} applied to $\Fix_G(\mathcal{T})$ and $\Fix_G(\mathcal{P})$ ensures the existence of a $W \in \mathcal{S}_{V_1}\lb l-1\rb$ such that $\Fix_G(\mathcal{T})\subseteq W\subseteq \Fix_G(\mathcal{P})\Fix_G(\mathcal{T})$. On the other hand, $\mathcal{T}'\subseteq \mathcal{P}$ which implies that $\Fix_G(\mathcal{P})\subseteq \Fix_G(\mathcal{T}')$. This proves the first condition.
		
		Next, we need to prove that $N_{G}(U, V)= \{g\in {G} : g^{-1}Vg\subseteq U\}$ is compact for every $V$ in the conjugacy class of an element of $\mathcal{S}_{V_1}$. Just as before, notice that $V=\Fix_G(\mathcal{T}')$ for some $\mathcal{T}'\in \mathfrak{T}_{V_1}\lb l'\rb$. Since $G$ satisfies the hypothesis \ref{Hypothese HV1} notice that
		\begin{equation*}
		\begin{split}
		N_{G}(U, V)&= \{g\in G : g^{-1}Vg\subseteq U\}= \{g\in G : g^{-1}\Fix_{G}(\mathcal{T}')g\subseteq \Fix_{G}(\mathcal{T})\}\\
		&= \{g\in G : \Fix_G(g^{-1}\mathcal{T}')\subseteq \Fix_G(\mathcal{T})\}= \{g\in G : g\mathcal{T}\subseteq \mathcal{T}'\}.
		\end{split}
		\end{equation*}	
		Since both $\mathcal{T}$ and $\mathcal{T}'$ are finite subtrees of $T$, this implies that $N_G(U,V)$ is a compact subset of $G$ which proves the second condition.
		
		Finally, we need to prove that for every $W$ in the conjugacy class of an element of $\mathcal{S}_{V_1}\lb l-1\rb$ such that $U\subseteq W$ we have
		\begin{equation*}
		W\subseteq N_G(U,U) =\{g\in G : g^{-1}Ug\subseteq U\}.
		\end{equation*} 
		The same reasoning as before ensures the existence of some $\mathcal{R}\in \mathfrak{T}_{V_1}$ such that $W = \Fix_G(\mathcal{R})$. On the other hand, since $U\subseteq W$ and since $G$ satisfies the hypothesis \ref{Hypothese HV1}, notice that $\mathcal{R}\subseteq \mathcal{T}$. We have several cases. If $l=1$, there exists vertices $v\in V_0$ and $w\in V_1$ such that $\mathcal{T}=B_T(v,2)$ and $\mathcal{R}=B_T(w,1)$. In particular, since $\mathcal{R}\subseteq \mathcal{T}$, this implies that $v\in \mathcal{R}$ and therefore that
		\begin{equation*}
		\begin{split}
		\Fix_G(\mathcal{R})&\subseteq \Fix_G(v)=\{h\in G : hB_T(v,2)\subseteq B_T(v,2)\}=\{h\in G: h\mathcal{T}\subseteq \mathcal{T}\}.
		\end{split}
		\end{equation*} 
		Similarly, if $l\geq 2$, notice that $Q_{\mathcal{T}}$ and $Q_{\mathcal{R}}$ are non-empty sets. Furthermore, since $\mathcal{R}\subseteq \mathcal{T}$, we have that $Q_{\mathcal{R}}\subseteq Q_{\mathcal{T}}$ and there exists a unique $v\in Q_{\mathcal{T}}-Q_{\mathcal{R}}$. On the other hand, Lemma \ref{alternative form of S}, implies that $d_T(v,Q_\mathcal{R})=2$ and since $\mathcal{R}=\bigcup_{w\in Q_{\mathcal{R}}}B_T(w,2)$ we observe that $\Fix_G(\mathcal{R})\subseteq \Fix_G(Q_{\mathcal{T}})$. Since $\mathcal{T}=\bigcup_{w\in Q_{\mathcal{T}}}B_T(w,2)$, this implies that
		\begin{equation*}
		\begin{split}
		\Fix_G(\mathcal{R})\subseteq \{h\in G: h\mathcal{T}\subseteq \mathcal{T}\}.
		\end{split}
		\end{equation*} 
		In both cases, since $G$ satisfies the hypothesis \ref{Hypothese HV1} we obtain that 
		\begin{equation*}
		\begin{split}
		\Fix_G(\mathcal{R})&\subseteq  \{h\in G : h\mathcal{T}\subseteq \mathcal{T}\}= \{h\in G : \Fix_G(\mathcal{T})\subseteq \Fix_G(h\mathcal{T})\}\\
		&=\{h\in G : h^{-1}\Fix_G(\mathcal{T})h\subseteq \Fix_G(\mathcal{T})\}= N_G(U,U)
		\end{split}
		\end{equation*} which proves the third condition.
	\end{proof}
	In particular, if $G\leq\Aut(T)^+$ is a closed non-discrete unimodular subgroup satisfying the hypothesis \ref{Hypothese HV1} and the property \ref{IPV1} Theorem \ref{la version paki du theorem de classification} provides a bijective correspondence between the equivalence classes of irreducible representations of $G$ at depth $l\geq1$ with seed $C\in \mathcal{F}_{\mathcal{S}_{V_1}}$ and the $\mathcal{S}_{V_1}$-standard representations of $\Aut_G(C)$. As a concrete example, the group $\Aut(T)^+$ of type preserving automorphisms of a $(d_0,d_1)$-semi-regular tree $T$ with $d_0,d_1\geq 3$ satisfies the hypothesis of Theorem \ref{theorem Ipv1 letter}. Other examples will be constructed in Section \ref{Application to universal groups of right-angled buildings}. 
	\subsection{Existence of $\mathcal{S}_{V_1}$-standard representations}\label{existence for IPV1}
	Let $T$ be a locally finite tree and let $V=V_0\sqcup V_1$ be a bipartition of $T$ such that every edge of $T$ contains exactly one vertex in each $V_i$. Let $\mathfrak{T}_{V_1}$ be the family of subtrees defined on page \pageref{page Tv1}, let $G$ be a closed non-discrete unimodular subgroup of $\Aut(T)^+$ and let $\mathcal{S}_{V_1}=\{\Fix_G(\mathcal{T}): \mathcal{T}\in \mathfrak{T}_{V_1}\}$. If $G\leq \Aut(T)^+$ is a closed unimodular subgroup  satisfying the hypothesis \ref{Hypothese HV1} and the property \ref{IPV1} we have shown that $\mathcal{S}_{V_1}$ is a generic filtration of $G$ that factorises$^+$ at all depths $l\geq 1$. In particular, Theorem \ref{la version paki du theorem de classification} ensures the existence of a bijective correspondence between the equivalence classes of irreducible representations of $G$ at depth $l\geq 1$ with seed $C\in \mathcal{F}_{\mathcal{S}_{V_1}}$ and the $\mathcal{S}_{V_1}$-standard representations of $\Aut_{G}(C)$. The following result treats the existence of $\mathcal{S}_{V_1}$-standard representations of $\Aut_{\mathcal{S}_{V_1}}(C)$ for all $C\in \mathcal{F}_{\mathcal{S}_{V_1}}$ at height $l\geq 1$ provided that $G$ satisfy some geometric property.
	\begin{proposition}
		Let $G\leq \Aut(T)^+$ be a closed non-discrete unimodular subgroup satisfying the hypothesis \ref{Hypothese HV1} and the property \ref{IPV1}, let $l\geq 1$ and let $\mathcal{T}\in \mathfrak{T}_{V_1}\lb l\rb$ be such that for every $\mathcal{R}\in \mathfrak{T}_{V_1} \lb l-1\rb$ with $\mathcal{R}\subseteq \mathcal{T}$ we have $$\Fix_G(\mathcal{R})\subsetneq\Stab_G(\mathcal{T})=\{g\in G: g\mathcal{T}\subseteq \mathcal{T}\}.$$
		Then, there exists a $\mathcal{S}_{V_1}$-standard representation of $\Aut_{\mathcal{S}_{V_1}}(\mathcal{C}(\Fix_G(\mathcal{T})))$.
	\end{proposition}
	\begin{proof}
		Let $C=\mathcal{C}(\Fix_G(\mathcal{T}))$. Lemma \ref{la forme des Sl pour IPV1} ensures that $C$ has height $l$ in $\mathcal{F}_{\mathcal{S}_{V_1}}$. Since $G$ satisfies the hypothesis \ref{Hypothese HV1} and as a consequence of Lemma \ref{la forme des Sl pour IPV1}, notice that $N_G(\Fix_G(\mathcal{T}),\Fix_G(\mathcal{T}))=\{g\in G : g\mathcal{T}\subseteq\mathcal{T}\}=\Stab_G(\mathcal{T})$, that $\Aut_{G}(C)\simeq \Stab_G(\mathcal{T})/\Fix_G(\mathcal{T})$ and that
		\begin{equation*}
		\begin{split}
		\tilde{\mathfrak{H}}_{\mathcal{S}_{V_1}}(\Fix_G(\mathcal{T}))&= \{W: \exists g\in G \mbox{ s.t. }gWg^{-1}\in \mathcal{S}_{V_1}\lb l-1 \rb \mbox{ and }\Fix_G(\mathcal{T}) \subseteq W \}\\
		&= \{\Fix_G(\mathcal{R}): \mathcal{R}\in \mathcal{S}_{V_1}\lb l-1\rb \mbox{ s.t. }\mathcal{R}\subseteq\mathcal{T}\}.
		\end{split}
		\end{equation*}
		Furthermore, the hypothesis on $G$ imply that $\Fix_G(\mathcal{T})\subsetneq\Fix_G(\mathcal{R})\subsetneq \Stab_G(\mathcal{T})$ for every $\mathcal{R}\in \mathfrak{T}_{V_1}\lb l-1\rb$ with $\mathcal{R}\subseteq \mathcal{T}$.
		
		We have two cases. If $\mathcal{T}\in \mathfrak{T}_{V_1}\lb 1\rb$, there exists $v\in V_0$ such that $\mathcal{T}=B_T(v,2)$ and every subtree $\mathcal{R}$ of $\mathcal{T}$ that belongs to $\mathcal{S}_{V_1}\lb 0\rb$ is of the form $B_T(w,1)$ for some $w\in \partial B_T(v,1)$. Let $\{\mathcal{T}_1,..., \mathcal{T}_d\}$ be the set of subtrees of $T$ of the form $\bigcup_{w\in \partial B_T(v,1)-\{u\}}B_T(w,1)$ for some $u\in \partial B_T(v,1)$. Notice that each element of $\Stab_G(\mathcal{T})=\Fix_G(v)$ permutes the vertices of $\partial B_T(v,1)$ and therefore permutes the elements of $\{\mathcal{T}_1,..., \mathcal{T}_d\}$. On the other hand, for every $i,j\in \{1,...,d\}$ with $i\not= j$, we have that $\mathcal{T}_i\cup  \mathcal{T}_j=\mathcal{T}$ and thanks to Proposition \ref{independence on trees}, $\Fix_G(\mathcal{T})\subsetneq \Fix_G(\mathcal{T}_i)\subsetneq \Stab_{G}(\mathcal{T})$. In particular, Proposition \ref{existence criterion} ensures the existence of an irreducible representation $\sigma$ of $\Aut_{G}(C)\simeq \Stab_G(\mathcal{T})/\Fix_G(\mathcal{T})$ without non-zero $p_{\Fix_G(\mathcal{T})}(\Fix_G(\mathcal{T}_i))$-invariant vector and therefore without non-zero $p_{\Fix_G(\mathcal{T})}(\Fix_G(\mathcal{R}))$-invariant vector for every subtree $\mathcal{R}\in \mathcal{S}_{V_1}\lb 0\rb$ such that $\mathcal{R}\subseteq \mathcal{T}$.
		
		If $\mathcal{T}\in \mathfrak{T}_{V_1}\lb l\rb$ for some $l\geq 2$, every subtree $\mathcal{R}\in \mathcal{S}_{V_1}\lb l-1\rb$ is such that $Q_{\mathcal{R}}\not=\es$. Let $\{\mathcal{T}_1,..., \mathcal{T}_d\}$ be the set of subtrees $\mathcal{R}$ of $\mathcal{S}_{V_1}\lb l-1\rb$ such that $\mathcal{R}\subseteq \mathcal{T}$ and notice that $Q_{\mathcal{T}_i}\subsetneq Q_{\mathcal{T}}$ $\forall i$. Furthermore, notice the elements of  $\Stab_G(\mathcal{T})$ permutes the elements of $Q_{\mathcal{T}}$ and therefore the elements of $\{\mathcal{T}_1,..., \mathcal{T}_d\}$. On the other hand, for every $i,j\in \{1,...,d\}$ with $i\not= j$, we have $\mathcal{T}_i\cup \mathcal{T}_j=\mathcal{T}$. In particular, Proposition \ref{existence criterion} ensures the existence of an irreducible representation $\sigma$ of $\Aut_{G}(C)\simeq \Stab_G(\mathcal{T})/\Fix_G(\mathcal{T})$ without non-zero $p_{\mathcal{T}}(\Fix_G(\mathcal{T}_i))$-invariant vectors and the result follows.
	\end{proof}
	\newpage
	
	\section{Universal groups of right-angled buildings}\label{Application to universal groups of right-angled buildings}
	The purpose of this section is to prove that the universal groups of certain semi-regular right-angled buildings embed as closed subgroups of the group $\Aut(T)^+$ of type-preserving automorphisms of a locally finite tree $T$ and that those subgroups satisfy the hypothesis of Theorem \ref{theorem Ipv1 letter} if the prescribed local action is $2$-transitive on panels. In particular, for every of such group we obtain a generic filtration factorising$^+$ at all depths and the machinery developed in Chapter \ref{Chapter Olshanskii's factor} applies.
	
	\subsection{Preliminaries}\label{section application right angle preliminaire}
	In this work, we realise buildings as $W$-metric spaces associated with a Coxeter system $(W,I)$. We now formalise these concepts and refer to \cite{AbramenkoBrown2008} for more details. We start with the notion of Coxeter group generalising the concept of finite group generated by reflexions in a real vector space. 
	\begin{definition}
		Let $I$ be a set and $\{m_{i,j}: i,j\in I\}$ be a set of elements $m_{i,j}\in \N_{\geq2}\cup \{\infty\}$ satisfying the following conditions:
		\begin{enumerate}[label=(\roman*)]
			\item $m_{i,i}=1$.
			\item $m_{i,j}\geq 2$ if $i\not=j$.
			\item $m_{i,j}=m_{j,i}$.
		\end{enumerate}
	 	The group $W=\langle I : (ij)^{m_{i,j}} \mbox{ for all finite }m_{i,j}\rangle$
	 	is called a \tg{Coxeter group}, $(W,I)$ is called a \tg{Coxeter system} and $\modu{I}$ is the rank of $W$. In addition, a Coxeter system is \tg{right-angled} if $m_{i,j}\in \{2,\infty\}$ for all $i,j\in I$.
	\end{definition}
	We now provide some important examples of Coxeter groups. 
	\begin{example}
		The dihedral group $D_{2m}$ with $m\in \N_{\geq 2}$. Consider a real vector space $V$ of dimension $2$. For each $m\in \N_{\geq 2}$, consider a pair $V_1,V_2$ of one dimensional subspaces with relative angle $\pi/m$ and the corresponding planar reflexions $\sigma_1,\sigma_2$ along these axes. The group generated by these reflexions is easily identified with the group of isometries of the plane stabilising a regular $m$-gone and is thus isomorphic to the dihedral group $D_{2m}$ of order $2m$. This proves that the dihedral group $D_{2m}$ is a Coxeter group of rank $2$ with presentation:
		$$D_{2m}\simeq \langle i_0,i_1: i_0^2, i_1^2,(i_0i_1)^m\rangle.$$ 
	\end{example}
	\begin{example}
		The infinite dihedral group $D_{\infty}$. Consider a real vector space $V$ of dimension $2$, $V_1$ and $V_2$ be two distinct and parallel lines in $V$ and let $\sigma_0,\sigma_1$ be the corresponding planar reflexions $\sigma_0,\sigma_1$ along these axes. The group generated by these reflections contains a subgroup of index $2$ isomorphic to $\Z$ generated by the translation of minimal step $i_0i_1$ and is therefore easily identified with the infinite dihedral group $D_{\infty}$. This proves that the infinite dihedral group $D_{\infty}$ is a right-angled Coxeter group of rank $2$ with presentation:
		$$D_{\infty}\simeq \langle i_0,i_1: i_0^2, i_1^2\rangle\simeq C_2*C_2.$$ 
	\end{example}
	\begin{example}\label{example group Sym n cox}
		The symmetric group $\Sym(n)$ with $n\geq 2$. Consider an orthonormal basis $\{e_1,...,e_{n}\}$ of $\R^{n}$. For each integer $l\in \{1,...,n-1\}$ we let $\sigma_l$ be the reflexion along the hyperplane perpendicular to $e_{l+1}-e_l$. The group generated by $\{\sigma_1,..., \sigma_{n-1}\}$ is easily seen as the group of permutations of the basis $\{e_1,...,e_{n}\}$ and is thus isomorphic to $\Sym(n)$. This proves that the symmetric group $\Sym(n)$ is a Coxeter group of rank $n-1$ with presentation:
		 $$\Sym(n) \simeq\langle i_1,i_1,..., i_{n-1}: i_l^2, (i_li_{l+1})^3, (i_{l}i_{l'})^2  \mbox{ if }\modu{l-l'}\geq2\rangle$$
	\end{example}
	Our next task is to define a building as a chamber complex equipped with a particular kind of ``metric'' ranging in a Coxeter group. In what follows, given a Coxeter system $(W,I)$ and an element $w\in W$, we denote by $l(w)$ the length of $w$ with respect to $I$. 
	\begin{definition}\label{building definition}
		Let $(W,I)$ be a Coxeter system. A \tg{building} $\Delta$ of type $(W,I)$ is a couple $(\Ch(\Delta), \delta)$ where $\Ch(\Delta)$ is a set called \tg{the set of chambers} of $\Delta$ and where $$\delta: \Ch(\Delta)\times \Ch(\Delta)\rightarrow W$$
		is a map satisfying the following conditions for all chambers $c,d\in \Ch(\Delta)$: 
		\begin{enumerate}[leftmargin=*, label=(\roman*)]
			\item $\delta(c,d)=1_W$ if and only if $c=d$.
			\item If $\delta(c,d)=w$ and $c'\in \Ch(\Delta)$ satisfies $\qq \delta(c',c)=i\in I$ we have that $\delta(c',d)\in \{w, iw\}$. If, in addition, $l(iw) = l(w) + 1$, then $\delta(c',d)=iw$.
			\item If $\delta(c,d)=w$, then for any $i\in I$, there is a chamber $c'\in \Ch(\Delta)$ such
			that $\delta(c',c)=i$ and $\delta(c',d)=iw$.
		\end{enumerate}
		The map  $\delta$ is called the \tg{Weyl distance} of $\Delta$. 
	\end{definition}
	\begin{definition}
		For every subset $J\subseteq I$ we let $W_J$ be the subgroup of $W$ generated by $J$ and we define the $J$-\textbf{residue} of $\Delta$ containing a chamber $c\in \Ch(\Delta)$ to be the set $$\mathcal{R}_J(c)=\{d\in \Ch(\Delta): \delta(c,d)\in W_J\}.$$ In addition, when $\mathcal{R}$ is a residue of $\Delta$, we denote by $\Ch(\mathcal{R})$ the set of chambers $c\in \Ch(\Delta)$ that belong to $\mathcal{R}$. An $\{i\}$-residue is called an $i$-\textbf{panel} and $\Delta$ is said to be \tg{thick} if each of its panels contains at least three chambers.
	\end{definition}
	 We now provide some examples. 
	\begin{example}
		Every Coxeter group $W$ with Coxeter system $(W,I)$ can be realised as a building of type $(W,I)$ and rank $\modu{I}$. To see this, consider the chamber set $\Ch(\Delta)=W$ and the Weyl distance
		$$\delta: W\times W \rightarrow W : (v,w) \mapsto v^{-1}w.$$
		We are left to show that $\delta$ satisfies the three conditions of the Definition \ref{building definition}. The first condition is trivially satisfied from the fact that $W$ is a group. Now, let $u,v,w\in W$ and notice that $$\delta(u,w)=u^{-1}w=u^{-1}vv^{-1}w=\delta(u,v)\delta (v,w).$$ The first part of the second condition follows directly and the second part follows from the properties of Coxeter groups. Finally the third condition follows from the fact that for any two chambers $v,w\in W$ and each $i\in I$, we have that $\delta(vi,v)=i$ and $\delta(vi,w)=iv^{-1}w=i\delta(v,w)$. This proves that $\delta$ is a Weyl distance and our claim follows. Now, notice from the definition of $\delta$ that two chambers $v,w\in W$ are $i$-adjacent if and only if $v=wi$. In particular, each panel of $\Delta$ is of the form $\{w,wi\}$ for some $w\in W$ and $i\in I$ so that this building is never thick. 
	\end{example}
	\begin{example}
		Every tree $T$ without leaves can be seen as a building $\Delta$ of type $(D_{\infty}, \{i_0,i_1\})$ with chamber set $\Ch(\Delta)=E$ and where panels are in bijections with vertices. To see this, consider a bipartition $V=V_0\sqcup V_1$ such that each edge of $T$ contains exactly one vertex in each of the $V_t$ with $t\in \{0,1\}$ and denote the map sending a vertex to its type by 
		$$\tau: V \rightarrow \{0,1\}; $$
		that is $\tau(v)=t$ if and only if $v\in V_t$. For every two edges $e,f\in E$, we consider the unique geodesic $\gamma=(v_0,...,v_n)$ in $T$ from $e$ to $f$ containing both edges. If $n=1$ we let $\delta(e,f)=1_{D_\infty}$ and if $n\geq 2$ we let $$\delta(e,f)= i_{\tau(v_1)}...i_{\tau(v_{n-1})}.$$
		It is clear from the definition of $\delta$ that for every edge $e\in E$, the set $$\mathcal{R}_{\{i_t\}}(e)=\{f\in E : \delta (e,f)\in W_{\{i_{t}\}}\}$$ is a set all the edges of $T$ containing the same vertex $v\in V$ of type $t$ as $e$. We are left to prove that $\delta: E\times E\rightarrow D_{\infty}$ satisfies the three conditions of the Definition \ref{building definition}. To prove the first condition, notice that two edges $e,f\in E$ coincide if and only if $n=1$ and hence $\delta(e,f)=1_{D_{\infty}}$. Now, let $e,f\in E$ be such that $\delta(e,f)=w\in D_{\infty}$ and let $e'\in E$ be such that $\delta(e',e)=i_t$ for some $t\in \{0,1\}$. Notice from the definition of $\delta$ that $e$ and $e'$ share a common vertex $v\in V_t$. There are two cases. Either, the geodesic containing both $v$ and $f$ contains $e$ or $e'$ and in that case, this edge is closer from $f$ than the other one so that $\delta(e',f)=i_tw$. Or this geodesic does not contain $e$ or $e'$ and the sequence of vertices $(v_1,...,v_{n-1})$ is the same for the couples $(e,f)$ and $(e',f) $ so that $\delta(e',f)=\delta(e,f)=w$. This implies that $\delta(e',f)\in \{w,i_tw\}$. Furthermore, if $l(i_tw)=l(w)+1$, notice that the geodesic from $v$ to $f$ contains $e$ and hence that $\delta(e',f)=i_tw$. The second condition follows. Finally, let $e,f\in E$ be such that $\delta(e,f)=w\in D_{\infty}$ and let $\{v_0,v_1\}$ be the set of vertices of $e$ (with $v_t\in V_t$). Since $T$ has no leaf, there exists an edge $e_t'$ containing $v_t$ and such that the distance in the tree between $e$ and $f$ is different from the distance between $e'_t$ and $f$. The definition of $\delta$ ensures that $\delta(e'_t,e)=i_t$ and that $\delta(e'_t,f)=i_t\delta(e,f)$ so that the third condition holds. Notice furthermore that this building is thick if and only if each vertex has degree at least $3$.
	\end{example}
	\begin{example}
		Let $\mathbb{K}$ be a finite field and consider a $\mathbb{K}$-vector space $V$ of dimension $n$. A flag in $V$ is a $(n+1)$-tuple $(V_0, V_1,V_2,...,V_n)$ of subspaces of $V$ such that each $V_i$ has dimension $i$ and $V_0\subsetneq V_1\subsetneq V_2\subsetneq ... \subsetneq V_n=V$. We recall from Example \ref{example group Sym n cox} that $\Sym(n)$ is a Coxeter group generated by the permutations $\sigma_i=(i,i+1)$ with $i\in \{1,...,n-1\}$. The set of flags of $V$ can be realised as the chamber set $\Ch(\Delta)$ of a building $\Delta$ of type $(\Sym(n), \{\sigma_1,..., \sigma_{n-1}\})$. The definition of the Weyl distance requires some preliminaries. We say that an ordered basis $(e_1,...,e_n)$ of $V$ represents the flag $F=(V_0, V_1,V_2,...,V_n)$ if $V_i=\langle e_1,...,e_i \rangle$ $\forall i\in \{1,...,n\}$. Now, let us prove that for any two flags $F$ and $F'$ there exists an ordered basis $(e_1,e_2,...,e_n)$ of $V$ and a permutation $\alpha \in \Sym(n)$ such that $(e_1,e_2,...,e_n)$ represents $F$ and $(e_{\alpha(1)},..., e_{\alpha(n)})$ represents $F'$. To see this, consider any ordered bases $(f_1,f_2,...,f_n)$ and $(f_1',f_2',...,f_n')$ representing respectively $F$ and $F'$ and consider the matrix $M$ whose $j$-th column consists of the components $f^i_j$ of the vector $f_j$ in the basis $\{f_1',...,f_n'\}$. Notice that we can multiply any column by an invertible element of $\mathbb{K}$ or replace a column by its sum with any linear combinations of the columns to its left without changing the flag $F$ that the matrix represents. Using this method, we change the ordered basis representing $F$ without changing $F$. Now, consider the first column whose $n$-th coordinate is non-zero renormalise it so that this coordinate equals $1$ and subtract it to each column to its right so that it is the only column of the matrix without a $0$ entry in $n$-th row. Ignoring this column for the rest of the algorithm we look at the first non-ignored column whose $(n-1)$-th coordinate is non-zero. Notice that such a column exists since $\{f_1,...,f_n\}$ is a basis of $V$. Renormalise this column so that this coordinate equals $1$ and subtract it to each non-ignored column to its right so that it becomes the only non-ignored column of the matrix without a $0$ entry in $(n-1)$-th row. We now ignore this column and iterate our procedure until each column is ignored. We end up with a new matrix $M'$. Now, let $(e_1,...,e_n)$ be the ordered basis of $V$ whose components in the basis $\{f_1',...,f_n'\}$ are given by the columns of $M'$. By construction $(e_1,...,e_n)$ represents the flag $F$. Furthermore, due to the form of the matrix, there exists a unique permutation $\alpha\in \Sym(n)$ such that the matrix $M''$ whose $i$-th column is the $\alpha(i)$-th column of $M'$ is upper triangular. In particular, $M''$ represents the same flag complex as $F'$. This proves as desired that $(e_1,e_2,...,e_n)$ represents $F$ and that $(e_{\alpha(1)},..., e_{\alpha(n)})$ represents $F'$. Now, notice that even if the ordered basis $(e_1,...,e_n)$ is in general not unique, the permutation $\alpha$ is only determined by $F$ and $F'$ since the matrix of coordinates of another basis $(f_1,...,f_n)$ representing $F$ in the basis $\{e_1,...,e_n\}$ must be upper triangular. In light of this remark we denote by $\alpha_{F,F'}$ the only permutation of $\Sym(n)$ with the above property. We define the $\Sym(n)$-metric of $\Delta$ as the map 
		$$\delta : \Ch(\Delta)\times \Ch(\Delta)\rightarrow \Sym(n): (F,F')\mapsto \alpha_{F,F'}^{-1}.$$
		In particular, notice that two flags $F=(V_0,V_1,...,V_n)$ and $F'=(V_0',V_1'...,V_n')$ are $\sigma_i$-adjacent for some $i\in \{1,...,n-1 \}$ if and only if $V_{t}=V_{t}'$ for all $t\not=i$ and $V_i\not=V'_i$. We are left to prove that $\delta$ satisfies the three conditions of the Definition \ref{building definition}. The first condition is clear from our definition of $\delta$. Now, let $F,G\in \Ch(\Delta)$ be such that $\delta(F,G)=\alpha^{-1}\in \Sym(n)$, $F'\in \Ch(\Delta)$ be such that $\delta(F',F)=\sigma_i$ and consider an ordered basis $(e_1,...,e_n)$ representing $F$ such that $(e_{\alpha(1)}, ... e_{\alpha(n)})$ represents $G$. We notice from our description of $\sigma_i$-adjacent chambers that $F'$ is represented by some ordered basis $(e_1,...,e_{i-1},e_{i+1}+qe_i,e_i,e_{i+2},...,e_n)$ for some scalar $q\in \mathbb{K}$. In particular, if $\alpha_{F,G}(i)>\alpha_{F,G}(i+1)$, notice that $\alpha_{F',G}=\alpha_{F,G}$ because we can change $e_{i+1}+qe_i$ to $e_{i+1}$ without changing the flag and, if $\alpha_{F,G}(i)<\alpha_{F,G}(i+1)$, notice that $\alpha_{F',G}=\alpha_{F,G}\sigma_i$. The second condition follows. For the third condition, let $F,G\in \Ch(\Delta)$ be such that $\delta(F,G)=\alpha\in \Sym(n)$ and consider any ordered basis $(e_1,...,e_n)$ representing $F$ such that $(e_{\alpha(1)}, ... e_{\alpha(n)})$ represents $G$. It is clear from the definition of $\delta$ that $F'=(e_1,...,e_{i-1}, e_{i+1}, e_i,..., e_n)$ satisfies $\delta(F',F)=\sigma_i$ and $\delta(F',G)=\sigma_i \alpha$. This proves that $\delta$ is a $\Sym(n)$-distance and our claim follows. 
	\end{example}
	In these notes, we will only be interested in right-angled buildings. We recall that a building is \tg{right-angled} if its type $(W,I)$ is given by a right-angled Coxeter system, that is, if for any two generators $i,j\in I$ we either have that $i$ and $j$ commute or generate a free product $C_2*C_2$ of two copies of the cyclic group of order two $C_2$. The following result ensures the existence of a rich family of regular buildings for any right-angled Coxeter system. 
	\begin{theorem}[{\cite[Proposition 1.2]{Haglund2003}}]
		Let $(W,I)$ be a right-angled Coxeter system and $(q_i)_{i\in I}$ be a set of positive integers with $q_i\geq 2$. There exists a
		right-angled building $\Delta$ of type $(W, I)$ such that for every $i\in I$, each $i$-panel of $\Delta$ has
		size $q_i$. This building is unique, up to isomorphism.
	\end{theorem}
	\begin{definition}
		The essentially unique building $\Delta$ given by the above theorem is the \tg{semi-regular building of prescribed thickness} $(q_i)_{i\in I}$.  
	\end{definition}
	Our next task is to recall one of the fundamental features of buildings that is the existence of combinatorial projections between residues. To this end, we make a series of definition.
	\begin{definition}
		Two chambers $c,d\in \Ch(\Delta)$ are said to be $i$\tg{-adjacent} for some $i\in I$ if $\delta(c,d)=i$. A \tg{gallery} in $\Delta$ between two chambers $c,d$ is a finite sequence $c_1,...,c_n$ of chambers such that $c_1=c$, $c_n=d$ and such that $c_t$ and $c_{t+1}$ are $i_t$-adjacent for all $t\in \{1,...,n-1\}$. In that case, the gallery is said to have \tg{length} $n-1$ and it is \tg{minimal} if there is no shorter gallery between $c$ and $d$. 
	\end{definition}
	This notion of gallery provides a discrete valued metric $$d_{\Delta}:\Ch(\Delta)\times \Ch(\Delta)\rightarrow \N$$ on the set of chambers $\Ch(\Delta)$ where $d_{\Delta}(c,d)$ is the length of a minimal gallery containing both $c$ and $d$. Now, given a chamber $c\in \Ch(\Delta)$ and a residue $\mathcal{R}$ in a building $\Delta$, the \tg{gate property} ensures the existence of a unique chamber $d\in \Ch(\mathcal{R})$ that is closest to $c$ for the chamber metric $d_{\Delta}$. This unique chamber is called the \tg{projection} of $c$ on $\mathcal{R}$ and is denoted by $\proj_{\mathcal{R}}(c)$. We refer to \cite{Tits1974} for more details about this notion and state some of its properties. To start with, we recall that for any two residues $\mathcal{R}$ and $\mathcal{R}'$ of $\Delta$, the set
	$$\{\proj_\mathcal{R}(c): c\in \Ch(\mathcal{R}')\}$$
	is the chamber-set of a residue of $\Delta$ contained in $\mathcal{R}$. Furthermore, as we recall below, if $\Delta$ is a right-angled building this notion of projection provides a way to partition the building into convex chamber sets. A subset $C\subseteq \Ch(\Delta)$ is said to be \tg{convex} if for any two chambers $c,d\in C$, $C$ contains each of the chambers appearing in any minimal gallery between $c$ and $d$. 
	\begin{definition}
		Let $\Delta$ be a right-angled building of type $(W,I)$, $J\subseteq I$ and $c\in \Ch(\Delta)$. The \textbf{$J$-wing} of $\Delta$ containing $c$ is the set $$X_J(c)=\{d\in \Ch(\Delta): \proj_{\mathcal{R}_J(c)}(d)=c\}.$$
		If $J = \{i\}$ is a singleton, we refer to this set as the \textbf{$i$-wing} $X_i(c)$ of $c$.
	\end{definition} It is shown in \cite{Caprace2014}, that  wings are convex chamber sets and that for every $J$-residue $\mathcal{R}$, $\Ch(\Delta)$ is partitioned by the $J$-wings $X_J(c)$ with $c\in \mathcal{R}$. We now recall some of their properties.
	\begin{lemma}[{\cite[Lemma $3.1$]{Caprace2014}}]\label{la J wing est l'intersection des j wing wesh wesh wesh}
		Let $\Delta$ be a right-angled building of type $(W,I)$, $J\subseteq I$ be a non-empty set and $c\in \Ch(\Delta)$. The $J$-wing $X_J(c)$ containing $c$ satisfies the following properties:
		\begin{enumerate}
			\item $X_J(c)=\bigcap_{i\in J}X_i(c)$. 
			\item $X_J(c) = X_J(c')$ for all $c'\in X_J (c) \cap \mathcal{R}_{J\cup J^\perp}(c)$.
		\end{enumerate}
	\end{lemma}
	\begin{lemma}\label{residues of different types intersects in a single chamber}
		Let $\Delta$ be a right-angled building of type $(W,I)$, $J,J'\subseteq I$ be two disjoint subsets and $c\in \Ch(\Delta)$. Then $\mathcal{R}_{J}(c)\subseteq X_{J'}(c)$.  
	\end{lemma}
	\begin{proof}
		The result follows directly from the fact the two residues $\mathcal{R}_{J}(c)$ and $\mathcal{R}_{J'}(c)$ contain $c$ and that the intersection of a $J$-residue and a $J'$-residue is a $J\cap J'$-residue. In particular, this proves that $\proj_{\mathcal{R}_{J'}(c)}(\mathcal{R}_{J}(c))=\{c\}$ and therefore that $\mathcal{R}_{J}(c)\subseteq X_{J'}(c)$. 
	\end{proof}
	\begin{lemma}[{\cite[Lemma $3.4$]{Caprace2014}}]\label{the inclusion of wings}
		Let $\Delta$ be a right-angled building of type $(W,I)$, $i,i'\in I$ be such that $i=i'$ or $m_{i,i'}=\infty$, and $c,c'\in \Ch(\Delta)$ be such that $c'\in X_i(c)$ but $c\not\in X_{i'}(c')$. Then, we have $X_{i'}(c')\subseteq X_i(c)$. 
	\end{lemma}
	
	Another feature of combinatorial projections is given by the relation of parallelism. Two residues $\mathcal{R}$ and $\mathcal{R}'$ in a building $\Delta$ are said to be \tg{parallel} if $\proj_{\mathcal{R}}(\mathcal{R}')=\mathcal{R}$ and $\proj_{\mathcal{R}'}(\mathcal{R})=\mathcal{R}'$. Notice that the chamber sets of parallel residues are in bijection under the respective projection maps and that two parallel residues have the same rank. Caprace showed that the relation of parallelism has a particular flavour in right-angled buildings.
	\begin{lemma}[{\cite[Corollary $2.9$]{Caprace2014}}]
		Let $\Delta$ be a right-angled building. Then the relation of parallelism of residues is an equivalence relation.
	\end{lemma}
	\begin{lemma}[{\cite[Proposition $2.8$]{Caprace2014}}]
		Let $\Delta$ be a right-angled building of type $(W,I)$ and let $J\subseteq I$. Then two $J$-residues $\mathcal{R}$ and $\mathcal{R}'$ are parallel if and only if they are both contained in a common $J\cup J^\perp$-residue where $J^\perp=\{i\in I: ij=ji\qq \forall j\in J\}$.
	\end{lemma}
	
	We now recall the notion of universal groups of semi-regular right-angled buildings introduced by Tom De Medts, Ana C. Silva and Koen Struyve in \cite{Universal2018} to generalise the concept of Burger and Mozes universal groups on trees \cite{BurgerMozes2000}. We start by recalling the definition of type-preserving automorphisms of buildings. 
	\begin{definition}
		Let $\Delta$ be a building of type $(W,I)$ and Weyl distance $\delta$. A \tg{type-preserving automorphism} of $\Delta$ is a bijection 
		$$g: \Ch(\Delta)\rightarrow \Ch(\Delta): c \mapsto gc$$
		preserving the Weyl distance in the sense that for all $ c,d\in \Ch(\Delta)$ one has
		$$\delta (gc,gd)=\delta (c,d).$$
		We denote by $\Aut(\Delta)$ the group of type-preserving automorphisms of $\Delta$.
	\end{definition}
	The definition of universal groups of right-angled buildings requires a notion of coloring of the building. Let $\Delta$ be a semi-regular right-angled building of type $(W, I)$ and prescribed thickness $(q_i)_{i\in I}$. 
	\begin{definition}
		For each $i\in I$, let $Y_i$ be a set of size $q_i$ which we will refer to as the set of $i$-colors of $\Delta$. A \tg{set of legal colorings} of $\Delta$ is a set $(h_i)_{i\in I}$ of maps $$h_i:\Ch(\Delta) \rightarrow Y_i$$ such that $\restr{h_i}{\Ch(\tau)}:\Ch(\tau) \rightarrow Y_i$ is a bijection on each $i$-panel $\tau$ of $\Delta$, and such that $h_i (c) = h_i (c' )$ for every $(I-\{i\})$-residue $\mathcal{R}$ and each $c, c' \in \Ch(\mathcal{R})$.
	\end{definition}
	Now, for each $i\in I$, let $G_i\leq \Sym(Y_i)$ be a transitive permutation group and let $(h_i)_{i\in I}$ be a set of legal colorings of $\Delta$.
	\begin{definition}
		 The \tg{universal group} $\mathcal{U}((h_i,G_i)_{i\in I})$ of $\Delta$ with respect to the set of legal colorings $(h_i)_{i\in I}$ and prescribed local action $(G_i)_{i\in I}$ is the subgroup of $\Aut(\Delta)$ defined by
		\begin{equation*}
		\{g\in \Aut(\Delta): (\restr{h_i}{g\tau})\circ g\circ (\restr{h_i}{\tau})^{-1}\in G_i \qq \forall i\in I \mbox{ and every }i-\mbox{panel }\tau \mbox{ of }\Delta\}
		\end{equation*}
	\end{definition}
	\noindent It appears that those groups satisfy a factorisation property similar the Tits independence property that we now introduce. 
	\begin{definition}
		Let $\Delta$ be a right-angled building of type $(W,I)$ and $J\subseteq I$. A subgroup $G\leq \Aut(\Delta)$ is said to satisfy the property ${\rm IP}_{J}$ if for all $J\cup J^\perp$-residues $\mathcal{R}$ of $\Delta$ we have that
		\begin{equation*}\tag{${\rm IP}_{J}$}\label{IPJ}
		\Fix_G(\mathcal{R})= \prod_{c\in \mathcal{R}} \Fix_G(V_{J}(c))
		\end{equation*} 
		where $V_{J}(c)=\{ d\in \Ch(\Delta): \proj_{\mathcal{R}}(d)\not=c\}$ is the complement of the $J$-wing containing $c$ and where $\Fix_G(\mathcal{R})=\{g\in G : gc=c\qq \forall c\in \Ch(\mathcal{R})\}$.
	\end{definition}
	The following proposition ensures that every universal group of right-angled building satisfies the property ${\rm IP}_{\{i\}}$ for every $i\in I$. Furthermore, Proposition \ref{independence IPJ for universal} below ensures that they also satisfy the property ${\rm IP}_{J}$ for every finite set $J\subseteq I$ such that $\{i\}\cup \{i\}^\perp=J$ $\forall i\in J$. 
	\begin{proposition}[{\cite[Proposition $3.16$]{Universal2018}}]\label{universal satify IP i de larticle universal}
		Let $G$ be a universal group of a semi-regular right-angled building $\Delta$. Then $G$ satisfies the property ${\rm IP}_{\{i\}}$ for every $i\in I$. Furthermore,  for all $i\in I$, $ g\in \Fix_G(\mathcal{R}_{\{i\}\cup \{i\}^\perp}(c))$ and for each $c\in \Ch(\Delta)$ the type-preserving automorphism 
		\begin{equation*}
		g^c : \Ch(\Delta)\rightarrow \Ch(\Delta): x \mapsto \begin{cases}
		\qq g x\qq&\mbox{ if }x \in X_{i}(c) \\
		\qq x \qq &\mbox{ if }x\in V_{i}(c)
		\end{cases} 
		\end{equation*}
		is an element of $G$.  
	\end{proposition}
	\begin{proposition}\label{independence IPJ for universal}
		Let $G\leq \Aut(\Delta)$ be a universal group of a semi-regular right-angled building $\Delta$. Then, $G$ satisfies the property ${\rm IP}_{J}$ for every finite set $J\subseteq I$ such that $\{i\}\cup \{i\}^\perp=J$  $\forall i\in J$. 
	\end{proposition}
	\begin{proof}
		If $\modu{J}=1$, the results follows directly from Proposition \ref{universal satify IP i de larticle universal}. Suppose therefore that $\modu{J}\geq 2$ and let us show that $$\Fix_G(\mathcal{R})= \prod_{c\in \mathcal{R}} \Fix_G(V_{J}(c))$$ for every $J$-residue $\mathcal{R}$. First, notice that $\Fix_G(V_{J}(c))$ is a subgroup of $\Fix_G(\mathcal{R})$ for every $c\in \mathcal{R}$. Lemma \ref{la J wing est l'intersection des j wing wesh wesh wesh} ensures that $X_J(c)=\bigcap_{i\in J}X_i(c)$ for every $c\in \Ch(\mathcal{R})$. In particular, taking the complement, we obtain that $V_J(c)=\bigcup_{i\in J}V_i(c)$ and therefore that $\Fix_G(V_J(c))=\bigcap_{i\in J}\Fix_G(V_i(c))$. Let $i\in J$ and remember from Proposition \ref{universal satify IP i de larticle universal} that $G$ satisfies the property ${\rm IP}_{\{i\}}$ which implies that $\Fix_G(V_i(c))\subseteq \Fix_G(\mathcal{R}_{\{i\}\cup \{i\}^{\perp}}(c))$. Furthermore, since  $J=\{i\}\cup \{i\}^\perp$, we obtain that $\mathcal{R}=\mathcal{R}_{\{i\}\cup \{i\}^{\perp}}(c)$ and therefore that $\Fix_G(V_J(c))\subseteq \Fix_G(\mathcal{R})$. Notice that for every two distinct $c,d\in \Ch(\mathcal{R})$, the supports of the elements of $\Fix_G(V_i(c))$ and $\Fix_G(V_i(d))$ are disjoint from one another which proves that $\prod_{c\in \mathcal{R}} \Fix_G(V_{J}(c))$ is a well defined subgroup of $G$. The above discussion proves that $\prod_{c\in \mathcal{R}} \Fix_G(V_{J}(c))\subseteq\Fix_G(\mathcal{R})$. To prove the other inclusion, let $g\in \Fix_G(\mathcal{R})$ and let $J=\{i_1,...,i_n\}$. For any $i\in J$ and $c\in \Ch(\mathcal{R})$, let $\mathcal{R}_i(c)$ be the $i$-panel containing $c$ in $\mathcal{R}$. Let us fix some chamber $c\in \Ch(\mathcal{R})$ and let 
		\begin{equation*}
		g_1^d : \Ch(\Delta)\rightarrow \Ch(\Delta): x \mapsto \begin{cases}
		\qq g x \qq&\mbox{ if }x \in X_{i_1}(d) \\
		\qq x \qq &\mbox{ if }x\in V_{i_1}(d)
		\end{cases} 
		\end{equation*}
		for every $d\in \mathcal{R}_{i_1}(c)$.  Proposition \ref{universal satify IP i de larticle universal} ensures that $g^d_1\in G$ $\forall d\in \mathcal{R}_{i_1}(c)$ and that $g=\prod_{d\in \mathcal{R}_{i_1}(c)}g^d_1$. On the other hand, for every $d\in \mathcal{R}_{i_1}(c)$, there exists a unique $i_2$-panel $\mathcal{R}_{i_2}(d)$ such that $d\in \mathcal{R}_{i_2}(d)$. Since $g^d_1\in \Fix_G(V_{i_1}(d))\subseteq \Fix_G(\mathcal{R})$, we can repeat the above argument and we obtain that $g^d_1 = \prod_{d'\in \mathcal{R}_{i_2}(d)}g^{d'}_2$ where 
		\begin{equation*}
		g^{d'}_2 : \Ch(\Delta)\rightarrow \Ch(\Delta): x \mapsto \begin{cases}
		\qq g^d_1 x \qq&\mbox{ if }x \in X_{i_2}(d') \\
		\qq x \qq &\mbox{ if }x\in V_{i_2}(d')
		\end{cases} 
		\end{equation*}
		for every $d'\in \mathcal{R}_{i_2}(d)$. Just as before, Proposition \ref{universal satify IP i de larticle universal} ensures that $g^{d'}_2\in G$ for every $d'\in \mathcal{R}_{i_2}(d)$. On the other hand, $\mathcal{R}_{i_1}(c)\cap \mathcal{R}_{i_2}(d)=\mathcal{R}_{i_1}(d)\cap \mathcal{R}_{i_2}(d)=\{d\}$. In particular, this implies that $\proj_{\mathcal{R}_{i_1}(c)}(\mathcal{R}_{i_2}(d))=\{d\}$ and $d'\in X_{i_1}(d)$. Since $\mathcal{R}$ is an $\{i_1\}\cup \{i_1\}^\perp$-residue, Lemma \ref{la J wing est l'intersection des j wing wesh wesh wesh} ensures that $X_{i_1}(d)=X_{i_1}(d')$. This proves that $g_2^{d'}$ has support in $X_{i_1}(d')\cap X_{i_2}(d')=X_{\{i_1, i_2\}}(d')$ and therefore that $g_{2}^{d'}\in \Fix_G(V_{\{i_1,i_2\}}(d'))$. Proceeding iteratively, for any of the constructed $g^d_k$ with $k\in \{2,...,n-1\}$, we set
		\begin{equation*}
		g^{d'}_{k+1} : \Ch(\Delta)\rightarrow \Ch(\Delta): x \mapsto \begin{cases}
		\qq g^d_k x \qq&\mbox{ if }x \in X_{i_{k+1}}(d') \\
		\qq x \qq &\mbox{ if }x\in V_{i_{k+1}}(d')
		\end{cases} 
		\end{equation*}
		for every $d'\in \mathcal{R}_{i_{k+1}}(d)$. Once more, Proposition \ref{universal satify IP i de larticle universal} ensures that $g^{d'}_{k+1}\in G$  $\forall d'\in \mathcal{R}_{i_{k+1}}(d)$ and that $g^d_k=\prod_{d'\in \mathcal{R}_{i_{k+1}}(d)}g^{d'}_{k+1}$. On the other hand, $\mathcal{R}_{i_l}(d)\cap \mathcal{R}_{i_{k+1}}(d)=\{d\}$ for every $l=1,...,k$ which implies that $d'\in X_{i_{l}}(d)$. Since $\mathcal{R}$ is an $\{i_l\}\cup \{i_l\}^\perp$-residue, Lemma \ref{la J wing est l'intersection des j wing wesh wesh wesh} ensures that $X_{i_l}(d)=X_{i_l}(d')$. Finally, since $g^d_k$ has support in $X_{i_1}(d)\cap X_{i_2}(d)\cap ...\cap X_{i_k}(d)$ this proves that $g_{k+1}^{d'}$ has support in $X_{i_1}(d') \cap ... X_{i_k}(d')\cap X_{i_{k+1}}(d')=X_{\{i_1,...,i_{k+1}\}}(d')$ and therefore that $g^{d'}_{k+1}\in \Fix_G(V_{\{i_1,i_2,..., i_{k+1}\}}(d'))$. 
	\end{proof}
	
	\subsection{Reduction to groups of automorphisms of trees}
	
	The purpose of this section is to show that the group of type preserving automorphisms $\Aut(\Delta)$ of certain semi-regular right-angled buildings $\Delta$ can be realised as closed subgroups of the group $\Aut(T)^+$ of type-preserving automorphisms of a locally finite tree $T$ in such a way that the universal groups of those buildings embed as closed subgroups $G\leq \Aut(T)^+$ satisfying the property \ref{IPV1}. This applies only to certain Coxeter types and motivates the following definition.
	\begin{definition}
		A right-angled Coxeter system $(W,I)$ is said to satisfy the hypothesis \ref{Hypothese cox star} if it is finitely generated and there exists $r\geq 2$ such that
		\begin{equation}\tag{$\star$}\label{Hypothese cox star}
		I=\bigsqcup_{k=1}^rI_k
		\end{equation}
		for some $I_k=\{i\}\cup \{i\}^\perp$  $\forall i\in I_k$ and $\forall k=1,...,r$.
	\end{definition} 
	\begin{remark}
		A right-angled Coxeter system satisfying the hypothesis \ref{Hypothese cox star} is isomorphic to a free product $W_1*W_2*...*W_r$ where each of the $W_k$ is a direct product of finitely many copies of the group of order $2$. In particular, $W$ is virtually free.
	\end{remark} 
	\noindent Let $(W,I)$ be a right-angled Coxeter system satisfying the hypothesis \ref{Hypothese cox star}, let $(q_i)_{i\in I}$ be a set of positive integers with $q_i\geq 2$ and let $\Delta$ be a semi-regular building of type $(W,I)$ and prescribed thickness $(q_i)_{i\in I}$. We associate a locally finite bipartite graph to $\Delta$ as follows. We let $V_0=\Ch(\Delta)$, $$V_1= \{ \mathcal{R}: \mathcal{R} \mbox{ is an }I_k-\mbox{residue of }\Delta \mbox{ for some }k\in \{1,...,r\}\}$$ and we define $T$ as the bipartite graph with vertex set $V= V_0 \sqcup V_1$ and where a chamber $c\in V_0$ is adjacent to a residue $\mathcal{R}\in V_1$ if $c\in \mathcal{R}$. 
	\begin{lemma}\label{the tree gamma associated with Delta}
		The graph $T$ is a locally finite tree. 
	\end{lemma}
	\begin{proof}
		The graph $T$ is locally finite since each chamber is contained in finitely many residues and since each $I_k$-residue is finite. The graph $T$ is path connected since every two chambers of $\Delta$ are connected by a gallery and since each such gallery corresponds naturally to a path in $T$. We now  show that $T$ does not contain any cycle. Suppose for a contradiction that there is a simple cycle in $T$, say
		\begin{equation*}
		c_1 - \mathcal{R}_1 - c_2 - ... - \mathcal{R}_n - c_1.
		\end{equation*}
		Since each chamber $c\in \Ch(\Delta)$ is contained in a unique residue $\mathcal{R}$ of type $I_t$ and since the cycle is simple, notice that $\mathcal{R}_1$ and $\mathcal{R}_n$ have different types. In particular, Lemma \ref{residues of different types intersects in a single chamber} ensures that $\mathcal{R}_n\subseteq X_{J_1}(c_1)$ where $J_1$ is the type of $\mathcal{R}_1$. On the other hand, as we show below, $\mathcal{R}_n\subseteq X_{J_1}(c_2)$. For now, we assume this inclusion that is $\mathcal{R}_n\subseteq X_{J_1}(c_1)\cap X_{J_1}(c_2)$ and we show this leads to a contradiction. Since the cycle is simple, we have that $c_1\not=c_2$. Furthermore, since $c_1,c_2\in\mathcal{R}_1$, there exists some $i\in J_1$ such that $c_1\not\in X_i(c_2)$. Hence, we have that $X_i(c_1)\cap X_i(c_2)=\es$ and therefore that $X_{J_1}(c_1)\cap X_{J_1}(c_2)=\es$. The desired contradiction follows from our inclusion.
		
		Now, let us prove that $\mathcal{R}_n\subseteq X_{J_1}(c_2)$. To this end, we show that $X_{J_{t+1}}(c_{t+2})\subseteq X_{J_t}(c_{t+1})$ for every $t\in \{1,..., n-2\}$ where $J_t$ is the type of $\mathcal{R}_t$. Since the cycle is simple notice that $\mathcal{R}_t$ and $\mathcal{R}_{t+1}$ have different types and that $\mathcal{R}_{t+1}\subseteq X_{J_t}(c_{t+1})$ for every $t=1,...,n-2$. On the other hand, since $c_{t+1}\not=c_{t+2}$, there exists some $i'\in J_{t+1}$ such that $c_{t+1}\not\in X_{i'}(c_{t+2})$. Notice for every $i\in J_t$ that $m_{i,i'}=\infty$ and that $c_{t+2}\in \mathcal{R}_{t+1}\subseteq X_{i}(c_{t+1})$. In particular, Lemma \ref{the inclusion of wings} implies that $X_{J_{t+1}}(c_{t+2})\subseteq X_{i'}(c_{t+2})\subseteq \bigcap_{i\in J_t}X_i(c_{t+1})= X_{J_t}(c_{t+1})$ which completes the induction. This proves as desired that $\mathcal{R}_n\subseteq X_{J_{n-1}}(c_{n})\subseteq ...\subseteq X_{J_1}(c_2)$.  
	\end{proof}
	
	Our next goal is to explicit an injective map $\alpha :\Aut(\Delta)\rightarrow\Aut(T)^+$ defining a homeomorphism on its image. Notice that any type-preserving automorphism $g\in \Aut(\Delta)$ is bijective on the set of chambers $\Ch(\Delta)$ but also on the $I_k$-residues of $\Delta$ for any fixed $k$. For every $g\in \Aut(\Delta)$ we define the map $\alpha(g): V\rightarrow V$ as follows: 
	\begin{itemize}
		\item If $v\in V_0$ then $v$ is a chamber $c\in \Ch(\Delta)$ and we define $\alpha(g)v=gc$.
		\item If $v\in V_1$ then $v$ is an $I_k$-residue $\mathcal{R}$ of $\Delta$ for some $k\in \{1,...,r\}$ and we define $\alpha(g)v= g\mathcal{R}$.
	\end{itemize} 
	The map $\alpha(g)$ clearly defines a type preserving bijection on $V$. In fact, $\alpha(g)$ is a tree automorphism of $T$ and $\alpha: \Aut(\Delta)\rightarrow \Aut(T)^+$ is a well defined group homomorphism  since  for every $g\in \Aut(\Delta)$, every residue $\mathcal{R}$ of $\Delta$ and every $ c\in \Ch(\Delta)$, we have that $c\in \mathcal{R}$ if and only if $gc\in g\mathcal{R}$.
	\begin{proposition}\label{the automrophism of buildings as automorphims of locally finite graphs}
		The map $\alpha :\Aut(\Delta)\rightarrow \Aut(T)^+$ is an injective group homomorphism; $\alpha(\Aut(\Delta))$ is a closed subgroup of $\Aut(T)^+$ and $\alpha$ defines a homeomorphism between $\Aut(\Delta)$ and $\alpha(\Aut(\Delta))$. 
	\end{proposition}
	\begin{proof}
		The homomorphism $\alpha$ is injective, since 
		\begin{equation*}
		\begin{split}
		\ker(\alpha)&=\{g\in \Aut(\Delta): \alpha(g)=1_{\Aut(T)}\}\\
		&\subseteq\{g\in \Aut(\Delta): gc=c\qq \forall c\in \Ch(\Delta)\}=\{1_{\Aut(\Delta)}\}.
		\end{split}
		\end{equation*}
		Remember that the sets $$U_{T}(F_{T})=\{g \in \Aut(T)^+: g v=v \qq\forall v\in F_T\}$$ 
		where $F_T\subsetneq V_0$ is finite form a basis of open neighbourhoods of the identity in $\Aut(T)^+$. On the other hand, an element $h\in \Aut(T)^+$ belongs to $\Aut(T)^+-\alpha(\Aut(\Delta))$ if and only if there exists $i\in I$ and two $i$-adjacent chambers $c,d\in V_0$ such that $hc$ and $hd$ are not $i$-adjacent. In particular, for every such automorphism $h$, the set $hU_T(\{c,d\})$ is an open neighbourhood of $h$ in $\Aut(T)^+-\alpha(\Aut(\Delta))$. This proves that the complement of $\alpha(\Aut(\Delta))$ is an open set and therefore that $\alpha(\Aut(\Delta))$ is a closed subgroup of $\Aut(T)^+$.
		
		Let $\Phi: \Ch(\Delta)\rightarrow V_0$ be the map sending a chamber of $\Delta$ to the corresponding vertex of $V_0\subseteq V$ and remember that the sets 
		$$U_\Delta(F_\Delta)=\{g\in \Aut(\Delta): g c=c \qq \forall c\in F_\Delta\}$$ 
		where $F_\Delta\subseteq \Ch(\Delta)$ is finite form a basis of open neighbourhoods of the identity in $\Aut(\Delta)$.
		In particular, notice that $\alpha : \Aut(\Delta)\rightarrow \alpha (\Aut(\Delta))$ is continuous since for every finite subset $F_T\subseteq V_0$ we have that $\alpha^{-1}(U_{T}(F_T)\cap \alpha(\Aut(\Delta))=U_\Delta(\Phi^{-1}(F_T))$. Finally, notice that  $\alpha : \Aut(\Delta)\rightarrow \alpha (\Aut(\Delta))$ is an open map since, for every finite set $F_\Delta\subsetneq\Ch(\Delta)$ we have that $$\alpha(U_{\Delta}(F_\Delta))=U_\Delta(\Phi(F_\Delta))\cap \alpha(\Aut(\Delta)).$$
	\end{proof}
	The following proposition shows that under this correspondence, the property \ref{IPV1} of groups of type-preserving automorphisms of trees is tightly related to the property \ref{IPJ} of groups of type-preserving automorphisms of right-angled buildings.  
	\begin{proposition}\label{lemme de on a IP IK alors IP V1}
		Let $G\leq \Aut(\Delta)$ be a closed subgroup satisfying the property ${\rm IP}_{I_k}$ for every $k=1,...,r$. Then, $\alpha(G)$ is a closed subgroup of $\Aut(T)^+$ satisfying the property ${\rm IP}_{V_1}$. 
	\end{proposition}
	\begin{proof}
		Proposition \ref{the automrophism of buildings as automorphims of locally finite graphs} ensures that $\alpha(G)$ is a closed subgroup of $\alpha(\Aut(\Delta))$ and therefore of $\Aut(T)^+$. Let $\Phi: \Ch(\Delta)\rightarrow V_0$ be the map sending a chamber of $\Delta$ to the corresponding vertex of $V_0\subseteq V$ and notice that:
		\begin{itemize}
			\item For every residue $\mathcal{R}\in V_1$, $\Phi(\mathcal{R})=\{v\in V_0: v\in \mathcal{R}\}=B_T(\mathcal{R},1)\cap V_0$. In particular, we have that
			\begin{equation*}
			\begin{split}
			\alpha(\Fix_{G}(\mathcal{R}))&=\Fix_{\alpha(G)}(\Phi(\mathcal{R}))=\Fix_{\alpha(G)}(B_T(\mathcal{R},1)\cap V_0)\\
			&=\Fix_{\alpha(G)}(B_T(\mathcal{R},1)).
			\end{split}
			\end{equation*} 
			\item For every $k\in \{1,...,r\}$ and every chamber $c\in \Ch(\Delta)$ we have
			\begin{equation*}
			\begin{split}
			\Phi(\Ch(\Delta)-X_{I_k}(c))&=\Phi(\{d\in \Ch(\Delta): \proj_{\mathcal{R}_{I_k}(c)}(d)\not=c\})\\
			&=\{v\in V_0: d_T(v,\mathcal{R}_{I_k}(c))< d_T(v,c)\}\\
			&=T( \mathcal{R}_{I_k}(c),c)\cap V_0.
			\end{split}
			\end{equation*}
			In particular, we obtain that
			\begin{equation*}
			\begin{split}
			\alpha (\Fix_G(\Ch(\Delta)-X_{I_k}(c)))&=\Fix_{\alpha(G)}(\Phi(\Ch(\Delta)-X_{I_k}(c)))\\
			&=\Fix_{\alpha(G)}(T(\mathcal{R}_{I_k}(c), c)).
			\end{split}
			\end{equation*}
		\end{itemize}
		The results follows from the definitions of the properties ${\rm IP}_{I_k}$ and ${\rm IP}_{V_1}$. 
	\end{proof}
	
	\noindent 
	The proof of Theorem \ref{Theorem E} requires one last preliminary. Let $(W,I)$ be a right-angled Coxeter system. Let $(q_i)_{i\in I}$ be a set of positive integers with $q_i\geq 2$ and let $\Delta$ be the semi-regular building of type $(W,I)$ and prescribed thickness $(q_i)_{i\in I}$. Finally, let $(h_i)_{i\in I}$ be a set of legal colorings of $\Delta$ with $i$-colors given by a set $Y_i$ of size $q_i$ and let $G_i\leq \Sym(Y_i)$ be transitive on $Y_i$. 
	\begin{proposition}\label{les groupes universelles sont unimodular}
		For every $J\subseteq I$, every $(g_{j})_{j\in J}\in \prod_{j\in J}^{}G_j$ and every $c\in \Ch(\Delta)$, there exists $g\in \mathcal{U}((h_i,G_i)_{i\in I})$ such that $g\mathcal{R}_J(c)=\mathcal{R}_{J}(c)$, $h_j\circ g=g_j\circ h_j$ for all $j\in J$ and $h_i\circ g=h_i$ for all $i\in I-J$.
	\end{proposition}
	\begin{proof}
		Let $(h_i')_{i\in I} $ be the set of legal colorings obtained from $(h_i)_{i\in I}$ by replacing $h_j$ by $g_j\circ h_j$ for all $j\in J$ and leaving the other colorings unchanged. Notice that $h_j'$ is still a legal coloring of $G$ for every $j\in J$ since for all $I-\{j\}$-residues $\mathcal{R}$ and for all $d,d'\in\mathcal{R}$ we have $h_j'(c)=g_j\circ h_j(d)=g_j\circ h_j(d')=h_j'(d')$. Now, let $c'$ be the chamber of $\mathcal{R}_{J}(c)$ with colors $h_j(c')=g_j\circ h_j(c)$ for every $j\in J$. Since $h_i'(c')=h_i(c)$ for every $i\in I$, \cite[Proposition $2.44$]{Universal2018} ensures the existence of an automorphism $g\in \Aut(\Delta)$ mapping $c$ to $c'$ and such that $h_i\circ g=h_i'$ for all $i\in I$. Since $c'\in \mathcal{R}_J(c)$ notice that $g$ stabilises $\mathcal{R}_J(c)$. Finally, notice that $g$ acts locally as the identity on $i$-panels for all $i\in I-J$ and as $g_j$ on $j$-panels for all $j\in J$. Hence, $g$ is as desired.
	\end{proof}
	The following result proves Theorem \ref{Theorem E}.
	\begin{theorem}\label{theorem pour les right angled buildings}
		Suppose that $(W,I)$ satisfies the hypothesis \ref{Hypothese cox star}.  The group $\alpha(\mathcal{U}((h_i,G_i)_{i\in I}))$ is a closed subgroup of $\Aut(T)^+$ satisfying the property \ref{IPV1}. In addition, if $G_i$ is \textit{2}-transitive on $Y_i$ for every $i\in I$, $\alpha(\mathcal{U}((h_i,G_i)_{i\in I}))$ is unimodular and satisfies the hypothesis \ref{Hypothese HV1}.
	\end{theorem}  
	\begin{proof}
		Let $G = \mathcal{U}((h_i,G_i)_{i\in I})$. The first part of theorem follows directly from Proposition \ref{independence IPJ for universal} and Proposition \ref{lemme de on a IP IK alors IP V1}. Now, suppose that $G_i$ is \textit{2}-transitive on $Y_i$ for every $i\in I$, let $\mathfrak{T}_{V_1}$ be the family of subtrees defined on page \pageref{page Tv1} and let $\mathcal{T}, \mathcal{T}'\in \mathfrak{T}_{V_1}$. If $\mathcal{T}\subseteq\mathcal{T}'$ we clearly have that $\Fix_{\alpha(G)}(\mathcal{T}')\subseteq \Fix_{\alpha(G)}(\mathcal{T})$. On the other hand, if $\mathcal{T}\not \subseteq \mathcal{T}'$, let $v_\mathcal{T}$ a vertex of $V_1\cap \mathcal{T}$ that is at maximal distance from $\mathcal{T}'$ and let $w_\mathcal{T}\in B_T(v_\mathcal{T},1)-\{v_\mathcal{T}\}$ be such that  $\mathcal{T}'\subseteq T(w_\mathcal{T}, v_\mathcal{T})$.  In particular, we have that $\Fix_{\alpha(G)}(T(w_\mathcal{T},v_\mathcal{T}))\subseteq \Fix_{\alpha(G)}(\mathcal{T}')$. Let $I_k$ denote the type of $v_\mathcal{T}$ seen as a residue of $\Delta$ and let $j\in I_k$. Consider an element $g_j\in G_i$ that is not trivial and such that $g_j\circ h_j(w_{\mathcal{T}})=h_j(w_\mathcal{T})$ and let $g_i={\rm id}_{Y_i}$ for every $i\in I-\{j\}$. Proposition \ref{les groupes universelles sont unimodular} ensures the existence of an element $g\in \mathcal{U}((h_i,G_i)_{i\in I})$ such that $g\mathcal{R}_J(c)=\mathcal{R}_{J}(c)$ and $h_i\circ g=g_i\circ h_i$ for every $i\in I$. Now, notice that there exists a unique vertex $v'_{\mathcal{
				T}}$ that is adjacent to $w_{\mathcal{T}}$ and such that $\mathcal{T}'\subseteq T(v'_\mathcal{T},w_\mathcal{T}) \cup \{w_\mathcal{T}\}$. Now, notice that $v'_\mathcal{T}$ is a residue of type $I_{k'}$ with $I_{k'}\not = I_k$ and we realise from the definition that $g$ fixes every chamber of $v'_{\mathcal{T}}$ or equivalently that $\alpha(g)$ fixes $B_T(v'_{\mathcal{T}},1)$ pointwise. Proposition \ref{lemme de on a IP IK alors IP V1} implies the existence of an element $h\in \Fix_{\alpha(G)}(B_T(v'_\mathcal{T}),1)\cap \Fix_G( T(v'_{\mathcal{T}},w_{\mathcal{T}}))$ such that $hv=\alpha(g)v$ for every $v\in T(w_{\mathcal{T}},v'_{\mathcal{T}})$. In particular, we have that $h\in \Fix_{\alpha(G)}(\mathcal{T}')$ but $h\not \in \Fix_{\alpha(G)}(\mathcal{T})$ since by the definition $g$ does not fix every chamber of $v_{\mathcal{T}}$. This proves as desired that $\Fix_{\alpha(G)}(\mathcal{T}')\subsetneq \Fix_{\alpha(G)}(\mathcal{T})$.
		
		To prove that $\alpha(G)$ is unimodular, we apply \cite[Corollary 5]{BaumgartnerBertrandWillis2007}. This result ensures that a group $G$ which acts $\delta$-$2$-transitively on the set of chambers of a locally finite building is unimodular. Choose a chamber $c\in \Ch(\Delta)$. Since $G$ is transitive on the chambers of $\Delta$, we need to show for any two chambers $d_1,d_2\in \Ch(\Delta)$ that are $W$-equidistant from $c$ that there exists an element $g\in \Fix_G(c)$ such that $gd_1=d_2$. First of all, notice from the hypothesis \ref{Hypothese cox star} and the solution of the word problem in Coxeter groups that every $w\in W$ admits a unique decomposition $w=w_1...w_n$ with $w_t\in W_{I_{k_t}}-\{1_W\}$ such that $I_{k_t}\not= I_{k_{t+1}}$ for every $t$. Suppose that $d_1,d_2\in \Ch(\Delta)$ have $W$-distance $w_1...w_n$ from $c$ with $w_t\in W_{I_{k_t}}$ and let us show the existence of $g$ by induction on $n$. If $n=1$, the result follows from Proposition \ref{theorem pour les right angled buildings} since for every $j\in I_{k_1}$, there exists an element $g_j\in G_j$ such that $g_j\circ h_j(c)=h_j(c)$ and $g_j\circ h_j(d_1)=h_j(d_2)$. Hence, there exists an element $g\in G$ such that $h_j\circ g_j=g\circ h_j$ for every $j\in J$. In particular, $gc=c$, $gd_1=d_2$ and the result follows. If $n\geq2$, we let $d_s'=\proj_{\mathcal{R}_{I_{k_n}}(d_s)}(c)$ and notice that $\delta(c,d_s')=w_1...w_{n-1}$. Our induction hypothesis therefore ensures the existence of a $g'\in  G$ such that $g'c=c$ and $g'd_1'=d_2'$. Now, notice that $\delta(d_2',d_2)=w_n$ and $\delta(d_2',g'd_1)=\delta (g'd_1',g'd_1)=\delta (d_1',d_1)=w_n$. In particular, Proposition \ref{theorem pour les right angled buildings} ensures the existences of an element $h\in G$ such that $hd=d$ for every $d\in \mathcal{R}_{I_{k_{(n-1)}}}(d_2')$ and $hg'd_1=d_2$. Since $G$ satisfies the property ${\rm IP}_{I_{k_{(n-1)}}}$ by Proposition \ref{independence IPJ for universal}, this implies the existence of an element $h'\in \Fix_{G}(V_{I_{k_{(n-1)}}}(d_2'))$ such that $h'b=hb$ for every $b\in X_{I_{k_{(n-1)}}}(d_2')$. Since $c\in V_{I_{k_{(n-1)}}}(d_2')$, the automorphism $h'g'\in G$ satisfies that $h'g'c=c$ and that $h'g'd_1=d_2$. The result follows.
	\end{proof}
	In particular, if $\mathcal{U}((h_i,G_i)_{i\in I})$ is non-discrete, Theorem \ref{theorem Ipv1 letter} applies to $\alpha(G)$ and Theorem \ref{la version paki du theorem de classification} provides a bijective correspondence between the equivalence classes of irreducible representations of $\alpha(G)$ at depth $l\geq 1$ with seed $C\in \mathcal{F}_{\mathcal{S}_{V_1}}$ and the $\mathcal{S}_{V_1}$-standard representations of $\Aut_{\alpha(G)}(C)$. We recall further that an existence criterion for those representations was given in Section \ref{existence for IPV1}. 
	
	Since $\alpha$ is a homeomorphism on its image, the same holds for the representations of $G$. Notice that under the correspondence given by $\alpha^{-1}$, the generic filtration $\mathcal{S}_{V_1}$ describes a generic filtration $\mathcal{S}_\Delta$ of $G$ that factorises$^+$ at all depths $l\geq 1$ and which can be interpreted as follows. We explicit this correspondence below. Let $\delta$ denote the $W$-distance of $\Delta$ and let consider the set $$\mathcal{R}'(c)=\{d\in \Ch(\Delta): \delta(c,d)=w\mbox{ s.t. } \exists k\in \{1,...,r\} \mbox{ for which }w\in W_{I_{k}}\}$$ 
	for every chamber $c\in \Ch(\Delta)$. By use of the correspondence $\Phi: \Ch(\Delta) \rightarrow V_0$ between chambers of $\Delta$ and vertices of $V_0$, notice, for every $c\in V_0$, that  $\Phi^{-1}(B_T(\Phi(c),2 )\cap V_0))=\mathcal{R}'(c)$. We define a family $\mathcal{T}_\Delta$ of subsets of $\Ch(\Delta)$ as follows:
	\begin{enumerate}
		\item $\mathcal{T}_\Delta\lb 0\rb =\{\mathcal{R}_{I_k}(c): k\in \{1,...,r\}, c\in \Ch(\Delta)\}$.
		\item For every $l$ such that $l\geq 0$, we define iteratively:
		\begin{equation*}
		\begin{split}
		\mathcal{T}_\Delta\lb l+1\rb =\{\mathcal{R}\subseteq \Ch(\Delta): \exists \mathcal{Q}\in \mathcal{T}_\Delta\lb l\rb, \qq \exists c\in &\mathcal{Q}\mbox{ s.t. }\mathcal{R}'(c)\not\subseteq \mathcal{Q}\\
		&\mbox{ and }\mathcal{R}=\mathcal{Q}\cup\mathcal{R}'(c)) \}.
		\end{split}
		\end{equation*}
		\item We set $\mathcal{T}_\Delta= \bigsqcup_{l\in \N}\mathcal{T}_\Delta\lb l\rb$.
	\end{enumerate}
	It is quite easy to realise that $\mathcal{S}_\Delta=\{\Fix_G(\mathcal{R}):\mathcal{R}\in \mathfrak{T}_\Delta \}$ is the generic filtration of $G$ corresponding to $\mathcal{S}_{V_1}$ under the correspondence given by $\alpha^{-1}$ and that 
	$$\mathcal{S}_\Delta\lb l\rb=\alpha^{-1}\big(\mathcal{S}_{V_1}\lb l\rb \big)=\{\Fix_G(\mathcal{R}): \mathcal{R}\in \mathcal{T}_\Delta\lb l\rb\}.$$

	\chapter{Nebbia's conjecture and Radu groups}\label{Chapter Radu groups}
	\section{Introduction and main results}
	
	Let $T$ be a $(d_0,d_1)$-semi-regular tree with $d_0,d_1\geq 4$. The closed subgroups $G\leq \Aut(T)$ acting $2$-transitively on the boundary $\partial T$ and whose local action at every vertex contains the alternating group have been extensively studied by N.~Radu in \cite{Radu2017}. Among other things, he obtained a classification and an explicit description of those groups when $d_0,d_1\geq 6$. We therefore refer to them as Radu groups. The purpose of this section is to classify the irreducible representations of Radu groups and to prove that these groups are uniformly admissible. This provides a substantial contribution to an important conjecture formulated in the 90's \cite{Nebbia1999} by C.~Nebbia. We start with some reminders on the related concepts. We recall from Chapter \ref{Chapter I preliminaries} that Type {\rm I} groups are in a sense the locally compact groups whose representation theory is well behaved with respect to irreducible representations since the uniqueness of the direct integral decomposition of a representation into irreducible representations is only guaranteed for Type {\rm I} groups and since the classification of the irreducible representation of a locally compact groups is known to be intractable unless the group is Type {\rm I}. In addition, we recall that CCR locally compact groups are of Type {\rm I} and that a locally compact group  $G$ is \textbf{CCR} if the operator $\pi(f)$ is compact for every irreducible representation $\pi$ of $G$ and each $f\in L^{1}(G)$. When the group $G$ is totally disconnected, this property is equivalent to the fact that the group is \tg{admissible} that is: for every compact open subgroup $K \leq G$ and each irreducible representation $\pi$ of $G$ the space $\Hr{\pi}^K$ of $K$-invariant vectors is finite dimensional \cite{Nebbia1999}. A totally disconnected locally compact group is said to be uniformly admissible if for every compact open subgroup $K \leq G$, there exists a positive integer $k_K$ depending on $K$ such that $\dim(\Hr{\pi}^K)\leq k_K$ for each irreducible representation $\pi$ of $G$. The relations between those notions is given by the following diagram for each totally disconnected locally compact group: 
	$$\mbox{\tg{Uniformly admissible}} \Rightarrow  \mbox{\tg{admissible}} \iff \mbox{\tg{CCR}} \Rightarrow \mbox{\tg{Type} \tg{\rm I}}.$$
	At the end of the $80$'s, the classification of the irreducible representations of the full automorphism group $\Aut(T)$ of any thick regular tree lead to the conclusion that these groups are all uniformly admissible \cite{Olshanskii1980}. A few years later, C.~Nebbia's work highlighted a surprising parallel between the regularity of the unitary dual of any group of automorphisms of a tree and the regularity of its action on the boundary. To be more precise, he proved that the regularity of the representations of a closed subgroups $G\leq \Aut(T)$ implies a regularity of its action on the boundary $\partial T$:
	\begin{theorem*}[{\cite{Nebbia1999}}]
		Let $T$ be a thick regular tree. Then, any closed unimodular CCR vertex-transitive subgroup $G\leq \Aut(T)$ acts transitively on the boundary $\partial T$.
	\end{theorem*} 
	Further progress going in that direction were recently achieved by C.~Houdayer and S.~Raum \cite{HoudayerRaum2019} and with even higher level of generality by P-E.~Caprace, M.~Kalantar and N.~Monod \cite{CapraceKalantarMonod2022}. Among other things, they showed that any closed non-amenable Type {\rm I} subgroup acting minimally on a locally finite tree $T$ must act $2$-transitively on the boundary $\partial T$ \cite[Corollary D]{CapraceKalantarMonod2022}. 
	On the other hand, supported by the classification of spherical and special representations holding for any closed subgroup acting $2$-transitively on the boundary, C.~Nebbia conjectured that the regularity of the action of the group on the boundary should also imply the regularity of its unitary dual \cite{Nebbia1999}. More precisely, he conjectured that every closed subgroup $G\leq \Aut(T)$ of automorphisms of a thick regular tree acting transitively on the boundary $\partial T$ is CCR. This conjecture extends naturally to the framework of locally finite tree. Furthermore, notice that the hypothesis that the tree is semi-regular is non-restrictive in the extension of the conjecture. Indeed, compact groups are CCR and Theorem \ref{Theorem Burger mozes acting on the boundary trans} ensures that the tree is semi-regular if there exists a closed non-compact subgroups $G\leq \Aut(T)$ acting transitively on the boundary $\partial T$. In particular, Nebbia's conjecture extends as follows:
	\begin{conjecture*}[Nebbia's conjecture]
		Let $T$ be a thick semi-regular tree. Every closed subgroup $G\leq \Aut(T)$ acting transitively on the boundary $\partial T$ is CCR. 
	\end{conjecture*}
	
	To the present days, this conjecture is supported by the fact that rank one semi-simple algebraic groups over local-fields \cite{Bernshtein1974}, \cite{Harish1970} and groups satisfying the Tits independence property \cite{Amann2003} are uniformly admissible. The purpose of this chapter is to prove that Radu groups are also uniformly admissible. We recall that those groups are the closed subgroups $G \leq \Aut(T)$ acting transitively on the boundary $\partial T$ and whose local action at every vertex contains the alternating group. As a reminder, the \tg{local action} $\underline{G}(v)$ of $G$ at $v$ is the subgroup of the group $\Sym(E(v))$ of permutations of the set $E(v)$ of edges containing $v$ obtained as the image under the natural projection map of the action of $\Fix_G(v)$ on $E(v)$. It is not hard to prove that Radu groups are of Type {\rm I}. Indeed, each Radu group $G$ contains a closed subgroup $H$ that is conjugate in $\Aut(T)$ to the semi-regular version of the universal group of Burger-Mozes $\Alt_{(i)}(T)^+$ see \cite[page 4]{Radu2017} and such that $G/H$ has a finite non-zero $G$-invariant-measure. Furthermore, this subgroup $H$ is of Type {\rm I} since it is $2$-transitive on the boundary $\partial T$ and satisfies the Tits independence property. It follows from \cite[Theorem 1]{KallmanI1973} that $G$ is of Type {\rm I}. The purpose of these notes is to go further. Among other things, we provide the following substantial contribution to Nebbia's conjecture.  
	\begin{theoremletter}\label{Theorem A}
		Let $T$ be a $(d_0,d_1)$-semi-regular tree with $d_0,d_1\geq 6$. Then, Radu groups are uniformly admissible and hence CCR.
	\end{theoremletter}
	To put this result into the perspective of Radu's paper, we recall that the local action $\underline{G}(v)\leq \Sym(E(v))$ at every vertex $v\in V$ of a closed subgroup $G\leq \Aut(T)$ acting $2$-transitively on the boundary is a $2$-transitive subgroup of $\Sym(E(v))$ see Theorem \ref{Theorem Burger mozes acting on the boundary trans}. On the other hand, \cite[Proposition B.1 and Corollary B.2]{Radu2017} ensure that
	\begin{equation*}
	\Theta= \{d \geq 6 : \mbox{ each }\mbox{2-transitive subgroup of }\Sym(d)\mbox{ contains Alt}(d)\}
	\end{equation*}
	is asymptotically dense in $\N$ in the sense that 
	$$\lim_{n\rightarrow \infty} \frac{\modu{\Theta\cap \{1,2,...,n\}}}{n}=1.$$ The ten smallest elements of $\Theta$ are $34$,
	$35$, $39$, $45$, $46$, $51$, $52$, $55$, $56$ and $58$. All together, this implies the following.
	\begin{theoremletter}\label{Thoerem letter ABCDE}
		Nebbia's conjecture is confirmed for any $(d_0,d_1)$-semi-regular tree with $d_0,d_1\in \Theta$ where $\Theta$ is the above asymptotically dense subset of $\N$.
	\end{theoremletter}
	We recall that each irreducible representation of a closed subgroup $G\leq \Aut(T)$ is either spherical, special or cuspidal. Furthermore, we recall from Sections \ref{Section spherical rep} and \ref{Section special rep} that the spherical and special representations are classified for any closed subgroup $G\leq \Aut(T)$ acting $2$-transitively on the boundary $\partial T$ and hence at the level of generality of the conjecture. On the other hand, taking advantage on the machinery developed in Chapter \ref{Chapter Olshanskii's factor} and of the explicit description of Radu groups, we provide a classification of the cuspidal representations of each simple Radu groups and deduce a description of the irreducible representations of any Radu group via the induction dynamic described in Section \ref{section induction depuis sous group indice 2}. This procedure starts by considering a family of groups $G_{(i)}^+(Y_0,Y_1)$ indexed by finite subsets $Y_0,Y_1\subseteq \N$ see Definition \ref{definition de GiXY} below. Those groups are abstractly simple when $T$ is a $(d_0,d_1)$-semi-regular tree with $d_0,d_1\geq 4$ and they exhaust the list of \textbf{simple Radu} groups when $d_0,d_1\geq 6$. Exploiting the fact that those groups are determined by suitable local conditions, we show that the generic filtration $\mathcal{S}_0$ of $G$ defined  in Section \ref{chapter generic filtration for Gi+YY} factorizes$^+$ at all strictly positive depths. Together with Theorem \ref{la version paki du theorem de classification}, this leads to the following result: 
	\begin{theoremletter}\label{theorem C}
		Let $T$ be a $(d_0,d_1)$-semi-regular tree with $d_0,d_1\geq 4$. Each cuspidal representation of $G_{(i)}^+(Y_0,Y_1)$ admits a non-zero matrix coefficient supported in a compact open subgroup. Furthermore, the cuspidal representations of $G_{(i)}^+(Y_0,Y_1)$ are in bijective correspondence given by processes of induction and inflation with a family of irreducible representations of finite groups. 
	\end{theoremletter}
	This bijective correspondence is explicitly given by Theorem \ref{la classification des reresentations cuspidale Radu} below. Among other things, this theorem proves that the cuspidal representations of $G_{(i)}^+(Y_0,Y_1)$ are square-integrable. In light of \cite[Corollary of Theorem 2]{Harish1970} and of the classification of spherical and special representations, this proves that the groups $G_{(i)}^+(Y_0,Y_1)$ are uniformly admissible and hence CCR (the details are gathered in Section \ref{section radu groups are Type I}). To put those results into the perspective of Radu's classification, we recall that every Radu group $G$ belongs to a finite chain $H_n\geq... \geq H_0$ with $n\in \{0,1,2,3\}$ such that $H_n=G$, $\lb H_t: H_{t-1} \rb =2$ for all $t$ and $H_0$ is conjugate in the group of type-preserving automorphisms $\Aut(T)^+$ to one of those $G_{(i)}^+(Y_0,Y_1)$  when $d_0,d_1\geq 6$. We therefore deduce a description of the irreducible representations of $G$ in terms of the irreducible representations of $H_0$ by the induction dynamic described in Section \ref{section induction depuis sous group indice 2}. In particular, when $d_0,d_1\geq 6$ we obtain a description of the cuspidal representations of any Radu group from the cuspidal representations of the simple Radu groups $G_{(i)}^+(Y_0,Y_1)$ and we deduce that all Radu groups are uniformly admissible. The author would like to underline that each Radu group $G$ satisfies the property \ref{IPk} for some $k\geq 1$ depending on $G$. In particular, Theorem \ref{thm B1} applies directly to these groups. However, unless $G$ satisfies the Tits independence property, that is $k=1$, this second approach never leads to a classification of all the cuspidal representations of $G$. 
	
	\subsection*{Structure of the chapter}
	We recall the classification of Radu groups in Section \ref{section Radu groups}. We describe the cuspidal representations of simple Radu groups in Section \ref{Section classification des rep de simple radu groups}; and we prove that Radu groups are all uniformly admissible in Section \ref{section radu groups are Type I}. 
	
	\section{The Radu groups}\label{section Radu groups}
	The purpose of this section is to introduce properly Radu groups and to recall some results concerning their classifications provided in \cite{Radu2017}. 
	Let $T$ be a $(d_0,d_1)$-semi-regular tree with $d_0,d_1\geq 4$ and with associated bipartition $V=V_0\sqcup V_1$. We denote by $\Aut(T)^+$  the group of type-preserving automorphisms of $T$ that is the set of automorphisms of $T$ leaving $V_0$ and $V_1$ invariants. Consider the following two sets
	\begin{equation*}
	\mathcal{H}_T=\{ G\leq\Aut(T): G\mbox{ is closed and } 2\mbox{-transitive on }\partial T \}
	\end{equation*}
	and 
	\begin{equation*}
	\mathcal{H}_T^+=\{ G\leq\Aut(T)^+: G\mbox{ is closed and } 2\mbox{-transitive on }\partial T \}.
	\end{equation*} 
	If $d_0 \not= d_1$, notice that every automorphisms of $T$ is type-preserving so that $\mathcal{H}_T^+= \mathcal{H}_T$. Now, let $G\leq \Aut(T)$ be a closed subgroup and notice that for each vertex $v\in V$, the stabiliser $\Fix_G(v)$ of $v$ acts on the set $E(v)$ of edges containing $v$. The image of $\Fix_G(v)$ in $\Sym(E(v))$ for the projection map given by this action form a group called the \textbf{local action} of $G$ at $v$ which we denote by  $\underline{G}(v)$. In the light of Theorem \ref{Theorem Burger mozes acting on the boundary trans}, every group $G\in \mathcal{H}_T^+$ is transitive on both $V_0$ and $V_1$. In particular, for each $t\in \{0,1\}$ there exists a subgroup $F_t\leq \Sym(d_t)$ such that each the local actions in $\{\underline{G}(v): v\in V_t\}$ is isomorphic to $F_t$. 
	\begin{definition}
		A \textbf{Radu group} $G$ is an element of $\mathcal{H}_T$ such that $\underline{G}(v)\simeq F_t\geq \Alt(d_t)$ for every vertex $v\in V_t$ and each $t\in \{0,1\}$. 
	\end{definition}
	Radu provided a classification as well as a concrete realisation of those groups. The statement of his result requires some preliminaries. For every vertex $v\in V$ and every positive integer $r\in \N$ let $$S(v,r)=\{w\in V: d(v,w)=r\}$$ be the set of vertices of $T$ at distance $r$ from $v$.
	\begin{definition}
		We define a \tg{legal coloring} $i:V\rightarrow \N$ of the tree $T$ to be the concatenation of a pair of maps $$i_0:V_0\rightarrow\{1,...,d_1\}\mbox{ and }i_1:V_1\rightarrow\{1,...,d_0\}$$ 
		such that $\restr{i_0}{S(v,1)}:S(v,1)\rightarrow \{1,...,d_1\}\mbox{ and }\restr{i_1}{S(w,1)}:S(w,1)\rightarrow \{1,...,d_0\}$ are bijections for all $v\in V_1$  and $w\in V_0$.
	\end{definition} 
	Given a legal coloring $i$ of $T$ and an automorphism $g\in \Aut(T)$, we define the \tg{local action} of $g$ at a vertex $v\in V$ as the permutation 
	\begin{equation*}
	\sigma_{(i)}(g,v)= \restr{i}{S(gv,1)}\circ g\circ \Big(\restr{i}{S(v,1)}\Big)^{-1}\in \begin{cases}
	\Sym(d_0)\q &\mbox{if }v\in V_0\\
	\Sym(d_1)\q &\mbox{if }v\in V_1
	\end{cases}.
	\end{equation*}  
	\begin{remark}
		If $d_0=d_1$, the tree $T$ is regular and these notions of legal coloring and of local action differ from the notions of legal coloring and local actions used  in \cite{BurgerMozes2000} to define the universals Burger Mozes group. Indeed, with our definition, the closed subgroup $G\leq \Aut(T)$ of all the automorphisms of trees $g\in G$ such that $\sigma_{(i)}(g,v)=\id$ $\forall v\in V$ is not vertex-transitive (not even transitive on $V_0$).
	\end{remark}
	Now, let $T$ be a $(d_0,d_1)$-semi-regular tree with $d_0,d_1\geq 4$ and $i$ be a legal coloring of $T$. For every vertex $v\in V$ and every finite subset $Y\subset \N$ let $$S_Y(v)=\bigcup_{r\in Y} S(v,r)$$ and for every set of vertices $B\subseteq V$ set   
	$$\Sgn_{(i)}(g,B)=\prod_{w\in B}^{}\sgn(\sigma_{(i)}(g,w))$$
	where $\sgn(\sigma_{(i)}(g,w))$ is the sign of the local action $\sigma_{(i)}(g,w)$ of the automorphism $g$ at $w$ for the legal coloring $i$. 
	\begin{definition}\label{definition de GiXY}
		For all (possibly empty) finite subsets $Y_0, Y_1$ of $\N$ and every legal coloring $i$ of $T$, we set
		\begin{equation*}
		G^+_{(i)}(Y_0,Y_1)=\left\{g\in \Aut(T)^+\ \middle\vert \begin{array}{l}
		\Sgn_{(i)}(g,S_{Y_0}(v))=1\qq \mbox{ for each }v\in V_{t_0},\\ 
		\Sgn_{(i)}(g,S_{Y_1}(v))=1\qq \mbox{ for each }v\in V_{t_1}
		\end{array}\right\},
		\end{equation*}
		where $t_0=\max(Y_0)\mod2$,  $t_1=(1+\max(Y_1))\mod2$ and $\max(\es)=0$. 
	\end{definition}
	\begin{remark}\label{remark les sommets ont des types opposers}
		The conditions on $t_0$ and $t_1$ are made in such a way that the vertices of $S_{Y_0}(v)$ with $v\in V_{t_0}$ at maximal distance from $v$ and the vertices of $S_{Y_1}(w)$ with $w\in V_{t_1}$ at maximal distance from $w$ have opposite types.
	\end{remark}
	Notice that $G_{(i)}^+(\es,\es)=\Aut(T)^+$ is the full group of type-preserving automorphisms and that $G^+_{(i)}(\{0\},\{0\})$ is a subgroup of each $G^+_{(i)}(Y_0,Y_1)$. Furthermore, if $T$ is a $d$-regular tree notice that $G^+_{(i)}(\{0\},\{0\})$ is conjugate to $U(\Alt(d))^+$ where $G^+=G\cap \Aut(T)^+$ and where $U(\Alt(d))$ is the universal Burger-Mozes group of the alternating group defined in \cite{BurgerMozes2000}. As we recall below, when $d_0,d_1\geq 6$, every simple Radu group is of the form $G_{(i)}^+(Y_0,Y_1)$ for some finite subsets $Y_0,Y_1\subset \N$ and some legal coloring $i$ of $T$. Furthermore, every Radu group $G$ belongs to a finite chain $H_n\geq... \geq H_0$ with $n\in \{0,1,2,3\}$ such that $H_n=G$, $\lb H_t: H_{t-1} \rb =2$ for all $t$ and $H_0$ is conjugate in $\Aut(T)^+$ to one of those $G_{(i)}^+(Y_0,Y_1)$. 

	We now recall some of the statements of \cite{Radu2017} as well as the explicit description of Radu groups following from those. This requires some preliminaries. In what follows, we denote the intersection of all normal cocompact closed subgroups of a locally compact $G$ by $G^{(\infty)}$. We recall from \cite[Proposition 3.1.2]{BurgerMozes2000} that $H^{(\infty)}$ belongs to $\mathcal{H}^+_T$ and is topologically simple for all $H\in \mathcal{H}^+_T$ (in our case these groups are even abstractly simple). Finally, we let $\mathcal{G}_{T}^+(i)$ be the set of groups $G^+_{i}(Y_0,Y_1)$ with non-empty finite subsets $Y_0, Y_1\subseteq\N$ such that $y=\max(Y_t)\mod2$ for each $y\in Y_t$ with $y\geq \max(Y_{1-t})$ ($t\in \{0,1\}$). 
	\begin{theorem*}[{\cite[Theorem A]{Radu2017}}]
		Let $T$ be a $(d_0,d_1)$-semi-regular tree with $d_0,d_1\geq 4$ and $i$ be a legal coloring of $T$. Then, for every finite subsets $Y_0, Y_1\subset \N$, the group $G^+_{(i)}(Y_0,Y_1)$ belongs to $\mathcal{H}^+_T$ and is abstractly simple.
	\end{theorem*}
	\begin{theorem*}[{\cite[Theorem B]{Radu2017}}]
		Let $T$ be a $(d_0,d_1)$-semi-regular tree with $d_0,d_1\geq 6$, $i$ be a legal coloring and $G\in \mathcal{H}^+_T$ be such that there exists two groups of permutations $\Alt(d_t)\leq F_t\leq \Sym(d_t)$ with $t\in \{0,1\}$ and such that $\underline{G}(v)\simeq F_t$ for each $v\in V_t$. Then, $G^{(\infty)}$ is conjugate in $\Aut(T)^+$ to an element of $\mathcal{G}_{T}^+(i)$ and $\lb G: G^{(\infty)}\rb \in \{1,2,4\}$.
	\end{theorem*}
	When $T$ is a $d$-regular tree, a similar result holds for all $G\in \mathcal{H}_T- \mathcal{H}_T^+$.
	\begin{theorem*}[{\cite[Corollary C]{Radu2017}}]
		Let $T$ be a $d$-regular tree with $d\geq 6$, $i$ be a legal coloring and $G\in \mathcal{H}_{T}- \mathcal{H}^+_T$ be such that there exists a permutation $\Alt(d)\leq F \leq \Sym(d)$ with $\underline{G}(v)\simeq F$ for each $v\in V$. Then, $\lb G: G^{(\infty)} \rb\in \{2,4,8\} $ and $G^{(\infty)}$ is conjugate to ${G}^+_{(i)}(Y,Y)$ for some finite subset $Y\subset\N$.
	\end{theorem*}
	In particular, when $d_0,d_1\geq 6$ each Radu group is conjugate in $\Aut(T)^+$ to one of the following groups with $X,Y_0,Y_1$ are finite subsets of $\N$ and $\epsilon_t\in \{-1,1\}$ for all $t\in \{0,1\}$.
	\begin{enumerate} 
		\item \begin{center}
			$G_{(i)}(\es,\es)= \Aut(T)$ 
		\end{center}
		\item \begin{center}
			$G_{(i)}^+(\es,\es)=\Aut(T)^+$
		\end{center}
		\item \begin{center}
			$G_{(i)}(X,X)=\{g\in \Aut(T)\vert \Sgn_{(i)}(g,S_X(v))=1\qq \mbox{for each }v\in V\}$
		\end{center}
		\item \begin{center}
			$G_{(i)}(X,X)^*=\{g\in \Aut(T)\vert \Sgn_{(i)}(g,S_X(v))=\Sgn_{(i)}(g,S_X(v')) \qq \forall v,v'\in V\}$
		\end{center}
		\item \begin{flushleft}
			$\q\q\q\q \q\q\q\q\q\q\q\q G_{(i)}(X^*,X^*)=\left\{ g\in \Aut(T) \qq \middle\vert \begin{array}{l}
			\Sgn_{(i)}(g,S_X(v_0))=\Sgn_{(i)}(g,S_X(v_0'))\qq\forall v_0, v_0'\in V_0, \\ \Sgn_{(i)}(g,S_X(v_1))=\Sgn_{(i)}(g,S_X(v_1')) \qq\forall v_1, v_1'\in V_1
			\end{array}\right\}$
		\end{flushleft}
		\item \begin{flushleft}
			$\q\q\q\q \q\q\q\q \q\q\q\q G_{(i)}'(X^*,X^*)=\left\{ g\in \Aut(T)\ \middle\vert \begin{array}{l}
			\qq \Sgn_{(i)}(g,S_X(v_0))= \epsilon_0 \qq\forall v_0\in V_0,\\
			\qq \Sgn_{(i)}(g,S_X(v_1))=\epsilon_1\qq\forall v_1\in V_1,\\
			\epsilon_0=\epsilon_1\mbox{ if and only if }g\in \Aut(T)^+
			\end{array}\right\}$
		\end{flushleft}
		\item \begin{flushleft}
			$\q\q\q\q \q\q\q\q \q\q\q\q G_{(i)}^+(Y_0,\es)=\left\{g\in \Aut(T)^+\ \middle\vert   \begin{array}{l}
			\Sgn_{(i)}(g,S_{Y_0}(v))=1\qq \forall v\in V_t,\\
			 \q \q   t=\max(Y_0)\qq \text{\rm mod}\qq 2\end{array} \right\}$
		\end{flushleft}
		\item \begin{flushleft}
			$\q\q\q\q \q\q\q\q \q\q\q\q G_{(i)}^+(\es,Y_1)=\left\{g\in \Aut(T)^+\ \middle\vert \begin{array}{l}
			\Sgn_{(i)}(g,S_{Y_1}(v))=1\qq \forall v\in V_t,\\ 
			\q\qq  t=\max(Y_1)+1\qq \text{\rm mod}\qq 2
			\end{array}\right\}$
		\end{flushleft}
		\item \begin{center}
			$G_{(i)}^+(Y_0,Y_1)=G_{(i)}^+(Y_0,\es)\cap G_{(i)}^+(\es,Y_1) $
		\end{center}
		\item \begin{flushleft}
			$\q\q\q\q \q\q\q\q \q\q\q\q G_{(i)}^+(Y_0^*,\es)=\left\{g\in \Aut(T)^+\ \middle\vert \begin{array}{l}
			\Sgn_{(i)}(g,S_{Y_0}(v))=\Sgn_{(i)}(g,S_{Y_0}(v'))\qq \forall v,v'\in V_t,\\
			  \q \q \q \q\q \q t=\max(Y_0)\qq \text{\rm mod}\qq 2
			\end{array}\right\}$
		\end{flushleft}
		\item \begin{flushleft}
			$\q\q\q\q \q\q\q\q \q\q\q\q G_{(i)}^+(\es,Y_1^*)=\left\{g\in \Aut(T)^+\ \middle\vert \begin{array}{l}
			\Sgn_{(i)}(g,S_{Y_1}(v))=\Sgn_{(i)}(g,S_{Y_1}(v'))\qq \forall v,v'\in V_t,\\
			\q \q \q \q\q \q t=\max(Y_1)+1\qq \text{\rm mod}\qq 2
			\end{array}\right\}$
		\end{flushleft}
		\item \begin{flushleft}
			$\q\q\q\q \q\q\q\q \q\q\q\q G_{(i)}^+(Y_0^*,Y_1^*)=\left\{g\in \Aut(T)\ \middle\vert \begin{array}{l}
			\Sgn_{(i)}(g,S_{Y_0}(v_0))=\Sgn_{(i)}(g,S_{Y_0}(v_0'))\qq \forall v_0,v_0'\in V_{t_0},\\
			\Sgn_{(i)}(g,S_{Y_1}(v_1))=\Sgn_{(i)}(g,S_{Y_1}(v_1'))\qq \forall v_1,v_1'\in V_{t_1},\\
			\q t_0=\max(Y_0)\qq \text{\rm mod}\qq 2,\qq t_1=\max(Y_1)+1\qq \text{\rm mod}\qq 2
			\end{array}\right\}$
		\end{flushleft}
	\end{enumerate}
	Furthermore, when $T$ is a $(d_0,d_1)$-semi-regular tree with $d_0,d_1\geq 4$, each of the above group is a closed subgroup $G\leq \Aut(T)$ acting $2$-transitively on the boundary $\partial T$ and whose local action at every vertex contains the alternating group of corresponding degree. The following statement follows directly from this description and provides the key step of the classification of irreducible representations of Radu groups presented in Section \ref{Section classification des rep de simple radu groups}. 
	\begin{theorem}\label{Corollary Radu simple then Radu}
		Let $T$ be a $(d_0,d_1)$-semi-regular tree with $d_0,d_1\geq 6$ and $i$ be a legal coloring of $T$. Every Radu group $G$ is contained in a finite chain $H_n\geq... \geq H_0$ with $n\in \{0,1,2,3\}$ such that $H_n=G$, $\lb H_t: H_{t-1} \rb =2$ for all $t$ and $H_0$ is conjugate in $\Aut(T)^+$ to some $G_{(i)}^+(Y_0,Y_1)$. 
	\end{theorem}
		\section{Cuspidal representations of the simple Radu groups}\label{Section classification des rep de simple radu groups}
	Let $T$ be a $(d_0,d_1)$-semi-regular tree with $d_0,d_1\geq 4$,  $V=V_0\sqcup V_1$ be the associated bipartition, $i$ be a legal coloring of $T$ and  $Y_0,Y_1\subseteq \N$ be two finite subsets. We recall from Section \ref{section Radu groups} that $G_{(i)}^+(Y_0,Y_1)$ is a simple Radu group and that every simple Radu group is of this form when $d_0,d_1\geq 6$. Our current purpose is to describe the irreducible representations of these groups and to prove that they are uniformly admissible. Since the group $G_{(i)}^+(Y_0,Y_1)$ acts $2$-transitively on the boundary, a classification of its spherical and special representations is given in Sections \ref{Section spherical rep} and \ref{Section special rep}. We now provide a classification of their cuspidal representations by relying on the explicit description of these groups and on the machinery developed in Chapter \ref{Chapter Olshanskii's factor}.
	
	\subsection{Generic filtration for the simple Radu groups}\label{chapter generic filtration for Gi+YY}
	
	Let $T$ be a $(d_0,d_1)$-semi-regular tree with $d_0,d_1\geq 4$, $V=V_0\sqcup V_1$ be the associated bipartition, $i$ be a legal coloring of $T$, consider two finite subsets $Y_0,Y_1\subseteq \N$. The purpose of this section is to explicit a generic filtration for $G_{(i)}^+(Y_0,Y_1)$. We adopt the same formalism as in Section \ref{application IPk}. Let $\mathfrak{T}_0$ be the family of subtrees of $T$ defined by
	$$\mathfrak{T}_{0}=\{B_T(v,r): v\in V, r\geq 1\}\sqcup\{B_T(e,r): e\in E, r\geq 0\}.$$
	Let $G\leq \Aut(T)$  be a closed subgroup and consider the basis of neighbourhoods of the identity given by the pointwise stabilisers of those trees
	\begin{equation*}
	\mathcal{S}_{0}=\{\Fix_G(\mathcal{T}): \mathcal{T}\in \mathfrak{T}_{0}\}.
	\end{equation*}
	We recall from Section \ref{application IPk} that $G$ satisfies the hypothesis \ref{Hypothese H0} if for all $\mathcal{T},\mathcal{T'}\in \mathfrak{T}_{0}$ we have that
		\begin{equation}\tag{$H_0$}\label{Hypothese H0}
		\Fix_G(\mathcal{T}')\subseteq \Fix_G(\mathcal{T})\mbox{ if and only if }\mathcal{T}\subseteq \mathcal{T}'. 
		\end{equation}
	In addition, the following lemma follows directly from Lemma \ref{la forme des Sl pour IPk}.
	\begin{lemma}\label{Lemme la forme des elements de Sl}
		Let $G\leq \Aut(T)$ be a closed non-discrete unimodular subgroup satisfying the hypothesis \ref{Hypothese H0}. Then, $\mathcal{S}_{0}$ is a generic filtration of $G$ and the sets $\mathcal{S}_{0}\lb l \rb$ can be described as follows:
		\begin{itemize}
			\item If $l$ is even $\mathcal{S}_{0}\lb l\rb = \{\Fix_G(B_T(e,\frac{l}{2})):e\in E\}.$
			\item If $l$ is odd $\mathcal{S}_{0}\lb l \rb=\{ \Fix_G(B_T(v,(\frac{l+1}{2}))): v\in V\}.$
		\end{itemize} 
	\end{lemma}
	\noindent We come back to our case  $G=G^+_{(i)}(Y_0,Y_1)$. 
	\begin{lemma}\label{lemma G satisfies hypothesis H0}
		The group $G^+_{(i)}(Y_0,Y_1)$ satisfies the hypothesis \ref{Hypothese H0}.
	\end{lemma}
	\begin{proof}
		For every set $X$, let $\Delta_X=\{(x,x): x\in X\}\subseteq X\times X$. We choose two functions 
		\begin{equation*}
		\psi_0 \fct{\{1,...,d_0\}\times \{1,...,d_0\}-\Delta_{ \{1,...,d_0\}}}{\Sym(d_0)}{(k,l)}{\psi_0(k,l)}
		\end{equation*}
		\begin{equation*}
		\psi_1 \fct{\{1,...,d_1\}\times\{1,...,d_1\}-\Delta_{ \{1,...,d_0\}}}{\Sym(d_1)}{(k,l)}{\psi_1(k,l)}
		\end{equation*}
		such that $\psi_t(k,l)$ is a non-trivial element of $\Alt(d_t)$ which fixes $k$ but not $l$. Notice that the existence of those functions is guaranteed from the fact that  $d_0,d_1\geq 4$. For brevity we denote by $G$ be the group $G^+_{(i)}(Y_0,Y_1)$. Let $\mathcal{T},\mathcal{T}'$ be two subtrees of $\mathfrak{T}_0$. If $\mathcal{T}\subseteq \mathcal{T}'$, we clearly have that $\Fix_G(\mathcal{T}')\subseteq \Fix_G(\mathcal{T})$. Now, suppose that $\mathcal{T}\not \subseteq \mathcal{T}'$. In order to prove that $G$ satisfies the hypothesis \ref{Hypothese H0}, we need to show that $\Fix_G(\mathcal{T}')\not \subseteq \Fix_G(\mathcal{T})$. Since $\mathcal{T}\not \subseteq \mathcal{T}'$, there exists a vertex $v\in \mathcal{T}$ that does not belong to $\mathcal{T}'$. Let $\gamma$ be the smallest geodesic from $v$ to $\mathcal{T}'$, let $v'$ be the vertex of $\gamma$ that is adjacent to $v$ and $t\in\{0,1\}$ be such that $v'\in V_t$. Let $w$ be the neighbour of $v'$ which is the closest to $\mathcal{T}'$. Notice that this vertex exists and is unique since $\mathcal{T}$ and $\mathcal{T}'$ are complete. The definition of $v'$ ensures that $\mathcal{T}'\subseteq T(v',v)=\{x\in V: d_T(x,v')<d_T(x,v)\}$. Now, notice the existence of an automorphism $g\in \Fix_{\Aut(T)^+}(T(v',v))$ such that $\sigma_{(i)}(g,v')=\psi_t(i(w),i(v))$ and $\sigma_{(i)}(g,x)$ is even for every $x\in V$. Notice that such an element belongs to $G^+_{(i)}(Y_0,Y_1)\cap \Fix_{\Aut(T)^+}(\mathcal{T}')$ but does not fix $v$ by construction which implies that $g\not\in \Fix_G(\mathcal{T})$.  		
	\end{proof}
	\noindent It follows from Lemma \ref{Lemme la forme des elements de Sl} that $\mathcal{S}_0$ is a generic filtration of $G^+_{(i)}(Y_0,Y_1)$. 
	\subsection{Factorisation of the generic filtration}\label{Section factorization}
	
	Let $T$ be a $(d_0,d_1)$-semi-regular tree with $d_0,d_1\geq 4$, $V=V_0\sqcup V_1$ be the associated bipartition, $i$ be a legal coloring of $T$, consider two finite subsets $Y_0,Y_1\subseteq \N$ and let $G=G_{(i)}^+(Y_0,Y_1)$. We have shown in Section \ref{chapter generic filtration for Gi+YY}, that $\mathcal{S}_0$ is a generic filtration $G$. The purpose of this section is to prove that this generic filtration factorises$^+$ at all depths $l\geq 1$.
	
	We start by defining some notations that will be used in the proof. For every two distinct vertices $v,w\in V$, let $\lb v,w\rb$ be the unique geodesic between $v$ and $w$. Suppose that $d(v,w)=n$, let $v=v_0,v_1,...,v_n=w$ be the sequence of vertices corresponding to $\lb v,w\rb$ in $T$, let 
	\begin{equation*}\label{definition Proj sur la geor}
	p_{\lb v,w\rb}\fct{\lb v,w\rb -\{v\}}{\lb v,w\rb}{v_i}{v_{i-1}}
	\end{equation*}
	and let
	\begin{equation*}
	\begin{split}
	T(v,w)&= \{x\in V: d_T(x,p_{\lb v,w\rb }(w))<d_T(x,w)\}\\
	&= \{x\in V: d_T(x,v_{n-1})<d_T(x,w)\}.
	\end{split}
	\end{equation*}
	\begin{figure}[H]\label{drawing00}\caption{The set $T(v,w)$}
		\begin{center}
			\includegraphics[scale=0.09]{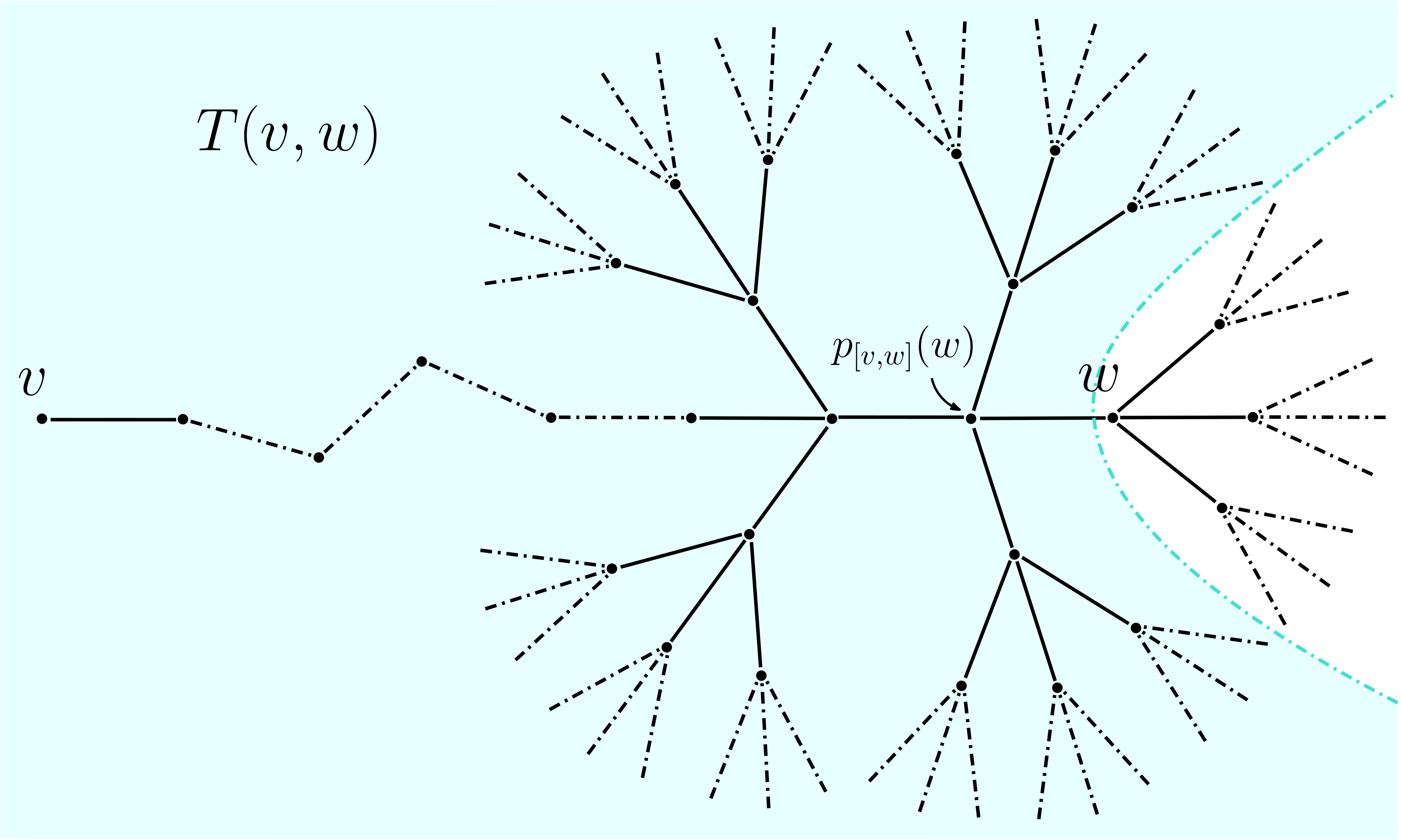}
		\end{center}
	\end{figure}
	\noindent The following intermediate result is the key ingredient required to prove the factorisation of the generic filtration $\mathcal{S}_0$ of $G_{(i)}^+(Y_0,Y_1)$ at all depths $l\geq 1$.
	\begin{proposition}\label{Proposition la premiere etape de la fcatorization pour GY1Y_2}
		For all $l,l'\in \N$ such that $l\geq 1$ and $l'\geq l$, for all $U$ in the conjugacy class of an element of $\mathcal{S}_0\lb l \rb$ and every $V$ in the conjugacy class of an element of $\mathcal{S}_0\lb l'\rb$ such that $V\not \subseteq U$, there exists a subgroup $W\in\mathcal{S}_0\lb l-1\rb$ such that $U\subseteq W\subseteq VU$. 
	\end{proposition}
	\begin{proof}
		To shorten the proof and for clarity of the argument, parts of the reasoning are proved in Lemmas \ref{Lemme Omega 0 est non vide} and \ref{Lemme Omega n est non vide} below. Since the proof is quite long and technical, we start by giving an idea of its structure. We begin the proof by identifying the group $W$ from $U$ and $V$. We then prove that each element of $W$ decomposes as a product of an element of $V$ and an element of $U$. The proof of this decomposition is where the technicalities come from. It is achieved by a compactness argument taking advantage from the fact that $G_{(i)}^+(Y_0,Y_1)$ is defined by local actions conditions.
		
		Let $G= G_{(i)}^+(Y_0,Y_1)$. As announced at the beginning of the proof, we start by identifying $W$. Notice that $\mathfrak{T}_0$ is stable under the action of $G$. Furthermore, for every $g\in \Aut(T)^+$ and for every subtree $\mathcal{T}$ of $T$, we have that $g\Fix_G(\mathcal{T})g^{-1}=\Fix_G(g\mathcal{T})$. In particular, there exist $\mathcal{T}$, $\mathcal{T}'\in \mathfrak{T}_0$ such that $U=\Fix_G(\mathcal{T})$ and $V=\Fix_G(\mathcal{T}')$. Since $V\not \subseteq U$, notice that $\mathcal{T}\not \subseteq \mathcal{T'}$. If $l$ is even, Lemma \ref{Lemme la forme des elements de Sl} ensures that $\mathcal{T}=B_T(e,\frac{l}{2})$ for some edge $e\in E$. Furthermore, since $\mathcal{T}\not \subseteq \mathcal{T'}$ and since $l'\geq l$, there exists a unique vertex $v\in e$ such that $\mathcal{T}'\subseteq T(v,w)\cup B_T(v,\frac{l}{2})$ where $w$ denotes the other vertex of $e$. In that case, we let $\mathcal{T}_W=B_T(v,\frac{l}{2})$. If on the other hand $l$ is odd, Lemma \ref{Lemme la forme des elements de Sl} ensures that $\mathcal{T}=B_T(w,\frac{l+1}{2})$ for some vertex $w\in V$. Furthermore, since $\mathcal{T}\not \subseteq \mathcal{T'}$ and since $l'\geq l$, there exists a unique vertex $v\in B_T(w,1)-\{w\}$ such that $\mathcal{T}'\subseteq T(v,w)\cup B_T(\{v,w\},\frac{l-1}{2})$. In that case, we let $\mathcal{T}_W=B_T(\{v,w\},\frac{l-1}{2})$. In both cases, we set $W=\Fix_G( \mathcal{T}_W)$. By construction, notice that $W\in \mathcal{S}_0\lb l-1\rb$ and that $U\subseteq W$ (since $\mathcal{T}_W\subseteq \mathcal{T}$). Our purpose is therefore to show that $W\subseteq V U$. Let $\alpha\in W$ and let us prove the existence of an element $\alpha_0\in U$ such that $\restr{\alpha}{\mathcal{T}'}=\restr{\alpha_0}{\mathcal{T}'}$. We start by explaining why the existence of $\alpha_0$ settles the proof. Indeed, if $\alpha_0$ exists, notice that the automorphism $\alpha_1=\alpha_0^{-1}\circ \alpha$ is an element of $G$ for which $\restr{\alpha_1}{\mathcal{T}'}=\restr{id}{\mathcal{T}'}$. In particular, we have that $\alpha_1\in \Fix_G(\mathcal{T}')$, $\alpha_0\in \Fix_G(\mathcal{T})$ and by construction $\alpha=\alpha_0 \circ \alpha_1$ which proves that $W\subseteq UV$. Applying the inverse map on both sides of the inclusion we obtain that $W\subseteq VU$ which settles the proof.  
		
		Now, let us prove the existence of $\alpha_0$. As announced, at the beginning of the proof, we are going to use a compactness argument taking advantage from the fact that $G_{(i)}^+(Y_0,Y_1)$ is defined by local actions conditions. To be more precise, we are going to define a descending chain of non-empty compact sets $\Omega_n\subseteq \Aut(T)^+$ and an increasing chain of finite subtrees $\mathcal{R}_n$ of $T$ such that $T=\bigcup_{n\in \N} \mathcal{R}_n$ and such that for all $h\in \Omega_n$ we have:
		\begin{enumerate}[label=(\roman*)]
			\item $h\in \Fix_G(\mathcal{T})$ and $\restr{h}{\mathcal{T}'}=\restr{\alpha}{\mathcal{T}'}$
			\item $\Sgn_{(i)}(h,S_{Y_0}(v))=1$ for all $v$ in $V_{t_0}\cap \mathcal{R}_n$.
			\item $\Sgn_{(i)}(h,S_{Y_1}(v))=1$ for all  $v$ in $V_{t_1}\cap \mathcal{R}_n$.
		\end{enumerate} 
		We recall that in the above $t_0=\max(Y_0)\mod2$,  $t_1=(1+\max(Y_1))\mod2$ and $\max(\es)=0$. Let us first show that this settles the existence of $\alpha_0$. Since the $\Omega_n$ form a descending chain of non-empty compact sets in a Hausdorff space we obtain $\bigcap_{n\in \N} \Omega_n\not =\es$. Let $\alpha_0 \in \bigcap_{n\in \N} \Omega_n$. Since $\alpha_0\in \Omega_0$, notice that $ \restr{\alpha_0}{\mathcal{T}}=\restr{\id}{\mathcal{T}}$, $\restr{\alpha_0}{\mathcal{T}'}=\restr{\alpha}{\mathcal{T}'}$. To see that $\alpha_0$ is as desired, we are left to show that $\alpha_0\in G_{(i)}^+(Y_0,Y_1)$. However, for every $v\in V_{t}(T)$, there exists a positive integer $n\in \N$ such that $v\in \mathcal{R}_n$ and since $\alpha_0\in \Omega_{n}$ we have that $\Sgn_{(i)}(\alpha_0,S_{Y_t}(v))=1$. This proves that $\alpha_0$ is as desired.
		
		We are left to define the descending chain of non-empty compact sets $\Omega_n\subseteq \Aut(T)^+$. Suppose that $\max(Y_0)\leq \max(Y_1)$ (the proof for $\max(Y_1)\leq \max(Y_0)$ is similar). Let $\gamma$ be the smallest geodesic of $T$ containing both the centre of $\mathcal{T}$, the centre of $\mathcal{T}'$ and oriented from $\mathcal{T}$ to $\mathcal{T}'$ (notice that the center is either a vertex or an edge depending on the values of $l$ and $l'$). Since $\mathcal{T}\not \subseteq \mathcal{T}'$ and since $l'\geq l$ notice that $\gamma $ contains at least two vertices. 
		\begin{figure}[H]\label{drawing01}\caption{The tree $\mathcal{T}_W$ and the geodesic $\gamma$}
			\begin{center}
				\includegraphics[scale=0.07]{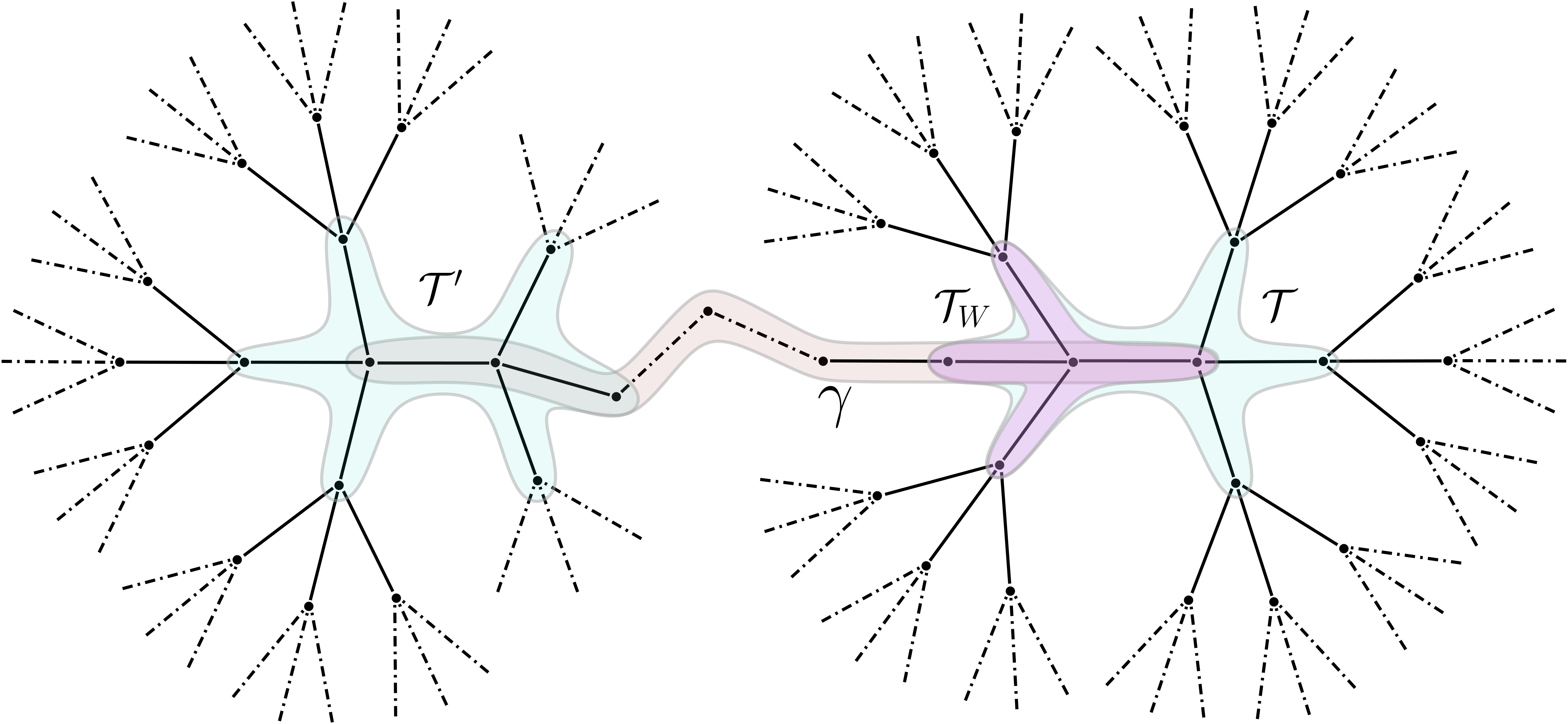}
			\end{center}
		\end{figure}
		The increasing chain of finite subtrees $\mathcal{R}_n$ of $T$ such that $T=\bigcup_{n\in \N} \mathcal{R}_n$ that we are going to use is $\mathcal{R}_n=B_T(\gamma,n)$. We let
		\begin{equation*}
		\Omega_{-1}\qq=\qq \left\{ h \in \Aut(T)^+ \middle\vert \begin{array}{l}\restr{h}{\mathcal{T}}=\restr{\id}{\mathcal{T}} \mbox{ and } \restr{h}{\mathcal{T}'}=\restr{\alpha}{\mathcal{T}'} 
		\end{array}\right\}.
		\end{equation*}
		Since $\alpha\in \Fix_G(\mathcal{T}_W)$ and since $\mathcal{T}_W$ contains every vertices of $\mathcal{T}\cap \mathcal{T}'$ notice that $\Omega_{-1}$ is not empty. Now, since $\max(Y_0)\leq \max(Y_1)$, notice that there exists a unique $r\in \N$ such that $\max(Y_0)+2r\leq \max(Y_1)\leq \max(Y_0)+2r +1$ (where one of this inequality is an equality). We let
		\begin{equation*}
		\Omega_0=\left\{h\in \Omega_{-1}\ \middle\vert \begin{array}{l}
		\Sgn_{(i)}(g,S_{Y_0}(v))=1\qq \mbox{ for each }v\in B_T(\gamma,2r)\cap V_{t_0},\\ 
		\Sgn_{(i)}(g,S_{Y_1}(v))=1\qq \mbox{ for each }v\in B_T(\gamma,0)\cap V_{t_1}
		\end{array}\right\}.
		\end{equation*} 
		Lemma \ref{Lemme Omega 0 est non vide} below ensures that this set is not empty. From there, we define the sets $\Omega_n$ by induction on $n$. For every $n\geq 1$, let $h_n$ an element of $\Omega_{n-1}$ and let
		\begin{equation*}
		\Omega_n\qq=\qq \left\{ h \in \Omega_{n-1}\qq \ \middle\vert \begin{array}{l}
		\q\q \restr{h}{B_T(\gamma,n+{\max(Y_1)})}=\restr{h_{n}}{B_T(\gamma,n+{\max(Y_1)})},\\
		\qq\Sgn_{(i)}(h,S_{Y_0}(w))=1\q \forall w\in B_T(\gamma,n+2r)\cap V_{t_0},\\
		\qq\Sgn_{(i)}(h,S_{Y_1}(w))=1\q \forall w\in B_T(\gamma,n)\cap V_{t_1}
		\end{array}\right\}.
		\end{equation*} 
		For this induction to make sense, it is important for $\Omega_{n}$ to be not empty for all $n\geq 1$. This is proved by Lemma \ref{Lemme Omega n est non vide} below which ensures that $\Omega_{n}$ is a non-empty compact set. The result follows.
	\end{proof}
	
	Our current purpose is to prove Lemmas \ref{Lemme Omega 0 est non vide} and \ref{Lemme Omega n est non vide}. To this end, we introduce some formalism that will be used in both proofs. For all $v\in V$, we are going to need an automorphism $h_{(v)}\in \Aut(T)^+$ that will be used to create an element of $\omega_{n+1}$ from an element of $\Omega_n$. We start by choosing four functions:
	\begin{equation*}
	\phi_0 \fct{\{1,...,d_0\}}{\Sym(d_0)}{k}{\phi_0(k)}
	\end{equation*}
	\begin{equation*}
	\phi_1 \fct{\{1,...,d_1\}}{\Sym(d_1)}{k}{\phi_1(k)}
	\end{equation*}
	\begin{equation*}
	\tilde{\phi}_0 \fct{\{1,...,d_0\}\times\{1,...,d_0\}}{\Sym(d_0)}{(k,l)}{\tilde{\phi}_0(k,l)}
	\end{equation*}
	\begin{equation*}
	\tilde{\phi}_1 \fct{\{1,...,d_1\}\times\{1,...,d_1\}}{\Sym(d_1)}{(k,l)}{\tilde{\phi}_1(k,l)}
	\end{equation*}
	such that $\phi_t(k)$ is an odd permutation of $\Sym(d_t)$ which fixes $k$ and $\tilde{\phi}_t(k,l)$ is an odd permutation of $\Sym(d_t)$ which fixes $k$ and $l$. 
	
	If $v\in V-\gamma$ we choose $w\in \gamma$ and let $h_{(v)}\in \Aut(T)^+$ be such that:
	\begin{enumerate}
		\item $h_{(v)}\in \Fix_{\Aut(T)^+}(T(p_{\lb w,v\rb }(v),v))$.
		\item $\sigma_{(i)}(h_{(v)},v)=\phi_t(i(p_{[w,v]}(v)))$ where $t\in\{0,1\}$ is such that $v\in V_t$.
	\end{enumerate}
	Notice that for all $v\in V-\gamma$ and every $w,w'\in \gamma$, $p_{\lb w,v\rb}(v) = p_{\lb w',v\rb}(v)$ (the definition of $p_{\lb w,v\rb}$ is on page \pageref{definition Proj sur la geor}) so that our choice of $w\in \gamma$ does not change the two properties that $h_{(v)}$ must satisfy. 
	\begin{figure}[H]\label{drawing02}\caption{The automorphism $h_{(v)}$}
		\begin{center}
			\includegraphics[scale=0.08]{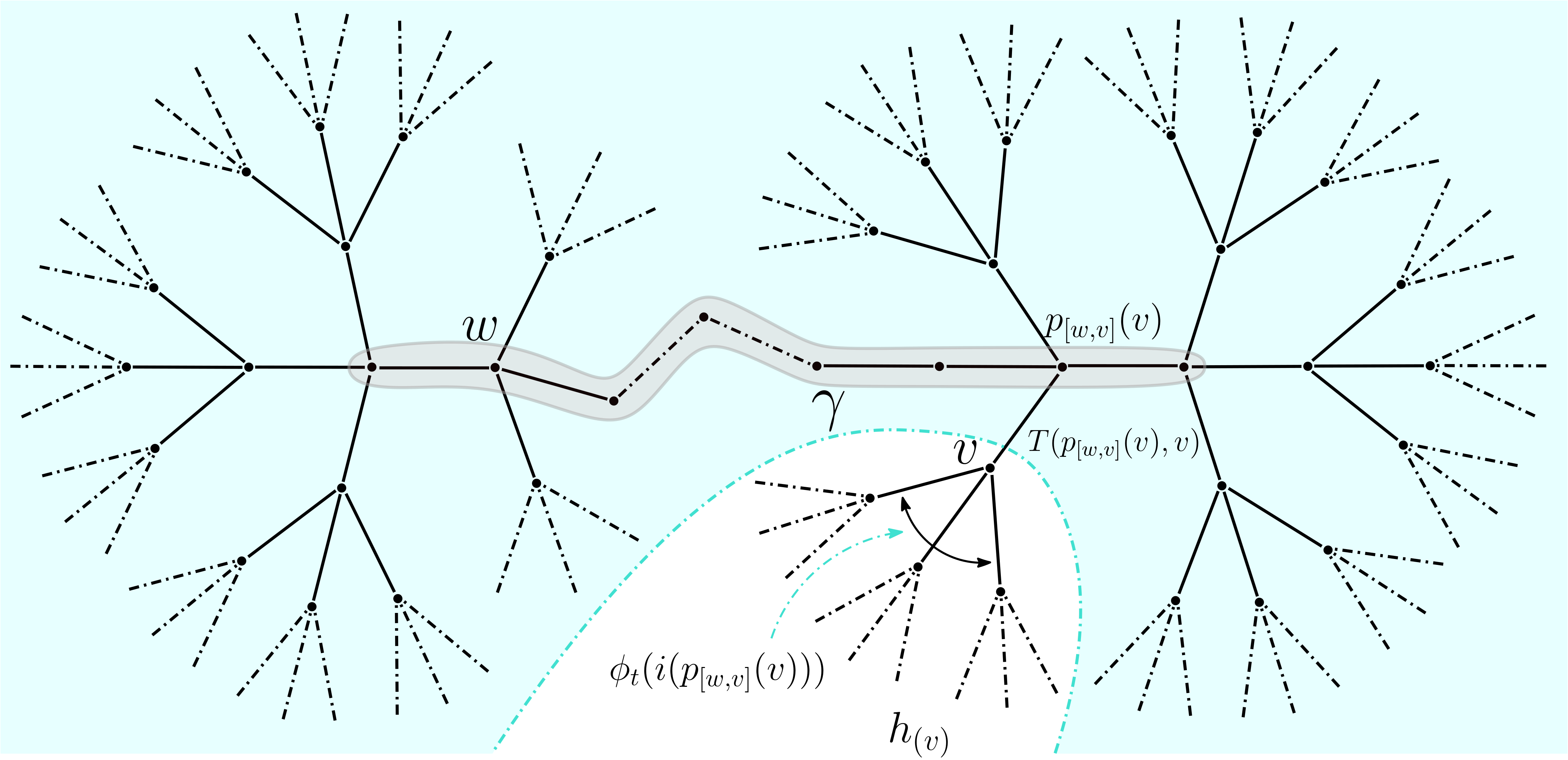}
		\end{center}
	\end{figure}
	
	If $v\in \gamma$, we have two cases. Remember that $\gamma$ has at least two vertices. If $v$ is an end of $\gamma$ let $w$ be the unique vertex of $\gamma$ that is adjacent to $v$ and choose an automorphism $h_{(v)}\in \Aut(T)^+$ such that:
	\begin{enumerate}
		\item $h_{(v)}\in \Fix_{\Aut(T)^+}(T(w,v))$.
		\item $\sigma_{(i)}(h_{(v)},v)=\phi_t(i(w))$ where $t\in\{0,1\}$ is such that $v\in V_t$.
	\end{enumerate}
	On the other hand, if $v$ is not an end of $\gamma$, let $w_1,w_2$ be the two neighbours of $v$ which belong to $\gamma$ and choose an automorphism $h_{(v)}\in \Aut(T)^+$ such that:
	\begin{enumerate}
		\item $h_{(v)}\in\Fix_{\Aut(T)^+}(T(w_1,v)\cup T(w_2,v))$.
		\item $\sigma_{(i)}(h_{(v)},v)=\tilde{\phi}_t(i(w_1),i(w_2))$ where $t\in\{0,1\}$ is such that $v\in V_t$.
	\end{enumerate}
	We are now ready to prove  Lemmas \ref{Lemme Omega 0 est non vide} and \ref{Lemme Omega n est non vide}.
	\begin{lemma}\label{Lemme Omega 0 est non vide}
		The set
		\begin{equation*}
		\Omega_0=\left\{h\in \Omega_{-1}\ \middle\vert \begin{array}{l}
		\Sgn_{(i)}(g,S_{Y_0}(v))=1\qq \mbox{ for each }v\in B_T(\gamma,2r)\cap V_{t_0},\\ 
		\Sgn_{(i)}(g,S_{Y_1}(v))=1\qq \mbox{ for each }v\in B_T(\gamma,0)\cap V_{t_1}
		\end{array}\right\}.
		\end{equation*} 
		is not empty.
	\end{lemma}
	\begin{proof}
		We recall that $\max(Y_0)\leq \max(Y_1)$, that $r\in \N$ is the unique integer such that $\max(Y_0)+2r\leq \max(Y_1)\leq \max(Y_0)+2r +1$, that $t_0=\max(Y_0) \mod 2$ and that $t_1=\max(Y_1)+1 \mod 2$. Remember from the proof of Proposition \ref{Proposition la premiere etape de la fcatorization pour GY1Y_2} that $\Omega_{-1}$ is not empty and let $h_0\in \Omega_{-1}$. We are going to modify the element $h_0$ with the automorphisms $h_{(v)}$ in order to obtain an element of $\Omega_{0}$. A concrete example of the procedure is given on a $4$-regular tree with $Y_0=\{0\}$ and $Y_1=\{1,2\}$ by figures \ref{drawing03}, \ref{drawing04} and \ref{drawing05}. In those figures:
		\begin{itemize}
			\item The hollow vertices are those concerned by the current and previous steps.
			\item The vertices circled in purple are the vertices for which we desire to change the sign $\Sgn_{(i)}(g,S_{Y_0}(v))$ or $\Sgn_{(i)}(g,S_{Y_1}(v))$ (depending on the step) without affecting the sign of other hollow vertices.
			\item The vertices circled in yellow are the vertices for which a change of the local action is applied in order to fulfil the desired change of sign (note that for our choice $Y_0=\{0\}$ those vertices are also the vertices circled in purple). 
		\end{itemize}
		Let $\{w_{0,0},...,w_{0,m_0}\}$ be the set of vertices $w\in B_T(\gamma,0)\cap V_{t_0}$ such that $\Sgn_{(i)}(h_0,S_{Y_0}(w))=-1$. For all $j=0,1,...,m_0$ we choose a vertex $$v_{0,j}\in \bigcap_{w\in \gamma-\{w_{0,j}\}}{}T(w_{0,j},w)$$ such that $d(v_{0,j},w_{0,j})=\max(Y_0)$. 	
		In particular, notice that $v_{0,j}\in S_{Y_0}(w_{0,j})$ but that $v_{0,j}\not \in S_{Y_0}(w)$ for every $w \in B_T(\gamma,0)\cap V_{t_0} $. Furthermore, since $\Sgn_{(i)}(h_0,S_{Y_0}(w_{0,j}))=-1$, $\restr{h_0}{\mathcal{T}}= \restr{\id}{\mathcal{T}}$, $\restr{h_0}{\mathcal{T}'}= \restr{\alpha}{\mathcal{T}'}$ and due to the form of $\mathcal{T}$ and $\mathcal{T}'$, the vertices $v_{0,j}$ must be such that the automorphisms $h_{(v_{0,0})},..., h_{(v_{0,m_0})}$ fix $\mathcal{T}\cup \mathcal{T}'$ pointwise. In particular, the automorphism 
		$$h_{0,0}= h_0 \circ h_{(v_{0,0})}\circ ...\circ h_{(v_{0,m_0})}$$ satisfies  $\restr{h_{0,0}}{\mathcal{T}}=\restr{h_0}{\mathcal{T}}=\restr{\id}{\mathcal{T}}$, $ \restr{h_{0,0}}{\mathcal{T}'}=\restr{h_0}{\mathcal{T}'}=\restr{\alpha}{\mathcal{T}'}$ and $$\Sgn_{(i)}(h_{0,0},S_{Y_0}(w))=1\qq \forall w\in \gamma\cap V_{t_0}.$$
		\begin{figure}[H]\caption{Step I of the proof of Lemma \ref{Lemme Omega 0 est non vide}} \label{drawing03}
			\includegraphics[scale=0.08]{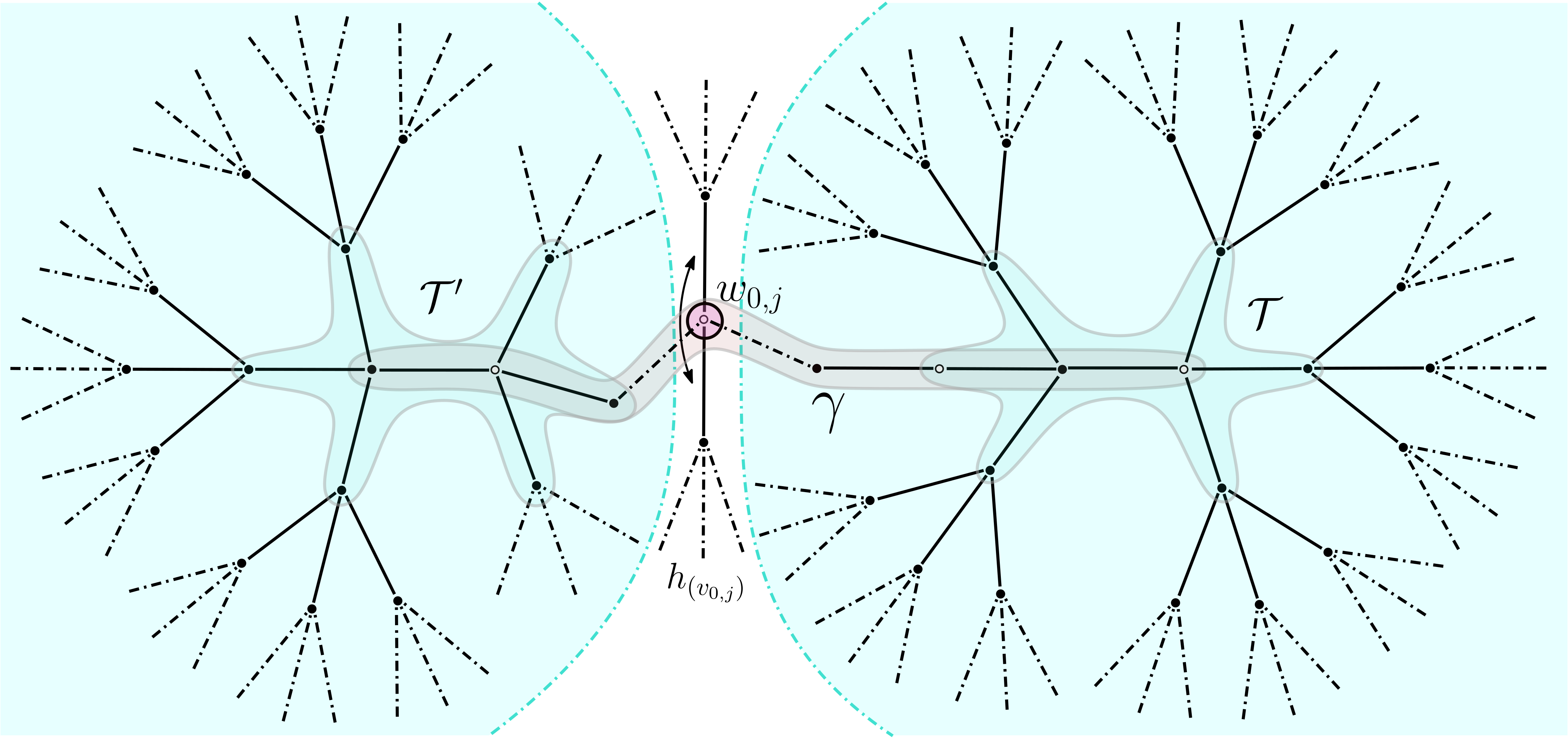}
		\end{figure}
		
		If $r\not=0$, we iterate this procedure. For every $1\leq \nu\leq 2r$, let $\{w_{\nu,0},...,w_{\nu,m_{\nu}}\}$ be the set of vertices $w\in B_T(\gamma,\nu)\cap V_{t_0}$ such that $$\Sgn_{(i)}(h_{\nu-1,0},S_{Y_0}(w))=-1.$$ 
		For all $j=0,1,...,m_\nu$ we choose a vertex $$v_{\nu,j}\in \bigcap_{w\in B_T(\gamma,\nu)\cap V_{t_0}-\{w_{\nu,j}\}}{}T(w_{\nu,j},w)$$ such that $d(v_{\nu,j},w_{\nu,j})=\max(Y_0)$. Hence, notice that $v_{\nu,j}\in S_{Y_0}(w_{\nu,j})$ but that $v_{\nu,j}\not \in S_{Y_0}(w)$ for every $w\in B_T(\gamma,\nu)\cap V_{t_0}-\{w_{\nu,j}\}$. Furthermore, notice that the automorphisms $h_{(v_{\nu,0})}, ..., h_{(v_{\nu,m_\nu})}$ fix $\mathcal{T}\cup \mathcal{T}'$ pointwise. In particular, the automorphism  
		$$h_{\nu,0}=h_{\nu-1,0} \circ h_{(v_{\nu,0})}\circ ...\circ h_{(v_{\nu,m_\nu})}$$ satisfies that  $\restr{h_{\nu,0}}{\mathcal{T}}=\restr{h_{\nu-1,0}}{\mathcal{T}}=\restr{\id}{\mathcal{T}}$, $ \restr{h_{\nu,0}}{\mathcal{T}'}=\restr{h_{\nu-1,0}}{\mathcal{T}'}=\restr{\alpha}{\mathcal{T}'}$ and $$\Sgn_{(i)}(h_{\nu,0},S_{Y_0}(w))=1\qq \forall w\in B_T(\gamma,\nu)\cap V_{t_0}.$$
		\begin{figure}[H]\caption{Step II of the proof of Lemma \ref{Lemme Omega 0 est non vide}} \label{drawing04}
			\includegraphics[scale=0.08]{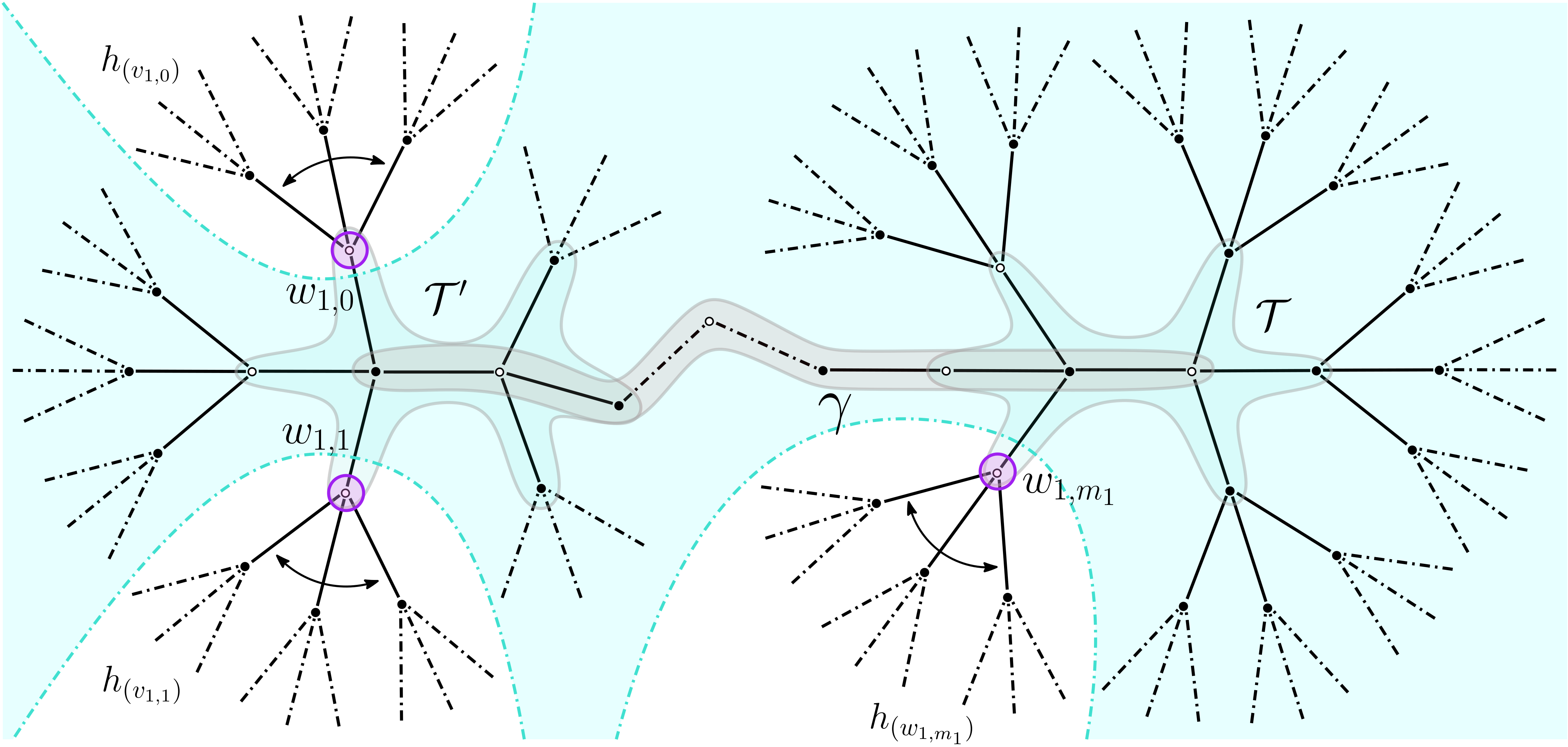}
		\end{figure}
		Consider the element $h_{2r,0}$ that we have just constructed. This element behaves as desired for the condition given by $Y_0$ on the vertices of $B_T(\gamma,2r)\cap V_{t_0}$. However, nothing ensures that the condition given by $Y_1$ on the vertices of $B_T(\gamma,0)\cap V_{t_1}$ is yet satisfied. We now take care of this task. Let $\{w_0,...,w_m\}$ be the set of vertices $w\in B_T( \gamma,0)\cap V_{t_1}$ satisfying that $\Sgn_{(i)}(h_{2r,0},S_{Y_1}(w))=-1$. For all $j=0,1,...,m$, let $v_j\in\bigcap_{w\in \gamma-\{w_j\}} T(w_j,w)$ such that $d(v_j,w_j)=\max(Y_1)$.
		In particular, notice that $v_{j}\in S_{Y_1}(w_{j})$ but that $v_{j}\not \in S_{Y_1}(w)$ for every $w\in B_T(\gamma,0) \cap V_{t_1}$. Furthermore, since $\max(Y_0)+2r\leq \max(Y_1)$ notice  from Remark \ref{remark les sommets ont des types opposers} that $v_j\not \in S_{Y_0}(v)$ for every $v\in B_T(\gamma,2r)\cap  V_{t_0}$. On the other hand, just as before, the automorphisms $h_{(v_{0})}, ..., h_{(v_{m})}$ fix $\mathcal{T}\cup \mathcal{T}'$ pointwise. In particular, $h_{1}= h_{2r,0}\circ h_{(v_{0})}\circ ...\circ h_{(v_{m})}$ satisfies:
		\begin{enumerate}[label=(\roman*)]
			\item $\restr{h_{1}}{\mathcal{T}}=\restr{h_{2r,0}}{\mathcal{T}}=\restr{\id}{\mathcal{T}}$ and $\restr{h_{1}}{\mathcal{T}'}=\restr{h_{2r,0}}{\mathcal{T}'}=\restr{\alpha}{\mathcal{T}'}$.
			\item $\Sgn_{(i)}(h_{1},S_{Y_0}(w))=1\qq \forall w\in B_T(\gamma,2r)\cap V_{t_0}.$
			\item $\Sgn_{(i)}(h_{1},S_{Y_1}(w))=1\qq \forall w\in B_T(\gamma,0)\cap V_{t_1}.$
		\end{enumerate} 
		\begin{figure}[H]\caption{Step III of the proof of Lemma \ref{Lemme Omega 0 est non vide}} \label{drawing05}
			\includegraphics[scale=0.08]{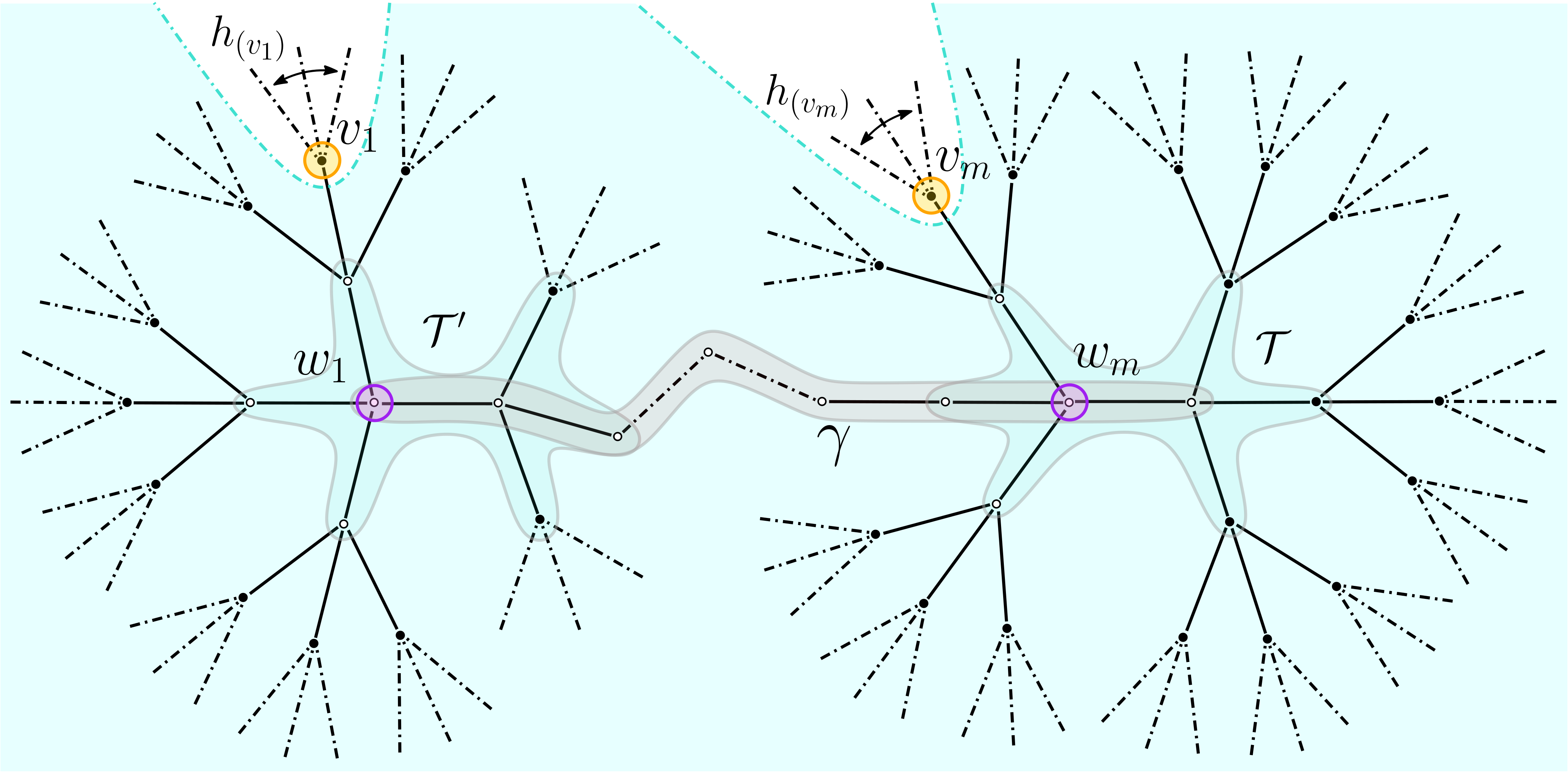}
		\end{figure}
		This proves that $h_1\in \Omega_0$ and therefore that $\Omega_0$ is not empty. 
	\end{proof}
	\begin{lemma}\label{Lemme Omega n est non vide}
		For all $n\geq 1$, the set	
		\begin{equation*}
		\Omega_n\qq=\qq \left\{ h \in \Omega_{n-1}\qq \ \middle\vert \begin{array}{l}
		\q\q \restr{h}{B_T(\gamma,n+{\max(Y_1)})}=\restr{h_{n}}{B_T(\gamma,n+{\max(Y_1)})},\\
		\qq\Sgn_{(i)}(h,S_{Y_0}(w))=1\q \forall w\in B_T(\gamma,n+2r)\cap V_{t_0},\\
		\qq\Sgn_{(i)}(h,S_{Y_1}(w))=1\q \forall w\in B_T(\gamma,n)\cap V_{t_1}
		\end{array}\right\}.
		\end{equation*}
		is a non-empty compact subsets of $\Aut(T)^+$. 
	\end{lemma}
	\begin{proof}
		We show that $\Omega_n$ is not empty by induction. Lemma \ref{Lemme Omega 0 est non vide} ensures that $\Omega_0$ is not empty. Suppose that $\Omega_{n-1}$ is not empty and let $h_n\in \Omega_{n-1}$ be the automorphism appearing in the definition of $\Omega_n$. Just as in the proof of Lemma \ref{Lemme Omega 0 est non vide}, we are going to modify $h_n$ with the automorphisms $h_{(v)}$ in order to obtain an element of $\Omega_{n}$. A concrete example of the procedure is given by figures \ref{drawing201} and \ref{drawing202} on a $4$-regular tree with $Y_0=\{0\}$ and even $\max{(Y_1)}$ (with the same conventions as before and where the vertices concerned by the current step are circled in pink). Let $\{\tilde{w}_0,...,\tilde{w}_{k}\}$ be the set of vertices $w\in B_T(\gamma,n+2r)\cap V_{t_0}$ such that $\Sgn_{(i)}(h_n,S_{Y_0}(w))=-1$. For each $j=0,1,...,k$ we choose a vertex $\tilde{v}_{j}\in \bigcap_{w\in B_T(\gamma,n+2r)\cap V_{t_0}-\{\tilde{w}_j\}}{}T(\tilde{w}_{j},w)$ such that $d(\tilde{v}_{j},\tilde{w}_{j})=\max(Y_0)$. Hence, notice that $\tilde{v}_j\not \in S_{Y_0}(w)$ for all $w\in B_T(\gamma,n+2r)\cap V_{t_0}-\{\tilde{w}_j\}$. Furthermore, since $\max(Y_1)-1\leq 2r+\max(Y_0)$ notice from Remark \ref{remark les sommets ont des types opposers} that $\tilde{v}_j\not \in S_{Y_1}(w)$ for all $w\in  B_T(\gamma,n-1)\cap V_{t_1}$. Furthermore, since $\Sgn_{(i)}(h_{n},S_{Y_0}(\tilde{w}_{j}))=-1$, $\restr{h_{n}}{\mathcal{T}}= \restr{\id}{\mathcal{T}}$, $\restr{h_{n}}{\mathcal{T}'}= \restr{\alpha}{\mathcal{T}'}$ and due to the form of $\mathcal{T}$ and $\mathcal{T}'$, notice that the automorphisms $h_{(\tilde{v}_{0})}, ..., h_{(\tilde{v}_{k})}$ fix $\mathcal{T}\cup \mathcal{T}'$ pointwise. In particular, $\tilde{h}_n=h_{n} \circ h_{(\tilde{v}_{0})}\circ ...\circ h_{(\tilde{v}_{k})}$ satisfies:
		\begin{enumerate}[label=(\roman*)]
			\item $\restr{\tilde{h}_n}{B_T(\gamma,n-1+\max{(Y_1)})}=\restr{h_{n}}{B_T(\gamma,n-1+\max(Y_1))}.$
			\item $\Sgn_{(i)}(h_{n},S_{Y_0}(w))=1\qq \forall w\in B_T(\gamma,n+2r)\cap V_{t_0}.$
			\item $\Sgn_{(i)}(h_{n},S_{Y_1}(w))=1\qq \forall w\in B_T(\gamma,n-1)\cap V_{t_1}.$
		\end{enumerate} 
		\begin{figure}[H]\caption{Step I of the proof of Lemma \ref{Lemme Omega n est non vide}}\label{drawing201}
			\includegraphics[scale=0.075]{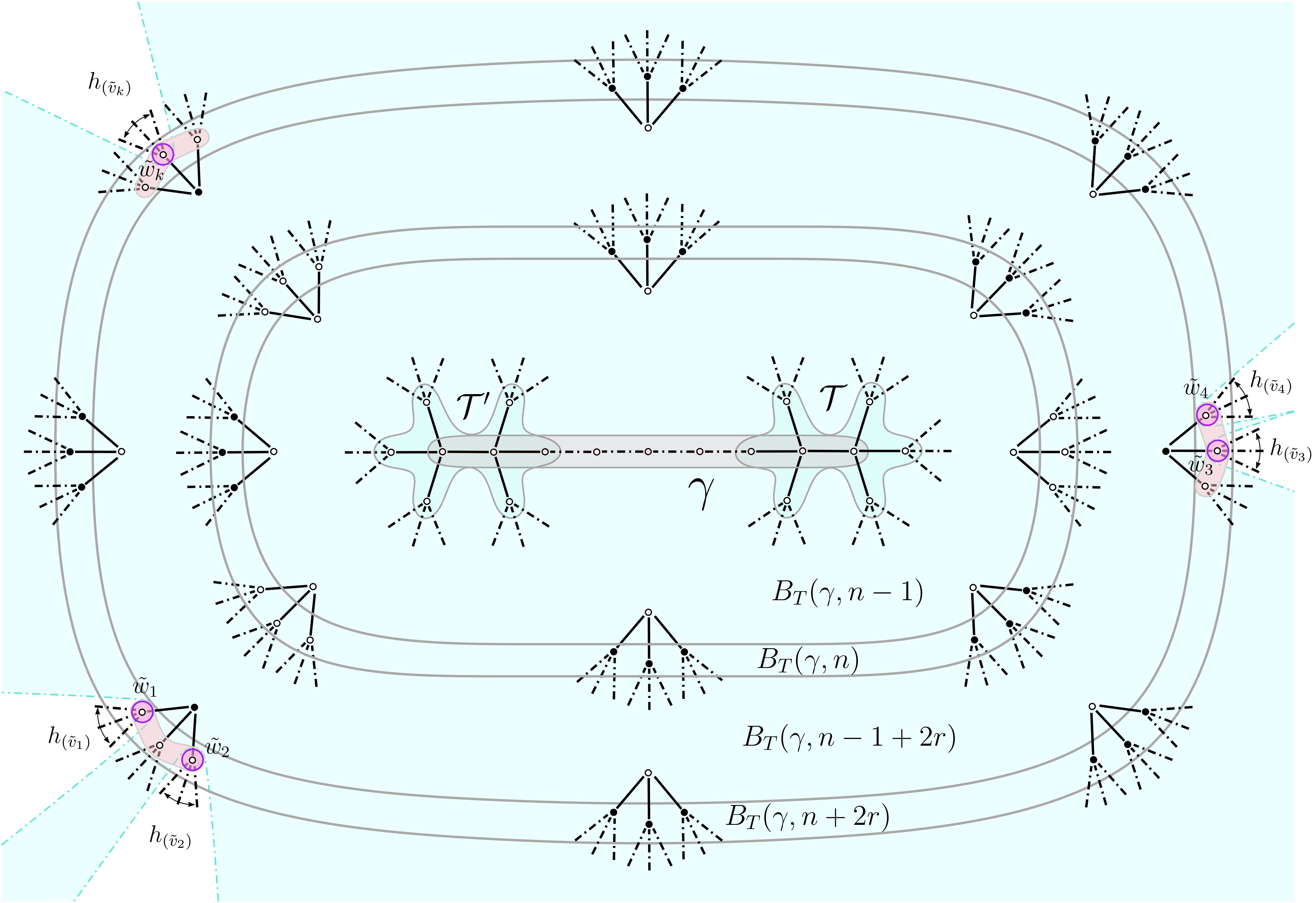}
		\end{figure}
		
		Now, let  $\{w_0,...,w_m\}$ be the set of vertices $w\in B_T(\gamma,n)\cap V_{t_1}$ such that $\Sgn_{(i)}(\tilde{h}_n,S_{Y_1}(w))=-1$. For each $j=0,1,...,m$, we choose a vertex $v_{j}\in \bigcap_{w\in B_T(\gamma,n)-\{w_j\}}{}T(w_{j},w)$ such that $d(v_{j},w_{j})=\max(Y_1)$. Since $\max(Y_0)+2r\leq\max(Y_1)$ notice from Remark \ref{remark les sommets ont des types opposers} that $v_j\not \in S_{Y_0}(w)$ for every $w\in B_T(\gamma,n+2r)\cap V_{t_0}$ and $v_j\not \in S_{Y_1}(w)$ for every $w\in B_T(\gamma,n)\cap V_{t_1}-\{w_j\}$. Just as before, notice that the automorphisms $h_{(v_{0})}, ...,h_{(v_{m})}$ fix $\mathcal{T}\cup \mathcal{T}'$ pointwise. In particular, $h_{n+1}=\tilde{h}_{n} \circ h_{(v_{0})}\circ ...\circ h_{(v_{m})}$ satisfies:
		\begin{enumerate}[label=(\roman*)]
			\item $\restr{h_{n+1}}{B_T(\gamma,n-1+{\max(Y_1)})}=\restr{h_{n}}{B_T(\gamma,n-1+{\max(Y_1)})}.$
			\item $\Sgn_{(i)}(h_{n+1},S_{Y_0}(w))=1\qq \forall w\in B(\gamma,n+2r)\cap V_{t_0}.$
			\item $\Sgn_{(i)}(h_{n+1},S_{Y_1}(w))=1\qq \forall w\in B(\gamma,n)\cap V_{t_1}.$
		\end{enumerate}
		\begin{figure}[H]\caption{Step II of the proof of Lemma \ref{Lemme Omega n est non vide}}\label{drawing202}
			\includegraphics[scale=0.075]{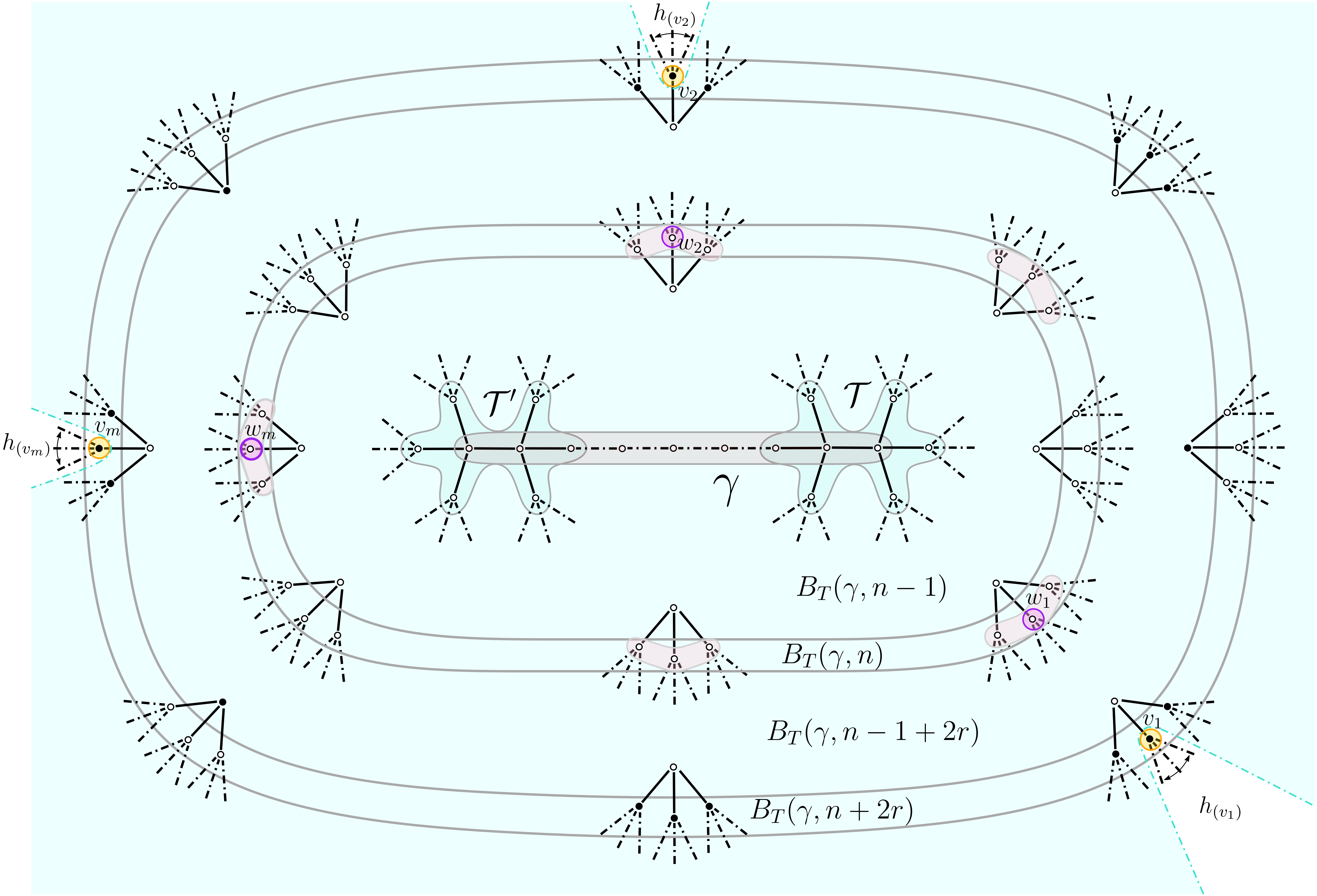}
		\end{figure}
		This proves that $\Omega_n$ is not empty. We now show that $\Omega_n$ is compact for every integer $n\geq 1$. To this end, notice that $\Omega_n$ is a closed subset of $\Aut(T)^+$ and that
		\begin{equation*}
		\Omega_n\qq \subseteq \qq h_{n}\qq \Fix_{\Aut(T)^+}(B_T(\gamma,n+\max(Y_1))).
		\end{equation*}
		Since the right hand side is a compact subset of $\Aut(T)^+$ the results follows.
	\end{proof}
	
	\noindent We are finally able to prove the result announced at the beginning of the section.
	\begin{theorem}\label{Theorem S 0 est fertile }
		The generic filtration $\mathcal{S}_{0}$ of $G^+_{(i)}(Y_0,Y_1)$ factorises$^+$ at all depths $l\geq 1$.
	\end{theorem}
	\begin{proof}
		For brevity, we let $G=G^+_{(i)}(Y_0,Y_1)$. To prove that $\mathcal{S}_0$ factorises$^+$ a depth $l\geq 1$, we shall successively verify the three conditions of the Definition \ref{definition olsh facto}.
		
		First, we need to prove that for every  $U$ in the conjugacy class of an element of $\mathcal{S}_0\lb l \rb$ and every subgroup $V$ in the conjugacy class of an element of $\mathcal{S}_0$ with $V \not\subseteq U$, there exists a $W$ in the conjugacy class of an element of $\mathcal{S}_0\lb l -1\rb$ such that $$U\subseteq W \subseteq V U.$$ Let $U$ and $V$ be as above. If $V$ is conjugate to an element of $\mathcal{S}_0\lb l'\rb$ for some $l'\geq l$, the result follows directly from Proposition \ref{Proposition la premiere etape de la fcatorization pour GY1Y_2}. Therefore, let us suppose that $l'<l$. By the definition of $\mathcal{S}_0$ and since $\mathfrak{T}_0$ is stable under the action of $G$, there exist two subtrees $\mathcal{T},\mathcal{T}'\in \mathfrak{T}_0$ such that $U=\Fix_G(\mathcal{T})$ and $V=\Fix_G(\mathcal{T}')$. There is two cases. Either, $\mathcal{T}'\subseteq \mathcal{T}$ and there exists a subtree $\mathcal{P}\in \mathfrak{T}_{0}$ such that $\mathcal{T}'\subseteq \mathcal{P}\subseteq \mathcal{T}$ and $\Fix_G(\mathcal{R})\in \mathcal{S}_0\lb l-1\rb$. In that case $$\Fix_G(\mathcal{T})\subseteq \Fix_G(\mathcal{P})\subseteq \Fix_G(\mathcal{T}')\subseteq \Fix_G(\mathcal{T}')\Fix_G(\mathcal{T})$$ and the result follows. Or else, $\mathcal{T}'\not \subseteq \mathcal{T}$ and since $l'< l$, there exists a subtree $\mathcal{P}\in \mathfrak{T}_0$ such that $\mathcal{T}'\subseteq \mathcal{P}\not= \mathcal{T}$ and $\Fix_G(\mathcal{P})\in \mathcal{S}_0\lb l\rb$. In particular, Proposition \ref{Proposition la premiere etape de la fcatorization pour GY1Y_2} ensures the existence of a subgroup $W\in \mathcal{S}_0\lb l-1\rb$ such that $U\subseteq W \subseteq \Fix_G(\mathcal{P})U\subseteq \Fix_G(\mathcal{T}')U$. This proves the first condition.
		
		Next, we need to prove that $N_{G}(U, V)= \{g\in {G} : g^{-1}Vg\subseteq U\}$ is compact for every $V$ in the conjugacy class of an element of $\mathcal{S}_0$. Just as before, notice that there exists a $\mathcal{T}'\in \mathfrak{T}_0$ such that $V = \Fix_G(\mathcal{T}')$. Since $G$ satisfies the hypothesis \ref{Hypothese H0}, notice that
		\begin{equation*}
		\begin{split}
		N_{G}(U, V)&= \{g\in G : g^{-1}Vg\subseteq U\}=\{g\in G : g^{-1}\Fix_G(\mathcal{T}')g\subseteq \Fix_G(\mathcal{T})\}\\ 
		&=\{g\in G : \Fix_G(g^{-1}\mathcal{T}')\subseteq \Fix_G(\mathcal{T})\}= \{g\in G : g\mathcal{T}\subseteq \mathcal{T}'\}.
		\end{split}
		\end{equation*}
		In particular, since both $\mathcal{T}$ and $\mathcal{T}'$ are finite subtrees of $T$, $N_G(U,V)$ is a compact subset of $G$ which proves the second condition.
		
		Finally, we need to prove for every $W$ in the conjugacy class of an element of $\mathcal{S}_0\lb l-1\rb$ with $U\subseteq W$ that 
		\begin{equation*}
		W\subseteq N_G(U,U) =\{g\in G : g^{-1}Ug\subseteq U\}.
		\end{equation*} 
		For the same reasons as before, there exists $\mathcal{R}\in \mathfrak{T}_0$ such that $W= \Fix_G(\mathcal{R})$. Furthermore, since $U\subseteq W$ and since $G$ satisfies the hypothesis \ref{Hypothese H0}, notice that $\mathcal{R}\subseteq \mathcal{T}$. Moreover, since $\Fix_G(\mathcal{R})$ has depth $l-1$, notice that $\mathcal{R}$ contains every interior vertex of $\mathcal{T}$. Since $G$ is unimodular and satisfies the hypothesis \ref{Hypothese H0}, this implies that
		
		\begin{equation*}
		\begin{split}
		\Fix_G(\mathcal{R})&\subseteq \{h\in G : h\mathcal{T}\subseteq \mathcal{T}\}= \{h\in G : \Fix_G(\mathcal{T})\subseteq\Fix_G(h\mathcal{T})\}\\
		&=\{h\in G : h^{-1}\Fix_G(\mathcal{T})h\subseteq\Fix_G(\mathcal{T})\}= N_G(U,U)
		\end{split}
		\end{equation*}
		which proves the third condition.
	\end{proof}

	\subsection{Description of cuspidal representations}\label{section cuspidal rep description Radu}
	
	The purpose of this section is to provide the explicit description of the cuspidal representations of $G_{(i)}^+(Y_0,Y_1)$ following from the factorisation$^+$ of the generic filtration $\mathcal{S}_0$ see Theorem \ref{la classification des reresentations cuspidale Radu} below. Let $T$ be a $(d_0,d_1)$-semi-regular tree with $d_0,d_1\geq 4$, $V=V_0\sqcup V_1$ be the associated bipartition, $i$ be a legal coloring of $T$, consider two finite subsets $Y_0,Y_1\subseteq \N$, let  $G=G_{(i)}^+(Y_0,Y_1)$ and 
	$$\mathfrak{T}_{0}=\{B_T(v,r): v\in V, r\geq 1\}\sqcup\{B_T(e,r): e\in E, r\geq 0\}.$$
	We have shown in Sections \ref{chapter generic filtration for Gi+YY} and \ref{Section factorization} that
	\begin{equation*}
	\mathcal{S}_{0}=\{\Fix_G(\mathcal{T}): \mathcal{T}\in \mathfrak{T}_{0}\}.
	\end{equation*}
	is a generic filtration of $G_{(i)}^+(Y_0,Y_1)$ that factorises$^+$ at all depths $l\geq 1$.  We adopt here the formalism developed in Section \ref{Section the bijective correspondence of theorem A}. Theorem \ref{la version paki du theorem de classification} provides a bijective correspondence between irreducible representations of $G_{(i)}^+(Y_0,Y_1)$ at depth $l$ with seed $C \in \mathcal{F}_{\mathcal{S}_0}$ and the $\mathcal{S}_0$-standard representations of $\Aut_G(C)$. We start by identifying those seeds and show that the cuspidal representations of $G$ are precisely the irreducible representations of $G$ at depth $l\geq 1$ with respect to $\mathcal{S}_0$. In light of Lemma \ref{Lemme la forme des elements de Sl} we consider the partition $\mathfrak{T}_0=\bigsqcup_{l\in \N}\mathfrak{T}_0\lb l\rb$ where:
	\begin{itemize}\label{de la merde one definit mathfrak T 0}
		\item $\mathfrak{T}_0\lb l\rb = \{B_T(e,\frac{l}{2}): e\in E\}$ if $l$ is even.
		\item $\mathfrak{T}_0\lb l \rb=\{ B_T\big(v,(\frac{l+1}{2})\big): v\in V\}$ if $l$ is odd.
	\end{itemize} 
	Notice that for all $l\in \N$, $\mathfrak{T}_0\lb l\rb$ is stable under the action of $G$. Furthermore, if $l$ is even or if $l$ is odd and $G$ is vertex-transitive the set $\mathfrak{T}_0\lb l\rb$ consists of a single $G$-orbit. On the other hand, if $l$ is odd and $G$ has two orbits of vertices the set $\mathfrak{T}_0\lb l\rb$ consists of the two $G$-orbits $\{B_T\big(v,(\frac{l+1}{2})\big): v\in V_0\}$ and $\{B_T\big(v,(\frac{l+1}{2})\big): v\in V_1\}$. In particular, in light of Lemma \ref{lemma G satisfies hypothesis H0}, there is either one or two elements of $\mathcal{F}_{\mathcal{S}_0}=\{\mathcal{C}(U): U\in \mathcal{S}_0\}$ at height $l$ and each such element is of the form
	$$C=\{\Fix_G(\mathcal{T}): \mathcal{T}\in \mathcal{O}\}$$
	for some $G$-orbit $\mathcal{O}$ of an element $\mathcal{T} \in\mathfrak{T}_0\lb l \rb$. It follows that $U\in \mathcal{S}_0\lb 0\rb$ if and only if $U\in \{\Fix_G(e):e\in E\}$. In particular, the cuspidal representations of $G$ are exactly the irreducible representations of $G$ at depth $l\geq 1$ with respect to $\mathcal{S}_0$.    Now, let $\pi$ be a cuspidal representation of $G$ with seed $C_\pi\in \mathcal{F}_{\mathcal{S}_0}$, let $U\in C_\pi$ and $\mathcal{T}\in \mathfrak{T}_0$ be such that $U=\Fix_G(\mathcal{T})$. Since $\mathcal{S}_0$ factorises$^+$ at all depths $l\geq 1$, Theorem \ref{la version paki du theorem de classification} ensures that $\pi$ is induced from an irreducible representation of $N_G(U)$ passing to the quotient $\Aut_G(C_\pi)\simeq N_G(U)/U$. Furthermore, since $G$ satisfies the hypothesis \ref{Hypothese H0}, notice that 
	\begin{equation*}
	\begin{split}
	N_G(U)= \{g\in G: gUg^{-1}=U\}
	&=\{g\in G: g\Fix_G(\mathcal{T})g^{-1}=\Fix_G(\mathcal{T})\}\\
	&=\{g\in G: \Fix_G(g\mathcal{T})=\Fix_G(\mathcal{T})\}\\
	&=\{g\in G: g\mathcal{T}=\mathcal{T}\}=\Stab_G(\mathcal{T})
	\end{split}
	\end{equation*}
	In particular, $\Aut_G(C_\pi)$ can be identified with the group of automorphisms of $\mathcal{T}$ obtained from the restriction of the action of $\Stab_G(\mathcal{T})$. Moreover, since $G$ satisfies the hypothesis \ref{Hypothese H0} notice that
	\begin{equation*}
	\tilde{\mathfrak{H}}_{\mathcal{S}_0}(U)=\{\Fix_G(\mathcal{R}): \mathcal{R}\in \mathfrak{T}_0, \qq \mathcal{R}\subsetneq \mathcal{T} \mbox{ and }\mathcal{R}\mbox{ is maximal for this property}\}
	\end{equation*}
	and
	$$\mathfrak{H}_{\mathcal{S}_0}(C_\pi)=\{ p_U(W): W\in \tilde{\mathfrak{H}}_{\mathcal{S}_0}(U)\}$$ 
	is the set of pointwise stabilisers (in $\Aut_G(C_\pi)$) of subtrees $\mathcal{R}\in \mathfrak{T}_0$ satisfying $\mathcal{R}\subsetneq \mathcal{T}$ and which are maximal for this property. In particular, the $\mathcal{S}_0$-\tg{standard} representations of $\Aut_G(C_\pi)$ are the irreducible representations which do not admit any non-zero invariant vector for the pointwise stabiliser of any subtree $\mathcal{R}$ of $\mathcal{T}$ which belongs to $\mathfrak{T}_0$ and is maximal for this property. This leads to the following description of the cuspidal representations of $G_{(i)}^+(Y_0,Y_1)$.
	
	\begin{theorem}\label{la classification des reresentations cuspidale Radu}
		Let $T$ be a $(d_0,d_1)$-semi-regular tree with $d_0,d_1\geq 4$, $i$ be a legal coloring of $T$, $Y_0,Y_1\subseteq \N$ be two finite subsets, $G=G_{(i)}^+(Y_0,Y_1)$. 
		Then, the cuspidal representations of $G$ are exactly the irreducible representations at depth $l\geq 1$ with respect to $\mathcal{S}_{0}$. Furthermore, if $\pi$ is a cuspidal representation at depth $l$ we have that:
		\begin{itemize}[leftmargin=*]
			\item $\Hr{\pi}^{\Fix_G(\mathcal{R})}=\{0\}$ for any $\mathcal{R}\in \bigsqcup_{r<l}\mathfrak{T}_0\lb r\rb$.
			\item There exists a unique $G$-orbit $\mathcal{O}$ of $\mathfrak{T}_0\lb l \rb$  such that $\pi$ admits a non-zero $\Fix_G(\mathcal{T})$-invariant vector for any $\mathcal{T}\in \mathcal{O}$. Furthermore, $\mathcal{O}$ is the only orbit in $\mathfrak{T}_0$ under the action of $G$ such that $C_\pi=\{\Fix_G(\mathcal{T}): \mathcal{T}\in \mathcal{O}\}$.
			\item If $\mathcal{O}$ is the unique $G$-orbit of $\mathfrak{T}_0\lb l \rb$ corresponding to $\pi$ and if $\mathcal{T}\in \mathcal{O}$, there exists a non zero function of positive type associated to $\pi$ supported in $\Stab_G(\mathcal{T})$. In particular, $\pi$ is square-integrable its equivalence class is isolated in the unitary dual $\widehat{G}$ for the Fell topology.	
		\end{itemize}
		Furthermore, for every $C\in \mathcal{F}_{\mathcal{S}_0}$ at height $l\geq 1$, there exists a bijective correspondence between the equivalence classes of irreducible representations of $G$ with seed $C$ and the equivalence classes of $\mathcal{S}_0$-standard representations of $\Aut_G(C)$. More precisely, denote by $\mathcal{O}$ the $G$-orbit in $\mathfrak{T}_0\lb l \rb$ such that $C=\{\Fix_G(\mathcal{T}): \mathcal{T}\in \mathcal{O}\}$. For every $\mathcal{T}\in \mathcal{O}$ the following holds:
		\begin{enumerate}[leftmargin=*, label=(\roman*)]
			\item If $\pi$ is a cuspidal representation of $G$ with seed $C$, $(\omega_\pi, \Hr{\pi}^{\Fix_G(\mathcal{T})})$ is a $\mathcal{S}_0$-standard representation of $\Aut_{G}(C)$ such that
			\begin{equation*}
			\pi \simeq T(\Fix_G(\mathcal{T}),\omega_\pi)=\Ind_{\Stab_G(\mathcal{T})}^{G}(\omega_\pi \circ p_{\Fix_G(\mathcal{T})}).
			\end{equation*}
			\item If $\omega$ is a $\mathcal{S}_0$-standard representation of $\Aut_{G}(C)$, the representation $T(\Fix_G(\mathcal{T}),\omega)$ is a cuspidal of $G$ with seed $C$. 
		\end{enumerate}
		Furthermore, if $\omega_1$ and $\omega_2$ are $\mathcal{S}_0$-standard representations of $\Aut_{G}(C)$, we have that $T(\Fix_G(\mathcal{T}),\omega_1)\simeq T(\Fix_G(\mathcal{T}),\omega_2)$ if and only if $\omega_1\simeq\omega_2$. In particular, the above two constructions are inverse of one another.
	\end{theorem}
	
	\subsection{Existence of cuspidal representations}\label{section existence des cuspidale}
	
	Let $T$ be a $(d_0,d_1)$-semi-regular tree with $d_0,d_1\geq 4$, $V=V_0\sqcup V_1$ be the associated bipartition, $i$ be a legal coloring of $T$, consider two finite subsets $Y_0,Y_1\subseteq \N$ and let $G=G_{(i)}^+(Y_0,Y_1)$. We adopt the same formalism as in Section \ref{section cuspidal rep description Radu}. Theorem \ref{la classification des reresentations cuspidale Radu} provides a description of the cuspidal representations of $G$. However, this result does not ensure the existence of those cuspidal representations. The purpose of this section is to prove the existence of a cuspidal representation with seed $C$ for each conjugacy class $C\in \mathcal{F}_{\mathcal{S}_0}$ at height $l\geq 1$. Equivalently we prove the following result.
	\begin{theorem}
		Let $G=G^+_{(i)}(Y_0,Y_1)$ and $C\in \mathcal{F}_{\mathcal{S}_0}$ be a conjugacy class at height $l\geq 1$. Then, there exists a ${\mathcal{S}_0}$-standard representations of $\Aut_G(C)$. 
	\end{theorem}
	The proof of this theorem is gathered in the following results. The first proposition treat the case of even height.
	\begin{proposition}\label{Proposition existence of standard for edges }
		For each $C\in \mathcal{F}_{\mathcal{S}_0}$ at even height $l\geq 1$, there exists a ${\mathcal{S}_0}$-standard representation of $\Aut_{G}(C)$.
	\end{proposition} 
	\begin{proof}
		Since $l$ is even, Lemma \ref{Lemme la forme des elements de Sl} ensures the existence of an edge $e\in E$ and an integer $r\geq 1$ such that $B_T(e,r)\in \mathfrak{T}_0$ and $C=\mathcal{C}(\Fix_{G}(B_T(e,r))$. For brevity, we let $\mathcal{T}=B_T(e,r)$. As observed in Section \ref{section cuspidal rep description Radu} we have that $\Aut_G(C)\simeq \Stab_G(\mathcal{T})/\Fix_G(\mathcal{T})$ and
		\begin{equation*}
		\begin{split}
		\mathfrak{H}_{\mathcal{S}_0}(\Fix_G(\mathcal{T}))&=\big\{\Fix_G(\mathcal{R})/\Fix_G(\mathcal{T}): \mathcal{R}\in \mathfrak{T}_0, \qq \mathcal{R}\subsetneq \mathcal{T} \\
		&\mbox{ }\mbox{ }\mbox{ }\mbox{ }\mbox{ }\mbox{ }\mbox{ }\mbox{ }\mbox{ }\mbox{ }\mbox{ }\mbox{ }\mbox{ }\mbox{ }\mbox{ }\mbox{ }\mbox{ }\mbox{ }\mbox{ }\mbox{ }\mbox{ }\mbox{ }\mbox{ }\mbox{ }\mbox{ }\mbox{ }\mbox{ }\mbox{ }\mbox{ }\mbox{ }\mbox{ }\mbox{ }\mbox{ }\mbox{ }\mbox{ and } \mathcal{R}\mbox{ is maximal for this property}\big\}\\
		&=\{\Fix_G(B_T(v,r))/\Fix_G(B_T(e,r)): v\in e\}.
		\end{split}
		\end{equation*} 
		Let $v_0,v_1$ be the two vertices of $e$, set $\mathcal{T}_i=B_T(v_i,r)$ and notice that $\mathcal{T}_0\cup \mathcal{T}_1=\mathcal{T}$. Moreover, notice that $\Stab_G(\mathcal{T})=\{g\in G : ge=e\}$ acts by permutations on the set $\{\mathcal{T}_0, \mathcal{T}_1\}$. Furthermore, since $G$ satisfies the hypothesis \ref{Hypothese H0} notice that $\Fix_G(\mathcal{T})\subsetneq \Fix_G(\mathcal{T}_i)\subsetneq \Stab_G(\mathcal{T})$. In particular, the hypothesis of Proposition \ref{existence criterion} are satisfied and the result follows.
	\end{proof}
	A similar reasoning leads to a proof of the existence of ${\mathcal{S}_0}$-standard representation of $\Aut_{G}(C)$ for all $C\in \mathcal{F}_{\mathcal{S}_0}$ with odd height $l> 1$.
	\begin{proposition}\label{existence of standard for odd depth Radu}
		For each $C\in \mathcal{F}_{\mathcal{S}_0}$ at odd height $l> 1$, there exists a ${\mathcal{S}_0}$-standard representation of $\Aut_{G}(C)$.
	\end{proposition}
	\begin{proof}
		Since $l$ is odd, Lemma \ref{Lemme la forme des elements de Sl} ensures the existence of a vertex $v\in V$ and an integer $r\geq 1$ such that $B_T(v,r+1)\in \mathfrak{T}_0$ and $$C=\mathcal{C}(\Fix_{G}(B_T(v,r+1)).$$ For brevity, we let $\mathcal{T}=B_T(v,r+1)$ of $T$. As observed in Section \ref{section cuspidal rep description Radu} we have that $\Aut_G(C)\simeq \Stab_G(\mathcal{T})/\Fix_G(\mathcal{T})$ and
		\begin{equation*}
		\begin{split}
		\mathfrak{H}_{\mathcal{S}_0}(\Fix_G(\mathcal{T}))&=\big\{\Fix_G(\mathcal{R})/\Fix_G(\mathcal{T}): \mathcal{R}\in \mathfrak{T}_0, \qq \mathcal{R}\subsetneq \mathcal{T} \\
		&\mbox{ }\mbox{ }\mbox{ }\mbox{ }\mbox{ }\mbox{ }\mbox{ }\mbox{ }\mbox{ }\mbox{ }\mbox{ }\mbox{ }\mbox{ }\mbox{ }\mbox{ }\mbox{ }\mbox{ }\mbox{ }\mbox{ }\mbox{ }\mbox{ }\mbox{ }\mbox{ }\mbox{ }\mbox{ }\mbox{ }\mbox{ }\mbox{ }\mbox{ }\mbox{ }\mbox{ }\mbox{ }\mbox{ }\mbox{ }\mbox{ and } \mathcal{R}\mbox{ is maximal for this property}\big\}\\
		&=\{\Fix_G(B_T(e,r-1))/\Fix_G(B_T(v,r)): e\in E(B_T(v,1))\}.
		\end{split}
		\end{equation*} 
		Let $\{w_1,...,w_d\}$ be the leaves of $B_T(v,1)$. For every $i=1,...,d$ let $\mathcal{T}_i= (B_T(v,1)-\{w_i\})^{(r-1)}$ where $\mathcal{R}^{(t)}=\{v\in V: \exists w\in \mathcal{R} \mbox{ s.t. }d_T(v,w)\leq t\}$. Notice that $\mathcal{T}_i\cup \mathcal{T}_j=\mathcal{T}$. On the other hand, in our case $$\Stab_G(\mathcal{T})=\{g\in G : gv=v\}=\Fix_G(v)$$ and $\Fix_G(v)$ acts by permutation on the set $\{\mathcal{T}_1,...,\mathcal{T}_d\}$. Finally, notice that $v\in \mathcal{T}_i$ and therefore that $\Fix_G(\mathcal{T}_i)\subseteq \Stab_G(\mathcal{T})$ for all $i=1,...,d$. Furthermore, since $G$ contains $G^{+}_{(i)}(\{0\},\{0\})$, notice that $\Fix_G(\mathcal{T})\subsetneq \Fix_G(\mathcal{T}_i)\subsetneq \Stab_G(\mathcal{T})$. In particular, Proposition \ref{existence criterion} ensures the existence of an irreducible representation of $\Aut_{G}(C)$ without non-zero $\Fix_G(\mathcal{T}_i)/\Fix_G(\mathcal{T})$-invariant vectors. The result follows from the fact that for every edge $e\in E(B_T(v,1))$ there exists some $i\in \{1,...,d\}$ such that $B_T(e,r)\subseteq \mathcal{T}_i$ and hence $\Fix_G(B_T(e,r))/\Fix_G(\mathcal{T})\subseteq \Fix_G(\mathcal{T}_i)/\Fix_G(\mathcal{T})$.
	\end{proof}
	The next lemma treats the remaining case $l=1$ where Proposition \ref{existence criterion} does not apply.
	\begin{lemma}\label{existence of standard for vertices Radu}
		For each $C\in \mathcal{F}_{\mathcal{S}_0}$ at height $1$, there exists a ${\mathcal{S}_0}$-standard representation of $\Aut_{G}(C)$. 
	\end{lemma}
	\begin{proof}
		Lemma \ref{Lemme la forme des elements de Sl} ensures the existence of a vertex $v\in V$ such that $C=\mathcal{C}(\Fix_{G}(B_T(v,1))$. For brevity, we let $\mathcal{T}=B_T(v,r)$. We recall as observed in Section \ref{section cuspidal rep description Radu} that $$\Aut_G(C)\simeq \Stab_G(\mathcal{T})/\Fix_G(\mathcal{T})$$ where $\Stab_G(\mathcal{T})=\{g\in G : g\mathcal{T}\subseteq \mathcal{T}\}= \{ g\in G : gv=v\}=\Fix_G(v)$. In particular, $\Aut_G(C)$ can be realised as the group of automorphisms of $B_T(v,1)$ obtained by restricting the action of $\Fix_G(v)$. Furthermore, we have that  
		\begin{equation*}
		\begin{split}
		\mathfrak{H}_{\mathcal{S}_0}(\Fix_G(\mathcal{T}))&=\big\{\Fix_G(\mathcal{R})/\Fix_G(\mathcal{T}): \mathcal{R}\in \mathfrak{T}_0, \qq \mathcal{R}\subsetneq \mathcal{T} \\
		&\mbox{ }\mbox{ }\mbox{ }\mbox{ }\mbox{ }\mbox{ }\mbox{ }\mbox{ }\mbox{ }\mbox{ }\mbox{ }\mbox{ }\mbox{ }\mbox{ }\mbox{ }\mbox{ }\mbox{ }\mbox{ }\mbox{ }\mbox{ }\mbox{ }\mbox{ }\mbox{ }\mbox{ }\mbox{ }\mbox{ }\mbox{ }\mbox{ }\mbox{ }\mbox{ }\mbox{ }\mbox{ }\mbox{ }\mbox{ }\mbox{ and } \mathcal{R}\mbox{ is maximal for this property}\big\}\\
		&=\{\Fix_G(f)/\Fix_G(B_T(v,1)): f\in E(B_T(v,1))\}.
		\end{split}
		\end{equation*} 
		Let $d$ be the degree of $v$ in $T$, let $X=E(B_T(v,1))$ and let $e\in X$. Since $G$ is closed subgroup of $\Aut(T)$ acting $2$-transitively on the boundary $\partial T$, Theorem \ref{Theorem Burger mozes acting on the boundary trans} ensures that $\underline{G}(v)$ is $2$-transitive. In particular, $\Aut_{G}(C)$ is $2$-transitive $X$ and  Lemma \ref{les rep de Qsur Qi} ensures the existence of an irreducible representation $\sigma$ of $\Aut_{G}(C)$ without non-zero $\Fix_{\Aut_{G}(C)}(e)$-invariant vector. Since $\Fix_G(v)$ is transitive on $E(B_T(v,1))$, this representation does not admit a non-zero $\Fix_{\Aut_{G}(C)}(f)$-invariant vector for any $f\in E(B_T(v,1))$. The lemma follows since $\Fix_{\Aut_{G}(C)}(f)=\Fix_G(f)/\Fix_G(B_T(v,1))$.
	\end{proof}
	\section{The Radu groups are uniformly admissible}\label{section radu groups are Type I}
	
	Let $T$ be a $(d_0,d_1)$-semi-regular tree with $d_0,d_1\geq 4$, $i$ be a legal coloring of $T$, consider two finite subsets $Y_0,Y_1\subseteq \N$ and let $G=G^+_{(i)}(Y_0,Y_1)$. The purpose of this section is to exploit the classification of the irreducible representations of $G$ to shows that this group is uniformly admissible and hence CCR. We recall that a totally disconnected locally compact group $G$ is uniformly admissible if for every compact open subgroup $K$, there exists an integer $k_K$ such that $\dim(\Hr{\pi}^K)<k_K$ for every irreducible representation $\pi$ of $G$. 
	
	The following result ensures the existence of a uniform bound for the spherical representations of any closed subgroup $G\leq \Aut(T)$ acting $2$-transitively on the boundary. 
	
	\begin{proposition}\label{Proposition les representations spheriques sont admissibles}
		Let $T$ be a thick semi-regular tree, $v\in V$ and $G\leq \Aut(T)$ be a closed subgroup acting $2$-transitively on the boundary. Then, for all $n\in \N_0$ there exists a positive integer $k_n\in \N$ such that ${\rm dim}(\Hr{\pi}^{\Fix_G(B_T(v,n))})\leq k_n$ for every spherical representation $\pi$ of $G$ admitting a non-zero $\Fix_G(v)$-invariant vector.
	\end{proposition}
	\begin{proof}
		
		For brevity, let $K=\Fix_G(v)$ and $K_n=\Fix_G(B_T(v,n))$ for all $ n\in \N$. Since $(G, \Fix_G(v))$ is a Gelfand pair by Lemma \ref{Lemma aut T fixator of a vertex is a Gelfand pair}, our hypothesis on $\pi$ and Theorem \ref{Theorem dimension plus petite que 1} ensure that ${\rm dim}(\Hr{\pi}^K)= 1.$ Now, let $\xi$ be a unit vector of $\Hr{\pi}^K$, $\eta$ be a unit vector of $\Hr{\pi}^{K_n}$ and $\mu$ be the Haar measure of $G$ renormalised in such a way that $\mu(K)=1$. Notice that 
		$$\varphi_{\xi, \eta}\fct{G}{\C}{g}{\prods{\pi(g)\xi}{\eta}}_{\Hr{\pi}}$$
		is $K$-right invariant and $K_n$-left invariant continuous function. On the other hand, since ${\rm dim}(\Hr{\pi}^K)= 1$, notice for all $g,h\in G$ that 
		\begin{equation*}
		\begin{split}
		\int_K\varphi_{\xi,\eta}(gkh)\qq\diff\mu(k)&=\prods{\int_K \pi(gkh)\xi\diff \mu( k)}{\eta}\\
		&=\prods{\pi(h)\xi}{\int_K\pi(k^{-1})\pi(g^{-1})\eta}\\
		&=\prods{\pi(h)\xi}{\alpha(\eta, g)\xi}=\overline{\alpha(\eta, g)}\varphi_{\xi,\xi}(h)
		\end{split}
		\end{equation*}
		for some $\alpha(\eta, g)\in \C$. However, $\varphi_{\xi,\xi}(1_G)=1$ and hence $\overline{\alpha(\eta, g)}=\varphi_{\xi,\eta}(g)$. This implies for all $g,h\in G$ that
		\begin{equation}\label{equation de semi sphericite}
		\int_K\varphi_{\xi,\eta}(gkh)\qq \diff \mu(k)=\varphi_{\xi,\eta}(g)\varphi_{\xi,\xi }(h).
		\end{equation}
		Since $\varphi_{\xi,\eta}$ is $K$-right invariant and $K_n$-left invariant notice that it can be realised as a function $\phi:Gv\rightarrow \C$ on the orbit $Gv$  and constant on the $K_n$-orbits of $v$. On the other hand, since $K_n$ is an open subgroup of the compact group $K$, the index of $K_n$ in $K$ is finite. As $K$ is transitive on the boundary of the tree, this implies that $K_n$ has finitely many orbit on $\partial T$. In particular, there exists an integer $N_n\geq \max\{2,n\}$ such that $\partial T(w,v)$ is contained in a single $K_n$-orbit for all $w\in \partial B_T(v,N_n)$ where $\partial T(w,v)$ is the set of ends of $T(w,v)=\{u\in V: d_T(u,w)<d_T(u,v)\}$ which are not vertices.
		Now, let $t\in G$ be such that $d_T(v,t v)=2$. Taking $k=t$ in \eqref{equation de semi sphericite} we obtain that the sum of values of $\phi$ on the vertices at distance $2$ from $gv$ is equal to $\phi(gv)\varphi_{\xi,\xi}(t)$. In particular, $\phi$ is entirely determined by the values it takes in $B_T(v,N_n)$. Hence, the space $\mathcal{L}_n$ of $K$-right invariant, $K_n$-left invariant functions $\varphi:G\rightarrow \C$ satisfying \eqref{equation de semi sphericite} has finite dimension bounded by the cardinal $k_n$ of $B_T(v,N_n)$. On the other hand, as $\pi$ is irreducible, the vector $\xi$ is cyclic which implies that the linear map $\Psi_n: \Hr{\pi}^{K_n}\rightarrow \mathcal{L}_n: \eta \rightarrow \overline{\varphi_{\xi,\eta}}$ is injective. It follows that $\dim(\Hr{\pi}^{K_n})\leq \dim(\mathcal{L}_n)\leq k_n<+\infty$.
	\end{proof}
	On the other hand, a uniform bound can be obtained for special and cuspidal representations via the following classical result. 
	\begin{theorem}[{\cite[Corollary of Theorem 2]{Harish1970}}]\label{Theorem Harsih chandra les reps de carré integrable sont CCR}
		Let $G$ be a locally compact group, $\pi$ be an irreducible square-integrable representation of $G$ and $K\leq G$ be a compact open subgroup. Then, there exists a positive integer $k_{K,\pi}$ depending on $K$ and $\pi$ such that $\dim(\Hr{\pi}^K)\leq k_{K,\pi}$. 
	\end{theorem}
	We are finally able to prove the announced result.
	\begin{theorem}\label{theorem simple radu groups are unifmroly admissible}
		$G^+_{(i)}(Y_0,Y_1)$ is uniformly admissible and hence CCR.
	\end{theorem}
	\begin{proof}
		For brevity, let $G=G^+_{(i)}(Y_0,Y_1)$. Let $K\leq G$ be a compact open subgroup, let $v\in V$ and for every $n\in \N$ let $K_n=\Fix_G(B_T(v,n))$. For every $n\in \N$, Theorem \ref{Proposition les representations spheriques sont admissibles} ensures the existence of a positive integer $k_{K_n}$ such that $\dim(\Hr{\pi}^{K_n})\leq k_{K_n}$ for every spherical representation $\pi$ of $G$. Now, for every $n\in \N$, we let $\Sigma_n$ be the subset of all equivalence classes of non-spherical irreducible representations of $G$ admitting non-zero $K_n$-invariant vectors. The classification of special and cuspidal representations of $G$ provided by Theorem \ref{thm la classification des speciale} and Theorem \ref{la classification des reresentations cuspidale Radu} ensure that this set $\Sigma_n$ is finite. Furthermore, for every $\sigma\in \Sigma_n$, Theorem \ref{Theorem Harsih chandra les reps de carré integrable sont CCR} ensures the existence of a constant $k_{\sigma, n}$ such that $\dim(\Hr{\pi}^{K_n})\leq k_{\sigma, n}$. In particular, for each $n\in \N$ the constant $k'_n=\max_{\sigma\in \Sigma_n} k_{\sigma, n}$ is finite and $\dim(\Hr{\pi}^{K_n})< k'_n$ for each special and cuspidal representation $\pi$ of $G$. It follows that $\dim(\Hr{\pi}^{K_n})\leq \max\{k_{K_n}, k'_{K_n}\}$ for all $\pi \in \widehat{G}$ and every $n\in \N$. On the other hand, since $(K_n)_{n\in  \N}$ is a basis of neighbourhoods of the identity there exists some $N\in \N$ such that $K_N\subseteq K$. It follows that 
		$$\dim(\Hr{\pi}^K)\leq \max\{k_{K_N}, k'_{K_N}\}\q \forall \pi \in \widehat{G}.$$
	\end{proof}
	The following result provides the last missing piece to deduce from Theorem \ref{theorem simple radu groups are unifmroly admissible} and Theorem \ref{Corollary Radu simple then Radu} that every Radu group on a $(d_0,d_1)$-semi-regular tree with $d_0,d_1\geq 6$ is uniformly admissible. 
	\begin{lemma}\label{lemma G admissible iff  H is admissible}
		Let $G$ be a totally disconnected locally compact group and $H$ be a closed subgroup of index $2$. Then, $G$ is uniformly admissible if and only if $H$ is uniformly admissible.
	\end{lemma}
	\begin{proof}
		Suppose that $G$ is uniformly admissible and let $K$ be a compact open subgroup of $H$. Since $H$ has index $2$ in $G$, it is a clopen subgroup of $G$ which implies that $K$ is a compact open subgroup of $G$. Since $G$ is uniformly admissible, there exists a positive integer $k_K$ such that $\dim(\Hr{\pi}^K)\leq k_K$ for every irreducible representation $\pi$ of $G$. Let $\sigma$ be an irreducible representation of $H$. Theorem \ref{les rep dun group loc compact par rapport a celle d'un de ses sous groupes} ensures that $\Ind_H^G(\sigma)$ is either irreducible or splits as a sum of two irreducible representations of $G$. On the other hand, Frobenius reciprocity (Theorem \ref{criterion of mackey weak frobenius reciprocity}) ensures that $${\rm I}(\sigma,\Res_H^G(\Ind_H^G(\sigma)))={\rm I}(\Ind_H^G(\sigma),\Ind_H^G(\sigma))\geq 1$$ so that $\sigma\leq \Res_H^G(\Ind_H^G(\sigma))$. All together, this proves that $$\dim(\Hr{\sigma}^K)\leq \dim\big(\Hr{\Ind_H^G(\sigma)}^K\big)\leq 2k_K$$ and $H$ is uniformly admissible. Suppose now that $H$ is uniformly admissible and let $K$ be a compact open subgroup of $G$. Since $H$ has index $2$ in $G$, it is a clopen subgroup of $G$ and $K\cap H$ is a compact open subgroup of $H$. Since $H$ is uniformly admissible, this implies the existence of a constant $k_{K\cap H}$ such that $\dim(\Hr{\sigma})\leq k_{{K\cap H}}$ for every irreducible representation $\sigma$ of $H$. Furthermore, Theorem \ref{les rep dun group loc compact par rapport a celle d'un de ses sous groupes} ensures that $\Res_H^G(\pi)$ is either an irreducible representation of $G$ or splits as a direct sum of $2$ irreducible representations of $H$. This implies that $$\dim(\Hr{\pi}^K)\leq \dim(\Hr{\pi}^{K\cap H})= \dim\big(\Hr{\Res_H^G(\pi)}^{K\cap H}\big)\leq 2k_{K\cap H}.$$
		Hence, $G$ is uniformly admissible. 
	\end{proof}

	\chapter{Fell topology of automorphism groups of tree}\label{Chapter Fell topology of Aut(T)}

	\section{Introduction and main results}
	
	Let $G$ be a second-countable locally compact group. One of the most fundamental aspect of representation theory is to determine its unitary dual $\widehat{G}$. We recall from Section \ref{section Fell topology} that this space admits a natural topology called the \tg{Fell topology} whose properties carries the information of whether the group is amenable, compact, Type {\rm I} or even CCR. In addition, this topology is well behaved with respect to natural operations such as restriction and induction.
	Now, let $T$ be a thick semi-regular tree and let $G$ be a closed subgroup $G\leq \Aut(T)$ acting $2$-transitively on the boundary $\partial T$. We recall from Definition \ref{definiton of spherical special cuspidal} that the irreducible representations of $G$ split in three disjoint subsets: spherical, special and cuspidal. We denote by $\sphe(G)$, $\spe(G)$ and $\cusp(G)$ their respective equivalence classes, so that $\widehat{G}=\sphe(G)\sqcup \spe(G)\sqcup\cusp(G)$. We recall that the classification of both spherical and special representations of $G$ can be obtained without further hypothesis and is a classical result. We refer to Sections \ref{Section spherical rep} and \ref{Section special rep} for a detailed exposition. Furthermore, we recall from Section \ref{cuspidal representations of the group of automorphism of a tree} that the cuspidal representations of $G$ have also been classified for certain groups such as the full group of automorphisms of a semi-regular tree. However, to the present days, no result so far describes the Fell topology of $\widehat{G}$ or even its restriction to the known part of the dual. The purpose of this chapter is to provide a description of the Fell topology of $\sphe(G)\sqcup \spe(G)$ for every closed subgroup $G\leq \Aut(T)$ that is $2$-transitive on the boundary based on the classical classification of spherical and special representations of $G$. The proper statement of our results requires some preliminaries. First of all, we recall in light of Proposition \ref{Proposition les representations spheriques sont admissibles}, Theorem \ref{thm la classification des speciale} and Theorem \ref{Theorem Harsih chandra les reps de carré integrable sont CCR} that each of the spherical or special representation of a closed subgroup $G\leq \Aut(T)$ acting $2$-transitively on the boundary is CCR. In particular, Lemma \ref{lemma Dixmier sur les rep CCR sont fermee} ensures that $\{\pi\}$ is closed in $\widehat{G}$ for all $\pi \in \sphe(G)\sqcup \spe(G)$. However, as we will see $ \sphe(G)\sqcup \spe(G)$ always contains pairs of non-Hausdorff points. We recall that two points $x$ and $y$ of a topological space $X$ are \tg{non-Hausdorff} if there does not exist open sets $V_x,V_y\subseteq X$ such that $x\in V_x$, $y\in V_y$ and $V_x\cap V_y=\es$. In particular, in order to describe properly the Fell topology on $\sphe(G)\sqcup \spe(G)$ we now define the interval with double endpoints. Consider the direct sum  $\mathcal{I}_-\sqcup \mathcal{I}_+$ of two copies of the real interval $\lb -1,1\rb$ say 
	$$\mathcal{I}_-=\{(t,-1)\in \R^2:t\in \lb -1,1\rb\}\mbox{ and }\mathcal{I}_+=\{(t,1)\in \R^2:t\in \lb -1,1\rb\}$$
	equipped with the topology induced from $\R^2$. We define an equivalence relation $\sim$ on this space as follows: two elements $(x_1,x_2),(y_1,y_2)\in \mathcal{I}_-\sqcup \mathcal{I}_+$ are equivalent for $\sim$ if $(x_1,x_2)=(y_1,y_2)$ or $x_1\not\in \{-1,1\}$ and $x_1=y_1$. The quotient topological space $\mathcal{I}_-\sqcup\mathcal{I}_+/\sim$ is called the \tg{interval with double endpoints} and is denoted by $\mathcal{I}$. Notice that this topological space has exactly two pairs of non-Hausdorff points $\{(-1,-1),(-1,1)\}$ and $\{(1,-1),(1,1)\}$. This is represented in the following figure by linking the non-Hausdorff points with dotted lines. 
	\begin{figure}[H]\caption{The interval with double endpoints}\label{figure 0}
		\includegraphics[scale=0.15]{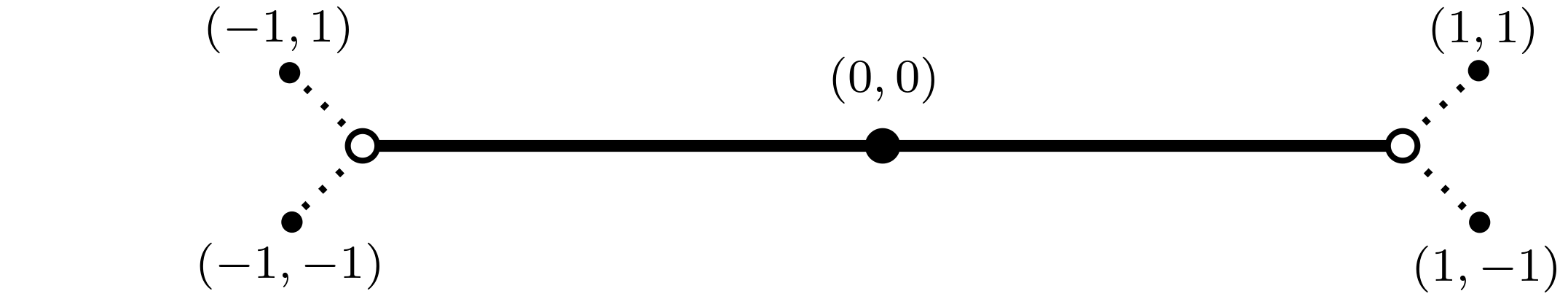}
	\end{figure}
	The following two results describe the Fell topology of $\sphe(G)\sqcup \spe(G)$ depending on the number of $G$-orbits of vertices and uses the same formalism as Theorems \ref{thm la classification des spherique cas trans}, \ref{thm la classification des spherique cas 2 orbites} and \ref{thm la classification des speciale}.
	\begin{theoremletter}\label{theorem la topo de fell de Sphe spe G cas transitif}
		Let $T$ be a thick regular tree,  $v\in V$ and $G\leq \Aut(T)$ be a vertex-transitive closed subgroup acting $2$-transitively on the boundary $\partial T$. Then, the following holds:
		\begin{enumerate}[label=(\roman*)]
			\item $\{\pi\}$ is closed in $\widehat{G}$ for all $\pi \in \sphe(G)\sqcup \spe(G)$ and the sets $\sphe(G)$ and $\sphe(G)\sqcup \spe(G)$ are open in $\widehat{G}$.
			\item $\phi_{v,G}: \sphe(G)\rightarrow \mathcal{I}_v=\lb -1,+1\rb$ is a homeomorphism.
			\item $\sphe(G)\sqcup \spe(G)$ is homeomorphic to the interval with double end with poles of non-Hausdorff points $\{\phi_{v,G}^{-1}(-1),\sigma^{+1}\}$ and $\{\phi_{v,G}^{-1}(1),\sigma^{-1}\}$.
		\end{enumerate}
	\end{theoremletter}
	\noindent The following figure represents the topology given by Theorem \ref{theorem la topo de fell de Sphe spe G cas transitif} with the same conventions as in Figure \ref{figure 0}
	\begin{figure}[H]\caption{$\sphe(G)\sqcup \spe(G)$ for vertex-transitive subgroups}\label{figure 1}
		\includegraphics[scale=0.15]{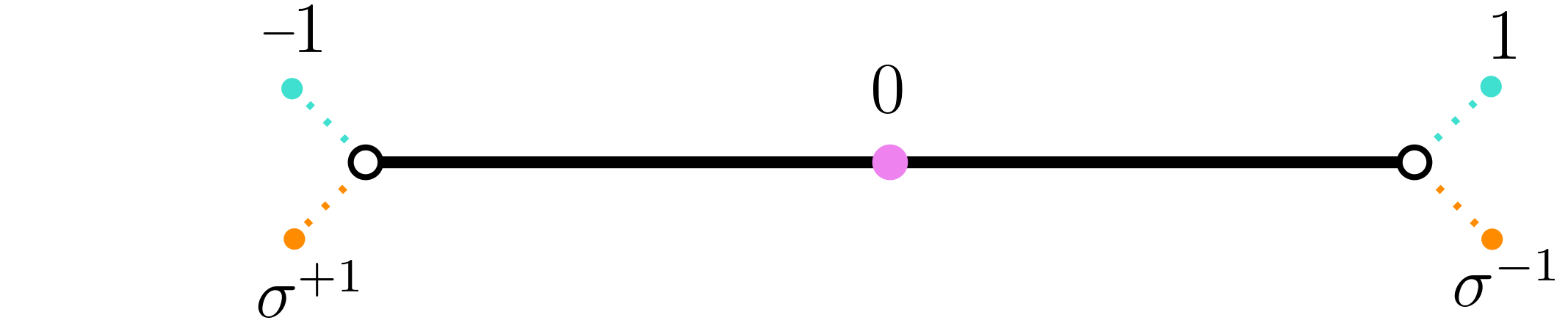}
	\end{figure}
	\begin{theoremletter}\label{theorem la topo de fell de Sphe spe G cas 2orbits}
		Let $T$ be a thick semi-regular tree, $v,v'\in V$ be two adjacent vertices with respective degree $d$ and $d'$ and let $G\leq \Aut(T)^+$ be a closed subgroup of type-preserving automorphisms acting $2$-transitively on the boundary. Then, the following holds:
		\begin{enumerate}[label=(\roman*)]
			\item $\{\pi\}$ is closed in $\widehat{G}$ for all $\pi \in \sphe(G)\sqcup \spe(G)$ and the sets $\sphe(G)$ and $\sphe(G)\sqcup \spe(G)$ are open in $\widehat{G}$.
			\item $\phi_{v,G}:\sphe(G)-\{\pi_{v'}\}\rightarrow  \mathcal{I}_v=\big\lb -\frac{1}{d'-1},1\big\rb$ is a homeomorphism.
			\item $\sphe(G)\sqcup \spe(G)$ is homeomorphic to the interval with double end with poles of non-Hausdorff points $\{\pi_v,\pi_{v'}\}$ and $\{\phi_{v,G}^{-1}(1),\sigma\}$.
		\end{enumerate}
	\end{theoremletter}
	The following figure represents the topology given by Theorem \ref{theorem la topo de fell de Sphe spe G cas 2orbits} with the same conventions as in Figure \ref{figure 0}.
	\begin{figure}[H]\caption{$\sphe(G)\sqcup \spe(G)$ for non-vertex-transitive subgroups}\label{figure 2}
		\includegraphics[scale=0.15]{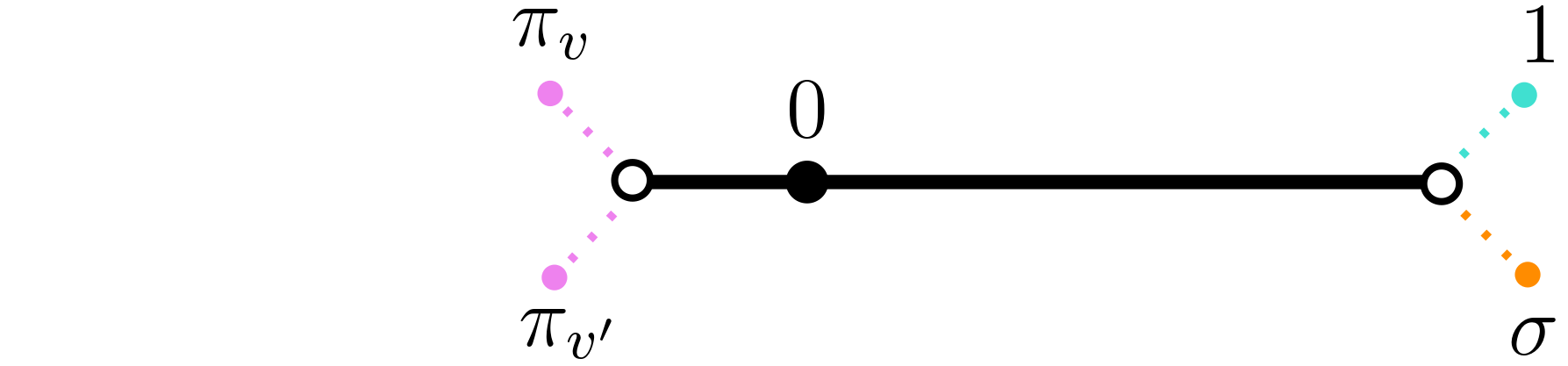}
	\end{figure}
	The proof of these results can be found in Sections \ref{Section preliminaries article Fell topology},  Sections \ref{Section topo of sphe} and \ref{Chapter Sphe et spe topo}. The purpose of Section \ref{Section preliminaries article Fell topology} is to provide the general preliminaries that will be used in Sections \ref{Section topo of sphe} and \ref{Chapter Sphe et spe topo} to describe respectively the Fell topology of $\sphe(G)$ and $\sphe(G)\sqcup \spe(G)$.
	
	Notice that every thick regular tree provides a natural setup to apply both Theorems \ref{theorem la topo de fell de Sphe spe G cas transitif} and \ref{theorem la topo de fell de Sphe spe G cas 2orbits} and to put them in relations in light of the induction-restriction dynamic provided by Theorem \ref{theoremletter induction restiction dynamic}. To be more precise, let $T$ be a thick regular tree and let $G\leq \Aut(T)$ be a closed vertex-transitive subgroup acting $2$-transitively on the boundary $\partial T$. The subgroup $G^+=G\cap \Aut(T)^+$ of type preserving automorphisms of $G$ is also a closed subgroup of $\Aut(T)$ acting $2$-transitively on the boundary $\partial T$ but is not vertex-transitive. In particular, Theorem \ref{theorem la topo de fell de Sphe spe G cas transitif} applies to $G$ while Theorem \ref{theorem la topo de fell de Sphe spe G cas 2orbits} applies to $G^+$. On the other hand, the operations of restriction and induction are continuous for the respective Fell topology and the dynamic of these of operations is described between $G^+$ and $G$ by Theorem \ref{theoremletter induction restiction dynamic}. For clarity of the exposition, we summarised its statement in Figure \ref{figure 4} below. In this figure, the representations corresponding to one another by induction and restriction are matched by two-sided arrows when there is no ambiguity and by colour matching when an ambiguity is possible.  
	\begin{figure}[H]\caption{Induction-restriction dynamic between $G$ and $G^+$}\label{figure 4}
		\includegraphics[scale=0.15]{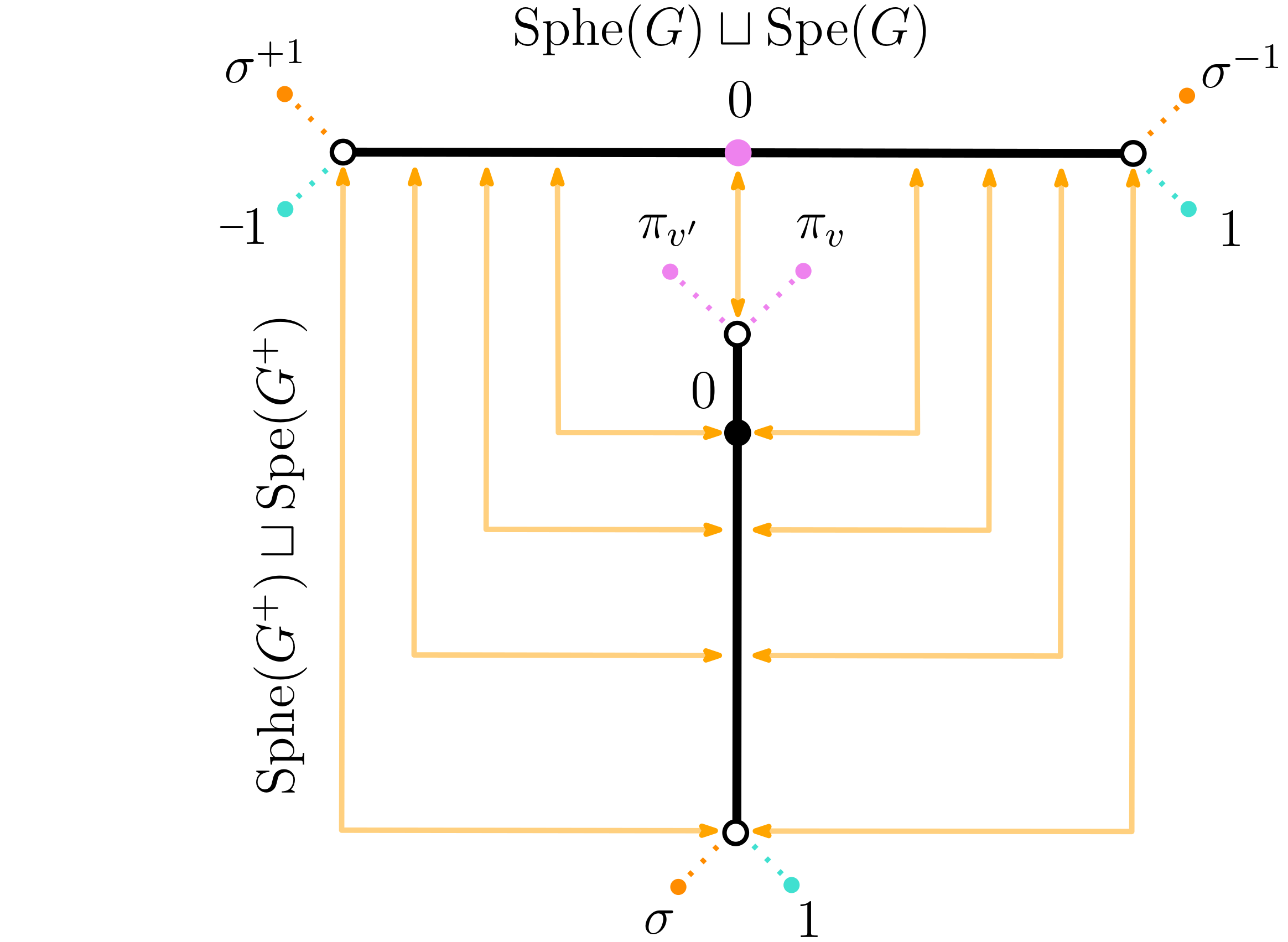}
	\end{figure}
	Notice for instance that the non-Hausdorff aspect of the exceptional spherical representations of $G^+$ could have easily been apprehended from this induction-restriction dynamic and from the description of the Fell topology of $G$. Indeed, the restriction of the representation $\phi_{v,G}^{-1}(0)$ splits as a direct sum 	$\Res_{G^+}^G(\phi_{v,G}^{-1}(0)) \simeq \pi_v \oplus \pi_{v'}$. On the other hand, the restriction $\Res_{G^+}^G(\phi_{v,G}^{-1}(\frac{1}{n}))$ is irreducible for all $n\in \N_0$ and the sequence $(\phi_{v,G}^{-1}(\frac{1}{n+1}))_{n\in \N}$ of spherical representations of $G$ converges to $\phi_{v,G}^{-1}(0)$ in $\widehat{G}$. It follows that the sequence $\big(\Res_{G^+}^G\big(\phi_{v,G}^{-1}(\frac{1}{n+1})\big)\big)_{n\in \N}$ of spherical representations of $G^+$ converges simultaneously to $\pi_v$ and $\pi_{v'}$ in $\widehat{G^+}$.

	Concerning the rest of the dual $\widehat{G}$ we recall that the set $\cusp(G)$ of all the cuspidal representations of closed subgroups $G\leq \Aut(T)$ acting $2$-transitively on the boundary has only been classified for certain families of groups, see Sections \ref{cuspidal representations of the group of automorphism of a tree} and \ref{section cuspidal rep description Radu}. For all of these families of groups, the resulting classification of cuspidal representations ensures that they are both integrable and square-integrable, leading to the conclusion that those groups are CCR and hence Type {\rm I}. In particular, \cite[Theorem 4.4.5]{Dixmier1977} ensures the existence of a dense open locally compact subset of $\widehat{G}$. 
	The following theorem is proved in Section \ref{Chapter cupsidal of radu and otehrs} and provides a description of the Fell topology on the rest of the dual of $\widehat{G}$ under the conjectural condition that the cuspidal representations of $G$ are both integrable and square integrable. In particular, for each such group $G$, we identify explicitly a dense locally compact open subset of $\widehat{G}$ and describes the cortex $\mbox{\rm Cor}(G)$ of $G$ (that is the set of all irreducible representations that are not Hausdorff-separated from the trivial representation $1_{\widehat{G}}$).
	\begin{theoremletter}\label{thm C}
		Let $T$ be a thick semi-regular tree, $G\leq \Aut(T)$ be a closed subgroup acting $2$-transitively on the boundary and suppose that each cuspidal representation of $G$ is both integrable and square-integrable. Then, the following holds:
		\begin{enumerate}[leftmargin=*, label=(\roman*)]
			\item The Fell topology on $\widehat{G}$ is $T_1$. 
			\item $\sphe(G)\sqcup \spe(G)$ is closed in $\widehat{G}$.
			\item $\cusp(G)$ is a countable discrete clopen subset of $\widehat{G}$.
			\item If $G$ is vertex-transitive, $\widehat{G}-\{\sigma^{-1},\sigma^{+1}\}$ is a dense locally compact open subset of $\widehat{G}$ and $\mbox{\rm Cor}(G)=\{{1},\sigma^{-1}\}$.
			\item If $G$ is not vertex-transitive, $\widehat{G}-\{\sigma,\pi_v,\pi_{v'}\}$ is a dense locally compact open subset of $\widehat{G}$ and $\mbox{\rm Cor}(G)=\{{1},\sigma\}$.
		\end{enumerate}
	\end{theoremletter}
	The Fell topology of $\widehat{G}$ described by Theorems \ref{theorem la topo de fell de Sphe spe G cas transitif}, \ref{theorem la topo de fell de Sphe spe G cas 2orbits} and \ref{thm C} is summarised in the Figure \ref{figure 10} below. 
	\begin{figure}[H]\caption{Fell topology of $\widehat{G}$}\label{figure 10}
		\includegraphics[scale=0.15]{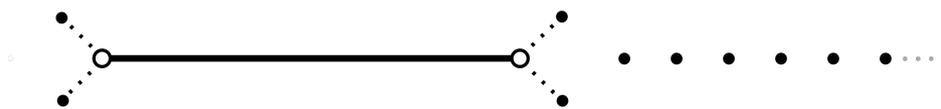}
	\end{figure}
	We notice the strong similarities with the Fell topology of the dual of $\PSl_2(\R)$ provided by the following picture.
	\begin{figure}[H]\caption{Fell topology on $\widehat{\PSl_2(\R)}$}\label{figure 11}
		\includegraphics[scale=0.15]{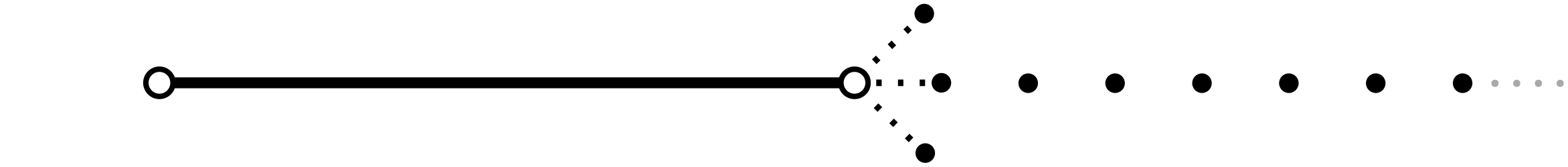}
	\end{figure}
	The description of the Fell topology of the dual of $\Sl_2(\R)$ is initially due to D.Mili\v{c}i\'{c} see \cite{Milicic1971} and an exposition of his result can be found in \cite[Section $7.6$ example $3$]{Folland2016}. Since $\PSl_2(\R)$ is the quotient of $\Sl_2(\R)$ by its centre, the dual space of $\PSl_2(\R)$ corresponds exactly to the irreducible representations of $\Sl_2(\R)$ whose action by $-\id_{\R^2}$ is trivial. In particular, since the Fell topology is described in terms of matrix coefficients, the Figure \ref{figure 11} is obtained from \cite[Figure $7.3$]{Folland2016} by restricting our attention to the trivial representation, the positive part $\pi_{it}^+$ of the spherical principal series, the complementary series $\kappa_s$ and the part of the discrete discrete series $\delta_n^\pm$ indexed by even integers.

	\section{Preliminaries}\label{Section preliminaries article Fell topology}
	
	The purpose of this section is to provide a proof of the first part of both Theorem \ref{theorem la topo de fell de Sphe spe G cas transitif} and Theorem \ref{theorem la topo de fell de Sphe spe G cas 2orbits} as well as some preliminaries required to end their proofs. We start by proving the following result.
	\begin{proposition}\label{proposition sphe and sphe spe are open}
		Let $T$ be a thick semi-regular tree and $G\leq \Aut(T)$ be a closed subgroup acting $2$-transitively on the boundary $\partial T$. Then, $\{\pi\}$ is closed in $\widehat{G}$ for all $\pi \in \sphe(G)\sqcup \spe(G)$ and $\sphe(G)$ and $\sphe(G)\sqcup \spe(G)$ are open subsets of $\widehat{G}$.
	\end{proposition}
	The proof of this result relies on the bijective correspondence between the equivalence classes of unitary representations of $G$ and of non-degenerate $*$-representations of  the maximal $C^*$-algebra $C^*(G)$ recalled in Section \ref{section Fell topology} and on the following result. 
	\begin{lemma}[{\cite[{Proposition $3.3.2$}]{Dixmier1977}}]\label{Lemma dixmier lower semi cont norm de pix}
		Let $A$ be a $C^*$-algebra and let $a\in A$. The function $\Psi_a: \tilde{\pi}\mapsto \norm{\tilde{\pi}(a)}{A}$ is lower semi-continuous in $\widehat{A}$ for the hull-kernel topology.
	\end{lemma}
	\begin{proof}[Proof of Proposition \ref{proposition sphe and sphe spe are open}]
		First of all, we recall in light of Proposition \ref{Proposition les representations spheriques sont admissibles}, Theorem \ref{thm la classification des speciale} and Theorem \ref{Theorem Harsih chandra les reps de carré integrable sont CCR} that each of the spherical or special representation of a closed subgroup $G\leq \Aut(T)$ acting $2$-transitively on the boundary is CCR. In particular, Lemma \ref{lemma Dixmier sur les rep CCR sont fermee} ensures that $\{\pi\}$ is closed in $\widehat{G}$ for all $\pi \in \sphe(G)\sqcup \spe(G)$.
		
		Now, let $K\leq G$ be a compact open subgroup, $\mu$ be the Haar measure of $G$ renormalised so that $\mu(K)=1$, $\mathds{1}_K$ be the characteristic function on $K$, let $e\in E$ be an edge with vertices $v,v'\in V$. Notice that $\mathds{1}_K\in L^1(G)\subseteq  C^*(G)$ and that $\tilde{\pi}(\mathds{1}_K)=\int_K\pi(k)\diff \mu (k)$ is the projection on the space $\hr{\pi}^K$ of $K$-invariant vectors of $\pi$ for every $\pi \in \widehat{G}$. In particular, notice that $ \norm{\tilde{\pi}(\mathds{1}_{K})}{C^*(G)}$ is equal to $1$ if $\pi$ admits a non-zero $K$-invariant vector and to $0$ otherwise. Now, notice that our hypothesis on $G$ ensures that $G$ is edge transitive. In particular, $\cusp(G)$ is exactly the set of irreducible representations of $G$ without non-zero $\Fix_G(e)$-invariant vectors. Therefore, Lemma \ref{Lemma dixmier lower semi cont norm de pix} ensures that $$\cusp(G)= \Psi_{\mathds{1}_{\Fix_G(e)}}^{-1}(0)=\Psi_{\mathds{1}_{\Fix_G(e)}}^{-1}(\lb 0,1/2\rb)$$ is closed in $\widehat{G}$. This proves as desired that $\sphe(G)\sqcup \spe(G)$ is open in $\widehat{G}$. Similarly, since $G$ is transitive the set of edges of $T$, the spherical representations of $G$ are exactly the irreducible representations admitting a non-zero $\Fix_G(v)$ or  $\Fix_G(v')$ invariant vector (possibly both). It follows from Lemma \ref{Lemma dixmier lower semi cont norm de pix} that 
		\begin{equation*}
		\begin{split}
		\spe(G)\sqcup \cusp(G)&=\Psi_{\mathds{1}_{\Fix_G(v)}}^{-1}(0)\cap \Psi_{\mathds{1}_{\Fix_G(v')}}^{-1}(0)\\
		&=\Psi_{\mathds{1}_{\Fix_G(v)}}^{-1}(\lb 0,1/2\rb)\cap \Psi_{\mathds{1}_{\Fix_G(v')}}^{-1}(\lb 0,1/2\rb)
		\end{split}
		\end{equation*} is closed in $\widehat{G}$. This proves as desired that $\sphe(G)$ is open in $\widehat{G}$. 
	\end{proof}
	
	Let $T$ be a thick semi-regular tree, $G\leq \Aut(T)$ be a closed subgroup acting $2$-transitively on the boundary $\partial T$ and let $e$ be an edge of $T$. We recall from the classifications given in Section \ref{Section special rep} that every special representation of $G$ has a one-dimensional subspace of $\Fix_G(e)$-invariant vectors. Proposition \ref{Proposition la dimension de invariant pour Fix Ge est plus petite que 2} below ensures that the space of $\Fix_G(e)$-invariant vectors of any irreducible representation of $G$ is not much larger. This will be the key to prove Theorems \ref{theorem la topo de fell de Sphe spe G cas transitif} an \ref{theorem la topo de fell de Sphe spe G cas 2orbits}. We start with the following two technical results.
	\begin{lemma}\label{lemma les focntions d coeff des spherique pour fix v et fix e sont dans L alpha}
		Let $G\leq \Aut(T)$ be a closed subgroup acting $2$-transitively on the boundary $\partial T$, $e$ be an edge of $T$ with set of vertices $\{v,v'\}$,  $\pi$ be a spherical representation of $G$ admitting a unit $\Fix_G(v)$-invariant vector $\xi\in \Hr{\pi}$, $\eta\in \Hr{\pi}$ be a $\Fix_G(e)$-invariant vector, $\tau$ be a translation of minimal step in $G$ such that $e\subseteq \lb v,\tau v\rb$ and $\varphi_\pi$ be the $\Fix_G(v)$-spherical function of positive type associated with $\pi$. Then, the function
		$$\varphi_{\xi,\eta}: G\rightarrow \C : g\mapsto \prods{\pi(g)\xi}{\eta}_{\Hr{\pi}}$$
		is $\Fix_G(e)$-left invariant , $\Fix_G(v)$-right invariant and satisfies that
		\begin{equation*}
		\int_{\Fix_G(v)}\varphi_{\xi,\eta}(g k \tau)\qq \diff \mu(k)= \varphi_\pi(\tau) \varphi_{\xi,\eta}(g)\qq\forall g\in G.
		\end{equation*}
	\end{lemma}
	\begin{proof}
		Since $\eta\in \Hr{\pi}$ is $\Fix_G(e)$-invariant and $\xi\in \Hr{\pi}$, $\eta\in \Hr{\pi}$ is $\Fix_G(e)$-invariant, the function $\varphi_{\xi,\eta}$ is both $\Fix_G(e)$-left invariant and $\Fix_G(v)$-right invariant. Now, let $\mu$ be the Haar measure of $G$ renormalised so that $\mu(\Fix_G(v))=1$. Since $\hr{\pi}^{\Fix_G(v)}$ is one-dimensional, there exists a function $\beta: G\rightarrow \C$ depending only on $\xi$ and $\eta$ such that
		$$\int_{\Fix_G(v)}\pi(k^{-1})\pi(g^{-1})\eta\diff \mu(k)= \beta(g)\xi.$$
		It follows for all $g,h\in G$ that
		\begin{equation*}
		\begin{split}
		\int_{\Fix_G(v)}\varphi_{\xi,\eta}(gkh)\qq \diff \mu(k)&=\prods{\pi(h)\xi}{\int_{\Fix_G(v)}\pi(k^{-1})\pi(g^{-1})\eta\diff \mu(k)}\\
		&=\prods{\pi(h)\xi}{\beta(g)\xi}=\overline{\beta(g)}\prods{\pi(h)\xi}{\xi}=\overline{\beta(g)} \varphi_{\pi}(h)
		\end{split}
		\end{equation*}
		Taking $h=1_G$ in this equality, we obtain that 
		\begin{equation*}
		\overline{\beta(g)}=\varphi_{\xi,\eta}(g)\qq \forall g\in G.
		\end{equation*}
		Therefore, taking $h=\tau$, we obtain as desired that
		\begin{equation*}
		\int_{\Fix_G(v)}\varphi_{\xi,\eta}(g k \tau)\qq \diff \mu(k)= \varphi_\pi(\tau)\varphi_{\xi,\eta}(g)\qq\forall g\in G.
		\end{equation*}
	\end{proof}
	 The following result ensures that such functions are entirely determined by the values they take on $1_G$ and $\tau$.
	\begin{proposition}\label{prop les fonction de L alpha sont polynomial and totally determined by the value on tau}
		Let $G\leq \Aut(T)$ be a closed subgroup acting $2$-transitively on the boundary $\partial T$, $e$ be an edge of $T$ with set of vertices $\{v,v'\}$, $\tau$ be a translation of minimal step in $G$ such that $e\subseteq \lb v,\tau v\rb$, $\alpha\in \C$ and let $\mathcal{L}_\alpha$ be the space of $\Fix_G(e)$-left, $\Fix_G(v)$-right invariant functions $\varphi: G\rightarrow \C$ satisfying that
		\begin{equation}\label{equation pour trouver les valeur de varphi xi eta pour sperique special}
		\int_{\Fix_G(v)}\varphi(g k \tau)\qq \diff \mu(k)= \alpha \varphi(g)\qq\forall g\in G.
		\end{equation}
		Then, each element $\varphi$ of $\mathcal{L}_\alpha$ is entirely determined by the values it takes on $\{\tau^m: m \in \Z\}$ and for each $m\in \Z$, there exists a polynomial $Q_m$ in three variables depending only on $m$, $G$, $T$ and $v$ such that 
		$$\varphi(\tau^m)=Q_m(\varphi(1_G), \varphi(\tau), \alpha).$$ 
	\end{proposition}
	\begin{proof}
		We recall from Lemma \ref{Lemma the factorization of $G$ into FixGe double cosets} that $G=\bigsqcup_{m\in \Z} \Fix_G(e)\tau^m \Fix_G(v)$ so that every function $\varphi$ of  $\mathcal{L}_\alpha$ is entirely determined by the values it takes on the elements $\tau^m$ with $m\in \Z$. We now prove the existence of $Q_m$ by induction on $m$ for $m\geq 0$ and $m\leq 0$ separately. The existence of $Q_m$ is obvious for $m=0$ and $m=1$. Now, let $m\geq 2$. Taking $h=\tau$ and $g=\tau^{m-1}$ in \eqref{equation pour trouver les valeur de varphi xi eta pour sperique special} we have that
		\begin{equation}\label{equation pour trouver les valeur de varphi xi eta pour sperique special vaprhi evaluer en tau}
		\alpha \varphi(\tau^{m-1})=\int_{\Fix_G(v)}\varphi(\tau^{m-1} k \tau)\qq \diff \mu(k).
		\end{equation}
		When $G$ is vertex-transitive, let $d$ be the degree of $v$, set $$U=\{k\in \Fix_G(v): k\tau v=\tau^{-1} v\}$$ and notice that 
		\begin{equation*}
		\begin{split}
		\alpha\varphi(\tau^{m-1})&=\int_{\Fix_G(v)} \varphi(\tau^{m-1}k\tau)\qq \diff\mu (k)\\
		&=\int_{\Fix_G(v)-U} \varphi(\tau^{m-1}k\tau)\qq \diff\mu (k)+\int_{U} \varphi(\tau^{m-1}k\tau)\qq \diff\mu (k)\\
		&=\frac{d-1}{d}\varphi(\tau^m)+\frac{1}{d}\varphi(\tau^{m-2}).
		\end{split}
		\end{equation*}
		In particular, this proves that
		$$\varphi(\tau^m)= \frac{d}{d-1}\alpha Q_{m-1}(\varphi(1_G),\varphi(\tau),\alpha)-\frac{1}{d-1}Q_{m-2}(\varphi(1_G),\varphi(\tau),\alpha).$$
		When $G$ has two orbits of vertices, let $d$ and $d'$ be the respective degree of $v$ and $v'$, set $U=\{k\in \Fix_G(v):  k\tau v=\tau^{-1}v\}$ and $$U'= \{k\in \Fix_G(v): d(k\tau v,v)\geq d(k\tau v,\tau^{-1} v)\},$$ and notice that
		\begin{equation*}
		\begin{split}
		\alpha\varphi(\tau^{m-1})&=\int_{\Fix_G(v)} \varphi(\tau^{m-1}k\tau)\qq \diff\mu (k)\\
		&=\int_{\Fix_G(v)-U'} \varphi(\tau^{m-1}k\tau)\qq \diff\mu (k)+\int_{U'-U} \varphi(\tau^{m-1}k\tau)\qq \diff\mu (k)\\
		&\q \q \q \q \q \q \q \q \q \q \q \q \q \q \q \q +\int_{U} \varphi(\tau^{m-1}k\tau)\qq \diff\mu (k) \\
		&=\frac{d-1}{d}\varphi(\tau^m)+\frac{d'-2}{d(d'-1)}\varphi(\tau^{m-1})+ \frac{1}{d(d'-1)}\varphi(\tau^{m-2}).
		\end{split}
		\end{equation*}
		In particular, this proves that
		\begin{equation*}
		\begin{split}
		\varphi(\tau^m)&= \bigg(\frac{d}{d-1}\alpha -\frac{d'-2}{(d'-1)(d-1)}\bigg) Q_{m-1}(\varphi(1_G),\varphi(\tau), \alpha)\\
		&\q \q \q \q \q\q \q \q \q -\frac{1}{(d'-1)(d-1)}Q_{m-2}(\varphi(1_G),\varphi(\tau), \alpha).
		\end{split}
		\end{equation*}
		In both cases, this proves, for all $m\geq 0$, that $\varphi(\tau^m)=Q_m(\varphi(1_G),\varphi(\tau), \alpha)$ for some polynomial $Q_m$ depending only on $m$, $G$, $T$ and $v$. Suppose now that $m\leq -1$. Taking $h=\tau$ and $g=\tau^{m+1}$ in \eqref{equation pour trouver les valeur de varphi xi eta pour sperique special} we have that
		\begin{equation}\label{equation pour trouver les valeur de varphi xi eta pour sperique special vaprhi evaluer en tau moins}
		\alpha \varphi(\tau^{m+1})=\int_{\Fix_G(v)}\varphi(\tau^{m+1} k \tau)\qq \diff \mu(k).
		\end{equation}
		For $m=-1$, this implies that
		\begin{equation*}
		\begin{split}
		\alpha\varphi(1_G)&=\int_{\Fix_G(v)} \varphi(k\tau)\qq \diff\mu (k)=\frac{d-1}{d}\varphi(\tau^{-1})+\frac{1}{d}\varphi(\tau).
		\end{split}
		\end{equation*}
		In particular, one has 
		$$\varphi(\tau^{-1})= \frac{d}{d-1}\alpha \varphi(1_G)- \frac{1}{d-1}\varphi(\tau)$$
		and the claim follows. For $m\leq -2$ we have two cases. If $G$ is vertex-transitive, \eqref{equation pour trouver les valeur de varphi xi eta pour sperique special vaprhi evaluer en tau moins} implies that
		\begin{equation*}
		\begin{split}
		\alpha\varphi(\tau^{m+1})&=\int_{\Fix_G(v)} \varphi(\tau^{m+1}k\tau)\qq \diff\mu (k)\\
		&=\int_{\Fix_G(v)-U} \varphi(\tau^{m+1}k\tau)\qq \diff\mu (k)+\int_{U} \varphi(\tau^{m+1}k\tau)\qq \diff\mu (k)\\
		&=\frac{d-1}{d}\varphi(\tau^m)+\frac{1}{d}\varphi(\tau^{m+2})
		\end{split}
		\end{equation*}
		where $v$ has degree $d$ and $U=\{k\in \Fix_G(v): k\tau v=\tau v\}$. In that case, this proves that
		$$\varphi(\tau^m)= \frac{d}{d-1}\alpha Q_{m+1}(\varphi(1_G),\varphi(\tau), \alpha)-\frac{1}{d-1}Q_{m+2}(\varphi(1_G),\varphi(\tau), \alpha).$$
		If on the other hand, $G$ has two orbits of vertices, let $d$ and $d'$ be the respective degrees of $v$ and $v'$, let $U=\{k\in \Fix_G(v):  k\tau v=\tau v\}$ and set $$U'= \{k\in \Fix_G(v): d(k\tau v,v)\geq d(k\tau v,\tau v)\}.$$ A  straightforward computation shows that
		\begin{equation*}
		\begin{split}
		\alpha\varphi(\tau^{m+1})&=\int_{\Fix_G(v)} \varphi(\tau^{m+1}k\tau)\qq \diff\mu (k)\\
		&=\int_{\Fix_G(v)-U'} \varphi(\tau^{m+1}k\tau)\qq \diff\mu (k)+\int_{U'-U} \varphi(\tau^{m+1}k\tau)\qq \diff\mu (k)\\
		&\q \q \q \q \q \q \q \q \q \q \q \q \q \q \q \q +\int_{U} \varphi(\tau^{m+1}k\tau)\qq \diff\mu (k) \\
		&=\frac{d-1}{d}\varphi(\tau^m)+\frac{d'-2}{d(d'-1)}\varphi(\tau^{m+1})+ \frac{1}{d(d'-1)}\varphi(\tau^{m+2})
		\end{split}
		\end{equation*}
		In particular, this proves that
		\begin{equation*}
		\begin{split}
		\varphi(\tau^m)&= \bigg(\frac{d}{d-1}\alpha -\frac{d'-2}{(d'-1)(d-1)}\bigg) Q_{m+1}(\varphi(1_G),\varphi(\tau), \alpha)\\
		&\q \q \q \q\q \q \q \q\q \q -\frac{1}{(d'-1)(d-1)}Q_{m+2}(\varphi(1_G),\varphi(\tau), \alpha).
		\end{split}
		\end{equation*}
	\end{proof}
	The following result follows directly.
	\begin{proposition}\label{Proposition la dimension de invariant pour Fix Ge est plus petite que 2}
		Let $G\leq \Aut(T)$ be a closed subgroup acting $2$-transitively on the boundary $\partial T$, $e$ be an edge of $T$ and $\pi\in \widehat{G}$. Then, $\dim\big(\hr{\pi}^{\Fix_G(e)}\big) \leq 2$.
	\end{proposition}
	\begin{proof}
		The result is trivial if $\pi$ is a cuspidal representation and follows directly from Theorem \ref{thm la classification des speciale} if $\pi$ is a special representation. So, suppose that $\pi$ is a spherical representation of $G$. Let $v$ be a vertex of $e$ for which $\pi$ admits a unit $\Fix_G(v)$-invariant vector $\xi \in \Hr{\pi}$, let $\varphi_\pi$ be the unique $\Fix_G$-spherical function of positive type associated with $\pi$ and let $\tau$ be a translation of minimal step in $G$ such that $e\subseteq \lb v,\tau v\rb$. We recall from Proposition \ref{prop les fonction de L alpha sont polynomial and totally determined by the value on tau} that each element of $\mathcal{L}_{\varphi_\pi(\tau)}$ is entirely determined by the values it takes on $1_G$ and $\tau$. In particular, $\mathcal{L}_{\varphi_\pi(\tau)}$ has dimension at most $2$. On the other hand, for every $\eta \in \Hr{\pi}^{\Fix_G(e)}$ Lemma \ref{lemma les focntions d coeff des spherique pour fix v et fix e sont dans L alpha} ensures that $$\overline{\varphi_{\xi,\eta}}: G\rightarrow \C: g\mapsto \prods{\eta}{\pi(g)\xi}$$ 
		belongs to $\mathcal{L}_{ \varphi_\pi(\tau)}$. This provides a linear map
		$$\hr{\pi}^{\Fix_G(e)}\rightarrow \mathcal{L}: \eta \mapsto\overline{\varphi_{\xi,\eta}}.$$
		Now, notice from the fact that $\pi$ is irreducible that $\xi$ is cyclic and hence that this map is injective. It follows that $\dim(\hr{\pi}^{\Fix_G(e)})\leq 2$. 
	\end{proof}
	\section{The Fell topology on spherical representations}\label{Section topo of sphe}
	Let $T$ be a thick semi-regular tree and let $G\leq \Aut(T)$ be a closed subgroup acting $2$-transitively on the boundary $\partial T$. The purpose of this section is to describe the restriction of the Fell topology on the open subset $\sphe(G)$ of $\widehat{G}$. We adopt the notations of Theorem \ref{thm la classification des spherique cas trans} and Theorem \ref{thm la classification des spherique cas 2 orbites}. We let $v\in V$ be a vertex, $\mathcal{I}_v$ be the corresponding real interval and $\phi_{v,G}$ be the bijective correspondence between the equivalence classes of spherical representations of $G$ admitting a $\Fix_G(v)$-invariant vector and $\mathcal{I}_v$. The following result ensures that this correspondence, is both open and continuous. 
	\begin{theorem}\label{Theorem convergence of spherical representations to spherical representations}
		Let $T$ be a thick semi-regular tree and let $G\leq \Aut(T)$ be a closed subgroup acting $2$-transitively on the boundary $\partial T$. Then, $\phi_{v,G}$ is a homeomorphism.
	\end{theorem}
	Since the domain of $\phi_{v,G}$ coincides with $\sphe(G)$ when $G$ is vertex-transitive and with $\sphe(G)-\{\pi_v\}$ when $G$ has two orbits of vertices, we obtain a complete description of the Fell topology on $\sphe(G)$ when $G$ is vertex-transitive and a description of the Fell topology on $\sphe(G)-\{\pi_v\}$ and $\sphe(G)-\{\pi_{v'}\}$ when $G$ has two orbits of vertices. On the other hand, we recall from  Proposition \ref{proposition sphe and sphe spe are open} that each spherical representation of $G$ is a closed point in $\sphe(G)$. In particular, the above Theorem leads to a complete description of the Fell topology on $\sphe(G)$ even when $G$ has two orbits of vertices (because $\{\pi_v\}$ and $\{\pi_{v'}\}$ are closed).

	The proof of Theorem \ref{Theorem convergence of spherical representations to spherical representations} requires some preliminaries.  Let $G$ be a locally compact group and $\mu$ be a left-invariant Haar measure on $G$. We recall that the convolution $\varphi*\psi$ of two functions $\varphi, \psi : G \rightarrow \C$ is defined (when it makes sense) as 
	\begin{equation*}
	\varphi * \psi (g):=\int_G  \varphi (h) \psi (h^{-1}g)  \diff \mu (h)\q \forall g\in G.
	\end{equation*}
	Furthermore, if $G$ is unimodular we have that
	\begin{equation*}
	\varphi * \psi (g) = \int_G \varphi (h) \psi (h^{-1}g) \qq \diff \mu(h)=\int_G \varphi (gh) \psi (h^{-1}) \qq \diff \mu(h)=\psi *\varphi(g^{-1}).
	\end{equation*}
	In particular, since the groups considered in this chapter are all unimodular and since these equalities simplify our proofs we will restrict our attention to this setup. 

	\begin{lemma}\label{Lemma la convergences sur les comapctes par la convolution}
		Let $G$ be a unimodular locally compact group, $K\subseteq G$ be a compact subset, $\varphi:G\rightarrow \C$ be a continuous function and $(\varphi_n)_{n\in \N}$ be a sequence of continuous functions converging uniformly on compacta to $\varphi$ when $n$ tends to infinity. Then, we have that 		
		\begin{equation*}
		\varphi_n*\mathds{1}_K\li{n}{\infty} \varphi*\mathds{1}_K 	\mbox{	and  }
		\mathds{1}_K*\varphi_n\li{n}{\infty} \mathds{1}_K*\varphi
		\end{equation*} 
		uniformly on compacta.
	\end{lemma}
	\begin{proof}
		Since $\phi * \psi(g)=\psi *\phi(g^{-1})$ for all $\phi,\psi\in C(G)$ and every $g\in G$ and since the inversion map $G\rightarrow G$ is a homeomorphism notice that it is enough to prove that
		\begin{equation*}
		\varphi_n*\mathds{1}_K\li{n}{\infty} \varphi*\mathds{1}_K
		\end{equation*}
		uniformly on compacta. To this end, we notice for all $\phi\in C(G)$ that
		\begin{equation*}
		\begin{split}
		\phi*\mathds{1}_K(g)=\int_G\phi(gh)\mathds{1}_K(h^{-1})\qq \diff \mu(h)= \int_{K^{-1}} \phi(gh)\qq \diff \mu(h)\q\forall g\in G.
		\end{split}
		\end{equation*}
		In particular, we have that
		\begin{equation*}
		\begin{split}
		\modu{\varphi_n*\mathds{1}_K(g)-\varphi *\mathds{1}_K(g)}\leq \int_{K^{-1}}\modu{ \varphi_n(gh)-\varphi(gh)}\diff \mu (h)\q \forall g\in G.
		\end{split}
		\end{equation*}
		For every compact subset $U\subseteq G$, notice that $UK^{-1}$ is compact and therefore that there exists positive integer $N_U\in \N$ such that
		\begin{equation*}
		\modu{ \varphi_n(g)-\varphi(g)}<\varepsilon \q \forall g\in UK^{-1},\qq \forall n\geq N_U.
		\end{equation*}
		In particular, for all $n\geq N_U$ and every $g\in U$ we have that
		\begin{equation*}
		\modu{\varphi_n*\mathds{1}_K(g)-\varphi *\mathds{1}_K(g)}\leq \int_{K^{-1}} \varepsilon\qq \diff \mu (h)=\mu(K^{-1})\varepsilon.
		\end{equation*}
		Since $K$ is compact notice that $\mu(K^{-1})<+\infty$. Hence, $\varphi_n*\mathds{1}_K$ converges to $\varphi*\mathds{1}_K$ uniformly on compacta as $n$ tends to infinity.
	\end{proof}
	\noindent The above lemma becomes particularly relevant in light of the following result. 
	\begin{lemma}\label{Lemma les fonctions de coeff convoluer devien des coeff de proj}
		Let $G$ be a unimodular locally compact group, $\pi$ be a representation of $G$, $\eta\in \hr{\pi}$, $K$ be a compact open subgroup of $G$ and let $P_K:\hr{\pi}\rightarrow \hr{\pi}^K$ be the projection on the space of $K$-invariant vectors. Then, we have that
		$$\mathds{1}_K*\prods{\pi(\cdot)\eta}{\eta}*\mathds{1}_K=\mu(K)^2\prods{\pi(\cdot)P_K\eta}{P_K\eta}.$$
	\end{lemma}
	\begin{proof}
		Since $K$ is a compact open subgroup of $G$, notice that $$P_K=\frac{1}{\mu(K)}\int_K \pi(k)\qq\diff \mu(k).$$  In particular, a  straightforward computation shows that
		\begin{equation*}
		\begin{split}
		\mathds{1}_K*\prods{\pi(\cdot)\eta}{\eta}*\mathds{1}_K(g)&=\int_G \mathds{1}_K(h)(\prods{\pi(\cdot)\eta}{\eta}*\mathds{1}_K)(h^{-1}g)\qq \diff \mu(h)\\
		&=\int_G\int_G \mathds{1}_K(h)\prods{\pi(h^{-1}gh')\eta}{\eta}\mathds{1}_K(h')\qq \diff \mu(h')\diff \mu(h)\\
		&=\int_K\int_K\prods{\pi(g)\pi(k')\eta}{\pi(k)\eta}\qq\diff \mu(k')\diff \mu(k)\\
		&=\mu(K)^2\prods{\pi(g)P_K\eta}{P_K\eta}.
		\end{split}
		\end{equation*}
	\end{proof}
	We are finally able to prove the main result of this section.
	\begin{proof}[Proof of Theorem \ref{Theorem convergence of spherical representations to spherical representations}] We adopt the same notations as in Theorems \ref{thm la classification des spherique cas trans} and \ref{thm la classification des spherique cas 2 orbites}. Furthermore, for brevity, we let $K=\Fix_G(v)$ and for every $\gamma \in \mathcal{I}_v$, we let $\pi_\gamma=\phi_{v,G}^{-1}(\gamma)$, $\varphi_\gamma$ be the unique $K$-spherical function of positive type associated to $\pi_{\gamma}$ and $\xi_\gamma\in \Hr{\pi_{\gamma}}$ be a unit $K$-invariant vector. Since both $\mathcal{I}_v$ and $\widehat{G}$ are second-countable topological space, $\phi_{v,G}$ is a homeomorphism if and only if for every sequence $\alpha_n$ of elements of $\mathcal{I}_v$ and for every $\alpha\in \mathcal{I}_v$, one has that $\alpha_n\li{n}{\infty}\alpha$ if and only if $\pi_{\alpha_n}\li{n}{\infty}\pi_{\alpha}$. 
	
	Suppose first that $\alpha_n$ converges to $\alpha$ and let $\tau\in G$ be a translation of minimal step in $G$ along a geodesic containing $v$. We recall from Lemma \ref{Lemma the factorization of $G$ into FixGe double cosets} that $G=\bigcup_{m\in \N} K \tau^m K$. In particular, since $K \tau^m K$ is open for all $m\in \N$ the sequence $\varphi_{\alpha_n}$ converges uniformly on compacta to $\varphi_\alpha$ if and only if $\varphi_{\alpha_n}(\tau^m)\li{n}{\infty}\varphi_\alpha(\tau^m)$ for all $m\in \N$. On the other hand, Lemma \ref{lemma calcul des fonctions psheriques} ensures the existence of a polynomial $Q_m$ depending only on $m$, $G$, $T$ and $v$ such that $$\varphi_{\gamma}(\tau^m)=Q_m(\gamma).$$ It follows that $\varphi_{\alpha_n}\li{n}{\infty}\varphi_\alpha$ uniformly on compacta which proves that $\pi_{\alpha_n}$ converges to $\pi_{\alpha}$ in the Fell topology as $n$ tends to infinity.

	To prove the other direction, suppose that $\pi_{\alpha_n}$ converges to $\pi_\alpha$ in the Fell topology and let us show that $\alpha_n$ converges to $\alpha$ as $n$ tends to infinity. The convergence in the Fell topology implies that the function 
	$\varphi_{\alpha}: G\rightarrow \C: g\mapsto \prods{\pi_\alpha(g)\xi_\alpha}{\xi_\alpha}$ is uniform on compacta limit of functions of positive type associated with $\pi_{\alpha_n}$. In particular, there exists a sequence of positive integers $l_n\in \N$ and a set of vectors $\zeta_n^k\in \hr{\pi_{\alpha_n}}$ with $k=1,..., l_n$ such that
	$$\s{k=1}{l_n}\prods{\pi_{\alpha_n}(\cdot)\zeta_n^k}{\zeta_n^k}\li{n}{\infty}\varphi_{\alpha}$$
	uniformly on compacta. Now, Lemma \ref{Lemma la convergences sur les comapctes par la convolution} ensures that
	$$\s{k=1}{l_n}\mathds{1}_K*\prods{\pi_{\alpha_n}(\cdot)\zeta_n^k}{\zeta_n^k}*\mathds{1}_K\li{n}{\infty}\mathds{1}_K*\varphi_{\alpha}*\mathds{1}_K$$
	uniformly on compacta and we obtain from Lemma \ref{Lemma les fonctions de coeff convoluer devien des coeff de proj} that
	\begin{equation}\label{equation la convergence dans les psherique pour les coefficient convoluer}
	\s{k=1}{l_n}\prods{\pi_{\alpha_n}(\cdot)P_K\zeta_n^k}{P_K\zeta_n^k}\li{n}{\infty}\varphi_{\alpha}
	\end{equation}
		uniformly on compacta. On the other hand, since $\pi_{\alpha_n}$ is a spherical representation, we have that $\dim(\hr{\pi_{\alpha_n}}^K)=1$. In particular for all $n\in \N$ and every $k=1,...,l_n$  we have that $P_K\zeta_n^k=\lambda^k_n \xi_{\alpha_n}$ for some complex numbers $\lambda_n^k$.Combining this observation with \eqref{equation la convergence dans les psherique pour les coefficient convoluer} we obtain a sequence $\lambda_n$ of positive real numbers such that
		$$\lambda_n\varphi_{\alpha_n}=\lambda_n\prods{\pi_{\alpha_n}(\cdot)\xi_{\alpha_n}}{\xi_{\alpha_n}}\li{n}{\infty} \varphi_{\alpha}$$
		uniformly on compacta. Since the pointwise convergence at $1_G$ implies that $\lambda_n\li{n}{\infty} 1$, it follows that 
		\begin{equation}\label{equation la convergence dans les psherique pour les coefficient ultime}
		\varphi_{\alpha_n}\li{n}{\infty} \varphi_{\alpha}
		\end{equation}
		uniformly on compacta. The pointwise convergence at $\tau$ in \eqref{equation la convergence dans les psherique pour les coefficient ultime} provides the desired convergence
		$$\alpha_n=\varphi_{\alpha_n}(\tau)\li{n}{\infty} \varphi_{\alpha}(\tau)=\alpha.$$
	\end{proof}
	As a direct corollary, the following result shows that if $G$ has two orbits of vertices, the two exceptional representations of $G$ cannot be $T_2$-separated by the Fell topology.
	\begin{lemma}\label{Lemma les non exceptionnelle converge vers les 2 en meme temps}
		Let $T$ be a  thick semi-regular tree and $G\leq \Aut(T)$ be a closed type-preserving subgroup acting $2$-transitively on the boundary $\partial T$. Then, a sequence of non-exceptional spherical representations $\pi_n\in \sphe(G)-\{\pi_{v},\pi_{v'}\}$ converges to $\pi_v$ if and only if it converges to $\pi_{v'}$.
	\end{lemma}
	\begin{proof}
		Theorem \ref{Theorem convergence of spherical representations to spherical representations} ensures that a sequence of non-exceptional spherical representations $\pi_n$ converges to $\pi_{v'}$ if and only if $\phi_{v,G}(\pi_n)\li{n}{\infty} -\frac{1}{d'-1}$ and converges to $\pi_{v}$ if and only if $\phi_{v',G}(\pi_n)\li{n}{\infty} -\frac{1}{d-1}$. On the other hand, we recall from Theorem \ref{thm la classification des spherique cas 2 orbites} that the two bijective correspondences $\phi_{v,G}:\sphe(G)-\{\pi_v\}\rightarrow \big\lbrack -\frac{1}{d'-1},1\big \rbrack$ and $\phi_{v',G}:\sphe(G)-\{\pi_{v'}\}\rightarrow \big\lbrack -\frac{1}{d-1},1\big \rbrack$ satisfy 
		\begin{equation*}
		\phi_{v',G}(\pi)=\frac{d(d'-1)}{d'(d-1)}\phi_{v,G}(\pi) + \frac{d-d'}{d'(d-1)} \q \forall \pi\in \sphe(G)-\{\pi_v,\pi_{v'}\}.\qedhere
		\end{equation*}
		In particular, this equality ensures as desired that $\phi_{v,G}(\alpha_n)$ converges to $-\frac{1}{d'-1}$ if and only if $\phi_{v',G}(\alpha_n)$ converges to $-\frac{1}{d-1}$.
	\end{proof}
	
	\section{The Fell topology of spherical and special representations}\label{Chapter Sphe et spe topo}
	Let $T$ be a thick semi-regular tree and let $G\leq \Aut(T)$ be a closed subgroup acting $2$-transitively on the boundary $\partial T$. The purpose of this section is to describe explicitly the Fell topology on the open set $\sphe(G)\sqcup \spe(G)$ of $\widehat{G}$. Since the Fell topology on $\sphe(G)$ was already described in Section \ref{Section topo of sphe}, we are left to understand the interactions between spherical and special representations of $G$ and the restriction of the Fell topology to $\spe(G)$. However, we recall from Theorem \ref{thm la classification des speciale} that $G$ has either one or two special representations depending on its number of orbits of vertices and these representations are closed points in $\sphe(G)\sqcup \spe(G)$ by Proposition \ref{proposition sphe and sphe spe are open}. This observation provides a complete description of the Fell topology on $\spe(G)$ and leaves us with the task of understanding when a sequence of spherical representations of $G$ converges to a special representation. We answer this question below depending on the number of $G$-orbits of vertices. 
	\subsection{Case 1 : Vertex-transitive subgroups}
	For the rest of this section, let $G\leq \Aut(T)$ be a closed vertex-transitive subgroup acting $2$-transitively on the boundary $\partial T$ and let us use the use the notations of Theorems \ref{thm la classification des spherique cas trans} and \ref{thm la classification des speciale} unless otherwise stated. Let $e$ be an edge of $T$ and  $\{v,v'\}\subseteq V$ be its corresponding set of vertices. We recall that $\mathcal{I}_v=\lbrack -1,1\rbrack$ and for brevity, for every $\gamma \in \mathcal{I}_v$, we let $\pi_\gamma=\phi_{v,G}^{-1}(\gamma)$, $\varphi_\gamma$ be the unique $\Fix_G(v)$-spherical function of positive type associated to $\pi_{\gamma}$ and $\xi_\gamma\in \Hr{\pi_{\gamma}}$ be a unit $\Fix_G(v)$-invariant vector. The purpose of this section is to prove the following theorem.
	\begin{theorem}\label{Theorem convergence des spheriques vers special pour transitive}
		Let $\epsilon\in \{-1,+1\}$ and $(\alpha_n)_{n\in \N}$ be a sequence of elements of $\mathcal{I}_v-\{-1,1\}$. Then, $\pi_{\alpha_n}$ converges to the special representations $\sigma^\epsilon$ in the Fell topology if and only if $\alpha_n$ converges to $-\epsilon$ as $n$ tends to infinity. 
	\end{theorem}
	The proof of this result takes advantage of the fact that the spherical representations of $G$ contain ``few'' $\Fix_G(e)$-invariant vectors. The following result provides an explicit description of a basis of the space of these $\Fix_G(e)$-invariant vectors in terms of the $\Fix_G(v)$-invariant vectors. 
	\begin{lemma}\label{Lemma Basis de H fix e pour les spherique trans}
		Let $\alpha \in \mathcal{I}_v-\{-1,1\}$ and $h\in G$ be an inversion of $e$. Then $$\{\xi_\alpha+\pi_\alpha(h)\xi_\alpha, \qq \xi_\alpha- \pi_\alpha(h)\xi_\alpha\}$$ is an orthogonal basis of $\hr{\pi_\alpha}^{\Fix_G(e)}$. 
	\end{lemma}
	\begin{proof}
		We recall from Proposition \ref{Proposition la dimension de invariant pour Fix Ge est plus petite que 2} that $\dim(\hr{\pi_\alpha}^{\Fix_G(e)}) \leq 2.$
		On the other hand, notice that $\prods{\pi_\alpha(h)\xi_\alpha}{\xi_\alpha}=\alpha \not=\pm1$ so that $\pi_\alpha(h)\xi_\alpha \not=\pm \xi_\alpha$. In particular, the vectors $\xi_\alpha +\pi_\alpha(h)\xi_\alpha$ and $\xi_\alpha -\pi_\alpha(h)\xi_\alpha$ are both non-zero. Furthermore, since $\xi_\alpha$ is $\Fix_G(v)$-invariant, notice that $\pi(h)\xi_\alpha$ is a $\Fix_G(v')$-invariant vector and thus that $\xi_\alpha +\pi_\alpha(h)\xi_\alpha$ and $\xi_\alpha -\pi_\alpha(h)\xi_\alpha$ are both $\Fix_G(e)$-invariant. On the other hand, as $\alpha$ is a real number, a direct computation shows that
		\begin{equation*}
		\begin{split}
		\prods{\xi_{\alpha}+\pi_{\alpha}(h)\xi_{\alpha}}{\xi_{\alpha}-\pi_{\alpha}(h)\xi_{\alpha}}&= \prods{\xi_{\alpha}}{\xi_{\alpha}}  + \prods{\pi_{\alpha}(h)\xi_{\alpha}}{\xi_{\alpha}} - \prods{\xi_{\alpha}}{\pi_{\alpha}(h)\xi_{\alpha}} \\
		&\q\q\q\q\q\q\q\q\q\q\q\q- \prods{\pi_{\alpha}(h)\xi_{\alpha}}{\pi_{\alpha}(h)\xi_{\alpha}}\\
		&= \prods{\xi_{\alpha}}{\xi_{\alpha}}+ \alpha - \overline{\alpha}- \prods{\xi_{\alpha}}{\xi_{\alpha}}=0.
		\end{split}
		\end{equation*}
		This proves that $\xi_\alpha +\pi_\alpha(h)\xi_\alpha$ and $\xi_\alpha -\pi_\alpha(h)\xi_\alpha$ are orthogonal and the result follows.
	\end{proof}
	We now prove the two implications of Theorem \ref{Theorem convergence des spheriques vers special pour transitive} separately.
	\begin{proposition}\label{Proposition convergence spherique special  cas transitif part I}
		If $\pi_{\alpha_n}$ converges to the special representations $\sigma^\epsilon$ in the Fell topology we have that $\alpha_n$ converges to $-\epsilon$ as $n$ tends to infinity. 
	\end{proposition}
	\begin{proof}
		For each  $\alpha\in \mathcal{I}_v-\{-1,+1\}$ and for each $\epsilon\in \{-1,+1\}$, we let $$\eta_{\alpha}^{\epsilon}= \frac{\xi_{\alpha}+ \epsilon\pi_{\alpha}(h)\xi_{\alpha}}{\norm{\xi_{\alpha}+ \epsilon\pi_{\alpha}(h)\xi_{\alpha}}{}}$$
		where $h$ is an inversion of the edge $e$. Lemma \ref{Lemma Basis de H fix e pour les spherique trans} ensures that $\{\eta_{\alpha}^{-1}, \eta_{\alpha}^{+1}\}$ is an orthonormal basis of $\hr{\pi_{\alpha}}^{\Fix_G(e)}$ and a  straightforward computation shows that
		\begin{equation}\label{lequation des prop de eta plus et eta moins des spehriques trans}
		\pi(h) \eta_{\alpha}^{\epsilon}= \epsilon\eta_{\alpha}^{\epsilon}.
		\end{equation}
		Now, let $\eta^\epsilon$ be a unit $\Fix_G(e)$-invariant vector of $\sigma^\epsilon$. We recall in light of Proposition \ref{proposition compute the special functions in both cases} that  $$\varphi_{\eta^\epsilon,\eta^\epsilon}: G\rightarrow \C: g\mapsto \prods{\sigma^{\epsilon}(g)\eta^\epsilon}{\eta^\epsilon}$$ is the unique function $\varphi: G\rightarrow \C$ satisfying simultaneously the three following conditions:
		\begin{enumerate}
			\item $\varphi$ is $\Fix_G(e)$-bi-invariant.
			\item $\varphi(1_G)=1$ and $\varphi(h)=\epsilon$.
			\item $\int_{\Fix_G(w)}\varphi(kg)\diff \mu(k)=0\q \forall g\in G, \qq \forall w\in\{v,v'\}.$
		\end{enumerate}
		In particular, $\varphi_{\eta^\epsilon,\eta^\epsilon}$ does not depends on our choice of $\eta^\epsilon\in \hr{\sigma^\epsilon}^{\Fix_G(e)}$. The convergence in the Fell topology implies that 
		$$\varphi_{\eta^\epsilon,\eta^\epsilon}: G\rightarrow \C: g\mapsto \prods{\sigma^{\epsilon}(g)\eta^\epsilon}{\eta^\epsilon}$$ is the uniform on compacta limit of functions of positive type associated with $\pi_{\alpha_n}$. In particular, there exists a sequence of positive integers $l_n\in \N$ and a set of vectors $\zeta_n^k\in \hr{\pi_{\alpha_n}}$ with $k=1,..., l_n$ such that 
		\begin{equation*}\label{Equation la convergence vers varphi sigma pluss general}
		\s{k=1}{l_n}\prods{\pi_{\alpha_n}(\cdot)\zeta_n^k}{\zeta_n^k}\li{n}{\infty}\varphi_{\eta^\epsilon,\eta^\epsilon}
		\end{equation*}
		uniformly on compacta. Let $B=\Fix_G(e)$ and notice from Lemma \ref{Lemma la convergences sur les comapctes par la convolution} that
		\begin{equation*}
		\s{k=1}{l_n}\mathds{1}_B*\prods{\pi_{\alpha_n}(\cdot)\zeta_n^k}{\zeta_n^k}*\mathds{1}_B\li{n}{\infty}\mathds{1}_B*\varphi_{\eta^\epsilon,\eta^\epsilon}*\mathds{1}_B
		\end{equation*}
		uniformly on compacta. Hence, Lemma \ref{Lemma les fonctions de coeff convoluer devien des coeff de proj} ensures that
		\begin{equation*}
		\s{k=1}{l_n}\prods{\pi_{\alpha_n}(\cdot)P_B\zeta_n^k}{P_B\zeta_n^k}\li{n}{\infty}\varphi_{P_B\eta^\epsilon,P_B\eta^\epsilon}
		\end{equation*}
		uniformly on compacta. Now, Lemma \ref{Lemma Basis de H fix e pour les spherique trans} ensures the existence of two sequences of complex numbers  $\beta_n^-$ and $\beta_n^+$ such that 
		\begin{equation}\label{equation conv unif des speheric sur special trans}
		\prods{\pi_{\alpha_n}(\cdot)(\beta_n^+\eta^{+1}_{\alpha_n}+\beta_n^-\eta_{\alpha_n}^{-1})}{\beta_n^+\eta^{+1}_{\alpha_n}+\beta_n^-\eta_{\alpha_n}^{-1}}\li{n}{\infty}\varphi_{\eta^\epsilon,\eta^\epsilon}
		\end{equation}
		uniformly on compacta. Since $\eta^{-1}_{\alpha_n}$ and $\eta^{+1}_{\alpha_n}$ are orthogonal and as a consequence of \eqref{lequation des prop de eta plus et eta moins des spehriques trans}, the pointwise convergence at $1_G$ and $h$ ensure that
		\begin{equation*}\label{equation 1 trans vers special}
		\modu{\beta_n^+}^2+\modu{\beta^-_n}^2\li{n}{\infty}\varphi_{\eta^\epsilon,\eta^\epsilon}(1_G)=1
		\end{equation*}
		and
		\begin{equation*}\label{equation 2 trans vers special}
		\modu{\beta_n^+}^2-\modu{\beta^-_n}^2\li{n}{\infty}\varphi_{\eta^\epsilon,\eta^\epsilon}(h)=\epsilon
		\end{equation*}
		It follows that:
		\begin{itemize}
			\item $\beta^+_n\li{n}{\infty}0$ and $\modu{\beta^-_n}\li{n}{\infty}1$ if $\epsilon=+1$.
			\item $\beta^-_n\li{n}{\infty}0$ and $\modu{\beta^+_n}\li{n}{\infty}1$ if $\epsilon=-1$.
		\end{itemize}
		In particular, \eqref{equation conv unif des speheric sur special trans} implies that 
		\begin{equation*}
		\prods{\pi_{\alpha_n}(\cdot)\eta^\epsilon_{\alpha_n}}{\eta^\epsilon_{\alpha_n}}\li{n}{\infty} \varphi_{\eta^\epsilon,\eta^\epsilon}
		\end{equation*}
		uniformly on compacta. Now, let $k\in \Fix_G(v)-\Fix_G(v')$. The pointwise convergence at $k$ implies that:
		\begin{equation}\label{Equation de convergence pour la spherique vers sigma Plus finale}
		\begin{split}
		\frac{\prods{\pi_{\alpha_n}(k)(\xi_{\alpha_n}+\epsilon\pi_{\alpha_n}(h)\xi_{\alpha_n})}{(\xi_{\alpha_n}+\epsilon\pi_{\alpha_n}(h)\xi_{\alpha_n})}}{\norm{\xi_{\alpha_n}+\epsilon\pi_{\alpha_n}(h)\xi_{\alpha_n}}{}^2}\li{n}{\infty}\varphi_{\eta^\epsilon,\eta^\epsilon}(k).
		\end{split}
		\end{equation}  
		On the other hand, since $\varphi_{\alpha_n}=\prods{\pi_{\alpha_n}(\cdot)\xi_{\alpha_n}}{\xi_{\alpha_n}}$ is real valued, notice that:
		\begin{equation}\label{equation norm de eta esplilon n trans}
		\begin{split}
		\norm{\xi_{\alpha_n}+\epsilon\pi_{\alpha_n}(h)\xi_{\alpha_n}}{}^2&= 2\prods{\xi_{\alpha_n}}{\xi_{\alpha_n}}+ \epsilon\prods{\pi(h)\xi_{\alpha_n}}{\xi_{\alpha_n}}+\epsilon \prods{\xi_{\alpha_n}}{\pi(h) \xi_{\alpha_n}}\\
		&= 2 + 2\epsilon\varphi_{\pi_{\alpha_n}}(h)= 2 ( 1 +\epsilon \alpha_n).
		\end{split}
		\end{equation}
		Furthermore, we have that
		\begin{equation*}
		\begin{split}
		\mbox{ }\mbox{ }\mbox{ }\mbox{ }\mbox{ }\mbox{ }\mbox{ }\mbox{ }&\mbox{ }\mbox{ }\mbox{ }\mbox{ }\mbox{ }\mbox{ }\mbox{ }\mbox{ }\mbox{ }\mbox{ }\mbox{ }\mbox{ }\prods{\pi_{\alpha_n}(k)(\xi_{\alpha_n}+\epsilon\pi_{\alpha_n}(h)\xi_{\alpha_n})}{(\xi_{\alpha_n}+\epsilon\pi_{\alpha_n}(h)\xi_{\alpha_n})}\\
		&=\prods{\pi_{\alpha_n}(k)\xi_{\alpha_n}}{\xi_{\alpha_n}}+ \epsilon\prods{\pi_{\alpha_n}(h^{-1}k)\xi_{\alpha_n}}{\xi_{\alpha_n}} +\epsilon \prods{\pi_{\alpha_n}(kh)\xi_{\alpha_n}}{\xi_{\alpha_n}}\\
		&\q\q\q\q\q\q\q\q\q\q\q\q\q\q\q\q\q\q\q\q\q\q+ \prods{\pi_{\alpha_n}(h^{-1}kh)\xi_{\alpha_n}}{\xi_{\alpha_n}}\\
		&= \varphi_{\alpha_n}(k)+ \epsilon\varphi_{\alpha_n}(h^{-1}k)+\epsilon \varphi_{\alpha_n}(kh)+ \varphi_{\alpha_n}(h^{-1}kh)\\
		&= \varphi_{\alpha_n}(1_G) + \epsilon\varphi_{\alpha_n}(\tau)+ \epsilon\varphi_{\alpha_n}(\tau)+ \varphi_{\alpha_n}(\tau^2)\\
		&= 1 + 2\epsilon\alpha_n + \frac{d}{d-1}\alpha_n^2-\frac{1}{d-1}= \frac{d-2}{d-1} + 2\epsilon\alpha_n+ \frac{d}{d-1}\alpha_n^2\\
		&= (1+\epsilon\alpha_n)\bigg(\frac{d-2}{d-1}+\frac{d}{d-1}\epsilon\alpha_n\bigg).
		\end{split}
		\end{equation*}
		Finally, since $\varphi_{\eta^\epsilon, \eta^\epsilon}(1_G)=1$ and since $\int_{\Fix_G(v)}\varphi_{\eta^\epsilon, \eta^\epsilon}(t)\diff \mu(t)=0$ notice that $\varphi_{\eta^\epsilon, \eta^\epsilon}(k)= \frac{-1}{d-1}$. In light of these computations, \eqref{Equation de convergence pour la spherique vers sigma Plus finale} becomes 
		\begin{equation*}
		\begin{split}
		\frac{1}{2}\bigg(\frac{d-2}{d-1}+\frac{d}{d-1}\epsilon\alpha_n\bigg)\li{n}{\infty}\frac{-1}{d-1}.
		\end{split}
		\end{equation*}
		This implies as desired that $\alpha_n$ converges to $\-\epsilon$ as $n$ tends to infinity.
	\end{proof}
	We now prove the other implication. 
	\begin{proposition}\label{Proposition si alphan converge cas trans anors on a pi alphan vers special}
		If $\alpha_n$ converges to $-\epsilon$, $\pi_{\alpha_n}$ converges to the special representations $\sigma^\epsilon$ in the Fell topology as $n$ tends to infinity. 
	\end{proposition}
	\begin{proof}
		We adopt the same notations as in the proof of Proposition \ref{Proposition convergence spherique special  cas transitif part I} and show that 
		$$\prods{\pi_{\alpha_n}(\cdot)\eta_{\alpha_n}^\epsilon}{\eta_{\alpha_n}^\epsilon}\li{n}{\infty} \varphi_{\eta^\epsilon,\eta^\epsilon}$$
		uniformly on compacta. We recall from Proposition \ref{proposition compute the special functions in both cases} that  $$\varphi_{\eta^\epsilon,\eta^\epsilon}: G\rightarrow \C: g\mapsto \prods{\sigma^{\epsilon}(g)\eta^\epsilon}{\eta^\epsilon}$$ is the unique function $\varphi: G\rightarrow \C$ satisfying simultaneously the three following conditions:
		\begin{enumerate}
			\item\label{prop 1 equation de merde} $\varphi$ is $\Fix_G(e)$-bi-invariant.
			\item\label{prop 2 equation de merde} $\varphi(1_G)=1$ and $\varphi(h)=\epsilon$.
			\item $\int_{\Fix_G(w)}\varphi(kg)\diff \mu(k)=0\q \forall g\in G, \qq \forall w\in\{v,v'\}.$
		\end{enumerate}
		Set $\psi_n^\epsilon=\prods{\pi_{\alpha_n}(\cdot)\eta_{\alpha_n}^\epsilon}{\eta_{\alpha_n}^\epsilon}$ and let us prove that $\psi_n^\epsilon$ converges uniformly on compacta to a function $\varphi: G\rightarrow \C$ satisfying the three above properties (the desired convergence thus follows). Since $\eta_{\alpha_n}^\epsilon$ is a unit $\Fix_G(e)$-bi-invariant vector, notice that $\psi_n^\epsilon$ is $\Fix_G(e)$-bi-invariant and satisfies $\psi_n^{\epsilon}(1_G)=1$. Now, let $\tau$ be a translation of step $1$ such that $\tau v= v'$ and let $\varphi_{\alpha_n}$ be the $\Fix_G(v)$-spherical function of positive type associated with $\pi_{\alpha_n}$. A  straightforward development of $\eta^\epsilon_{\alpha_n}$ shows for all $g\in G$ that
		$$ \psi_n^\epsilon(g)= \frac{1}{\norm{\xi_{\alpha_n}+ \epsilon\pi_{\alpha_n}(h)\xi_{\alpha_n}}{}^2}(\tilde{\varphi}^\epsilon_n(g)+ \epsilon\tilde{\varphi}^\epsilon_n(gh) )$$
		where $$\tilde{\varphi}^\epsilon_n: G\rightarrow \C : g\mapsto \prods{\pi_{\alpha_n}(g)\xi_{\alpha_n}}{\xi_{\alpha_n}+ \epsilon\pi_{\alpha_n}(h)\xi_{\alpha_n}}_{\Hr{ \pi}}.$$ 
		Now, we recall from \eqref{equation norm de eta esplilon n trans} that
		$$\norm{\xi_{\alpha_n}+ \epsilon\pi_{\alpha_n}(h)\xi_{\alpha_n}}{}^2= 2 (1+\epsilon \alpha_n).$$
		On the other hand, since $\xi_{\alpha_n}$ is a unit $\Fix_G(v)$-invariant vector of $\hr{\pi}$ and since $\xi_{\alpha}+ \epsilon\pi_{\alpha}(h)\xi_{\alpha}$ is a $\Fix_G(e)$-invariant vector of $\Hr{\pi}$, Lemma \ref{lemma les focntions d coeff des spherique pour fix v et fix e sont dans L alpha} and Proposition \ref{prop les fonction de L alpha sont polynomial and totally determined by the value on tau} ensure for all $g\in G$ the existence of a polynomial in three variable $Q_g$ depending only on $G$, $T$, $v$ and $g$ such that 
		$$ \tilde{\varphi}^\epsilon_n(g)=Q_g(\tilde{\varphi}^\epsilon_n(1_G), \tilde{\varphi}^\epsilon_n(\tau),\varphi_{\alpha_n}(\tau)).$$
		Now, notice that $\tilde{\varphi}^\epsilon_n(1_G)$, $\tilde{\varphi}^\epsilon_n(\tau)$ and $\varphi_{\alpha_n}(\tau)$ are all polynomial expressions of $\alpha_n$ depending only on $T$, $G$, $v$ and $\epsilon$. Indeed, a  straightforward computation shows that 
		\begin{equation*}
		\begin{split}
		\tilde{\varphi}^\epsilon_n(1_G)=\prods{\xi_{\alpha_n}}{\xi_{\alpha_n}+ \epsilon\pi_{\alpha_n}(h)\xi_{\alpha_n}}_{\Hr{ \pi}}= 1 + \epsilon \varphi_{\alpha_n}(h^{-1})= 1+ \epsilon\alpha_n,
		\end{split}
		\end{equation*}
		 \begin{equation*}
		 \begin{split}
		 \tilde{\varphi}^\epsilon_n(\tau)=\prods{\pi_{\alpha_n}(\tau)\xi_{\alpha_n}}{\xi_{\alpha_n}+ \epsilon\pi_{\alpha_n}(h)\xi_{\alpha_n}}_{\Hr{ \pi}}= \varphi_{\alpha_n}(\tau) + \epsilon \varphi_{\alpha_n}(h^{-1}\tau)= \alpha_n+\epsilon.
		 \end{split}
		 \end{equation*}
		 and  the definition of $\varphi_{\alpha_n}$ ensures that $\varphi_{\alpha_n}(\tau)=\alpha_n$. In particular, for each $g\in G$, there exists a polynomial in one variable $\tilde{Q}_g^\epsilon$ depending only on $T$, $G$, $v$, $g$ and $\epsilon$ such that 
		 \begin{equation}\label{equation qui montre que psin epslon est une expression polynomial de alphan }
		 	\psi^\epsilon_n(g)= \frac{(Q_g+Q_{gh})(\tilde{\varphi}^\epsilon_n(1_G), \tilde{\varphi}^\epsilon_n(\tau),\varphi_{\alpha_n}(\tau))}{\norm{\xi_{\alpha_n}+ \epsilon\pi_{\alpha_n}(h)\xi_{\alpha_n}}{}^2}= \frac{\tilde{Q}_g^\epsilon(\alpha_n)}{2(1+\alpha_n)}.
		 \end{equation}
		Since $\alpha_n$ converges to $-\epsilon$, notice that the denominator converges to zero as $n$ tend to infinity. However, for all $g\in G$, one has that $$ \modu{\psi^\epsilon_n(g)}=\modu{\prods{\pi_{\alpha_n}(g)\eta_{\alpha_n}^\epsilon}{\eta_{\alpha_n}^\epsilon}}\leq \norm{\eta_{\alpha_n}^\epsilon}{\Hr{\pi}}^2 = 1.$$ 
		This is only possible if $1-x$ divides $\tilde{Q}_g^\epsilon(x)$. It follows from \eqref{equation qui montre que psin epslon est une expression polynomial de alphan } that $\psi_n^\epsilon(g)$ is a polynomial expression of $\alpha_n$ depending only on $T$, $G$, $v$ and $g$. In particular $\psi_n^\epsilon$ converges pointwise to a $\Fix_G(e)$-invariant function $\varphi$ as $n$ tend to infinity. The limit $\varphi$ is clearly $\Fix_G(e)$-bi-invariant. Furthermore, since $\Fix_G(e)$ is both compact and open, $\varphi$ is continuous and the convergence is uniform on compacta. On the other hand, notice form a  straightforward computation that $\varphi(1_G)= \lim_{n\rightarrow \infty}\psi_n^\epsilon(1_G)=1$ and that
		\begin{equation*}
		\psi_n^{\epsilon}(h)= \prods{\pi_{\alpha_n}(h)\eta_{\alpha_n}^\epsilon}{\eta_{\alpha_n}^\epsilon}=\epsilon\prods{\eta_{\alpha_n}^\epsilon}{\eta_{\alpha_n}^\epsilon}=\epsilon.
		\end{equation*}
		This proves that $\varphi$ satisfies the conditions \ref{prop 1 equation de merde} and \ref{prop 2 equation de merde}. To prove that it satisfies the last one, notice from the definition that 
		\begin{equation*}
		\begin{split}
		&\q\q\q\q\q\q\q\q\q\q\q2(1+\epsilon\alpha_n)\int_{\Fix_G(v)}\psi_n^{\epsilon}(kg)\diff \mu(k)\\
		&= \int_{\Fix_G(v)} \varphi_{\alpha_n}(kg)+ \epsilon\varphi_{\alpha_n}(kgh)+ \epsilon\varphi_{\alpha_n}(h^{-1}kg)+ \varphi_{\alpha_n}(h^{-1}kgh)\diff\mu(k)\q \forall g\in G.
		\end{split}
		\end{equation*}
		Let $\mu$ be renormalised in such a way that $\mu(\Fix_G(v))=1$ and remember from the fact that $\varphi_{\alpha_n}$ is a $\Fix_G(v)$-spherical function that $$\int_{\Fix_G(v)}\varphi_{\alpha_n}(g_1kg_2)\diff \mu(k)= \varphi_{\alpha_n}(g_1)\varphi_{\alpha_n}(g_2)\q \forall g_1,g_2\in G.$$ Therefore, we obtain for all $g\in G$ that
		\begin{equation*}
		\begin{split}
		\q\q\q\q&\q\qq  2(1+\epsilon\alpha_n)\int_{\Fix_G(v)}\psi_n^{\epsilon}(kg)\diff \mu(k)\\
		&= \varphi_{\alpha_n}(g)+ \epsilon\varphi_{\alpha_n}(gh)+\epsilon \varphi_{\alpha_n}(h^{-1})\varphi_{\alpha_n}(g)+ \varphi_{\alpha_n}(h^{-1})\varphi_{\alpha_n}(gh)\\
		&= \varphi_{\alpha_n}(g)+ \varphi_{\alpha_n}(gh)+ \alpha_n\varphi_{\alpha_n}(g)+ \alpha_n\varphi_{\alpha_n}(gh)\\
		&= (1+\epsilon\alpha_n)\varphi_{\alpha_n}(g)+ (1+\epsilon\alpha_n)\epsilon\varphi_{\alpha_n}(gh)\\
		&=(1+\epsilon\alpha_n)(\varphi_{\alpha_n}(g)+\epsilon\varphi_{\alpha_n}(gh)).
		\end{split}
		\end{equation*}
		This proves that $\int_{\Fix_G(v)}\psi_n(kg)\diff \mu(k)=\frac{\varphi_{\alpha_n}(g)+\epsilon\varphi_{\alpha_n}(gh)}{2}$. On the other hand, since $\alpha_n$ converges to $-\epsilon$ as $n$ tends to infinity, notice that $$\lim_{n \rightarrow \infty}\varphi_{\alpha_n}(g)=-\epsilon\lim_{n\rightarrow \infty}\varphi_{\alpha_n}(gh).$$ It follows that
		$$\int_{\Fix_G(v)}\psi_n(kg)\diff \mu(k)=\frac{\varphi_{\alpha_n}(g)+\epsilon\varphi_{\alpha_n}(gh)}{2}\li{n}{\infty}0.$$
		Similarly, for every $g\in G$ we have that
		\begin{equation*}
		\begin{split}
		&\q\q\q\q\q\q\q\q  2(1+\epsilon\alpha_n)\int_{\Fix_G(v')}\psi_n^{\epsilon}(kg)\diff \mu(k)\\
		&= \int_{\Fix_G(v')} \varphi_{\alpha_n}(kg)+ \epsilon\varphi_{\alpha_n}(kgh)+ \epsilon\varphi_{\alpha_n}(h^{-1}kg)+ \varphi_{\alpha_n}(h^{-1}kgh)\diff\mu(k).
		\end{split}
		\end{equation*}
		Hence, since $hv=v'$ we obtain that
		\begin{equation*}
		\begin{split}
		&\q \q\q\q \q\q\q \q 2(1+\epsilon\alpha_n) \int_{\Fix_G(v')}\psi_n^{\epsilon}(kg)\diff \mu(k)\\
		&=  \varphi_{\alpha_n}(h)\varphi_{\alpha_n}(h^{-1}g)+ \epsilon\varphi_{\alpha_n}(h) \varphi_{\alpha_n}(h^{-1}gh)+ \epsilon\varphi_{\alpha_n}(h^{-1}g)+ \varphi_{\alpha_n}(h^{-1}gh)\\
		&= \alpha_n\varphi_{\alpha_n}(h^{-1}g)+\alpha_n\epsilon \varphi_{\alpha_n}(h^{-1}gh)+\epsilon\varphi_{\alpha_n}(h^{-1}g)+ \varphi_{\alpha_n}(h^{-1}gh)\\
		&= (1+\epsilon\alpha_n)(\varphi_{\alpha_n}(h^{-1}gh)+\epsilon \varphi_{\alpha_n}(h^{-1}g)).
		\end{split}
		\end{equation*}
		This proves that $\int_{\Fix_G(v')}\psi_n^\epsilon(kg)\diff \mu(k)=\frac{\varphi_{\alpha_n}(h^{-1}gh)+\epsilon\varphi_{\alpha_n}(h^{-1}g)}{2}$. Just as before, since $\alpha_n$ converges to $-\epsilon$ as $n$ tends to  infinity, notice that $$\lim_{n \rightarrow \infty}\varphi_{\alpha_n}(h^{-1}g)=-\epsilon\lim_{n \rightarrow \infty}\varphi_{\alpha_n}(h^{-1}gh)$$ and therefore that
		$$\int_{\Fix_G(v')}\psi_n^{\epsilon}(kg)\diff \mu(k)=\frac{\varphi_{\alpha_n}(h^{-1}gh)+\epsilon \varphi_{\alpha_n}(h^{-1}g)}{2}\li{n}{\infty}0.$$
		The result follows.
	\end{proof}
	
	\subsection{Case 2 : Groups with two orbits of vertices}
	For the rest of this section, let $G\leq \Aut(T)$ be a closed type-preserving subgroup acting $2$-transitively on the boundary $\partial T$ and let us use the notations of Theorems \ref{thm la classification des spherique cas 2 orbites} and \ref{thm la classification des speciale} unless otherwise stated. Let $e$ be an edge of $T$, $\{v,v'\}\subseteq V$ be its corresponding set of vertices and let $d$ and $d'$ be their respective degree in $T$. We recall that $\mathcal{I}_v=\lbrack -\frac{1}{d'-1},1\rbrack$ and for brevity, for every $\gamma\in \mathcal{I}_v$, we let $\pi_\gamma=\phi_{v,G}^{-1}(\gamma)$, $\varphi_\gamma$ be the unique $\Fix_G(v)$-spherical function of positive type associated to $\pi_\gamma$ and $\xi_\gamma\in \hr{\pi_\gamma}$ be a unit $\Fix_G(v)$-invariant vector. The purpose of this section is to prove the following theorem.
	\begin{theorem}\label{Theorem convergence des spheriques vers special pour 2 orb}
		Let $(\alpha_n)_{n\in \N}$ be a sequence of elements of $\mathcal{I}_v-\big\{\frac{-1}{d'-1},1\big\}$. Then, $\pi_{\alpha_n}$ converges to the special representations $\sigma$ in the Fell topology if and only if $\alpha_n$ converges to $1$ as $n$ tends to infinity. 
	\end{theorem}
	Just as in the vertex transitive case, the proof of this result takes advantage of the fact that the spherical representations contain ``few'' $\Fix_G(e)$-invariant. The following result provides a description of a basis of the space of these $\Fix_G(e)$-invariant vectors in terms of the $\Fix_G(v)$-invariant vectors. 
	\begin{lemma}\label{Lemma G has 2 orbits then dimension of fige invariant is 2}
		Let $\pi$ be a non-exceptional non-trivial spherical representation of $G$ and let $\alpha =\phi_{v,G}(\pi)$. Then, the set
		$$\bigg\{\xi_\alpha,\qq  \int_{\Fix_G(v')}\pi(k)\xi_\alpha\diff \mu(k)\bigg\}$$ is a basis of $\hr{\pi}^{\Fix_G(e)}$.
	\end{lemma}
	\begin{proof}
		Notice that $\alpha \in \mathcal{I}_v-\big\{\frac{-1}{d'-1},1\big\}$ since $\pi$ is both non-exceptional and non-trivial. We recall that $\varphi_\alpha= \prods{\pi(\cdot)\xi_\alpha}{\xi_\alpha}$ is a $\Fix_G(v)$-spherical function satisfying $\varphi_\alpha(1_G)=1$ and $\varphi_\alpha(\tau)=\alpha$ for all $ \tau\in G$ such that $d(\tau v,v)=2$. Now, let $\tau$ be such an element of $G$ and renormalise the Haar measure $\mu$ of $G$ so that $\mu (\Fix_G(v))=1$ and notice from a  straightforward computation shows that
		\begin{equation}\label{equation du produit scalaire de eta' avec lui meme}
		\begin{split}
		\prods{\int_{\Fix_G(v')}\pi (k)&\xi_\alpha \diff \mu(k)}{\int_{\Fix_G(v')}\pi (k')\xi_\alpha \diff \mu(k')}\\
		&=\int_{\Fix_G(v')}\int_{\Fix_G(v')}\prods{\pi (k'^{-1}k)\xi_\alpha }{\xi_\alpha }\diff \mu(k)\diff\mu(k')\\
		&=\mu(\Fix_G(v')) \int_{\Fix_G(v')}\prods{\pi (k)\xi_\alpha }{\xi_\alpha }\diff \mu(k)\\
		&=\mu(\Fix_G(v'))\int_{\Fix_G(v')} \varphi_\alpha(k)\diff \mu(k)\\
		&=\frac{d'}{d\qq}\bigg(\int_{\Fix_G(e)} \varphi_\alpha(k)\diff \mu(k)+\int_{\Fix_G(v')-\Fix_G(e)} \varphi_\alpha(k)\diff \mu(k)\bigg)\\
		&=\frac{d'}{d\qq}\Big(\mu(\Fix_G(e)) \varphi(1_G) + (\mu(\Fix_G(v'))-\mu(\Fix_G(e)))\varphi(\tau)\Big)\\
		&= \frac{d'}{d\qq}\bigg(\frac{1}{d}+ \bigg(\frac{d'}{d\mbox{ }}-\frac{1}{d}\bigg)\alpha\bigg)=\frac{d'}{d^2}\big(1+ (d'-1)\alpha\big).
		\end{split}
		\end{equation}
		In particular, since $\alpha\not=\frac{-1}{d'-1}$, notice that $\int_{\Fix_G(v')}\pi(k)\xi_\alpha \diff \mu(k)$ is not zero. Set $$\xi_\alpha'=\frac{\int_{\Fix_G(v')}\pi(k)\xi_\alpha\diff \mu(k)}{\big\lVert\int_{\Fix_G(v')}\pi(k)\xi_\alpha\diff \mu(k)\big\rVert}$$ and notice from another computation that 
		\begin{equation}\label{equation wumpa doubi chevre}
		\begin{split}
		\Big\lVert\int_{\Fix_G(v')}\pi(k)\xi_\alpha\diff \mu(k)\Big\rVert^2 &= \prods{\int_{\Fix_G(v')}\pi (k)\xi_\alpha \diff \mu(k)}{\int_{\Fix_G(v')}\pi (k')\xi_\alpha \diff \mu(k')}\\
		&=\int_{\Fix_G(v')}\int_{\Fix_G(v')}\prods{\xi_\alpha }{\pi (k^{-1}k')\xi_\alpha }\diff\mu(k')\diff \mu(k)\\
		&=\frac{d'}{d\qq }\int_{\Fix_G(v')}\prods{\xi_\alpha }{\pi (k')\xi_\alpha }\diff \mu(k').
		\end{split}
		\end{equation}
		In particular, in light of \eqref{equation du produit scalaire de eta' avec lui meme} notice that 
		\begin{equation}\label{ta mere est une chevre 2000}
		\begin{split}
		\prods{\xi_\alpha }{\xi_\alpha '} &= \frac{d}{d'}\sqrt{\prods{\int_{\Fix_G(v')}\pi (k)\xi_\alpha \diff \mu(k)}{\int_{\Fix_G(v')}\pi (k')\xi_\alpha \diff \mu(k')}}\\
		&\q= \sqrt{\frac{1+ (d'-1)\alpha}{d'}}.
		\end{split}
		\end{equation}
		Thus, we obtain that $\prods{\xi_\alpha }{\xi'_\alpha}=1$ if and only if $\alpha= 1$, and since both vector are renormalised, this proves that $\xi_\alpha $ and $\xi_\alpha '$ are linearly independent for our choice of $\alpha$. The result follows from Proposition \ref{Proposition la dimension de invariant pour Fix Ge est plus petite que 2}. 
	\end{proof}	
	We now prove the two implications of Theorem \ref{Theorem convergence des spheriques vers special pour 2 orb} separately.
	\begin{proposition}\label{Proposition 2 orbites et les spheriques tendent vers la speciale}
		If $\pi_{\alpha_n}$ converges to the special representations $\sigma$ in the Fell topology we have that $\alpha_n$ converges to $1$ as $n$ tends to infinity. 
	\end{proposition}
	\begin{proof}
		For each $\alpha\in \mathcal{I}_v-\big\{\frac{-1}{d'-1},1\big\}$, let $\xi_{\alpha}\in \hr{\pi_{\alpha}}$ be a unit $\Fix_G(v)$-invariant vector of $\pi_{\alpha}$ and set $$\xi_{\alpha}':=\frac{\int_{\Fix_G(v')} \pi_{\alpha}(k)\xi_{\alpha} \diff \mu(k)}{\norm{\int_{\Fix_G(v')} \pi_{\alpha}(k)\xi_{\alpha_n} \diff \mu(k)}{}}$$ and recall from Lemma \ref{Lemma G has 2 orbits then dimension of fige invariant is 2} that $\{\xi_{\alpha}, \xi'_{\alpha}\}$ is a basis of $\hr{\pi_{\alpha}}^{\Fix_G(e)}$. Now, let $\eta$ be a unit $\Fix_G(e)$-invariant vector of $\hr{\sigma}$. We recall from Proposition \ref{proposition compute the special functions in both cases} that the function of positive type $\varphi_{\eta,\eta}=\prods{\sigma(\cdot)\eta}{\eta}$ is the unique function $\varphi: G\rightarrow \C$ satisfying the three following properties:
			\begin{enumerate}
				\item $\varphi$ is $\Fix_G(e)$-bi-invariant.
				\item $\varphi(1_G)=1$.
				\item $\int_{\Fix_G(w)}\varphi(kg)\diff \mu(k)=0\q \forall g\in G, \qq \forall w\in\{v,v'\}.$
			\end{enumerate}
		On the other hand, the convergence in the Fell topology implies that
		$$\varphi_{\eta,\eta}: G\rightarrow \C: g\mapsto \prods{\sigma(g)\eta}{\eta}$$ is uniform on compacta limit of functions of positive type associated with $\pi_{\alpha_n}$. In particular, there exists a sequence of positive integer $l_n\in \N$ and a set of vectors $\zeta_n^k\in \hr{\pi_{\alpha_n}}$ with $k=1,..., l_n$ such that:
		\begin{equation}\label{Equation 2 orb la convergence vers varphi sigma moins general}
		\s{k=1}{l_n}\prods{\pi_{\alpha_n}(\cdot)\zeta_n^k}{\zeta_n^k}\li{n}{\infty}\varphi_{\eta,\eta}
		\end{equation}
		uniformly on compacta. Let $B=\Fix_G(e)$ and notice from Lemma \ref{Lemma la convergences sur les comapctes par la convolution} that
		\begin{equation*}
		\s{k=1}{l_n}\mathds{1}_B*\prods{\pi_{\alpha_n}(\cdot)\zeta_n^k}{\zeta_n^k}*\mathds{1}_B\li{n}{\infty}\mathds{1}_B*\varphi_{\eta,\eta}*\mathds{1}_B
		\end{equation*}
		uniformly on compacta. Since $\xi_{\alpha_n}$ and $\xi_{\alpha_n}'$ form a basis of $\hr{\pi_{\alpha_n}}^{\Fix_G(e)}$ the above convergence ensures, in light of Lemma \ref{Lemma les fonctions de coeff convoluer devien des coeff de proj}, the existence of two sequence $\beta_n$, $\beta_n'\in \C$ such that
		\begin{equation}\label{equation 2 orbites convergence apres ka convolution par l'indicatrice de Fix Ge}
		\prods{\pi_{\alpha_n}(\cdot)\beta_n\xi_{n}+\beta_n'\xi_{n}'}{\beta_n\xi_{n}+\beta_n'\xi_{n}'}\li{n}{\infty}\varphi_{\eta,\eta}
		\end{equation}
		uniformly on compacta. The pointwise convergence at $1_G$ implies that
		\begin{equation}\label{Equation 2 orb equation A}
		\modu{\beta_n}^2+2\mbox{\rm Re}(\beta_n\overline{\beta_n'}\prods{\xi_n}{\xi_{n}'})+\modu{\beta_n'}^2\li{n}{\infty}1.
		\end{equation}
		On the other hand, convoluting both side of \eqref{equation 2 orbites convergence apres ka convolution par l'indicatrice de Fix Ge} by $\mathds{1}_{K}$ with $K=\Fix_G(v)$, Lemmas \ref{Lemma la convergences sur les comapctes par la convolution} and \ref{Lemma les fonctions de coeff convoluer devien des coeff de proj} ensure that
		\begin{equation}\label{Equation 2 orb equation B premise}
		\prods{\pi_{\alpha_n}(\cdot)\beta_n\xi_{\alpha_n}+\beta_n'P_{K}\xi_{\alpha_n}'}{\beta_n\xi_{\alpha_n}+\beta_n'P_{K}\xi_{\alpha_n}'}\li{n}{\infty}0
		\end{equation}
		uniformly on compacta. Since $P_{K}=\int_{\Fix_G(v)}\pi_{\alpha_n}(k)\diff \mu(k)$ we have that
		\begin{equation*}
		\prods{\xi_{\alpha_n}}{P_{K}\xi_{\alpha_n}'}=\prods{P_{K}\xi_{\alpha_n}}{\xi'_{\alpha_n}}=\prods{\xi_{\alpha_n}}{\xi_{\alpha_n}'}.
		\end{equation*}
		On the other hand, since $\xi_{\alpha_n}'$ is a unit $\Fix_G(v')$-invariant vector, notice that
		\begin{equation*}
		\prods{P_{K}\xi_{\alpha_n}'}{\xi'_{\alpha_n}}=\int_{K} \prods{\pi_{\alpha_n}(k)\xi'_{\alpha_n}}{\xi'_{\alpha_n}}=\int_{K}\varphi_{\xi'_{\alpha_n},\xi'_{\alpha_n}}(k)\diff \mu(k)
		\end{equation*}
		where $\varphi_{\xi'_{\alpha_n},\xi'_{\alpha_n}}$ is the $\Fix_G(v')$-spherical function associated with $\pi_{\alpha_n}$. In particular, we obtain from a direct computation that
		\begin{equation*}
		\begin{split}
		\prods{P_{K}\xi'_{\alpha_n}}{\xi'_{\alpha_n}}&=\int_{\Fix_G(e)} \varphi_{\xi'_{\alpha_n},\xi'_{\alpha_n}}(k)\diff \mu(k)+ \int_{\Fix_G(v)-\Fix_G(e)}\varphi_{\xi'_{\alpha_n},\xi'_{\alpha_n}}(k)\diff \mu(k)\\
		&=\frac{1}{d}\varphi_{\xi'_{\alpha_n},\xi'_{\alpha_n}}(1_G)+ \frac{d-1}{d}\varphi_{\xi'_{\alpha_n},\xi'_{\alpha_n}}(\tau_{v'})\\
		&= \frac{1}{d}+ \frac{d-1}{d} \bigg(\frac{d(d'-1)}{d'(d-1)}\alpha_n+\frac{d-d'}{d'(d-1)}\bigg)=\frac{1+(d'-1)\alpha_n}{d'}.
		\end{split}
		\end{equation*}
		In light of these computations, the pointwise convergence of \eqref{Equation 2 orb equation B premise} at $1_G$ becomes
		\begin{equation}\label{Equation 2 orb equation B}
		\modu{\beta_n}^2+ 2\mbox{\rm Re}(\beta_n\beta_n'\prods{\xi_{\alpha_n}}{\xi_{\alpha_n}'})+\frac{1+(d'-1)\alpha_n}{d'} \modu{\beta_n'}^2 \li{n}{\infty}0.
		\end{equation}
		Together, \eqref{Equation 2 orb equation A} and \eqref{Equation 2 orb equation B} ensure that $\beta_n'$ does not converge to $0$ and that
		\begin{equation*}
		\bigg(\frac{1+(d'-1)\alpha_n}{d'}-1\bigg) \modu{\beta_n'}^2\li{n}{\infty}0.
		\end{equation*}
		This proves as desired  that $\alpha_n\li{n}{\infty} 1$. 
	\end{proof}
	We now prove the other implication.
	\begin{proposition}
		If $\alpha_n$ converges to $1$, $\pi_{\alpha_n}$ converges to the special representation $\sigma$ as $n$ tends to infinity. 
	\end{proposition}
	\begin{proof}
		We adopt the same notations as in the proof of Proposition \ref{Proposition 2 orbites et les spheriques tendent vers la speciale}. For brevity, let $K=\Fix_G(v)$, $K'=\Fix_G(v')$, $\eta_{{\alpha_n}}= \frac{\xi_{\alpha_n}- \xi_{\alpha_n}'}{\norm{\xi_{\alpha_n}- \xi_{\alpha_n}'}{}}$ and $\psi_n=\prods{\pi_{\alpha_n}(\cdot)\eta_{\alpha_n}}{\eta_{\alpha_n}}$. To prove the convergence in the Fell topology, we show that 
		$$\psi_n\li{n}{\infty} \varphi_{\eta,\eta}$$
		uniformly on compacta. We recall from Proposition \ref{proposition compute the special functions in both cases} that $\varphi_{\eta,\eta}$ is the unique function $\varphi:G\rightarrow\C$ satisfying the three following properties:
		\begin{enumerate}
			\item $\varphi$ is $\Fix_G(e)$-bi-invariant,
			\item $\varphi(1_G)=1$,
			\item $\int_{\Fix_G(w)}\varphi(kg)\diff \mu(k)=0\q \forall g\in G, \qq \forall w\in \{v,v'\}.$
		\end{enumerate}
		To prove the desired convergence, we show that $\psi_n$ converges uniformly on compacta to a function $\varphi: G\rightarrow \C$ satisfying the three above conditions. Now, notice for all $g\in G$ that 
		\begin{equation*}
		\begin{split}
		\psi_n(g)&= \frac{1}{\norm{\xi_{\alpha_n}- \xi_{\alpha_n}'}{}^2} (\tilde{\varphi}_n(g)-\tilde{\varphi}_n'(g) )
		\end{split}
		\end{equation*}
		where
		$$\tilde{\varphi}_n: G\rightarrow \C : g\rightarrow \prods{\pi_{\alpha_n}(g)\xi_{\alpha_n}}{\xi_{\alpha_n}-\xi_{\alpha_n}'}_{\hr{\pi_{\alpha_n}}}$$
		and 
		$$\tilde{\varphi}'_n: G\rightarrow \C : g\rightarrow \prods{\pi_{\alpha_n}(g)\xi_{\alpha_n}'}{\xi_{\alpha_n}-\xi_{\alpha_n}'}_{\hr{\pi_{\alpha_n}}}.$$
		Observe that $\xi_{\alpha_n}$ is a $\Fix_G(e)$-invariant vector of $\Hr{\pi_\alpha}$ and let $\tau\in G$ be translation of step $2$ such that $e\subseteq \lb v,\tau v\rb$ and $e\subseteq \lb v,\tau v\rb$ and hence such that $e\subseteq \lb v', \tau^{-1}v'\rb$. 
		Since $\xi_{\alpha_n}$ is a unit $\Fix_G(v)$-invariant vector and since $\xi_{\alpha_n}'$ is a $\Fix_G(v')$-invariant vector,  Lemma \ref{lemma les focntions d coeff des spherique pour fix v et fix e sont dans L alpha} and Proposition \ref{prop les fonction de L alpha sont polynomial and totally determined by the value on tau} ensure the existence of polynomials in two variables  $Q_g$ and $Q'_g$ depending only on $T$, $G$, $v$, $v'$ and $g$ such that 
		\begin{equation}\label{prout de chevre le prequel 1}
		\tilde{\varphi}_n(g)= Q_g(\tilde{\varphi}_n(1_G),\tilde{\varphi}_n(\tau), \alpha_n)
		\end{equation}
		and 
		\begin{equation}\label{prout de chevre le prequel 2}
		\tilde{\varphi}_n'(g)= Q_g'(\tilde{\varphi}_n'(1_G),\tilde{\varphi}_n'(\tau^{-1}), \prods{\pi_{\alpha_n}(\tau)\xi_{\alpha_n}'}{\xi_{\alpha_n}'}).
		\end{equation}
		Now, let $\gamma_n= \prods{\xi_{\alpha_n}}{\xi_{\alpha_n}'}_{\Hr{\alpha_n}}$ and let us show that $\alpha_n$, $\tilde{\varphi}_n(1_G)$, $\tilde{\varphi}_n(\tau)$, $\tilde{\varphi}_n'(1_G)$ and $\tilde{\varphi}_n'(\tau^{-1})$ are all polynomial expressions of $\gamma_n$ depending only on $T$, $G$, $v$ and $v'$. Indeed, we recall from \eqref{ta mere est une chevre 2000} that $\gamma_n$ is a positive real number and that
		\begin{equation}\label{Prout de chevre 1}
		\begin{split}
			\alpha_n = \frac{d'\gamma_n^2-1}{d'-1}.
		\end{split}
		\end{equation}
		Furthermore, some  straightforward computations show that
		\begin{equation}\label{Prout de chevre 2}
		\begin{split}
		\tilde{\varphi}_n(1_G)= \prods{\xi_{\alpha_n}}{\xi_{\alpha_n}-\xi_{\alpha_n}'}= \prods{\xi_{\alpha_n}}{\xi_{\alpha_n}}-\prods{\xi_{\alpha_n}}{\xi_{\alpha_n}'}= 1- \gamma_n
		\end{split}
		\end{equation}
		and that
		\begin{equation}\label{Prout de chevre 3}
		\begin{split}
		\tilde{\varphi}_n'(1_G)=  \prods{\xi_{\alpha_n}'}{\xi_{\alpha_n}-\xi_{\alpha_n}'}= \prods{\xi_{\alpha_n}'}{\xi_{\alpha_n}}-\prods{\xi_{\alpha_n}'}{\xi_{\alpha_n}'}= \overline{\gamma_n}-1= \gamma_n-1.
		\end{split}
		\end{equation}
		Now, since $d(v,v')=d(\tau v,v')$, notice the existence of some $k'\in \Fix_G(v')$ such that $k'\tau v=v$. In particular, since $\xi_{\alpha_n}'$ is $\Fix_G(v')$-invariant and since $k'\tau\in \Fix_G(v)$ this implies that
		\begin{equation}\label{ouf le wumpa lumpa 1}
		\begin{split}
		\prods{\pi_{\alpha_n}(\tau)\xi_{\alpha_n}}{\xi_{\alpha_n}'}=\prods{\pi_{\alpha_n}(k'\tau)\xi_{\alpha_n}}{\xi_{\alpha_n}'}=\prods{\xi_{\alpha_n}}{\xi_{\alpha_n}'}=\gamma_n
		\end{split}
		\end{equation}
		It follows that 
		\begin{equation}\label{Prout de chevre 4}
		\begin{split}
		\tilde{\varphi}_n(\tau)= \prods{\pi_{\alpha_n}(\tau)\xi_{\alpha_n}}{\xi_{\alpha_n}-\xi_{\alpha_n}'}&= \prods{\pi_{\alpha_n}(\tau)\xi_{\alpha_n}}{\xi_{\alpha_n}}- \prods{\pi_{\alpha_n}(\tau)\xi_{\alpha_n}}{\xi_{\alpha_n}'}\\
		&\q\q\q= \alpha_n - \gamma_n=  \frac{d'\gamma_n^2-1}{d'-1}- \gamma_n.
		\end{split}
		\end{equation}
		Finally, for our last expression we have that
		\begin{equation}\label{Prout de chevre 5 BB}
		\begin{split}
		\tilde{\varphi}_n'(\tau^{-1})=\prods{\pi_{\alpha_n}(\tau^{-1})\xi_{\alpha_n}'}{\xi_{\alpha_n}-\xi_{\alpha_n}'}= \prods{\pi_{\alpha_n}(\tau^{-1})\xi_{\alpha_n}'}{\xi_{\alpha_n}}-\prods{\pi_{\alpha_n}(\tau^{-1})\xi_{\alpha_n}'}{\xi_{\alpha_n}'}.
		\end{split}
		\end{equation}
		Now, notice in light of \eqref{ouf le wumpa lumpa 1} that 
		$$\prods{\pi_{\alpha_n}(\tau^{-1})\xi_{\alpha_n}'}{\xi_{\alpha_n}}= \overline{\prods{\pi_{\alpha_n}(\tau)\xi_{\alpha_n}}{\xi_{\alpha_n}'}}=\overline{\gamma_n}=\gamma_n $$
		On the other hand, since $d(v',\tau^{-1}v')= d(v', \tau v')$, there exists an element $k'\in \Fix_G(v')$ such that $k'\tau^{-1}v'=\tau v'$. In particular, since $\xi_{\alpha_n}'$ is $\Fix_G(v')$-invariant and since $k'\tau\in \Fix_G(v)$ this implies that
		$$\prods{\pi_{\alpha_n}(\tau^{-1})\xi_{\alpha_n}'}{\xi_{\alpha_n}'}=\prods{\pi_{\alpha_n}(k'\tau^{-1})\xi_{\alpha_n}'}{\xi_{\alpha_n}'}=\prods{\pi_{\alpha_n}(\tau)\xi_{\alpha_n}'}{\xi_{\alpha_n}'}$$
		Now, let $$U'=\{k'\in \Fix_G(v'): k'\tau v= v\}$$ and renormalise the Haar measure $\mu$ of $G$ so that $\mu(\Fix_G(v'))=1$. We observe from the definition of $\xi_{\alpha_n}'$, from \eqref{equation wumpa doubi chevre} and from Lemma \ref{lemma calcul des fonctions psheriques} that
		\begin{equation}\label{Prout de chevre 5}
		\begin{split}
		\prods{\pi_{\alpha_n}&(\tau)\xi_{\alpha_n}'}{\xi_{\alpha_n}'}\\
		&=\bigg(\frac{d}{d'\gamma_n}\bigg)^2\int_{\Fix_G(v')}\int_{\Fix_G(v')}\varphi_{\alpha_n}(k'\tau k)\diff \mu(k)\diff \mu(k')\\
		&=\bigg(\frac{d}{d'\gamma_n}\bigg)^2\int_{\Fix_G(v')}\bigg(\int_{\Fix_G(e)}\varphi_{\alpha_n}(k'\tau k)\diff \mu(k)\\
		&\q\q\q\q\q\q\q\q\q\q\q\q+ \int_{\Fix_G(v')-\Fix_G(e)}\varphi_{\alpha_n}(k'\tau k)\diff \mu(k)\bigg)\diff \mu(k')\\
		&=\bigg(\frac{d}{d'\gamma_n}\bigg)^2\bigg(\frac{1}{d'}\int_{\Fix_G(v')}\varphi_{\alpha_n}(k'\tau)\diff \mu(k')+\frac{d'-1}{d'}\int_{\Fix_G(v')}\varphi_{\alpha_n}(k'\tau^2)\diff \mu(k')\bigg)\\
		&=\bigg(\frac{d}{d'\gamma_n}\bigg)^2 \bigg(\frac{1}{d'}\int_{U'}\varphi_{\alpha_n}(k'\tau)\diff \mu(k')+\frac{1}{d'}\int_{\Fix_G(v')-U'}\varphi_{\alpha_n}(k'\tau)\diff \mu(k')\\
		&\qq\q \q \q + \frac{d'-1}{d'}\int_{U'}\varphi_{\alpha_n}(k'\tau^2)\diff \mu(k')+  \frac{d'-1}{d'}\int_{\Fix_G(v')-U'}\varphi_{\alpha_n}(k'\tau^2)\diff \mu(k')\bigg)\\
		&= \bigg(\frac{d}{d'\gamma_n}\bigg)^2\bigg(\frac{1}{d'^2}\varphi_{\alpha_n}(1_G)+ \frac{d'-1}{d'^2}\varphi_{\alpha_n}(\tau) +  \frac{d'-1}{d'^2}\varphi_{\alpha_n}(\tau) +  \frac{(d'-1)^2}{d'^2} \varphi_{\alpha_n}(\tau^2) \bigg)\\
		&=   \bigg(\frac{d}{d'^2\gamma_n}\bigg)^2\bigg(1+ 2(d'-1)\alpha_n \\
		&\q \q \q \q \q \q\q + (d'-1)^2\bigg(\frac{d}{d-1}\alpha_n^2 -\frac{d'-2}{(d'-1)(d-1)}\alpha_n-\frac{1}{(d'-1)(d-1)}\bigg)\bigg)\\
		&= \bigg(\frac{d}{d'^2\gamma_n}\bigg)^2 \frac{(1+(d'-1)\alpha_n)((d'-1)d\alpha_n+(d-d'))}{(d-1)}\\
		&= \bigg(\frac{d}{d'^2\gamma_n}\bigg)^2 \frac{(d'\gamma_n)^2(d\gamma_n^2-1)}{d-1} = \frac{d^2(d\gamma_n^2-1)}{d'^2(d-1)}
		\end{split}
		\end{equation}
		It follows from \eqref{Prout de chevre 5 BB} that 
		\begin{equation}\label{Prout de chevre 6}
		\begin{split}
		\tilde{\varphi}_n'(\tau^{-1})= \gamma_n- \frac{d^2(d\gamma_n^2-1)}{d'^2(d-1)}.
		\end{split}
		\end{equation}
		Together, \eqref{Prout de chevre 1}, \eqref{Prout de chevre 2}, \eqref{Prout de chevre 3}, \eqref{Prout de chevre 4}, \eqref{Prout de chevre 5} and \eqref{Prout de chevre 6} prove that  $\alpha_n$, $\tilde{\varphi}_n(1_G)$, $\tilde{\varphi}_n(\tau)$, $\tilde{\varphi}_n'(1_G)$ and $\tilde{\varphi}_n'(\tau^{-1})$  are all polynomial expressions of $\gamma_n$ depending only on $T$, $G$, $v$ and $v'$. In particular, for each $g\in G$ it follows from \eqref{prout de chevre le prequel 1} and \eqref{prout de chevre le prequel 2} that there exists a polynomial in one variable $\tilde{Q}_g$ depending only on $T$, $G$, $v$, $v'$ and $g$ such that 
		\begin{equation}\label{equation qui montre que psin epslon est une expression polynomial de alphan PROUT }
		\psi_n(g)=\frac{\tilde{Q}_g(\gamma_n)}{\norm{\xi_{\alpha_n}-\xi_{\alpha_n}'}{}^2}= \frac{\tilde{Q}_g(\gamma_n)}{2(1-\gamma_n)}.
		\end{equation}
		However, since $\alpha_n$ converges to $1$, so does $\gamma_n$. In particular, the denominator converges to zero as $n$ tend to infinity. However, for all $g\in G$, one has that $$ \modu{\psi_n(g)}=\modu{\prods{\pi_{\alpha_n}(g)\eta_{\alpha_n}}{\eta_{\alpha_n}}}\leq \norm{\eta_{\alpha_n}}{\Hr{\pi}}^2 = 1.$$ 
		This is only possible if $1-x$ divides $\tilde{Q}_g(x)$. It follows from \eqref{equation qui montre que psin epslon est une expression polynomial de alphan PROUT } that $\psi_n(g)$ is a polynomial expression of $\gamma_n$ depending only on $T$, $G$, $v$, $v'$ and $g$. In particular $\psi_n$ converges as desired pointwise to a $\Fix_G(e)$-invariant function $\varphi$ as $n$ tend to infinity. On the other hand, since $\eta_{\alpha_n}$ is a unit $\Fix_G(e)$-invariant vector, notice that $\psi_n$ is $\Fix_G(e)$-invariant and satisfies $\psi_n(1_G)=1$. Hence, the limit function $\varphi$ is $\Fix_G(e)$-invariant and satisfies $\varphi(1_G)=1$. Now, consider the function $$\phi_n: G\rightarrow \C : g\mapsto \prods{\pi_{\alpha_n}(g)\xi_{\alpha_n}}{\xi'_{\alpha_n}}.$$ Since $\dim(\hr{\pi_{\alpha_n}}^{\Fix_G(v)})=1$ we notice from a  straightforward computation that 
		\begin{equation*}
		\begin{split}
		&\q\qq 2(1-\gamma_n)\int_K \psi_n(kg) \diff \mu(k)\\
		&= \int_{\Fix_G(v)} \prods{\pi_{\alpha_n}(kg)\xi_{\alpha_n}}{\xi_{\alpha_n}}\diff \mu(k) -\int_{\Fix_G(v)} \prods{\pi_{\alpha_n}(kg)\xi_{\alpha_n}}{\xi_{\alpha_n}'}\diff \mu(k)\\
		&\q\q\q\q\qq-\int_{\Fix_G(v)} \prods{\pi_{\alpha_n}(kg)\xi_{\alpha_n}'}{\xi_{\alpha_n}}\diff \mu(k) +\int_{\Fix_G(v)} \prods{\pi_{\alpha_n}(kg)\xi_{\alpha_n}'}{\xi_{\alpha_n}'}\diff \mu(k)\\
		&= \prods{\pi_{\alpha_n}(g)\xi_{\alpha_n}}{\xi_{\alpha_n}}-\prods{\pi_{\alpha_n}(g)\xi_{\alpha_n}}{P_{K}\xi_{\alpha_n}'}- \prods{\pi_{\alpha_n}(g)\xi_{\alpha_n}'}{\xi_{\alpha_n}}+\prods{\pi_{\alpha_n}(g)\xi_{\alpha_n}'}{P_{K}\xi_{\alpha_n}'}\\
		&= \varphi_{\alpha_n}(g)- \prods{\xi_{\alpha_n}}{\xi_{\alpha_n}'}\varphi_{\alpha_n}(g) -\overline{\phi_n(g^{-1})} +\prods{\xi_{\alpha_n}}{\xi_{\alpha_n}'}\overline{\phi_n(g^{-1})}\\
		&= (1-\gamma_n)\big(\varphi_{\alpha_n}(g) -\overline{\phi_n(g^{-1})}\big).
		\end{split}
		\end{equation*} 
		In particular, we have  that
		$$\int_{\Fix_G(v)} \psi_n(kg) \diff \mu(k)= \frac{\varphi_{\alpha_n}(g) -\overline{\phi_n(g^{-1})}}{2}.$$
		On the other hand, as $\alpha_n\li{n}{\infty}1$ notice that 
		$$\norm{\xi_{\alpha_n}-\xi_{\alpha_n}'}{}^2=2(1-\gamma_n)\li{n}{\infty}0.$$
		This implies that 
		$$\varphi_{\alpha_n}(g)\li{n}{\infty}1\mbox{  and  }\phi_{n}(g)\li{n}{\infty}1 \q \forall g\in G.$$ 
		This proves, as desired, that
		$$\int_{\Fix_G(v)} \psi_n(kg) \diff \mu(k)\li{n}{\infty}0.$$
		By symmetry of the role played by $\xi_{\alpha_n}$ and $\xi_{\alpha_n}'$, we obtain similarly that
		$$\int_{{\Fix_G(v')}} \psi_n(k'g) \diff \mu(k')\li{n}{\infty}0$$
		and the result follows.
	\end{proof}

	\section{A proof of Theorem \ref{thm C}}\label{Chapter cupsidal of radu and otehrs}
	
	Let $T$ be a thick semi-regular tree. Up to this point, we have described the restriction of the Fell topology to both $\sphe(G)$ and $\sphe(G)\sqcup \spe(G)$ for any closed subgroup $G\leq \Aut(T)$ acting $2$-transitively on the boundary $\partial T$ see Theorem \ref{theorem la topo de fell de Sphe spe G cas transitif} and Theorem \ref{theorem la topo de fell de Sphe spe G cas 2orbits}. However, the rest of the dual $\widehat{G}$ has not yet been classified at this level of generality. On the other hand, the cuspidal representations have been classified for certain families of closed subgroups $G\leq \Aut(T)$ acting $2$-transitively on the boundary $\partial T$ see Sections \ref{cuspidal representations of the group of automorphism of a tree} and \ref{section cuspidal rep description Radu}. For all of these families, the resulting classification ensures that each cuspidal representation is both integrable and square-integrable, leading to the conclusion that these groups are CCR and hence Type {\rm I} see Section \ref{section radu groups are Type I}. Furthermore, \cite[Theorem 4.4.5]{Dixmier1977} ensures the existence of a dense open locally compact subset of $\widehat{G}$ for any Type {\rm I} group $G$. 
	The purpose of this section is to prove Theorem \ref{thm C} which provides a description of the Fell topology of the entire dual $\widehat{G}$, identifies explicitly a dense locally compact open subset of $\widehat{G}$ and describes the cortex of any closed subgroup $G\leq \Aut(T)$ acting $2$-transitively on the boundary $\partial T$ under the conjectural condition that the cuspidal representations of $G$ are both integrable and square integrable. We recall that the cortex $\mbox{\rm Cor}(G)$ of a locally compact group $G$ is the set of irreducible representations of $G$ that are not Hausdorff-separated from the trivial representation $1_{\widehat{G}}$. We recall that this theorem can be stated as follows. 
	\begin{theorem*}
		Let $T$ be a thick semi-regular tree, $G\leq \Aut(T)$ be a closed subgroup acting $2$-transitively on the boundary and suppose that each cuspidal representation of $G$ is both integrable and square-integrable. Then, the following result holds:
		\begin{enumerate}[leftmargin=*, label=(\roman*)]
			\item\label{item theorem letter fell topo 1} The Fell topology on $\widehat{G}$ is $T_1$. 
			\item\label{item theorem letter fell topo 2} $\sphe(G)\sqcup \spe(G)$ is closed in $\widehat{G}$.
			\item\label{item theorem letter fell topo 3} $\cusp(G)$ is a countable discrete clopen subset of $\widehat{G}$.
			\item\label{item theorem letter fell topo 4} If $G$ is vertex-transitive, $\widehat{G}-\{\sigma^{-1},\sigma^{+1}\}$ is a dense locally compact open subset of $\widehat{G}$ and $\mbox{\rm Cor}(G)=\{{1},\sigma^{-1}\}$.
			\item\label{item theorem letter fell topo 5} If $G$ is not vertex-transitive, $\widehat{G}-\{\sigma,\pi_v,\pi_{v'}\}$ is a dense locally compact open subset of $\widehat{G}$ and $\mbox{\rm Cor}(G)=\{{1},\sigma\}$.
		\end{enumerate}
	\end{theorem*}
	\begin{proof}
		To prove \ref{item theorem letter fell topo 1}, we recall from Theorem \ref{Theorem les differentes implication topo fell} that the Fell topology on $\widehat{G}$ is $T_1$ if and only if the $G$ is CCR. On the other hand, a totally disconnected locally compact group is CCR if and only if every of its irreducible representation $\pi$ of $G$ is admissible in the sense that for every compact open subgroup $K \leq G$ the space $\Hr{\pi}^K$ of $K$-invariant vectors is finite dimensional \cite{Nebbia1999}. Since  $G$ is a closed subgroup acting $2$-transitively on the boundary $\partial T$, Proposition \ref{Proposition les representations spheriques sont admissibles} ensures that each spherical representation of $G$ is admissible. On the other hand, our hypothesis and the classification of special representations ensures that every other irreducible representation of $G$ is square-integrable. In particular, Theorem \ref{Theorem Harsih chandra les reps de carré integrable sont CCR} ensures that these irreducible representations are also admissible. It follows that $G$ is CCR. Now, to prove \ref{item theorem letter fell topo 2} and \ref{item theorem letter fell topo 3} we recall that \cite[Corollary 1 pg.223 ]{DufloMoore1976} ensures that the equivalence class of an irreducible representation that is both integrable and square-integrable is open in the unitary dual $\widehat{G}$ for the Fell topology. It follows directly from our hypotheses that $\cusp(G)$ is open in $\widehat{G}$ and hence that $\sphe(G)\sqcup \spe(G)$ is closed. On the other hand, Lemma \ref{lemma dixmier sur les reps de carre integrable} ensures that each irreducible square-integrable representation of $G$ appears as subrepresentation of the regular representation. Since $G$ is separable, our hypotheses implies that $G$ has at most countably many cuspidal representation. It follows that $\cusp(G)$ is a countable discrete open subset of $\widehat{G}$. On the other hand, $\cusp(G)$ is closed in $\widehat{G}$ by Proposition \ref{proposition sphe and sphe spe are open}. This proves \ref{item theorem letter fell topo 3}. Now, \ref{item theorem letter fell topo 4} and \ref{item theorem letter fell topo 5} follow from the above discussion and from Theorem \ref{theorem la topo de fell de Sphe spe G cas transitif} and Theorem \ref{theorem la topo de fell de Sphe spe G cas 2orbits}.
	\end{proof}

\clearpage
\newpage
\bibliographystyle{alpha}
\bibliography{bibliography}	
\addcontentsline{toc}{chapter}{Bibliography}
\end{document}